\documentclass[oneside,11pt,a4paper]{article}
\usepackage{amsmath}
\usepackage{amssymb}
\usepackage{amsthm}
\usepackage[alphabetic]{amsrefs}
\usepackage{array}
\usepackage{fancyhdr}
\usepackage{euscript}
\usepackage{graphics,graphicx}
\usepackage{cancel}
\usepackage{fancybox}
\usepackage{verbatim}
\usepackage{tikz}
\usepackage{tikz-cd}
\usetikzlibrary{arrows.meta,positioning}
\usetikzlibrary{arrows.meta,patterns,calc,decorations.pathreplacing,3d}
\definecolor{bubbleC}{HTML}{2166AC}   \definecolor{domC}{HTML}{E08214}      \definecolor{cuspC}{HTML}{B2182B}     \definecolor{covC}{HTML}{1B7837}      \definecolor{holeC}{HTML}{762A83}     \tikzset{
  bnd/.style={gray!70!black,line width=0.9pt},          bub/.style={bubbleC,line width=0.8pt},
  bubfill/.style={fill=bubbleC!12},
  dom/.style={fill=domC!35,draw=domC!80!black,line width=0.8pt},
  lat/.style={fill=black},
  cusp/.style={fill=cuspC},
  cov/.style={covC,line width=0.8pt,dashed},
  cellE/.style={black!55,line width=0.3pt},
  lab/.style={font=\footnotesize},
dnode/.style={circle,draw,fill=white,minimum size=7pt,inner sep=0pt,line width=0.7pt},
  dedge/.style={line width=0.7pt},
  dlab/.style={font=\scriptsize,inner sep=2pt},
  dm/.style={font=\scriptsize,above,inner sep=2pt},
}
\newcommand{\dstep}{1.3cm}  \newcommand{\dyndouble}[2]{\draw[dedge] ([yshift=1.6pt]#1.east)--([yshift=1.6pt]#2.west);
  \draw[dedge] ([yshift=-1.6pt]#1.east)--([yshift=-1.6pt]#2.west);
  \draw[dedge] ($(#1)!0.5!(#2)+(-2.5pt,4pt)$)--($(#1)!0.5!(#2)+(2.5pt,0pt)$)--($(#1)!0.5!(#2)+(-2.5pt,-4pt)$);}
\newcommand{\dotat}[2][lat]{\fill[#1] (#2) circle (1.6pt);}

\usepackage{xcolor}
\usepackage{longtable}
\usepackage{graphicx}
\usepackage[tableposition=below]{caption}
\usepackage[colorlinks=false,hyperindex]{hyperref}

\hypersetup{
	colorlinks,
	citecolor=violet,
	linkcolor=violet,
	urlcolor=blue,
	pdftitle={The Basic Theory of Clifford-Bianchi Groups for Hyperbolic n-Space},
	pdfauthor={Taylor Dupuy, Anton Hilado, Colin Ingalls, Adam Logan}}

\newcommand{\agentskill}[1]{\href{https://github.com/tdupu/agent-skills-public/blob/8d13efb5e66306aeae6e3b948bc463e7e2db1b65/skills/#1/SKILL.md}{\texttt{#1}}}

\DeclareFontFamily{OT1}{rsfs}{}
\DeclareFontShape{OT1}{rsfs}{n}{it}{<-> rsfs10}{}
\DeclareMathAlphabet{\mathscr}{OT1}{rsfs}{n}{it}

\newtheoremstyle {custom}{}{}{}{}{}{.}{ }{\thmname{}\thmnumber{}\thmnote{\bfseries #3}}

\newtheoremstyle {Theorem}{}{}{\itshape}{}{}{.}{ }{\thmname{\bfseries #1}\thmnumber{\;\bfseries #2}\thmnote{\;(\bfseries #3)}}

\theoremstyle{Theorem}
\newtheorem{theorem}{Theorem}[subsection]

\newtheorem{corollary}[theorem]{Corollary}
\newtheorem{lemma}[theorem]{Lemma}
\newtheorem{proposition}[theorem]{Proposition}

\theoremstyle{definition}
\newtheorem{definition}[theorem]{Definition}
\newtheorem{notation}[theorem]{Notation}
\newtheorem{example}[theorem]{Example}
\newtheorem{remark}[theorem]{Remark}
\newtheorem{algorithm}[theorem]{Algorithm}

\newtheorem{conjecture}[theorem]{Conjecture}
\newtheorem{question}[theorem]{Question}

\newcommand{\quat}[2]{\left( \frac{#1}{#2}\right) }
\newcommand{\frakp}{{\mathfrak p}}
\newcommand{\frakP}{{\mathfrak P}}
\newcommand{\GL}{\operatorname{GL}}
\newcommand{\SL}{\operatorname{SL}}
\newcommand{\PSL}{\operatorname{PSL}}

\newcommand{\code}{C}
\newcommand{\ZZ}{\mathbb{Z}}
\let\Z=\ZZ
\newcommand{\QQ}{\mathbb{Q}}
\let\Q=\QQ
\newcommand{\lat}{{\Lambda}}
\newcommand{\NN}{\mathbb{N}}

\newcommand{\RR}{\mathbb{R}}
\newcommand{\CC}{\mathbb{C}}
\newcommand{\HH}{\mathbb{H}}
\newcommand{\FF}{\mathbb{F}}

\newcommand{\mon}{\vartriangleright} \newcommand{\Hcal}{\mathcal{H}}

\renewcommand{\t}{\operatorname{t}}

\newcommand{\Aut}{\operatorname{Aut}}
\newcommand{\Br}{\operatorname{Br}}

\newcommand{\Clf}{\operatorname{Clf}}
\newcommand{\adr}{\operatorname{adr}}
\renewcommand{\Vec}{\operatorname{Vec}}

\newcommand{\Ring}{\mathsf{Ring}}

\newcommand{\Res}{\operatorname{Res}}

\newcommand{\sh}{\operatorname{sh}}
\newcommand{\fppf}{\operatorname{fppf}}
\newcommand{\et}{\operatorname{\acute{e}t}}

\newcommand{\Pic}{\operatorname{Pic}}
\newcommand{\Disc}{\operatorname{Disc}}
\newcommand{\Spec}{\operatorname{Spec}}

\newcommand{\Pin}{\operatorname{Pin}}
\newcommand{\GSpin}{\operatorname{GSpin}}
\newcommand{\Spin}{\operatorname{Spin}}
\newcommand{\PSpin}{\operatorname{PSpin}}
\newcommand{\SO}{\operatorname{SO}}
\renewcommand{\O}{\operatorname{O}}

\newcommand{\SV}{\operatorname{SV}}
\newcommand{\GG}{\mathbb{G}}

\newcommand{\conj}{\operatorname{conj}}
\newcommand{\orb}{\operatorname{orb}}

\newcommand{\height}{\operatorname{ht}}
\newcommand{\tr}{\operatorname{tr}}
\newcommand{\Tr}{\operatorname{Tr}}
\newcommand{\Isom}{\operatorname{Isom}} 
\newcommand{\SU}{\operatorname{SU}}

\newcommand{\PSU}{\operatorname{PSU}}
\newcommand{\GM}{\operatorname{GM}}
\newcommand{\Stab}{\operatorname{Stab}}
\newcommand{\Cl}{\operatorname{Cl}}

\newcommand{\PP}{\mathbb{P}}

\newcommand{\Ocal}{\mathcal{O}}

\renewcommand{\Im}{\operatorname{Im}}
\renewcommand{\Re}{\operatorname{Re}}

\newcommand{\id}{\operatorname{id}}

\newcommand{\KO}{\operatorname{KO}}
\newcommand{\BO}{\operatorname{BO}}
\newcommand{\KSp}{\operatorname{KSp}}
\newcommand{\Mod}{\mathsf{Mod}}
\newcommand{\Hom}{\operatorname{Hom}}
\newcommand{\KU}{\operatorname{KU}}

\newcommand{\BSp}{\operatorname{BSp}}
\newcommand{\U}{\operatorname{U}}
\newcommand{\Sp}{\operatorname{Sp}}
\newcommand{\Hilb}{\operatorname{Hilb}}
\newcommand{\wt}{\operatorname{wt}}
\newcommand{\diag}{\operatorname{diag}}
\newcommand{\Weyl}{\operatorname{Weyl}}
\newcommand{\Dynkin}{\operatorname{Dynkin}}

\newcommand{\quo}{\mathsf{quo}}
\newcommand{\rem}{\mathsf{rem}}
\renewcommand{\L}{\mathsf{L}}
\newcommand{\R}{\mathsf{R}}
\renewcommand{\gcd}{\mathsf{gcd}}
\newcommand{\gcdcoeffs}{\mathsf{gcdcoeffs}}

\newcommand{\nrd}{\operatorname{nrd}}
\newcommand{\Nrd}{\operatorname{Nrd}}
\newcommand{\trd}{\operatorname{trd}}
\newcommand{\ns}{\operatorname{Ns}}
\newcommand{\n}{\operatorname{N}}
\renewcommand{\tr}{\operatorname{Tr}}
\newcommand{\herm}{\operatorname{herm}}
\newcommand{\pos}{\operatorname{pos}}
\newcommand{\BS}{\operatorname{BS}}

\newcommand{\Sat}{\operatorname{Sat}}

\newcommand{\Xcal}{\mathcal{X}}
\newcommand{\Ycal}{\mathcal{Y}}

\newenvironment{psmallmatrix}
{\left(\begin{smallmatrix}}
	{\end{smallmatrix}\right)}

\newcommand{\Orb}{\mathsf{Orb}}
\newcommand{\Man}{\mathsf{Man}}

\newcommand{\Cusps}{\operatorname{Cusps}}
\newcommand{\GrpSch}{\mathsf{GrpSch}}
\newcommand{\PSO}{\operatorname{PSO}}
\newcommand{\PO}{\operatorname{PO}}

\newcommand{\Adj}{\operatorname{Adj}}

\newcommand{\Span}{\operatorname{Span}}

\newcounter{Chapcounter}

\newcommand{\chapter}[1] 
{ {\centering          
		\addtocounter{Chapcounter}{1} \Large \underline{\textbf{ \color{blue} Chapter \theChapcounter: ~#1}} }   
	\addcontentsline{toc}{section}{ \color{blue} Chapter:~\theChapcounter~~ #1}    
}

\begin{document}

\title{
	The Basic Theory of Clifford-Bianchi Groups for Hyperbolic $n$-Space
}
\author{Taylor Dupuy, Anton Hilado, Colin Ingalls, Adam Logan }

\maketitle

\begin{abstract}
	Let $K$ be a $\QQ$-Clifford algebra associated to an $(n-1)$-ary positive definite quadratic form and let $\Ocal$ be a maximal order in $K$. 
	For each such order we define Clifford-Bianchi groups $\SL_2(\Ocal)$ which generalize Bianchi groups associated to orders in imaginary quadratic fields.
The present paper is about the arithmetic of maximal orders $\Ocal$ and the actions of $\SL_2(\Ocal)$ on hyperbolic space $\Hcal^{n+1}$ via M\"{o}bius transformations $x\mapsto (ax+b)(cx+d)^{-1}$. This  generalizes the well known theory for hyperbolic 3-space to all dimensions.
	
	We develop the general theory of orders and, via experiments, exhibit explicit orders in low dimensions that have exotic connections to lattice theory. 
	These include, for example, higher-dimensional analogues of the Hurwitz order one of which has a lattice of Clifford vectors equal to the $E_8$ lattice.
	We develop the abstract and computational theory for determining their fundamental domains for their Clifford-Bianchi groups and as well as explicit generators and relations for these groups.

	The existence of the fundamental domain is a ``proof by moduli interpretation", following a generalized version Swan's ``Bianchi-Humbert Theory'', and the domains we obtain do not follow from the general considerations involving Siegel sets.
	We make connections to the classical literature on symmetric spaces and arithmetic groups and prove that these groups are $\ZZ$-points of a $\ZZ$-group scheme and are arithmetic subgroups of $\SO_{1,n+1}(\RR)^{\circ}$ with their action via fractional linear transformations.
	
	We report on our findings concerning certain Clifford-Bianchi groups acting on $\Hcal^4$, $\Hcal^5$, and $\Hcal^6$.

Agents using Claude Code (\texttt{claude-opus-5-5}) and OpenAI Codex (model not
recorded) assisted with proofreading, Lean formalization, and figures; their
use is described in \S\ref{sec:ai}.

\end{abstract}

\tableofcontents

\section{Introduction}
The present paper is about \emph{Clifford-Bianchi groups} $\SL_2(\Ocal)$ for $\Ocal$, an order in a Clifford algebra $\Clf(q)$ associated to a positive definite integral quadratic $(n-1)$-form $q$ (general background on Clifford algebras is given in \S\ref{sec:clifford-algebras}). 
These groups act on hyperbolic $(n+1)$-space $\Hcal^{n+1}$ and simultaneously generalize the modular group action on hyperbolic $2$-space and Bianchi group actions on hyperbolic $3$-space.

These groups are not new but, in our opinion, are very neglected.
The groups were very briefly introduced by Vahlen \cite{Vahlen1902}, then picked up by \cite{Maass1949} and only repopularized by Ahlfors in the 1980s with \cite{Ahlfors1984}.
Following this there were some important contributions by Elstrodt-Grunewald-Mennicke \cite{Elstrodt1987, Elstrodt1988,Elstrodt1990},as well as \cite{Waterman1993} \cite{Maclachlan1989} \cite{Vulakh1993,Vulakh1995,Vulakh1999}. We postpone a more detailed account of these contributions until \S\ref{sec:history-intro}.

To place our manuscript into context it is perhaps beneficial to compare our manuscript to Swan's paper \cite{Swan1971} on Bianchi Groups. 
Swan's paper gives an accounting of fundamental domains and their groups and generators for Bianchi groups in a similar fashion to what we do for the Clifford-Bianchi setting. 
This is intentional as \cite{Swan1971} is foundational for various aspects of computational number theory surrounding Bianchi groups and Bianchi modular forms.
See for example  \cite{Cremona1984,Whitley1990,Cremona1994,Bygott1998,Lingham2005,Rahm2010,Sengun2012,Rahm2013}. 

Like, Swan's paper, our paper provides an algorithm for computing fundamental domains, tables of generators and relations for these groups, and a description of fundamental domains in terms of moduli of quadratic forms. 

Unlike Swan's paper, we do not have a theory of orders in Clifford algebras analogous to the theory of orders in imaginary quadratic number fields available to us. Similarly, Swan had available to him a theory of M\"{o}bius transformations on hyperbolic three space that was more developed than the analogous theory in the Clifford-Bianchi setting. 
In fact, the words ``Clifford-Bianchi setting'' is new.

In \S\ref{sec:clifford-uniformization-intro} a review of the classical theory $\SL_2(\CC_n)$ acting on $\Hcal^{n+1}$ is given highlighting some of its remarkable properties. 

In \S\ref{sec:classification} we discuss some orders $\Ocal$ that we found in low dimensions. These have strange lattices associated with them and via a search, we found, for example there is an order $\Ocal_{E_8} \subset \quat{-1,-1,-1,-1,-1,-1,-1}{\QQ}$ whose Clifford vectors are an $E_8$ root lattice and where the group $\PSL_2(\Ocal_{E_8})$ acts on $\Hcal^9$.

In \S\ref{sec:clifford-euclideanity-intro} we give an introduction to the notion of Clifford-Euclideanity. Many of the orders that we found experimentally have this property.

In \S\ref{sec:fundamental-domains-intro} we say more about arithmeticity and our $\ZZ$-group schemes, and in \S\ref{sec:bianchi-humbert} we discuss our ``Bianchi-Humbert theory'' (terminology coming from Swan) which tells us that our spaces correspond to certain ``quadratic forms'' and that the fundamental domain are those quadratic forms which attain a certain minimal value. 
See in particular Remark~\ref{rem:abstract-fundamental-domain} which contrasts our fundamental domain construction with the fundamental domain obtained from the abstract theory of symmetric spaces and arithmetic groups.

In \S\ref{sec:algorithms-intro} we give a census of our algorithmic contributions.

In \S\ref{sec:history-intro} we give a history of the Clifford uniformization and previous contributions to this area.

\subsection{The Clifford Uniformization of Hyperbolic Space}\label{sec:clifford-uniformization-intro}

When our quadratic form is $q=q(y_1,\ldots,y_m)=d_1y_1^2+d_2y_2^2+\cdots + d_m y_m^2$ over a commutative ring $R$ we use the notation 
$$\Clf(q) = \quat{-d_1,-d_2,\ldots,-d_m}{R}$$ 
for the associated Clifford algebra.Typically, $R$ will be $\RR,\QQ$, or $\ZZ$.  Although it is not required for the truth of the theorems, when we discuss specific examples over $\ZZ$ we normally assume that the $d_i$ are all squarefree.

Let $n\in \ZZ_{\geq 1}$. 
The \emph{Clifford numbers} $\CC_n$ are the Clifford algebra of the quadratic form $f_{n-1}$ on $\RR^{n-1}$ given by
$$\CC_n = \Clf(\RR^{n-1},f_{n-1}), \qquad f_{n-1}=y_1^2+y_2^2 + \cdots + y_{n-1}^2,$$ 
Explicitly, they are the associative algebra with the presentation
$$\CC_n = \RR[i_1,i_2,\ldots,i_{n-1}], \qquad i_a^2=-1, \quad 1\leq a \leq n-1, \qquad i_ai_b = -i_bi_a, \quad 1\leq a< b\leq n-1.$$ 
We have $\CC_1\cong\RR$, $\CC_2 \cong\CC$, and $\CC_3 \cong \HH$, Hamilton's quaternions. 
For every $n\geq 4$ there exists some $m< n$ such that either $\CC_n=\CC_m \oplus \CC_m$ or there exists some $r\geq 2$ such that $\CC_n = M_r(\CC_m)$. 
The Bott Periodicity Theorem for Clifford algebras tells us that for every $n\geq 1$ we have $\CC_{n+8} = M_{16}(\CC_n)$ (\S\ref{sec:arithmetic-bott}).
For what follows, $V_n \subset \CC_n$ denotes the $n$-dimensional vector subspace of \emph{Clifford vectors} $x=x_0 + x_1 i_1 + \cdots + x_{n-1} i_{n-1},$ where $x_j \in \RR$ for $0 \leq j <n$.

What is remarkable is that there exist groups $\SL_2(\CC_n) \subset M_2(\CC_n)$ of
$2\times 2$ matrices which act on $S^n$ and $\Hcal^{n+1}$ via 
$$ \begin{pmatrix} a & b \\ c & d \end{pmatrix} \cdot x = (ax+b)(cx+d)^{-1}. $$
Here we have let $S^n = V_n \cup \lbrace \infty \rbrace$ denote the $n$-sphere viewed as a one-point compactification of $V_n$, and $\Hcal^{n+1}$ is given as a subset of the Clifford vectors.

\begin{definition}
	By the \emph{Clifford uniformization of hyperbolic $(n+1)$-space} we mean hyperbolic $(n+1)$-space in its half space model as the Riemannian manifold
	\begin{equation}\label{eqn:clifford-uniformization}
		\Hcal^{n+1} = \lbrace x \in V_{n+1} \colon x_n > 0 \rbrace, \qquad ds^2 = (dx_0^2 + \cdots +dx_n^2)/x_n^2.
	\end{equation}
	given explicitly as a subset of $V_{n+1}$. 
\end{definition}
This presentation is different from the models as halfspaces in $\RR^{n+1}\cong V_{n+1}$ in that it comes equipped with a $\SL_2(\CC_n)$ representation via M\"{o}bius transformations.

In the case $n=1$ this is the usual action of $\SL_2(\RR)$ on $\Hcal^2$, and in the case $n=2$ this is the Bianchi action of $\SL_2(\CC)$ on $\Hcal^3$.
For higher $\Hcal^{n+1}$ these actions are less well-known.

Before defining our groups we recall the definition of the Clifford group. 
This subgroup of the group of units in the associative algebra plays such an important role and will be so much more commonly used than the group of units that we give it the notation $\CC_n^{\times}$.
\begin{definition}\label{def:clifford-group-anisotropic}
The \emph{Clifford group} $\CC_n^{\times}$, is the subgroup of $U(\CC_n)$ (the group of invertible elements of $\CC_n$) generated by Clifford vectors. 
\end{definition}

We now define one of our main groups. 
In what follows the symbol $*$ denotes the transpose involution which is defined in \S\ref{sec:involutions}.
We write $\CC_n^{\mon}=\CC_n^{\times}\cup\{0\}$ for the \emph{Clifford monoid}, the multiplicative monoid generated by the Clifford vectors.
For $g=\begin{psmallmatrix} a & b \\ c & d\end{psmallmatrix} \in M_2(\CC_n)$ the element $\Delta(g) = ad^*-bc^*$ is called the \emph{pseudo-determinant} of $g$.
\begin{definition}
The \emph{$2\times 2$ Clifford special linear group} $\SL_2(\CC_n)$ is the set of $g=\begin{psmallmatrix} a & b \\ c & d\end{psmallmatrix} \in M_2(\CC_n)$ satisfying the following conditions
\begin{equation}\label{eqn:clifford-algebraic}
a,b,c,d\in\CC_n^{\times}\cup\lbrace 0\rbrace, \quad \Delta(g) = ad^*-bc^*=1, \quad ab^*, cd^* \in V_n.
\end{equation}
\end{definition}
The conditions $c^*a, d^*b\in V_n$ follow from \eqref{eqn:clifford-algebraic} (see the proof of Theorem~\ref{thm:gl2-short}), and $\SL_2(\CC_n)$ is the group of Definition~\ref{def:gl2} for $C=\CC_n$ and $\Delta=1$ (Theorem~\ref{thm:gl2-short} with $R=\RR$).
These conditions \eqref{eqn:clifford-algebraic} are \emph{Clifford-algebraic}, i.e., they can be described in terms of the ring operations together with the various involutions of the Clifford algebra. 
We can convert them into genuine algebraic conditions using a Weil restriction process, which gives rise to an affine group scheme (\S\ref{sec:weil-restriction}). 
For the $\ZZ$-group scheme structures we require slightly more Clifford-algebraic conditions than those in \eqref{eqn:clifford-algebraic}.

We pause to highlight three very remarkable and unusual properties of the group $\SL_2(\CC_n)$.
\begin{enumerate}
	\item $\SL_2(\CC_n)$ is closed under the usual multiplication of matrices 
	\begin{equation}\label{eqn:group-law}
	\begin{pmatrix}
	a_1 & b_1\\
	c_1 & d_1 
	\end{pmatrix}  \begin{pmatrix}
	a_2 & b_2\\
	c_2 & d_2 
	\end{pmatrix} = \begin{pmatrix}
	a_1a_2 + b_1 c_2 &  a_1 b_2 + b_1 d_2\\
	c_1 a_2 + d_1 c_2 & c_1 b_2 + d_1 d_2 
	\end{pmatrix}.
	\end{equation}
	As stated in \eqref{eqn:clifford-algebraic}, the entries of the matrices are members of the monoid $\CC_n^{\mon}$ under the ring-multiplication of $\CC_n$.
	This is surprising because $\CC_n^{\mon}$ is not closed under addition of
	elements of the Clifford algebra, yet the group law for $\SL_2(\CC_n)$ uses
	addition.
	\item $\SL_2(\CC_n)$ acts on $\Hcal^{n+1}$ by the
	usual formula for the M\"{o}bius representation.
	If $\begin{psmallmatrix}a&b \\ c & d \end{psmallmatrix} \in \SL_2(\CC_n)$ and $x \in V_n$,
	then $(ax+b)(cx+d)^{-1} \in V_n \cup \{\infty\}$.
	We are combining $a,b,c,d \in \CC_n$ ($\dim_{\RR}(\CC_n) = 2^{n-1}$) with $x \in V_n$ ($\dim_{\RR}(V_n)=n$) to get back into the subspace $V_n$. 
	General Clifford-algebraic combinations of elements of $\CC_n$ with elements from $V_{n}$ are not in $V_n$.
	
	\item One can naively generalize the action of M\"obius tranformations on hyperbolic space via fractional linear transformations with Clifford-algebraic operations.
	While there exist theories of M\"obius transformations and inversive geometry in higher dimensional hyperbolic spaces, they make heavy use of matrices. 
	The Clifford-algebraic theory genuinely
uses a fractional linear representation of the form $x\mapsto (ax+b)(cx+d)^{-1}$.
	This makes generalizations of results from the theory of arithmetic hyperbolic 3-manifolds much more tractable. 
\end{enumerate}

We show how $\GL_2(\CC_n)$ can be arrived at as the set of $\begin{psmallmatrix} a &b \\ c & d\end{psmallmatrix}\in M_2(\CC_n)$ such that for every $m\geq n$ the map $S^m\to S^m$ defined by $x\mapsto (ax+b)(cx+d)^{-1}$ is a homeomorphism (\S\ref{sec:mobius}); $\SL_2(\CC_n)$ is its subgroup defined by $\Delta=1$.
This is where the conditions \eqref{eqn:clifford-algebraic} come from.
We give a full proof in \S\ref{subsec:clifford-unif}.

Now we turn to arithmetic.
Let $q= d_1 y_1^2 + d_2 y_2^2 + \cdots + d_{n-1} y_{n-1}^2$ be a primitive positive definite integral quadratic form on $\ZZ^{n-1}$. 
Let 
$$K = \Clf(\QQ^{n-1},q) = \quat{-d_1,-d_2,\ldots,-d_{n-1}}{\QQ} = \QQ[\sqrt{d_1}i_1, \sqrt{d_2}i_2, \ldots, \sqrt{d_{n-1}}i_{n-1}]$$ be the rational Clifford algebra for this quadratic form. 
(Note that in the case $n=2$ these $K$ are just imaginary quadratic fields.)
\begin{definition}\label{def:clifford-bianchi-groups}
	A \emph{Clifford-Bianchi group} is a group of the form $\SL_2(\Ocal)=\SL_2(\CC_n) \cap M_2(\Ocal)$ where $M_2(\Ocal)$ denotes the $2\times 2$-matrices where $\Ocal$ is an order in a Clifford algebra $K$ associated to a positive definite rational quadratic form which is closed under the reversal involution $x\mapsto x^*$.
\end{definition}
In \S\ref{sec:clifford-bianchi-groups} we give another definition that allows us to view this group as the points of a $\ZZ$-group scheme.
This is important as it gives us access to the fppf topology and fixes some issues in the literature.

\subsubsection{Arithmetic Manifolds and Satake Abelian Varieties}
Let $K$ be a positive definite rational Clifford algebra.
Let  $\Ocal$, an order in $K$ closed under reversal.
The Clifford vectors for $K$ and $\Ocal$ are defined by 
$$\Vec(K) = K \cap V_n, \qquad \Vec(\Ocal) = \Ocal \cap V_n.$$ 

Simultaneously generalizing the constructions for $\Hcal^2$ and $\Hcal^3$ we can define the partial Satake compactification of hyperbolic $(n+1)$-space.
\begin{definition}
The \emph{partial Satake compactification} of $\Hcal^{n+1}$ is the space 
$$\Hcal^{n+1,\Sat} = \Hcal^{n+1} \cup \Vec(K) \cup \lbrace \infty \rbrace,$$ with the topology of Definition~\ref{def:satake}. 
The elements $\Vec(K) \cup \lbrace \infty \rbrace$ are called the \emph{cusps}. 
\end{definition}

We remark that in contrast to the classical situation of the Bianchi groups $\SL_2(\Ocal)$, where $\Ocal$ is the ring of integers in a quadratic field $F$ and the compactification is constructed by adding
points in bijection with $\PP^1(F)$, it is not appropriate to compactify by
adding points in bijection with $K \cup \{\infty\}$, essentially because,
for $n\geq 3$, $2^{n-1} = \dim_\Q K$ is greater than the dimension of the boundary of
$\Hcal^{n+1}$, which is $n$.

\begin{definition}
Let $\Gamma$ a finite index subgroup of $\SL_2(\Ocal)$. 
The \emph{Clifford-Bianchi orbifold} and \emph{Satake compactified Clifford-Bianchi orbifold} are defined by  
$$\Ycal(\Gamma)=[\Gamma \backslash \Hcal^{n+1}], \quad \Xcal(\Gamma) = [\Gamma \backslash \Hcal^{n+1,\Sat}].$$
We will let $Y(\Gamma)=\vert \Ycal(\Gamma) \vert$ and $X(\Gamma) = \vert \Xcal(\Gamma) \vert$ denote their underlying coarse spaces.
\end{definition}
Some discussion of their categorical status is found following Definition~\ref{def:pBS}.	
These naive-looking compactifications coincide with partial Satake compactifications in the symmetric space literature \cite{Ji2002,Borel2006} (\S\ref{sec:modular-spaces}) .  

Remarkably, due to a construction of Satake \cite{Satake1966}, the points of $\Hcal^{n+1}$, with a choice of auxiliary parameters, give abelian varieties with ``even Clifford multiplication" (\S\ref{sec:satake-abelian}). This generalizes Shimura curves associated to indefinite quaternion algebras over $\Q$ and the abelian varieties in this construction are a generalization of the famous Kuga-Satake construction associated to K3 surfaces.

\begin{remark}[\S\ref{sec:satake-abelian}]\label{thm:parametrize-ab-vars}
	Let $\Gamma = \PSL_2(\Ocal)$, for $\Ocal$ an order closed under reversal in the Clifford algebra associated to a 
	positive definite integral quadratic $(n-1)$-form. 
	We expect that the orbifolds $\Ycal(\Gamma) = [\Gamma \backslash \Hcal^{n+1}]$ (with a choice of auxiliary data) parametrize abelian varieties of dimension $2^{n}$ with $\Ocal_+$-multiplication, where $\Ocal_+$ is the  even subalgebra of $\Ocal$.
	What \S\ref{sec:satake-abelian} shows is that each $z\in\Hcal^{n+1}$ gives a polarized abelian variety $A_z$ of dimension $2^{n}$.
\end{remark}

\subsection{Classification in Low Dimension}\label{sec:classification}
We have found many new and interesting orders in Clifford algebras in low dimensions which generalize well-known orders like the Gaussian integers in $\QQ[i] = \quat{-1}{\QQ}$ and the Hurwitz order in $\QQ[i,j] = \quat{-1,-1}{\QQ}$.
We give generators for $\SL_2(\Ocal)$ and describe their fundamental polyhedra.

At the end of this manuscript (\S\ref{sec:classical}--\S\ref{sec:1-1-1-1}), similar to Swan \cite{Swan1971}, we have included descriptions of fundamental domains as well as generators and relations for $\SL_2(\Ocal)$ for several explicit orders $\Ocal$; the sections are organized by rational Clifford algebra $K = \quat{-d_1,-d_2,\ldots,-d_{n-1}}{\QQ} \subset \CC_n$ while the subsections are organized by orders $\Ocal \subset K$. 

A table of the orders investigated is given in Table~\ref{tab:orders}.
Figure~\ref{fig:orders} shows how these orders contain one another.
\begin{figure}[htbp!]
\centering
\begin{tikzpicture}[lab,scale=0.92,o/.style={draw,rounded corners,inner sep=2pt,fill=white},
                    o2/.style={o,fill=cuspC!15},>=Stealth]
 \node[o] (L) at (1.0,0) {$\ZZ[i,j]$}; \node[o] (H) at (1.0,1.4) {$\Ocal_3$ Hurwitz};
 \draw (L)--(H) node[midway,right]{2};
 \node at (1.0,-0.7){$\quat{-1,-1}{\QQ}$};
 \node[o] (C13) at (5.6,0) {$\ZZ[i,\sqrt3 j]$};
 \node[o] (S1) at (4.5,1.4) {$\Ocal(-1,-3)_1$}; \node[o] (S2) at (6.7,1.4) {$\Ocal(-1,-3)_2$};
 \draw (C13)--(S1) (C13)--(S2); \draw[dashed,rounded corners] (3.4,1.1) rectangle (7.8,1.7);
 \node at (5.6,-0.7){$\quat{-1,-3}{\QQ}$};
 \node[o2] (Z3) at (10.6,0) {$\ZZ[i_1,i_2,i_3]$};
 \node[o] (Z3z) at (10.6,1.1) {$\ZZ[i_1,i_2,i_3,\zeta]$};
 \node[o] (O4) at (10.6,2.2) {$\Ocal_4$};
 \draw (Z3)--(Z3z) (Z3z)--(O4) node[midway,right]{4};
 \node at (10.6,-0.7){$\quat{-1,-1,-1}{\QQ}$};
\end{tikzpicture}\\[10pt]
\begin{tikzpicture}[lab,scale=0.92,o/.style={draw,rounded corners,inner sep=2pt,fill=white},o2/.style={o,fill=cuspC!15}]
 \node[o] (Z) at (2,0) {$\ZZ[i_1,i_2,a_3]$};
 \node[o2] (A) at (0,1.4) {$A(-1,-1,-3)$};
 \foreach \j/\x in {0/2.2,1/3.6,2/5.0}{\node[o] (B\j) at (\x,1.4) {$B_{\j}$}; \draw (Z)--(B\j);}
 \draw (Z)--(A) node[midway,left]{16}; \node[right] at (3.0,0.6){16 each};
 \draw[dashed,rounded corners] (1.7,1.1) rectangle (5.5,1.7);
 \node at (2,-0.7){$\quat{-1,-1,-3}{\QQ}$};
 \node[o] (Z4) at (9,0) {$\ZZ[i_1,\dots,i_4]$};
 \node[o2] (Ob) at (7.2,1.4) {$\Ocal_{5,!}$};
 \foreach \j/\x in {0/8.5,1/9.45,2/10.4,3/11.35,4/12.3}{\node[o,inner sep=1.2pt,font=\scriptsize] (P\j) at (\x,1.4) {$\Ocal_{5,\j}$}; \draw (Z4)--(P\j);}
 \draw (Z4)--(Ob); \draw[dashed,rounded corners] (7.95,1.1) rectangle (12.85,1.7);
 \node at (9.5,-0.7){$\quat{-1,-1,-1,-1}{\QQ}$};
\end{tikzpicture}
\caption{The orders of Table~\ref{tab:orders}. In each algebra the Clifford order is at the bottom and the orders above it contain it; numbers are indices stated in the text. Dashed boxes enclose Clifford-conjugate orders. Red fill: $\SL_2(\Ocal)$ has $2$ cusps; white: $1$ cusp. Here $\zeta=(1+i_1+i_2+i_3)/2$, $a_3=\sqrt3 i_3$, $B_j=B(-1,-1,-3)_j$.}\label{fig:orders}
\end{figure}
For each $K$, as above, we call  
$$\quat{-d_1,-d_2,\ldots,-d_{n-1}}{\ZZ} \subset K$$ 
the \emph{Clifford order}. 
We primarily study maximal orders $\Ocal$ containing the Clifford order.
Let $\Gamma=\SL_2(\Ocal)$ and $\Gamma_{\infty}$ 
be the stabilizer of $\infty$. 
In the subsection corresponding to a particular $\Ocal$ we give an explicit description of the Clifford group $\Ocal^{\times}$, a fundamental domain $F \subset V_n$ for $\Gamma_{\infty}$, an explicit description of the sides of a closed convex fundamental polyhedron $\overline{D} \subset \Hcal^{n+1}$ for  $\Gamma$ (including cusps of $\Xcal(\Gamma))$, and generators and relations for $\Gamma$
for the rows up to and including $\Ocal_4$.

\begin{table}[h]
  \begin{center}
    \small\setlength{\tabcolsep}{4pt}\begin{tabular}{|c|c|c|c|c|c|} \hline
      $\CC_n$ &  Clifford Algebra & Clifford Class of Order & Name & \#Cusps & Reference \\ \hline \hline
      $\CC_1$& $\QQ$ & $\ZZ$&  integers & 1 & \S\ref{sec:psl-2-z} \\ \hline 
      $\CC_2$ & $\QQ[i] = \quat{-1}{\QQ}$ & $\ZZ[i] = \quat{-1}{\ZZ}$ &  Gaussian integers &1 & Ex.~\ref{ex:psl-2-z-i}\\ \hline 
      & $\QQ[\sqrt{-3}] = \quat{-3}{\QQ}$ & $\ZZ[\zeta_3]$&  Eisenstein integers &1 & \cite[(2.9)]{FGT2010} \\ \hline
      & $\QQ[\sqrt{-19}] = \quat{-19}{\QQ}$ & $\ZZ[\frac{1+\sqrt{-19}}{2}]$&   &1 & \S\ref{S:sqrt19} \\ \hline
      $\CC_3$ & $\quat{-1,-1}{\QQ}$ & $\Ocal_3$ & Hurwitz order & 1 & \S\ref{S:hurwitz} \\ \hline 
      & & $\quat{-1,-1}{\ZZ}$ & Lipschitz order& 1 & \S\ref{sec:lipschitz} \\ \hline
& $\quat{-1,-3}{\QQ}$ &  $\Ocal(-1,-3)_i$, $i=1,2$ & stained glass order & 1 &\S\ref{S:quat13} \\ \hline
      $\CC_4$ & $\quat{-1,-1,-1}{\QQ}$ & $\Ocal_4$ & triality order & 1 & \S\ref{sec:maximal-c4} \\ \hline
      & & $\quat{-1,-1,-1}{\ZZ}$ & Clifford order & 2 & \S\ref{sec:1-1-1-order} \\ \hline 
      & $\quat{-1,-1,-3}{\QQ}$ & $B(-1,-1,-3)_i$, $i=0,1,2$ &  & 1 & \S\ref{sec:113-first} \\ \hline
      &  & $A(-1,-1,-3)$ &  & 2 & \S\ref{sec:113-second}\\ \hline
      $\CC_5$ & $\quat{-1,-1,-1,-1}{\QQ}$ & $\Ocal_{5,!}$ & oddball order&  2 & \S\ref{sec:c5-oddball} \\ \hline
      &  & $\Ocal_{5,i}$, $i=0,1,2,3,4$  & [5,1,4]-code order & 1 & \S\ref{sec:c5-reasonable}  \\ \hline 
    \end{tabular}
  \end{center}
  \caption{A table of interesting orders $\Ocal$ where $\SL_2(\Ocal)$ acts on certain $\Hcal^{m}$ with $m\leq 6$.  We use $\#\mbox{Cusps}$ to denote
    $\#(X(\Gamma) \setminus Y(\Gamma))$.}
  \label{tab:orders}
\end{table}

The maximal orders $\Ocal$ were found using \texttt{magma}.
This algorithm combines algorithms for $p$-maximal orders and discriminant considerations (Algorithm~\ref{alg:maximal-orders-algo}).
In these orders we discovered that, while all of them are isomorphic as algebras if the rank of the underlying lattice is even and greater than or equal to 4 (Proposition~\ref{prop:maximal-orders-unique}), the orders would clump into classes with isomorphic lattices $\Vec(\Ocal)$ and isomorphic Clifford groups $\Ocal^{\times}$. 
\begin{example}\label{ex:113}
 In the Clifford algebra $\quat{-1,-1,-3}{\QQ}$ we found four maximal orders 
$$A(-1,-1,-3), \quad B(-1,-1,-3)_0, \quad B(-1,-1,-3)_1, \quad B(-1,-1,-3)_2$$ 
and $\vert A(-1,-1,-3)^{\times} \vert = \vert B(-1,-1,-3)_i ^{\times}\vert = 24$ for $i=0,1,2$, but these groups are not isomorphic: $A(-1,-1,-3)^{\times}\cong\SL_2(\FF_3)$, while $B(-1,-1,-3)_i^{\times}$ is the dicyclic group of order $24$, which has elements of order $12$. 
This is an example where there are actually two isomorphism classes of orders which also correspond to the Clifford conjugacy classes (Proposition~\ref{prop:no-zeta-12}). 

\end{example}
Two orders are called \emph{Clifford conjugate} if they are conjugate by an element of the Clifford monoid, for instance a rational \emph{Clifford vector} (Definition~\ref{def:clifford-conjugacy});
in Example~\ref{ex:113} the $B(-1,-1,-3)_i$ for $i=0,1,2$ are one Clifford conjugacy class while $A(-1,-1,-3)$ is in its own Clifford conjugacy class.

We remark that even Clifford conjugacy doesn't capture all the information about the group. 
\begin{example}
  Let $\HH = \quat{-1,-1}{\RR}$, and let the imaginary units be $i,j$.
In the quaternion algebra $\quat{-1,-3}{\QQ} = \QQ[i,\sqrt{3}j]\subset \HH = \RR[i,j]$ there are two Clifford conjugate orders, $\Ocal(-1,-3)_1$ and $\Ocal(-1,-3)_2$, where one contains 
the 6th root of unity $(1+\sqrt{3}j)/2$ while the other contains $(1+\sqrt{3}ij)/2$.
Note that $ (1+\sqrt{3}j)/2$ is a Clifford vector while $(1+\sqrt{3}ij)/2$ is not. 
\end{example}

We will now report on how Gaussian integers generalize; this is the case of maximal orders in $\quat{-1,-1,\ldots, -1}{\QQ}$. 
\begin{example}
\begin{enumerate}
\item As is well-known, the Gaussian integers $\Ocal_2=\ZZ[i_1]$ are the unique maximal order in $\QQ[i_1]$.
\item In the quaternion case $\QQ[i_1,i_2]$, the Clifford order $\ZZ[i_1,i_2]$ is called the Lipschitz order and it is 
contained in a unique maximal order $\Ocal_3=\ZZ[i_1,i_2][\frac{1+i_1+i_2+i_1i_2}{2}]$.
We observe that $\Vec(\Ocal_3)=\Vec(\ZZ[i_1,i_2]) \cong \ZZ^3$. 
\item In the case $\QQ[i_1,i_2,i_3]$, we report that there is a unique maximal order containing $\ZZ[i_1,i_2,i_3]$ 
 $$\Ocal_4 = \ZZ[i_1,i_2,i_3][\frac{1+i_1+i_2+i_3}{2},\frac{1+i_1i_2i_3}{2}], \quad \Vec(\Ocal_4) \cong D_4^*,$$ 
where $D_4^* = \ZZ^4 + (1/2,1/2,1/2,1/2)\ZZ$ is the dual of the checkerboard lattice $D_4$; it is similar to $D_4$. 
The lattice $\Vec(\Ocal_4)$ has a $D_4$-root system whose Dynkin diagram has a triality automorphism. 
In \S\ref{sec:maximal-c4} we observe that one finds a non-Clifford vector generator of $\Ocal_4^{\times}$ that acts on this lattice by the triality automorphism. 
The general relationship governing this is unclear. 
\item One might suspect that in $\QQ[i_1,i_2,i_3,i_4]$ there is a unique maximal order containing the Clifford order, but this is no longer the case. 
We compute six maximal orders, five of which are Clifford conjugate, $\Ocal_{5,0}, \Ocal_{5,1}, \Ocal_{5,2}, \Ocal_{5,3}, \Ocal_{5,4}$, and one, $\Ocal_{5,!}$, which we have taken to calling ``the oddball''. 

The lattices of these are different. 
The lattice of $\Ocal_{5,!}$ is the standard cubic lattice while the lattices of $\Ocal_{5,i}$ have one extra generator from the cubic lattices given by 
$$\frac{i_1+i_2+i_3+i_4}{2}, \frac{i_0+i_2+i_3+i_4}{2},\frac{i_0+i_1+i_3+i_4}{2},\frac{i_0+i_1+i_2+i_4}{2},\frac{i_0+i_1+i_2+i_3}{2},$$
respectively; here we have let $i_0=1$.
Each of these corresponds to one of the $5$ doubly-even $[5,1,4]$ binary linear codes. 
The $5$ indicates the ambient space, the $1$ corresponds to the dimension, and the $4$ indicates the minimal Hamming weight of a nonzero element---in our case the one-dimensional $\FF_2$-vector spaces in $\FF_2^5$ are generated by $01111$, $10111$, $11011$, $11101$, and $11110$ respectively. \end{enumerate}
\end{example}

Generalizing this to higher dimensions it appears that for each of these maximal orders $\Ocal \subset \quat{(-1)^{n-1}}{\QQ} \subset \CC_n$ there exists a doubly even binary code $C \subset \FF_2^n$ related to $\Vec(\Ocal)$. 
In what follows we will identify $c\in C$ with $c = (c_0,c_1,\ldots,c_{n-1}) \in \ZZ^n$ consisting of vectors of zeros and ones so that its dot product with $I = (i_0,i_1,\ldots,i_{n-1})$ makes sense. 
Let $\Lambda_C=\lbrace x\in\ZZ^n \colon x \bmod 2\in C\rbrace$ be the inverse image of $C$ under $\ZZ^n\to\FF_2^n$; we call it the \emph{lattice associated to the code}, and we identify $x\in\QQ^n$ with $x\cdot I$.
For each $\Ocal$ there exists a doubly even binary code $C$ such that 
 $$ \Vec(\Ocal) = \frac{1}{2}\Lambda_C = \Span_{\ZZ}\left(\lbrace i_0,i_1,\ldots,i_{n-1}\rbrace\cup\lbrace \frac{c\cdot I}{2} \colon c \in C \rbrace\right)$$
(Lemma~\ref{lem:order-to-code}).

We conjecture the following.
\begin{conjecture} 
	For every maximal doubly even binary code $C \subset \FF_2^n$ there exists a unique maximal order $\Ocal_C \subset \quat{(-1)^{n-1}}{\QQ} \subset \CC_n$ containing the Clifford order such that 
 \begin{equation}\label{eqn:code-order}
 \Ocal_C  \supset \ZZ[\frac{c\cdot I}{2}\colon c\in C].
 \end{equation}
\end{conjecture}
Note that it is not the case that every order $\Ocal$ is equal to $\ZZ[\Vec(\Ocal)]$ as an algebra---the Hurwitz order is a counter-example to this statement. 

The existence of orders $\Ocal_C$ has been checked up to codes of length
$10$ using the list of doubly even codes at \cite[\href{https://rlmill.github.io/de_codes/}{Doubly-Even Codes}]{Miller2024}.  In fact, we have
checked that for every code $C$ up to length $10$ there exists an order
whose code is $C$.

In particular, there exists an order $\Ocal_{H(8,4)}\subset \CC_8$ associated to the Hamming code $C=H(8,4)$ whose lattice $\Vec(\Ocal)=\frac{1}{2}\Lambda_{H(8,4)}$ is a rescaled copy of the $E_8$-lattice.

\begin{remark}
Computational costs prevent the authors from implementing this construction in the case of the Golay code $G_{24}$, whose associated lattice $\Lambda_{G_{24}}$ is $\sqrt{2}$ times the Niemeier lattice with root system $24A_1$ (not the Leech lattice); similarly for the other Niemeier lattices that come from doubly even self-dual codes of length $24$.
A crude approximation tells us it would roughly take 200,000 years with our current implementation. 
The existence and investigation of these orders, and in particular their relationship to Bott periodicity phenomena, is an important area for future investigations.\footnote{Note the periodicity of $8$ in both phenomena} 
\end{remark}

\subsubsection{Clifford-Euclideanity}\label{sec:clifford-euclideanity-intro}

We now return to the general setting of $K$, a rational Clifford algebra for a positive definite form, and $\Ocal \subset K$, an order. 
In this manuscript we develop a theory of  \emph{Clifford-Euclidean} orders (\S\ref{sec:clifford-euclidean}). 
Just as in the classical theory (see \cite{Cremona1984,Whitley1990,Bygott1998} in the Bianchi case) this is an important simplifying hypothesis one can impose.
\begin{definition}
	Let $\Ocal$ be a $*$-stable order (Definition~\ref{def:clifford-stable}) and let $\Ocal^{\mon}$ be its Clifford monoid (Definition~\ref{defn:CliffordGroup}).
	We say that $\Ocal$ is \emph{right Clifford-Euclidean} if there is a norm function $N\colon \Ocal^{\mon}\to \NN$ such that $N(x)=0$ if and only if $x$ is a zero-divisor, and such that for all $x,y\in \Ocal^{\mon}$ with $N(x)>0$ and $xy^*\in \Vec(\Ocal)$ there exist some $q\in \Vec(\Ocal)$ and some $r\in \Ocal^{\mon}$ with $y=xq+r$ and $N(r)<N(x)$
	(Definition~\ref{def:clifford-euclidean}).
\end{definition}
At the heart of Clifford-Euclideanity is the ability to approximate elements of $\Vec(K)$ well by elements of $\Vec(\Ocal)$. 
In particular Clifford-Euclideanity is implied by the lattice $\Vec(\Ocal)$ having a small covering radius (\S\ref{sec:covering-radius}).

The hypothesis of having two Clifford group elements which have ratios which are vectors is ubiquitous in this theory. This is perhaps because evaluating Mobius transformations at vectors generates these examples.

We develop the theory so that, for positive definite $q$, right Clifford-Euclidean implies that the order is right cuspidally principal, and similarly on the left (Corollary~\ref{cor:Euclidean-implies-unimodular}). Naturally we make use of algorithms for left division and left greatest common divisor (Algorithm~\ref{alg:gcd}). 
Like the Bianchi setting, when $\Gamma=\SL_2(\Ocal)$ for $\Ocal$ Clifford-Euclidean, we prove that $\Xcal(\Gamma)$ has a single cusp.  
Unlike in the Bianchi setting, where the cusps of $\PSL_2(\Ocal_K)$ are in bijection with the ideal class group of $K$, the natural map from orbits of $\Vec(K) \cup \lbrace \infty \rbrace$ to the ideal class set $\Cl(\Ocal)$ (one needs to work with ideal class \emph{sets} now because our rings are non-commutative) need not be injective, even for maximal orders.
The right target is a finer \emph{cuspidal class set} $\Cl^{\mathrm{cusp}}(\Ocal)$, in which left ideals are identified only up to right multiplication by the Clifford group: the cusps map onto $\Cl^{\mathrm{cusp}}(\Ocal)$, injectively on the cusps that we call completable, and for every maximal order in Table~\ref{tab:orders} all cusps are completable, so that the cusps are in bijection with $\Cl^{\mathrm{cusp}}(\Ocal)$ (Theorem~\ref{thm:finitely-many-cusp-classes} and Remark~\ref{rem:cusp-map-not-injective}).
We also show that there are finitely many cusps in the Satake compactifications (Corollary~\ref{cor:class-set-finite}).

\subsection{Fundamental Domains and Finite Presentations of $\SL_2(\Ocal)$ }\label{sec:fundamental-domains-intro}

Just as the end of Swan's paper \cite{Swan1971} builds towards explicit examples of fundamental domains and finite presentations for $\SL_2(\Ocal_K)$ for various imaginary quadratic fields, we also build towards explicit examples of fundamental domains and finite presentations of $\SL_2(\Ocal)$, but now for $\Ocal$ various orders in Clifford algebras:
\begin{itemize}
\item In \S\ref{sec:classical} we review some modular and Bianchi groups $\PSL_2(\ZZ)$, $\PSL_2(\ZZ[i])$,
and $\SL_2(\Ocal_{\QQ(\sqrt{-19})})$, now from the viewpoint of Clifford-Bianchi groups. 
	\item In \S\ref{S:hamilton} we give the fundamental domains and presentations for the $\SL_2(\Ocal)$ for $\Ocal=\Ocal_3$, the Hurwitz order, and $\Ocal= (-1,-1/\ZZ)$, the Lipschitz order. These give the first interesting examples for fundamental domains acting on $\Hcal^4$. 	
	\item In \S\ref{S:quat13} we investigate $\SL_2(\Ocal)$ for the various orders in $\quat{-1,-3}{\QQ}$.
	This is the first example where the orders are not unique and there exist several inequivalent orders.  
	These give a striking example of how much more complicated the theory of complex multiplication can be in this setting. One particular order, which we call the stained glass order, has a fundamental domain who stabilizer of infinity happens to share symmetries with a stained glass window pictured in Figure~\ref{fig:stained-glass}. These appear to be new. 
	\item In \S\ref{sec:1-1-1} we study the $\SL_2(\Ocal_4)$ acting on $\Hcal^5$ here $\Ocal_4$ is the unique generalization of the Gaussian integers in this setting. 
	The stabilizer of infinity $\Vec(\Ocal_4)\rtimes (\Ocal_4^{\times}/\lbrace \pm 1\rbrace)$ in $\PSL_2(\Ocal_4)$ is an extension of an index-$2$ subgroup (of order $288$) of the group of orientation preserving automorphisms of the $D_4$ root system by a $D_4$ root lattice. 
	In particular the stabilizer of infinity contains the famous triality automorphism which permutes the legs of the $D_4$-root system.
	
\begin{figure}[htbp!]
	\centering
	\begin{minipage}{0.58\textwidth}\centering\begin{tikzpicture}[scale=2.6, line cap=round, line join=round,
    dot/.style={circle, fill, inner sep=1.3pt}]
  \def\r{0.57735}\draw[thin] (0,0) circle (1);
\foreach \a in {0,30,...,330} \draw[very thin] (0,0) -- (\a:1);
\draw[thin] (0:1) -- (60:1) -- (120:1) -- (180:1) -- (240:1) -- (300:1) -- cycle;
  \draw[thin] (0:1) -- (120:1) -- (240:1) -- cycle;
  \draw[thin] (60:1) -- (180:1) -- (300:1) -- cycle;
\draw[semithick] (30:\r) -- (90:\r) -- (150:\r) -- (210:\r) -- (270:\r) -- (330:\r) -- cycle;
\fill[black!12] (0,0) -- (90:\r) -- (30:\r) -- cycle;
  \draw[line width=1.6pt] (0,0) -- (90:\r) -- (30:\r) -- cycle;
  \draw[line width=1.6pt] (0,0) -- (60:0.5);
\draw[->, thin] (1,0) -- (1.18,0) node[right] {$i$};
  \draw[->, thin] (0,1) -- (0,1.18) node[above] {$j$};
\node[dot] at (0,0) {}; \node[below left=1pt, fill=white, inner sep=1pt] at (0,0) {$0$};
  \node[dot] at (0:1) {};
  \node[dot] at (90:1) {};
  \node[dot] at (60:1) {}; \node[above right] at (60:1) {$J=\tfrac12 i+\tfrac{\sqrt3}{2}j$};
  \node[dot] at (120:1) {}; \node[above left] at (120:1) {$J-i$};
  \node[dot] at (30:1) {}; \node[right=2pt] at (30:1) {$\tfrac{\sqrt3}{2}i+\tfrac12 j$};
  \node[dot] at (90:\r) {};
  \node[dot] at (60:0.5) {};
  \node[dot] at (30:\r) {};
\draw[very thin, gray] (90:\r) -- (-0.32,0.66) node[left, black, fill=white, inner sep=1pt] {$\tfrac{1}{\sqrt3}j$};
  \draw[very thin, gray] (60:0.5) -- (0.40,0.70) node[right, black, fill=white, inner sep=1pt] {$\tfrac12 J$};
  \draw[very thin, gray] (30:\r) -- (1.10,0.27) node[right, black, fill=white, inner sep=1pt] {$\tfrac12 i+\tfrac{1}{2\sqrt3}j$};
\end{tikzpicture}\end{minipage}\hfill
	\begin{minipage}{0.38\textwidth}\centering\includegraphics[width=\linewidth]{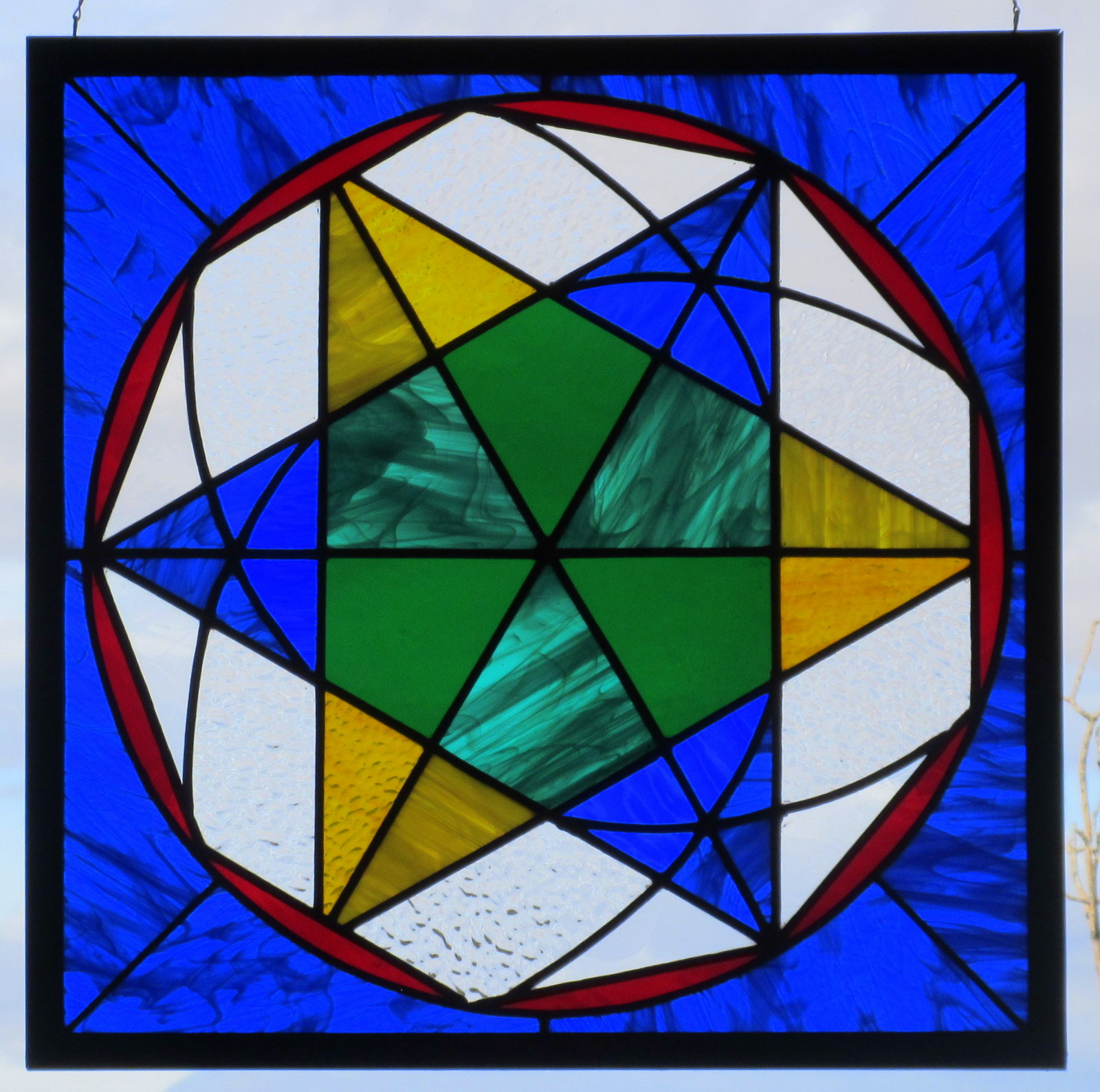}\end{minipage}
	\caption{Left: the cross section, in the plane spanned by $i$ and $j$, of the orbits of the fundamental domain for the action of $\SL_2(\Ocal(-1,-3)_2)_{\infty}$ on the boundary of $\Hcal^4$.
		A fundamental domain is the union of the two shaded triangles, drawn with heavy lines.
		Right: \emph{Enneagon Inscribed in a Circle} by Paul Powers of Power Squared Gallery in Santa Fe, New Mexico, reproduced with the artist's permission (\url{https://powersquaredglassworks.com/gallery}).
		The symmetry of the stained glass piece is closely related to that of $\SL_2(\Ocal(-1,-3)_2)_{\infty}$.
	}\label{fig:stained-glass}
\end{figure}
	\item In \S\ref{sec:1-1-3} we study the maximal orders $\Ocal$ in $\quat{-1,-1,-3}{\QQ}$. 
	These can be thought of as higher dimensional examples of the Eisenstein integers. 
	\item Finally in \S\ref{sec:1-1-1-1} we give presentations $\SL_2(\Ocal)$ for $\Ocal$ maximal orders in $\quat{-1,-1,-1,-1}{\QQ}$. These groups act on $\Hcal^6$.
	This is the first case where there is not a unique generalization of the Gaussian integers and Hurwitz order. 
\end{itemize}

\begin{remark}
	We note that finite presentations for the Lipschitz order appear in our \S\ref{sec:lipschitz}.
	The higher-dimensional Clifford order $\ZZ[i_1,i_2,i_3]$ is discussed in \S\ref{sec:1-1-1-order}, and we are able to analyze it due to $\SL_2(\ZZ[i_1,i_2,i_3])$ having finite index in $\SL_2(\Ocal_4)$.
	
	It is interesting to note that the $\SL_2(\ZZ[i_1,i_2,\ldots,i_n])$ become more poorly behaved as $n\to\infty$. 
	After our investigations we now understand this to be largely due to the ``porcupine nature" of the hypercube in larger dimensions or, equivalently, that $\ZZ[i_1,i_2,\ldots,i_n]$ lives ``deep inside" maximal orders.
	We show that as soon as we get to $\Hcal^5$ the group $\SL_2(\ZZ[i_1,i_2,i_3])$ starts to develop cusps as a consequence. 
	The distance lemma (\S\ref{sec:distance-lemma}) was one of our earlier observations of this behavior. 
	One should compare this to \cite[Theorem 11]{Maclachlan1989} where they indicate that they were also aware of this ``bad behavior" of the  $\SL_2(\ZZ[i_1,\ldots,i_n])$. 
\end{remark}

In \S\ref{sec:gens-and-rels} we detail how to extract finite presentations for these groups. 
This crucially depends on a good description of their fundamental domains with their side-pairings (elements which map a fundamental tile to an adjacent tile and some information about which tiles get identified). We use the tessellation of $\Hcal^{n+1}$ by the fundamental domain and walks on the graphs where vertices are tiles and edges are drawn between two tiles if they share a facet. 

In \S\ref{sec:gens-and-rels} we also explain the principle, going back to Poincar\'{e}, that a finite presentation of a group acting discretely on a space can be read off from the geometry of a fundamental domain for the action.

In \S\ref{sec:algorithms} we describe the algorithms we used to compute fundamental domains for $\SL_2(\Ocal)$, where $\Ocal$ is a Clifford-stable order in a Clifford algebra over $\QQ$; our debt to \cite{Bygott1998,Lingham2005,Rahm2010} is explained there.

\subsection{Fundamental Domains and Presentations: Technical Underpinnings}

We now describe our theoretical contributions to the construction of these fundamental domains. 
Sections \S\ref{sec:arithmeticity},  \S\ref{sec:representability}, and \S\ref{sec:fundamental-domains} are basic infrastructure for the theory. 

\begin{remark}\label{rem:abstract-fundamental-domain}
	The reader familiar with Borel's book \cite{Borel2019} might note that the theory of Siegel sets and reduction theory for arithmetic groups inside real semisimple Lie groups implies the existence of a fundamental domain (\cite[17.8]{Borel2019}) so once we have shown that the groups are arithmetic will have shown the existence of a fundamental domain. 
	
	This theory is inadequate for generating a side-pairings --- a necessary part of determining finite presentations for fundamental domains.
	Indeed, this is the approach taken by Elstrodt-Grunewald-Mennicke in \cite{Elstrodt1988}.

	From general considerations they know that $\vert \Gamma \backslash \Cusps(G)\vert < \infty$ where $G=\SL_2(\CC_n)$ is the real point of some algebraic group and $\Gamma$ is an arithmetic subgroup. 
	They then take $\gamma_1,\ldots,\gamma_h$ representatives of $\Gamma \backslash \Cusps(G)$ and for each $\gamma_j$ define $\Lambda_j \subset V_n$ by letting $\Lambda_j$ be the unique lattice such that $\gamma_j^{-1} \Gamma \gamma_j \cap U(\RR) = \lbrace \tau_{\lambda} \colon \lambda \in \Lambda_j \rbrace$.
	Here, $U(\RR)$ is the unipotent radical and is given as an element of $\SL_2(\CC_n)$ consisting of matrices of the form $\begin{psmallmatrix} 1 & b \\ 0 & 1 \end{psmallmatrix}$ for $b$, a Clifford vector.
	They then use fundamental domains $F_1,\ldots,F_h$ contained in $V_n$ for $\Lambda_1,\ldots,\Lambda_h$, and consider cylinders above the fundamental domain bounded by a given fixed $r>0$, so sets of the form $F_j \times \RR_{\geq r}$.
	Then for \emph{some} compact set $\Omega$ they know that the fundamental domain $D$ satisfies
	\begin{equation}\label{eqn:fundamental-domain-borel}
		D \subseteq \Omega \cup \bigcup_{j=1}^h \gamma_j(F_j \times \RR_{\geq r}).
	\end{equation}
	This is not enough to describe a side-pairing and hence generators and relations for $\SL_2(\Ocal)$.
\end{remark}

\subsubsection{Fundamental Domains in the Style of Swan}\label{sec:bianchi-humbert}
Readers familiar with \cite{Swan1971} will recognize \S\ref{sec:fundamental-domains} are our version of his ``Bianchi-Humbert Theory''. This is fundamental for determining the polyhedral structure and side-pairing of the fundamental domains. 

Recall that \cite{Swan1971} establishes a bijection between points in hyperbolic three space modulo $\SL_2(\Ocal)$ for $\Ocal$ an order in an imaginary quadratic field and Hermitian forms. 
Swan then proves that members of the fundamental domain correspond to those Hermitian forms which attain its minimal unimodular value at a specified point. 

In our setting we establish a correspondence between classes of ``positive Clifford-Hermitian forms'' and $\Hcal^{n+1}$ as $\Gamma$-sets, but unlike the Hermitian forms in the positive case, they don't seem to correspond to anything familiar and, at least from the authors' perspective, are completely contrived for the purpose of giving a ``proof by moduli interpretation'' for the purposes of copying Swan's construction of the fundamental domain.\footnote{After developing this theory  we found that \cite{Vulakh1993} has developed a number of connections in this direction, building on \cite{Elstrodt1988}, complementing our computations in connection to Lagrange Spectra, Markov Spectra, and work of Margulis.}

We now describe this theory.
Let $\Gamma = \PSL_2(\Ocal)$ act on $\Hcal^{n+1}$, where $\Ocal$ is a Clifford-stable order (Definition~\ref{def:clifford-stable}).
Let $\Gamma_{\infty}$ be the stabilizer of $\infty$. Let $F \subset V_n$ be a fundamental domain for $\Gamma_{\infty}$ containing $0\in V_n$. 
We say that $(\mu,\nu) \in \Ocal^2$ is \emph{unimodular} if and only if $\exists \begin{psmallmatrix} * & * \\ \mu & \nu \end{psmallmatrix}\in \SL_2(\Ocal)$. 
This is analogous to a pair of integers being coprime. \begin{definition}
	The \emph{bubble domain} is the set 
	$$B=\lbrace x \in \Hcal^{n+1} \colon \forall (\lambda,\mu),  \vert x -\mu^{-1}\lambda\vert \geq 1/\vert\mu\vert \rbrace$$
	where $(\lambda,\mu)$ run over unimodular pairs with $\mu\neq 0$.
\end{definition}
Using the bubble domain $B$ and $F$ we can now describe the fundamental domain $D$ for $\Gamma$.
\begin{theorem}[Theorem~\ref{T:fundamental-domain}]\label{thm:bubble-domain}
	Let $F$ be an open fundamental domain for $\Gamma_{\infty}$ acting on $V_n$ which is the interior of a compact convex polytope.
	The set
	$$D = \lbrace x_0+x_1i_1+\cdots + x_ni_n \in B^{\circ} \colon x_0+x_1 i_1+\cdots + x_{n-1} i_{n-1} \in F  \rbrace \subset \Hcal^{n+1},$$
	where $B^{\circ}$ is the interior of $B$, is an open fundamental domain for $\Gamma$.
\end{theorem}

The points $x\in \Hcal^{n+1}$ correspond to a class of positive definite Clifford-Hermitian forms (\S\ref{sec:hermitian}).
Two positive definite Clifford-Hermitian forms are considered equivalent if, and only if, they differ by some element of $\RR_{>0}^{\times}$.
This bijection between $\Hcal^{n+1}$ and classes of forms is $\SL_2(\CC_n)$ equivariant (Definition~\ref{def:matrix-to-form}, \S\ref{sec:action-hermitian}). 
One can then evaluate our forms at unimodular points to get unimodular values. 
Of all the unimodular values there is a minimal one. 
The notion of a unimodular input at which a $\RR_{>0}^{\times}$-class of forms attains its minimal unimodular value is well-defined. 
Those classes with minimal unimodular value at $(0,1)$ are those forms which correspond to points $x \in H^{n+1}$ which are contained in $B$ (\S\ref{sec:bubble-domain}).

Breaking down this condition is what then gives rise to the boundary spheres in the definition of the bubble domain $B$.
The fact that such a theory can even exist in the Clifford setting is surprising.

\subsubsection{Arithmeticity}

In \S\ref{sec:fundamental-domains} we use the fact that our Clifford-Bianchi groups act properly discontinuously and with finite covolume (See the proofs of Theorem~\ref{T:fundamental-domain} and Theorem~\ref{T:fundamental-domain-finiteness} for example).

Sections \S\ref{sec:representability} and \S\ref{sec:arithmeticity} establish discreteness and finite covolumeness via the following sequence of very technical results about algebraic groups. 

\begin{theorem}[Representability, \S\ref{sec:weil-restriction}]\label{thm:representability-main}
	Let $q$ be an integral positive definite quadratic form.
	\begin{enumerate}
		\item There exists a $\ZZ$-ring scheme $\Clf_q$ such that for every commutative ring $R$ we have  $\Clf_q(R)=\Clf(q_R),$ where $q_R$ is the quadratic form base changed to $R$.
		\item There exists a $\ZZ$-group scheme $\SL_2(\Clf_q)$ such that for every commutative ring $R$ we have $\SL_2(\Clf_q)(R) = \SL_2(\Clf(q_R)).$
	\end{enumerate}
\end{theorem}

\begin{theorem}[Spin Isomorphisms, \S\ref{sec:our-group-schemes}]\label{thm:mcinroy-correction}
Let $q$ be a positive definite integral quadratic form in $n-1$ variables and let $Q=x^2-yz+q$, a $\ZZ$-form of the real quadratic form with one negative and $n+1$ positive squares. We use the signature convention $(1,n+1)$ for one negative and $n+1$ positive squares; replacing $Q$ by $-Q$ does not change $\Spin_Q$.
Then $\SL_2(\Clf_q) \cong \Spin_Q$ as $\ZZ$-group schemes.
	Let $M+M^{T}$ be the matrix of the polar form $b_q(x,y)=q(x+y)-q(x)-q(y)$ of $q$.
	For a prime $p$, the quadratic module $Q\otimes\ZZ_{(p)}$ is regular (for $n$ even) or semiregular (for $n$ odd) in the sense of Knus \cite[Ch.~I, (5.3.5) and Ch.~IV, (3.1.4)]{Knus1991} if and only if $p\nmid\det(M+M^{T})$; for diagonal $q=d_1y_1^2+\cdots+d_{n-1}y_{n-1}^2$ with $n\geq 2$ these are the primes not dividing $2d_1\cdots d_{n-1}$.
	In particular, we have
	$$\SL_2(\CC_n) \cong \SL_2(\Clf_q)(\RR) \cong \Spin_Q(\RR) \cong \Spin_{1,n+1}(\RR).$$
\end{theorem}

\begin{theorem}[Arithmeticity \S\ref{sec:arithmeticity}]\label{thm:arithmeticity-main}
	For $\Ocal$ an order in $\quat{-d_1,-d_2,\ldots,-d_{n-1}}{\QQ}$ which is closed under $*$, the groups $\PSL_2(\Ocal)$ are arithmetic. More precisely, $\PSL_2(\Ocal)$ can be identified with an arithmetic subgroup of $\SO_{1,n+1}(\RR)^{\circ}$ acting on $\Hcal^{n+1}$, and the action is given by M\"obius transformations.
	
	As a consequence the group $\Gamma = \PSL_2(\Ocal) \subset \PSL_2(\CC_n)$ with its subspace topology is discrete, its orbits on $\Hcal^{n+1}$ are discrete in the subspace topology, and $\Gamma \backslash \Hcal^{n+1}$ has finite volume. \end{theorem}

\begin{remark}
	\begin{enumerate}
		\item 	The Spin Isomorphisms Theorem (Theorem~\ref{thm:mcinroy-correction}) requires a Bott Periodicity Theorem for Clifford algebras over a general base (see \S\ref{sec:arithmetic-bott}). 
		We could not find a reference for this theorem in the literature which works over $\ZZ$ and hence provided one. 
		\item The statement Spin Isomorphism Theorem appears implicitly in  \cite{McInroy2016}, but the definition and proof there are incorrect.
		\item Over $\QQ$, the Spin Isomorphism Theorem is obtained in \cite{Elstrodt1987} but this is inadequate for working with the fppf topology on $\Spec(\ZZ)$. 
		\item Finally, we want to emphasize that the Representability Theorem is not just for functoriality but also for representability of the functor as a scheme---that is, that its defining functor comes from polynomial equations. This is what allows us to derive exact sequences involving $\ZZ$-points of group schemes.
		
		The yoga of spin sequences and group schemes is difficult to isolate from the rest of the paper. For example the last place where Spin sequences are used is Question~\ref{prop:known-groups} and the remark after it, where we make some connections between one of our Clifford-Bianchi groups and a group that appears in the hyperbolic geometry literature).
	\end{enumerate}
\end{remark}

Discreteness and finite covolumeness follows from an application of the Borel--Harish-Chandra Theorem (\S\ref{sec:borel-harish-chandra}). This is important to justify  our later computations.

\begin{remark}\label{rem:thesis}
	Our proof of arithmeticity (Theorem~\ref{thm:arithmeticity}) uses the fppf topology: the sequence $1\to\mu_2\to\Spin_Q\to\PSpin_Q\to 1$ of Theorem~\ref{thm:short-exact-sequence} is exact over $\ZZ$ in the fppf topology but not in the \'etale topology at $p=2$. The group scheme $\PSpin_Q$ differs from $\SO_Q$ at $p=2$ and at the primes dividing $d_1\cdots d_{n-1}$ (Remark~\ref{rem:spin-sequence-at-2}).
	For example, see \cite[page 20]{Auel2009} for the issue with the \'etale topology and $p=2$ in this setting.
\end{remark}

\begin{remark}
We want to take a moment to compare what happens in our manuscript with what already exists in the literature. 
	\begin{enumerate}
		\item Arithmeticity of $\SL_2(\ZZ[i_1,\ldots,i_n])$ appears in \cite{Maclachlan1989} and for general $\SL_2(\Ocal)$, with $\Ocal$ $*$-stable, appears in \cite{Elstrodt1988}.
		The argument in \cite{Maclachlan1989} is given in terms of what we would call ``Weil restriction methods''.
		
		Remark~\ref{rem:MWW-arithmeticity} gives a detailed discussion of what \cite{Maclachlan1989} does and what we do in terms of ``Weil restriction methods''.
		Briefly, the argument of loc. cit. makes no connection to the orthogonal representation or the spin representation or spin exact sequences.
		\item In the manuscript the Representability Theorem follows directly from Theorem~\ref{thm:representability} and Definition~\ref{def:gl2}. 
		We remark that Definition~\ref{def:gl2} is adapted from a definition of Elstrodt-Grunewald-Mennicke \cite{Elstrodt1987} which in loc. cit. was used only over fields.
		We remark that Definition~\ref{def:gl2} was chosen to fit the hypotheses of Theorem~\ref{thm:representability} and the correct definition is not clear. 
		
		Interestingly, in \cite{Elstrodt1988} they define $\SL_2(\Ocal)$ as $\SL_2(\CC_n) \cap M_2(\Ocal)$. 
		For arithmeticity they simply state that these groups are stabilizers of some lattices and hence are arithmetic.
	\end{enumerate}
	
\end{remark}

\subsection{Algorithms}\label{sec:algorithms-intro}

We compute fundamental domains and finite presentations of $\PSL_2(\Ocal)$ for orders $\Ocal$ in rational Clifford algebras of positive definite forms, using \texttt{magma}.
Our algorithms are inspired by the work on Bianchi groups of Cremona \cite{Cremona1984,Cremona1994}, Whitley \cite{Whitley1990}, Bygott \cite{Bygott1998}, Lingham \cite{Lingham2005} and Rahm \cite{Rahm2010}, but they are not direct generalizations of any of these.
The algorithms we use are the following.
\begin{itemize}
	\item \emph{Maximal orders.} Algorithm~\ref{alg:maximal-orders-algo} computes the maximal orders containing the Clifford order, by applying Algorithm~\ref{alg:p-maximal} at each prime whose square divides the discriminant.
	\item \emph{Greatest common divisors.} For a Clifford-Euclidean order, Algorithms~\ref{alg:gcd} and~\ref{alg:gcd-coefs} compute a generator $\gamma$ of $\Ocal\alpha+\Ocal\beta$ and elements $\lambda,\mu$ with $\lambda\alpha+\mu\beta=\gamma$. Consequently every cusp is equivalent to $\infty$ (Corollary~\ref{cor:Euclidean-implies-unimodular}).
	\item \emph{Units.} The group $\Ocal^\times$ is finite, and Algorithm~\ref{alg:units} lists its elements by a short-vector search in a lattice. The fundamental domain $F$ of the stabilizer $\Gamma_\infty$ of $\infty$ is then a fundamental domain for $\Ocal^\times/\{\pm1\}$ acting on the Voronoi cell of $0$ in $\Vec(\Ocal)$.
	\item \emph{Bubbles.} Algorithm~\ref{alg:under-spheres} decides by linear programming whether a hemisphere lies under the union of finitely many others. We use it to discard the bubbles that do not bound the fundamental domain $D$. Remark~\ref{rem:stopping} gives a test that certifies, after the fact, that no bubble is missing.
	\item \emph{Generators and relations.} Theorem~\ref{thm:generators-ocal} gives generators of $\PSL_2(\Ocal)$ from the sides of $D$, once a matrix crossing each bubble side is known. Relations come from cycles in the tessellation by translates of $D$ (\S\ref{sec:finite-relations}).
\end{itemize}

Several parts of this are worth pointing out.
The stabilizer of $\infty$ can have rich geometry: for the order $\Ocal_4$ of \S\ref{sec:maximal-c4}, $F$ is made of $4$ of the $1152$ characteristic simplices of the $24$-cell, the Voronoi cell of $\Vec(\Ocal_4)$.
Crossing a bubble means completing a pair $(c,d)$ to a matrix in $\SL_2(\Ocal)$; \S\ref{sec:tidy} does this explicitly for ``tidy'' cusps, and the general case is open (Question~\ref{qstn:crossing-bubbles}).
In every Clifford-Euclidean example we computed, $B(0)$ is the only bubble up to $\Gamma_\infty$ (Conjecture~\ref{conj:euclidean-bubbles}).
We also treat orders that are not Clifford-Euclidean, such as $\ZZ[i_1,i_2,i_3]$, which has two cusps (\S\ref{sec:1-1-1-order}), and the order $\Ocal_{5,!}$, whose fundamental domain needs two bubbles (\S\ref{sec:c5-oddball}).
Finally, we know no a priori bound on the size of the bubbles needed (Question~\ref{qstn:bubble-bound}).

\subsection{Previous Work On Clifford-Bianchi Groups}\label{sec:history-intro}
	There are many contributions to this subject in the literature.
	The following is a history of the subject and how it relates to the present manuscript.

	We start with Ahlfors \cite{Ahlfors1984}, who, in his manuscript, gives a history of the Clifford uniformization up to 1984.  
	According to Ahlfors, the manuscript \cite{Vahlen1902} introduced the $\PSL_2(\CC_n)$ representations of $\Isom(\Hcal^{n+1})^{\circ}$ (the paper \cite{Vahlen1902} has 9 pages).\footnote{
   
	Some authors call $\SL_2(\CC_n)$ ``Vahlen groups'' and use the notation $\SV_2(\CC_n)$ or $\SV_n$.
	We have avoided this notation for clarity and Nazi affiliations (see \cite{Segal2003}).
	}
	This method was ignored until a paper of Maass in 1949 \cite{Maass1949}.
	In \cite{Ahlfors1984} Ahlfors writes
	\begin{quote}
		In a Comptes Rendus note of 1926 R. Fueter [Fu] showed that the transition from $M(C)$ to $M(U)$ can be easily and elegantly expressed in terms of quaternions. It seemed odd that this discovery should come so late, when quaternions were already quite unpopular, and sure enough a search of the literature by D. Hejhal turned up a paper from 1902 by K. Th. Vahlen [Va] where the same thing had been done, not only with quaternions but more generally, in any dimension, with Clifford numbers. It is strange that this paper passed almost unnoticed except for an unfavorable mention in an encyclopedia article by E. Cartan and E. Study. Vahlen was finally vindicated in 1949 when H. Maass \cite{Maass1949} rediscovered and used his paper. Meanwhile the theory of Clifford algebras had taken a different course due to applications in modern physics, and Vahlen was again forgotten.
	\end{quote}
	The paper \cite{Ahlfors1984} gives a proof of Mostow rigidity in this Clifford setting.
	Since Ahlfors, it appears to have been picked up by two separate groups of hyperbolic geometers in the late 1980s: Elstrodt, Grunewald, and Mennicke (EGM) \cite{Elstrodt1990}, and Maclachlan, Waterman, and Wielenberg (MWW)\cite{Maclachlan1989}.
	
	The EGM group published three papers \cite{Elstrodt1987}, \cite{Elstrodt1988}, \cite{Elstrodt1990}. 
	The EGM manuscript \cite{Elstrodt1990} proves lower bounds on the smallest positive eigenvalue of the Laplacian on $Y(\Gamma)$ for congruence groups $\Gamma$ acting on $\Hcal^{n+1}$.
	Let $\lambda_1^{\Gamma}$ be the smallest positive eigenvalue of the Laplacian on $Y(\Gamma)$; for non-cocompact $\Gamma$ the continuous spectrum is $[n^2/4,\infty)$.
	For $n\geq 2$, in loc. cit. they prove that $\lambda_1^{\Gamma}\geq (2n-1)/4$.
	The analogous bound $\lambda_1^{\Gamma}\geq 3/16$ for congruence subgroups of $\SL_2(\ZZ)$ acting on $\Hcal^2$ is due to Selberg \cite{Selberg1965}, and this generalization was an open problem around that time.
	This result was also proved by Li, Piatetski-Shapiro, and Sarnak \cite{Li1987} using different methods.
	
	In \cite{Elstrodt1988} there are some developments in Bianchi-Humbert Theory---in particular there is a correspondence between $\Hcal^{n+1}$ and Clifford-Hermitian forms but no presentation of the fundamental domain in this setting. 
	They use the general theory of arithmetic groups for their fundamental domains.
	Later, in \cite{Vulakh1993}, the theory is further developed, but computation of various $D$, as in \cite{Swan1971}, is not pursued. 
	(The idea of doing this is stated briefly in a remark on page 960 but there is an accidental conflation of $B$ and $D$).
	The excellent paper \cite{Vulakh1993} develops the theory of the Markov spectrum in this setting.
	
	Around the same time Maclachlan-Waterman-Wielenberg published a single joint paper \cite{Maclachlan1989} and Waterman published \cite{Waterman1993}, which is in the same spirit as \cite{Ahlfors1984}. 
	The paper \cite{Waterman1993} contains basic information of the sort that would be found in a chapter on M\"{o}bius transformations in a complex analysis textbook. 
	The paper \cite{Maclachlan1989} gives finite presentations for $\SL_2(\ZZ[i_1,i_2])$ in terms of graph amalgamation products.
    
	The literature has become more sparse since the late 80's, when
        these papers were published.
	We essentially know of four groups of authors which have worked with the Clifford uniformization: Vulakh \cite{Vulakh1993,Vulakh1995,Vulakh1999}, McInroy \cite{McInroy2016} (which was already mentioned), Krau{\ss}har et al \cite{Krausshar2004,Bulla2010,Constales2013,Grob2015}\footnote{Krau{\ss}har has many other papers but those can be found by following this thread}, as well as a book of Shimura \cite{Shimura2004} which appears not to take any of the aforementioned papers into account.
	
	As stated previously, \cite{Vulakh1993} develops the theory of integral Clifford forms and Markov spectra.
	The manuscripts \cite{Vulakh1995, Vulakh1999} are a continuation of this.
	Developing connections and running experiments here is a very interesting area for future investigation, especially in connection to the Satake construction.
	
	After this there are the papers of Krau{\ss}har, well summarized in the book \cite{Krausshar2004}, and the paper by McInroy \cite{McInroy2016} which was already mentioned.
	Krau{\ss}har's work is Clifford-analytic in nature and deals with Clifford-analytic analogs of classical automorphic forms mostly for the groups spanned by translations $\tau_{i_a}$ for $0\leq a \leq m$ for some $m\leq n$ and inversion $x \mapsto -x^{-1}$.
	The computation in \S\ref{sec:lipschitz-generators} shows that for $n=3$ the group generated by the translations $\tau_1,\tau_{i_1},\tau_{i_2}$ and the inversion is all of $\PSL_2(\ZZ[i_1,i_2])$; we do not know whether the analogous statement holds for $n\geq 4$.

	Finally, Shimura's AMS Monograph \cite{Shimura2004} contains some ideas related to Clifford uniformization. The bibliography contains 26 citations, 11 of which are his own papers, and none from the above references.
	The case of orthogonal groups of signature $(1,m)$ is studied briefly, and we can see in \S14.4 a variant of the Clifford uniformization. 
	We can also see some applications of Ahlfors's Magic Formula (Theorem~\ref{thm:magic}), for example, in \S14.8. 
	In the subsequent chapter the author then reverts to an adelic formalism.
	Developing connections to the work of Shimura is an interesting area for future investigation.

\subsection{AI and Software Discussion}\label{sec:ai}
Agents using Claude Code~\cite{AnthropicClaudeCode2026}
(\texttt{claude-opus-5-5}) and OpenAI Codex~\cite{OpenAICodex2026}
(model not recorded) assisted with proofreading, mathematical and expository
corrections, figure preparation, and a Lean formalization using
Lean~4~\cite{Lean4}, Mathlib~\cite{Mathlib2020}, and
Tau~Ceti~\cite{TauCeti2026}. The agent reviews and the Lean port exposed
errors that this revision corrects. Agents drafted parts of the manuscript,
including proof text. The authors reviewed and revised the AI-assisted
material incorporated in this version, including a substantial rewrite of
the maximal-orders argument. The formalization is incomplete, and the
arguments in this manuscript do not depend on it. The workflows used are
available in the public \texttt{agent-skills-public}
repository~\cite[\href{https://github.com/tdupu/agent-skills-public/blob/8d13efb5e66306aeae6e3b948bc463e7e2db1b65/README.md}{README}]{AgentSkillsPublic}.
The authors take responsibility for the submitted text. Appendix~\ref{app:ai-methods}
describes the methods and records.

\subsection*{Acknowledgements}
Dupuy is supported by the National Science Foundation DMS-2401570.
Logan was supported by Simons Foundation grant $\#550023$ in the context of the Simons Collaboration on Algebraic Geometry, Number Theory, and Computation while this research was being done; he would also like to thank ICERM for its hospitality and the Tutte Institute for Mathematics and Computation for its support and encouragement of his external research.
The authors would like to thank Asher Auel, Eran Assaf, Spencer Backman, David Dummit, Daniel Martin, Veronika Potter, and John Voight for useful conversations. 
The authors also thank Paul Powers of Power Squared Gallery in Santa Fe for allowing us to reproduce an image in Figure~\ref{fig:stained-glass}.

\section{Background on Clifford Algebras}\label{sec:clifford-algebras}

A basic reference for a general theory of Clifford algebras over any commutative ring is \cite{Hahn2004}.

In this paper we will work with Clifford algebras over an arbitrary
commutative unital ring $R$
without assuming that $2 \ne 0$ in $R$ or even that $2$ is invertible.
In order to do so, we use a definition of quadratic forms that does
not require them to arise from bilinear forms.  We begin the chapter by
presenting the definition.

\begin{definition}\label{def:quadratic-form}
  (cf.~\cite[Definition 7.1]{EKM2008})
  Let $R$ be a commutative ring and $W$ an $R$-module.  If $q: W \to R$
  satisfies $q(tw) = t^2 q(w)$ and the function $\phi: W \times W \to R$
  defined by $\phi(u,v) = q(u+v)-q(u)-q(v)$ has the property that
  $\phi(u,tv+w) = t\phi(u,v) + \phi(u,w)$ for all $u, v, w \in W$ and
  $t \in R$, 
  then $q$ is a {\em quadratic form} on $W$, and the pair $(W,q)$ is a
  {\em quadratic space}.
\end{definition}

  As pointed out in \cite{NLabQF}, we may replace the use of $\phi$ by
  its definition, thus obtaining a definition expressed purely in terms of
  $q$.  The conditions on $q$, in addition to $q(tw) = t^2 q(w)$, are that
  $q(tv+w) + tq(v)+tq(w) = tq(v+w) + t^2q(v) + q(w)$ and 
  $q(u+v+w)+q(u)+q(v)+q(w) = q(u+v)+q(u+w)+q(v+w)$ for all $u, v, w \in W$
  and $t \in R$.

It is common to abstract the properties of the function $\phi$ as follows.
\begin{definition}\label{def:bilinear-forms} (cf.~\cite[Definition 1.1]{EKM2008})
  A {\em symmetric bilinear form} is a function $B: W \times W \to R$
  satisfying $B(u,v) = B(v,u)$ and $B(u,tv+w) = tB(u,v)+B(u,w)$ for
  $u, v, w \in W$ and $t \in R$.  Given a symmetric bilinear form, we
  obtain a quadratic form $q(w) = B(w,w)$, and we say that $q$
  {\em arises from} $B$.  
\end{definition}

  If $q$ is a quadratic form, then $2q$ always arises from a symmetric
  bilinear form, namely $q(u+v)-q(u)-q(v)$.  In particular,
  if $2$ is invertible in $R$,
  then every quadratic form $q$ over $R$ can be written as $2 \cdot \frac{q}{2}$
  and hence arises from a symmetric bilinear form.
  Conversely, if $2$ is not
  invertible in $R$, then the form $q(x,y) = xy$ on the $R$-module
  $R^2$ does not arise from a symmetric bilinear form.

In most applications $W$ will be a free module of finite rank.
This isn't necessary for the next few sections but we will assume for simplicity that $W$ is projective of finite constant rank.

\begin{definition}\label{defn:CliffordAlgebra}
  Let $R$ be a commutative ring and let $(W,q)$ be a quadratic space over $R$. 

  The \emph{Clifford algebra}
  associated to $q$ is the algebra
  $$\Clf(W,q):=T(W)/I_{q}$$
  where $T(W)$ is the tensor algebra and $I_{q}$ is the ideal
  generated by $v^{2}+q(v)$ for all $v\in W$.
\end{definition}

We will also use the notations $\Clf(q)$, $\Clf(V)$ for $\Clf(V,q)$ to denote the Clifford algebra associated to $q$. If $q$ is a quadratic form over $R$ and $R'$ is any $R$-algebra, we may use the notation $\Clf(q,R')$ to denote $\Clf(q_{R'})\simeq \Clf(q)\otimes_R R'$ where $q_{R'}$ denotes the base extension of $q$ to $R'$.

Let $R$ be a ring and let $(W,q)$ be a quadratic space over $R$ where $W$ is free of rank $n$ with basis $\gamma_1,\ldots,\gamma_n$.
Then the Clifford algebra $\Clf(q)$, as an $R$-module, is free of rank $2^{n}$ \cite[Chapter IV]{Knus1991} (for $W$ only projective of rank $n$, it is projective of rank $2^n$; we use the basis $\gamma_S$ only when $W$ is free) with basis 
$\lbrace  \gamma_S : S \subset \lbrace 1,2,\ldots, n \rbrace \rbrace$ where if $S = \lbrace s_1,\ldots,s_r \rbrace$ and $s_1< s_2<\ldots <s_r$ we have $\gamma_S = \gamma_{s_1} \gamma_{s_2} \cdots \gamma_{s_r}$.
We use the notation that $\gamma_{\emptyset} =1$.

\subsection{Clifford Vectors, the Universal Property, and Involutions}\label{sec:involutions}
For every quadratic space $(W,q)$ there is a natural inclusion of $R$-modules $W \to \Clf(W,q)$.
This embedding allows us to define many things. 
First, we define the Clifford vectors. 
\begin{definition}
	\begin{enumerate}
		\item The space of \emph{Clifford vectors} $\Vec(W,q)$ is the sub-$R$-module of $\Clf(W,q)$ generated by $W$ and $1$.
		We also use the notation $\Vec(C)$ for $C$ a Clifford algebra to denote this space.
		\item The space of \emph{imaginary Clifford vectors} is the subspace of $\Vec(W,q)$, which is the image of the natural inclusion $W \to \Clf(W,q)$. We identify $W$ with this space. We also use the notation $\Im(\Vec(C))$ for $C$ a Clifford algebra to denote this space.
	\end{enumerate} 
\end{definition}

\begin{remark}\label{rem:different-names}
  Our terminology differs from that of \cite{McInroy2016} and
  \cite{Auel2009}.  What we call ``imaginary Clifford vectors'', they call
  ``Clifford vectors'', and our ``Clifford vectors'' are their ``paravectors''.
\end{remark}

We give two statements of the universal property of Clifford algebras. 
The first statement is useful for constructing maps, and the second statement is useful for proving that a certain algebra is isomorphic to the Clifford algebra of some subring.  
\begin{lemma}[Universal Property]\label{lem:universal-property}
	Let $R$ be a commutative ring. 
	Let $(W,q)$ be a quadratic space over $R$. 
	Let $A$ be an $R$-algebra. 
	If $f:W\to A$ is a morphism of $R$-modules such that $f(w)^2=-q(w)$, then there is a unique map 
	$\widetilde{f}: \Clf(W,q) \to A$ of $R$-algebras such that we have the following commutative diagram
	$$\begin{tikzcd}
	W \arrow[r, "f"] \ar[hook, d] & A \\
	\Clf(W,q) \arrow[ur,"\widetilde{f}"] & 
	\end{tikzcd}.$$
\end{lemma}
\begin{proof}
	The morphism $f$ induces a morphism from the tensor algebra $T(W) \to A$. 
	The relation $f(w)^2=-q(w)$ implies that the morphism above factors through $T(W)/I_q = \Clf(q)$. 
\end{proof}

Another way to look at this is to define the category of $q$-algebras.

Fix an $R$-module $W$ and a quadratic form $q$ on it.
The objects of this category are pairs $(A,f)$ consisting of an $R$-algebra $A$ and a morphism of $R$-modules $f:W\to A$ such that $f(w)^2=-q(w)$. 
A morphism $(A,f) \to (B,g)$ is a morphism $\varphi:A\to B$ of $R$-algebras such that 
the following diagram commutes 
$$\begin{tikzcd}
A \arrow[rr, "\varphi"] & & B \\
& W\ar[ul,"f"] \ar[ur,"g"] &
\end{tikzcd}.$$
Then Lemma~\ref{lem:universal-property} says that $\Clf(q)$ is the initial object in the category of $q$-algebras. 

This allows us to define the conjugations.
\begin{definition}[Involutions]\label{def:involutions}
	Let $C = \Clf(W,q)$ be a Clifford algebra associated to a quadratic space $(W,q)$ over a commutative ring $R$.
	Each involution is defined as a $R$-linear map from $C  $ to $C$. 
	\begin{enumerate}
		\item (Parity Involution)
		The \emph{parity involution} is the unique involution $c\mapsto c'$ extending the linear map $\alpha:W\to W$ given by $\alpha(v) =  -v$.
		\item (Transpose Involution)
		The \emph{transpose} involution is the morphism $c\mapsto c^*$ induced by the map on $T(W)$ given by $v_1 \otimes v_2 \otimes \cdots \otimes v_r \mapsto v_r \otimes v_{r-1} \otimes \cdots \otimes v_1.$
		\item(Clifford Conjugation)
		The \emph{Clifford conjugation} is $a\mapsto \overline{a} := (a^{*})' = (a')^*.$ 
	\end{enumerate}
\end{definition}
These morphisms behave on a Clifford algebra $C$ as $(ab)'=a'b'$, $(ab)^* = b^*a^*$, $\overline{ab}=\overline{b}\overline{a}$. So the transpose and Clifford conjugation are algebra anti-isomorphisms, and the parity involution is an algebra isomorphism.

\begin{remark}
Satake \cite{Satake1966} and McInroy \cite{McInroy2016} use $a^t$ for the transpose map, and
    Sheydvasser \cite{Sheydvasser2019} denotes it by $a^{\ddag}$.
\end{remark}

\subsection{Examples}

\begin{example}[Clifford Numbers]
	Let $f_{n-1}=x_{1}^{2}+\ldots +x_{n-1}^{2}$. Throughout this paper we will use $\CC_{n}$ to denote 
	$$\CC_n=\Clf(\RR^{n-1},f_{n-1}).$$
	This is sometimes called the ring of \emph{Clifford numbers}. 
	As an algebra, $\CC_n$ is generated by elements $i_1,i_2,\ldots,i_{n-1}$ which satisfy $i_j^2=-1$ for $0<j<n$ and $i_ji_k=-i_ki_j$ for $0< j<k <n$.  
The Clifford vectors in this space will be denoted by $V_n$, which is an $n$-dimensional real vector space with basis $1,i_1,\ldots,i_{n-1}$. 
This example generalizes the reals, complexes, and quaternions simultaneously as we have
	$$ \CC_1 \cong \RR, \quad \CC_2 \cong \CC, \quad \CC_3 \cong \HH.$$
	Note that the Clifford vectors in $\CC_3 \simeq \HH$ form a three-dimensional vector space. A history of Clifford numbers and their relation to hyperbolic geometry is found in \cite[Section7]{Ahlfors1984}.
\end{example}

\begin{example}[Imaginary Quadratic Fields]
	For a unary quadratic form $q(x) = dx^2$ where $d$ is an integer, we have $\Clf(\ZZ,q) \cong \ZZ[x]/(x^2+d)$, and when $-d$ is not a square we further have $\Clf(\ZZ,q) \cong \ZZ[\sqrt{-d}]$
	and $\Clf(\QQ,q_{\QQ}) \cong \QQ(\sqrt{-d})$. 
\end{example}

\begin{example}[Hyperbolic Plane over $\ZZ$]  Let $V$ be the free $\Z$-module generated by $\alpha^1,\
	\alpha^2$ and let $W$ be the free $\Z$-module generated by $\beta^1,\beta^2$, with dual bases $\alpha_1,\alpha_2$ and $\beta_1,\beta_2$. We define quadratic forms $q_V = \alpha_1^2-\alpha_2^2$ and 
	$q_W = \beta_1\beta_2$, 
	where we think of these expressions in the dual basis vectors as an operation on functions pointwise.  
	We define a map $\iota: V \to W$ of quadratic spaces over $\Z$ by $$\iota(x\alpha^1+y\alpha^2) = (x+y)\beta^1+(x-y)\beta^2.$$
	There are isomorphisms
	$$\Clf(q_V) \xrightarrow{\sim} \left\{ \begin{pmatrix} a & b \\ c & d \end{pmatrix} \in M_2(\Z)  \colon  a \equiv d, b \equiv c \pmod{2} \right\},   \quad \alpha^1 \mapsto  \begin{pmatrix} 0 & 1 \\ -1 & 0 \end{pmatrix}, \quad  \alpha^2 \mapsto  \begin{pmatrix} 0 & 1 \\ 1 & 0 \end{pmatrix}$$
	$$\quad \Clf(q_W) \xrightarrow{\sim} M_2(\Z), \quad  \beta^1 \mapsto  \begin{pmatrix} 0 & 1 \\ 0 & 0 \end{pmatrix}, \quad \beta^2 \mapsto  \begin{pmatrix} 0 & 0 \\ -1 & 0 \end{pmatrix}$$ 
	that describe the induced inclusion $\iota: \Clf(q_V) \to \Clf(q_W).$
\end{example}

\begin{example}[Indefinite Clifford Numbers]
	For integers $p,q\geq 0$ (here $q$ is not a quadratic form) we use $C_{p,q} = \Clf(\RR^{p+q}, x_1^2+\cdots +x_p^2 -y_1^2-\cdots-y_q^2)$.
	These are going to appear in our ``Bott Periodicity'' relations in Lemmas~\ref{lem:bott2} and~\ref{thm:bott1}.
\end{example}

The next example generalizes quaternionic Hilbert symbols to Clifford algebras.
\begin{example}[Clifford Algebra Symbols] Let $d_1,\ldots,d_{n-1}$ be elements of a ring $R$ and consider the $(n-1)$-ary quadratic form given by $ q(x_1,x_2,\ldots,x_{n-1}) = d_1 x_1^2 + d_2 x_2^2 + \cdots + d_{n-1} x_{n-1}^2.$
	The \emph{symbol}  is 
	$$ \left ( \frac{-d_1,\ldots, -d_{n-1}}{R} \right) := \Clf(R^{n-1},q).$$
	In the case where $R=\ZZ$ and $d_1,d_2,\ldots, d_{n-1}$ are non-zero positive integers for $0<i<n$, the inclusion of Clifford algebras 
	$$\left ( \frac{-d_1,\ldots, -d_{n-1}}{\ZZ} \right) \subset \left ( \frac{-d_1,\ldots, -d_{n-1}}{\QQ} \right)\subset \quat{-d_1,\ldots, -d_{n-1}}{\RR}$$
	$$\ZZ[\sqrt{d_1}i_1,\ldots,\sqrt{d_{n-1}}i_{n-1}] \subset \QQ[\sqrt{d_1}i_1,\ldots,\sqrt{d_{n-1}}i_{n-1}] \subset  \CC_n,$$
	provides a nice generalization of the inclusion of a ring of integers of an imaginary quadratic field into the imaginary quadratic field, which can be viewed inside $\CC$:
	$$\ZZ[\sqrt{-d}] \subset \QQ(\sqrt{-d}) \subset \CC. $$
	The last identification $ \left ( \frac{-d_1,\ldots, -d_{n-1}}{\RR} \right) \cong \CC_n$ uses that the $\gamma_j := \sqrt{d_j} i_j$ for $0<j<n$ generate $\left ( \frac{-d_1,\ldots, -d_{n-1}}{\QQ} \right)$ as a $\QQ$-algebra. 
	The order $\Clf(\ZZ^{n-1},q)$ is called the \emph{Clifford order}.
\end{example}

\subsection{Conjugation Actions and Norms}

Let $C=\Clf(W,q)$ be a Clifford algebra.
\begin{definition}\label{def:form-from-mult}
  Let $\beta(v,w) = -\{v,w\} = -(vw+wv)$.
\end{definition}
We observe that $\beta(v,w) \in R$, because
$\beta(v,w) = v^2 + w^2 - (v+w)^2 = q(v+w) - q(v) - q(w) $.
If $q$ arises from a bilinear form $B$, we have $\beta(v,w) = 2B(v,w)$.
However, we do not require $2$ to be invertible in $R$, so we cannot assume
this.
Note that the definition implies that
\begin{equation}\label{eqn:form-from-comm}
\lbrace v',w \rbrace =v'w+ wv' = \beta(v,w), \quad v,w\in W.
\end{equation}
When $q(v)$ is invertible, this tells us that the action by an imaginary Clifford vector by conjugation is just reflection:
$$v'wv^{-1} = (v'w+wv'-wv')v^{-1} = \beta(v,w)v^{-1}-wv'v^{-1} = w-(\beta(v,w)/q(v))v.$$
Note that conjugation is the reflection across the hyperplane orthogonal to $Rv$, and is given by 
\begin{equation}\label{eqn:reflection}
r_{v}(w) = w- (\beta(v,w)/q(v) )v.
\end{equation}
Observe that $r_v(v)=-v$ and if $w$ is orthogonal to $v$ then $r_v(w)=w$. 

Suppose now that $W$ has rank $n$ with basis $\gamma_1,\ldots, \gamma_n$. 
Definition~\ref{def:form-from-mult} implies, for example, that $\gamma_i\gamma_j = -\gamma_j\gamma_i-\beta(\gamma_i,\gamma_j),$ for $1\leq i<j\leq n.$
In the case that $q$ arises from a symmetric bilinear form $B$ with
an orthogonal basis,
i.e., a basis such that $B(\gamma_j,\gamma_k) =0$ for $j\neq k$, we have
$\gamma_i\gamma_j = -\gamma_j\gamma_i$.
However, this does not hold in examples such as the quadratic form $xy$
over $\ZZ$.

The Clifford norm, Clifford trace, and spinor norm are defined as 
$$ \nrd(x) = x\overline{x}, \quad \trd(x) = x+\overline{x}, \quad \ns(x) = xx^*.$$
There is another norm and trace defined via the left regular representation $\psi_x(a) = xa$ for $a \in \Clf(W,q)$. 
On the Clifford monoid, these norms will be related.

\subsection{Even Clifford Algebra}
Let $(W,q)$ be a quadratic space of rank $n$ over $R$. 
Recall that $\gamma_S = \gamma_{s_1} \cdots \gamma_{s_r}$ where $S= \lbrace s_1,\ldots, s_r\rbrace \subset \lbrace 1,\ldots, n \rbrace$ and $s_1<s_2<\ldots< s_r$ is a basis.

Every $\gamma_S$ is written as a product of an even or odd number of imaginary Clifford vectors.
When $q$ is diagonal the Clifford algebra $\Clf(W,q)$ is a $(\ZZ/2\ZZ)^n$-graded algebra where $S \subset \lbrace 1,\ldots,n\rbrace$ is put in bijection with $(\mathbb{Z}/2)^n$ so that $\gamma_S$ has the appropriate degree.
This gives the Clifford algebra $\Clf(W,q)$ a decomposition as an $R$-module in the following way
\begin{equation}\label{eqn:weight}
\Clf(W,q) = \bigoplus_{d=0}^{n} \Clf(W,q)_d, \quad \Clf(W,q)_d = \bigoplus_{S \colon \vert S \vert = d} R \gamma_S.
\end{equation}
Note that the decomposition \eqref{eqn:weight} is given by weight of the elements of $(\ZZ/2\ZZ)^n$ and is not a true grading. 
The map $(\ZZ/2\ZZ)^n \to \ZZ/2\ZZ$ given by $(c_1,\ldots,c_n)\mapsto c_1+\cdots +c_n$ induces an $\ZZ/2\ZZ$-grading, and we define the even and odd parts of the Clifford algebra, denoted by $\Clf(W,q)_+$ and $\Clf(W,q)_-$ respectively, by 
$$\Clf(W,q)_+ = \bigoplus_{d=0}^{\lfloor n/2 \rfloor} \Clf(W,q)_{2d}, \quad \Clf(W,q)_{-} = \bigoplus_{d=0}^{\lfloor (n-1)/2 \rfloor} \Clf(W,q)_{2d+1}.$$ 
We will use the fact that $\Clf(W,q)_+$ is a subring.
We remark that $\Clf(W,q)_-$ is a $\Clf(W,q)_+$-bimodule. 
When $q$ is not diagonal only the $\ZZ/2\ZZ$-grading makes sense.

\subsection{The Clifford Monoid and Clifford Groups}

In what follows, we have chosen our notation to be consistent with that of \cite{Waterman1993,Ahlfors1984,Elstrodt1987}, which differs from other places in the literature (say \cite{Knus1991}) that discusses spin groups and their connection to Clifford vectors.
Our definitions come from several desiderata: 
1) we wanted our theory to be consistent with the papers of \cite{Waterman1993,Ahlfors1984,Elstrodt1987}, this means that Clifford vectors needed to follow Ahlfors' convention and that Clifford groups over fields of characteristic zero needed to be generated by Clifford vectors for $q$ strongly anisotropic (see Definition~\ref{def:strongly-anisotropic});
2) for arithmeticity conditions we needed to have definitions which worked over $\ZZ$ and gave rise to group schemes which would allow us to apply Spin exact sequences and integral versions of these groups as subgroups of the rational versions of these groups for various integral domains.

The following Lemma can be skipped, but we keep it because it is useful to readers digging in to the various normalizations of conjugations in the literature.
They vary widely and sorting through all of them was tedious.
\begin{lemma}\label{lem:conjugations}
	Let $(W,q)$ be a quadratic space over a commutative ring $R$.
	Let $a \in \Clf(W,q)$.
	\begin{enumerate}
		\item \label{item:units-nrd} If $\nrd(a) \in R^{\times}$ then $a^{-1} = \overline{a} u$ for some $u\in R^{\times}$. 
		\item \label{item:units-ns} If $\ns(a) \in R^{\times}$ then $a^{-1} = a^* u$ for some $u\in R^{\times}$. 
		\item \label{item:conj-ns} If $\ns(a) \in R^{\times}$ then we have the following equivalences of conditions
		$$ a W a^* \subset W \quad \iff \quad a W a^{-1} \subset W, $$
		$$ a \Vec(q) a^* \subset \Vec(q) \quad \iff \quad a \Vec(q) a^{-1} \subset \Vec(q). $$
		\item \label{item:conj-nrd} If $\nrd(a) \in R^{\times}$ then we have the following equivalences of conditions
			$$ a W \overline{a} \subset W \quad \iff \quad a W a^{-1} \subset W, $$
	$$ a \Vec(q) \overline{a} \subset \Vec(q) \quad \iff \quad a \Vec(q) a^{-1} \subset \Vec(q). $$
	\item \label{item:conj-spin} If $a \in \Clf(q)_+$ then $\nrd(a) = \ns(a)$ and the condition $\nrd(a) \in R^{\times}$ (or equivalently $\ns(a) \in R^{\times}$) implies the following equivalences of conditions
		$$ a W a^* \subset W \quad \iff \quad  a W \overline{a} \subset W \quad \iff \quad a W a^{-1} \subset W, $$
	$$ a \Vec(q) a^* \subset \Vec(q) \quad \iff \quad a \Vec(q) \overline{a} \subset \Vec(q) \quad \iff \quad a \Vec(q) a^{-1} \subset \Vec(q). $$
	\item \label{item:ns-to-vec} If $a W a^* \subset W$ and $\ns(a) \in R$ then $a \Vec(q) a^* \subset \Vec(q)$.  
	\item \label{item:nrd-to-vec} If $a W \overline{a} \subset W$ and $\nrd(a) \in R$ then $a \Vec(q) \overline{a} \subset \Vec(q)$.
	\item \label{item:unit} If $a \in U(\Clf(q))$ and $\nrd(a) \in R$ then $\nrd(a)\in R^{\times}$. The converse also holds.
	\end{enumerate}
\end{lemma}
\begin{proof}
	The proofs of the second and first assertions are similar so we will only prove the first. If $\nrd(a) \in R^{\times}$ we have $a \overline{a} = u$ for some unit $u\in R^{\times}$. 
	This implies $a \overline{a} u^{-1} = 1$ which implies $a^{-1} = \overline{a} u^{-1}$; this right inverse is two-sided, because left multiplication by $a$ is then a surjective endomorphism of the finitely generated $R$-module $\Clf(q)$, hence injective, and $a(\overline{a}u^{-1}a-1)=0$. The assertion for $\ns(a) = a a^* \in R^{\times}$ is exactly the same with $a^*$ replacing $\overline a$. 
	
	To prove the third assertion we will assume $aa^* \in R^{\times}$ so that $a^{-1} = a^* u $ for some $u\in R^{\times}$. 
	The statement $aWa^{-1} \subset W$ then is $a W a^*u \subset W$ which implies $aWa^* \subset u^{-1} W = W$, so the forward direction is proved. 
	This argument is reversible and you can start with $aWa^* \subset W$ and then multiply both sides by $u$ on the right to get the converse.
	
	The argument for $\Vec(q)$ is similar and the argument for $a\overline{a} \in R^{\times}$ is the same. 
	
	For the last two assertions the condition that $a \in \Clf(q)_+$ implies that $\overline{a} = (a')^* = a^*$ since $a'=a$ for all elements of the even part of the Clifford algebra. 
	This makes the spinor norm and Clifford norm the same on the even subspace.

	Assertions \eqref{item:ns-to-vec} and \eqref{item:nrd-to-vec} follow from $\Vec(q) = R + W$, since $a\cdot 1\cdot a^* = \ns(a) \in R$ and $a \cdot 1\cdot \overline{a} = \nrd(a)\in R$.
	
	To prove \eqref{item:unit} we suppose $x \in \Clf(q)$ is a unit with $\nrd(x)\in R$. 
	Let $y$ be its inverse and $r = \nrd(x)$.
	Then $r\,\overline{y}y = x\overline{x}\,\overline{y}y = x\,\overline{yx}\,y = xy = 1$, so $r$ is a unit of $\Clf(q)$ lying in $R$.
	Since $R$ is an $R$-module direct summand of $\Clf(q)$, there is an $R$-linear map $\rho:\Clf(q)\to R$ with $\rho(1)=1$ (for $W$ free, the $\gamma_\emptyset$-coordinate), and $1=\rho(r\,\overline{y}y)=r\,\rho(\overline{y}y)$. Hence $\nrd(x)\in R^{\times}$.
	The converse is \eqref{item:units-nrd}.
\end{proof}
When going through the literature we keep the conjugation lemma in mind and come up with the following definitions.

\begin{definition}\label{defn:CliffordGroup}
  Let $\Clf(W,q)$ be a Clifford algebra over $R$.
  The \emph{Clifford monoid} of $\Clf(W,q)$ is defined to be
  \begin{equation}\label{eqn:clifford-monoid}
 \Clf(W,q)^\mon = \lbrace x \in \Clf(W,q) \colon x\overline{x} \in R , \ x \Vec(W,q) x^* \subset \Vec(W,q) \rbrace.
 \end{equation}

The \emph{Clifford group}, denoted $\Clf(W,q)^\times$,
is the subset of the Clifford monoid consisting of invertible elements.
It is a group: if $x\in\Clf(W,q)^\mon$ is invertible, then $\overline{x}x = x^{-1}(x\overline{x})x = x\overline{x}$ is a unit of $R$ by Lemma~\ref{lem:conjugations}\eqref{item:unit}, and $v\mapsto xvx^*$ is an injective endomorphism of $\Vec(W,q)$ that stays injective modulo every maximal ideal of $R$ (as $\Vec(W,q)$ is a direct summand of $\Clf(W,q)$), hence is surjective by Nakayama's lemma; so $x^{-1}\Vec(W,q)(x^{-1})^*=\Vec(W,q)$ and $x^{-1}\in\Clf(W,q)^\mon$.
We denote the full group of units of $\Clf(W,q)$ by $U(\Clf(W,q))$, though
this notion will be used infrequently. More generally, if
$S \subset \Clf(W,q)$ is an $R$-subalgebra of a Clifford algebra, we use
$S^\mon$ and $U(S)$ to denote the intersection of
$\Clf(W,q)^\mon$ with $S$ and the unit group of $S$, and we call $S^\times = U(S)\cap\Clf(W,q)^\times$ the group of \emph{Clifford units} of $S$.

Note that, even if $R$ is a field, the Clifford monoid
is not equal to the Clifford group since the monoid contains $0$ and the
group does not.
\end{definition}

\subsection{Conjugation, Norms, and the Orthogonal Group, and Strong Anisotropy}
The \emph{Clifford numbers} are the Clifford algebra $\CC_n = \Clf(\RR^{n-1}, y_1^2+\cdots+y_{n-1}^2)$. 
The associative algebra $\CC_n$ is generated by $i_1,\ldots,i_{n-1}$ satisfying $i_r^2=-1$ for $1\leq r \leq n-1$ and $i_ri_s=-i_si_r$ for any $1\leq r < s\leq n-1$. 
This makes $\CC_n$ an associative $\RR$-algebra of real vector space dimension $2^{n-1}$.
We can write elements $a\in \CC_n$ as $a= \sum_S a_S i_S$, where $a_S\in \RR$ and $S$ runs over subsets of $\lbrace 1,\ldots, n-1\rbrace$ and $i_S$ is the product of the $i_j$ for $j\in S$ written in increasing order. 
For example, $i_{\lbrace 1,3,4\rbrace} = i_1i_3i_4$. 
In later sections we often omit the braces and use notation like $i_{134}$ to denote $i_1i_3i_4$.
We give $\CC_n$ the standard Euclidean norm by the usual formula $\vert a \vert^2 = \sum_{S} a_S^2$.

Just as with all Clifford algebras, $\CC_n$ comes with the involutions $a\mapsto a',a^*,\overline{a}=(a')^*=(a^*)'$ (Definition~\ref{def:involutions}), the last of which is related to the absolute value by $\vert a \vert^2 = \Re(a \overline{a})$.

Finally, the center of $\CC_n$, which we denote by $Z(\CC_n)$, is $\RR$ if $n$ is odd and $\RR[I]$ where $I=i_1i_2\cdots i_{n-1}$ if $n$ is even. 
We note that 
$$ I^2 = (-1)^{n(n-1)/2} = \begin{cases}
1,& n = 0,1 \mod 4\\
-1,& n = 2,3 \mod 4 
\end{cases}.$$

The Clifford algebra $\CC_n$ is an exceptionally nice Clifford algebra. 
Most of its nice properties arise from the fact that the quadratic form $f_{n-1}=y_1^2 + \cdots + y_{n-1}^2$ over $\RR^{n-1}$ has a property called \emph{strong anisotropy} which we now define. 

First we need an auxiliary quadratic form.
\begin{definition}
	Let $R$ be a commutative ring and let $W$ be a free $R$-module of rank $n$. 
	Let $(W,q)$ be a quadratic space over $R$.  Let $\{\gamma_S : S \subset \{1,\ldots,n\}\}$ be the usual basis of $\Clf(W,q)$.
	For $a \in \Clf(W,q)$ let $r_\emptyset(a)$ be the $\gamma_\emptyset$ component of $\nrd(a)$:
	$$\nrd(a) = a\overline{a} =r_{\emptyset} \gamma_{\emptyset} + \sum_{S \neq \emptyset} r_S \gamma_S.$$
	for some $r_S \in R.$ 
	Note that each $r_S=r_S(a) $ is a function of $a$ alone.
	We define the \emph{big form} $\hat{q}$ on the Clifford algebra $\Clf(W,q)$ by
	$$ \hat{q}(a) = r_{\emptyset}(a).$$
	This is just the $\gamma_{\emptyset}=1$-component (or ``real component'') of $a\overline{a}$.
\end{definition}

This defines a new quadratic space from an old one. 
We now arrive at the following definition.
\begin{definition}\label{def:strongly-anisotropic}
	We say that $(W,q)$ is {\it strongly anisotropic} if the quadratic space $(\Clf(W,q),\hat{q})$ is anisotropic.\
	We will call a Clifford algebra \emph{strongly anisotropic} if its defining form is.
\end{definition}

We will see in Proposition~\ref{prop:big-form-02} that $f_{n-1}$ is strongly anisotropic, and more generally so is $d_1y_1^2+\cdots+d_{n-1}y_{n-1}^2$ with all $d_j>0$, the big form being computed in this basis.

\begin{remark}
We note that the definition of the form $\widehat{q}$ depends on the fact that the algebra basis being generated by a basis of vectors and freeness of the module. 
Here is a example that shows that this definition is not well defined if we allow any basis since the notion of ``scalar part'' can be ambiguous.

  Consider the form $ax^2+bxy+cy^2$ over the ring $\ZZ[a,b,c]$.
 Let $a = \alpha + \beta e_1+ \gamma e_2 + \delta e_1e_2.$
 So the ``scalar part'' of $a$ is $\alpha$, but if we reorder our basis then
 $a =  \alpha + \beta e_1+ \gamma e_2 + \delta(b-e_2e_1)$
 which now has ``scalar part'' $\alpha +\delta b$.  
 This issue also occurs when we take the ``scalar part'' of elements of the form $a\overline{a}$.
 
  The simplest fix is to use the trace to split the inclusion $R \to \Clf$, but we need $1/2 \in R$ to do this. We then can define $\widehat{q}(x) = \frac{1}{2^n}\Tr(x \overline{x})$ where $\Tr$ is the trace of the left regular representation. 
  
  We use a basis of the free module to generate the standard basis of the Clifford algebra.
\end{remark}

The geometric content and utility of $\CC_n$ and $\quat{-d_1,\ldots,d_{n-1}}{\QQ} \subset \CC_n$ for $d_1,\ldots,d_n \in \NN$ come from its interaction with its Clifford vectors. 
In $\CC_n$ we use the special notation $V_n=\Vec(\CC_n)$ for the Clifford vectors. 
First, note that for $a \in V_n$ nonzero we have $a\bar{a} = \vert a \vert^2$.
This implies that $a^{-1} = \bar{a}/\vert a \vert^2$ for $a\in V_n$.
\begin{definition}
	We let $\Gamma_n$ denote the subgroup of $\CC_n^{\times}$ generated by the nonzero Clifford vectors under multiplication.
	(Similarly, if $C=\Clf(W,q)$ with $q$ a nondegenerate quadratic form on a finite-dimensional vector space $W$ over a field of characteristic $\neq 2$, let $\Gamma(C)$ be the subgroup of $C^\times$ generated by the invertible Clifford vectors; if $(W,q)$ is strongly anisotropic these are exactly the nonzero Clifford vectors.)
\end{definition}
We will see in a moment that $\Gamma_n$ (resp.~$\Gamma(C)$) is actually the
Clifford group $\CC_n^\times$ (resp.~$C^\times$).

Note that if $a \in \Gamma_n$ the anti-commutative property of $a\mapsto \overline{a}$ implies that $a\overline{a} = \vert a \vert^2$ (resp.~$a \overline a = \widehat{q}(a)$).

There is some interesting geometry related to
$$\pi_a: \Gamma_n \to \O(V_n), \quad \pi_a(x) = axa^*/\vert a\vert^2  =ax(a')^{-1}.$$

This is a rotation, and it is a classical fact that this map is surjective onto the special orthogonal group \cite[Theorem 3]{Waterman1993}.
The proof generalizes in a straighforward manner to the strongly anisotropic case.
Surjectivity is a consequence of the fact that the orthogonal group is generated by reflections.
In fact, \cite[Theorem 2]{Waterman1993} tells us that this transformation for $a\in V_n$ is the reflection $r_1$ followed by the reflection $r_a$ where for $b\in V_n$ the reflection $r_b$ is the reflection in the plane perpendicular to $b$.
The composition is, then, a rotation around the plane spanned by $1$ and $a$ in the counterclockwise direction. 
For example, $\pi_{i_j}(x) = i_jx(i_j')^{-1} = i_j x i_j$, which is a $180$-degree rotation in the $x_0x_j$-plane where $x=x_0+x_1i_1+\cdots + x_{n-1}i_{n-1}$.
We also make the observation that $\pi_a(x) = axa^*$ preserves $V_n$ for $a\in \CC_n^\times$ because $i_rxi_r$ does for each $r$.
Similar statements can be proved in the strongly anisotropic case (albeit without the direct geometric interpretations in Euclidean space).

\begin{lemma}\label{lem:generated-by-clifford-vectors}
	The group $\Gamma_n$ (resp.~$\Gamma(C)$, for $C$ as above) generated by the nonzero (resp.~invertible) Clifford vectors of $\CC_n$ (resp.~$C$) under multiplication is equal to $\CC_n^\times$ (resp.~$C^{\times}$), the Clifford group. 
\end{lemma}
\begin{proof}
	We treat both cases at once; let $K$ be the base field ($K=\RR$ for $\CC_n$).
	For $v=v_0+w\in\Vec(C)$ with $v_0\in K$, $w\in W$, put $N(v)=v\overline{v}=v_0^2+q(w)$; this is a nondegenerate quadratic form on $\Vec(C)$ (for $\CC_n$ it is $\vert v\vert^2$), with polar form $B(u,v)=u\overline{v}+v\overline{u}$, and $v$ is invertible if and only if $N(v)\neq 0$, with $v^{-1}=\overline{v}/N(v)$.
	For $b$ with $N(b)\neq 0$ let $r_b(v)=v-(B(b,v)/N(b))\,b$ be the reflection of $(\Vec(C),N)$ in $b$ (cf.~\eqref{eqn:reflection}); note $r_1(v)=-\overline{v}$.
	It is clear that $\Gamma(C)\subset C^\times$.
	Let $x\in C^\times$ and $r=x\overline{x}\in K^\times$. Then $\overline{x}=rx^{-1}$, so $\overline{x}x=r$, and applying the parity involution gives $x^*x'=r$.
	By the definition of the Clifford monoid, $T_x(v)=r^{-1}xvx^*$ maps $\Vec(C)$ to itself, and
	\begin{samepage}
	$$N(T_x(v)) = r^{-2}\,xv(x^*x')\overline{v}\,\overline{x} = N(v),$$
	so $x\mapsto T_x$ is a homomorphism $C^\times\to\O(\Vec(C),N)$.
	\end{samepage}
	For an invertible Clifford vector $a$ we have $a^*=a$, and $r_a(r_1(v)) = -\overline{v}+(av+\overline{v}\,\overline{a})a/N(a) = ava/N(a)$, so $T_a=r_a\circ r_1$.
	By the Cartan--Dieudonn\'e theorem every element of $\SO(\Vec(C),N)$ is a product of an even number of reflections, and $r_br_c=(r_br_1)(r_cr_1)^{-1}=T_bT_c^{-1}$; hence $\SO(\Vec(C),N)\subseteq T(\Gamma(C))$.
	Conversely, $T(C^\times)\subseteq\SO(\Vec(C),N)$: otherwise $r_1=T_x^{-1}(T_xr_1)$ would equal $T_y$ for some $y\in C^\times$; then $yy^*=-y\overline{y}$ gives $y^*=-\overline{y}$, so $y$ is odd, and $T_y(w)=w$ gives $yw=-wy$ for all $w\in W$; expanding $y$ in the basis $\gamma_S$ attached to an orthogonal basis of $W$ shows that no nonzero odd element anticommutes with all of $W$.
	Finally, if $T_x=\mathrm{id}$ then $xx^*=r=x\overline{x}$, so $x^*=\overline{x}$ and $x$ is even, and $xw=wx$ for all $w\in W$; the same expansion shows that an even element commuting with $W$ lies in $K$. So the kernel of $T$ is $K^\times\subset\Gamma(C)$, and therefore $C^\times=\Gamma(C)$.
\end{proof}

The following Lemma is very useful so we call it the ``Useful Lemma''.
\begin{lemma}[The Useful Lemma]\label{lem:useful}
	Let $a,c\in \CC_n^{\times}$. 
	We have $a^*c \in V_n$ if and only if $ac^{-1}\in V_n$.
	A similar statement holds for $C$ a strongly anisotropic Clifford algebra over a field of characteristic not equal to $2$. 
\end{lemma}
\begin{proof}
	For $b \in \CC_n^{\times}$, we have $\nrd(b) \cdot b^{-1} = \overline{b}$, because this is true for $b \in V_n \setminus \{0\}$, and if it holds for $b_1, b_2 \in \CC_n^{\times}$ it holds for $b_1 b_2$ as well. Since $\overline{c} = \nrd(c)c^{-1}$ with $\nrd(c)\in\RR^\times$, we have $ac^{-1} \in V_n$ if and only if $a\overline{c} \in V_n$. Since $c^*\in\CC_n^\times$, the map $x\mapsto c^*x(c^*)^*=c^*xc$ and its inverse preserve $V_n$, so this holds if and only if $c^*(a\overline{c})c = \nrd(c)\,c^*a \in V_n$, if and only if $c^* a \in V_n$, if and only if $a^* c \in V_n$ (because $V_n$ is closed under $*$).
\end{proof}

Here are some results connected with Clifford norms and strong anisotropy.

\begin{lemma}\label{lem:norms-of-vectors}
  Let $v$ be a Clifford vector that is not an element of $R$.
  Then $v^2 - (\trd v) v + \nrd v = 0$. If $R$ is an integral domain, then
  $x^2 - (\trd v) x + \nrd v$ is the minimal polynomial of $v$ over the fraction field of $R$.
\end{lemma}

\begin{proof} Writing $v = v_0 + v_W$ where $v_0 \in R$ and $v_W$ is an
  imaginary Clifford vector, we have $\bar v = v_0 - v_W$ and so
  $\trd v = 2v_0, \nrd v = v_0^2 - v_W^2$.  Thus
  $v^2 - (\trd v) v + \nrd v = v_0^2 + 2v_0 v_W + v_W^2 - 2v_0(v_0 + v_W) + v_0^2 - v_W^2 = 0$.
  If $R$ is an integral domain with fraction field $K$, then $v_W\neq0$ in
  $W\otimes_RK$, because $v\notin R$ and $W$ is projective; so $v\notin K$, the
  minimal polynomial of $v$ over $K$ has degree $2$, and it is
  $x^2-(\trd v)x+\nrd v$.
\end{proof}

\begin{lemma}\label{lem:char-poly-vec}
  Let $v$ be a Clifford vector in $\Clf(q)$ that is not in $R$.
  Then the characteristic polynomial of the $R$-module endomorphism of
  $\Clf q$ given by (left or right) multiplication by $v$ is
  $(x^2-(\trd v) x + \nrd v)^{2^{n-1}}$, where $n$ is the rank of $W$.
\end{lemma}

\begin{proof}
  Since characteristic polynomials are compatible with specialization, it
  suffices to do this in a generic example.  Thus we take $q=\sum_{1\le i\le j\le n}a_{ij}x_ix_j$ to be a generic
  quadratic form in $n$ variables and $v=y_0+\sum_{j=1}^n y_j\gamma_j$ a generic
  Clifford vector; the ring is then the integral domain
  $R = \ZZ[a_{ij},y_0,\dots,y_n]$, with fraction field $K$.
  Write $v=v_0+v_W$, with $v_0\in R$ and $v_W\in W$.
  The polynomial $x^2-(\trd v)x+\nrd v$ annihilates $v$, hence annihilates
  multiplication by $v$, and its discriminant
  $(\trd v)^2-4\nrd v=4v_W^2=-4q(v_W)$ is not a square in $K$, since it has
  degree $1$ in $a_{11}$. So this polynomial is irreducible over $K$, and it is
  the minimal polynomial of multiplication by $v$ on $\Clf(q_K)$.

  It follows that the characteristic polynomial is a power of this
  polynomial.  Comparing degrees ($\Clf(q)$ has rank $2^n$) gives the desired result.
\end{proof}

  \begin{proposition}\label{prop:big-form-02} 
	Let $R$ be an integral domain with fraction field $K$ and let $W$ be a free $R$-module of rank $n$. 
	Suppose $R$ does not have characteristic 2.
	Let $(W,q)$ be a quadratic space. 
        Let $a \in \Clf(W,q)$.
	\begin{enumerate}
	\item  If $\nrd(a) \in R$ then $\nrd(a) = \hat{q}(a)$.
		\item If $R \subseteq \RR$ and the basis $\gamma_1,\ldots,\gamma_n$ of $W$ used to define $\hat{q}$ is orthogonal, with $q(\gamma_i)>0$ for all $i$, then $(W,q)$ is strongly anisotropic.
		\item If $(W,q)$ is strongly anisotropic then $\nrd(a) =0$ if and only if $a = 0.$
		\item If $N$ is the algebra norm defined by
$N(a) = \det(x \mapsto ax)$ on $\Clf(q)$, $q_K$ is nondegenerate and $n\geq 1$, we have $N(a) = \nrd(a)^{2^{n-1}}$ if $a$ is in the Clifford monoid. \end{enumerate}
\end{proposition}\label{prop:anisotropic}
\begin{proof} The first and third statements are immediate. For the second, by hypothesis
	$q = d_1y_1^2+\cdots+d_ny_n^2$ in the basis $\gamma_1,\ldots,\gamma_n$, where all $d_i > 0$.
	Let $a = \sum_{S} a_S \gamma_S \in \Clf(W,q)$.
	Since the basis is orthogonal, $\gamma_S\overline{\gamma_T}$ is a scalar multiple of $\gamma_{S\triangle T}$, so
	the only products that contribute to $\hat{q}(a)$ are of the form
	\begin{equation}\begin{split}
	\gamma_{s_1}\cdots\gamma_{s_r}\overline{ \gamma_{s_1}\cdots\gamma_{s_r}}
	& =   (-1)^r \gamma_{s_1}\cdots\gamma_{s_r} \gamma_{s_r}\cdots\gamma_{s_1}\\
	& =  (-1)^r (-d_{s_1})\cdots (-d_{s_r})\\
	& =  d_{s_1} \cdots d_{s_r}.
	\end{split}\end{equation}
	So we see that $\hat{q}(a) = \sum_{S} a_S^2 \prod_{s \in S} d_s \geq 0$, with equality if and only if all $a_S = 0$.
	
	For the fourth, let $a$ be in the Clifford monoid and $r=\nrd(a)\in R$.
	If $r=0$ and $a\neq 0$, then $\overline{a}\neq 0$ and $a\overline{a}=0$, so $x\mapsto ax$ is not injective and $N(a)=0=r^{2^{n-1}}$.
	If $r\neq 0$, then $a$ is invertible in $\Clf(q_K)$, with inverse $\overline{a}/r$, so it lies in the Clifford group of $\Clf(q_K)$,
	and by Lemma~\ref{lem:generated-by-clifford-vectors} we can write $a=v_1\cdots v_m$ with invertible $v_i\in\Vec(q_K)$.
	Both $N$ and $\nrd$ are multiplicative on such products, and $N(v)=\nrd(v)^{2^{n-1}}$ for Clifford vectors $v$:
	for $v\notin K$ this is Lemma~\ref{lem:norms-of-vectors} and Lemma~\ref{lem:char-poly-vec}, and for $v\in K$ both sides equal $v^{2^n}$.
	Hence $N(a)=\nrd(a)^{2^{n-1}}$.
\end{proof}

\begin{lemma}\label{lem:same-norms-clifford-monoid} 
	On the Clifford monoid of $\CC_n=\quat{-1,-1,\dots,-1}{\RR}$, 
	the Clifford norm $\nrd(a) = a \bar a$
	is real and coincides with the Euclidean norm $|a|^2 = \sum_i a_i^2$.
\end{lemma}

\begin{proof} By the definition of the Clifford monoid, $\nrd(a) = a\bar a \in \RR$.
	Since $|a|^2 = \Re(a\bar a)$, it follows that $\nrd(a) = |a|^2$.
\end{proof}

This also lets us see the behavior on basis elements. 
\begin{corollary}\label{cor:norms-over-q} 
	Let $d_1, \dots, d_n$ be positive rational numbers and let
	$C=\quat{-d_1,\dots,-d_n}{\QQ}$. On the Clifford monoid of $C$, the
	Clifford norm $\nrd(a)=a\bar a$ coincides with the scaled Euclidean norm
	for which the set
	\[\left\{\gamma_S: S \subseteq \{1,2,\dots,n\}\right\}\]
	is orthogonal and $\gamma_S$ has norm $\prod_{i\in S}d_i$.
	Here $\gamma_1,\ldots,\gamma_n$ are the generators with $\gamma_j^2=-d_j$.
\end{corollary}

\begin{proof} The embedding of $\quat{-d_1,\dots,-d_n}{\QQ}$
	into $\quat{-1,\dots,-1}{\RR}$, taking the generators to
	$\sqrt{d_j} i_j$ for $1 \le j \le n$, preserves both of these, so the
	result follows from Lemma~\ref{lem:same-norms-clifford-monoid}.
\end{proof}

We note that the Clifford norm of a Clifford algebra does not coincide with its reduced norm as an order.
\begin{example}
  Consider the Clifford algebra of a non-degenerate quadratic form $q$ on $\CC^{10}$ over $\CC$.  Then $C = \Clf(\CC^{10},q) \simeq \CC^{2^5 \times 2^5 }$.
  Let $a \in \CC.$  Let $\Nrd$ be the reduced norm of $C$ as an order over $\CC$, as defined in \cite[pg 122]{Reiner1975}.
  Then
$$
    \nrd(a)  =  a^2, \quad N(a)  = a^{2^{10}}, \quad \Nrd(a)  = a^{2^5}.$$
    \end{example}

We have the following relationships.
\begin{corollary}
  Let $R$ be an integral domain of characteristic $\neq 2$ with fraction field $K$, and let $W$ be a free $R$ module of even rank $n$ with a quadratic form $q$ such that $q_K$ is non-degenerate.
  \begin{enumerate}
    \item For any $a \in \Clf(W,q)$ we have $N(a) = \Nrd(a)^{2^{n/2}}.$
    \item For any $a$ in the Clifford monoid we have $N(a) = \nrd(a)^{2^{n-1}}.$
  \end{enumerate}
  \end{corollary}
\begin{proof}
  By the Corollary in \S\ref{sec:factoring}, $\Clf(q_K)$ is central simple over $K$ of degree $2^{n/2}$, and for a central simple algebra of degree $m$ the norm of the left regular representation is $\Nrd^m$ \cite[\S 9]{Reiner1975}; this gives (1). Part (2) is Proposition~\ref{prop:big-form-02}(4).
\end{proof}
\subsection{The Pin and Spin Groups}

\begin{definition}\label{defn:spin}
	We define the \emph{Spin group} to be 
	$$\Spin(W,q) = \lbrace a \in \Clf(q)_+^\mon \colon \ns(a)=1, aWa^* \subset W \rbrace.$$
\end{definition}
For degenerate $q$ this group can be strictly larger than the classical spin group: for $q=0$ on $\QQ^2$ the element $1+e_1e_2$ lies in $\Spin(W,q)$. For $q$ nondegenerate over a field and representing a unit it agrees with the classical spin group, and the condition $a\in\Clf(q)^{\mon}$ is then redundant. 
Note that the condition $a \in \Clf(q)_+$ implies that $\ns(a) = \nrd(a)$. 
This means that the condition $aWa^* \subset W$ could be replaced by $a W a^{-1} \subset W$ or $a W \overline{a} \subset W$.
The condition that $a \in \Clf(q)^{\mon}$ also implies that $a \Vec(q) a^* \subset \Vec(q)$ for $a \in \Spin(W,q)$.
For the purpose of explaining how this Spin group matches up with other Spin groups from the literature that some readers may be more familiar with (and to cite theorems from these papers), we compare the above definition of the Spin group with an alternative definition. 
To state this Lemma we need to define the \emph{imaginary Clifford}, \emph{general spin}, and \emph{pin} groups: 
$$ \widetilde{\Clf}(W,q)^{\times} = \lbrace x \in \Clf(W,q)^{\times} \colon xW(x')^{-1} \subset W \rbrace, \quad \GSpin(W,q) = \widetilde{\Clf}(W,q)^\times_+,$$
The argument after Definition~\ref{defn:CliffordGroup}, applied to $W$ and $w\mapsto xw(x')^{-1}$, shows that $\widetilde{\Clf}(W,q)^\times$ is closed under inverses.
$$\Pin(W,q) = \lbrace u \in \widetilde{\Clf}(W,q)^\times \colon u \bar{u} =1 \rbrace.$$ 
The literature sometimes defines
  Spin as $\Pin(W,q) \cap \GSpin(W,q)$.
We show that this is a consequence of our definitions.
\begin{lemma} 
	\begin{enumerate}
		\item $\Spin(W,q) = \Pin(W,q) \cap \GSpin(W,q)$.
		\item $\GSpin(W,q) = \lbrace a \in \Clf(q)_+^\mon \colon \nrd(a)\in R^\times, aWa^* \subset W \rbrace $
	\end{enumerate}
\end{lemma}
\begin{proof}
	We prove the first assertion. 
	Suppose that $x \in \Pin(q) \cap \GSpin(q)$.
	Then $x \overline{x} = 1$ implies that $xx^* =1$.
	Also, note that in this situation $x^*=x^{-1}$ and $x'=x$.
	By the property that $x \in \widetilde{\Clf}(q)^\times$ we get $xW(x')^{-1} \subset W$ but $(x')^{-1} = x^*$ and we are done.
	
	Conversely, suppose that $a\in \Spin(q)$.
	The condition $\ns(a)=1$ implies $a\overline{a}=1$.
	The condition $aWa^*\subset W$ implies that $aW(a')^{-1} \subset W$.
	Since $a\in \Clf(q)_+^\times$ is invertible we have $a\in \widetilde{\Clf}(q)$ and hence $a\in \Pin(q)$.
	Since $a$ is even we have $a \in \GSpin(q)$, which proves the result.
	
	Using our conjugation lemma and the definition of $\Clf(q)^{\mon}$ we can check that 
	\begin{eqnarray*}
          \GSpin(q) &=&  \lbrace a \in \widetilde{\Clf}(q)_+ \colon \ns(a)\in R^{\times}, a W a^* \subset W \rbrace \\
          &=& \lbrace a \in \Clf(q)_+ \colon a W a^{-1} \subset W, a \Vec(q) a^* \subset \Vec(q) , \nrd(a)=\ns(a) \in R^{\times} \rbrace \\
		&=& \lbrace a \in \Clf(q)_+^{\mon} \colon \nrd(a) \in R^{\times}, aW a^{-1} \subset W \rbrace \\
	  &=& \lbrace a \in \Clf(q)_+^{\times} \colon a W (a')^{-1} \subset W \rbrace
	\end{eqnarray*}
These equalities are just compositions of facts from Lemma~\ref{lem:conjugations}.
One unravels both definitions to arrive at the long definition in the middle of these equalities.
\end{proof}

\begin{remark}
  We caution the reader again that notations vary from source to source.
  The Clifford monoid is denoted as $\mathcal{PT}(W)$ in \cite[\S 6]{McInroy2016}, where it is called the ``paravector Clifford group''.
  What we call the imaginary Clifford group is called the ``Clifford group'' in both \cite{Auel2009} and \cite{McInroy2016}.
The definition of $\Pin(W,q)$ in \cite{Auel2009} (following \cite{Knus1991}) is different from the one here: in particular, the condition
$\ns(u) = uu^*=1$ is imposed.
In \cite[p. 220]{Lounesto2001} there is yet another definition:
elements of the Pin group are required to have $u\overline{u}= \pm 1$
(and \cite{Lounesto2001} only works with real Lie groups).
	Since the conjugation $x\mapsto x'$ is trivial on the even part of Clifford algebras, our Spin groups coincide with those of \cite{Knus1991} and \cite{Auel2009}.
	For the experts we remark that this allows us to ignore distinctions between the ``naive orthogonal group'' and Knus' ``fancy orthogonal group'' \cite{Knus1991}.
	
	We also record that $NC_0(M,q)$ in \cite{McInroy2016} is $\GSpin$.
\end{remark}

\subsection{Clifford $\GL_2$ and $\SL_2$}\label{sec:clifford-bianchi-groups}
Let $C = \Clf(W,q)$ for $(W,q)$, a quadratic space over a commutative ring $R$ where $W$ is a projective $R$-module. 

\begin{definition}
Let $g = \begin{psmallmatrix} a & b  \\ c & d \end{psmallmatrix}$ be a 2 by 2 matrix with entries in $C$. 
The \emph{pseudodeterminant} is defined to be $\Delta(g) = ad^* - bc^*.$
\end{definition}

We now define the Clifford version of $\GL_2$.
\begin{definition}\label{def:gl2}We define the \emph{Clifford general linear group}
	$\GL_2(C)$ to be the set of matrices $ \begin{psmallmatrix}
	 a & b \\ 
	 c & d  
	 \end{psmallmatrix}$
	where 
	\begin{enumerate}
		\item \label{item:det-cond} $ad^*-bc^* \in R^\times$.
	\item \label{item:abcd-cond} $ab^* = ba^*$ and $cd^* = dc^*$.
		\item \label{item:norm-cond} $a\overline{a},b\overline{b},c\overline{c},d\overline{d} \in R$.
		\item \label{item:ratios} $a\overline{c},b\overline{d} \in \Vec(C)$.

		\item \label{item:five} if $x \in \Vec(C)$,
then  $ax\overline{b}+b\overline{x}\,\overline{a},cx\overline{d}+d\overline{x}\,\overline{c} \in R$. 
        \item \label{item:six} if $x \in \Vec(C)$, then $ ax\overline{d} +b \overline{x}\, \overline{c} \in \Vec(C).$
	\end{enumerate}
	The \emph{Clifford special linear group} $\SL_2(C)$ is defined to be $\{g \in \GL_2(C): \Delta(g) = 1\}$.
\end{definition}
After the proof of Theorem~\ref{lem:pre-exceptional}
we will justify the name by showing that $\GL_2(C)$ and $\SL_2(C)$ are
in fact groups.  This long definition is needed to work in the generality where $R$ is a commutative ring with no other hypothesis.  This was first used for fields in \cite[Definition 3.1]{Elstrodt1987} and our observation is that this is the correct definition for repairing the integral spin isomorphism of \cite{McInroy2016}.

\begin{lemma}\label{lem:int-domain-module}
  Let $R$ be an integral domain and $M$ a free $R$-module with a direct sum
  decomposition $M = P \oplus Q$.  Suppose that $m \in M, r \in R$ are such
  that $rm \in P$.  Then $m \in P$ or $r = 0$.
\end{lemma}

\begin{proof}
  Let $m = p + q$ with $p \in P$ and $q \in Q$.  Then $rm = rp + rq$, so
  $rm \in P$ if and only if $rq \in P$.  But $rq \in Q$, so this is equivalent
  to $rq = 0$.  Because $R$ is an integral domain, this is equivalent to
  $r = 0$ or $q = 0$; the conclusion follows, since $q = 0$ if and only if
  $m \in P$.
\end{proof}

\begin{lemma} 
  Let $R$ be an integral domain of characteristic $\neq 2$, $W$ be a free $R$-module, and $(W,q)$ a strongly anisotropic quadratic space.  Then the Clifford monoid $C^\mon$ is closed under transposition.
\end{lemma}
\begin{proof}  We let $K$ be the fraction field of $R$ and we consider the quadratic form $(W_K,q_K) = (W\otimes_RK),q_K)$.  
  Let $x \in \Clf(q)^\mon\setminus \{0\}$. By strong anisotropy $\nrd(x)\neq0$, so $x$ is a unit of $\Clf(q_K)$ and lies in $\Clf(q_K)^\times$.
  We use the fact (Lemma~\ref{lem:generated-by-clifford-vectors})
  that the Clifford group is generated by Clifford vectors.  
	
	So there are vectors $v_1,\ldots,v_m \in
  \Vec(q_K)\setminus \{0\}$ such that $x = v_1\cdots v_m$.  Each $v_i^*\in \Vec(q_K)\setminus \{0\}$ and so $x^* =v_m^*\cdots v_1^* \in \Clf(q_K)^{\times}$; that is, $x^*x' = \nrd(x^*) \in K$ and $x^*\Vec(q_K)x \subseteq \Vec(q_K)$.
  Since $W$ is free, $R = R\gamma_\emptyset$ and $\Vec(q)$ are direct summands of $\Clf(q)$, so $K\cap\Clf(q) = R$ and $\Vec(q_K)\cap\Clf(q)=\Vec(q)$.
  As $x^*\in\Clf(q)$, we get $x^*x'\in R$ and $x^*\Vec(q)x\subseteq\Vec(q)$, i.e. $x^*\in\Clf(q)^\mon$.
  \end{proof}

\begin{theorem}\label{thm:gl2-short}

	Let $R$ be an integral domain of characteristic $\neq 2$, let $W$ be a free $R$-module, and $(W,q)$ a strongly anisotropic quadratic space.  Then $\GL_2(C)$ is described by the formula 
	\begin{equation}\label{eqn:simple-defn}
	\GL_2(C) = \left \lbrace
	\begin{pmatrix}
	a & b \\
	c & d 
	\end{pmatrix}
	\colon a,b,c,d \in C^\mon,  ad^*-bc^* \in R^{\times}, ab^*,dc^* \in \Vec(C)\right \rbrace.
	\end{equation}
        The matrix $\frac{1}{\Delta}\begin{psmallmatrix}d^*&-b^*\\-c^*&a^*\end{psmallmatrix}$ is a two-sided inverse of $\begin{psmallmatrix}a&b\\c&d\end{psmallmatrix}$.
\end{theorem}

\begin{remark}\label{rem:mcinroy} Our proof is inspired by
  \cite[Thm.~6.1 (1)]{McInroy2016}, with some added details and corrections.
  In particular, the proof of \cite[Thm.~6.1 (1)]{McInroy2016} accidentally assumes
  that for $a \in C$ we have $\nrd(a) = 0$ if and only if $a = 0$, which
  requires some additional hypothesis, for example 
  that the form $(W,q)$ is strongly anisotropic.  
  Its definition also accidentally assumes that entries of the paravector version of $\SV_2$ for $\Clf(q_R)$ of a general ring $R$ (the analog we are interested in) must have invertible entries. 
  This precludes, for example, translation matrices like $\begin{psmallmatrix} 1 & 2 \\ 0 & 1 \end{psmallmatrix}$ being elements of this group.
  
  We also note that Theorem~\ref{thm:gl2-short} 
  was previously proved in \cite[Prop.~3.7]{Elstrodt1987} in the case where
  $R$ is a field of characteristic not equal to $2$.  
  The definition in this particular form in the case of $\Clf(q) = \CC_n$ comes from earlier work of Ahlfors (see \cite{Ahlfors1984}).
\end{remark}

\begin{proof}
  We begin by assuming that $a,b,c,d$ satisfy the conditions in equation~\eqref{eqn:simple-defn}.  
  The conditions \eqref{item:det-cond} and \eqref{item:norm-cond} of Definition~\ref{def:gl2} are immediate. Since $ab^*,cd^* \in \Vec(C)$ we have that $ab^* = (ab^*)^* = ba^*$ and $cd^* =dc^*$ giving \eqref{item:abcd-cond}. 

  We begin by showing that $\overline{a}b,\overline{c}d \in \Vec(C)$.  If $b =0$ then clearly $\overline{a}b \in \Vec(C)$, so we assume that $b \neq 0$.
  Note that $ab^* \in \Vec(C)$ so $b'\overline{a} \in \Vec(C)$.
  So $b^*(b'\overline{a})b = \nrd(b^*)\overline{a}b$.  Since $\nrd(b^*) \neq 0$ by Proposition~\ref{prop:big-form-02}(3),
  it follows from Lemma~\ref{lem:int-domain-module}
  that $\overline{a}b \in \Vec(C)$.

Similarly, $\overline{d}c \in \Vec(C)$.
So we have established that
\begin{equation}\label{eqn:small} a,b,c,d \mbox{ satisfy \eqref{eqn:simple-defn}}
  \implies \overline{a}b ,\overline{d}c \in \Vec(C), \end{equation}
and it follows that $\overline{b}a,\overline{c}d \in \Vec(C)$.

  Next, we show that $a^*c,b^*d \in \Vec(C)$.  Since $C^\mon$ is closed under transpose, we know that $\nrd(c^*) \in R$.
  Let $X = \begin{psmallmatrix} a & b \\ c & d \end{psmallmatrix}$ be
  a matrix satisfying the conditions in~\eqref{eqn:simple-defn} and let
  $\Delta = ad^*-bc^*$.  Since $\Delta \in R^{\times}$,
  it follows that $ad^*-bc^* = d'\overline{a}-c'\overline{b}$.  We multiply by $c^*$ on the left and $a$ on the right to obtain
$$c^*\Delta a = c^*(d'\overline{a}-c'\overline{b})a = c^*d'a\overline{a} - c^*c'\overline{b}a = (\overline{d}c)^*\nrd(a) - \nrd(c^*)\overline{b}a.$$
Since $\overline{b}a,\overline{d}c \in \Vec(C)$ we get $\Delta c^*a \in \Vec(C)$ and so
  $c^*a \in \Vec(C)$, by Lemma~\ref{lem:int-domain-module}.
  Hence $a^*c \in \Vec(C)$ and similarly $b^*d \in \Vec(C)$.

  Now we show that 
  $ Y = \frac{1}{\Delta} \begin{psmallmatrix}
    d^* & - b^* \\ -c^* & a^* \end{psmallmatrix}$
  is a two-sided inverse of $X$.
  Showing that $XY = 1$ is a simple calculation, so we consider
  $$U  = \Delta YX =
  \begin{pmatrix} d^*a - b^*c & d^*b-b^*d \\ a^*c - c^*a & a^*d-c^*b \end{pmatrix}.$$
  
  We first note that $a^*c, d^*b \in \Vec(C)$ so $U$ is diagonal.
  If $d^*a -b^*c \in R^{\times}$ we have that
  $a^*d - c^*b \in R^{\times}.$  Hence it remains to show that
  $d^*a -b^*c = ad^*-bc^*$.
We use that $ab^*, a^*c \in \Vec(C)$ and $\Delta \in R$ in the following calculation; since a Clifford vector is fixed by $*$, we have $ab^*=ba^*$ and $a^*c=c^*a$:
  \begin{eqnarray*}
    \overline{a}a (d^*a-b^*c) & = & \overline{a} (ad^*)a - \overline{a}ab^*c \\
    & = & \overline{a} (ad^*)a - \overline{a}ba^*c \\
    & = & \overline{a} (ad^*)a - \overline{a}bc^*a \\
    & = & \overline{a} (ad^*- bc^*)a \\
  & = & \nrd(a) \Delta \end{eqnarray*}
So if $\nrd(a) \neq 0$ then $Y$ is an inverse for $X$.  If $\nrd(a) =0$, then $a=0$, so $\Delta = -bc^* \in R^\times$.
  Hence $b$ has the right inverse $-\Delta^{-1}c^*$, which is a two-sided inverse because $C$ is a finitely generated $R$-module \cite[Tag 05G8]{stacks-project}.
  So $c = -\Delta (b^*)^{-1}$ and $b^*c = -\Delta = bc^*$, hence $d^*a-b^*c=\Delta$ and $Y$ is an inverse for $X$.
Next, we note that $Y \in \GL_2(C).$  Since $C^\mon$ is closed under transpose, we have that $a^*,b^*,c^*,d^* \in C^\mon$.  Since $d^*a-b^*c =ad^*-bc^* \in R^\times$ we see that the pseudo-determinant of $Y$ is in $R^\times$.  The last condition in the statement is that $d^*b, a^*c \in \Vec(C)$ which we have already established.  Hence we conclude that $\GL_2(C)$ is closed under inverses.

If we apply equation~\eqref{eqn:small} to $Y$ we see that
$\overline{d}^*b^*=d'b^*\in\Vec(C)$, hence $b\overline{d}=(d'b^*)^*\in\Vec(C)$, and similarly $a\overline{c} \in \Vec(C)$, establishing item~\eqref{item:ratios}.

We now check \eqref{item:five}, that $ax\overline{b}+b\overline{x}\,\overline{a} \in R$ for all $x \in \Vec(C)$. 
The condition $cx\overline{d}+d\overline{x}\,\overline{c} \in R$ is proved similarly. 
First suppose that $u,v  \in \Vec(C)$ with $u = r+p$ and $v=s+q$ with $r,s \in R$ and $p,q \in \Im \Vec(C)$.  Then $w=uv$ satisfies
  \begin{multline} \label{eqn:yybarinR} w+\overline{w} = (r+p)(s+q)+\overline{(r+p)(s+q)}
=  (r+p)(s+q)+\overline{(s+q)}\,\overline{(r+p)}
   \\ =  rs+sp+rq+pq+rs-sp-rq+qp \in R.  \end{multline}
  Let $y = ax \overline{b}$. 
  If $b=0$ we are done. So consider
  $$\nrd(b^*)ax\overline{b}=ab^*b'x\overline{b}=(ab^*)(bx'b^*)'$$
  Let $u = ab^*$ and $v=(bx'b^*)'$, so  $\nrd(b^*)y = uv$ is the product of two elements of $\Vec(C)$ and so
  $\nrd(b^*)(y+\overline{y}) \in R$, giving condition \eqref{item:five}.

  Finally we establish condition \eqref{item:six}.
  First, let us assume that $a=0$.  In this case we have that
  $ad^*-bc^* = -bc^* \in R^\times$.  We also know that $b \neq 0$ and $b^* \in C^\mon$ so $0 \neq \nrd(b^*) \in R$.  Let $x \in \Vec{C}$ and consider
  $$ \nrd(b^*) b \overline{x}\, \overline{c} =
  b \overline{x} b^* b'\overline{c} = (b \overline{x} b^*) (bc^*)' \in \Vec(C)$$
  So $ax\overline{d}+b\overline{x}\,\overline{c} \in \Vec{C}$ establishing \eqref{item:six} when $a=0$.
 Now we assume $x \in \Vec(C)$ and we will show that $axa^* \in \Vec(C)$ if and only if $ax\overline{d}+b\overline{x}\,\overline{c} \in \Vec(C)$ when $a \neq 0$.  Note that
  \begin{eqnarray*} \nrd(a^*)(ax\overline{d}+b\overline{x}\,\overline{c})
  & = & axa^*a'\overline{d} + \nrd(a^*)b\overline{x}\,\overline{c}\\
    & = & axa^*(ad^*-bc^*)'+axa^*b'\overline{c}+\nrd(a^*)b\overline{x}\,\overline{c}.\end{eqnarray*}
    If $b=0$ we are done, so now we assume $b \neq 0$.  Since $x, a^*b' \in \Vec(C)$ we have that $r= x(a^*b') + \overline{a^*b'}\overline{x} \in R$ by~\eqref{eqn:yybarinR}.
    So $x(a^*b') = r - b^*a'\overline{x}.$
    Continuing the calculation above yields
   \begin{eqnarray*} \nrd(a^*)(ax\overline{d}+b\overline{x}\,\overline{c})
     & = & axa^*(ad^*-bc^*)'+a(r-b^*a'\overline{x})\overline{c}+\nrd(a^*)b\overline{x}\,\overline{c}\\
     & = & axa^*(ad^*-bc^*)'+ra\overline{c}-ab^*a'\overline{x}\,\overline{c}
     +a^*a'b\overline{x}\,\overline{c} \\
     & = & axa^*(ad^*-bc^*)'+ra\overline{c}-ab^*a'\overline{x}\,\overline{c}
     +ba^*a'\overline{x}\,\overline{c} \\
     & = & axa^*(ad^*-bc^*)'+ra\overline{c}
     ,\end{eqnarray*}
   where we used $ab^* = ba^*$ and $\nrd(a^*) \in R$ since $ab^* \in \Vec{C}$ and $a^* \in C^\mon$.
   Since $a\overline{c} \in \Vec{C}$, we now have that $axa^* \in \Vec(C)$
   if and only if  $ax\overline{d}+b\overline{x}\,\overline{c} \in \Vec(C)$ when $a \neq 0$
   and we have shown that the short definition in the theorem implies the long definition of Definition~\ref{def:gl2}.
We will use this contrapositive of this step to help establish the converse.
  
  Conversely, let $X=\begin{psmallmatrix} a & b \\ c & d \end{psmallmatrix}\in\GL_2(C)$ and $\Delta = ad^*-bc^*$. 
  We need to show that $a,b,c,d\in C^{\mon}$, $\Delta\in R^{\times}$ and that $ab^*, cd^* \in \Vec(C)$.
  The condition $\Delta\in R^\times$ is \eqref{item:det-cond}.

  We first use that $\GL_2(C)$ is a group (Corollary~\ref{lem:gl-is-group}).
  By \eqref{item:det-cond} and \eqref{item:abcd-cond} we have $XY=1$ for $Y=\frac{1}{\Delta}\begin{psmallmatrix} d^* & -b^* \\ -c^* & a^* \end{psmallmatrix}$, so $Y=X^{-1}\in\GL_2(C)$.
  Condition~\eqref{item:norm-cond} for $Y$ gives $\nrd(a^*),\nrd(b^*),\nrd(c^*),\nrd(d^*)\in R$, and condition~\eqref{item:ratios} for $Y$ gives $d^*c',b^*a' \in \Vec(C)$.
  Applying $*$, and then $'$, both of which preserve $\Vec(C)$, we obtain $\overline{d}c,\overline{a}b \in \Vec(C)$ and $a^*b'\in\Vec(C)$.

  Next we show that $a\in C^\mon$. By \eqref{item:norm-cond}, $\nrd(a)\in R$. Let $x\in\Vec(C)$; we show that $axa^*\in\Vec(C)$, which is clear if $a=0$.
  If $a\neq 0$, the computation in the proof of \eqref{item:six} above used only $ab^*=ba^*$, $\nrd(a^*)\in R$ and $a^*b'\in\Vec(C)$, which we now have, and it gives
  $$\nrd(a^*)(ax\overline{d}+b\overline{x}\,\overline{c}) = axa^*\Delta+ra\overline{c}$$
  with $r\in R$ (and $r=0$ if $b=0$). By \eqref{item:six} and \eqref{item:ratios}, and since $\Delta\in R^\times$, we get $axa^*\in\Vec(C)$. So $a\in C^\mon$.
  Since $S=\begin{psmallmatrix} 0 & -1 \\ 1 & 0 \end{psmallmatrix}\in\GL_2(C)$, the matrices $XS$, $SX$ and $SXS$ lie in $\GL_2(C)$; their upper left entries are $b$, $-c$ and $-d$, so $b,c,d\in C^\mon$ by the same argument.

  Lastly, since $d\in C^\mon$ and $\overline{d}c\in\Vec(C)$, we have $d(\overline{d}c)d^* = \nrd(d)cd^* \in \Vec(C)$.
  If $d\neq 0$ then $\nrd(d)\neq 0$, so $cd^* \in \Vec(C)$ by Lemma~\ref{lem:int-domain-module}; if $d=0$ this is clear.
  Similarly $a(\overline{a}b)a^*=\nrd(a)ba^*\in\Vec(C)$ gives $ba^*\in\Vec(C)$, hence $ab^*=(ba^*)^*\in\Vec(C)$.
 
\end{proof}

\begin{remark}
		Let $A$ be an order in a Clifford algebra $C$.  
		Suppose that $A$ is not closed under involutions.
		One can look at matrices with entries in $A$ satisfying Equation \eqref{eqn:simple-defn}.
		This is a monoid but it is not clearly a group. 
		This is the reason we assume our orders are closed under the involutions.
\end{remark}

\begin{corollary}\label{cor:stars-match}
 If $\begin{psmallmatrix}a&b\\c&d \end{psmallmatrix} \in \SL_2(\CC_n)$
 then $c^*a=a^*c,b^*d = d^*b, ab^*=ba^*,cd^*=dc^*$. 
\end{corollary}
\begin{proof} This follows from the formula for the inverse in
  Theorem~\ref{thm:gl2-short} and in particular from the fact that it is a
  two-sided inverse.
\end{proof}

\subsection{Factoring Clifford Algebras}\label{sec:factoring}

In this section we will prove the Decomposition Lemma
(Lemma~\ref{lem:decomposition}), showing that if $(W,q)$ is a quadratic space
with a submodule $U$ of degree $2$ with a complement $V$,
then under certain conditions the
Clifford algebra $\Clf(q)$ is isomorphic to
$\Clf(q\vert_U) \otimes \Clf(\delta q\vert_V)$, where $\delta = \Disc(U,q\vert_U)$.  In order to state and prove our
result, we first recall the definition of the discriminant. 
\begin{definition}
  Suppose that $(W,q)$ is a quadratic space where $W$ is a
  free module of finite rank $n$, choose a basis
  $\{e_i\}$, and let $\phi(u,v)=q(u+v)-q(u)-q(v)$ be the polar form, so that $\phi(e_i,e_i)=2q(e_i)$.
  The {\em discriminant} of $q$ is the element
  $\Disc q \in R/{R^\times}^2$ represented by $(-1)^{n(n-1)/2}\det \phi(e_i,e_j)$ if $n$ is even,
  and by $(-1)^{n(n-1)/2}\frac12\det \phi(e_i,e_j)$ if $n$ is odd.
\end{definition}
For $n$ odd, $\det\phi(e_i,e_j)$ is $2$ times a polynomial with integer coefficients in the $q(e_i)$ and the
$\phi(e_i,e_j)$ ($i<j$), and $\frac12\det\phi(e_i,e_j)$ denotes the value of that polynomial.
A change of basis by $P\in\GL_n(R)$ multiplies $\det\phi(e_i,e_j)$ by $\det(P)^2$, so $\Disc q$ is well defined.
For example, if $U$ has basis $e_1,e_2$, then $\Disc(U,q\vert_U) = \phi(e_1,e_2)^2-4q(e_1)q(e_2)$.

Now we give the definition of, and a basic fact about,
complements in quadratic spaces.

\begin{definition}\label{def:orth-complement}
  Let $(W,q)$ be a quadratic space and let $U$ be a submodule of $W$.
  The {\em orthogonal complement} $U^\perp$ of $U$ in $W$ is the set
  $\{w \in W: q(u+w) = q(u) + q(w) \mbox{ for all } u \in U\}$.
\end{definition}

\begin{proposition}\label{prop:complement-is-sub}
  The orthogonal complement of $U$ is a submodule of $W$.
\end{proposition}

\begin{proof}
  This follows straightforwardly from Definition~\ref{def:quadratic-form}.
  Indeed, suppose that $v, v' \in U^\perp$.
  \begin{samepage}
  For all $u \in U$ we have
  \begin{align*}
    q(u+v+v') &= q(u+v)+q(u+v')+q(v+v')-q(u)-q(v)-q(v')\\
    &= q(u)+q(v)+q(u)+q(v')+q(v+v')-q(u)-q(v)-q(v')\\
    &= q(u)+q(v+v'),\\
  \end{align*}
  which shows that $v+v' \in U^\perp$.
  \end{samepage}
  Similarly, for $t \in R$ we have
  $$q(tv+u) = tq(v+u)+t^2q(v)+q(u)-tq(u)-tq(v) = t^2q(v)+q(u) = q(tv)+q(u)$$
  and it follows that $tv \in U^\perp$ as well.
\end{proof}

\begin{definition}\label{def:decomp}
  Suppose that $W$ is written as a direct sum $U \oplus V$ where
  $V \subseteq U^\perp$.  We say that $W = U \oplus V$ is a
  {\em decomposition} of $W$.
\end{definition}

The following is an integral version of a well-known Lemma over fields for which we could not find an adequate reference. See, for example, Lemma 2.8 of \cite{Conrad2012} for the statement over fields.
\begin{lemma}[Decomposition Lemma]\label{lem:decomposition}
  Let $(W,q)$ be a quadratic space over $R$.  Suppose given a decomposition
  $W = U \oplus V$ where $U$ is free of rank $2$ and $V$ is finitely generated
  and such that $q(u+v) = q(u)+q(v)$ for all $u \in U, v \in V$.  If
  $\delta = \Disc(U,q\vert_U)$ is invertible, then
  $$\Clf(q) \cong \Clf(q\vert_U) \otimes_R \Clf(\delta q\vert_V).$$
\end{lemma}
\begin{proof}
  We fix a basis $e_1, e_2$ for $U$ and use the same notation for the
  generators of $\Clf(q\vert_U)$ and $\Clf(q)$.
  Given $v \in V$, let $f_v, g_v$ be the
  corresponding elements of $\Clf(\delta q\vert_V)$ and $\Clf(q)$.
  Let $a = q(e_1), c = q(e_2), b = q(e_1+e_2)-a-c$, so that $e_1e_2+e_2e_1=-b$ and $\delta = b^2-4ac$, and let $d = 2e_1e_2+b = e_1e_2-e_2e_1$.
  We show that the map
  $\phi: \Clf(q\vert_U) \otimes_R \Clf(\delta q\vert_V) \to \Clf(q)$
  given by $\phi(e_i \otimes 1) = e_i$ for $i = 1, 2$ and
  $\phi(1 \otimes f_v) = d g_v$ for $v \in V$ is an isomorphism.  There are
  three things to check.

  First we show that the relations among the generators of
  $\Clf(q \vert_U) \otimes \Clf(\delta q \vert_V)$ are satisfied by their
  images in $\Clf(q)$.  For the $e_i \otimes 1$, this is clear.
  The element $d$ anticommutes with $e_1$: we have
  $de_1 = e_1e_2e_1 - e_2e_1e_1 = e_1e_2e_1 + ae_2$ and
  $e_1d = e_1e_1e_2 - e_1e_2e_1 = -ae_2 - e_1e_2e_1$;
  similarly $d$ anticommutes with $e_2$.
  Since $g_v$ anticommutes with $e_1$ and $e_2$, it follows that
  $\phi(1\otimes f_v) = dg_v$ commutes with $e_1$ and $e_2$, and that $d$ commutes with $g_v$.
  Using $(e_1e_2)^2 = e_1(-b-e_1e_2)e_2 = -be_1e_2-ac$, $(e_2e_1)^2 = -be_2e_1-ac$
  and $e_1e_2e_2e_1 = e_2e_1e_1e_2 = ac$, we get
  $$d^2 = (e_1e_2)^2+(e_2e_1)^2-e_1e_2e_2e_1-e_2e_1e_1e_2 = (b^2-2ac)-2ac = \delta.$$
  Hence
  $$\phi(f_v)\phi(f_{v'})+\phi(f_{v'})\phi(f_v) = d^2(g_vg_{v'}+g_{v'}g_v) = -\delta\bigl(q(v+v')-q(v)-q(v')\bigr),$$
  and in particular $\phi(f_v)^2 = -\delta q(v)$; these are the defining relations of
  $\Clf(\delta q\vert_V)$.  This completes the verification.

  Second, we need to prove that $\phi$ is surjective.  This is clear,
  since $\Clf(q)$ is generated by the $e_i$ and the $g_v$, and $\delta$ is a
  unit, so $d$ is invertible with $d^{-1}=\delta^{-1}d$ and $g_v = \delta^{-1}d\,(dg_v)$ lies in the image.

  Finally, to show that $\phi$ is injective, it suffices to note that an
  inverse is given by $\phi^{-1}(e_i) = e_i \otimes 1$ and
  $\phi^{-1}(g_v) = \delta^{-1}d \otimes f_v$; these elements satisfy the defining relations of
  $\Clf(q)$ because $d$ anticommutes with $e_1,e_2$ and $(\delta^{-1}d\otimes f_v)^2 = \delta^{-2}d^2\otimes f_v^2 = -q(v)$.
\end{proof}

\begin{remark}
  When the nondegenerate quadratic form has odd rank $n$ over a field of characteristic $\neq 2$, the product $J = \gamma_1\cdots \gamma_n$ of an orthogonal basis generates the center. 
\end{remark}

\begin{example}
	Consider $A=\quat{-2,-3,-5}{\QQ}$. 
	Then we have $\quat{-2,-3,-5}{\QQ} \cong \quat{-2,-3}{\QQ}\otimes_{\QQ} \quat{30}{\QQ} \cong \quat{-2,-3}{\QQ}\otimes_{\QQ} \QQ(\sqrt{30}) \cong \quat{-2,-3}{\QQ(\sqrt{30})}.$
	This implies that $A$ is not central over $\QQ$: its center is
        $\QQ(\sqrt{30})$.
	If the basis vectors are $\gamma_1,\gamma_2,\gamma_3$ with $\gamma_1^2=-2, \gamma_2^2 = -3, \gamma_3^2=-5$ with $\gamma_i\gamma_j = -\gamma_j\gamma_i$ for $i\neq j$, one can check that the central element $J:= \gamma_1\gamma_2\gamma_3$ satisfies $J^2 = 30$.
	This all has to do with the rank of the quadratic form being odd. Note that since the rank of $A$ is not a square, $A$ cannot be a central simple algebra over $\Q$.
\end{example}

\begin{corollary}\label{cor:clifford-tensor-product}
	Let $(W,q)$ be a nondegenerate quadratic space of dimension $n$ over a field $K$ of characteristic $\neq 2$, let $C=\Clf(W,q)$, and let $J=\gamma_1\cdots\gamma_n$ for an orthogonal basis $\gamma_1,\ldots,\gamma_n$ of $W$.
	\begin{enumerate}
		\item If $n$ is odd, then $Z(C) = K[J]\cong K[t]/(t^2+\Disc q)$, $C\cong C_+\otimes_K Z(C)$, and $C_+$ is a tensor product of $(n-1)/2$ quaternion algebras over $K$, hence central simple over $K$. Thus $C$ is central simple over the field $Z(C)$ if $-\Disc q$ is not a square in $K$, and $C\cong C_+\times C_+$ otherwise.

\item If $n\geq 2$ is even, 
then $C$ is a tensor product of $n/2$ quaternion algebras over $K$, hence central simple over $K$, and $Z(C_+)=K[J]\cong K[t]/(t^2-\Disc q)$; $C_+$ is central simple over the field $Z(C_+)$ if $\Disc q$ is not a square in $K$, and a product of two central simple $K$-algebras otherwise.
	\end{enumerate}
\end{corollary}
\begin{proof}
	For orthogonal $\gamma_i,\gamma_j$ we have $\Disc(\langle\gamma_i,\gamma_j\rangle,q)=-4q(\gamma_i)q(\gamma_j)\neq 0$, so repeated application of Lemma~\ref{lem:decomposition} writes $C$ as a tensor product of $\lfloor n/2\rfloor$ quaternion algebras, tensored, if $n$ is odd, with $\Clf(\langle c\rangle)\cong K[t]/(t^2+c)$ for some $c\in K^\times$.
	For $C_+$ use the isomorphism $\Clf\bigl(q(\gamma_n)\,q\vert_{\langle\gamma_1,\ldots,\gamma_{n-1}\rangle}\bigr)\to C_+$, $\gamma_i\mapsto\gamma_i\gamma_n$, which maps a basis to a basis.
	Finally, $J$ is central in $C$ (if $n$ is odd) or in $C_+$ (if $n$ is even), and $J^2=(-1)^{n(n+1)/2}q(\gamma_1)\cdots q(\gamma_n)$, which is $-\Disc q$ (if $n$ is odd) or $\Disc q$ (if $n$ is even) modulo $(K^\times)^2$.
\end{proof}

\subsection{Orthogonal Representations of Clifford-Bianchi Groups}\label{sec:naive}
Recall that there is a representation of $\SL_2(\CC)$ into the Lorentz transformations $\O_{1,3}(\RR)^{\circ}$ induced by acting on the augmented Pauli matrices by conjugation.
In order to get a map $\SL_2(\CC_n) \to \O_{1,n+1}(\RR)$ we generalize this construction in a naive fashion.

Let $q = \sum_{j=1}^{n-1} d_j y_j^2$ be a positive definite quadratic form in $n-1$ variables with $d_j$ positive squarefree integers. 
Let $\Ocal = \Clf(\ZZ^{n-1},q)$ and let $K = \Ocal\otimes_{\ZZ}\QQ$.
Let $\gamma_j$ be the generators of $\Ocal$ for $1\leq j \leq n-1$ and embed $\Ocal$ into $\CC_n$ using $\gamma_j = \sqrt{d_j}i_j$ for $1 \leq j \leq n-1$.
The Clifford-Hermitian matrices $\tau_j$ defined by
$$ \tau_{n+1} = \begin{pmatrix} 1 & 0 \\ 0 & 1 
\end{pmatrix}, \quad \tau_{0} = \begin{pmatrix} 0 & 1 \\ 1 & 0 
\end{pmatrix}, \quad \tau_n = \begin{pmatrix} 1 & 0 \\0 & -1 
\end{pmatrix}, \quad \tau_j = \begin{pmatrix} 0 & -\gamma_j \\ \gamma_j & 0 
\end{pmatrix}, \qquad j=1,\ldots, n-1
$$
where $\gamma_j = \sqrt{d_j}i_j$, will be called the \emph{Pauli matrices}.
The $\RR$-span of these matrices is the collection of Clifford-Hermitian matrices $M_2(\CC_n)_{\herm}$ (this will be rigorously defined in \S\ref{sec:hermitian}).
In the case $n=2$ and $d_1=1$, so that $\gamma_1 = i = \sqrt{-1}$, the matrices $\tau_0,\tau_1,\tau_2$ are the classical Pauli matrices $\sigma_1,\sigma_2,\sigma_3$, and $\tau_3$ is the identity matrix.

We now generalize the well-known representation of $\SL_2(\CC)$ into the group of Lorentz transformations. 
A general element of $M_2(\CC_n)_{\herm}$ can be written as $Y = \sum_{j=0}^{n+1} y_j \tau_j$ where $y_j \in \RR$. 
We have $ Y = \begin{psmallmatrix}
y_{n+1} + y_n & \overline{y} \\
y & y_{n+1}-y_n
\end{psmallmatrix}$
where $y=y_0+y_1 \gamma_1 + \cdots + y_{n-1} \gamma_{n-1}$.
Using the determinant of Clifford-Hermitian matrices from \S\ref{sec:hermitian}, we find
$$\det(Y) = y_{n+1}^2 - y_n^2-y_0^2 - d_1 y_1^2 - \cdots - d_{n-1}y_{n-1}^2.$$
This is a new quadratic form $Q = y_{n+1}^2 -y_n^2-y_0^2-q$.
Up to sign it is the form of Theorem~\ref{thm:mcinroy-correction}: on the lattice $L$ below, $-Q=b\overline{b}-ac$, which is $x^2-yz+q$ with $x=b_0$, $y=a$ and $z=c$; replacing $Q$ by $-Q$ changes neither $\SO_Q$ nor $\Spin_Q$.
Let $L\subset M_2(\CC_n)_{\herm}$ be the lattice of matrices $\begin{psmallmatrix} a & \overline{b}\\ b & c\end{psmallmatrix}$ with $a,c\in\ZZ$ and $b\in\Vec(\Ocal)$; on $L$ we have $Q = ac-b\overline{b}$, an integral quadratic form, and we write $\SO_Q(\ZZ)$ for the group of determinant-one isometries of $(L,Q)$.
This gives us an action of $\SL_2(\Clf(q))$ on $M_2(\CC_n)_{\herm}$ via $g\cdot Y = gYg^{\dag}$, where $g^{\dag}$ is the Clifford adjoint of \S\ref{sec:hermitian}, and this gives rise to a representation.
The map given by $(g,Y) \mapsto gYg^{\dag}$ gives a group homomorphism 
\begin{equation}\label{eqn:lorentz-rep}
\varphi:\SL_2(\Clf_q(\ZZ)) \to \SO_Q(\ZZ),
\end{equation}
which generalizes the famous representation of $\SL_2(\CC)$ into the Lorentz group using the $n=2$ Pauli matrices.
It is a homomorphism because $(gh)^{\dag}=h^{\dag}g^{\dag}$; it preserves $L$ because $g$ satisfies items \eqref{item:norm-cond}--\eqref{item:six} of Definition~\ref{def:gl2}; and it preserves $Q$ because
$$Q(\rho_g y)=\det(gYg^{\dag}) = \det(Y)= Q(y)$$
for $g\in\SL_2(\CC_n)$ \cite[\S 3]{Elstrodt1988}. The image has determinant one because $\SL_2(\CC_n)$ acts on $\Hcal^{n+1}$ by orientation-preserving isometries (\S\ref{sec:action-hermitian}).
One can now ask the following question.
\begin{question}
	Is the image of $\varphi$ from \eqref{eqn:lorentz-rep} surjective? Is the image an arithmetic group?
\end{question}
Analyzing this map by hand is extremely difficult. 
We set up an exact sequence of group schemes in the fppf topology to attack this problem.

\subsection{Arithmetic Bott Periodicity}\label{sec:arithmetic-bott}

This section gives a treatment of periodicity for Clifford algebras, which will be used in our study of arithmetic groups. 
We need to give an integral version of the statement that 
\begin{equation}\label{eqn:real-bott}
C_{p,q} \cong (C_{p+1,q})_+
\end{equation}
i.e., that basic real Clifford algebras can be related to even parts of the real Clifford algebra. 
In the special case of Clifford algebras over $\RR$, this theorem can be found in \cite[Chapter 7]{Porteous1995}.
In our application, we then have $M_2(\CC_n) \cong C_{n,1} \cong (C_{n+1,1})_+.$
Under this isomorphism we have $\SL_2(\CC_n) \cong \Spin_{1,n+1}(\RR)$.
One of our goals is to strengthen this to an isomorphism of $\ZZ$-group schemes, and this section builds the arithmetic periodicity variant of real Bott periodicity \eqref{eqn:real-bott} (see also Figure~\ref{fig:bott-clock}) to do this.
\begin{figure}[htbp!]
\centering
\begin{tikzpicture}[lab,scale=1.15]
 \draw[gray!60] (0,0) circle (2.3); \draw[gray!40] (0,0) circle (3.5);
 \foreach \n/\A/\X [count=\k from 0] in {1/\RR/{\O}, 2/\CC/{\O/\U}, 3/\HH/{\U/\Sp}, 4/{\HH\oplus\HH}/{\BSp\times\ZZ}, 5/{\HH(2)}/{\Sp}, 6/{\CC(4)}/{\Sp/\U}, 7/{\RR(8)}/{\U/\O}, 8/{\RR(8)\oplus\RR(8)}/{\BO\times\ZZ}}{
   \pgfmathsetmacro\ang{90-45*\k}
   \node[fill=white,inner sep=1pt] at (\ang:2.3) {$\CC_{\n}=\A$};
   \node[domC!70!black,fill=white,inner sep=1pt] at (\ang:3.5) {$\X$};}
 \draw[-{Stealth},thick,bubbleC] (60:0.9) arc[start angle=60,end angle=-240,radius=0.9];
 \node[bubbleC,align=center] at (0,0){$\CC_{n+8}$\\$=\CC_n(16)$};
\end{tikzpicture}
\caption{Bott periodicity as a clock, from Table~\ref{tab:bott-periodicity} in Appendix~\ref{A:classical-bott}. Inner ring: $\CC_n$ for $n\equiv 1,\dots,8\pmod 8$, with $R(m)=M_m(R)$. Outer ring (orange): the spaces $X_n^\infty$; at $n=1$ the table has no entry, and $\O$ is the value at $n=9$. One turn of the clock is $\CC_{n+8}=\CC_n(16)$.}\label{fig:bott-clock}
\end{figure}

Let $q$ be a quadratic form over a ring $R$ in the variables $x_1,\ldots,x_n$, and consider the form $q+x^2$, which is an orthogonal direct sum.  We will write the generators of $\Clf(q)$ as $e_1,\ldots,e_n$ and the generators of $\Clf(q+x^2)$ as
$e_1,\ldots,e_n,e$. Write $\overline{\otimes}$ for the graded tensor product.
Note first that $\Clf(q)$ is a subalgebra of $\Clf(q+x^2)$ under the inclusion
$e_i \mapsto e_i$, which is compatible with the natural inclusion
$\Clf(q) \simeq \Clf(q) \otimes 1 \subseteq \Clf(q) \overline{\otimes} \Clf(x^2).$

Similarly, consider the orthogonal direct sum $q-yz$, and let the generators
of $\Clf(\ZZ^{n+2},q-yz)$ be $e_1,\ldots,e_n,f,g$.  As before we identify
the $e_i \in \Clf(\ZZ^n,q)$ with their images in $\Clf(\ZZ^{n+2},q-yz)$ under
the obvious injective homorphism.  and we similarly identify $f,g$ with their  images from $\Clf(\ZZ^2,-yz).$
We have that
$$e_jf = -fe_j,  \quad e_j g = -ge_j, \quad f^2=g^2 = 0, \quad fg+gf =1.$$

In $\Clf(\ZZ^{n+1},q+x^2)$ we also have $e^2=-1$ and, since the sum is orthogonal, $(e_j+e)^2 = e_j^2+e^2$, so $e_je=-ee_j$.
These identities imply the following rules:
	 $$ xf=fx, \quad xg = gx, \quad xe=ex, \qquad x\in \Clf(q)_+$$
	 $$ yf=-fy, \quad yg=-gy, \quad ye=-ey, \qquad y\in \Clf(q)_-$$
	 $$ zf=fz', \quad zg = gz', \quad ze=ez', \qquad z \in \Clf(q). $$

\begin{definition}
	On $\Clf(q+x^2)_+$ we define the \emph{Satake involutions} $x\mapsto x^s$ and $x\mapsto x^{\alpha}$ by 
	 \begin{equation}\label{eqn:satake-involutions}
	 x^{\alpha} = -exe, \quad x^s = x^{*\alpha} = x^{\alpha *}.
	 \end{equation}
\end{definition}
There are closely related to the Cartan involution on $\operatorname{Lie}(\Spin_Q(\RR))$. 

\begin{lemma}[First Bott Periodicity]\label{lem:bott2}
	Let $v_++v_- \in \Clf(q)$ be a decomposition into graded components.
	Then there is an isomorphism of algebras 
	\begin{equation}\label{eqn:bott2}
	\phi: \Clf(q) \to \Clf(q+x^2_0)_+, \quad \phi(v_++v_-) = v_++v_-e.
	\end{equation}
	Furthermore we have 
	\begin{equation}\label{eqn:phi-involutions}
		\phi(v')=\phi(v^*)^*=\phi(v)^{\alpha}, \qquad \phi(v^*) = \phi(v')^* =\phi(v)^s, \qquad \phi(\overline{v})=\phi(v)^*.
	\end{equation}
\end{lemma}
\begin{proof}

	Let $W$ be the $R$-module associated to $q$, and let $e$ be the generator of $\Clf(q+x_0^2)$ associated to $x_0$, so that $e^2=-1$ and $we=-ew$ for $w\in W$.
	Consider the map $W\to\Clf(q+x_0^2)_+$ given by $w\mapsto we$. Since
	$$ (we)^2 = wewe = -w^2e^2 = -q(w), $$
	Lemma~\ref{lem:universal-property} gives an algebra homomorphism $\phi:\Clf(q)\to\Clf(q+x_0^2)_+$ with $\phi(w)=we$.
	Since $(e_ie)(e_je)=-e_ie_je^2=e_ie_j$, on basis vectors $e_S$ for $S \subset \lbrace 1,\ldots, n\rbrace$ we have
	\begin{equation}\label{eqn:beta2}
	\phi(e_S) =\begin{cases}
	e_S, & \vert S \vert \equiv 0 \mod 2 \\
	e_S e, & \vert S \vert \equiv 1 \mod 2
	\end{cases},
	\end{equation}
	so $\phi(v_++v_-)=v_++v_-e$. The elements on the right of \eqref{eqn:beta2} form a basis of $\Clf(q+x_0^2)_+$, so $\phi$ is bijective.

	The first part of the first identity follows from $\phi(v^*)^* = (v_+^* + v_-^* e)^* = v_++ev_- = v_+-v_-e = \phi(v')$.
	The second part of the first identity follows from
	$\phi(x)^{\alpha}= -e(x_++x_- e)e = -ex_+e -ex_-e^2 = x_++ex_- = x_+-x_- e = \phi(x').$

	The third identity is the first applied to $v^*$, since $\overline{v}=(v^*)'$.
	The second identity follows: $\phi(x')^*=\phi(x^*)$ is the $*$ of the first identity, and $\phi(x)^s = \phi(x)^{*\alpha} = \phi(\overline{x})^{\alpha} = \phi(\overline{x}') = \phi(x^*)$.

\end{proof}

\begin{example}
We have $\CC\cong \HH_+$ since 
$\Clf(\RR,x^2)=\CC \cong \Clf(\RR^2,x^2+y^2)_+ =\HH_+ = \RR[k] $ where $k=ij$. 
\end{example}

Let $q$ be a quadratic form on a free $R$-module of finite rank.
Consider the quadratic form $Q=q-yz$, where $q$ is orthogonal to $-yz$.  We write the generators of $\Clf(-yz)$ as $f,g$, which again is naturally a subalgebra of $\Clf(Q).$
In what follows we use $M_2(\Clf(q)) \simeq  \Clf(q) \overline{\otimes} \Clf(-yz)$.

\begin{definition}
Let $A=\begin{psmallmatrix} a & b \\ c & d \end{psmallmatrix} \in M_2(\Clf(q))$. 
We define the \emph{Clifford adjugate} $\Adj(A) = A^{a}$, the \emph{transported transpose} $A^{\tau}$ and \emph{transported parity involutions} $A^{\sigma}$ by 
\begin{equation}\label{eqn:matrix-involutions-i}
A^{a} = \begin{pmatrix}
d^* & -b^* \\ -c^* & a^* \end{pmatrix}, \qquad A^{\tau} = \begin{pmatrix} \overline{d}& \overline{b} \\ \overline{c} & \overline{a} \end{pmatrix}, \qquad A^{\sigma} = \begin{pmatrix} a' & -b' \\ -c' & d' \end{pmatrix}.
\end{equation}
\begin{equation}\label{eqn:matrix-involutions-ii}
 A^{\underline{\tau}} = \begin{pmatrix} \overline{d} & b^* \\ c^* & \overline{a} \end{pmatrix}, \qquad A^{\underline{\sigma}} = \begin{pmatrix} a & -b' \\ -c' & d \end{pmatrix}.
\end{equation}
These transformations satisfy $A^a= A^{\tau\sigma} = A^{\sigma \tau}$.

\end{definition}
In the theorem below we will see that $A^a$, $A^{\sigma}$ and $A^{\tau}$ are cooked up to correspond to Clifford conjugation, the sign changing transformation $x\mapsto x'$ (parity), and the Clifford transpose $x\mapsto x^*$.

\begin{lemma}[Second Bott Periodicity] \label{thm:bott1}
	There is an isomorphism of $R$-algebras
	\begin{equation}\label{eqn:bott1}
	\iota: M_2(\Clf(q))
	\to \Clf(q-yz), \quad A=\begin{pmatrix} a & b \\ c & d \end{pmatrix}
	\mapsto a fg+ bf + c'g +d'gf,
	\end{equation}
	where $a,b,c,d \in \Clf(q)$. 
	We have
	\begin{equation}\label{eqn:iota-involutions}
	\iota(A^{\tau}) = \iota(A)^*, \qquad \iota(A^{\sigma}) = \iota(A)', \qquad  \iota(\Adj(A)) = \overline{\iota(A)}.
	\end{equation}
\end{lemma}
\begin{proof}
Let $W$ be the module of $q$. The $R$-linear map $W\oplus Rf\oplus Rg\to M_2(\Clf(q))$, $w+sf+tg\mapsto \begin{psmallmatrix} w & s\\ t & -w\end{psmallmatrix}$, satisfies
$\begin{psmallmatrix} w & s\\ t & -w\end{psmallmatrix}^2=(w^2+st)1_2=-(q(w)-st)1_2$,
so by the universal property it extends to an $R$-algebra homomorphism $\Theta:\Clf(q-yz)\to M_2(\Clf(q))$ with
$\Theta(w)=\begin{psmallmatrix} w&0\\0&-w\end{psmallmatrix}$, $\Theta(f)=\begin{psmallmatrix} 0&1\\0&0\end{psmallmatrix}$ and $\Theta(g)=\begin{psmallmatrix} 0&0\\1&0\end{psmallmatrix}$.
For $x\in\Clf(q)$ we have $\Theta(x)=\begin{psmallmatrix} x&0\\0&x'\end{psmallmatrix}$, since both sides are algebra homomorphisms that agree on $W$.
Hence $\Theta(xfg)=\begin{psmallmatrix} x&0\\0&0\end{psmallmatrix}$, $\Theta(xf)=\begin{psmallmatrix} 0&x\\0&0\end{psmallmatrix}$, $\Theta(xg)=\begin{psmallmatrix} 0&0\\x'&0\end{psmallmatrix}$ and $\Theta(xgf)=\begin{psmallmatrix} 0&0\\0&x'\end{psmallmatrix}$, so $\Theta(\iota(A))=A$ for every $A$.
If $e_1,\ldots,e_n$ is a basis of $W$, the elements $e_Sfg, e_Sf, e_Sg, e_Sgf$ for $S\subseteq\{1,\ldots,n\}$ form a basis of $\Clf(q-yz)$, since they are obtained from the basis $e_S, e_Sf, e_Sg, e_Sfg$ by the substitution $gf=1-fg$.
The map $\Theta$ sends them to $e_S$, $\pm e_S$ times the four matrix units, which form a basis of $M_2(\Clf(q))$.
So $\Theta$ is bijective, and $\iota=\Theta^{-1}$ is an isomorphism of $R$-algebras.

The map $A\mapsto A^{\sigma}$ is an automorphism of $M_2(\Clf(q))$ (apply $x\mapsto x'$ to the entries, then conjugate by $\begin{psmallmatrix}1&0\\0&-1\end{psmallmatrix}$), and $A\mapsto A^{\tau}$ is an anti-automorphism (apply $x\mapsto\overline{x}$ to the entries, transpose, then conjugate by $\begin{psmallmatrix}0&1\\1&0\end{psmallmatrix}$).
The homomorphisms $z\mapsto\Theta(z')$ and $z\mapsto\Theta(z)^{\sigma}$ agree on $W$, $f$ and $g$, so they are equal; likewise the anti-homomorphisms $z\mapsto\Theta(z^*)$ and $z\mapsto\Theta(z)^{\tau}$ agree on $W$, $f$ and $g$, so they are equal.
Applying $\iota$ gives the first two identities in \eqref{eqn:iota-involutions}; the third follows from them, since $\overline{x}=(x')^*$ and $\Adj(A)=(A^{\sigma})^{\tau}$.
Writing $x=x_0+x_1$ for the decomposition of $x\in\Clf(q)$ into even and odd parts, we have $c'=c_0-c_1$ and $d'=d_0-d_1$, so
\begin{equation}\label{eqn:iota-even-odd}
\iota(A) = \underbrace{(a_0fg+b_1f-c_1g+d_0gf)}_{\mbox{even}} + \underbrace{(a_1fg+b_0f+c_0g-d_1gf)}_{\mbox{odd}}.
\end{equation}
\end{proof}

We now give an arithmetic version of the Bott periodicity statement \eqref{eqn:real-bott}.
Note that 
$$\begin{pmatrix} a & b \\ c & d \end{pmatrix}\begin{pmatrix}
	d^* & -b^* \\ -c^* & a^* \end{pmatrix}=\begin{pmatrix} ad^*-bc^* & -ab^*+ba^* \\
	cd^* -dc^* & -cb^*+da^* \end{pmatrix}.$$ 
If $g=\begin{psmallmatrix} a & b \\ c & d \end{psmallmatrix} \in \SL_2(\CC_n)$, we have $g\Adj(g) = \Adj(g) g = \Delta(g) 1_2$ where $1_2$ is the $2\times 2$ identity matrix. 
\begin{theorem}\label{lem:pre-exceptional}
	Let $q$ be a quadratic form over a commutative ring $R$ on a free module $W$ of finite rank.
	Let $Q =q-yz+x^2$.
	The composition of the first and second Bott periodicity maps $\psi = \phi \circ \iota: M_2(\Clf(q) ) \to \Clf(Q)_+$ is an isomorphism of associative algebras given by 
	\begin{equation}\label{eqn:psi-map}
	 \psi\begin{pmatrix} a & b \\ 
	c & d \end{pmatrix} = (a_0 fg+b_1 f - c_1g + d_0 gf)+(a_1fg+b_0f+c_0g-d_1gf)e 
	\end{equation}
	where $a=a_0 + a_1, b=b_0+b_1, c=c_0+c_1, d=d_0+d_1$ corresponds to the $\ZZ/2\ZZ$-grading $\Clf(q) = \Clf(q)_+ \oplus \Clf(q)_-$. 
	Also $e,f,g$ satisfy $e^2=-1$ and $fg=1-gf$, and correspond to the basis for the quadratic space associated to the form $-yz+x^2$.

	The map satisfies 
	 $$\psi(\Adj(A)) = \psi(A)^* $$
	and hence $\psi(A\Adj(A)) =\ns(\psi(A))$; in formulas  
	$$\psi \left( \begin{pmatrix} ad^*-bc^* & -ab^*+ba^* \\
		cd^* -dc^* & -cb^*+da^* \end{pmatrix}\right )= \psi(A)\psi(A)^*.$$
	If $A \in \GL_2(\Clf(q))$, then $\ns(\psi(A))=\Delta(A)$.

        Furthermore, $\psi$ restricts to multiplicative bijections
        $$ \SL_2(\Clf(q)) \to \Spin(Q), \qquad \GL_2(\Clf(q)) \to  \GSpin(Q),$$
        which are group isomorphisms by Corollary~\ref{lem:gl-is-group}.

\end{theorem}
\begin{proof}
By Lemma~\ref{thm:bott1} and Lemma~\ref{lem:bott2} (applied to $q-yz$), $\psi$ is an isomorphism of $R$-algebras, and \eqref{eqn:psi-map} follows from \eqref{eqn:iota-even-odd} and the recipe for $\phi$, which multiplies the odd part on the right by $e$.
By \eqref{eqn:iota-involutions} and \eqref{eqn:phi-involutions}, using
$x^{s}=-ex^*e$ and $e^{-1}=-e$, we obtain
$$\psi(\Adj(A)) = \phi(\overline{\iota(A)}) = \phi(\iota(A))^* = \psi(A)^*, \qquad \psi(A^{\tau}) = \phi(\iota(A)^*) = \psi(A)^{s} = e^{-1}\psi(A)^*e.$$
In particular $\ns(\psi(A))=\psi(A\Adj(A))$, and $\psi(r1_2)=r$ for $r\in R$.

Let $W_Q = W\oplus Rf\oplus Rg\oplus Re$ be the module of $Q$, and let $z=\psi(A)$.
Since $z$ is even we have $\nrd(z)=\ns(z)$, and $z\Vec(Q)z^*\subseteq\Vec(Q)$ follows from $\ns(z)\in R$ and $zW_Qz^*\subseteq W_Q$.
So, by the Lemma following Definition~\ref{defn:spin}, $z\in\GSpin(Q)$ if and only if (i) $\ns(z)\in R^\times$ and (ii) $zW_Qz^*\subseteq W_Q$; and $z\in\Spin(Q)$ if and only if moreover $\ns(z)=1$.

Condition (i). Since $\psi$ is injective, $\ns(z)=\psi(A\Adj(A))$ lies in $R^\times$ if and only if $A\Adj(A)=u1_2$ with $u\in R^\times$.
As $A\Adj(A)=\begin{psmallmatrix} \Delta(A) & ba^*-ab^*\\ cd^*-dc^* & \Delta(A)^*\end{psmallmatrix}$, this holds if and only if conditions \eqref{item:det-cond} and \eqref{item:abcd-cond} of Definition~\ref{def:gl2} hold, and then $\ns(z)=\Delta(A)$.

Condition (ii). Right multiplication by $e$ is injective and $z^*e=e\psi(A^{\tau})$, so (ii) is equivalent to $z\,(W_Qe)\,\psi(A^{\tau})\subseteq W_Qe$.
We have $\psi^{-1}(we)=\begin{psmallmatrix} w&0\\0&-w\end{psmallmatrix}$ for $w\in W$, $\psi^{-1}(fe)=\begin{psmallmatrix} 0&1\\0&0\end{psmallmatrix}$, $\psi^{-1}(ge)=\begin{psmallmatrix} 0&0\\1&0\end{psmallmatrix}$ and $\psi^{-1}(ee)=-1_2$. Hence 
$$  X :=\psi^{-1}(W_Qe) = \left\{ \begin{pmatrix} p & v \\ s & \overline{p} \end{pmatrix}
  \colon p\in \Vec(q),\ v,s\in R \right\},$$
and (ii) holds if and only if $AMA^{\tau}\in X$ for all $M\in X$.
For $M=\begin{psmallmatrix} p & v \\ s & \overline{p} \end{psmallmatrix}$ we compute
$$AMA^{\tau} = \begin{pmatrix}
ap\overline{d}+b\overline{p}\,\overline{c}+v\,a\overline{c}+s\,b\overline{d} & ap\overline{b}+b\overline{p}\,\overline{a}+v\,a\overline{a}+s\,b\overline{b}\\
cp\overline{d}+d\overline{p}\,\overline{c}+v\,c\overline{c}+s\,d\overline{d} & cp\overline{b}+d\overline{p}\,\overline{a}+v\,c\overline{a}+s\,d\overline{b}
\end{pmatrix},$$
and the lower right entry is always the conjugate of the upper left entry.
Taking $(p,v,s)$ equal to $(0,1,0)$, $(0,0,1)$ and $(p,0,0)$ shows that $AXA^{\tau}\subseteq X$ if and only if conditions \eqref{item:norm-cond}, \eqref{item:ratios}, \eqref{item:five} and \eqref{item:six} of Definition~\ref{def:gl2} hold.

Therefore $\psi(A)\in\GSpin(Q)$ if and only if $A\in\GL_2(\Clf(q))$, and $\psi(A)\in\Spin(Q)$ if and only if $A\in\SL_2(\Clf(q))$.
Since $\psi$ is bijective and multiplicative, it restricts to multiplicative bijections $\GL_2(\Clf(q))\to\GSpin(Q)$ and $\SL_2(\Clf(q))\to\Spin(Q)$.
\end{proof}

\begin{remark}
	A version of this over fields is proved in \cite[Proposition 4.1]{Elstrodt1987}.
	There is also a version of this in \cite{McInroy2016}, but that paper unfortunately has an error where the definition of the paravector group (closely related to our group) is defined so that $\SL_2(\ZZ)$ would not contain the element $\begin{psmallmatrix} 1 & 2 \\ 0 & 1 \end{psmallmatrix}$. 
	The entries of elements of
        the $\operatorname{PSV}$ are required to be invertible.
\end{remark}

\begin{corollary}\label{lem:gl-is-group}
	Let $C=\Clf(W,q)$, where $W$ is a free module of finite rank over a commutative ring $R$.
	\begin{enumerate}
	\item  The sets $\GL_2(C)$ and $\SL_2(C)$ are groups, and the inverse of $g\in\GL_2(C)$ is $\Delta(g)^{-1}\Adj(g)$.
        \item $\Delta$ is a homomorphism from $\GL_2(C)$ to $R^\times$.
	\end{enumerate}
	In particular, the bijections of Theorem~\ref{lem:pre-exceptional} are group isomorphisms $\SL_2(\Clf(q))\cong\Spin(Q)$ and $\GL_2(\Clf(q))\cong\GSpin(Q)$.
\end{corollary}

\begin{proof}
  Let $X=\{\begin{psmallmatrix} p & v \\ s & \overline{p}\end{psmallmatrix} : p\in\Vec(C),\ v,s\in R\}\subseteq M_2(C)$.
  The proof of Theorem~\ref{lem:pre-exceptional} shows that $g\in M_2(C)$ lies in $\GL_2(C)$ if and only if (a) $g\Adj(g)=u1_2$ for some $u\in R^\times$, in which case $u=\Delta(g)$, and (b) $gXg^{\tau}\subseteq X$.
  The maps $A\mapsto A^{\tau}$ and $A\mapsto \Adj(A)=(A^{\sigma})^{\tau}$ are anti-automorphisms of $M_2(C)$.

  If $g,h\in\GL_2(C)$, then $gh\Adj(gh)=g\,h\Adj(h)\Adj(g)=\Delta(g)\Delta(h)1_2$ and $(gh)X(gh)^{\tau}=g(hXh^{\tau})g^{\tau}\subseteq X$.
  So $gh\in\GL_2(C)$ and $\Delta(gh)=\Delta(g)\Delta(h)$. Clearly $r1_2\in\GL_2(C)$ for $r\in R^\times$.

  Let $g\in\GL_2(C)$ and $u = \Delta(g)$.
  Left multiplication by $g$ is a surjective endomorphism of the finitely generated $R$-module $M_2(C)$, hence injective \cite[Tag 05G8]{stacks-project}; so the right inverse $u^{-1}\Adj(g)$ of $g$ is a two-sided inverse, and $\Adj(g)g=u1_2$.
  As $\Adj(\Adj(g))=g$, this says that $\Adj(g)$ satisfies (a).
  The $R$-linear map $T\colon X\to X$, $M\mapsto gMg^{\tau}$, is injective, because $g$ and $g^{\tau}$ are units of $M_2(C)$.
  The module $X\cong \Vec(C)\oplus R^2$ is free of finite rank and is a direct summand of $M_2(C)$, compatibly with base change.
  So for each maximal ideal $\mathfrak m$ of $R$, $T\otimes k(\mathfrak m)$ is the map $M\mapsto \bar gM\bar g^{\tau}$ on $X\otimes k(\mathfrak m)\subseteq M_2(C\otimes k(\mathfrak m))$, where $\bar g$ is the image of $g$, a unit; hence $T\otimes k(\mathfrak m)$ is injective and therefore bijective.
  By Nakayama's lemma \cite[Tag 00DV]{stacks-project}, $T$ is surjective, so $g^{-1}X(g^{-1})^{\tau}\subseteq X$ and $\Adj(g)X\Adj(g)^{\tau}=u^2g^{-1}X(g^{-1})^{\tau}\subseteq X$.
  Thus $\Adj(g)\in\GL_2(C)$, and $g^{-1}=u^{-1}\Adj(g)\in\GL_2(C)$.
  This proves (1) for $\GL_2(C)$ and (2); since $\SL_2(C)$ is the kernel of $\Delta$, it is a subgroup.
\end{proof}
\section{Orders in Rational Clifford Algebras}\label{sec:maximal-orders}

Let $F$ be a field. An $F$-algebra $A$ is \emph{separable} if $A$ is semisimple,
the center $Z(A)$ is an \'etale $F$-algebra, and $\dim_F(A) < \infty$.  
Let $F$ be a number field. We recall that an {\em order} in a separable $F$-algebra $A$ is a subring
$\Ocal \subset A$ that is a finitely generated $\ZZ$-module and spans $A$
as a $\QQ$-vector space.

\subsection{Maximal Orders}\label{subsec:maximal-orders}
There is some discussion in \cite[Chapter 10]{Voight2020} that readers may find helpful.  
Our approach is more explicit and computational. 
We begin this section with examples that indicate that the involutions on a
Clifford algebra do not preserve individual maximal orders (it is clear that
the set of maximal orders is preserved by each involution).
Following that, we explain when maximal orders in Clifford algebras are
unique up to conjugation and how to enumerate them. 

\begin{example}
	The Clifford algebra $\CC_n$ has an
	obvious order $\ZZ[i_1,i_2,\ldots.i_{n-1}]= \Clf(\Z^{n-1}, x_1^2+\cdots+x_{n-1}^2)$ generated as a $\Z$-algebra by $i_1, \dots, i_{n-1}$ and
	as a $\Z$-module by the products of these taken without repetition.  
	We call this order the \emph{Clifford order}.
	
	For $n = 2$ it is well-known that this is the unique maximal order---the Gaussian integers; for $n = 3$, these are the Lipschitz quaternions, which are not a maximal
	order.  
\end{example}

The order above is contained in a unique maximal order, namely
the Hurwitz quaternions, which are obtained by adjoining
$(1+i_1+i_2+i_1i_2)/2$.
\begin{example}\label{ex:orders-not-closed}
	We consider the usual quaternion algebra $\CC_3$. Their
	most familiar maximal order is the {\em Hurwitz order} $\Ocal_3$
	\cite[Section 11.1]{Voight2020}, generated as a group by
	$1, i, j, (1+i+j+k)/2$.  

	The order $\Ocal_3$ is closed under
	the standard involutions of $\CC_3$.  
	
	On the other hand, the element
	$\alpha = (3j+4k)/5$ has minimal polynomial $x^2+1$; it is a consequence
	of the Skolem-Noether theorem
	\cite[Main Theorem 7.7.1, Corollary 7.7.3]{Voight2020} that it is conjugate
	to $i$.  
	It is therefore contained in a maximal order $\Ocal'$.
	However, $\Ocal'$ does not contain $\alpha^* = (3j-4k)/5$, for
	if it did it would certainly contain $\alpha + \alpha^* = 6j/5$, which is
	not integral.

	Every order in $\CC_3$ is closed under the conjugation $v \mapsto \bar v$, because
	$\bar v = (v+\bar v)-v$ and $$v+\bar v = \trd(v) \in \Z,$$ for every integral $v \in \CC_3$.
	Since $v^* = \overline{v'}$, an order in $\CC_3$ is closed under the parity involution
	if and only if it is closed under $*$; in particular, the maximal order $\Ocal'$ above is
	closed under neither. Similarly to the above, in $\CC_4$ consider
	the element $\beta = (3i_1 + 4i_1 i_2)/5$; although $\beta$ is integral (indeed $\beta^2=-1$) and
	hence contained in a maximal order, we have $\beta' = (-3i_1 + 4i_1 i_2)/5$
	and so no order contains both $\beta$ and $\beta'$, because $\beta+\beta' = 8i_1i_2/5$ is not integral.
\end{example}
Experiments show that there are many interesting maximal orders in $\quat{-1,-1,\ldots,-1}{\QQ}$ generalizing the Hurwitz order. 
In $\CC_4$ there is another unique one, and in $\CC_5$ there are six of them; these are all conjugate (Proposition~\ref{prop:maximal-orders-unique}), but under conjugation by the Clifford monoid they break into two classes. 
We will return to this later, in $\S$\ref{sec:1-1-1-1}.

We will often need to make the hypothesis that our order $\Ocal$ is closed under $*$. 
This hypothesis appears in \cite{Elstrodt1990} under the term ``compatible''.
This is needed, for example, for $\SL_2(\Ocal)$ to make sense.
\begin{definition}\label{def:clifford-stable}
Let $C$ be a Clifford algebra over $K$.  An order $\Ocal \subset C$ is
{\em Clifford-stable} if $\Ocal^* = \bar \Ocal = \Ocal$.  If $\Ocal^* = \Ocal$ we say that
  $\Ocal$ is \emph{$*$-stable}.
\end{definition}
Every explicit order considered in \S\S\ref{sec:classical}--\ref{sec:1-1-1-1} is Clifford-stable (see Remark~\ref{rem:star-vs-clifford-stable}).
Experimentally, we found maximal orders in the quaternion algebra $\quat{-2,-13}{\QQ}$ which are not $*$-stable.
\begin{remark}\label{rem:star-vs-clifford-stable}
  Since $\bar x = (x^*)'$, a $*$-stable order is Clifford-stable if and only if it is stable under $x\mapsto x'$.
  \begin{enumerate}
  \item If $\Ocal$ is generated as a ring by a set $S$ with $S^*\subseteq \Ocal$ and $S'\subseteq\Ocal$, then $\Ocal$ is Clifford-stable, because $\{x\in\Ocal \colon x^*\in\Ocal\}$ and $\{x\in\Ocal\colon x'\in\Ocal\}$ are subrings of $\Ocal$. This applies to $S\subseteq\Vec(\Ocal)$: for $v\in\Vec(\Ocal)$ we have $v^*=v$ and $v'=\bar v = 2\Re(v)-v$ with $2\Re(v)\in\ZZ$.
  \item In a quaternion algebra $\bar x=\trd(x)-x$, so every order is stable under $x\mapsto \bar x$, and every $*$-stable order is Clifford-stable.
  \item If $n-1$ is even and $J = i_1 i_2\cdots i_{n-1}$ is a unit of $\Ocal$, then $x' = JxJ^{-1}$, so $\Ocal$ is stable under $x\mapsto x'$. This holds for every order of $\CC_n$ containing $\ZZ[i_1,\ldots,i_{n-1}]$ when $n$ is odd.
  \item In general, $*$-stability does not imply Clifford-stability, even for orders containing the Clifford order: in $\quat{-1,-1,-5}{\QQ}$, with $a_3=\sqrt{5}\,i_3$, the order generated over $\quat{-1,-1,-5}{\ZZ}$ by 
  $$\alpha = (1+i_1+i_2+i_1i_2+a_3+3i_1a_3+i_2a_3+3i_1i_2a_3)/4$$ is $*$-stable but does not contain $\alpha'$.
  \end{enumerate}
\end{remark}

We now define the discriminant of an order.
\begin{definition} Let $K$ be a number field and $A$ a $K$-algebra containing
	an order $\Ocal$ as above.  Let $a_1, \dots, a_n$ be a $\Z$-basis for
	$\Ocal$.  The {\em $\Z$-algebra discriminant} $D_\Ocal$ of $\Ocal$
	uses $\tr$, the trace of an element in the left regular representation
	(i.e., $\tr x$ is the trace of the matrix of the linear transformation
	$y \mapsto xy$ on $A$), and is defined by
	$$D_{\Ocal} = \det \left( \tr(a_i a_j) \right)_{i,j}.$$
	The same formula, applied to a $\ZZ$-basis of any full
	$\ZZ$-lattice $L \subset A$ (a finitely generated subgroup spanning $A$ over $\QQ$),
	defines $D_L$.  We say that $\Ocal$ is {\em maximal} if it is not properly contained in any other order.
\end{definition}

\begin{remark} An alternative approach to discriminants would begin from
	the observation that $\Ocal_{K,A} = \{x \in K: x\Ocal \subseteq \Ocal\}$
	is an order in $K$: we could then consider $\Ocal$ as a module for this ring.
	If $\Ocal$ is a free module, we can choose a basis and write formally the
	same definition as above. However, since $\Ocal_{K,A}$ need not be the
	maximal order, modules over it need not be locally free, which leads to
	problems. The present definition is simpler and is sufficient for our
	purposes. Alternatively, we could restrict to $\Ocal_K$-algebras.
\end{remark}

We now give some basic properties of discriminants.
\begin{lemma}\label{lem:discriminant-basics}
	\begin{enumerate}
		\item \label{item:disc-int} $D_\Ocal \in \Z$.
		\item \label{item:disc-ind} If $L \subseteq L'$ are full $\ZZ$-lattices in $A$ with $[L':L] = p$, then $D_{L'} = D_{L}/p^2$.
		\item \label{item:disc-nonzero} If $D_\Ocal \ne 0$ then the same holds for every order of $A$.
		\item If $D_\Ocal \ne 0$ then every order of $A$ is contained in at
		least one maximal order, and in only finitely many.
	\end{enumerate}
\end{lemma}

\begin{proof}\begin{enumerate}
		\item 
		For $x\in\Ocal$, the matrix of $y\mapsto xy$ in a $\ZZ$-basis of $\Ocal$ is integral, so its trace is integral, and
		so $D_\Ocal$ is the determinant of an integral matrix.
		\item 
		We may choose a basis $a_1, \dots, a_n$ of $L'$ such that
		$pa_1, a_2, \dots, a_n$ is a basis of $L$. Then the matrix used
		to compute $D_L$ is obtained from that for $D_{L'}$ by multiplying
		the first row and column by $p$.
		\item 
		Let $\Ocal'$ be an order of $A$. By basic results on finitely
		generated abelian groups, we know that $\Ocal \cap \Ocal'$ has finite
		index in both $\Ocal$ and $\Ocal'$. Thus it is an order
		and the result follows from (\ref{item:disc-ind}).
		\item Let $\Ocal'$ be an order of $A$; then $D_{\Ocal'}\neq 0$ by (\ref{item:disc-nonzero}).
If $\Ocal''\supseteq\Ocal'$ is an order, then applying (\ref{item:disc-ind}) along a
		composition series of the finite group $\Ocal''/\Ocal'$ gives
		$D_{\Ocal'}=[\Ocal'':\Ocal']^2D_{\Ocal''}$, and $D_{\Ocal''}\in\ZZ$ by (\ref{item:disc-int}).
		So $m=[\Ocal'':\Ocal']$ satisfies $m^2\mid D_{\Ocal'}$ and $\Ocal'\subseteq\Ocal''\subseteq\frac1m\Ocal'$.
		Hence the orders containing $\Ocal'$ have bounded index, and there are only finitely
		many of them, each corresponding to a subgroup of one of the finite groups $\frac1m\Ocal'/\Ocal'$.
		Both parts of the statement are now immediate.
	\end{enumerate}
\end{proof}

\begin{definition}
  An order $\Ocal$ is \emph{$p$-maximal} if its index in every order
  containing it is not divisible by $p$. Equivalently, an order $\Ocal$ is
  $p$-maximal if $v_p(D_{\Ocal}) = v_p(D_{\Ocal'})$ for all $\Ocal' \supseteq \Ocal$.
\end{definition}

\begin{remark} Before proceeding to our algorithms, we remark on the
  difference between them and previous work in the subject, starting
  from \cite{IvanyosRonyai} and including the \texttt{MaximalOrder} command
  in \texttt{magma}.  This work is focused on finding a single maximal
  order that contains a given order, whereas we insist on enumerating
  all of them.  Unfortunately this means that it is difficult to estimate
  the running time, because we do not know in advance how many maximal
  orders there might be.  For example, if the number of maximal orders
  containing an order of discriminant $D$ in a fixed $\Q$-algebra
  can grow faster than polynomially in $D$, then it is unreasonable
  to expect a polynomial-time algorithm to list them all.
\end{remark}

\begin{algorithm}\label{alg:p-maximal}
	Let $A$ be a finite-dimensional semisimple algebra over $\QQ$.
	Let $R$ be an order in $A$ of $\ZZ$-discriminant $D$ and let $p$ be prime.
	We may list all $p$-maximal orders in $A$ containing $R$ with index a power of $p$
	by the following algorithm.
	We begin with a queue consisting of the pair
	$(R,D)$ and an empty list of $p$-maximal orders.
	\begin{enumerate}
		\item \label{item:top} Take the top element $(\Ocal,D_\Ocal)$ from the queue.
		\item \label{item:all-sub-p} Construct all submodules $M$ of
		$A$ containing $\Ocal$ with index $p$.  If $p^2 \nmid D_\Ocal$, none of
		them generates an order.
		\item \label{item:enlarge} For each such submodule $M$, determine whether the algebra $\Ocal_M$
		it generates is an order. (To do so, start with the submodule $M$, and
		recursively enlarge $M$ by products of generators until either all
		products are in $M$, in which case we have an order, or the discriminant of the
		enlarged module is not an integer, in which case we do not.)
		If none of these is an order,
		then add $\Ocal$ to the list of $p$-maximal orders.
		\item For each order among the $\Ocal_M$, calculate the index $p^{n_M}$ with which
		it contains $R$. Then the discriminant of $\Ocal_M$ is $D_R/p^{2n_M}$.
		If $(\Ocal_M,D_R/p^{2n_M})$ is not in the queue, add it.
		\item Stop when the queue is empty.
	\end{enumerate}
$\vartriangle$
\end{algorithm}

\begin{proposition} The above algorithm terminates and constructs all
	$p$-maximal orders containing $R$ with index a power of $p$, and only those.
\end{proposition}
\begin{proof} 
	The assumption on $R$ and $A$ implies that $D_R$ is not zero.
	First, at every application of (\ref{item:all-sub-p})
	only finitely many submodules $M$ are constructed, and by
	Lemma~\ref{lem:discriminant-basics} (\ref{item:disc-ind}), no chain of
	submodules has length greater than $v_p(D_\Ocal)/2$.
        In addition, every time we enlarge $M$ in Step
        \ref{item:enlarge}, we divide the discriminant by $p^{2k}$ for some
        $k \ge 1$, and we stop if it is not an integer,
        so this is a finite process. This proves the
	termination. Note that the reasoning in this paragraph is only valid
        under the semisimplicity hypothesis: otherwise we have $D_\Ocal = 0$ and
        cannot conclude anything.  

        It is clear from the construction that the algorithm
	cannot end with any non-$p$-maximal orders on the queue. Finally,
	for correctness, let $\Ocal^p$ be a $p$-maximal order containing $R$ with $[\Ocal^p:R]$ a power of $p$,
	and let $R = \Ocal_0 \subseteq \Ocal_1 \subseteq \dots \Ocal_n = \Ocal^p$
	be a maximal chain of orders. On the first pass through the algorithm,
	we obtain a submodule $M_1 \subseteq \Ocal_1$ that contains $R$ with index $p$ and generates $\Ocal_1$
	(it exists because $\Ocal_1/R$ is a nontrivial $p$-group; if it did not generate $\Ocal_1$, we could insert another order into the chain). Thus
	$(\Ocal_1,D_{\Ocal_1})$ appears in the queue. By induction, it follows that
	all $(\Ocal_i,D_{\Ocal_i})$ appear at some stage, and in particular
	$(\Ocal^p,D_{\Ocal^p})$.
\end{proof}

We now show how to use this algorithm to construct maximal orders rather
than just $p$-maximal orders.

\begin{lemma}\label{lem:sum-orders}
	Let $\Ocal$ be an order contained in the orders
	$\Ocal_1, \Ocal_2$ with indices
	$n_1, n_2$ respectively, and suppose that $(n_1,n_2) = 1$. Then the
	group $\Ocal_1 + \Ocal_2$ is an algebra.
\end{lemma}

\begin{proof} It suffices to show that this set is closed under multiplication.
	Let $a_i \in \Ocal_i$ for $i = 1, 2$.  Then $(n_1 a_1) a_2 \in \Ocal_2$ and
	$a_1 (n_2 a_2) \in \Ocal_1$. Thus for every integer $k$ of the form
	$b_1 n_1 + b_2 n_2$
	for $b_1, b_2 \in \Z$ we have $k(a_1 a_2) \in \Ocal_1 + \Ocal_2$, and in
	particular for $k = 1$.
\end{proof}

\begin{lemma}\label{lem:orders-local}
	Let $\Ocal'$ be an order containing $\Ocal$ with index $n$, and let
	$n = \prod_{i=1}^m {p_i}^{a_i}$ be the prime factorization of $n$. Then
	there are unique orders $\Ocal_{p_1}, \dots, \Ocal_{p_m}$ containing
	$\Ocal$ with index ${p_i}^{a_i}$ and contained in $\Ocal'$,
	and $\Ocal_{p_i}\cap\Ocal_{p_j}=\Ocal$ for $i\neq j$. More generally, the algebra
	$\Ocal_S$ generated
	by the $\Ocal_{p_i}$ for $i \in S$ is the unique intermediate order
	containing $\Ocal$ with index $\prod_{i \in S} {p_i}^{a_i}$.
\end{lemma}

\begin{proof} Let $Q$ be the abelian group $\Ocal'/\Ocal$, let
	$p \in \{p_1,\dots,p_m\}$, and let $Q_p$ be the $p$-Sylow subgroup of $Q$.
	Let $\Ocal_p$ be the inverse image of $Q_p$ in $\Ocal'$; it
	is the only subgroup of
	$\Ocal'$ containing $\Ocal$ with index $p^a$, so it is the only candidate
	for such an algebra. To see that it is an algebra, let $b,c \in \Ocal_p$
	and let $r$ be the smallest positive integer such that $rbc \in \Ocal$.
	Since $p^a b, p^a c \in \Ocal$, we must have $r|p^{2a}$. Thus, the image
	of $bc$ in $\Ocal'/\Ocal$ has an order with a power of $p$ and is contained in 
	$Q_p$, so $bc \in \Ocal_p$ as desired. Since distinct Sylow subgroups of $Q$
	meet trivially, $\Ocal_{p_i}\cap\Ocal_{p_j}=\Ocal$ for $i\neq j$. The more general statement now
	follows from Lemma~\ref{lem:sum-orders}.
\end{proof}

\begin{corollary}\label{cor:assemble-orders}
  Let $\Ocal$ be an order and let $S = \{p_1, \dots, p_m\}$
	be the primes whose squares divide $D_\Ocal$.  Then there is a 
	bijection between maximal orders containing $\Ocal$ and $m$-tuples
	$(\Ocal_{p_1}, \dots, \Ocal_{p_m})$ of $p_i$-maximal orders containing
	$\Ocal$ with $[\Ocal_{p_i}:\Ocal]$ a power of $p_i$. The tuple maps to
	$\sum_i\Ocal_{p_i}$; the inverse sends a maximal order $M$ to the inverse
	images in $M$ of the $p_i$-primary subgroups of $M/\Ocal$.
\end{corollary}

\begin{proof} Let $M$ be a maximal order and let $\Ocal_{p_i}$ be the inverse
	image in $M$ of the $p_i$-primary subgroup of $M/\Ocal$ (with $\Ocal_{p_i}=\Ocal$ if $p_i\nmid[M:\Ocal]$). We claim that
	$\Ocal_{p_i}$ is $p_i$-maximal. If not, it is contained in $\Ocal^+_{p_i}$
	with $[\Ocal^+_{p_i}:\Ocal_{p_i}]$ a power of $p_i$. Let
	$S' = S \setminus \{p_i\}$ and consider the order $\Ocal_{S'}$ of
	Lemma~\ref{lem:orders-local}. By Lemma~\ref{lem:sum-orders}, the group
	$\Ocal^+_{p_i} + \Ocal_{S'}$ is an order; it properly contains $M$,
	contradiction.
	
	Conversely, if all the $\Ocal_{p_i}$ are $p_i$-maximal orders of $p_i$-power index over $\Ocal$, their sum $M$ is
	an order by Lemma~\ref{lem:sum-orders}. If $M'$ is an order containing $M$,
	then $[M':M]$ is supported on the $S$ by Lemma~\ref{lem:discriminant-basics}.
	In particular, if $p_i$ divides the index, then $\Ocal_{p_i}(M')$
	properly contains $\Ocal_{p_i}$, a contradiction since $\Ocal_{p_i}$ was
	assumed to be $p_i$-maximal.
\end{proof}

We can now use Algorithm~\ref{alg:p-maximal} to find all maximal orders.

\begin{algorithm}\label{alg:maximal-orders-algo}
  Let $\Ocal$ be the Clifford order of the Clifford algebra of a nondegenerate quadratic form.
  The following algorithm finds all of the maximal orders containing $\Ocal$.
  \begin{enumerate}
  \item Factor the discriminant of $\Ocal$. Let $p_1, \dots, p_n$ be the primes whose squares divide it.
  \item For each $p_i$, use Algorithm~\ref{alg:p-maximal} to determine the set $M_{p_i}$ of $p_i$-maximal orders
    containing $\Ocal$ with index a power of $p_i$.
  \item For each tuple $(\Ocal_{p_i})_i\in\prod_i M_{p_i}$, take the order
    $\sum_i\Ocal_{p_i}$.
  \end{enumerate}
\end{algorithm}
Corollary~\ref{cor:assemble-orders} ensures the correctness of the algorithm.
We implemented this algorithm in \texttt{magma} to obtain our examples. 
See for example the case of $\quat{-1,-1,-1}{\QQ}$ in \S\ref{sec:1-1-1}.

\begin{remark}\label{rem:small-fixed-or-not}
  In Example~\ref{ex:orders-not-closed} we showed that not all maximal orders
  of a Clifford algebra are preserved by the involutions. Perhaps this is
  not surprising, since the involutions themselves are not preserved by
  conjugation (for example, $(aba^{-1})^* = (a^{-1})^* b^* a^*$, not
  $a b^* a^{-1}$); thus, a conjugate of an order preserved by the involutions
  should not be expected to be preserved itself.
  To get a clearer sense of the situation, we might want
  to restrict to maximal orders containing a fixed order that is preserved
  by the involutions.
  
  In particular, it turns out that each of the six maximal orders containing
  $\ZZ[i_1,i_2,i_3,i_4]$ is preserved by the standard involutions. However, it
  is not clear that every maximal order containing $\ZZ[i_1,i_2,i_3,i_4,i_5]$ is
  preserved by the standard involutions.  
  See Conjecture~\ref{conj:fixed-order-exists} and
  Question~\ref{q:most-orders-not-fixed} in \S\ref{sec:open-problems}.

  On the other hand, for other Clifford algebras it can be shown that no maximal order is
  preserved by the standard involutions. 
In forthcoming work we will prove by an elementary argument that this is true for the quaternion algebra $\quat{-2,-13}{\QQ}$, for example. Examining this in the present article will lead us too far afield.
\end{remark}

In general, it follows from Lemma~\ref{lem:discriminant-basics}
that any order in any rational Clifford algebra associated to a nondegenerate
quadratic form
is contained in only finitely many maximal orders. In a situation where
Bott periodicity applies, we obtain the structure of one maximal order.

It is well-known that the left ideal classes of a maximal order
in a quaternion algebra over
$\QQ$ ramified only at $p$ and $\infty$ are in bijection with the
supersingular $j$-invariants over $\overline{\FF_p}$, so that the number of
these is roughly $p/12$, while the conjugacy classes of maximal orders
correspond to the orbits $\{j,j^p\}$, so that there are roughly $p/24$ of them \cite[Theorem 25.3.15; Section 42.3.8]{Voight2020}.
However, this turns out not to give good intuition
for the maximal orders in Clifford algebras. Before stating the next result, we recall
some definitions and results from \cite{Reiner1975}.

\begin{definition}\label{def:reiner}
  Let $K$ be a number field and let $A$ be a central simple $K$-algebra. Given a full
  $\Ocal_K$-lattice $M \subseteq A$, 
  
  \begin{itemize}
  	\item   Let $\Ocal_\ell(M) = \{x \in A: xM \subseteq M\}$ and
  	$\Ocal_r(M) = \{x \in A: Mx \subseteq M\}$ be the \emph{left order} and \emph{right order} for $M$ respectively. 
  	\item   If these are both maximal orders, we call $M$ a {\em normal ideal} (clearly $M$ is a left ideal in $\Ocal_\ell(M)$ and a
  	right ideal in $\Ocal_r(M)$).
  	\item   Given a place $\frakp$ of $K$, let the {\em local index} $m_\frakp$
  	\cite[Section 32]{Reiner1975} be the order of $A \otimes_K K_\frakp$ in the Brauer group of $K_\frakp$.
  	\footnote{recall that $\Br K_\frakp \cong \QQ/\ZZ$ if $\frakp$ is nonarchimedean, of order $2$ if $\frakp$ is real, and trivial
  	if $\frakp$ is complex.}
  	\item   Let $\Cl_A(\Ocal_K)$ be the {\em ray class group} of $A$
  	\cite[Section 35]{Reiner1975}, which is the group  of invertible fractional ideals of $\Ocal_K$ modulo the principal ideals
  	generated by $x \in K^\times$ that are positive at all real places $\frakp$ of $K$ of local index $2$.
  	\item If $[A:K]=n^2$, define $H_A\subseteq\Cl_A(\Ocal_K)$ to be the subgroup
  	generated by the classes of $\frakp^{n/m_\frakp}$ for finite primes $\frakp$.
  	\item Let $M, N$ be full $\Ocal_K$-lattices in $A$. If $\Ocal_r(M) = \Ocal_\ell(N)$, we say that $MN$ is a {\em proper product}
  	(\cite[p. 183]{Reiner1975}).   This definition is extended inductively to products of arbitrary length.
  \end{itemize}
\end{definition}

\begin{lemma}\label{lem:max-ords} Let $A$ be a central simple algebra of degree $n$
  over a number field $K$. Let
  $\Lambda_1, \Lambda_2$ be maximal orders of $A$.  Let $I(\Lambda_1)$ be the group of two-sided fractional ideals
  of $\Lambda_1$.
  \begin{enumerate}
  \item\label{item:normal-ideal}
    There exists a normal ideal with left order $\Lambda_1$ and right order $\Lambda_2$, and any two such ideals
    $M, N$ satisfy $NM^{-1} \in I(\Lambda_1)$.
  \item\label{item:ideal-group} $I(\Lambda_1)$ is a free abelian group on the prime ideals of $\Lambda_1$.
  \item\label{item:reduced-norm} Let $\frakP$ be a prime ideal of $\Lambda_1$ and let $\frakp = \frakP \cap \Ocal_K$.
    The reduced norm of $\frakP$ is
    $\frakp^{n/m_\frakp}$ (whose class belongs to $H_A$), and reduced norm is multiplicative on proper products.
  \item\label{item:distance}
    Let $d(\Lambda_1,\Lambda_2) = [\Nrd(M)] \in \Cl_A(\Ocal_K)/H_A$, where $M$ is a normal ideal with left and
    right order $\Lambda_1, \Lambda_2$, respectively. Then $d(\Lambda_1,\Lambda_2)$ is well-defined (i.e., it does not
    depend on the choice of $M$) and, for any maximal order $\Lambda_3$,
    $d(\Lambda_1,\Lambda_3) = d(\Lambda_1,\Lambda_2) + d(\Lambda_2,\Lambda_3)$.
  \item\label{item:dist-conj-invariant} For all $x\in A^\times$ we have $d(\Lambda_1,x^{-1}\Lambda_1 x) = 0$, and
    $d(\Lambda_1,\Lambda_2)$ depends only on the conjugacy class of $\Lambda_2$.
  \end{enumerate}
\end{lemma}

\begin{proof}
  All of these follow readily from results of \cite{Reiner1975}.  To be precise:
  \begin{enumerate}
  \item One such ideal is $\Lambda_1\Lambda_2$.  For the second statement, see Theorem 22.21 and its proof.
  \item Theorem 22.10.
  \item See (24.6), Theorem 24.11, and the local calculation in the proof of
    Theorem 24.13 of \cite{Reiner1975}.
  \item The first statement follows from (\ref{item:normal-ideal}) and (\ref{item:reduced-norm}).
    For the second statement, let $M_1, M_2$ be such that the left order of
    $M_i$ is $\Lambda_i$ and the right order is $\Lambda_{i+1}$.  Then the product $M_1M_2$ is proper and the result follows.
  \item The first statement follows by choosing $M = \Lambda_1 x$, since
    $\Nrd(\Lambda_1x)=\Nrd(x)\Ocal_K$ \cite[Corollary 24.12]{Reiner1975}.
    At every ramified real place the reduced norm is positive \cite[Theorem 33.15]{Reiner1975}.
    The second statement follows from \eqref{item:distance}.
  \end{enumerate}
\end{proof}

\begin{lemma}\label{lem:distance-injective}
  Suppose that $d(\Lambda,\Lambda') = 0$ and $[A:K] > 4$. Then $\Lambda, \Lambda'$ are conjugate.
\end{lemma}

\begin{proof}
  As above, choose a normal ideal $M$ with left and right order $\Lambda, \Lambda'$ respectively.  Let $X \in I(\Lambda)$
  with $[\Nrd(M)] = [\Nrd(X)]$ in $\Cl_A(\Ocal_K)$. Since the product $X^{-1}M$ is normal, the ideal $N = X^{-1}M$ is a
  normal ideal with left and right order $\Lambda, \Lambda'$. The equality in the ray class group ensures that
  $\Nrd(N) = a\Ocal_K$ with $a$ positive at every real place ramified in $A$.  By Eichler's theorem
  \cite[Theorem 34.9]{Reiner1975}, it follows that $N = \Lambda x$ for some $x \in A^\times$ (the hypothesis $[A:K] > 4$
  ensures that the Eichler condition \cite[Definition 34.3, Remark 34.4]{Reiner1975} is satisfied). Thus the right order of $N$
  is $x^{-1}\Lambda x$.
\end{proof}

\begin{lemma}\label{lem:distance-surjective}
  The map $d$ from the set of maximal orders to $\Cl_A(\Ocal_K)/H_A$ is surjective.
\end{lemma}

\begin{proof}
  Fix a maximal order $\Lambda$ and a class in $\Cl_A(\Ocal_K)$, and choose an integral ideal $I$ representing it (multiply an arbitrary fractional ideal in
  the class by a positive integer).  For each $\frakp|I$, there exists $x_\frakp \in A_\frakp^\times$ such that
  $\Nrd(x_\frakp) \Ocal_{K,\frakp} = I\Ocal_{K,\frakp}$ by \cite[Theorem 33.4]{Reiner1975}.  Define $x_\frakp = 1$ for
  $\frakp \nmid I$.  Now \cite[Theorem 5.3]{Reiner1975} proves that there is a left $\Lambda$-lattice $L \subset A$ of full rank
  whose completion at every $\frakp$ is $\Lambda_\frakp x_\frakp$.  In particular, the right orders of these local lattices are
  $x_\frakp^{-1} \Lambda_\frakp x_\frakp$ and hence are maximal orders. Since the completion of the right order is the right order of the
  completion, it follows from \cite[Corollary 11.6]{Reiner1975} that the right order $\Lambda'$ of $L$ is maximal.  From the
  construction, $d(\Lambda,\Lambda') = [I]$.
\end{proof}

\begin{corollary}\label{lem:distance-bijective}
  If $[A:K]>4$, then, for a fixed maximal order $\Lambda$, the map
  $\Lambda'\mapsto d(\Lambda,\Lambda')$ induces a bijection from the conjugacy
  classes of maximal orders to $\Cl_A(\Ocal_K)/H_A$.
\end{corollary}
\begin{proof}
  This follows from Lemmas~\ref{lem:max-ords}, \ref{lem:distance-injective},
  and \ref{lem:distance-surjective}.
\end{proof}

We now introduce some notation for rational Clifford algebras and reformulate the target of $d$.
  For the rest of this section, fix $q = d_1x_1^2 + \dots + d_rx_r^2$, a nondegenerate quadratic form over
  $\QQ$ of rank $r > 3$.  
  Let $C=\Clf(q)$ be the rational Clifford algebra associated to $q$, with generators $i_1,\dots,i_r$, so that $i_j^2=-d_j$.
 
  Write $C_+=\Clf(q)_+$ for its even subalgebra.
  Put $J=i_1\cdots i_r$ and 
  $$\delta=J^2=(-1)^{r(r+1)/2}d_1\cdots d_r.$$
  Let $s_+$ and $s_-$ be the numbers of positive and of negative coefficients $d_j$.  If $r$ is odd and $\delta$ is not a
  square in $\QQ$, then let $F = \QQ(J) \cong \QQ(\sqrt \delta)$ and let $m = 2^{(r-1)/2}$.
  Multiplication identifies $C$ with $C_+\otimes_\QQ F$ in this case. Let $\Cl_C(F)$ be the
  narrow class group of $F$ if $s_+ - s_- \equiv 3 \bmod 8$, and the class group otherwise.
  Let $S$ be the set of primes $\frakp$ of $\Ocal_F$ lying over rational primes $p$ that split in $F$
  for which $C_+\otimes_\QQ\QQ_p$ is not a matrix algebra over $\QQ_p$. Let
  $$H'_C = \big\langle\, \Cl_C(F)^m,\ [\frakp]^{m/2} \colon \frakp\in S \,\big\rangle$$
  and $G_C = \Cl_C(F)/H'_C$.

\begin{lemma}\label{lem:clc-k}
  Suppose $r$ is odd and $\delta$ is not a square in $\QQ$. Then
  \[\Cl_C(\Ocal_F)=\Cl_C(F),\qquad \Cl_C(\Ocal_F)/H_C=G_C.\]
\end{lemma}

\begin{proof}
  If $F$ is an imaginary quadratic field, there is nothing to do; there is no distinction between narrow class group and class group
  and no infinite primes are ramified, so $\Cl_C(\Ocal_F) = \Cl_C(F) = \Cl(F)$.  If $F$ is real, then 
  the parity involution $x\mapsto x'$ sends $J$ to $-J$, so it is semilinear for the nontrivial automorphism of $F$ and
  identifies the completions of $C$ at the two real places of $F$.
  So both real places ramify in $C$ exactly when $C\otimes_\QQ\RR$ is a product of two matrix algebras over $\HH$.
  Since $i_j^2=-d_j$, the algebra $C\otimes_\QQ\RR$ is the real Clifford algebra with $s_-$ generators of square $+1$ and $s_+$
  of square $-1$, and by the classification of real Clifford algebras \cite[\S 17.6, p.~228]{Lounesto2001} this is a product of two matrix algebras over
  $\HH$ exactly when $s_--s_+\equiv 5\bmod 8$, that is, when $s_+-s_-\equiv 3\bmod 8$. 
Hence $\Cl_C(\Ocal_F)=\Cl_C(F)$.
  Reversion gives $C_+\cong C_+^{\mathrm{op}}$, so the Brauer class of $C_+$
  has order at most $2$. At a finite prime $\frakp$ of $F$ over $p$, the
  invariant of $C\otimes_F F_\frakp$ is
  $[F_\frakp:\QQ_p]$ times the invariant of $C_+\otimes_\QQ\QQ_p$
  \cite[Theorems 31.4 and 31.9]{Reiner1975}. Thus the local index is $2$
  exactly at the primes in $S$: at a nonsplit $p$ the local degree is $2$
  and kills the invariant, while at a split $p$ the local degree is $1$.
  Since the fractional ideals are generated by prime ideals, $H_C$ is generated by
  the $m$th powers in $\Cl_C(F)$ and $[\frakp]^{m/2}$ for $\frakp\in S$.
  Thus $\Cl_C(\Ocal_F)/H_C=G_C$.
\end{proof}

\begin{proposition}\label{prop:maximal-orders-unique}
  \begin{enumerate}
  \item \label{item:max-orders-conjugate}
    If $r$ is even, or if $r$ is odd and $\delta$ is a square, then all maximal orders of $C$ are conjugate.
  \item \label{item:max-orders-classes}
    Suppose that $r$ is odd and $\delta$ is not a square in $\QQ$. Fix a maximal
    order $\Lambda$ of $C$. Then the map
    $\Lambda'\mapsto [\Nrd(\Lambda\Lambda')]$ induces a bijection from the
    conjugacy classes of maximal orders of $C$ onto $G_C$. In particular,
    all maximal orders are conjugate if and only if $\#\Cl_C(F)$ is odd.
  \item \label{item:max-orders-isomorphism}
    In the situation of (\ref{item:max-orders-classes}), there are $\frac12 \left(|G_C|+|G_C[2]|\right)$
    isomorphism classes of maximal orders of~$C$, where $G_C[2]$ is the kernel of multiplication by $2$. In particular, isomorphism
    of maximal orders is equivalent to conjugacy if and only if $G_C$ has exponent $1$ or $2$.
  \end{enumerate}
\end{proposition}

\begin{proof} 
  \begin{enumerate}
  \item If $r$ is even, then $C$ is central simple over $\QQ$ of degree
    $2^{r/2}>2$. The ray class group over $\QQ$ is trivial, so
    Corollary~\ref{lem:distance-bijective} gives the result.

    In the second case a further argument is needed, since the centre of $C$ is $\QQ \oplus \QQ$ and so $C$ is not central simple
    over $\QQ$. Choose $t\in\QQ$ with $t^2=\delta$. The map
    $\QQ[J]\cong\QQ\times\QQ$ sending $J$ to $(t,-t)$ gives
    $C\cong C_+\otimes_\QQ\QQ[J]\cong C_+\times C_+$.
Now, by \cite[Theorem 10.5]{Reiner1975} the maximal orders of a product are just the products of the orders of the factors,
    and conjugacy of a product is the same as conjugacy of all factors.  Let $\Lambda_1, \Lambda_2$ be maximal orders of $C$.
    Each factor is a maximal order in $C_+$, a central simple algebra over $\QQ$
    of degree $2^{(r-1)/2}>2$, and hence its maximal orders are conjugate by
    Corollary~\ref{lem:distance-bijective}. Thus $\Lambda_1$ and $\Lambda_2$
    are conjugate in $C$.
  \item A maximal $\ZZ$-order in $C$ is the same as a maximal
    $\Ocal_F$-order \cite[Theorem 10.5(iii)]{Reiner1975}. By
    Lemma~\ref{lem:max-ords}, $d(\Lambda,\Lambda')=[\Nrd(\Lambda\Lambda')]$.
    Thus the first
    statement follows from Corollary~\ref{lem:distance-bijective} and
    Lemma~\ref{lem:clc-k}. For the second, we
    need $r > 3$. This implies that $m = 2^{(r-1)/2} \ge 4$, and so that
    $H'_C \subseteq \Cl_C(F)^2$. Conversely, if $\Cl_C(F)$ has odd order,
    its $m$th-power subgroup is all of $\Cl_C(F)$, so $H'_C=\Cl_C(F)$.
  \item Let $\Lambda_1, \Lambda_2$ be abstractly isomorphic orders.  Extend the isomorphism to
    $\phi \otimes \QQ: \Lambda_1 \otimes \QQ \to \Lambda_2 \otimes \QQ$, i.e., to an automorphism of $C$.
    It fixes $F$, the centre of $C$, as a set, so either it fixes $F$
    pointwise or it acts by the nontrivial automorphism.  If it fixes $F$ pointwise, it is an automorphism of $C$ as an $F$-algebra,
    so it is an inner automorphism by the Skolem-Noether theorem \cite[Theorem 7.21]{Reiner1975}.  If not, let $\alpha$ be the
    parity involution; then $\alpha \circ (\phi \otimes \QQ)$ is inner.  Thus two maximal orders $\Lambda_1, \Lambda_2$
    are abstractly isomorphic if and only if $\Lambda_2$ is conjugate to either $\Lambda_1$ or $\alpha(\Lambda_1)$.
    
    Let $\sigma$ be the nontrivial automorphism of $F$.
    Note that if $M$ is a normal ideal with left and right orders $\Lambda_1,  \Lambda_2$, then $\alpha(M)$ is a normal ideal
    with left and right orders $\alpha(\Lambda_1), \alpha(\Lambda_2)$ and that $\Nrd(\alpha(M)) = \sigma(\Nrd(M))$.
    Further, the action of $\sigma$ on $G_C$ is by inversion.

    Choose a maximal order $\Lambda_0$ and let $c = d(\Lambda_0,\alpha(\Lambda_0))$.  Then
    $d(\Lambda_0,\alpha(\Lambda_1)) = d(\Lambda_0,\alpha(\Lambda_0)) + d(\alpha(\Lambda_0),\alpha(\Lambda_1)) = c - d(\Lambda_0,\Lambda_1)$
    and so the isomorphism classes are in bijection with the orbits of $x \to c-x$ on $G_C$.  To prove the statement, it suffices
    to show that $\Lambda_0$ can be chosen with $c = 0$, i.e., such that $\Lambda_0$ is conjugate to $\alpha(\Lambda_0)$. If this is true,
    the result follows immediately.

    To find such a $\Lambda$, choose a maximal order $\Lambda_+$ of $C_+$.
    For each prime $p$ write $F_p=F\otimes_\QQ\QQ_p$, so that $C\otimes_\QQ\QQ_p=C_{+,p}\otimes_{\QQ_p}F_p$ and $\alpha=\id\otimes\sigma$,
    and choose a maximal order $\Gamma_p$ of $C\otimes_\QQ\QQ_p$ as follows.
    If $p$ splits in $F$, or $C_{+,p}$ is a matrix algebra over $\QQ_p$, let $\Gamma_p=\Lambda_{+,p}\otimes_{\ZZ_p}\Ocal_{F,p}$;
    this is a maximal order, isomorphic to $\Lambda_{+,p}\times\Lambda_{+,p}$ or to $M_m(\Ocal_{F,p})$.
    If $p$ ramifies in $F$ and $C_{+,p}$ is not split, then, since the class of $C_+$ has order at most $2$,
    its local index at $p$ is $2$ \cite[Theorems 7.4 and 31.4]{Reiner1975},
    so $C_{+,p}\cong M_{m/2}(D)$ with $D$ a quaternion division algebra;
    then $D$ contains the unramified quadratic extension $W$ of $\QQ_p$, and the ring of integers $\Ocal_W$ lies in the maximal order of $D$ \cite[Theorem 12.8 and the discussion after Theorem 14.3]{Reiner1975}.
    Then $\Ocal_W\otimes_{\ZZ_p}\Ocal_{F,p}$ is the ring of integers of the unramified quadratic extension $W F_p$ of $F_p$, and exactly one maximal order $\Gamma_0$ of $D\otimes_{\QQ_p}F_p\cong M_2(F_p)$ contains it, because all its fractional ideals are homothetic. As $\alpha$ preserves $\Ocal_W\otimes_{\ZZ_p}\Ocal_{F,p}$, it preserves $\Gamma_0$; let $\Gamma_p=M_{m/2}(\Gamma_0)$.
    If $p$ is inert in $F$ and $C_{+,p}$ is not split, let $\Gamma_p$ be any maximal order (one exists by \cite[Corollary 10.4]{Reiner1975}). 
    In the first two cases $\alpha(\Gamma_p)=\Gamma_p$.
    Since $C_+$ splits at almost all $p$ \cite[Theorem 32.1(i)]{Reiner1975}, we have $\Gamma_p=\Lambda_{+,p}\otimes_{\ZZ_p}\Ocal_{F,p}$ for almost all $p$,
    and so there is a maximal order $\Lambda$ of $C$ with these completions \cite[Theorem 5.3, Corollary 11.6]{Reiner1975}.
    The normal ideal $\Lambda\,\alpha(\Lambda)$ has trivial completions except at the finitely many inert primes $p$ of the last case, where its reduced norm is a power of $p\Ocal_{F,p}$.
    In particular, since the $p$ are inert in $F$, the reduced norm ideal is generated by a positive rational number.  It follows that
    $c = 0$ for this maximal order as desired.
  \end{enumerate}
\end{proof}

\begin{remark}\label{rem:need-n-gt-3}
	If $r = 3$ this does not work. The Clifford algebra in this
	case is a quaternion algebra over its centre, or a product of two quaternion
	algebras over $\QQ$, and it need not satisfy the Eichler condition: for example,
	$\quat{-1,-1,-3}{\QQ}$ is a totally definite quaternion algebra over $\QQ(\sqrt 3)$.
	We will see examples in \S\ref{sec:1-1-3}
        that prove that nonisomorphic maximal orders do indeed exist.
\end{remark}

\begin{remark}
	Given an isomorphism of a Clifford algebra with a matrix algebra over
	$\Q$, a quadratic field, or a quaternion algebra, or with a sum of two such
	rings, it is easy to write down a single maximal order. For $r>3$,
	Proposition~\ref{prop:maximal-orders-unique} gives one conjugacy class in
	case~\eqref{item:max-orders-conjugate} and $|G_C|$ classes in
	case~\eqref{item:max-orders-classes}. For $\CC_n$ with $n>4$ there is only one class: either
	$n-1$ is even, or $\delta=\pm1$ and $\QQ(i)$ has class number $1$. On the other hand,
	for $q=x_1^2+x_2^2+x_3^2+x_4^2+5x_5^2$ the centre is $\QQ(\sqrt{-5})$ and
	$G_C\cong\ZZ/2\ZZ$, so there are two conjugacy classes, and they are not isomorphic.
	For $q=x_1^2+\cdots+x_5^2-x_6^2-3x_7^2$ the centre $\QQ(\sqrt 3)$ has class number $1$
	but narrow class number $2$, and again there are two conjugacy classes.
\end{remark}

\subsection{Units and Zero Divisors}\label{sec:units}

The following remark explains how zero-divisors can appear unusual in this non-commutative setting.
\begin{remark} 
	In a Clifford algebra, an element $x$ is a left zero-divisor
	if and only if it is a right zero-divisor, because both are equivalent to
	the constant coefficient of the minimal polynomial being $0$.  
	However,
	the concept of ``zero-divisor'' requires caution in
	noncommutative rings. Let $R$ be a ring with nonzero central elements
	$s, t$ satisfying $st = 0$. Then we have $sxt = 0$ for all $x \in R$,
	even if $x$ is a unit.  
	
	This phenomenon occurs in some Clifford algebras:
	for example, in $\CC_4$ the elements $s = (1+i_{123})/2, t = (1-i_{123})/2$
	are commuting orthogonal idempotents (we know by Bott periodicity that
	such elements must exist) and so $svt = 0$ even if $v$ is a nonzero vector.
\end{remark}

The following comes in handy in the discussion of Euclidean orders.
\begin{corollary}\label{cor:no-zero-divisors}
	Let $R$ be an integral domain, let $W$ be a free $R$-module, and let
	$(W,q)$ be a strongly anisotropic quadratic space (Definition~\ref{def:strongly-anisotropic}).
	Then no nonzero element of $\Clf(q)^\mon$ is a left or right zero-divisor.
\end{corollary}

\begin{proof} Let $v \ne 0$ be an element of $\Clf(q)^\mon$.  Then
	$v \bar v \in R$ by the definition of the Clifford monoid, so $v \bar v = \hat q(v)$, and $\hat q(v) \neq 0$ because $\hat q$ is anisotropic.
	Since $R$ is a commutative integral domain and $\Clf(q)$ is a free $R$-module, multiplication by $v \bar v$ is injective on $\Clf(q)$.
	Thus, if $av = 0$ then $av \bar v = 0$, contradiction, so $v$ is not a
	right zero-divisor. Over the fraction field of $R$ the element $v$ has the inverse $\bar v/(v\bar v)$, so $\bar v v = v \bar v$; using $\bar v v$ in the same way, we see that
	$v$ is not a left zero-divisor either.
\end{proof}

The following fact about zero divisors is fundamental, but several of the graduate algebra textbooks we checked didn't contain it so we include it. 
\begin{lemma}\label{lem:left-implies-right}
	Let $R$ be a commutative integral domain with fraction field $K$, let
	$A$ be a finite-dimensional $K$-algebra, and let $\Ocal$ be an $R$-order in
	$A$. Let $x \in \Ocal$. Then $x$ is a left zero-divisor if, and only if, it
	is a right zero-divisor.
\end{lemma}

\begin{proof}
	If $x$ is a left zero-divisor in $\Ocal$, then it is not a left unit in $A$
	and hence not a right unit in $A$ either (here we use that $A$ is
        finite-dimensional over $K$). It is therefore a right
	zero-divisor in $A$: say $wx = 0$ with $w \ne 0$. Because $\Ocal$ is
	an $R$-order in $A$, we can find $r \ne 0 \in R$ with $rw \in \Ocal$.
	Since $r$ is a unit in $A$, we must have $rw \ne 0$, so the equation
	$rwx = 0$ implies that $x$ is a right zero-divisor. Similarly in the
	opposite direction.
\end{proof}

\begin{lemma}\label{lem:units-norm-1}
Let $d_1,\dots,d_n$ be positive integers and let $\Ocal$ be an order in $(-d_1,\ldots,-d_n/\QQ)$. Let $x \in \Ocal^\times$. Then $\nrd(x) = 1$.
\end{lemma}

\begin{proof}
	Let $x \in \Ocal^\times$ and let $c_x$ be the characteristic polynomial of
	multiplication by $x$ on the algebra and $m_x$ the minimal polynomial
	of $x$. If $x$ belongs to an order, then $m_x$ must
	be integral, and for $x$ to be a unit of $\Ocal$ the constant term of
	$m_x$ must be $\pm 1$ (otherwise $m_{1/x}$ would not be
	integral). It follows that the constant term of
	every factor of $m_x$ is $\pm 1$. Since $c_x$ is a product of powers of
	factors of $m_x$, the same holds for $c_x$. The constant term of $c_x$
	is the algebra norm $N_A(x)$ (the sign is correct because $\dim \CC_n$ is
	even).
	
	Thus $\nrd(x)^{2^{n-1}} = N_A(x) = 1$ (Proposition~\ref{prop:big-form-02}), so $\nrd(x) = \pm 1$. But $\nrd(x)=|x|^2\geq 0$ in $\CC_{n+1}$ (Corollary~\ref{cor:norms-over-q}).
\end{proof}

\begin{proposition}\label{prop:units-finite}
  Let $d_1, \dots, d_n$ be positive integers. Then 
  the group of Clifford units $\Ocal^\times$ of any order in
  $\quat{-d_1,\ldots,-d_n}{\QQ}$ is finite. 
\end{proposition}

\begin{proof} This group is contained in the intersection of the discrete
	set $\Ocal$ with the compact set defined by $|x|^2 = \nrd(x) = 1$ in $\CC_{n+1}$ (Lemma~\ref{lem:units-norm-1}).
\end{proof}

\begin{remark} We will see later (Proposition~\ref{rem:not-invariant})
	that different maximal orders of the same Clifford algebra, even if
	conjugate, may have different groups of Clifford units.
\end{remark}

We describe an algorithm to determine the group of units of an order
$\Ocal \subset \quat{-d_1,\ldots,-d_n}{\QQ}$.
This implements the proof of Proposition~\ref{prop:units-finite}, which is essentially that elements of $\Ocal^{\times}$ need to have norm 1 and, conversely, that if $\overline{\Ocal}=\Ocal$ then any element $x$ of norm 1 in the Clifford monoid of $\Ocal$ lies in $\Ocal^{\times}$, since $x^{-1}=\bar x$.

\begin{algorithm}[Units]\label{alg:units}
Let $\Ocal \subset \quat{-d_1,\ldots,-d_n}{\QQ}$ be a Clifford-stable order. Let 
  $$a_0=1,\quad  a_1 = \sqrt{d_1}\, i_1, \quad a_2 = \sqrt{d_2}\, i_2, \quad a_{12} = a_1a_2, \ldots$$
  be the standard basis of the algebra, so that $a_j^2=-d_j$ and $a_S\overline{a_S}=\prod_{s\in S}d_s$.

  \begin{enumerate}
  \item The quadratic form for which the $a_S$ are orthogonal and $q(a_S) = a_S \overline{a_S}$
    is positive definite and coincides with the reduced norm on Clifford monoid
    elements. 
  \item Let $\Lambda$ be the lattice with this quadratic form.
    Let $d$ be the LCM of the denominators of the coefficients of
    the basis elements of $\Ocal$ with respect to the standard basis.
    Then $d \Ocal$ defines a sublattice of $\Lambda$ in which all units
    have norm $d^2$.  
  
  \item Using standard lattice methods, we enumerate the elements of
    $d \Ocal$ of norm $d^2$, retaining only those that belong to the Clifford
    monoid.   $\vartriangle$
  \end{enumerate}
\end{algorithm}

\begin{remark}
  The running time of this algorithm tends to be inordinately long for
  $n > 5$. As a practical compromise, we generally proceed as follows.
 
  \begin{enumerate}
  \item Consider all sums $\sum_{i \in S} s_i b_i$, where $b_i$ is the basis of
    the order, the $s_i$ are $\pm 1$, and $S$ is a multiset of size at most
    $k$ for some $k$. 
  
  \item Again, determine which of these are Clifford monoid elements of
    norm $1$ and retain them as units. Let the set of these units be $U_k$.
  \end{enumerate}
  
  Typically, if $k$ is not too small then $U_k$ will generate
  $\Ocal^\times$.  We may increase our confidence by computing
  $U_k, U_{k+1}, \dots, U_{k+\ell}$; if they are all equal it is likely that they
  coincide with the full group of units.
\end{remark}

\begin{remark}
	If $\Ocal$ is $*$-stable and we want to find the order of the group of units but not list all
	the units, this may be done by means of the action $v\mapsto uvu^*$ of the units on the
	Clifford vectors $\Vec(\Ocal)$. This gives us a representation to $\GL_{n+1}(\Z)$ with
	kernel $\pm 1$.  
	
	As a practical matter, finding the order of the image
	in \texttt{magma} can be slow, and it is often better to reduce modulo a
	small prime such as $3$ (it is well known that the reduction map is
	injective on torsion groups) and use the \texttt{CompositionTree}
	functionality to determine the order of the reduction.
\end{remark}

\subsection{Unimodularity/Coprimality}\label{sec:unimodularity}

We will need to know what it means for $a,b\in \Ocal$ to be ``coprime''.
To be consistent with the other literature we call this property ``unimodularity''. 

\begin{definition}\label{def:right-unimodular}
	Let $\Ocal$ be a $*$-stable order.
	We say that $(\mu,\nu) \in \Ocal^2$ is \emph{right unimodular} if and only if there exists a matrix $\begin{psmallmatrix} * & * \\ \mu & \nu \end{psmallmatrix}\in \SL_2(\Ocal)$.  If there is $\begin{psmallmatrix} \mu & * \\ \nu & * \end{psmallmatrix} \in \SL_2(\Ocal)$, then $(\mu,\nu)$ is \emph{left unimodular}.        
\end{definition}

We want to prove that this condition is equivalent to several other conditions so we can work with it fluidly. 
One thing we will want to do is take adjoints of our matrices.
\begin{definition}
	The \emph{Clifford adjoint} of an $m\times n$ matrix $A \in M_{m,n}(\CC_n)$ is the $n\times m$ matrix $A^{\dag} \in M_{n,m}(\CC_n)$ given by taking the conjugate transpose:
	$$ A^{\dag} = (\overline{A})^{\t} = \overline{(A^{\t})}.$$
\end{definition}

\begin{lemma}\label{L:closed-under-conj}
	Both $\GL_2(\CC_n)$ and $\SL_2(\CC_n)$ are closed under $g\mapsto g^{\dag}$.
	When $\Ocal$ is Clifford-stable (Definition~\ref{def:clifford-stable}) then $\SL_2(\Ocal)$ is closed under $g\mapsto g^{\dag}$.
\end{lemma}
\begin{proof}
	Let $g = \begin{psmallmatrix} a & b \\ c & d \end{psmallmatrix}$ be in $\GL_2(\CC_n)$. 
	One can check that $\Delta(g^{\dag}) = \overline{\Delta(g^{-1})}$.
	To see this, note that $\Delta(g^{-1}) = d^*a-b^*c$ and $\Delta(g^{\dag}) = (a^*)'d'-(c^*)'b'$. 
	We apply $*$ then $'$ to the $\Delta(g^{-1})$ to get $\Delta(g^{\dag}) \in \RR$. 
	
	The usual conditions $ab^*, cd^*, c^*a, d^*b \in V_n$ turn into $\bar{a} \bar{c}^*, \bar{b} \bar{d}^*, \bar{b}^*\bar{a}, \bar{d}^*\bar{c} \in V_n$.
	These are obtained from the original conditions as a set by applying the Clifford conjugation.
	If moreover $\Ocal$ is Clifford-stable and $g\in\SL_2(\Ocal)$, then the entries $\bar a,\bar b,\bar c,\bar d$ of $g^{\dag}$ lie in $\overline{\Ocal}=\Ocal$, so $g^{\dag}\in\SL_2(\Ocal)$.
\end{proof}

The following lemma is used later when we do our Bianchi-Humbert theory for Clifford-Hermitian forms.
All of the conditions in this lemma are equivalent to saying that the pair $(\mu,\nu)$ is right unimodular.
\begin{lemma}\label{L:proper-inputs}
	Suppose that $\Ocal$ is Clifford-stable (Definition~\ref{def:clifford-stable}).
	The following are equivalent for a pair of elements in $(\mu,\nu)$. 
	\begin{enumerate}
		\item \label{I:left}$\exists \begin{pmatrix} \bar{\mu} & * \\ \bar{\nu} & * \end{pmatrix} \in \SL_2(\Ocal)$
		\item \label{I:right}$\exists \begin{pmatrix} *& \bar{\mu} \\ *&\bar{\nu} \end{pmatrix} \in \SL_2(\Ocal)$
		\item \label{I:top} $\exists \begin{pmatrix} \mu & \nu \\ * & * \end{pmatrix}\in \SL_2(\Ocal)$
		\item \label{I:bottom}$\exists \begin{pmatrix} * & * \\ \mu & \nu \end{pmatrix}\in \SL_2(\Ocal)$; i.e. $(\mu,\nu)$ is right unimodular. 
		\item \label{I:invert}$\exists \begin{pmatrix} * & * \\ -\nu & \mu \end{pmatrix}\in \SL_2(\Ocal)$
		\item \label{I:units} $\forall u \in \Ocal^{\times}, \exists \begin{pmatrix} u\mu & u\nu \\ * & * \end{pmatrix}\in \SL_2(\Ocal)$.
	\end{enumerate}
\end{lemma}
\begin{proof}
	These follow from simple matrix identities.
	Let $S= \begin{psmallmatrix} 0 & -1 \\ 1 & 0 \end{psmallmatrix}$, $g = \begin{psmallmatrix} a & b \\ c & d \end{psmallmatrix}.$
	\begin{itemize}
		\item We have $gS = \begin{psmallmatrix} b & -a \\ d & -c\end{psmallmatrix}$.
		Similarly, $gS \begin{psmallmatrix}
		-1& 0 \\
		0 & -1 
		\end{psmallmatrix} = \begin{psmallmatrix} -b & a \\ -d & c \end{psmallmatrix}$.
		This gives the equivalence between \eqref{I:left} and \eqref{I:right}.
		\item For the fact that \eqref{I:left} and \eqref{I:right} are equivalent to \eqref{I:top} and \eqref{I:bottom} we use Lemma~\ref{L:closed-under-conj} which shows that the group is closed under Clifford adjunction. 
\item Assertion \eqref{I:invert} follows from the identity $g(-S) = \begin{psmallmatrix} -b & a \\ -d & c \end{psmallmatrix}$, whose bottom row is $(-d,c)$. 
\item The assertion about units comes from the identity $\begin{psmallmatrix} u & 0 \\
		0 & u^{*-1} \end{psmallmatrix} g = \begin{psmallmatrix} ua & ub \\ u^{*-1}c & u^{*-1} d \end{psmallmatrix}.$
	\end{itemize}
\end{proof}

\subsection{Clifford-Euclidean Rings} \label{sec:clifford-euclidean}

Recall that $R^\mon$ denotes the Clifford monoid of an order of a Clifford algebra $R$
(Definition~\ref{defn:CliffordGroup}).

\begin{definition}\label{def:clifford-euclidean}
Let $R$ be a $*$-stable order (Definition \ref{def:clifford-stable})
  in the rational Clifford algebra of a positive definite quadratic form; then $0$ is the only zero divisor in $R^{\mon}$, and the two descriptions of $R^{\mon}$ (via Definition~\ref{defn:CliffordGroup}, and as $R\cap(K^{\times}\cup\lbrace 0\rbrace)$ with $K^{\times}$ the Clifford group of Definition~\ref{def:clifford-group-anisotropic}) agree. We say that $R$ is {\em right Clifford-Euclidean} if there is a norm function $N: R^\mon \to \NN$ such that
	\begin{enumerate}
		\item $N(x) = 0$ if and only if $x$ is a zero-divisor.
		\item \label{I:division-algorithm} For all $x,y \in R^\mon$ with $N(x) > 0$ and $xy^* \in \Vec(R)$, there exists some $q\in \Vec(R)$ and some $r\in R^\mon$ such that 
		\begin{equation}\label{eqn:left-division}
		y=xq+r
		\end{equation}
		where $N(r)<N(x)$.
	\end{enumerate}
	If condition~\ref{I:division-algorithm} holds under the further assumption $xR + yR = R$, we say that $R$ is {\em weakly Clifford-Euclidean}. 
\end{definition}

We now make a remark on the terminology of ``left'' vs ``right'' Euclidean and the connection to the theory of quasi-Euclideanity.
\begin{remark}\label{rem:left-vs-right-division}
	\begin{enumerate}
	\item Note that it is well established that
	$x^{-1} y$ is ``right division by $x$'' and not ``left division by $x$''.
	The survey by Ahmadi, Jain, Lam, and Leroy on noncommutative Euclidean rings \cite{Alahmadi2014}, Brung's paper \emph{Left Euclidean Rings} \cite{Brungs1973}, and Voight's book \cite{Voight2020} all use this convention. 

	This unfortunately means ``left multiplication by the inverse of $x$'' is ``right division by $x$'' and ``right division by $x$'' is reversed by ``left multiplication by $x$''. 
	On the other hand, the right division algorithm will end up giving generators for right ideals. 
	This correspondence between left ideals and left division is perhaps a good mnemonic for remembering this convention.
	\item 	It is possible to develop the theory of Clifford-Euclidean rings by bootstrapping from the well-developed theory of quasi-Euclidean rings which is reviewed in \cite{Alahmadi2014}.
	We state the definition for the interested reader. 
	Let $R$ be an associative ring.
	A pair $(a,b)\in R^2$ is \emph{right quasi-Euclidean} if and only if there exists a sequence of elements $(q_1,r_1), (q_2,r_2), \ldots, (q_{n+1},r_{n+1})$ such that $r_{i-1} = r_i q_{i+1} + r_{i+1}$ for $0\leq i\leq n$ with $r_{-1}=a$, $r_0=b$, and $r_{n+1}=0$.
	If all pairs $(a,b)$ are right quasi-Euclidean, then we say $R$ is \emph{right quasi-Euclidean}. 
Assume that $0$ is the only zero divisor in $R^{\mon}$ (for example, $q$ positive definite), the hypothesis of Algorithm~\ref{alg:gcd}.
If $R$ is a right Clifford-Euclidean order in a rational Clifford algebra $K$, then every $(a,b)\in (R^{\mon})^2$ with $b^{-1}a \in \Vec(K)$ is right quasi-Euclidean. 
	If $R$ is weakly Clifford-Euclidean, then every right unimodular pair $(a,b)$ is right quasi-Euclidean.  
	Algorithm~\ref{alg:gcd} proves these statements. In neither case is
	the converse obvious, because the definition of right quasi-Euclidean allows
	the use of any elements of the ring as partial quotients and remainders,
	rather than only vectors and monoid elements respectively required by the definition
	of Clifford-Euclidean.
	\end{enumerate}
\end{remark}

In the formula for right division we use the notation \eqref{eqn:left-division}
$$ \quo_{\R}(y,x)= q, \quad  \rem_{\R}(y,x)=r.$$
There is also a notion of \emph{left division} where we are concerned with approximating $yx^{-1}$ rather than $x^{-1}y$.  In this situation we will have 
$ y=qx+r$ and use the notation $\quo_{\L}(y,x)=q$, $\rem_{\L}(y,x) = r.$

A general principle is that the involution $*$ interchanges left and right for operations such as quotient, remainder, GCD, and coefficients of GCD.
For example, we have  $\quo_{\L}(y,x)^* = \quo_{\R}(y^*,x^*),$ and $\rem_{\L}(y,x)^*=\rem_{\R}(y^*,x^*) $
where we use $(yx^{-1})^*=(x^*)^{-1}y^*$ (it may be useful to recall that taking inverse commutes with all of the involutions).
This implies that if $\Ocal$ is right Clifford-Euclidean then $\Ocal^*$ is left Clifford-Euclidean. 
Since Definition~\ref{def:clifford-euclidean} requires $\Ocal^*=\Ocal$, a right Clifford-Euclidean order is also left Clifford-Euclidean, for the norm $x\mapsto N(x^*)$, so we simply call it Clifford-Euclidean.

In what follows we will use that the algorithm for greatest common divisors makes sense for right and left division. 
\begin{algorithm}[GCD]\label{alg:gcd}
	Let $R$ be Clifford-Euclidean, and assume that $0$ is the only zero-divisor in $R^{\mon}$ (for example, $q$ positive definite), so that $N(r)=0$ forces $r=0$.
Since $R=R^*$ (and $R^{\mon}$ is $*$-stable, as when $q$ is positive definite), applying $*$ to the right division of Definition~\ref{def:clifford-euclidean} gives the left division $y=qx+r$ used below.
Given $a,b\in R^{\mon}$ with $a^*b$ a vector\footnote{By the Useful Lemma, this is equivalent to $ab^{-1}$ being a vector.} we compute the element $g\in R$ such that $Ra + Rb=Rg$. 
	Set $r_{-1}=a$ and $r_0=b$ and define inductively $q_j$ and $r_j$ for $j\geq 1$ by 
	\begin{equation}\label{eqn:gcd-computation}
	r_{j-2} = q_j r_{j-1} + r_j.
	\end{equation}
	(At $j=1$ we get $a=q_1b+r_1$.)
	Eventually this algorithm produces some $n$ with $r_{n+1}=0$ so that $r_n$ is an element such that $Rr_n=Ra+Rb$.
	We call this element the \emph{left gcd} and denote it by $\gcd_{\L}(a,b)$. 
	$$ \Ocal a + \Ocal b = \Ocal r_n,  \qquad \gcd_{\L}(a,b) = r_n.$$ 
	It is useful to note that in \eqref{eqn:gcd-computation} we always have $q_j \in \Vec(K)$ for $j=1,\ldots,n$ and $r_jr_{j-1}^{-1} \in \Vec(K)$ for $j=0,\ldots,n+1$.
	$\vartriangle$
\end{algorithm}

There is similarly a right GCD which we will denote by $\gcd_{\R}(a,b)$. 
The two algorithms and outputs are related by $ \gcd_{\R}(a,b)^* = \gcd_{\L}(a^*,b^*).$

\begin{algorithm}[GCD Coefficients]\label{alg:gcd-coefs}
Let $R$ be Clifford-Euclidean, under the hypotheses of Algorithm~\ref{alg:gcd}, and assume that every nonzero element of $\Vec(K)$ is invertible (for example, $q$ positive definite).
We introduce an algorithm
  $\gcdcoeffs_{\L}$ that accepts $a,b\in R^{\mon}$ with $ab^{-1}$ a vector
  and returns Clifford monoid elements $c,d$ such that
  $c a+ d b = \gcd_{\L}(a,b).$
  To do this we use \eqref{eqn:gcd-computation} to write the last nonzero remainder $r_n$ in terms of $a$ and $b$ in a way very similar to the classical algorithm.

	If $n=0$, that is, if $r_1=0$, then $\gcd_{\L}(a,b)=b$ and we take $c=0$ and $d=1$. Suppose now that $n\geq 1$.
	We introduce variables $c_j$ and $d_j$ to write
	\begin{equation}\label{eqn:coefficient-recurrence}
	r_n = c_j r_{n-j-1}+d_j r_{n-j}.
	\end{equation}
	When $j=1$ we have $r_n = r_{n-2} - q_nr_{n-1}$,
	which gives us $c_1=1$ and $d_1 = -q_n$.
	We also get a recurrence. The equation \eqref{eqn:gcd-computation} written with the index $n-j$ gives $r_{n-j} = r_{n-j-2} - q_{n-j} r_{n-j-1}$.
	Substituting this into \eqref{eqn:coefficient-recurrence} gives $ r_n = c_jr_{n-j-1} + d_j(r_{n-j-2}-q_{n-j}r_{n-j-1}) = d_jr_{n-j-2} + (c_j-d_jq_{n-j})r_{n-j-1}$,
	which gives
	\begin{equation}
	c_{j+1} = d_j, \quad d_{j+1} = c_j -d_jq_{n-j}, \qquad 1 \leq j \leq n-1.
	\end{equation}
	Continuing with the case $j=n$ gives $\gcd_{\L}(a,b)=r_n = c_nr_{-1}+d_n r_0=c_n a+d_n b.$
	So $c=c_n$, $d=d_n$ and the algorithm $\gcdcoeffs_{\L}$ gives outputs defined by the expression below
	$$ \gcdcoeffs_{\L}(a,b) = (c,d), \qquad ca+db = \gcd_{\L}(a,b).$$
	$\vartriangle$
\end{algorithm}

Similarly, we can define an algorithm $\gcdcoeffs_{\R}$ that accepts
$a, b \in R^\mon$ with $a^{-1}b \in \Vec(K)$ and returns $c, d \in R^\mon$
such that $ac+bd = \gcd_{\R}(a,b)$.
As usual we have $\gcdcoeffs_{\R}(a,b)^* = \gcdcoeffs_{\L}(a^*,b^*).$

One of the things that we are interested in is determining when an equation $ad^* + bc^* =1 $
can be lifted to a matrix $\begin{psmallmatrix} a & b \\ c & d \end{psmallmatrix}$.
Lemma~\ref{eqn:ratios-of-coeffs} is useful for this, as we will need to check when appropriate ratios of columns and rows are vectors.

The definition of $\GL_2(\CC_n)$ says
that we only need to check $g^{-1}(0),g^{-1}(\infty) \in V_n$ and that $g(0),g(\infty)\in V_n$ come for free.
This is a simple but non-obvious algebraic manipulation that is worth knowing.

\begin{lemma}\label{lem:just-need-two}
	Let $\Ocal$ be a $*$-stable order and let
	$g=\binom{a \ \ b}{c \ \ d} \in M_2(\Ocal)$ with $a,b,c,d\in \Ocal^{\mon}$, $\Delta=ad^*-bc^*=1$.
	If $a,c \neq 0$ and $a^{-1}b, c^{-1}d\in V_n$ then $g \in \SL_2(\Ocal)$.
\end{lemma}
\begin{proof}
	By Theorem~\ref{thm:gl2-short} it suffices to show that $ab^*, dc^* \in V_n$.
	If $b=0$ then $ab^*=0$.
	Otherwise $a^*$ and $b^*$ are invertible and $(a^{-1}b)^* = b^*(a^*)^{-1} \in V_n$, so the Useful Lemma~\ref{lem:useful} gives $ba^* \in V_n$, hence $ab^*=(ba^*)^* \in V_n$.
	In the same way $c^{-1}d \in V_n$ gives $cd^* \in V_n$, hence $dc^* = (cd^*)^* \in V_n$.
\end{proof}

\begin{lemma}\label{eqn:ratios-of-coeffs}
	For $c$ and $d$, as above with $ca+db=\gcd_{\L}(a,b)$, we have $c,d \in R^{\mon}$ and $cd^* \in \Vec(K)$. In particular, $c^{-1}d$ is a vector if $c \neq 0$.
\end{lemma}
\begin{proof}
	If $n=0$ then $(c,d)=(0,1)$ and there is nothing to prove, so let $n \geq 1$.
	Recall that the nonzero elements of $R^{\mon}$ are invertible in $K$ (Corollary~\ref{cor:no-zero-divisors}), that $\Vec(K) \subset K^{\mon}$, and that $K^{\mon}$ is closed under multiplication.
	We show by induction on $j$ that $c_j,d_j \in R^{\mon}$ and $c_jd_j^* \in \Vec(K)$.
	For $j=1$ we have $c_1=1$, $d_1=-q_n \in \Vec(R)$ and $c_1d_1^*=-q_n^* \in \Vec(K)$.
	Suppose that the claim holds for $j$. Then $c_{j+1}=d_j \in R^{\mon}$ and
	$$c_{j+1}d_{j+1}^* = d_j(c_j-d_jq_{n-j})^* = (c_jd_j^*)^* - d_jq_{n-j}^*d_j^*,$$
	which is a vector because $\Vec(K)$ is closed under $*$ and $d_j \Vec(K) d_j^* \subset \Vec(K)$ by the definition of the Clifford monoid.
	If $d_j=0$ then $d_{j+1}=c_j \in R^{\mon}$.
	Otherwise $d_j$ is invertible and $d_{j+1}=d_j(d_j^{-1}c_j-q_{n-j})$, where $d_j^{-1}c_j \in \Vec(K)$: this is clear if $c_j=0$, and otherwise $(d_j^{-1}c_j)^*=c_j^*(d_j^*)^{-1}$ is a vector by the Useful Lemma~\ref{lem:useful}, because $(c_j^*)^*d_j^*=c_jd_j^*$ is a vector.
	Hence $d_{j+1}$ is the product of an element of $R^{\mon}$ and a vector, so $d_{j+1} \in R \cap K^{\mon} = R^{\mon}$.
	Finally, if $c \neq 0$ and $d \neq 0$ then $c^{-1}d \in \Vec(K)$ by the same application of the Useful Lemma, and if $d=0$ then $c^{-1}d=0$.
\end{proof}

\begin{remark}\label{rem:still-vec}
	Recall from Corollary~\ref{cor:no-zero-divisors} that, in the case of
	most interest to us where $R$ is an order in $\quat{-d_1,\dots,-d_n}{\Q}$
	(where the $d_i$ are positive integers), the only right zero-divisor in
	$R^\mon$ is $0$. We also point out that \eqref{I:division-algorithm} imposes
	no condition
	if $xy^* \notin \Vec(R)$ or if $N(x)=0$. If $xy^* \in \Vec(R)$ and
	such $q, r$ do exist, then it follows that $xr^* \in \Vec(R)$, because
	$xr^* = x(y-xq)^* = xy^* - xqx^*$ and $(xqx^*)^*$ is a vector when $q$ is.
\end{remark}

\begin{proposition}\label{prop:weak-euclidean-implies-weak-unimodular}
Suppose $\Ocal$ is a Clifford-Euclidean order stable under $*$.  Then $(\mu,\nu)$ is right unimodular (there exists an element of $\begin{psmallmatrix} * & * \\
\mu & \nu \end{psmallmatrix} \in \SL_2(\Ocal)$) if and only if $\mu^{-1}\nu$ is a Clifford vector, $\mu,\nu \in \Ocal^{\mon}$ and $\mu \Ocal + \nu \Ocal = \Ocal$ as right ideals. 
\end{proposition}
\begin{proof}
		Suppose $(\mu,\nu)$ is right unimodular. Then by the matrix condition $\mu,\nu \in \Ocal^{\mon}$, $\mu^{-1}\nu $ is a Clifford vector since elements of $\SL_2(\Ocal)$ induce M\"{o}bius transformations, and $a\nu^*-b\mu^*=1$ for some $a,b\in \Ocal$. Taking $*$ gives the implication. 
	
	Conversely, suppose that $\mu,\nu \in \Ocal^{\mon}$, that $\mu^{-1}\nu$ is a Clifford vector and that $\mu\Ocal + \nu\Ocal = \Ocal$.
	By the Useful Lemma $\mu\nu^*$ is a vector.
	Since $\Ocal = \Ocal^*$, the order $\Ocal$ is right Clifford-Euclidean (Definition~\ref{def:clifford-euclidean}); let $N$ be its norm.
The defining form $q$ is positive definite (Definition~\ref{def:clifford-euclidean}); in an orthogonal basis of its $\QQ$-vector space it has positive diagonal coefficients, so it is strongly anisotropic over $\QQ$ (Proposition~\ref{prop:big-form-02}(2), with $R=\QQ$), so by Corollary~\ref{cor:no-zero-divisors} the only $x \in \Ocal^{\mon}$ with $N(x)=0$ is $x=0$, and every nonzero element of $\Ocal^{\mon}$ is invertible in $K$.
For $x \in \Vec(\Ocal)$ let $T_x = \begin{psmallmatrix} 1 & x \\ 0 & 1 \end{psmallmatrix}$ and let $S = \begin{psmallmatrix} 0 & -1 \\ 1 & 0 \end{psmallmatrix}$; by Theorem~\ref{thm:gl2-short} both lie in $\SL_2(\Ocal)$.
	We run the right division algorithm on the row $(x,y)=(\mu,\nu)$ using these matrices.
	If $x \neq 0$, then $xy^* \in \Vec(\Ocal)$, so there are $q \in \Vec(\Ocal)$ and $r \in \Ocal^{\mon}$ with $y = xq+r$ and $N(r)<N(x)$, and
	$$(x,y)\,T_{-q}S = (x,r)S = (r,-x).$$
	By Remark~\ref{rem:still-vec} we have $xr^* \in \Vec(\Ocal)$, so the new row again satisfies $r(-x)^* = -(xr^*)^* \in \Vec(\Ocal)$.
	Since $N$ decreases, after finitely many steps we reach $x=0$, so that $(\mu,\nu)M = (0,u)$ for some $M \in \SL_2(\Ocal)$ and some $u \in \Ocal^{\mon}$.
	Then $(\mu,\nu) = (0,u)M^{-1}$, so $\Ocal = \mu\Ocal+\nu\Ocal \subseteq u\Ocal$ and $u$ is a unit of $\Ocal$.
	As $\Ocal$ is $*$-stable, $(u^*)^{-1} = (u^{-1})^* \in \Ocal$, and $D = \begin{psmallmatrix} (u^*)^{-1} & 0 \\ 0 & u \end{psmallmatrix} \in \SL_2(\Ocal)$ by Theorem~\ref{thm:gl2-short}.
	The matrix $g = DM^{-1} \in \SL_2(\Ocal)$ has bottom row $(0,u)M^{-1} = (\mu,\nu)$.
\end{proof}

If the order is compatible and Clifford-Euclidean, we can use the same method as in the classical setting to prove that all cusps are unimodular.

\subsection{Clifford-Principal Ideal Rings}
In this section we sort out some consequences of Clifford-Euclideanity.

\begin{definition}\label{def:cpir}
  Let $R$ be an order in a Clifford algebra $K$. We say
	that $R$ is {\em right cuspidally principal} if for all $v \in \Vec(K)$ there exists some
	$\binom{a \ \ b}{c \ \ d} \in \SL_2(R)$
	with $c^{-1}d = v$. 
\end{definition}
We will just abusively say that $R$ is \emph{cuspidally principal} or \emph{principal} for convenience. 
Note that this is saying that all of the cusps $v\in \Vec(K) \cup \lbrace \infty \rbrace$ are principal (in the sense that they are in the orbit of $\infty$ under $\PSL_2(R)$). 
The unimodular pair for $\infty$ is by convention $(c,d) = (0,1)$. 

\begin{corollary}\label{cor:Euclidean-implies-unimodular}
If $R$ is a right Clifford-Euclidean order (so $q$ is positive definite and $0$ is the only zero divisor in $R^{\mon}$, Definition~\ref{def:clifford-euclidean}) then it is right cuspidally principal. 
	A similar statement holds for left Clifford-Euclidean orders. 
\end{corollary}

\begin{proof}[Proof]
Let $w\in\Vec(K)$; we find $\binom{a\ \ b}{c\ \ d}\in\SL_2(R)$ with $c^{-1}d=w$. If $w=0$, then $\begin{psmallmatrix}0&-1\\1&0\end{psmallmatrix}$ will do. Otherwise choose a positive integer $m$ with $mw\in R$ and put $x=m$ and $y=mw$. Then $x,y\in R^{\mon}$ ($y$ is a nonzero vector and $q$ is positive definite), $xy^*=m^2w$ is a vector, and $x^{-1}y=w$.
	Given such $x, y \in R^\mon$ not both $0$ with $xy^*$ a vector, we start by
finding $v \in R$ such that $xR + yR = vR$. To do so, we use the
	Clifford-Euclidean algorithm. Choose $q_0$ such
	that $y = xq_0 + r_0$ with $N(r_0) < N(x)$; now let $r_{-2} = y, r_{-1} = x$.
	Inductively define $q_i, r_i$ by choosing $q_i$ such that
	$r_{i-2} = r_{i-1} q_i + r_i$ with $N(r_i) < N(r_{i-1})$ until $r_i = 0$
	is reached. At each step we have
	$$r_{i-2}R + r_{i-1}R = (r_{i-1}q_i + r_i)R + r_{i-1}R = r_{i-1}R + r_iR,$$
	so at the end we have $yR + xR = r_{i-1}R$ and we choose $v = r_{i-1}$.
	Then $x, y \in vR$, so writing $x = vc, y = vd$ we have
	$x^{-1}y = c^{-1}v^{-1}vd = c^{-1}d$, and since $xR + yR = vR$ we may write
	$v = xe + yf = vce + vdf$. Because $v$ has nonzero norm it is not a
	zero-divisor and so $ce + df = 1$. This does not quite give the
	desired $M \in \SL_2(\Ocal)$, because $e$ and $f$ need not belong to the Clifford monoid,
but it shows that $cR + dR = R$. Moreover $c=v^{-1}x$ and $d=v^{-1}y$ lie in $R^{\mon}$ because $v$ is a monoid element, and $c^{-1}d=x^{-1}y$ is a vector because $xy^*$ is (Useful Lemma), so we have reduced to Proposition~\ref{prop:weak-euclidean-implies-weak-unimodular}. 
	
\end{proof}

\begin{definition} Let $R$ be an order in a Clifford algebra.
	If all right ideals of $R$ generated by elements of
	$R^\mon$ are generated by a single element of $R^\mon$
	we say that $R$ is {\em Clifford-principal}.
\end{definition}

In \S\ref{sec:maximal-c4} we give an example, in $\quat{-1,-1,-1}{\QQ}$, of a Clifford-Euclidean (hence cuspidally principal) order that is not Clifford-principal.
We do not know whether there are any Clifford-principal
orders in Clifford algebras with $3$ or more imaginary units, nor whether
a Clifford-principal order is always cuspidally principal
(a negative answer to the first of these implies a positive answer to the
second).

If an order is not Clifford-Euclidean, we cannot conclude that
there is more than one equivalence class of cusps: just as commutative
principal ideal domains are not Euclidean in general, cuspidally principal
orders need not be Clifford-Euclidean.
Indeed, for $R$ an order in an imaginary quadratic
field, all of the Clifford conditions reduce to the ordinary conditions, since
all elements are Clifford vectors; it is well-known that the maximal orders
of the quadratic fields of discriminant $-19, -43, -67, -163$ are
principal but not Euclidean. We announce that it is also possible for a
maximal order in a quaternion algebra to be principal but not Euclidean
(we expect the same method of proof to be applicable
in larger Clifford algebras as well but do not have any examples).
This is presented in a forthcoming manuscript \cite{DHIL2}.

We have been careful to work out examples that prove that no two of the definitions are equivalent. 
 
\begin{example}
	\begin{enumerate}
		\item We show that $\ZZ[i,j]$ and the Hurwitz order $\Ocal_3$ in $\quat{-1,-1}{\QQ}$ are both Clifford-Euclidean.
		\item We have proved that no maximal order in $\quat{-2,-13}{\QQ}$ is Clifford-Euclidean. This will be presented in a subsequent paper \cite{DHIL2}.
		\item In $\quat{-1,-1,-1,-1}{\QQ}$ we have two maximal orders which are conjugate but whose underlying lattice and Clifford groups are not the same.
	\end{enumerate}
\end{example}
The property of being Clifford-Euclidean should not be expected to be isomorphism-invariant because the geometry of the set of integral Clifford vectors need not be preserved under conjugation. 
In \S\ref{sec:1-1-3} and \S\ref{sec:1-1-1-1}, on orders in $\quat{-1,-1,-3}{\QQ}$ and $\quat{-1,-1,-1,-1}{\QQ}$ respectively, we show that there exist isomorphic orders $\Ocal$ and $\Ocal'$ with lattices $\Vec(\Ocal)$ and $\Vec(\Ocal')$ which are not isometric (see Proposition~\ref{rem:not-invariant}). 

\subsubsection{The Covering Radius and Clifford-Euclideanity}\label{sec:covering-radius}

\begin{definition}\label{def:covering-radius}
The \emph{covering radius} (\cite[1.2]{SPLAG}) of a lattice $\Lambda \subset \RR^n$ is the maximal distance from a point in $\RR^n$ to a lattice point: $ \rho(\Lambda) = \sup_{x\in \RR^{n}} \min \lbrace \vert x-\lambda\vert \colon \lambda \in \Lambda \rbrace. $
A \emph{hole} is a local maximum of the function $x\mapsto \min_{\lambda \in \Lambda} \vert x - \lambda \vert$, and a \emph{deep hole} is a global maximum of the function.
\end{definition}
For example, for the standard cubic lattice $\ZZ^n$ the vector $(1/2,1/2,\cdots,1/2)$ is the unique hole up to translation and it is deep. For typical lattices there are holes
that are not deep and multiple translation orbits of deep hole.

\begin{theorem}\label{thm:small-covering-radius-euclidean}
	Let $K \subset \CC_n$ be a rational Clifford algebra with a $*$-stable order $\Ocal \subset K$ associated to a positive definite quadratic form such that for $x \in \Vec(K)$ we have $N(x) = \vert x \vert^2$.
	Let $\Lambda = \Vec(\Ocal)$ be the integral Clifford vectors. 
	If the covering radius of $\Vec(\Ocal)$ is less than $1$, then $\Ocal$ is Clifford-Euclidean for the norm $N$.
\end{theorem}
\begin{proof}
	Suppose that $a,b\in \Ocal^{\mon}$ such that $ba^*\in \Vec(K)$. 
	Let $q \in \Vec(\Ocal)$ be an element closest to $a^{-1} b$.
	Then $a^{-1}b-q =r_0 \in \Vec(K)$ has norm less than one. 
	This implies that $b=aq+r$ with $r=ar_0$. 
	We have $\nrd(r)=\nrd(ar_0) = \nrd(a)\nrd(r_0)<N(a)$.
\end{proof}

\subsection{Codes and Lattices in $\quat{-1,-1,\ldots,-1}{\QQ}$}\label{sec:codes}
A \emph{binary code} is just a  vector space $C \subset \FF_2^n$. 
Elements of this space are called \emph{codewords}. 
The dimension of the code is $\dim_{\FF_2}(C)$.
The \emph{length} of the code is the space of the ambient dimension $n$.
We will write codewords as binary strings $b_1\cdots b_n$ where $b_j \in \lbrace 0,1\rbrace$.
In this way we also identify codes with subsets of $\lbrace 0, 1 \rbrace^n$ and allow ourselves to talk about $c \cdot v$ for $v\in \ZZ^n$.

The \emph{weight} $\wt(w)$ of a word is the number of nonzero entries, and the \emph{Hamming distance} between two codewords $v,w$ is $\wt(v-w)$. 
The \emph{minimal distance} of a code $C$ is then $\wt(v-w)$ where $v$ and $w$ range over elements of $C$. 
Binary codes are often classified by $[n,k,d]$ where $n$ is the length, $k$ is the dimension, and $d$ is the minimal distance (see \cite[Ch3, \S 2]{SPLAG}).
A code is \emph{even} if every codeword has even weight. 
The simplest such example that is not trivial is the code spanned by $11$ in $\FF_2^2$. 
A code is \emph{doubly even} if every codeword has weight divisible by $4$.

Orders in $\quat{-1,-1,\ldots,-1}{\QQ}$ correspond to various codes, and this section elaborates on this relationship. 
This relationship between codes and orders in Clifford algebras also appears in \cite[Section 6]{Iga2021}. 

Let $\Ocal$ be an order such that $\ZZ[i_1,i_2,\ldots,i_{n-1}] \subset \Ocal \subset \QQ[i_1,i_2,\ldots,i_{n-1}].$ 
Let $I = (i_0,i_1,\ldots,i_{n-1})$.
For such $\Ocal$ there exists a doubly even code $C$ with
$\Vec(\Ocal) = \frac{1}{2}\lat_C$, where $\lat_C$ is the lattice defined below;
the construction is described in Lemma~\ref{lem:order-to-code}.
The order $\Ocal$ contains $\ZZ[\Vec(\Ocal)]$, but it need not be equal to it (consider the Hurwitz order).
The lattice $\lat_C$ associated to the binary code $\code$ is the inverse image of $\code$ under the natural map $\ZZ^n \to \FF_2^n$ and $\Vec(\Ocal) = \frac{1}{2}\lat_C$. 
Figure~\ref{fig:codes-lattices}(a) shows this construction for $n=2$.
\begin{figure}[htbp!]
\centering
\begin{tikzpicture}[scale=0.7,lab]
 \foreach \x in {-2,...,2}{\foreach \y in {-2,...,2}{
   \pgfmathparse{int(mod(\x+\y+10,2))}
   \ifnum\pgfmathresult=0 \fill (\x,\y) circle (3pt); \else \draw (\x,\y) circle (3pt); \fi}}
 \draw[gray,->] (-2.6,0)--(2.6,0); \draw[gray,->] (0,-2.6)--(0,2.6);
 \node at (0,-3.2){(a) $C=\{00,11\}\subset\FF_2^2$: $\Lambda_C$ = filled points};
\end{tikzpicture}\qquad
\begin{tikzpicture}[scale=2.0,lab,x={(1cm,0cm)},y={(0.45cm,0.35cm)},z={(0cm,1cm)}]
 \foreach \a/\b/\c in {0/0/0,1/0/0,0/1/0,1/1/0,0/0/1,1/0/1,0/1/1,1/1/1}{\coordinate (v\a\b\c) at (\a,\b,\c);}
 \fill[domC!35] (v000)--(v100)--(v110)--(v010)--cycle;
 \fill[bubbleC!25] (v000)--(v110)--(v111)--(v001)--cycle;
 \fill[holeC!20] (v000)--(v110)--(v101)--cycle; \fill[holeC!20] (v000)--(v011)--(v101)--cycle;
 \draw[holeC,thick] (v000)--(v110)--(v101)--(v000)--(v011)--(v110) (v011)--(v101);
 \foreach \p/\q in {000/100,100/110,110/010,010/000,001/101,101/111,111/011,011/001,000/001,100/101,110/111,010/011}{\draw (v\p)--(v\q);}
 \foreach \p in {000,100,010,110,001,101,011,111}{\fill (v\p) circle (1.2pt);}
 \node[below left] at (v000){$000$}; \node[above right] at (v111){$111$};
 \node at (0.5,0,-0.45){(b) $H(8,4)$ on $\FF_2^3$};
\end{tikzpicture}
\caption{(a) A binary code $C\subset\FF_2^n$ gives the lattice $\Lambda_C\subset\ZZ^n$, the preimage of $C$; here $n=2$ and $C=\{00,11\}$, so $\Lambda_C$ is the checkerboard lattice, and Lemma~\ref{lem:order-to-code} gives $\Vec(\Ocal)=\frac12\Lambda_C$. (b) Index the $8$ coordinates of $H(8,4)$ by the points of $\FF_2^3$. Its $14$ words of weight $4$ are the $14$ affine planes (since $H(8,4)$ is the Reed--Muller code $\mathrm{RM}(1,3)$, \cite[Ch.~3]{SPLAG}): $6$ faces (orange shows one), $6$ diagonal rectangles (blue shows one), and the $2$ tetrahedra $x_1+x_2+x_3=0,1$ (purple shows one).}\label{fig:codes-lattices}
\end{figure}
The lattice $\Lambda = \Vec(\Ocal)$ is then also spanned by a doubly even code in the sense that $\Lambda$ is a cubic lattice coming from $\ZZ[i_1,\ldots,i_{n-1}]$ together with vectors of the form $(c \cdot I)/2$.
Direct sum of codes $C_1\oplus C_2$ corresponds to direct sum of lattices $\Lambda_{C_1} \oplus \Lambda_{C_2}$, and the basic objects $\Lambda_C$ in this theory are $A_1$ (empty column), $D_{2n}$ (make $n$ pairs of columns, then consider the code whose words are $1$ on an even number of pairs), $E_7$ (the $[7,3,4]$ simplex code), and $E_8$ (extended Hamming code $H(8,4)$), each up to scaling by $\sqrt 2$.

\begin{lemma}\label{lem:order-to-code}
	For every order $\Ocal$ in $\QQ[i_1,\ldots,i_{n-1}]$ containing all the $i_j$, we have $\Vec(\Ocal) = \frac{1}{2}\Lambda_C$ for some doubly even binary code of length $n$. 
	More precisely, if $v\in \Vec(\Ocal)$ then $v$ takes the form 
	$$v = \frac{c_0 + c_1 i_1 + \cdots + c_{n-1} i_{n-1}}{2}$$
	where $c_j \in \ZZ$ for $j=0,\ldots,n-1$.
	Moreover, $\Vec(\Ocal)$ is generated by $\ZZ+\ZZ i_1+\cdots+\ZZ i_{n-1}$ together with the vectors $(c\cdot I)/2$ for $c \in C$, where $C\subseteq \lbrace 0, 1 \rbrace^n = \FF_2^n$ is a doubly even code.
	These make the codewords of our code.
\end{lemma}
\begin{proof}
	Indeed, if $a_j$ is the $i_j$-coefficient of $v$ (so $c_j=2a_j$),
	then $a_j$ is the $1$-coefficient of $w=\overline{i_j}v=-i_j v$. Since $\overline{u}v+\overline{v}u=2B(u,v)$ for Clifford vectors $u,v$, where $B$ is the polarization of $\nrd$, we have $w+\overline{w}=2a_j$ and $w\overline{w}=\overline{w}w=\nrd(v)$, so $w\in\Ocal$ satisfies $x^2-2a_jx+\nrd(v)$.
	Since $w$ is integral, $2a_j\in\ZZ$ (if $w\notin\QQ$ this is its minimal polynomial; if $w\in\QQ$ then $w=a_j\in\ZZ$).  
	
	By translation by elements of the standard cubic lattice $\Lambda_0=\ZZ + \ZZ i_1 + \cdots + \ZZ i_{n-1}$ we can assume that any $v\in \Vec(\Ocal)$ is congruent modulo $\Lambda_0$ to an element $(1/2)\sum_{j=0}^{n-1} c_j i_j$ with $ c_j \in \lbrace -1,0,1 \rbrace$.  
	By integrality, $(c_0^2+c_1^2 + \cdots + c_{n-1}^2 )/4 \in \ZZ$ which means $4$ divides $\wt(c)$ where $c$ is the codeword determined by squaring all of the coefficients. 
\end{proof}

An example below shows that this is not a bijection. 

\begin{theorem}\label{thm:clifford-euclidean-orders}
	If $\Ocal$ is $*$-stable, $\Vec(\Ocal) = \frac{1}{2}\Lambda_C$ and  $\rho(\Lambda_C)<2$, then $\Ocal$ is Euclidean.
\end{theorem}

\begin{proof}
  This follows from Theorem~\ref{thm:small-covering-radius-euclidean}.
\end{proof}

In some ways, this construction is not well-behaved. Given a code $\code$
we might define $\Ocal_\code = \Z[\frac{I \cdot c}{2}]$. However, if
$\code$ is the code of an order $\Ocal$, it does not follow that
$\Ocal = \Ocal_\code$ (consider the Hurwitz order).
In fact, it is not even clear
that $\Ocal_\code$ is an order at all, although we conjecture that it is
in Conjecture~\ref{conj:generate-order}.  See also Example~\ref{ex:golay-order}.

Here is a list of important examples that can be found elsewhere in the paper. 
\begin{example}
The lattice $D_4^*=\ZZ^4+(1/2,1/2,1/2,1/2)\ZZ$, the dual of $D_4$ (which is similar to $D_4$), is associated to the code generated by $1111$ and is the lattice of the maximal order $\Ocal_4 \subset \quat{-1,-1,-1}{\QQ}$. One calculates that this is the unique maximal order containing the Clifford order of
	$\quat{-1,-1,-1}{\QQ}$. The map from orders to codes of Lemma
	\ref{lem:order-to-code} is not injective: the Lipschitz and Hurwitz orders both have the code $\{000\}$, and $\Ocal_4$ and its suborder $\ZZ[i_1,i_2,i_3,(1+i_1+i_2+i_3)/2]$ of index $4$ both have the code generated by $1111$.
	We do not know whether it is surjective in general (compare Conjecture~\ref{conj:generate-order}).
	This order is explored in \S\ref{sec:1-1-1}.
\end{example}

\begin{example}
	There are five $[5,1,4]$ codes which appear for orders of
	$\quat{-1,-1,-1,-1}{\QQ}$. 
	These are the one-dimensional $\FF_2$-subspaces of $\FF_2^5$ generated by a word of weight $4$.
	These orders are all isomorphic; they are explored in \S\ref{sec:c5-reasonable}.
\end{example}

There are some interesting examples related to the extended Hamming code $C=H(8,4)$ and its associated lattice $\Lambda_C\cong\sqrt2 E_8$. 
This construction is related to the largest possible $n$ such that we can find an order in $\quat{(-1)^n}{\QQ}.$
\begin{example}
	Let $H(8,4)$ be the extended Hamming code (\cite[Chapter 3, 2.4.2]{SPLAG}).
	This code is a doubly even code $C=H(8,4)$ of dimension $4$ and length $8$.
Figure~\ref{fig:codes-lattices}(b) shows its words of weight $4$.
	The associated lattice $\Lambda_C$ is isometric to $\sqrt 2E_8$, so it has covering radius $\sqrt 2$ (the covering radius of $E_8$ is $1$).
	Similarly, the lattice associated to the code of length
	$9$, obtained from $H(8,4)$ by adding $0$ to the end of
	every word, has covering radius $\sqrt{3}$. 
	
	The direct sum of two copies of this code $H(8,4)\oplus H(8,4)$ has length $16$ and radius $2$. 
If $\Ocal = \ZZ[\frac{c\cdot I}{2} \colon c \in C]$ is an order (see Conjecture~\ref{conj:generate-order}), then it is almost a Clifford-Euclidean order, since the bound is not strictly less than two.
There does not exist a doubly even code $C$ such that $\Lambda_C$ has covering radius less than 2 if $n>16$ as the vector $(1/2,\dots,1/2)$ is at distance greater than $2$ from all integral vectors.
	
	If we require the covering radius to be strictly less than $2$, then
	the bound is $n \le 15$, but in fact we know of no examples of such codes
	with $n > 9$.  
	Our padded $H(8,4)$ serves as an example in dimension $9$.
\end{example}

\begin{example}\label{ex:golay-order}
	The extended Golay code $G_{24}$ of type $[24,12,8]$ is doubly even, and its associated lattice $\Lambda_{G_{24}}$ is $\sqrt{2}$ times the Niemeier lattice with root system $24A_1$; it is not the Leech lattice. One can consider the ring $\Ocal=\ZZ[\frac{c\cdot I}{2} \colon c \in G_{24}] \subset \quat{(-1)^{23}}{\QQ}$. 
	We do not know if $\Ocal$ is an order.
	In other words, if we start at $\CC_{24}$ and adjoin these
	halves of codewords, do we stop at something containing the standard order
	with finite index? It is not obvious that the denominators of
        products of codewords do not grow without limit.  
	If there is such an order, then $\Vec(\Ocal) = \frac{1}{2}\Lambda_{G_{24}}$, since the Golay code is maximal among doubly even codes of its length.
	A crude extrapolation tells us that it would take $200000$ years to define
	$\CC_{24}$ in \texttt{magma} with our current implementation, and hence checking this computationally is beyond our capabilities.
\end{example}

\begin{remark}
	There do not exist lattices $\Lambda_C$ of doubly even codes $C$ with $\rho(\Lambda_C)<2$ for $n>16$. 
	Every doubly even code is contained in one of maximal dimension, and there is a list of maximal-dimensional ones,
	so we can check Robert L. Miller's database of Doubly-Even Codes \cite[\href{https://rlmill.github.io/de_codes/}{Doubly-Even Codes}]{Miller2024} in small dimension.

	When analyzing the covering radius, we are talking about lattices $\Lambda_C$ that are spanned by $2\ZZ^n$ and the words of the code read as integral
	vectors.
	The covering radius of $\Lambda_C$ is not determined by the covering radius of the code (for $C=\FF_2^n$ the code has covering radius $0$, while $\Lambda_C=\ZZ^n$ has covering radius $\sqrt n/2$), but for $y\in\ZZ^n$ the squared distance from $y$ to $\Lambda_C$ is the Hamming distance from $y \bmod 2$ to $C$, so $\rho(\Lambda_C)^2$ is at least the covering radius of the code. (The covering radius of
	a lattice is the maximal distance from a point in the ambient space to
	a lattice point; the covering radius of a code is the maximum of the
	minimum Hamming weight of a coset.)  
	From Fact 3 in \cite{AssmusPless1983},
	the covering radius of a self-dual code is at least $d/2$ where $d$ is
	the minimum weight, with equality if and only if the code is an extended
	perfect code. But essentially the only perfect codes are Hamming codes and
	Golay codes, and the extended Hamming codes aren't doubly even beyond
	length $8$. 
	This suggests that there are no more easy examples of Euclidean orders coming from straightforward covering radius considerations.
\end{remark}

\section{Weil Restriction and Representability of $\SL_2(\Ocal)$}\label{sec:representability}

In this section we prove that $\SL_2(\Ocal)$ is the group of $\ZZ$-points of a $\ZZ$-group scheme. 
This is later used to conclude arithmeticity of this group as a subgroup of $\SO_{1,n+1}(\RR)$ in \S\ref{sec:arithmeticity}.

In \S\ref{sec:weil-restriction} we give a systematic treatment of Weil restriction for Clifford algebras, which establishes that for integral quadratic forms $q$ the group schemes $G=\SL_2(\Clf_q)$ are $\ZZ$-group schemes, and that $\SL_2(\Ocal)$ for $\Ocal = \Clf_q(\ZZ)$ are indeed $\ZZ$-points of group schemes.
From using our integral Bott periodicity developed in \S \ref{sec:arithmetic-bott} we show that for positive definite quadratic forms $q$ there exists an integral quadratic form $Q$ such that $\SL_2(\Clf_q) \cong \Spin_Q$ as $\ZZ$-group schemes (see Theorem~\ref{thm:mcinroy}).
After this, we use the Spin exact sequence $1\to\mu_2\to\Spin_Q\to\PSpin_Q\to 1$ in the fppf topology of $\Spec(\ZZ)$, where $\PSpin_Q=\Spin_Q/\mu_2$ is a $\ZZ$-group scheme mapping to $\SO_Q$.  Thus we are able to show that $\PSL_2(\Clf_q(\ZZ))$ coincides with an index-$2$ subgroup of $\PSpin_{Q}(\ZZ)$, whose image in $\SO_Q(\RR) =\SO_{1,n+1}(\RR)$ has finite index in $\SO_Q(\ZZ)$, and that $\PSL_2(\Ocal)$ is commensurable with $\PSL_2(\Clf_q(\ZZ))$ if $\Ocal$ is an order closed under Clifford involutions.  This establishes that $\PSL_2(\Ocal)$ is an arithmetic group. The morphism $\PSpin_Q\to\SO_Q$ is not an isomorphism over $\ZZ$ (Remark~\ref{rem:spin-sequence-at-2}).

Our methods and results have a nonempty intersection with the papers of Maclachlan-Waterman-Wielenberg \cite{Maclachlan1989} and Elstrodt-Grunewald-Mennicke in \cite{Elstrodt1988}, and we now make some remarks on this. 
\begin{remark}\label{rem:MWW-arithmeticity}
	First, we would like to make some clarifying remarks on the proof of arithmeticity of 
	$$\Gamma = \PSL_2(\ZZ[i_1,\ldots,i_{n-1}])$$ 
	in \cite[Theorem 4 and its Corollary on page 745]{Maclachlan1989}.
	In addition to the common practice of taking $\ZZ$-points of $\QQ$-group schemes (which is deprecated for reasons discussed in Appendix \ref{sec:group-schemes}), they view $\Gamma$ as a subgroup of two different Lie groups and use conclusions about arithmeticity in one to deduce arithmeticity for the other. 
	One of the subgroups is $\SO_{1,n+1}(\RR)^{\circ}$, which is related to the geometry of M\"obius transformations, and the second is a closed subgroup $G(\RR)$ of $\GL_{2^n}(\RR)$ (again they do not pay much attention to the ring of definition).
	The claim in loc.~cit. is that $\Gamma$ is an arithmetic subgroup of $\SO_{1,n+1}(\RR)^{\circ}$ via a map $A \mapsto \widetilde{A}$ given in section 3 of loc.~cit.
	
	Via what we would call a ``Weil restriction argument'', they map $\SL_2(\CC_n)$ to some $G(\RR)$ and $G$ as a closed $\QQ$-subgroup $G$ of $\GL_{2^n,\QQ}$. 
	This map is done by mapping $\CC_n$ via the left regular representation to $M_{2^{n-1}}(\RR)$, and then taking $2\times 2$ matrices with entries being matrices in $M_{2^{n-1}}(\RR)$. 
	The Clifford group elements $\CC^{\times}_n$ map to elements in $\GL_{2^{n-1},\QQ}(\RR)$, which then give a map from $\SL_2(\CC_n)$ or $\GL_2(\CC_n)$ to $\GL_{2^n}(\RR)$.
	The proof of \cite[Theorem 3]{Maclachlan1989} at the top of page 745 then states that this representation identifies $\GL_2(\CC_n)$ with the $\RR$-points of a $\QQ$-algebraic subgroup, and that the image of $\Gamma=\SL_2(\ZZ[i_1,...,i_{n-1}])$ is $G(\ZZ)$. 
	We emphasize that they are discussing an arithmetic subgroup of $G$. 
	From this, in the proof of \cite[Corollary, page 745]{Maclachlan1989}, they
	deduce ``by Theorems 2 and 3'' that $\Gamma$, by the map $A \mapsto \widetilde{A}$, has finite covolume in $\SO_{1,n+1}(\RR)$. 
	This uses the arithmeticity of the inclusion of $\Gamma$ in $G(\RR)$ to conclude the consequences of Borel--Harish-Chandra (their Theorem 2) for the group $\SO_{1,n+1}(\RR)$.
	This step is valid: the map $A \mapsto \widetilde{A}$ comes from a continuous surjective homomorphism of Lie groups $G(\RR)\cong\SL_2(\CC_n) \to \SO_{1,n+1}(\RR)^{\circ}$ with finite kernel, and such a homomorphism carries a discrete subgroup of finite covolume to a discrete subgroup of finite covolume.
	
	This is not the approach taken in this section of the manuscript, and we give a different proof based on exact sequences of group schemes.
\end{remark}

\begin{remark}\label{rem:EGM-arithmeticity}
	In \cite[Definition 6.1]{Elstrodt1988}, they define $\SL_2(\Ocal)$ for $*$-stable $\Ocal$ as $\SL_2(\CC_n) \cap M_2(\Ocal)$.
	After this they state that $\SL_2(\CC_n)$ acts on $\CC_n\times \CC_n$ via usual matrix multiplication, and that via this action one can see $\SL_2(\Ocal)$ as an arithmetic subgroup since it is the stabilizer of the lattice $\Ocal\times \Ocal \subset \CC_n \times \CC_n$.
\end{remark}

\subsection{Weil Restriction for Clifford Algebras}\label{sec:weil-restriction}

Let $R$ be a commutative ring. 
Let $(W,q)$ be a quadratic module which is locally free of rank $n$. 
Let $A = \Clf(W,q)$ be the associated Clifford algebra. 
We develop a formalism to turn ``Clifford algebra functors'' into schemes (or group schemes or ring schemes). 
An example is the functor which takes an $R$-algebra $R'$ to the set of Clifford algebra elements in $C_{R'}=\Clf(W\otimes_R R', q_{R'})$ given by $ \lbrace x \in C_{R'} \colon x =\overline{x} \rbrace$ for $C$ a Clifford algebra defined over $R$.

The idea is clear to people familiar with Weil restriction techniques \cite{Bosch1990}.

\begin{definition}\label{def:adr}
We define the category of \emph{associative difference rings (with
			one involution and one anti-involution which commute)} $\Ring_R^{\adr}$ 
		to be the category of 
		associative $R$-algebras together with three involutions $x\mapsto
		x'$, $x\mapsto x^*$, $x\mapsto \overline{x}$ where $(x')^* =
		(x^*)' = \overline{x}$. We also require $(xy)'=x'y'$ and $(xy)^*
		=y^*x^*$. Morphisms in the category are associative $R$-algebra
		morphisms which commute with the involutions.
\end{definition}

It is a nontrivial problem to characterize Clifford algebras within difference algebras and find the correct notion of a category of associative algebras, which includes Clifford algebras and their orders but not being too hard or too soft, which also admits the correct functorial properties that allows us, for example, to keep track of involutions and Clifford vectors under isomorphisms given by Bott periodicity.
  
The following is an attempt at a
preliminary definition that is more rigid than the one given above.
\begin{definition}
  Let $R$ be a commutative ring with total quotient algebra $K$.
  We define the category of {\em abstract Clifford algebras}
  $\Ring_R^{\Clf}$ to be the category of associative $R$-algebras $A$
  with commuting involutions $x\mapsto x'$, $x\mapsto x^*$ and
  $x\mapsto\bar x=(x')^*=(x^*)'$, where $'$ is multiplicative and $*$ is
  anti-multiplicative, and an
  $R$-saturated submodule $V$ of $A$ containing $1$, 
  generating $A \otimes_R K$ as a $K$-algebra, and such that $v\overline{v} \in R$
  for all $v \in V$.
  Morphisms $(A,',*,\bar{\phantom{x}},V) \to (B,',*,\bar{\phantom{x}},W)$
  are morphisms
  $\phi: A \to B$ of associative difference rings with $\phi(V) \subseteq W$.
  A {\em Clifford ideal} of an abstract Clifford algebra is a two-sided
  ideal stable under the three involutions.
\end{definition}

One verifies immediately that the kernel of a Clifford algebra morphism
is a Clifford ideal, and that if $A$ is an abstract Clifford algebra and
$I$ a Clifford ideal, then the usual quotient $A/I$ inherits an abstract
Clifford algebra structure over $R/(I \cap R)$. Indeed, this algebra
has the usual universal property of a quotient ring, namely that all
Clifford algebra morphisms $\phi: A \to B$ with $\phi(I) = \{0\}$ factor
through $A \to A/I$. 

The following example shows that not all abstract Clifford algebras are Clifford algebras. 
\begin{example}\label{rem:not-all-abs-are-cliff}	
  Let $A_0$ be the Clifford algebra over
  $\ZZ$ associated to the diagonal quadratic form in $3$ variables with
  diagonal matrix $(1,1,2)$.  
  Let $I$ be the ideal generated by
  $i_1i_2a_3$ where $i_1^2=-1$, $i_2^2=-1$, and $a_3^2=-2$.
  This is a Clifford ideal, so the quotient $A_0/I$ is
  an abstract Clifford algebra. We have $I \cap \ZZ = 2 \ZZ$, so one might
  expect that $A_0/I$ is the Clifford algebra $A'$ over $\ZZ/2\ZZ$
  associated to the reduction mod $2$ of the same form.
  However, that is not true because $i_1i_2a_3 = 0$ in $A_0/I$ but not in $A'$.
\end{example}

\begin{example} Every Clifford-stable order $\Ocal$ (\ref{def:clifford-stable})
  in a Clifford algebra is an abstract Clifford algebra, taking the
  involutions to be those on the Clifford algebra and the distinguished
  submodule to be the set of Clifford vectors in $\Ocal$.
\end{example}

\begin{example}
  For $n \in \NN$ let $\CC_n$ be the Clifford algebra associated to the form
  $q(x) = \sum_{i=1}^{n-1}x_i^2$. The ring homomorphism $\CC_2 \to \CC_6$ taking
  $i_1$ to $i_1i_2i_3i_4i_5$ is a map of associative difference rings in this
  sense, but it is not a map of abstract Clifford algebras because the
  image of the Clifford vector $i_1$ is not a Clifford vector.
\end{example}

Note that kernels of morphisms of difference algebras $A\to A'$ are two-sided ideals $I \subset A$ which are closed under the involutions, and conversely, every morphism of difference algebras factors through the quotient by such an ideal. 

Let $(A,\sigma_1,\sigma_2)$ be an associative difference ring. 
For any associative $A$-algebra $B$ 
there exists a ring $B^1$ which is infinitely generated as a noncommutative ring by the symbols $b'$ and $b^*$ for $b\in B$ modulo the two-sided ideal generated by the relations $a' = \sigma_1(a)$, $a^*=\sigma_2(a)$ for $a\in A$; $(b_1b_2+b_3)' = b_1'b_2'+b_3'$ and $(b_1b_2+b_3)^* = b_2^*b_1^*+b_3^*$ for $b_1,b_2,b_3\in B$; $(b')^* = (b^*)'$ for $b\in B$. 
The ring $B^1$ comes with the natural structure of a difference algebra where we understand that $b''=b$ and $b^{**}=b$ for every $B$. 
This ring is called the \emph{jet or prolongation ring} in differential algebra \cite{Moosa2010}. 

With $A$ still a difference ring we consider the ring $A[z_1,\ldots,z_n]$ of noncommutative polynomials (this is just a free algebra).
The ring $A[z_1,\ldots,z_n]^1$ is called the ring of \emph{associative difference polynomials} (or just \emph{difference polynomials}).  

Let $R$ be a commutative ring. 
Let $R[z_1,\ldots,z_n]$ be the ring of associative difference polynomials.
Let $I \subset R[z_1,\ldots,z_n]^1$ be a difference ideal. 
A \emph{(noncommutative) involution scheme for two involutions with one anticommutative} (or simply a \emph{difference scheme}) is a functor from the category of $R$-difference algebras to sets given by $A\mapsto \Hom(R[z_1,\ldots,z_n]^1/I, A)$ where $\Hom$ is a morphism in the category of $R$-difference algebras.
Informally, the homomorphisms are thought of as solutions to the equations defining the difference ideal $I$.  

\begin{theorem}\label{thm:representability}
Let $F$ be a field. 
	Let $R$ be a ring with fraction field $F$. 
	Let $q$ be a quadratic form over $R$. 
	Let $\Ocal \subset \Clf(F^d, q)$ be an $R$-order which is closed under $'$ and $*$ and is free of finite rank as an $R$-module (for example, any such order when $R$ is $\ZZ$ or a principal ideal domain), and let
	$X$ be a difference scheme over $\Ocal$.
	Then there exists a scheme $\Res_{\Ocal/R}(X)$ over $\Spec(R)$ that represents
	the functor on $R$-algebras given by $R' \to X(\Ocal\otimes_R R')$.
\end{theorem}
\begin{proof}
	Let $I$ be the difference ideal in $\Ocal[z_1,\ldots,z_m]^1$ defining $X$. 
	We need to show that $X(\Ocal\otimes_R R')$ are the solutions of a set of polynomial equations in $R'$ over $R$. 
	Let $\gamma_s$ be a basis for $\Ocal$ as an $R$-module where $s$ runs over a finite index set $A$. 
	Write $\gamma_s \gamma_t = \sum_{r} c^{s,t}_r \gamma_r$, $\gamma_s' = \sum_{r} e_{s,r} \gamma_r$ and $\gamma_s^* = \sum_{r} f_{s,r} \gamma_r$ with $c^{s,t}_r, e_{s,r}, f_{s,r} \in R$; this uses that $\Ocal$ is closed under $'$ and $*$. 
	Let $f(z_1,\ldots,z_m) \in I$.
	We will write $z_j = \sum_s x_{j,s} \gamma_s$ with $\lbrace x_{j,s}\colon 1\leq j \leq m, s\in A \rbrace$, a collection of indeterminates for a commutative polynomial ring. 
	Then we may write
	$$ f( \sum_s x_{1,s} \gamma_s, \ldots, \sum_s x_{m,s} \gamma_s ) = \sum_s P_{f,s} \gamma_s $$
	where $P_{f,s} \in R[ x_{j,s} \colon 1\leq j \leq m, s\in A]$ are (commutative) polynomials with coefficients in the commutative ring $R$ in $m\cdot \vert A \vert$ variables $x_{k,s}$. 
	We then define $\widetilde{I}$ to be the ideal in this polynomial ring generated by $P_{f,s}$ for $f \in I$ and $s\in A$, and then take 
	$$ \Res_{\Ocal/R}(X) = \Spec R[ x_{j,s} \colon 1\leq j \leq m, s\in A]/\widetilde{I}. $$
	Here $z_j'$, $z_j^*$ and $\overline{z_j}$ are expanded using the constants $e_{s,r}$ and $f_{s,r}$. Since the $\gamma_s\otimes 1$ form an $R'$-basis of $\Ocal\otimes_R R'$, an $R$-algebra homomorphism $R[ x_{j,s}]/\widetilde{I}\to R'$ is the same as a point of $X(\Ocal\otimes_R R')$, and this bijection is compatible with $R$-algebra homomorphisms $R'\to R''$.
\end{proof}

\begin{definition}
We call the scheme $\Res_{\Ocal/R}(X)$ in the above theorem the \emph{Weil restriction} of the difference scheme $X$. 
\end{definition}

We remark that if $\Ocal$ is an order in a quadratic field, we recover the usual Weil restriction.
The special case of $\underline{\Ocal}$ is the Weil restriction of the zero ideal in the noncommutative polynomial ring $R[z]$.
We have $\underline{\Ocal}(R) = \Ocal$ and $\underline{\Ocal}(R') = \Ocal\otimes_R R'$ for $R'$ an $R$-algebra.

\subsection{Clifford, Pin, and Spin Group Schemes}\label{sec:our-group-schemes}
Let $q$ be a quadratic form over $\ZZ$. 
By the Weil restriction technique we can show that the functors
$$R \mapsto \Clf(q_R),\Clf(q_R)^\mon,\Clf(q_R)^\times,\GL_2(\Clf(q_R)),\SL_2(\Clf(q_R)),\Spin(q_R),$$
which take a commutative ring $R$ to any of the various groups and monoids, are schemes.
In addition, we can replace $\Clf(q_R)$ in these constructions with $\underline{\Ocal}$ when $\Ocal$ is an order in $\Clf(q_{\QQ})$, which is Clifford-stable (Definition~\ref{def:clifford-stable}), hence closed under the involutions $'$ and $*$.

We will denote these functors by 
$$\Clf_q, \quad \Clf_q^\mon,\quad  \Clf_q^\times, \quad \GL_2(\Clf_q), \quad \SL_2(\Clf_q), \quad \Spin_q.$$

We will have similarly defined functors, with $\underline{\Ocal}$ replacing $\Clf_q$, when $\Ocal$ is an order in $\Clf(q_{\QQ})$, which is Clifford-stable. For $\GL_2$ and $\SL_2$ we apply Definition~\ref{def:gl2} to $C=\Ocal\otimes_{\ZZ}R$, reading ``$\in R$'' as ``$\in \ZZ\cdot 1\otimes R$'' and $\Vec(C)$ as $(\Ocal\cap\Vec(q_{\QQ}))\otimes R$. For $(\Ocal\otimes_\ZZ R)^{\mon}$, $(\Ocal\otimes_\ZZ R)^{\times}$ and $\Spin$ we read $W$ as $(\Ocal\cap W_{\QQ})\otimes R$, $\Clf_+$ as $(\Ocal\cap\Clf(q_{\QQ})_+)\otimes R$ and $\Vec$ as $(\Ocal\cap\Vec(q_{\QQ}))\otimes R$. 

\begin{remark}
	We do not use a group scheme $\PSL_2(\Clf_q)$. When $q$ is positive definite, so that $\Clf(q_\RR)\cong \CC_n$, we write $\PSL_2(\Ocal) = \SL_2(\Ocal)/\lbrace \pm 1 \rbrace$ and $\PSL_2(\CC_n) = \SL_2(\CC_n)/\lbrace \pm 1 \rbrace$. These are abstract groups, and $\PSL_2(\Ocal) \subset \PSL_2(\CC_n)$.
\end{remark}

All of the Clifford algebra constructions described in Definitions~\ref{defn:CliffordGroup},~\ref{defn:spin}
are actually schemes.
\begin{theorem}\label{thm:group-schemes}
	Let $q$ be an $n$-ary quadratic form over $\ZZ$. 
	The functor $\Clf_q$ is an associative-algebra scheme over $\ZZ$.
The functor $\Clf_q^\mon$ is a monoid scheme over $\ZZ$. 
The functors $\Clf_q^\times$, $\GL_2(\Clf_q)$, $\SL_2(\Clf_q)$, and $\Spin_q$ are affine group schemes over $\ZZ$.
	A similar statement holds for $\Clf_q$ replaced by $\underline{\Ocal}$ if $\Ocal \subset \Clf(q_{\QQ})$ is a Clifford-stable order (Definition~\ref{def:clifford-stable}).
\end{theorem}
\begin{proof}
	This follows from Weil restriction.
	For $\GL_2(\Clf_q)$ and $\SL_2(\Clf_q)$ the group axioms hold for every commutative ring $R$ by Corollary~\ref{lem:gl-is-group}, applied to $C=\Clf(q_R)=\Clf(q)\otimes_{\ZZ}R$.
For $\underline{\Ocal}$, the submodules $\ZZ\cdot 1$ and $\Ocal\cap\Vec(q_{\QQ})$ are direct summands of $\Ocal$, because the quotients are torsion-free (for $\ZZ\cdot1$ this uses $\Ocal\cap\QQ=\ZZ$, as elements of an order are integral over $\ZZ$; it also makes $R\to\Ocal\otimes_\ZZ R$ injective, so that $\Delta$ can be read as a unit of $R$). The same saturation fact makes $\Vec(C)$ a direct summand of $C$ compatibly with base change, which is what the Nakayama argument for Clifford units after Definition~\ref{defn:CliffordGroup} needs. 
The proof of Corollary~\ref{lem:gl-is-group} uses only the direct-summand property, that $*$ and $x\mapsto\overline{x}$ are anti-automorphisms of $C$, and that $C$ is a finitely generated $R$-module, so it applies verbatim to $C=\Ocal\otimes_{\ZZ}R$.
\end{proof}

The following Theorem we will not attempt for $\underline{\Ocal}$ as it requires tracking through the arithmetic Bott periodicity theorems with an order replacing a Clifford algebra for two reasons: 1) doing it for $\Clf_q$ is sufficient for our application; 2) tracing through Bott periodicity for a general $\underline{\Ocal}$ would require a tremendous amount of work.

\begin{theorem}[Exceptional Isomorphism]\label{thm:mcinroy}
	Let $q$ be an $(n-1)$-ary quadratic form and let $Q = x^2-yz+q$.
	There is an isomorphism $\SL_2(\Clf_q) \xrightarrow{\sim} \Spin_Q$ as group schemes over $\Spec(\ZZ)$. 
\end{theorem}
\begin{proof}
	For every commutative ring $R$, Theorem~\ref{lem:pre-exceptional} applied to $q_R$ gives a multiplicative bijection $\psi_R\colon\SL_2(\Clf(q_R))\to\Spin(Q_R)$, which is a group isomorphism by Corollary~\ref{lem:gl-is-group}.
	The map $\psi_R$ is given by the formula \eqref{eqn:psi-map}, which is compatible with every ring homomorphism $R\to R'$, so the maps $\psi_R$ form an isomorphism of group functors $\SL_2(\Clf_q)\to\Spin_Q$.
	By the Yoneda lemma it is an isomorphism of the representing group schemes.
\end{proof}

We remind the reader that schemes are automatically sheaves on the \'etale and fppf sites. 
Because we want to do sheaf cohomology on the fppf site of $\Spec(\ZZ)$, we define morphisms of group schemes to be surjective if they are epimorphisms of fppf sheaves.

For a regular or semiregular form $q$, Bass \cite{Bass1974} relates $\Spin_q$ to $\SO_q$ by a long exact sequence involving the discriminant module; see also \cite[\S C.4]{Conrad2014}. 
The forms $Q$ used in \S\ref{sec:arithmeticity} are not of this kind at $2$ and at the primes dividing $d_1\cdots d_{n-1}$, and there $\Spin_Q\to\SO_Q$ is not surjective as a morphism of fppf sheaves (Remark~\ref{rem:spin-sequence-at-2}). 
We therefore work with the quotient $\Spin_q/\mu_2$, which makes sense for every $q$.
For a quadratic form $q$ on a finite free $\ZZ$-module $W$, let $\O_q$ be
the affine $\ZZ$-group scheme of all isometries of $q$: for a commutative
ring $R$, its points are the $R$-linear automorphisms $g$ of $W\otimes R$
with $q_R\circ g=q_R$. No nondegeneracy is assumed.
The second roots of unity $\mu_2\subset\Spin_q$ form a central closed
subgroup scheme. The fppf quotient sheaf $\PSpin_q=\Spin_q/\mu_2$ is
represented by an affine group scheme over $\ZZ$.

\begin{samepage}
\begin{theorem}\label{thm:short-exact-sequence}
Let $q$ be a quadratic form over $\ZZ$ on a free module $W$ of finite rank\footnote{The construction uses a basis of $W$, and the result does not depend on that choice up to canonical isomorphism.}.
	We have the following exact sequence of sheaves of groups on the big fppf site of $\Spec(\ZZ)$:
	\begin{equation}
	1 \to \mu_2 \to \Spin_q \to \PSpin_q \to 1.
	\end{equation}
	The morphism $\pi:\Spin_q\to\O_q$, $u\mapsto \pi_u$, where $\pi_u(x) = uxu^*$ on imaginary Clifford vectors $x$, factors through a morphism of $\ZZ$-group schemes $\overline{\pi}:\PSpin_q\to\O_q$.
\end{theorem}
\end{samepage}
\begin{proof}
	For $r\in\mu_2(R)$ the scalar $r$ is even, $\ns(r)=r^2=1$ and $rWr^*=r^2W=W$, so $\mu_2\subset\Spin_q$. It is central, and it is the closed subscheme of $\Spin_q$ on which the coordinates of all non-scalar monomials vanish.
	The group scheme $\mu_2=\Spec(\ZZ[t]/(t^2-1))$ is finite and free of rank $2$ over $\ZZ$, and it acts freely by multiplication on the affine scheme $\Spin_q$ (Theorem~\ref{thm:group-schemes}).
	So the groupoid $\mu_2\times\Spin_q\rightrightarrows\Spin_q$ has finite locally free source and target and is an equivalence relation, and by \cite[Tag 03BM]{stacks-project} the fppf quotient sheaf $\Spin_q/\mu_2$ is represented by an affine scheme $\PSpin_q$, with $\Spin_q\times_{\PSpin_q}\Spin_q=\mu_2\times\Spin_q$.
	Since $\mu_2$ is normal, $\PSpin_q$ is a group scheme, and the sequence is exact: $\Spin_q\to\PSpin_q$ is an epimorphism of fppf sheaves by construction, with kernel $\mu_2$. 
	Since $\pi$ is trivial on $\mu_2$, it factors through $\PSpin_q$.
\end{proof}

\begin{corollary}\label{cor:pspin-gspin-quotient}
	The inclusion $\Spin_q\subset\GSpin_q$ induces an isomorphism
	$\PSpin_q\cong\GSpin_q/\GG_m$, where the quotient is an fppf sheaf.
	For every commutative ring $R$ with $\Pic(R)=0$, we have
	$\PSpin_q(R)=\GSpin_q(R)/R^{\times}$.
	In particular, $\PSpin_q(\ZZ)=\GSpin_q(\ZZ)/\lbrace\pm1\rbrace$.
\end{corollary}
\begin{proof}
	The induced morphism $\Spin_q/\mu_2\to\GSpin_q/\GG_m$ is injective, because a scalar $r\in R^{\times}$ lies in $\Spin_q(R)$ exactly when $r^2=1$. 
	It is surjective: if $g\in\GSpin_q(R)$ and $c=\ns(g)$, then $t^{-1}g\in\Spin_q(R[t]/(t^2-c))$, and $R[t]/(t^2-c)$ is an fppf cover of $R$.
	Finally, the fibre of $\GSpin_q\to\GSpin_q/\GG_m$ over a point $x\in(\GSpin_q/\GG_m)(R)$ is a $\GG_m$-torsor over $\Spec(R)$, and these are classified by $H^1(\Spec(R)_{\fppf},\GG_m)=\Pic(R)$ \cite[Tags 03AJ, 03P8]{stacks-project}. If $\Pic(R)=0$ this torsor has a section, so $x$ lifts to $\GSpin_q(R)$, uniquely up to $R^{\times}$.
\end{proof}

\begin{remark}\label{rem:spin-sequence-at-2}
	Let $q=d_1x_1^2+\cdots+d_{n-1}x_{n-1}^2$, $Q=x^2-yz+q$ and $N=2d_1\cdots d_{n-1}$, and let $\SO_Q\subset\SL_{n+2}$ be the closed $\ZZ$-subgroup scheme of matrices preserving $Q$.
	Over $\ZZ[1/N]$ the form $Q$ is non-degenerate in the sense of \cite[\S C.1]{Conrad2014}, and $\overline{\pi}$ is an isomorphism $\PSpin_Q\to\SO_Q$ over $\ZZ[1/N]$ by \cite[Lemma C.4.1 and Proposition C.4.10]{Conrad2014}.
	For $n\geq 2$ this fails over $\ZZ$.
	If $d_i$ is odd, then $xx_i$ is a nontrivial element of $\Spin_Q(\FF_2)=\PSpin_Q(\FF_2)$ which acts trivially on vectors.
	If $q=x_1^2$, then $g=s_xs_{x+x_1}$ lies in $\SO_Q(\ZZ)$ but not in $\overline{\pi}(\PSpin_Q(\ZZ))$: its lifts to $\GSpin_Q(\QQ)$ are $\lambda x(x+x_1)$ with $\lambda\in\QQ^{\times}$, of spinor norm $2\lambda^2\neq\pm1$. If $q=3x_1^2$, the same holds for $g=s_xs_{x_1}$, whose lifts $\lambda xx_1$ have spinor norm $3\lambda^2$.
	In particular $1\to\mu_2\to\Spin_Q\to\SO_Q\to 1$ is not exact over $\ZZ$, nor over $\ZZ[1/2]$.
\end{remark}

\section{The Clifford Uniformization of Hyperbolic Space}

This section gives an overview of the theory of M\"obius transformations in the Clifford setting. These transformations of the form $x\mapsto (ax+b)(cx+d)^{-1}$ allow us to generalize many of the formulas of complex analysis to the Clifford algebra context.

For example, the Cayley transformation to the unit ball in this setting from $\Hcal^{n+1}$ to $B^{n+1}$ is given by $C(x) = ((1-i_n)x-(1+i_n))((1-i_n)x+(1+i_n))^{-1}$, the M\"obius transformation of $\begin{psmallmatrix} 1-i_n & -(1+i_n) \\ 1-i_n & 1+i_n \end{psmallmatrix}\in\GL_2(\CC_{n+1})$, which has $\Delta=4$. For $n=1$ this is the classical map $(x-i)(x+i)^{-1}$.
It has the property that $C(i_n)=0$, $C(0) = -1$, $C(\infty) = 1$, $C(t i_n) = (t-1)(t+1)^{-1}$ for $t \in \RR$.
If $x$ is real we have  $C(x)= (x-i_n)^2/(x^2+1) = ( (x^2 -1)- 2xi_n )/(x^2+1)$
whose components $( (x^2-1)/(x^2+1), -2x/(x^2+1))$ constitute
the famous rational parametrization of the unit circle. 
Also, in this setting one can see for example that the conformal bijections from the
ball to itself are exactly the maps $ g(x) = (x-u)(1-\overline{u}x)^{-1}$ composed with orthogonal transformations of $V_{n+1}$,
for $u \in V_{n+1}$ with $\vert u \vert < 1$.

\subsection{The Positive and Special Orthogonal Groups}
In preparation for our work with hyperbolic space, which will be defined in
\S\ref{subsec:clifford-unif}, we will discuss the orthogonal group
of a real quadratic form of signature $(1,m)$.
  We will let $\O(1,m)=\O_{1,m}(\RR)$ denote the real points of the $\ZZ$-group scheme $\O_{1,m}$. 
The real Lie group $\O(1,m)$ has four connected components. The $m$-dimensional
hyperboloid defined by $x^2 - y_1^2 -\dots - y_m^2 = 1$ has two connected components,
one where $x \geq 0$ and one where $x\leq 0$ (there are no points with $|x| < 1$, so
these are genuinely distinct components).

\begin{definition}\label{def:so-po-pso}
  Following \cite[p. 58]{Ratcliffe2019}, we define $\SO(1,m)=\SO_{1,m}(\RR)$ to be
  the subgroup of $\O(1,m)$ where the determinant is positive; 
  $\PO(1,m)$ to be the subgroup acting trivially on the set of components;
  and $\PSO(1,m)=\O_{1,m}(\RR)^{\circ} = \PO(1,m) \cap \SO(1,m)$.	
\end{definition}

\begin{remark} The group $\SO$ contains $\operatorname{diag}(-1,-1,1,\ldots,1)$,
  while $\PO$ does not. Thus $\SO\not\subset\PO$, and there is no induced
  action of $\SO(1,m)$ on $\Hcal^{m}$.
  The corresponding isometry groups are
  \[
    \PO(1,m)=\Isom(\Hcal^m),\qquad
    \PSO(1,m)\cong\Isom(\Hcal^m)^{\circ}.
  \]
\end{remark}

\subsection{Symmetric Spaces}

For those unfamiliar with Riemannian manifolds, a Riemannian metric on a manifold $M$
can informally be described as a smoothly varying family of inner products
on the tangent spaces $TM_x$ at points of $M$: see \cite{Ratcliffe2019} for formal definitions.

Recall that a Riemannian manifold $X$ is \emph{homogeneous} if its group of isometries $\Isom(X)$ acts transitively: for all $x,y \in X$ there exists some $\phi \in \Isom(X)$ such that $\phi(x)=y$ \cite[1.1.1]{Morris2015}.
A homogeneous space $X$ is a \emph{symmetric space} if, and only if, it is connected and there exists some nontrivial $\phi \in \Isom(X)$ such that $\phi^2 =\id_X$ and $\phi$ has an isolated fixed point. 
Note that by homogeneity this implies that every $x\in X$ admits some involutive isometry $\phi$ having $x$ as an isolated fixed point \cite[1.1.5]{Morris2015}.

There is a particular description in terms of connected Lie groups that might be more familiar to readers. 
If $X$ is a connected homogeneous space then $G=\Isom(X)^{\circ}$ is a Lie group which acts transitively, and $X \cong G/K$ as Riemannian manifolds where $K$ is the stabilizer of some point.
Here $G$ is given its $G$-invariant Riemannian metric, and $K$ being compact assures us that the metric descends to $G/K$. (This doesn't hold if we just assume that $K$ is a closed subgroup \cite[\S 1.2]{Morris2015}.)
This isomorphism works in the converse direction: if $G$ is a connected semisimple Lie group with finite center and $K$ is a maximal compact subgroup, then $G/K$, with a $G$-invariant metric, is a symmetric space. 

A complete Riemannian manifold $M$ is a \emph{locally symmetric space} if and only if its universal cover $X$ is a symmetric space. 
This means there is a group of isometries $\Gamma \subset \Isom(X)$ so that $\Gamma$ acts freely and properly discontinuously on $X$ and $M\cong \Gamma \setminus X$ as Riemannian manifolds.

\subsection{M\"obius Transformations}\label{sec:mobius}

This subsection gives an account of M\"obius transformations in our setting following Ahlfors \cite{Ahlfors1984} \cite[\S 2.2]{Ahlfors1985}.
M\"obius transformations on $\Hcal^n$ for $n\geq 2$ are well-studied
and not new.
Their presentation as fractional linear transformations using Clifford groups and Clifford vectors is also not new, but not well-studied. 
For this reason we include a summary of this theory. 
It will be helpful to recall that if $j = i_1i_2\cdots i_m$ then $j^2 = (-1)^{{m+1\choose 2}}$, $i_aj = (-1)^{m-1}ji_a$ for $1\leq a\leq m$, and $j^*=(-1)^{ {m\choose 2} }j$. Also, the center $Z(\CC_n)$ of $\CC_n$ is given by  
\begin{equation}\label{eqn:center}
Z(\CC_n) = \begin{cases}
\RR, & \mbox{ $n$ odd } \\
\RR[J], & \mbox{ $n$ even} 
\end{cases}, \quad J^2 = (-1)^{ {n \choose 2} }, \quad J=i_1i_2\cdots i_{n-1}
\end{equation}

Just as in complex analysis, we consider M\"obius tranformations which act on
the extended plane. 

\begin{definition}\label{def:one-point-comp}
	We let $S^n = V_n \cup \{\infty\}$ and give it the topology of the one-point compactification.
\end{definition}
Our notation is
justified by the well-known fact that $S^n$ is homeomorphic to the $n$-sphere.

\begin{definition}\label{def:moebius}
	A \emph{M\"obius transformation} is a homeomorphism $g:S^n\to S^n$ of the form $g(x) = (ax+b)(cx+d)^{-1}$ for $\begin{psmallmatrix} a  & b\\  c & d\end{psmallmatrix} \in M_2(\CC_n)$ such that $x\mapsto (ax+b)(cx+d)^{-1}$ induces a homeomorphism $S^m \to S^m$ for all $m\geq n$. 
	For $m\geq n$ the action goes through the inclusion $\CC_n \subset \CC_m$.
	The group of M\"obius transformations will be denoted $\GM(n)$.
\end{definition}

Following Ahlfors, we will often conflate M\"obius transformations $g(x)=(ax+b)(cx+d)^{-1}$ and the matrix $\begin{psmallmatrix}a&b\\c&d\end{psmallmatrix}$ that induces it. 
This section justifies this procedure.
The aim of this section is to show that if $g = \begin{psmallmatrix}a&b\\c&d\end{psmallmatrix}\in M_2(\CC_n)$ satisfies the condition of Definition~\ref{def:moebius}, then $g \in \GL_2(\CC_n)$.

\begin{theorem}\label{thm:gm-to-gl2}
	Let $\begin{psmallmatrix} a&b\\c&d\end{psmallmatrix}\in M_2(\CC_n)$ be a matrix such that $x\mapsto(ax+b)(cx+d)^{-1}$ induces a homeomorphism of $S^m$ for every $m\geq n$. Then $\begin{psmallmatrix} a&b\\c&d\end{psmallmatrix}\in\GL_2(\CC_n)$; in particular each of $a,b,c,d$ is $0$ or lies in the Clifford group $\CC_n^\times$.
\end{theorem}

We will prove this theorem by means of a sequence of lemmas.

\begin{lemma}\label{lem:identity-scalar}
	If for all $x\in S^n = V_n \cup \lbrace \infty \rbrace$ we have $(ax+b)(cx+d)^{-1}=x$, then $\begin{psmallmatrix} a & b \\ c & d \end{psmallmatrix}= \begin{psmallmatrix} z & 0 \\ 0 & z \end{psmallmatrix}$
	for some $z\in Z(\CC_n)$.
\end{lemma}

\begin{proof} Applying $\begin{psmallmatrix} a & b \\ c & d \end{psmallmatrix}$ to $0, \infty, 1$ (recall that $1 \in V_n$)
	we obtain successively $b = 0, c = 0, a = d$.
	The first two conditions tell us that $\begin{psmallmatrix} a & b \\ c & d \end{psmallmatrix}$ is diagonal. For all
	Clifford vectors $v$ we must have $ava^{-1} = v$, so $av = va$. The Clifford
	vectors generate $\CC_n$, so this implies $a \in Z(\CC_n)$.
\end{proof}

In the rest of this subsection $g=\begin{psmallmatrix}a&b\\c&d\end{psmallmatrix}\in M_2(\CC_n)$ is a matrix as in Definition~\ref{def:moebius}.
For $m\geq n$ we write $D_m=\lbrace x\in V_m\colon cx+d\in U(\CC_m)\rbrace$ and $g_m\colon S^m\to S^m$ for the induced homeomorphism.
We read Definition~\ref{def:moebius} in the weakest sense: $D_n\neq\emptyset$, and $g_m(x)=(ax+b)(cx+d)^{-1}$ for every $x\in D_m$; in particular $(ax+b)(cx+d)^{-1}\in V_m$ for $x\in D_m$.
An element of $\CC_n$ is invertible in $\CC_m$ if and only if it is invertible in $\CC_n$, because $\CC_m$ is a free left $\CC_n$-module, so $D_n=D_m\cap V_n$.
The set $D_m$ is the complement in $V_m$ of the zero set of the polynomial $x\mapsto\det(L_{cx+d})$, where $L_y$ is left multiplication by $y$ on $\CC_m$. This polynomial does not vanish on $D_n$, so $D_m$ is open and dense in $V_m$, and $g_m$ restricts to $g_n$ on $S^n$.
We put $g_1=\begin{psmallmatrix}d^*&-b^*\\-c^*&a^*\end{psmallmatrix}$ and recall that $x^*=x$ for $x\in V_m$.

\begin{lemma}\label{lem:inverses}
	There is $\Delta\in\RR^\times$ with $g_1g=gg_1=\Delta I_2$. Equivalently,
	$c^*a=a^*c$, $b^*d=d^*b$, $ab^*=ba^*$, $cd^*=dc^*$ and $d^*a-b^*c=ad^*-bc^*=\Delta$.
	In particular $\Delta(g)=\Delta(g_1)=\Delta$, and $\Delta^{-1}g_1$ is the inverse of $g$ in $M_2(\CC_n)$.
\end{lemma}

\begin{proof}
	Let $m\geq n$ and $x\in D_m$, and put $y=g_m(x)\in V_m$ and $w=cx+d$.
	Applying $*$ to $yw=ax+b$ gives $(xc^*+d^*)y=xa^*+b^*$. Substituting $y=(ax+b)w^{-1}$ gives
	\begin{equation}\label{eqn:star-identity}
	(xc^*+d^*)(ax+b)=(xa^*+b^*)(cx+d).
	\end{equation}
	Both sides are polynomial in $x$, so \eqref{eqn:star-identity} holds for all $x\in V_m$, because $D_m$ is dense.
	Put $P=d^*a-b^*c$, so that $c^*b-a^*d=-P^*$. Replacing $x$ by $tx$ with $t\in\RR$ and comparing powers of $t$ gives
	$x(c^*a-a^*c)x=0$, $d^*b=b^*d$ and $Px=xP^*$ for all $x\in V_m$.
	Taking $x=1$ gives $c^*a=a^*c$ and $P=P^*$, so $P$ commutes with every element of $V_m$.
	For $m=n+1$ this puts $P$ in $Z(\CC_{n+1})\cap\CC_n$, which is $\RR$ by \eqref{eqn:center}, since $i_1\cdots i_n\notin\CC_n$.
	Expanding $g_1g$ now gives $g_1g=PI_2$.

	Suppose $P=0$. For $x\in D_n$ we get $0=g_1g\begin{psmallmatrix}x\\1\end{psmallmatrix}=g_1\begin{psmallmatrix}y\\1\end{psmallmatrix}w$, so $d^*y=b^*$ and $c^*y=a^*$ for all $y$ in the nonempty open set $g_n(D_n)\subset V_n$.
	An affine map that vanishes on a nonempty open set is zero, so $a=b=c=d=0$, which contradicts $D_n\neq\emptyset$.
	Hence $\Delta:=P\in\RR^\times$, and $\Delta^{-1}g_1$ is a left inverse of $g$ in the finite-dimensional algebra $M_2(\CC_n)$, hence a two-sided inverse.
	So $gg_1=\Delta I_2$ as well, and expanding $gg_1$ gives the remaining relations.
\end{proof}

\begin{lemma}\label{lem:poles}
	If $x\in V_n$ and $cx+d$ is not invertible, then $g_n(x)=\infty$. In particular, if $g_n(\infty)=\infty$, then $D_n=V_n$.
\end{lemma}

\begin{proof}
	For $x\in D_n$ with $y=g_n(x)$, the second row of $g_1g\begin{psmallmatrix}x\\1\end{psmallmatrix}=\Delta\begin{psmallmatrix}x\\1\end{psmallmatrix}$ reads $(a^*-c^*y)(cx+d)=\Delta$.
	Let $x_0\in V_n\setminus D_n$ and suppose that $y_0=g_n(x_0)\in V_n$.
	Choose $x_k\in D_n$ with $x_k\to x_0$. Then $g_n(x_k)\to y_0$, and in the limit $(a^*-c^*y_0)(cx_0+d)=\Delta\neq0$.
	So $cx_0+d$ has a left inverse in $\CC_n$ and is therefore invertible, a contradiction.
	The second statement follows because $g_n$ is injective.
\end{proof}

\begin{lemma}\label{lem:affine-moebius}
	If $c=0$, then $a,d\in\CC_n^\times$, $b\in\CC_n^\times\cup\lbrace0\rbrace$ and $g\in\GL_2(\CC_n)$.
\end{lemma}

\begin{proof}
	Since $D_n\neq\emptyset$, $d$ is invertible, so $D_m=V_m$ for all $m\geq n$, and $v:=g_n(0)=bd^{-1}\in V_n$.
	By Lemma~\ref{lem:inverses}, $d^*a=\Delta$, so $a=\Delta f$ with $f=(d^*)^{-1}$, and $d^{-1}=f^*$.
	Hence $g_m(x)-v=\Delta\, fxf^*$ for $x\in V_m$, and so $fV_mf^*\subset V_m$ for all $m\geq n$.
	For $m=n$ this is the second condition in \eqref{eqn:clifford-monoid}.
	For $m=n+1$ and $x=i_n$, which anticommutes with $i_1,\ldots,i_{n-1}$, we get $fi_nf^*=i_nf'f^*\in V_{n+1}$.
	Since $\CC_{n+1}=\CC_n\oplus i_n\CC_n$ and $V_{n+1}\cap i_n\CC_n=\RR i_n$, this forces $f'f^*\in\RR$. Applying $'$ gives $f\overline{f}\in\RR$.
	So $f\in\CC_n^\mon$, and since $f$ is invertible, $f\in\CC_n^\times$ (Definition~\ref{defn:CliffordGroup} and Lemma~\ref{lem:generated-by-clifford-vectors}).
	Hence $a=\Delta f$ and $d=(f^{-1})^*$ lie in $\CC_n^\times$, and $b=vd$ is $0$ or lies in $\CC_n^\times$.
	Finally $ad^*-bc^*=\Delta\in\RR^\times$, $ab^*=\Delta fd^*v=\Delta v\in V_n$ and $dc^*=0$.
	So $g\in\GL_2(\CC_n)$ by Theorem~\ref{thm:gl2-short} and Corollary~\ref{lem:gl-is-group}.
\end{proof}

\begin{proof}[Proof of Theorem~\ref{thm:gm-to-gl2}]
	For $v\in V_n$ let $T_v=\begin{psmallmatrix}1&v\\0&1\end{psmallmatrix}$ and $\sigma=\begin{psmallmatrix}0&1\\1&0\end{psmallmatrix}$. These lie in $\GL_2(\CC_n)$ by Theorem~\ref{thm:gl2-short}.

	\emph{Case 1: $g_n(\infty)=\infty$.}
	By Lemma~\ref{lem:poles}, $D_n=V_n$. In particular $d$ is invertible and $v=g_n(0)=bd^{-1}\in V_n$.
	The matrix $h=T_{-v}g=\begin{psmallmatrix}a-vc&0\\c&d\end{psmallmatrix}$ satisfies Definition~\ref{def:moebius}: it has the same sets $D_m$ and induces $x\mapsto g_m(x)-v$. It fixes $0$ and $\infty$.
	By Lemma~\ref{lem:inverses} applied to $h$, $a-vc=\Delta_h(d^*)^{-1}$ with $\Delta_h\in\RR^\times$.
	For $x\in V_n\setminus\lbrace0\rbrace$ the vector $h(x)$ is nonzero, because $h$ is injective, and
	$$h(x)^{-1}=(cx+d)x^{-1}(a-vc)^{-1}=\Delta_h^{-1}\left(cd^*+dx^{-1}d^*\right)$$
	(compare \cite[(7.3)]{Ahlfors1984}).
	As $x\to\infty$ we have $h(x)\to\infty$, so $h(x)^{-1}\to0$, while the right-hand side tends to $\Delta_h^{-1}cd^*$.
	Hence $cd^*=0$, so $c=0$, and Lemma~\ref{lem:affine-moebius} gives $g\in\GL_2(\CC_n)$.

	\emph{Case 2: $p:=g_n^{-1}(\infty)\in V_n$.}
	Let $h=gT_p\sigma=\begin{psmallmatrix}ap+b&a\\cp+d&c\end{psmallmatrix}$.
	For $x\in V_m\setminus\lbrace0\rbrace$ we have $(cp+d)x+c=\bigl(c(p+x^{-1})+d\bigr)x$ and $(ap+b)x+a=\bigl(a(p+x^{-1})+b\bigr)x$.
	So wherever $(cp+d)x+c$ is invertible, the formula of $h$ agrees with the homeomorphism $x\mapsto g_m(p+x^{-1})$ of $S^m$ (where $0\mapsto g_m(\infty)$ and $\infty\mapsto\infty$); at $x=0$ this follows by continuity.
	The domain of $h$ in $V_n$ contains $(y-p)^{-1}$ for every $y\in D_n$ (note that $p\notin D_n$), so it is nonempty.
	Thus $h$ satisfies Definition~\ref{def:moebius}, and it fixes $\infty$. By Case 1, $h\in\GL_2(\CC_n)$, so $g=h\sigma T_{-p}\in\GL_2(\CC_n)$ by Corollary~\ref{lem:gl-is-group}.

	In both cases the entries of $g$ lie in $\CC_n^\mon=\CC_n^\times\cup\lbrace0\rbrace$ by \eqref{eqn:simple-defn}.
\end{proof}

\begin{proposition}\label{prop:vn-equalities}
	If $g=\begin{psmallmatrix}a&b\\c&d\end{psmallmatrix}\in M_2(\CC_n)$ induces a homeomorphism of $S^m$ for every $m\geq n$ via $x\mapsto(ax+b)(cx+d)^{-1}$, then $\Delta(g)=\Delta(g_1)\in\RR^\times$, and
	\begin{equation}\label{E:conditions}
	c^*a, \ b^*d, \ ab^*, \ cd^*\in V_n.
	\end{equation}
\end{proposition}

\begin{proof}
	The first statement is Lemma~\ref{lem:inverses}.
	Since $g\in\GL_2(\CC_n)$, \eqref{eqn:simple-defn} gives $ab^*\in V_n$ and $cd^*=dc^*\in V_n$.
	Applying \eqref{eqn:simple-defn} to $g^{-1}=\Delta^{-1}g_1\in\GL_2(\CC_n)$ gives $d^*b,a^*c\in V_n$, and $V_n$ is closed under $*$.
\end{proof}

\begin{remark}\label{rem:m-clause}
	The proof uses the condition of Definition~\ref{def:moebius} only for $m=n$ and $m=n+1$.
	The case $m=n+1$ cannot be dropped. For $n=4$ the element $u=1+\frac12 i_1i_2i_3$ is central and invertible, so $\begin{psmallmatrix}u&0\\0&u\end{psmallmatrix}$ induces the identity of $S^4$.
	But $u\overline{u}=\frac54+i_1i_2i_3\notin\RR$, and this matrix sends $i_4\in V_5$ to $\frac53 i_4+\frac43 i_1i_2i_3i_4\notin V_5$.
	For $n=2$, $\begin{psmallmatrix}i_1&0\\0&1\end{psmallmatrix}$ preserves $V_2$ but has $\Delta=i_1\notin\RR$.
	For this reason \cite[Theorem 6]{Waterman1993}, whose definition asks only for a bijection of one sphere, does not apply as stated; compare \cite[p.~246]{Lounesto2001}.
\end{remark}

By scaling we may assume that $\Delta(g)=\pm1$.

We now state a result so useful that we refer to it, and some of
its consequences, as the {\em Magic Formula}.
\begin{theorem}[Ahlfors's Magic Formula; see {\cite[(7.2)]{Ahlfors1984}} and {\cite[Theorem 7(i)]{Waterman1993}}]\label{thm:magic}
	If $g = \binom{ a  \ \ b}{  c \ \ d} \in \GL_2(\CC_n)$, $m\geq n$, and $x, y\in V_m$ are such that $cx+d$ and $cy+d$ are invertible, then
	\begin{equation}\label{eqn:magic}
	g(x) - g(y)^* = \Delta(g) (yc^* +d^*)^{-1}(x-y)(cx+d)^{-1}.
	\end{equation}
\end{theorem}

This theorem implies the following useful corollaries. 
\begin{corollary}\label{cor:magic-corollaries}
	\begin{enumerate}
		\item $g(y)=g(y)^*$.
		\item\label{item:derivative-formula} The formula for $g'(x):V_n \to V_n$ if $g(x)\neq \infty$ is given by $g'(x) \cdot u = \Delta(g) (xc^*+d^*)^{-1}u(cx+d)^{-1}$ for all $u \in V_n$.
Note that, up to the factor $\Delta(g)$, this is the map $v\mapsto ava^*$ with $a=(xc^*+d^*)^{-1}$, which is $\vert a\vert^2$ times the rotation $\pi_a$. \item \label{I:conformal} $\vert g'(x)\vert = \vert \Delta(g) \vert \vert cx+d\vert^{-2}$.
		\item \label{I:magic-heights} $g(x)_n/x_n = \frac{\Delta(g)}{\vert \Delta(g) \vert } \vert g'(x) \vert $.
	\end{enumerate}
\end{corollary}

\begin{proof}
	For the first one, specialize the Magic Formula in Theorem~\ref{thm:magic} to $x=y$.
	To prove the second, consider $g(x+tu)-g(x)$ for $t\in \RR$ and let $t\to 0$. 
	The third follows by taking absolute values of the previous. 
We obtain the fourth by applying the Magic Formula to $x$ and to its reflection $\hat{x}=x-2x_ni_n$ in $V_n$: since the entries of $g$ lie in $\CC_n$, we have $g(x)-g(\hat{x})=2g(x)_ni_n$, and this equals $\Delta(g)(\hat{x}c^*+d^*)^{-1}(2x_ni_n)(cx+d)^{-1}$; comparing $n$th components and using item~\eqref{I:conformal} gives the claim.
\end{proof}

This Corollary has geometric consequences. 
For example, the fact that $g \in \GL_2(\CC_n)$ with $\Delta(g)>0$ preserves $\Hcal^{n+1}$ follows from
the behavior of the $n$th coordinate (Siegel height):
$$ g(x)_n = \vert g'(x) \vert x_n = \Delta(g)\, x_n/\vert cx+d\vert^2.$$
If $\Delta(g)<0$, then $g$ maps $\Hcal^{n+1}$ to the lower half-space; for example $\begin{psmallmatrix}1&0\\0&-1\end{psmallmatrix}$ acts by $x\mapsto -x$.
Note that $g(x)=\infty$ is not possible, because it would imply that 
$x=-c^{-1}d\in V_n$, which is not an element of $\Hcal^{n+1}$.
Also, the statement (\ref{I:conformal}) about the magnitude of $g'(x)$ says
that the mapping is conformal. 
In the next section we will go on to prove that
elements $g\in\GL_2(\CC_n)$ with $\Delta(g)>0$ induce isometries of $\Hcal^{n+1}$.

\subsection{The Clifford Uniformization of Hyperbolic Space}\label{subsec:clifford-unif}

The result of this section is that M\"{o}bius transformations in the Clifford sense act as orientation preserving isometries, and conversely that every orientation preserving isometry comes from such a transformation. 
This has been done several times in the literature, and we collect it here in a form that is convenient for future reference (see \cite[Section 3]{Maclachlan1989}, \cite[Theorem 2.3]{Elstrodt1987}, \cite[\S7.5, Theorem]{Ahlfors1984}, \cite[Theorem 5(ii)]{Waterman1993}). 

We remind the reader that in this section the notation $\CC_n$ is used only for $n\geq 1$.

\begin{definition}
	The \emph{Clifford uniformization of Hyperbolic Space} $\Hcal^{n+1}$ is the set $\lbrace x \in V_{n+1} \colon x_n >0 \rbrace$
together with its structure of a Riemannian manifold with metric $ds^2 = (dx_0^2 + \cdots +dx_n^2)/x_n^2$, which we denote by $ds=\vert dx\vert/x_n$.
	Here the coordinates $(x_0,\dots,x_n)$ of $V_{n+1}$ identify $\Hcal^{n+1}$ with $\RR^n\times\RR_{>0}$, the height being $x_n$ (the standard upper half-space model). \end{definition}
For most of the manuscript we identify $\Hcal^{n+1}$ with this particular Riemannian manifold. 
The volume form for this manifold is $dx_0\cdots dx_n/x_n^{n+1}$. 

\begin{remark}
  Hyperbolic $m$-space can also be described by means of a metric on one sheet of a real quadric or
on a ball. These constructions are reviewed in \cite[Chapters 3-4]{Ratcliffe2019}.
\end{remark}

Figure~\ref{fig:mobius-basic-maps} shows the basic maps of part~\eqref{item:factorization-in-gm} acting on $V_2$.
\begin{figure}[htbp!]
\centering
\begin{tikzpicture}[scale=0.62,lab]
\input{tikz/new-mobius-data}
 \foreach \nm/\ttl/\sx [count=\k from 0] in {id/{$x$}/0, trans/{$x+\mu$}/4.4, inv/{$-x^{-1}$}/8.8, dil/{$\lambda^2x$}/13.2, rot/{$uxu^*$}/17.6}{
  \begin{scope}[xshift=\sx cm]
   \draw[gray!50,->] (-2,0)--(2,0); \draw[gray!50,->] (0,-2)--(0,2);
   \draw[gray!60,dotted] (0,0) circle (1);
   \begin{scope}\clip (-2,-2) rectangle (2,2);
    \ifnum\k=0 \fill[domC!40,draw=domC!80!black] \flagid; \fi
    \ifnum\k=1 \fill[domC!15] \flagid; \fill[domC!40,draw=domC!80!black] \flagtrans; \fi
    \ifnum\k=2 \fill[domC!15] \flagid; \fill[domC!40,draw=domC!80!black] \flaginv; \fi
    \ifnum\k=3 \fill[domC!15] \flagid; \fill[domC!40,draw=domC!80!black] \flagdil; \fi
    \ifnum\k=4 \fill[domC!15] \flagid; \fill[domC!40,draw=domC!80!black] \flagrot; \fi
   \end{scope}
   \node at (0,-2.5){\ttl};
  \end{scope}}
\end{tikzpicture}
\caption{The basic M\"obius transformations of Theorem~\ref{thm:isometries-and-mobius}\eqref{item:factorization-in-gm}, acting on $V_2\cong\CC$. Light orange: a test region; dark orange: its image. From left to right: the identity, translation by $\mu=-1.2-0.4i_1$, the inversion $x\mapsto -x^{-1}$ (dotted: the unit circle), dilatation by $\lambda^2=1.5$, and $x\mapsto uxu^*$ with $u=e^{i_1\pi/6}$, a rotation by $\pi/3$. Computed exactly and plotted.}\label{fig:mobius-basic-maps}
\end{figure}
\begin{theorem}\label{thm:isometries-and-mobius}
	\begin{enumerate}
		\item If $g = \begin{psmallmatrix} a & b \\ c & d \end{psmallmatrix}\in\GL_2(\CC_n)$ and $\Delta(g)>0$, then $x\mapsto (ax+b)(cx+d)^{-1}$ restricts to an isometry of $\Hcal^{n+1}$. 
		\item\label{item:factorization-in-gm} Every element of $\GM(n)$ is a composition of translations, inversions, dilatations, special orthogonal transformations, and the reflection $x\mapsto -x$. 
		\item\label{item:mobius-equals-isom} If $n$ is even, every element of $\GM(n)$ is induced by some $g\in\GL_2(\CC_n)$ with $\Delta(g)>0$. If $n$ is odd, the sign of $\Delta(g)$ depends only on the element of $\GM(n)$ induced by $g$.
		\item If $\begin{psmallmatrix} a & b \\ c & d \end{psmallmatrix} \in \GL_2(\CC_n)$, then $g(x) = (ax+b)(cx+d)^{-1}$ is an element of $\GM(n)$. 
		In particular, the map $\GL_2(\CC_n) \to \GM(n)$ is well-defined.
\item  We have the following isomorphisms of groups: $\PSL_2(\CC_n) \cong \SO_{1,n+1}(\RR)^{\circ} \cong \Isom(\Hcal^{n+1})^{\circ}=\Isom^+(\Hcal^{n+1})$, the group of orientation-preserving isometries. \end{enumerate}
\end{theorem}
\begin{proof}
	\begin{enumerate}
		\item Let $g = \begin{psmallmatrix} a & b \\ c & d \end{psmallmatrix} \in \GL_2(\CC_n)$ with $\Delta(g)>0$. Since $\GL_2(\CC_n)$ acts on $V_{n+1}\cup\lbrace\infty\rbrace$ \cite[\S7.5, Theorem, part II]{Ahlfors1984}, Corollary~\ref{cor:magic-corollaries} applies to $x\in\Hcal^{n+1}\subset V_{n+1}$, and the height formula above shows that $g$ maps $\Hcal^{n+1}$ to itself.
		The metric for $\Hcal^{n+1}$ in Clifford form is $\vert dx \vert /x_n$. 
		We have 
		$
		\frac{\vert dg(x)\vert}{g(x)_n} = \frac{\vert g'(x)dx \vert }{g(x)_n} = \frac{\vert g'(x) \vert \vert dx \vert}{x_n\vert g'(x) \vert } = \frac{\vert dx \vert }{x_n}.$
This proves that the differential of $g$ preserves the metric. Moreover $g$ is smooth on $\Hcal^{n+1}$, and it is a bijection of $\Hcal^{n+1}$ onto itself, with inverse induced by $\Delta(g)^{-1}\begin{psmallmatrix} d^* & -b^*\\ -c^* & a^*\end{psmallmatrix}$, which again has positive pseudo-determinant; so $g$ acts on $\Hcal^{n+1}$ by isometries. \item This is \cite[Lemma 10]{Waterman1993}; it is very similar to the decomposition of classical M\"{o}bius transformations, and there is a discussion of this following the proof.
		Waterman's factorization is an identity of formulas on $D_n=\lbrace x\in V_n\colon cx+d\in\CC_n^\times\rbrace$. Each basic map is a homeomorphism of $S^n=V_n\cup\lbrace\infty\rbrace$ (the inversion $x\mapsto -x^{-1}$ swapping $0$ and $\infty$), and $D_n$ is $V_n$ minus at most the one point $-c^{-1}d$, hence dense in $S^n$; so the continuous maps $g$ and the composite agree on all of $S^n$. \item By Lemma~\ref{lem:identity-scalar}, two matrices in $\GL_2(\CC_n)$ that induce the same element of $\GM(n)$ differ by a factor $z\in Z(\CC_n)$, and $\Delta(zg)=zz^*\Delta(g)$.
		If $n$ is odd, then $z\in\RR^\times$ by \eqref{eqn:center}, so $zz^*>0$.
		If $n$ is even, then $J=i_1\cdots i_{n-1}$ is central, $JJ^*=-1$, and $Jg\in\GL_2(\CC_n)$ by Theorem~\ref{thm:gl2-short}; so if $\Delta(g)<0$, then $Jg$ induces the same element of $\GM(n)$ as $g$ and $\Delta(Jg)=-\Delta(g)>0$.
		\item 
		Consider the function $V_n\setminus \lbrace -c^{-1}d\rbrace \to \CC_n$ given by $g(x) = (ax+b)(cx+d)^{-1}$. 
		Note that for every $x\neq -c^{-1}d$ we have $(cx+d)= c(x+c^{-1}d)$, since nonzero Clifford vectors are invertible. 
\begin{samepage}
		If $d=0$, then $\Delta(g)=-bc^*$, so $b=-\Delta(g)(c^*)^{-1}$ and
		\[
		  g(x)=ac^{-1}+bx^{-1}c^{-1}
		      =ac^{-1}-\Delta(g)(cxc^*)^{-1}.
		\]
		This is a sum of Clifford vectors. Henceforth assume $d\neq0$.
		\end{samepage}
		We know $g(0)\in V_n$ by the hypotheses on $\GL_2(\CC_n)$ and the useful formula Lemma~\ref{lem:useful}. 
		The Magic Formula holds for elements in $\begin{psmallmatrix} a & b \\ c & d \end{psmallmatrix} \in \GL_2(\CC_n)$, and we see that 
$g(x)-g(0) = \Delta(g) (d^*)^{-1}x(cx+d)^{-1} = \Delta(g)[ cd^* + dx^{-1}d^*]^{-1}.$
		Both terms on the right-hand side are Clifford vectors, so $g(x)-g(0) \in V_n$. 
Since $g(0)\in V_n$ this proves $g(x)\in V_n$.
		The induced map is a homeomorphism of $S^m$ for every $m\geq n$: for $c\neq 0$, $g=T_{ac^{-1}}\circ(\text{scaling})\circ(x\mapsto c'xc^{-1})\circ(x\mapsto -x^{-1})\circ T_{c^{-1}d}$, where $T_\mu$ is the translation by $\mu$ (and similarly for $c=0$); each factor is a homeomorphism of $S^m$, and $\GL_2(\CC_n)\subset\GL_2(\CC_m)$. Hence the map $x\mapsto g(x)$ gives an element of $\GM(n)$. Surjectivity of the map comes from Theorem~\ref{thm:gm-to-gl2}.
\item For the first isomorphism, let $Q$ be the form of Theorem~\ref{thm:mcinroy-correction}, which has signature $(1,n+1)$ over $\RR$, and let $^{\circ}$ denote the identity component in the matrix topology. The isomorphism $\psi\colon\SL_2(\CC_n)\cong\Spin_Q(\RR)$ of Theorem~\ref{thm:mcinroy-correction}, followed by the vector representation, has kernel $\lbrace\pm1\rbrace$ and image the kernel of the spinor norm, which is $\SO_{1,n+1}(\RR)^{\circ}$: the time-orientation function $h\mapsto\langle u_0,hu_0\rangle$, for a time-like vector $u_0$, separates the two cosets of $\SO_{1,n+1}(\RR)^{\circ}$ in $\SO_Q(\RR)$. The second isomorphism is \cite[Corollary 2, pg 64]{Ratcliffe2019}---it follows from the hyperboloid model of $\Hcal^{n+1}$ inside $\RR^{n+2}$.
		Note that $\SO_{1,n+1}(\RR)^{\circ}=\PSO(1,n+1)=\PO(1,n+1)\cap\SO(1,n+1)$ in the notation of Definition~\ref{def:so-po-pso}, while Ratcliffe's $\operatorname{PO}(1,n+1)$ is the full group $\Isom(\Hcal^{n+1})$. 
		
		The action of item~1 identifies $\PSL_2(\CC_n)$ with the group of orientation-preserving isometries of $\Hcal^{n+1}$.
Elements of $\SL_2(\CC_n)$ preserve $\Hcal^{n+1}$, and they preserve orientation: by Corollary~\ref{cor:magic-corollaries}, with $w=cx+d$, we have $dg(x)u=\Delta(g)(w^*)^{-1}uw^{-1}=\Delta(g)\vert w\vert^{-2}\,w'uw^{-1}$, a positive multiple of the rotation $u\mapsto w'uw^{-1}$, which has determinant $+1$. Every orientation-preserving M\"obius transformation of $S^n$, hence by Poincar\'e extension every orientation-preserving isometry of $\Hcal^{n+1}$, is induced by an element of $\SL_2(\CC_n)$ \cite[\S7.5, Theorem, parts II and III]{Ahlfors1984}; see also \cite[Theorem 5(ii) and p.~96]{Waterman1993}, whose $\CC_n$ is our $\CC_{n+1}$, and the classification of isometries of the upper half-space as compositions of similarities and inversions \cite[\S4.6]{Ratcliffe2019}. 
Concretely, the generators used are the translations, the dilatations, the maps $x\mapsto uxu$ for unit vectors $u\in V_n$, and $x\mapsto -x^{-1}$, induced by $S=\begin{psmallmatrix}0&-1\\1&0\end{psmallmatrix}$. The kernel of the action is $\lbrace \pm 1\rbrace$: an element acting trivially on $\Hcal^{n+1}$ fixes $S^n$ pointwise, so it is $zI_2$ with $z\in Z(\CC_n)$ by Lemma~\ref{lem:identity-scalar}; it also fixes $i_n$, so $z$ commutes with $i_n$, which forces $z\in\RR$ by \eqref{eqn:center}; finally $\Delta(zI_2)=z^2=1$.
		For the action on $S^n$ alone the kernel can be larger: if $n$ is even, then $J=i_1\cdots i_{n-1}$ is central and $\begin{psmallmatrix}J&0\\0&J\end{psmallmatrix}\in\GL_2(\CC_n)$ acts trivially on $S^n$, but it sends $i_n$ to $-i_n$ and has $\Delta=JJ^*=-1$, so it does not lie in $\SL_2(\CC_n)$.
	\end{enumerate}
\end{proof}

We now make some remarks about the factorization of M\"obius transformations in item~\ref{item:factorization-in-gm}. 
As stated, this is done by Waterman in  \cite[Lemma 10]{Waterman1993}. He shows that every M\"obius transformation $g$ is the composition of basic transformations: translations $x\mapsto x+\mu$ for $\mu \in V_n$; inversion $x\mapsto -x^{-1}$; dilatation $x\mapsto \lambda^2 x$, for $\lambda \in \RR_{>0}$; trivial maps induced by $\begin{psmallmatrix} \lambda & 0 \\ 0 & \lambda \end{psmallmatrix}$ where $\lambda \in \RR$ is nonzero;
special orthogonal transformations (which we call rotations) $x\mapsto ax(a')^{-1}$ for $a \in \CC_n^\times$; and reflections $x\mapsto -x$. The map $x \mapsto axa^*$ equals $\nrd(a)=a\bar a>0$ times the rotation $x\mapsto ax(a')^{-1}$, so it is a rotation composed with a dilatation, and special orthogonal only when $\nrd(a)=1$. As stated above, this is done in the language of Riemannian geometry in \cite[equation (4.3.1), p.~112]{Ratcliffe2019}, without Clifford-M\"{o}bius transformations.

\subsection{Hyperbolic Space as a Symmetric Space}\label{sec:symmetric-space}

The general theory of hyperbolic space as a symmetric space \cite[Section 1.2]{Morris2015} tells us that 
\begin{equation}\label{E:symmetric-space-classic}
\Hcal^{n+1} \cong \O_{1,n+1}(\RR)^{\circ}/\O_{n+1}(\RR)^{\circ}.
\end{equation}

In dimensions $2$ and $3$ we have the famous presentations $\Hcal^2 \cong \SL_2(\RR)/\SO_2(\RR)$ and $\Hcal^3 \cong \SL_2(\CC)/\SU(2).$
In the first case the map $\SL_2(\RR) \to \Hcal^2$ takes some $g\in \SL_2(\RR)$ and maps it to $g(i)$. 
The stabilizer of $i$ is $\SO_2(\RR)$. 
For 3-space, $\Hcal^3$ is usually presented as the set of $(z,\zeta) \in \CC\times \RR$ with $\zeta>0$. 
In this case, the map is given by $g\mapsto g(0,1)$ and $\SU(2)$ is the stabilizer of $(0,1)$ under $\SL_2(\CC)$. 
We now generalize this presentation to higher dimensions.

The first step is to define groups which will be the stabilizers of $i_n$ in $\SL_2(\CC_n)$. 
\begin{definition}
The \emph{Clifford special unitary group} is
\[
  \SU_2(\CC_n)
  =\left\{\begin{pmatrix}a&b\\-b'&a'\end{pmatrix}\in\SL_2(\CC_n)\right\}.
\]
The \emph{Clifford projective special unitary group} is
\[\PSU_2(\CC_n)=\SU_2(\CC_n)/\lbrace\pm1\rbrace.\]
\end{definition}
The groups $\SU_2(\CC_n)$ and $\PSU_2(\CC_n)$
are maximal compact subgroups of $\SL_2(\CC_n)$ and $\PSL_2(\CC_n)$,
respectively, for the topology induced from $M_2(\CC_n)$.
Indeed, let $K\supseteq\SU_2(\CC_n)$ be a compact subgroup of $\SL_2(\CC_n)$ and $g\in K$. There is $u=\begin{psmallmatrix} s & ts^{-1}\\ 0 & s^{-1}\end{psmallmatrix}\in\SL_2(\CC_n)$ with $s\in\RR_{>0}$ and $t\in V_n$, acting by $x\mapsto s^2x+t$, such that $u(i_n)=g(i_n)$; then $u^{-1}g$ fixes $i_n$, so it lies in $\SU_2(\CC_n)\subseteq K$, and hence $u\in K$. The powers of $u$ are unbounded unless $s=1$ and $t=0$, that is, unless $g(i_n)=i_n$; so $g\in\SU_2(\CC_n)$. The same argument applies in $\PSL_2(\CC_n)$. The proof that they are indeed stabilizers of $i_n$ is straightforward and found in \cite[p. 97]{Waterman1993}).

\begin{theorem}\label{T:symmetric-space-clifford}
	As symmetric spaces we have $\Hcal^{n+1}\cong \SL_2(\CC_n)/\SU_2(\CC_n)$.
\end{theorem}
\begin{proof}
This follows from the classical description of $\Hcal^{n+1}$ given as a symmetric space given in \eqref{E:symmetric-space-classic}, together with the description of $\PSL_2(\CC_n) \cong \O_{1,n+1}(\RR)^{\circ}$ and its maximal compact subgroup being $\PSU_2(\CC_n)$, together with the fact that $\SL_2(\CC_n) \to \PSL_2(\CC_n)$ and $\SU_2(\CC_n)\to \PSU_2(\CC_n)$ are maps of degree $2$ with the same kernel.
\end{proof}

Here, as in the case for $\Hcal^2$, there is a simple description of Siegel heights (traditionally defined in terms of the Iwasawa decomposition $G= K \cdot S \cdot U$ where $S$ is a maximal torus and $U$ is the unipotent radical of a minimal parabolic group). 
For $x \in \Hcal^{n+1}$ the \emph{Siegel height} (or just \emph{height}) of $x$ is $\height(x)=x_n$. 

\begin{remark} Let us consider the dimension $d_n$ of the Clifford monoid in
	$\CC_n$ (which is equal to the Clifford group $\CC_n^{\times}$).  
	The expected dimension of $\SL_2(\CC_n)$ in terms of $d_n$ is
	$4d_n - d_n - 2(d_n - n)$, since an element of $\SL_2$ is obtained by
	choosing $4$ elements $a,b,c,d$ of $\CC_n^\mon$ and imposing the conditions
	that the determinant is $1$ and that $ab^*, dc^*$ are vectors (i.e., that
	they lie in an $n$-dimensional subspace). Since
	$\SL_2(\CC_n) \cong \SO_{1,n+1}(\RR)^\circ$ has dimension $\binom{n+2}{2}$,
	we predict from this that $d_n = \binom{n}{2}+1$. 

	This is indeed correct. To see this, note that there is an exact sequence $1 \to \RR^{\times} \to \CC_n^{\times} \to \SO_n(\RR) \to 1$
	given by $a \mapsto \rho_a$ where $\rho_a(x) = ax(a')^{-1} = \pi_a(x)/\vert a \vert^2$ (\cite[Theorem 2]{Waterman1993}).
	Since $\dim(\SO_n(\RR)) = {n \choose 2}$ and $\dim(\RR^{\times})=1$, we have $d_n = {n \choose 2 } +1$.
\end{remark}

In what follows, it is good to keep in mind that $\Spin_{1,n+1}(\RR)$ is connected for $n\geq 1$.
\begin{theorem}\label{thm:h-as-sym-space}
	$\Hcal^{n+1} \cong \Spin_{1,n+1}(\RR)/K$ where $K=\Spin_{n+1}(\RR)$ is a maximal compact subgroup of $\Spin_{1,n+1}(\RR)$.
\end{theorem}
\begin{proof}
Give $\Spin_{1,n+1}(\RR)$ (Definition~\ref{defn:spin}) the topology induced from the real Clifford algebra $\Clf(Q)\otimes\RR$, and let $K=\Spin_{n+1}(\RR)$ be the subgroup corresponding to $\SU_2(\CC_n)$ under the isomorphism $\Spin_{1,n+1}(\RR)\cong\SL_2(\CC_n)$ of Theorem~\ref{thm:mcinroy-correction}. Since $\SU_2(\CC_n)$ is a maximal compact subgroup of $\SL_2(\CC_n)$ (see the discussion before Theorem~\ref{T:symmetric-space-clifford}), $K$ is a maximal compact subgroup, and the result follows from Theorem~\ref{T:symmetric-space-clifford}. \end{proof}

\subsection{Locally Symmetric Spaces and Clifford-Bianchi Modular Spaces}\label{sec:modular-spaces}
We first define our modular spaces $Y(\Gamma)$ and $\Ycal(\Gamma)$ associated to $\Gamma = \PSL_2(\Ocal)$ for $\Ocal$, an order in a $B = \Clf(q)_{\QQ}$ where $q$ is a positive definite integral quadratic form. 
We will make no hypotheses on torsion in our group $\Gamma$. 
We define $\Ycal(\Gamma)$, the \emph{open modular orbifold of level $\Gamma$}, to be the quotient orbifold and $Y(\Gamma)$, \emph{the open modular manifold of level $\Gamma$}, to be its associated coarse space
$$\Ycal(\Gamma) = [\Gamma \backslash \Hcal^{n+1}], \quad  Y(\Gamma) = \vert \Ycal(\Gamma) \vert.$$ 
The object $\Ycal(\Gamma)$ is to be regarded as an object of $\Orb$, the $2$-category of orbifolds, and $Y(\Gamma)$ is a Hausdorff topological space, which is an object of $\Man$, the category of $\mathcal{C}^{\infty}$-manifolds, when $\Gamma$ is torsion-free. 
The stabilizers of $\Ycal(\Gamma)$ are finite since the intersection of a compact subgroup and a discrete subgroup are finite.
It follows that $\Ycal(\Gamma)$ is a Deligne-Mumford orbifold.

Since $\Hcal^{n+1}$ is contractible and comes with a $\Gamma$-action with quotient
$\Ycal(\Gamma)$, we conclude that $\Ycal(\Gamma)$ is a $B\Gamma$ (also called a $K(\Gamma,1)$), i.e., has
$\pi_1 \cong \Gamma$ and all other homotopy groups trivial. 
A rigorous treatment of this for $\Gamma$ not torsion-free requires homotopy sequences for orbifolds and the \emph{Borel construction}, which relates the quotient $\Ycal(\Gamma)$ to fibration of topological spaces $\Gamma \to EG \to BG$, where $\Gamma$ in our case is given the discrete topology.
The definition of homotopy groups for orbifolds is given on page 25 of \cite{Adem2007}, and the Borel construction and homotopy groups of quotient orbifolds is discussed around Proposition 1.51 on page 26 \cite{Adem2007}.

Now we move on to the compactifications. 
We will let $\overline{\Hcal}^{n+1} = \Hcal^{n+1} \cup V_n \cup \lbrace \infty \rbrace$ where $V_n \cup \lbrace \infty \rbrace$ is given the usual topology of the $n$-sphere $S^n$. 
We will call this the \emph{full compactification}.
In the model of hyperbolic space as the unit open ball, this corresponds to taking the closed unit ball. 

We now define the partial Satake compactification, which is the analog of the partial Satake compactifications $\Hcal^2 \cup \QQ \cup \lbrace \infty \rbrace$ and $\Hcal^3 \cup \QQ(\sqrt{-D}) \cup \lbrace \infty \rbrace$ in the classical and Bianchi settings.

\begin{definition}\label{def:satake}
	Let $\Gamma \subset \SL_2(\Ocal)$ be a Clifford-Bianchi group with $\Ocal \subset \CC_n$, an order in a rational Clifford algebra of a positive definite quadratic form.
	Let $B= \Ocal\otimes \QQ$ be the rational Clifford algebra. 
	The \emph{partial Satake compactification} is a topological space
        with underlying set
	$$ \Hcal^{n+1,\Sat} = \Hcal^{n+1} \cup \Vec(B) \cup \lbrace \infty \rbrace. $$
        To define the topology at $\infty$, let the 
        $U(R)= \lbrace x \in \Hcal^{n+1} \colon x_n> R \rbrace \cup \lbrace \infty \rbrace $
        for positive real $R$ be basic open sets. Near $b \in \Vec(B)$, let
        $A$ be an element of $\SL_2(B)$ that takes $b$ to $\infty$. 
        We declare the collection of $A^{-1}(U(R)) $ to be basic open sets near $b$. (This does not depend on the
        choice of $A$, since two different elements differ by an element of
        the stabilizer of $\infty$, which maps each $U(R)$ to some $U(R')$ and so preserves this family of basic open sets.)
\end{definition}
Figure~\ref{fig:satake} shows these horoballs for $n=1$ and a picture of $Y(\PSL_2(\ZZ))$ with its cusp.
\begin{figure}[htbp!]
\centering
\begin{tikzpicture}[scale=2.2,lab]
 \fill[bubbleC!12] (-1.3,1.25) rectangle (1.3,1.6); \draw[bub] (-1.3,1.25)--(1.3,1.25) node[right]{$x_n=R$};
 \foreach \c/\d in {0/0.8,1/0.8,-1/0.8,0.5/0.2,-0.5/0.2,0.333/0.089,-0.333/0.089,0.667/0.089,-0.667/0.089}{
   \draw[bub,fill=bubbleC!12] (\c,\d/2) circle (\d/2);}
 \draw[bnd] (-1.3,0)--(1.3,0);
 \foreach \c in {0,1,-1,0.5,-0.5,0.333,-0.333,0.667,-0.667}{\fill (\c,0) circle (0.9pt);}
 \node at (0,1.45){$U(R)$}; \node at (0,1.55-1.3){};
 \node at (0,-0.25){(a) horoballs in $\Hcal^{2,\Sat}$};
\end{tikzpicture}\qquad
\begin{tikzpicture}[scale=1.0,lab]
 \draw[thick] (-0.25,3.2) .. controls (-0.3,2) and (-1.5,1.2) .. (-1.4,0.4) .. controls (-1.3,-0.4) and (1.3,-0.4) .. (1.4,0.4) .. controls (1.5,1.2) and (0.3,2) .. (0.25,3.2);
 \draw (-0.25,3.2) arc[start angle=180,end angle=0,x radius=0.25,y radius=0.06];
 \draw[dashed] (0,3.25)--(0,3.9); \fill[cuspC] (0,3.95) circle (2.2pt) node[right]{cusp $\infty$ added};
 \fill (-0.7,0.25) circle (1.5pt) node[above]{$i$}; \node at (-0.7,-0.55){order 2};
 \fill (0.7,0.25) circle (1.5pt) node[above]{$\zeta_3$}; \node at (0.7,-0.55){order 3};
 \node at (1.6,2.0){$Y(\Gamma)$}; \node at (0,-1.1){(b) $\Gamma=\PSL_2(\ZZ)$ (cartoon)};
\end{tikzpicture}
\caption{(a) The topology of $\Hcal^{n+1,\Sat}$ in Definition~\ref{def:satake}, drawn for $n=1$: the horoball $U(R)=\{x_n>R\}\cup\{\infty\}$ and its images $A^{-1}U(R)$, horoballs tangent to the boundary at cusps. (b) A cartoon of $Y(\Gamma)=|\Gamma\backslash\Hcal^2|$ for $\Gamma=\PSL_2(\ZZ)$: a sphere with a cusp and two cone points; the Satake compactification adds one point at the cusp. Not to scale.}\label{fig:satake}
\end{figure}

The basic open sets containing elements
  of $\Vec(B) \cup \{\infty\}$ are called {\em horoballs}  
  and their boundaries are {\em horospheres}.  
The horospheres are orthogonal to all geodesics with limit at the corresponding point of
  the boundary: for this, see \cite[pages 127, 132]{Ratcliffe2019}. 
Horospheres based at a point $b\in V_n$ are Euclidean spheres tangent to $\partial \overline{\Hcal}^{n+1}$ at $b$; those based at $\infty$ are the horizontal hyperplanes.
We call the elements of 
  $$\Cusps(B):=\Vec(B) \cup \lbrace \infty\rbrace$$ 
  the \emph{cusps} and those that are in the orbit of $\infty$ under $\Gamma$ we call \emph{principal cusps}.

Here is a theoretical way to arrive at the partial Satake compactification. 
Let $D\subset \Hcal^{n+1}$ be an open fundamental domain for $\Gamma$. 
Let $\overline{D} \subset \overline{\Hcal}^{n+1}$ be the closure of $D$ in the full compactification and let $\partial_\infty D = \overline{D} \cap \partial \Hcal^{n+1}$. 
We then define the boundary to be $\bigcup_{\delta \in \partial_\infty D}\orb_{\Gamma}(\delta)$; this is just adding ideal points of geodesics which are involved in the boundary of our fundamental domain. 
See the discussion in Borel-Ji \cite[end of \S I2]{Borel2006} and \cite[Chapter 12]{Ratcliffe2019} and \cite[Introduction]{Ji2002}.
In the case that all of our cusps of $\Gamma$ are principal, i.e., $\Vec(B) \cup \lbrace \infty\rbrace = \orb_{\Gamma}(\infty)$, it is easy to see that $\bigcup_{\delta \in \partial_\infty D}\orb_{\Gamma}(\delta) = \Cusps(B)$.

For the purpose of comparing our compactifications to what appears in the arithmetic groups literature, we introduce some group schemes.

For simplicity we will informally refer to  group schemes, Lie groups, and groups
as simply ``groups''.
Let $G = \SL_2(\underline{\Ocal})$ be our group scheme defined over $\ZZ$ from Theorem~\ref{thm:group-schemes} (for this we assume, as there, that $\Ocal$ is Clifford-stable). 
Here $\underline{\Ocal}$ is the affine ring scheme defined via Weil restriction such that $\underline{\Ocal}(\ZZ)=\Ocal$. 

We have $G(\ZZ) = \SL_2(\Ocal), G(\QQ) = \SL_2(B)$ and $G(\RR) = \SL_2(\CC_n)$. 
We define the $\ZZ$-group schemes $S,M,P$ and $U$ by their functors of points (c.f. \cite[Definition 5.2]{Elstrodt1988} ) where a ring $R$ maps to
\begin{align*}
S(R) &= \left \lbrace \begin{pmatrix}a & 0 \\ 0 & a^{-1} \end{pmatrix} \colon a \in R^{\times} \right\rbrace, & M(R) &=\left  \lbrace \begin{pmatrix}a & 0 \\ 0 & (a^*)^{-1} \end{pmatrix} \in G(R) \colon a\overline{a}=1 \right \rbrace, \\
P(R) &= \left \lbrace \begin{pmatrix}a & b \\ 0 & (a^*)^{-1} \end{pmatrix} \in G(R)\right \rbrace =\Stab_G(\infty), & U(R) &= \left \lbrace \begin{pmatrix}1 & b \\ 0 & 1 \end{pmatrix} \in G(R)\right \rbrace.
\end{align*}

The group $S$ is the maximal split torus over $\ZZ$ ($S_{\QQ}$ is a maximal split torus of $G_\QQ$, and $S_{\RR}$ is a maximal split torus of $G_{\RR}$). 
The group $M$ is the maximal connected anisotropic subgroup of the centralizer of $S$.

As group schemes we have $P \cong \Vec(\underline{\Ocal})\rtimes \underline{\Ocal}^{
	\times} \cong U\rtimes SM $; this follows from the corresponding fact for groups of matrices over rings. 
This will appear later in equation \eqref{eqn:horocycle-decomposition}.

The set of rational parabolic subgroups of $G_{\QQ}$ is in bijection with $G(\QQ)/P(\QQ)$, and abstractly \cite[Ch.~7]{Borel2019} defines the cusps of $G_{\QQ}$ to be 
$$\Cusps(G) := G(\QQ)/P(\QQ).$$
We can see that these group-theoretically defined cusps are in bijection with our cusps. 
\begin{lemma}
	$\Cusps(G) \cong \Cusps(B)$ as left $\SL_2(B)$-sets. 
\end{lemma}
\begin{proof}
  This is the orbit-stabilizer theorem. The group $G(\Q) = \SL_2(B)$ acts on
  $\Vec(B) \cup \{\infty\}$ in the usual way, and the stabilizer of $\infty$
  is the set of matrices with $c = 0$, which is $P(\QQ)$. So there is a
  bijection of $\SL_2(B)$-sets between the set of cosets $\Cusps(G)$ and the
  orbit of $\infty$ in $\Cusps(B)$. However, the action of $G(\Q)$ on
  $\Cusps(B)$ is transitive, since for all $v \in \Vec(B)$
  the matrix $\begin{psmallmatrix}v & v-1 \\ 1 & 1 \end{psmallmatrix}$ takes
  $\infty$ to $v$.

\end{proof}
In section~\ref{sec:cusps} we prove that $\Gamma\setminus \Cusps(B)$ is finite. 

We wish now to describe some group-theoretic compactifications for the purpose of comparing them to our compactifications. 
To describe the geometry of these compactifications we need to discuss some geometry of geodesics.
For each of our cusps $c$ we can give local coordinates in terms of geodesics. 
One collection of parameters will allow us to select which geodesic we wish to follow and the other coordinate moves along the said geodesic.

Before proceeding to some terse group theory, we review the familiar case of $\overline{\Hcal}^2$.
The case of $\overline{\Hcal}^2$ is shown in Figure~\ref{fig:horocycle-plane}.
There will be horospheres, which in the case of $\overline{\Hcal}^2$ are horocircles and are the Euclidean circles in Figure~\ref{fig:horocycle-plane} tangent to the boundary $\partial \overline{\Hcal}^2 \cong \RR$. 
These horospheres are orthogonal to the geodesics and hence choosing a point on the circle amounts to choosing a geodesic terminating at the cusp.
In the general case the horoballs will be generalizations of the filled-in horocircles.
\begin{figure}[htbp!]
	\centering
	\begin{tikzpicture}[scale=1.15, line join=round]
\foreach \r in {0.3,0.6,0.9,1.2,1.5}
    \draw[thick, black!55] (0,\r) circle (\r);
\begin{scope}
    \clip (-3.3,0) rectangle (3.3,3.4);
    \draw[thin] (0,0) -- (0,3.4);
    \foreach \c in {0.4,0.8,1.3,1.9,2.7}
      { \draw[thin] (\c,0) circle (\c); \draw[thin] (-\c,0) circle (\c); }
  \end{scope}
  \draw[very thick] (-3.4,0) -- (3.4,0) node[right] {$V_n$};
  \fill (0,0) circle (1.8pt) node[below=2pt] {$x$};
\draw[<-, thin] (1.4,2.37) -- (1.75,3.3) node[right] {geodesics ending at $x$};
  \draw[<-, thin] ({-1.5*sin(50)},{1.5+1.5*cos(50)}) -- (-2.0,3.55) node[left] {horocycles at $x$};
\end{tikzpicture}
	\caption{The horocycle decomposition at a cusp $x$ on the boundary, the image of Figure~\ref{fig:horocycle-infty} under an element $g$ with $g(\infty)=x$.
		Pictured in $\Hcal^2$: the horocycles at $x$ are the circles tangent to the boundary at $x$, and the geodesics ending at $x$ are the vertical line through $x$ and the semicircles through $x$ centred on the boundary.
		Each geodesic meets each horocycle orthogonally.
	}\label{fig:horocycle-plane}
\end{figure}

The geometry of horoballs and horospheres is given in terms of the Iwasawa decomposition, which we now review.
We have $S(\RR) \cong \RR^{\times}$, $S(\RR)^{\circ} = \RR_{> 0}$, $U(\ZZ) \cong \Vec(\Ocal)$, $U(\QQ) \cong \Vec(B)$, $U(\RR) \cong V_n$.
We have $M(\ZZ) \cong \Ocal^{\times}$, $M(\QQ) = B^{\times,1}$, and $M(\RR) = \CC_n^{\times,1}$ where $C^{\times,1} = \lbrace a \in C^{\times} \colon \vert a \vert =1 \rbrace$ for $C \subset \CC_n$, a subring.
We will let 
$$M=M(\RR), \quad  A=S(\RR)^{\circ}, \quad N = U(\RR),$$
so that $P=MAN$ is the Langlands decomposition.

All of these groups can either be defined starting from the parabolic group $P$ or the split torus $S$, and in the literature they are often given subscripts $P$ so $A=A_P = A_{P_{\infty}}=A_{\infty}$ and $N = N_P =N_{P_{\infty}}=N_{\infty}$. 
When we change the cusp $\infty$ to $c$ (or equivalently when we change the rational parabolic subgroup) ,we get $P_c$, $A_c$ and $N_c$, respectively.

Continuing with this decomposition at $\infty$, note that for $t \in A = \RR_{>0}^{\times}$ we write $\pi_t(x) = t^2x$ for the action of $\begin{psmallmatrix}t&0\\0&t^{-1}\end{psmallmatrix}$ (this is not the rotation $\pi_a$ of a Clifford group element), and for $v\in N$ we have $\tau_v(x) = x+v$. 
The \emph{horocycle decomposition} of $\Hcal^{n+1}$ is the map
\begin{equation}\label{eqn:horocycle-decomposition}
  A\times N \to \Hcal^{n+1}, \quad (t,v) \mapsto \tau_v(\pi_t(i_n)),
 \end{equation}
which is a specific case of a general construction in \cite{Borel2006}. 
Note that $v$ parametrizes the horocycles, and $t$ is the parameter of the geodesic flow. 
Also, $\lim_{t\to\infty} \tau_v(\pi_t(i_n))=\infty$, which is the cusp that we started with.
If $g\in \SL_2(B)$ has $g(\infty)=c \in \Cusps(B)$, the horocycle decomposition moves with it, and Figure~\ref{fig:horocycle-infty} becomes Figure~\ref{fig:horocycle-plane} where the planes (horospheres) move to the circles tangent to the boundary at $c$, and 
the vertical geodesics move to the semicircles through $c$ centred on the boundary, which meet the horocycles at $c$ orthogonally.
\begin{figure}[htbp!]
	\centering
	\begin{tikzpicture}[x=1cm, y=1cm, line join=round]
\def\plane#1#2{\filldraw[draw=black, fill=#2, fill opacity=0.85]
      (0,#1) -- ++(4,0) -- ++(1.6,1.1) -- ++(-4,0) -- cycle;}
\def\pts{0.13/0.3, 0.17/0.7, 0.57/0.2, 0.61/0.6, 0.89/0.4}
  \def\seg#1#2{\foreach \u/\v in \pts
      \draw[thick] ({4*\u+1.6*\v},{1.1*\v+#1}) -- ({4*\u+1.6*\v},{1.1*\v+#2});}
  \plane{0}{black!15}
  \foreach \u/\v in \pts \fill ({4*\u+1.6*\v},{1.1*\v}) circle (1.3pt);
  \seg{0}{1}
  \plane{1}{white}
  \seg{1}{1.9}
  \plane{1.9}{white}
  \seg{1.9}{2.7}
  \plane{2.7}{white}
  \foreach \u/\v in \pts
      \draw[thick, ->] ({4*\u+1.6*\v},{1.1*\v+2.7}) -- ({4*\u+1.6*\v},{1.1*\v+3.6});
\node[right] at (5.7,0.55) {boundary $V_n\subset\partial\Hcal^{n+1}$};
  \node[right] (hl) at (5.9,2.1) {horospheres $\cong V_n$};
  \foreach \h in {1,1.9,2.7} \draw[very thin, gray] (hl.west) -- (5.45,{\h+0.95});
  \node[above] at (3.3,4.45) {geodesics, toward the cusp $\infty$};
\draw[->, thick, black!60] (-0.5,1.2) -- (-0.5,3.0) node[midway, left] {$t\in A$};
  \draw[->, thick, black!60] (0.3,-0.35) -- (3.7,-0.35) node[midway, below] {$v\in V_n$};
\end{tikzpicture}
	\caption{The horocycle decomposition at $\infty$, pictured in the upper half-space model.
		The horizontal planes are horospheres based at $\infty$; the shaded plane is the boundary $V_n$.
		The vertical lines are geodesics ending at $\infty$.
		Acting by $t\in A=A_{P_{\infty}}$ moves along the geodesics, and acting by $v\in V_n$ moves within a horosphere.
	}\label{fig:horocycle-infty}
\end{figure}

In 
\cite{Borel1973}, Borel and Serre define a compactification $H^{\BS}$ for symmetric spaces $H$ such that for torsion-free arithmetic groups $\Gamma$ the inclusion $\Gamma \setminus H \to \Gamma \setminus H^{\BS}$ is a homotopy equivalence. 
The Borel-Serre compactification of $\Hcal^{n+1}$ comes from the horocycle decomposition associated to a boundary point $c\in \Cusps(B)$.
The idea is to compactify the geodesic flow by letting $\lim_{t\to\infty} \tau_v(\pi_t(i_n))$ be an element of $\overline{A}\times N$. Here we use $s = 1/t$ so that $s=0$ corresponds to $t=\infty$. The chart at infinity $H_{\infty}$ of the Borel-Serre compactification is given by 
\begin{equation}
H_{\infty}=\overline{A_{\infty}} \times N_{\infty},
\end{equation} 
and, more generally, $H_c = \overline{A_c} \times N_c$. 
At $\infty$ we have $H_{\infty} \cong \RR_{\geq 0} \times V_n.$ 
In this chart $\RR_{>0} \times V_n \cong \Hcal^{n+1} $ via $(s,v) \mapsto v+s^{-2}i_n$ (since $\tau_v(\pi_t(i_n))=v+t^2i_n$ and $s=1/t$), so we have compactified by adding all of the directions $v$ in which a geodesic can approach $\infty$. 
We do the same thing for every $c\in \Cusps(B)$ and Figure~\ref{fig:horocycle-plane} becomes very striking. 

\begin{definition}\label{def:pBS}
The \emph{partial Borel-Serre compactification} of $\Hcal^{n+1}$ relative to $\Gamma$ is defined to be 
$\Hcal^{n+1,\BS} = \bigcup_{c\in \Cusps(B)} H_c$.
In other words, it is 
$(\Cusps(B)\times \RR_{\geq 0} \times V_n)/\sim$ where $(c_0,t_0,v_0) \sim (c_1,t_1,v_1)$ if and only if $g_{t_0,v_0}^{c_0}(i_n) = g_{t_1,v_1}^{c_1}(i_n)$ when $t_0\neq 0$ and $t_1 \neq 0$. 
Here $(c,t,v) \mapsto g^c_{t,v}$ takes the tuple to its corresponding element $g^c_{t,v} \in G(\RR)$ given by the Iwasawa decomposition at the cusp $c$.
\end{definition}
The charts $\varphi_c:H_c \xrightarrow{\sim} \RR_{\geq 0} \times \RR^n$ make $\Hcal^{n+1,\BS}$ a manifold with corners. In other words, every point has a neighborhood
isomorphic to $\RR^{n+1,j} = \RR_{\geq 0}^j \times \RR^{n+1-j}$ for some $j\in\lbrace 0,1\rbrace$
(the rational rank is one, so there are no corners of higher codimension).

The prototypical manifold with corners is the closed positive ``quadrant'' $\RR_{\geq 0}^n$ where the non-negative parts of the coordinate hyperplanes are included.
The basic properties of the category of manifolds with corners $\Man^c$ is developed in \cite{Joyce2012} and \cite[\S8.5-8.9]{Joyce2014}.
In particular, there is a 2-category of orbifolds with corners $\Orb^c$ which can be developed from the perspective of charts, \'etale proper groupoids, and representable stacks in the category $\Man^c$ where the Grothendieck topology is given by open coverings. 
The theory is well-summarized in \cite[example 3.4]{Solomon2020} where the theory of quotient orbifolds is treated. 
The \emph{Borel-Serre compactification} of $\Ycal(\Gamma)$ and $Y(\Gamma)$ are then defined to be objects of $\Orb^c$ and $\Man^c$, respectively, given by 
	 $$ \Xcal^{\BS}(\Gamma) = [\Gamma \setminus \Hcal^{n+1,\BS} ], \quad X^{\BS}(\Gamma) = \vert \Xcal^{\BS}(\Gamma) \vert.$$
We can now see that our orbifolds are more complicated to deal with because $\Xcal(\Gamma) :=[\Gamma \setminus \Hcal^{n+1,\Sat}]$ and $X(\Gamma) = \vert \Xcal(\Gamma) \vert$; as we now explain, they are not orbifolds or manifolds with corners in the above sense. 
We see that the charts $H_c$ of $\Hcal^{n+1,\BS}$ 
are now blown down as shown in Figure~\ref{fig:cone}(b). 
This can be pictured as a 
successive collapsing of cones as pictured in Figure~\ref{fig:cone}(a), which readers may be more familiar with from Hatcher \cite{Hatcher2002}.

\begin{figure}[htbp!]
	\centering
	\begin{tikzpicture}[line join=round]
\begin{scope}
    \draw[dashed] (-1,0) -- (-1,2.2); \draw[dashed] (1,0) -- (1,2.2);
    \foreach \x in {-0.6,-0.2,0.2,0.6} \draw[thin] (\x,0) -- (\x,2.2);
    \draw[very thick] (-1,0) -- (1,0);
    \node[below] at (0,-0.1) {$I\times[0,\infty)$};
    \node[left] at (-1.05,1.1) {$\scriptstyle -b$};
    \node[right] at (1.05,1.1) {$\scriptstyle b$};
    \draw[->, thick] (1.8,1.1) -- (3.2,1.1) node[midway, above] {quotient};
    \begin{scope}[shift={(4.6,0)}]
      \draw[dashed] (0,0) -- (110:2.35); \draw[dashed] (0,0) -- (70:2.35);
      \foreach \a in {78,86,94,102} \draw[thin] (0,0) -- (\a:2.3);
      \fill (0,0) circle (2pt);
      \node[below] at (0,-0.1) {$\bigl(I\times[0,\infty)\bigr)\big/\bigl(I\times\{0\}\bigr)$};
    \end{scope}
    \node at (2.3,-1.0) {(a)};
  \end{scope}
\begin{scope}[shift={(7.6,0)}]
    \foreach \x in {-1.3,-0.9,...,1.3} \draw[thin] (\x,0) -- (\x,2.2);
    \draw[very thick] (-1.5,0) -- (1.5,0);
    \node[below] at (0,-0.1) {$\RR_{\geq 0}\times\RR$};
    \draw[->, thick] (1.9,1.1) -- (3.3,1.1) node[midway, above] {quotient};
    \begin{scope}[shift={(4.9,0)}]
      \foreach \a in {12,24,...,168} \draw[thin] (0,0) -- (\a:2.0);
      \draw[very thin, gray, dotted] (-2.1,0) -- (2.1,0);
      \fill (0,0) circle (2pt);
      \draw[->, thin] (0.7,0) arc (0:50:0.7) node[midway, right] {$\theta$};
      \node[below] at (0,-0.1) {$(\RR_{\geq 0}\times\RR)\big/(\{0\}\times\RR)$};
    \end{scope}
    \node at (2.6,-1.0) {(b)};
  \end{scope}
\end{tikzpicture}
	\caption{(a) Pinching the boundary $I\times\{0\}$ of the strip $I\times[0,\infty)$, where $I=(-b,b)$, to a point gives a cone.
		(b) Pinching the boundary $\{0\}\times\RR$ of $\RR_{\geq 0}\times\RR$ to a point; the result is parametrized by $(r,\theta)$ with $r\geq 0$ and $\theta\in(0,\pi)$.
	}\label{fig:cone}
\end{figure}

We can see that the collapsed chart in two dimensions 
(Figure~\ref{fig:cone}(b)) can be parametrized by $(r,\theta)$ with $\theta \in (0,\pi)$ and $r\geq 0$ where all of $(0,\theta)$ are identified. One can then see that these are equivalent to using charts which are locally homeomorphic to open subset of $\RR_{>0}\times \RR_{>0} \cup \lbrace (0,0)\rbrace$, which do not fall within the realm of manifolds with corners and are something else. 
To our knowledge, there is no theory of manifolds and orbifolds with ``$o$-minimal corners'' and the development of this theory is well beyond the scope of this paper. 

We want to close this subsection with a comparison between the compactifications given and others existing in the literature. 
For the purposes of exposition, let us call the previously introduced the Satake compactifications the ``easy Satake'' compactifications. 
We will call the Satake compactification as used at the end of \cite[Chapter 17]{Borel2019} the ``hard Satake'' compactification. 
We will refer to the compactification used by Borel in \cite[Chapter 17]{Borel2019} and by EGM in \cite{Elstrodt1988} as the ``Borel compactification''. 

In Borel's construction from \cite[Chapter 17]{Borel2019}, he selects a finite set of representatives for $\Cusps(G) = G(\QQ)/P(\QQ)$, which we will call $\lbrace c_1,\ldots,c_h\rbrace$. 
Let $\Gamma^c=c^{-1} \Gamma c$, and $\Gamma^c_{\infty}$ is an extension of $\Ocal_c^{\times}$ by $\Lambda_c$ 
 where $\Ocal_c^{\times} = \Gamma_{\infty}^c\cap M(\QQ)$ (contrary to what may be suggested
by the notation, we do not know whether there is actually an order whose group of units
is $\Ocal_c^{\times}$).

Over every cusp $c\in \Cusps(G)$, one has a component 
$$E_c = \Lambda_c \backslash M(\RR)U(\RR)/K\cap M(\RR)$$ 
on which 
$$L_c = (\Gamma^c \cap P(\QQ))/(\Gamma^c \cap U(\QQ)) \cong \Gamma_{\infty}^c/\Lambda_c \cong \Ocal_c^{\times}$$ acts. 
Using the semidirect product structure, we see that 
$$E_c \cong (\Lambda_c \backslash U(\RR)) \times (M(\RR)/M(\RR)\cap K) \cong (V_n/\Lambda_c)\times 1.$$ 
Generally, there is some fibration $\sigma: E_c \to B_c$ with $B_c = M(\RR)/(M(\RR)\cap K)$ which is nontrivial, but in the rank one case this is always trivial, as can be verified directly ($M(\RR)$ is compact so $K\cap M(\RR)=M(\RR)$). 
So, in summary, the easy Satake compactification is the blow-down of Borel-Serre.
This language of ``blow-down'' is the language used by \cite{Borel1973} and means that, in the map from Borel-Serre to Satake, the fibers over the cusps collapse to a point. 
It is also used to indicate that the cusp in the Satake compactification is replaced by many points, so that each geodesic terminating in the cusp in the Satake compactification now terminates in different points, each  identified with the asymptotic direction in which the geodesics approached the cusp.

On the level of symmetric spaces, $\Hcal^{n+1,\BS}$ is obtained from $\Hcal^{n+1}$ by adding the face $\lbrace 0\rbrace\times N_c \cong V_n$ of the chart $H_c$ at each $c \in \Cusps(B)$.

\subsection{Satake Abelian Varieties}\label{sec:satake-abelian}
Satake gives a way to parametrize abelian varieties with certain specified endomorphism algebra \cite[Proposition 3 and material before]{Satake1966}.
Using his theory, we can show that there is a construction depending on certain parameters where every $z \in \Hcal^{n+1}$ corresponds to abelian varieties with ``Clifford multiplication''. 
Satake's construction generalizes both the Kuga-Satake construction, which associates to K3 surfaces abelian varieties of large dimension, and the theory of Shimura curves $[\Gamma \backslash \Hcal^2]$, which parametrize abelian surfaces with quaternionic multiplication (see \cite[Ch.~43]{Voight2020} for an overview).
Note in particular that the Shimura construction is closely related to integral ternary quadratic forms $Q$ via its even Clifford algebra $\Clf(Q)_+$; in the Satake generalization, we will use an indefinite form $Q$ for our abelian varieties. 

Given a quadratic form $Q$ we consider the symmetric space $H=\Spin_Q(\RR)/K$ where  
$K$ is its maximal compact subgroup. 
Given an order $\Ocal \subset B=\Clf(Q)_{\QQ}$, we let $\Lambda = \Ocal_+ = \Ocal \cap B_+$ and $V = (B_+)_{\RR}$.
In the case that $Q$ has signature $(1,n+1)$ we have $H\cong \Hcal^{n+1}$.

\begin{theorem}[{\cite{Satake1966}}]
  The torus $V/\Lambda$ admits a family of complex structures $J_z: V \to V$ parametrized by $z \in H$, satisfying $J_z^2=-\id$ for
  every $z\in H$ and admitting a Riemann form. 
	
	In the case of signature $(1,n+1)$ a complex structure $J_0$ is determined by $b_1$ and $b_2$ described by the conditions in equations \eqref{eqn:satake01} and \eqref{eqn:satake02}; it varies according to conjugation $g J_0 g^{-1}$ for $g \in \Spin_Q(\RR)$. 
	A Riemann form $E$ is given by a suitable choice of $\mu \in \Ocal_+$ described in equation \eqref{eqn:satake03}.
\end{theorem}

In the case that the signature of $Q$ is $(1,n+1)$, the complex structure $J$, which depends on parameters $b_1$ and $b_2$, is given by 
\begin{equation}\label{eqn:satake01} J(x) = xb_1 + i x b_2, \quad i = i_1i_2\ldots i_{n+1}, \quad  b_1 \in \Ocal_+, \quad b_2 \in \begin{cases} 
\Ocal_+, & n+1 \mbox{ even }\\
\Ocal_-, & n+1 \mbox{ odd } 
\end{cases}.
\end{equation}
\begin{equation}\label{eqn:satake02}
 b_1^2 + (-1)^{{n+2 \choose 2}} b_2^2 =-1, \quad b_1b_2+b_2b_1=0.
 \end{equation}
This makes $V=\CC_{n+2}[j]_+$ into a complex vector space where $j$ is the generator of $C_{n+1,1}$ such that $j^2=1$.  
Give the lattice $\Ocal_+$ a Riemann-form $E:\Lambda\times \Lambda\to \ZZ$ given by
\begin{equation}\label{eqn:satake03}
 E(x,y) = \trd(\mu x^* y), \qquad \mu \in \Ocal_+, \quad \mu^* = -\mu.
\end{equation}

The general theory of Riemann forms is summarized succinctly in \cite[43.4.9,43.4.10]{Voight2020}. 
Recall that this means $E: \Lambda\times \Lambda \to \ZZ$ is alternating, and $\ZZ$-bilinear with $\RR$-linear extension $E_{\RR}: V\times V \to \RR$ such that $E_{\RR}(Jx,Jy) = E_{\RR}(x,y)$ for all $x,y\in V$ and the map $(x,y) \to E_{\RR}(Jx,y)$ is a symmetric positive definite $\RR$-bilinear form on $V$. 
These two conditions are equivalent to $H:V\times V \to \RR[J]$ defined by $H(x,y) = E_{\RR}(Jx,y) + J E_{\RR}(x,y)$ being a positive definite Hermitian form on $V$ with $\Im(H)=E_{\RR}$; so, given $H$ which is Hermitian positive definite on $V$ with $\Im H(\Lambda) \subset \ZZ$, then $\Im H\vert_{\Lambda}$ is a Riemann form for $(V,\Lambda)$. 

By taking $[g] \in \Spin_Q(\RR)/K$ one can vary the complex structure by taking $J^g = gJg^{-1}$, where $g$ acts on $V$ by left multiplication. 
In our application we have $\Spin_{Q}(\RR)/K \cong \SL_2(\CC_n)/\SU_2(\CC_n) \cong \Hcal^{n+1}$ with $[g] \mapsto g(i_n) = z$ so $J_z = J^g$, and this map is well-defined. 
The $Q$ is the $Q$ coming from Bott periodicity, so $Q=x^2-yz+q$ where $q$ is our positive definite quadratic form of rank $n-1$. 
This means $Q$ has real signature $(1,n+1)$ and we have $\Spin_Q(\RR)/K \cong \Spin_{1,n+1}(\RR)/K \cong \SL_2(\CC_n)/\SU_2(\CC_n) \cong \Hcal^{n+1}$. 
So to every $z\in \Hcal^{n+1}$ we get a complex abelian variety $A_z$ such that $A(\CC)\cong (\Clf(Q)_{\RR})_+/\Lambda$ where the complex structure and polarization are specified. 
Note that the rank of $Q$ is $m=3+(n-1)=n+2$, so $\dim_{\RR} V = 2^{m-1}=2^{n+1}$. 
This implies the abelian variety $A_z$ has complex dimension $2^{m-2} = 2^{n}$; hence via this Satake construction of abelian varieties (depending on choices $\mu$, $b_1,b_2$) the symmetric space $\Hcal^{n+1}$ actually gives polarized abelian varieties $A_z$ of dimension $2^{n}$, and for every $a\in\Ocal_+$ with $ab_1=b_1a$ and $ab_2=b_2a$ the right multiplication $x\mapsto xa$ is an endomorphism of $A_z$ for every $z\in \Hcal^{n+1}$, since it preserves $\Lambda$ and commutes with $J_z$. 

Note that quotients of $\Hcal^2$ parametrize complex abelian surfaces with a fixed quaternionic multiplication, while the points of $\Hcal^3$ give complex abelian fourfolds $A_z$ with $A_z(\CC)\cong V/\Ocal_+$ and $\dim_{\QQ}(\Ocal_+\otimes_{\ZZ}\QQ)=8$ for some fixed order $\Ocal$. 

\section{Arithmeticity of $\SL_2(\Ocal)$}\label{sec:arithmeticity}

The aim of this section is to prove that $\SL_2(\Ocal)$ can be viewed as an arithmetic subgroup of $\SO_{1,n+1}(\RR)$. 
This is done in Theorem~\ref{thm:arithmeticity}.
This has consequences for our fundamental domains because the Borel--Harish-Chandra Theorem then implies that the group action is discrete and finite covolume. 
This has abstract consequences for the construction of our fundamental domain constructions in section~\ref{sec:fundamental-domains} (although we don't make use of the general constructions of Fundamental domains of Borel in \cite{Borel2019}, which are difficult to work with in practice and do not mirror the classical constructions in $\Hcal^2$ and $\Hcal^3$).

Discussions of arithmeticity of the groups $\PSL_2(\Ocal)$ for $\Ocal=\ZZ[i_1,\ldots,i_{n-1}]$ exist in the literature but notation and terminology varies. 
Because of this, in section \ref{sec:arithmetic-groups} we give a definition of arithmetic subgroups of semisimple real Lie groups in terms of group schemes over the integers.
We are also careful to emphasize that this notion is a property of a \emph{subgroup} $\Gamma$ of a real Lie group and not an abstract property of a group $\Gamma$.

\subsection{Arithmetic Groups}\label{sec:arithmetic-groups}

We take \cite[Definition 5.1.19]{Morris2015} as a basis for our definition of arithmetic group---the only change (other than style) is that we are careful not to take $\ZZ$-points of $\QQ$-group schemes (see Appendix~\ref{sec:group-schemes} for a discussion of $\ZZ$-points of $\QQ$-group schemes and why this is deprecated).

\begin{definition}\label{D:arithmetic}
	Let $H$ be a semisimple Lie group and $\Gamma$ be a subgroup. 
	We say that $\Gamma<H$ is an \emph{arithmetic subgroup} or $(\Gamma,H)$ is an \emph{arithmetic pair} if and only if one of the following holds:
	\begin{enumerate}
		\item (Base Case: Integer Points of Group Schemes) There exists an integer $n$ and a closed $\ZZ$-subgroup scheme $G \subset \SL_n$ such that there exists an isomorphism of Lie groups $\phi: H \stackrel{\cong} \to G(\RR)$ such that $\phi(\Gamma) = G(\ZZ)$. 
		\item (Taking Commensurable Groups) The group $\Gamma<H$ is commensurable to $\Gamma'<H$ for $(\Gamma',H)$, an arithmetic pair.
		\item (Taking Connected Components) The pair $(\Gamma,H)$ comes from taking connected components of a known arithmetic pair $(\Gamma',H')$, i.e., $(\Gamma,H) \cong (\Gamma' \cap (H')^0, (H')^0).$
		\item (Taking Quotients by Compacts) The pair $(\Gamma,H)$ comes from taking the quotient of a known arithmetic pair $(\Gamma',H')$ by a compact normal subgroup of $K'<H'$, i.e. $(\Gamma,H) \cong (\Gamma'/(\Gamma'\cap K'),H'/K').$
	\end{enumerate}
\end{definition}

\begin{remark}
Note that $K'$ can be taken to be a finite group above. 
Thus if $(\Gamma',H')$ is an arithmetic pair and $K'<H'$ is a finite normal subgroup, then $\Gamma'/(\Gamma'\cap K')$ is an arithmetic subgroup of $H'/K'$.
\end{remark}

\begin{remark}\label{rem:arithmetic-iso-invariant}
If $(\Gamma,H)$ is an arithmetic pair and $\alpha:H\to H'$ is an isomorphism of Lie groups, then $(\alpha(\Gamma),H')$ is an arithmetic pair. Indeed, each of the clauses (1)--(4) of Definition~\ref{D:arithmetic} is preserved by composing the isomorphism appearing in it with $\alpha$, so this follows by induction on the number of clauses used.
\end{remark}

\begin{remark}
 When $(\Gamma,H)$ is an arithmetic pair we sometimes call $\Gamma$ simply an arithmetic group. 
	While this terminology, which omits the ambient group, is often used in the literature, we emphasize that being a \emph{subgroup} of a particular group is an important part of the definition.
    We point out some variations.
	In \cite{Krausshar2004}, Krau\ss har takes ``arithmetic group'' to mean discrete subgroup (without any reference to an algebraic group).
	The authors of \cite{Maclachlan1989} and \cite{Maclachlan2003} define an ``arithmetic group'' to be commensurable with the $\ZZ$-points of some group scheme over $\ZZ$ (the emphasis here is on the $\exists$ quantifier on the group). 
	In \cite{Morris2015} (and most other references on arithmetic groups), an arithmetic group is taken to be the $\ZZ$-points of a group scheme defined over $\QQ$.
\end{remark}

\begin{theorem}[Arithmeticity]\label{thm:arithmeticity}
	Let $\varphi:\PSL_2(\CC_n) \to \O_{1,n+1}(\RR)^{\circ}$ be the isomorphism of groups given first by acting via M\"obius transformations on $\Hcal^{n+1}$ and then by identifying orientation-preserving isometries of $\Hcal^{n+1}$ with $\O_{1,n+1}(\RR)^{\circ}$.
	Let $q=d_1x_1^2+\cdots + d_{n-1}x_{n-1}^2$ be a positive definite quadratic form over $\ZZ$, let $\Ocal$ be an order in $\quat{-d_1,\ldots,-d_{n-1}}{\QQ}$ which is closed under $*$ (for example $\Ocal = \quat{-d_1,\ldots,-d_{n-1}}{\ZZ}$), and let $\Gamma = \PSL_2(\Ocal)$ be the image of $\SL_2(\Ocal)$ in $\PSL_2(\CC_n)$.
	The group $\varphi(\Gamma)$ is an arithmetic subgroup of $\O_Q(\RR)^{\circ}$ where $Q=x^2-yz +q$.
\end{theorem}
\begin{proof}
	Write $\Lambda=\quat{-d_1,\ldots,-d_{n-1}}{\ZZ}=\Clf_q(\ZZ)$.
	In this proof we work with two isomorphisms of Lie groups $\PSL_2(\CC_n)\to\O_Q(\RR)^{\circ}$: the first is the representation $\varphi$ coming from the M\"obius action on $\Hcal^{n+1}$, and the second, which we will call $\psi$, is induced on real points by the morphism $\SL_2(\Clf_q) \to \Spin_Q \to \O_Q$ of
   Theorem~\ref{thm:mcinroy}, which we will call ``the Spin isomorphism''.
   Of $\psi$ we only use that it is a homomorphism of $\QQ$-groups $\SL_2(\Clf_q)_\QQ\to\SO_{Q,\QQ}$ whose real points map $\SL_2(\CC_n)$ onto $\O_Q(\RR)^{\circ}$ with kernel $\{\pm 1\}$.

  \emph{Transport.} Let $\sigma = \varphi\psi^{-1}$, an automorphism of the Lie group $\O_Q(\RR)^{\circ}$. Then $\sigma(\psi(\Gamma))=\varphi(\Gamma)$, so by Remark~\ref{rem:arithmetic-iso-invariant} it suffices to show that $\psi(\Gamma)$ is an arithmetic subgroup of $\O_Q(\RR)^{\circ}$.

  \emph{Reduction to $\Ocal=\Lambda$.} The orders $\Ocal$ and $\Lambda$ are $\ZZ$-lattices of full rank in $\quat{-d_1,\ldots,-d_{n-1}}{\QQ}$, so writing a $\ZZ$-basis of each in terms of a $\ZZ$-basis of the other gives an integer $N\geq 1$ with $N\Lambda\subset\Ocal$ and $N\Ocal\subset\Lambda$.
  Since $\Ocal$ is closed under $*$, the group $\SL_2(\Ocal)$ lies in $M_2(\Ocal)^\times$, and reduction modulo $N$ is a homomorphism $\SL_2(\Ocal)\to M_2(\Ocal/N\Ocal)^\times$ to a finite group.
  Its kernel has finite index in $\SL_2(\Ocal)$ and is contained in $1+NM_2(\Ocal)\subset M_2(\Lambda)$, hence in $\SL_2(\Lambda)$.
  Exchanging the roles of $\Ocal$ and $\Lambda$, we see that $\SL_2(\Ocal)\cap\SL_2(\Lambda)$ has finite index in both $\SL_2(\Ocal)$ and $\SL_2(\Lambda)$ (compare \cite[Proposition 5.4.5]{Morris2015}).
  Hence $\psi(\PSL_2(\Ocal))$ and $\psi(\PSL_2(\Lambda))$ are commensurable, and by part (2) of Definition~\ref{D:arithmetic} it suffices to treat $\Ocal=\Lambda$.

  \emph{The case $\Ocal=\Lambda$.} Let $\SO_Q\subset\SL_{n+2}$ be the closed $\ZZ$-subgroup scheme of matrices preserving $Q$. By parts (1) and (3) of Definition~\ref{D:arithmetic}, $\SO_Q(\ZZ)\cap\O_Q(\RR)^{\circ}$ is an arithmetic subgroup of $\O_Q(\RR)^{\circ}$.
  Embed $\SL_2(\CC_n)$ in $\GL(M_2(\CC_n))$ by left multiplication.
  Its integral points with respect to the lattice $M_2(\Lambda)$ form
  \[
    \SL_2(\Lambda)=\SL_2(\Clf_q)(\ZZ)\cong\Spin_Q(\ZZ).
  \]
Since $\psi$ is a surjective homomorphism defined over $\QQ$, \cite[Proposition 9.4.5]{Morris2015} (see also \cite[Theorem 4.1]{Platonov1994}) 
shows that $\psi(\SL_2(\Lambda))$ is commensurable with $\SO_Q(\ZZ)$, and hence with $\SO_Q(\ZZ)\cap\O_Q(\RR)^{\circ}$, which contains it.
  By part (2) of Definition~\ref{D:arithmetic}, $\psi(\PSL_2(\Lambda))$ is an arithmetic subgroup of $\O_Q(\RR)^{\circ}$.

  \begin{samepage}
  We can say more, using only $\ZZ$-points of $\ZZ$-group schemes. The short exact sequence of fppf sheaves of groups over $\Spec(\ZZ)$ of Theorem~\ref{thm:short-exact-sequence} gives the exact sequence
	\[
	  1 \to \mu_2(\ZZ) \to \Spin_Q(\ZZ) \to \PSpin_Q(\ZZ)
	  \xrightarrow{\delta} H^1(\Spec(\ZZ)_{\fppf},\mu_2),
	\]
  where $\delta(x)$ is the class of the $\mu_2$-torsor $\Spin_Q\times_{\PSpin_Q}\Spec(\ZZ)$, the fibre over $x$.
  \end{samepage}
  Since $\mu_2$ is central, multiplication in $\Spin_Q$ makes $\delta$ a homomorphism, so the image $\PSL_2(\Lambda)$ of $\Spin_Q(\ZZ)$ is the kernel of $\delta$.
	We will now compute $\vert H^1(\Spec(\ZZ)_{\fppf},\mu_2)\vert =2$, which shows that $\PSL_2(\Lambda)$ has index at most $2$ in $\PSpin_Q(\ZZ)$. The index is exactly $2$: $u=x(y+z)\in\GSpin_Q(\ZZ)$ has $\ns(u)=-1$, so its class in $\PSpin_Q(\ZZ)=\GSpin_Q(\ZZ)/\lbrace\pm1\rbrace$ (Corollary~\ref{cor:pspin-gspin-quotient}) is not in the image of $\Spin_Q(\ZZ)$.
	The group $\PSpin_Q(\ZZ)$ maps to $\SO_Q(\ZZ)$ by $\overline{\pi}$; it contains $\psi(\PSL_2(\Lambda))$ as a subgroup of index $2$, so its image has finite index in $\SO_Q(\ZZ)$, but this image is in general a proper subgroup (Remark~\ref{rem:spin-sequence-at-2}).

	To see this, we use the Kummer sequence for multiplication by 2 in the fppf topology  on $\Spec(\ZZ)$,
	$$ 1 \to \mu_2 \to \GG_m \xrightarrow{[2]} \GG_m \to 1.$$
	The long exact cohomology sequence of a short exact sequence (\cite[III 2.6.1]{Knus1991}) gives 
	$$ 1\to \mu_2(\ZZ)  \to \GG_m(\ZZ)  \xrightarrow{[2]}\GG_m(\ZZ) \to H^1(\Spec(\ZZ)_{\fppf},\mu_2) \to H^1(\Spec(\ZZ)_{\fppf}, \GG_m). $$
	We have $\mu_2(\ZZ)=\GG_m(\ZZ) = \GG_m(\ZZ)=\lbrace \pm 1\rbrace$, and the image of the map $[2]:\GG_m(\ZZ) \to \GG_m(\ZZ)$ is the map sending $-1$ to $1$. 
	Also, by \cite[III 2.7.4]{Knus1991} we can compare the fppf topology to the \'etale topology to conclude that $H^1(\Spec(\ZZ)_{\fppf},\GG_m) = \Pic(\ZZ)=1$.
	This proves that $\lbrace \pm 1 \rbrace = \GG_m(\ZZ)\to H^1(\Spec(\ZZ)_{\fppf},\mu_2)$ is surjective with trivial kernel and hence an isomorphism.
\end{proof}

\begin{remark}
	It is necessary to move to the fppf topology here as the sequences are not exact in the \'etale topology. 
	See for example \cite[III.2.7.3 part (3)]{Knus1991}.
 Write $\ZZ_{(p)}^{\sh}$ for the strict Henselization of the localization of $\ZZ$ at the prime $p$.
 Above primes $p\neq 2$ and the generic point, the morphism 
$\GG_m(\ZZ_{(p)}^{\sh})  \xrightarrow{[2]}\GG_m(\ZZ_{(p)}^{\sh})$ is surjective,
so exactness of the Kummer and Spin sequences can be checked on $\Spec(\ZZ)_{\et}$
		by \cite[Proposition 1.41]{Auel2009}.
		We cannot use this result for $p=2$ because the morphism $\GG_m(\ZZ_{(2)}^{\sh}) \to \GG_m(\ZZ_{(2)}^{\sh})$ given by $t\mapsto t^2$ is not surjective, and hence the $2$-Kummer sequence is not exact in the \'etale topology. 
		Knowing this, one might naively try to show that the morphism is exact on the fppf stalks as can be done for proofs using the \'etale topology away from~$2$.
		This is extremely difficult, however, as the stalks in the fppf topology are not easy to describe: see
		\cite[\href{https://mathoverflow.net/questions/42258/what-are-the-tau-local-rings-for-a-subcanonical-grothendieck-topology-tau-on-the-cat}{MO42258}]{mathoverflow-grothendieck-topology}.	
\end{remark}

\subsection{Borel--Harish-Chandra}\label{sec:borel-harish-chandra}

\begin{definition}
	Let $H$ be a Lie group. 
	A \emph{lattice} $\Gamma<H$ is a subgroup which is discrete and has finite covolume with respect to the Haar measure on $H$.
\end{definition}

\begin{theorem}[Borel--Harish-Chandra]\label{thm:borel-harish-chandra}
	Let $H$ be a semisimple Lie group. 
	If $\Gamma < H$ is an arithmetic group, then $\Gamma$ is a lattice.
\end{theorem}
\begin{proof}
	See \cite[Major Theorem 5.1.11]{Morris2015}.
See also \cite[Theorem 4.14]{Platonov1994} (cf.\ \cite[Theorem 4.13]{Platonov1994}).
\end{proof}

\begin{corollary}
	If $\Gamma<H$ is an arithmetic group, then it acts properly discontinuously by left multiplication on $H$. 
	It will also act properly discontinuously on $H/K$ for $K$ a maximal compact subgroup. 
\end{corollary}
\begin{proof}
  See \cite[p. 96]{Boothby1986}. Only the discreteness of $\Gamma$ is used here. For $\Gamma=\SL_2(\Ocal)$ acting on $\Hcal^{n+1}\cong\SL_2(\CC_n)/\SU_2(\CC_n)$ discreteness is elementary and does not need Theorem~\ref{thm:arithmeticity}: $\SL_2(\Ocal)\subset M_2(\Ocal)$, and $\Ocal$ is a lattice in $\CC_n$. 
  See also \cite[\href{https://math.stackexchange.com/questions/3186183/left-translation-on-lie-group-of-a-discrete-subgroup-is-properly-discontinuous}{MSE186183}]{math-stack-exchange-left-translation}.
\end{proof}

\section{Fundamental Domains for $\PSL_2(\Ocal)$}\label{sec:fundamental-domains}

Let $\Ocal$ be a Clifford-stable order (Definition~\ref{def:clifford-stable}) in $\quat{-d_1,\ldots,-d_{n-1}}{\QQ} \subset \CC_n$ for $d_j \in \ZZ$ squarefree coprime integers, for example $\Ocal = \quat{-d_1,\ldots,-d_{n-1}}{\ZZ}$.
Let $\Gamma = \PSL_2(\Ocal)$. 
In this section we give a construction of the open fundamental domain $D \subset \Hcal^{n+1}$ for $\Gamma$.
Interestingly, the proof is a ``proof by moduli interpretation''. 
Here, we first show that every point in $\Hcal^{n+1}$ is in $\SL_2(\CC_n)$-bijection with a homothety class of positive definite Clifford-Hermitian matrices.
We then classify the elements $x$ of a $\Gamma$-orbit in $\Hcal^{n+1}$ with maximal $x_n$ as those which correspond to homothety classes of Hermitian forms achieving a minimal value at a certain standard point. 

This mirrors the theory of quadratic forms in $\SL_n(\RR)$ and the theory of Hermitian forms in $\SL_2(\CC)$ but with some stranger definitions.
This section further develops Bianchi-Humbert Theory and is inspired by section 3 of \cite{Swan1971}.

\subsection{Special $\Ocal$ and Stabilizers of $\infty$}
At times we will need to talk about stabilizers of $\infty$ and the Clifford group of $\Ocal$. 
In this subsection we only use that $\Ocal$ is $*$-stable, i.e., that
$\Ocal^* = \Ocal$; the Clifford-stability assumed at the start of this section
is used through Lemma~\ref{L:proper-inputs} and Lemma~\ref{L:ideals}\eqref{I:vertical}.

\begin{definition}
  Recall that $\Ocal^\times$, the group of Clifford units of $\Ocal$, is
  defined to be $U(\Ocal) \cap \CC_n^\times$.
  We say that $\Ocal$ is \emph{special} if $\Ocal^\times \neq \lbrace \pm 1\rbrace$. In this case we will also say
  that $\Gamma = \PSL_2(\Ocal)$ and $\Gamma\backslash\Hcal^{n+1}$ are special.
\end{definition}

\begin{lemma}\label{lem:stabilizer-of-infinity}
	The stabilizer group $\SL_2(\Ocal)_{\infty}$ is generated by matrices of the form
	$ \tau_s = \begin{psmallmatrix} 1 & s \\
	0 & 1 
	\end{psmallmatrix}$, and $\sigma_t = \begin{psmallmatrix} t & 0 \\ 0 & t^{*-1} \end{psmallmatrix},$ for  $s \in \Vec(\Ocal)$, and $t \in \Ocal^{\times}.$
In fact $\SL_2(\Ocal)_{\infty} \cong \Vec(\Ocal) \rtimes \Ocal^{\times}$, where $t\in\Ocal^{\times}$ acts on $\Vec(\Ocal)$ by $s\mapsto tst^*$. 
\end{lemma}
\begin{proof}
	We know that they must be upper triangular of the form $g=\begin{psmallmatrix}a & b \\ 0 & d
	\end{psmallmatrix} \in \SL_2(\Ocal)$.
	The condition $\Delta(g) = ad^* = 1$ implies that $d^*=a^{-1}$ or $d=a^{*-1}$.
	We can factor these matrices as $\begin{psmallmatrix}
	a & b \\
	0 & a^{*-1}
	\end{psmallmatrix} = \begin{psmallmatrix} a & 0 \\
	0 & (a^*)^{-1} \end{psmallmatrix} \begin{psmallmatrix} 1 & a^{-1}b \\ 0 & 1 \end{psmallmatrix},$
	so the matrices of the form 
	$ \tau_s = \begin{psmallmatrix} 1 & s \\
	0 & 1 
	\end{psmallmatrix}$ and $\sigma_t = \begin{psmallmatrix} t & 0 \\ 0 & t^{*-1} \end{psmallmatrix}$
	generate $\SL_2(\Ocal)_{\infty}$. Here $a \in \Ocal^{\times}$, since $a^{-1}=d^*\in \Ocal$ by $*$-stability, and $a^{-1}b \in \Vec(\Ocal)$, since $\sigma_a^{-1}g = \begin{psmallmatrix} 1 & a^{-1}b \\ 0 & 1 \end{psmallmatrix}$ lies in $\SL_2(\Ocal)$. 
\end{proof}

\begin{corollary}
	The group $\Gamma = \PSL_2(\Ocal)$ is special if and only if $\Gamma_{\infty}$ is larger than $\Lambda = \Vec(\Ocal)$. 
\end{corollary}
\begin{proof}
By Lemma~\ref{lem:stabilizer-of-infinity}, $\Gamma_\infty$ is the image of $\Vec(\Ocal)\rtimes\Ocal^\times$ in $\PSL_2(\Ocal)$. Since $\tau_s\neq\pm1$ for $s\neq0$ and $\sigma_{-1}=-1$, the subgroup $\Lambda$ embeds in $\Gamma_\infty$ and $\Gamma_\infty/\Lambda\cong\Ocal^\times/\{\pm1\}$. Hence $\Gamma_\infty\neq\Lambda$ if and only if $\Ocal^\times\neq\{\pm1\}$.
\end{proof}

\subsection{Clifford-Hermitian Matrices}\label{sec:hermitian}
Recall the \emph{Clifford adjoint} $A^{\dag} = (\overline{A})^{\t}$ of a matrix $A\in M_{m,n}(\CC_n)$, defined before Lemma~\ref{L:closed-under-conj}.

\begin{definition}
	We say $A \in M_2(\CC_n)$ is \emph{Clifford-Hermitian} (or simply \emph{Hermitian} for convenience) if and only if it has Clifford vector entries and $A = A^{\dag}$. We denote the collection of $2\times 2$ Clifford-Hermitian matrices by $M_2(\CC_n)_{\herm}$. 
\end{definition}
Note that if $A$ is Hermitian, then it has the form
$$ A = \begin{pmatrix} a & b\\ \overline{b} & c \end{pmatrix}, \qquad a,c \in \RR, \quad b\in V_n.$$
(In \S\ref{sec:naive} the off-diagonal entries were written $\overline{y}$ and $y$; this is the same form with $b=\overline{y}$.)
We now define the discriminant, positive definiteness, and homothety (cf.  \cite[Section 3]{Elstrodt1988}).
\begin{definition}
	Let $A = \begin{psmallmatrix} a & b  \\ \overline{b} & c \end{psmallmatrix}$ be Hermitian.
	\begin{enumerate}
		\item  Its \emph{discriminant} (or \emph{naive determinant}) is defined to be the real number $\det(A) = ac-\vert b\vert^2.$
		\item We say $A$ is \emph{positive definite} if and only if $a,c>0$ and $\det(A)>0$.
		We will denote the cone of positive definite Hermitian matrices by $M_2(\CC_n)_{\herm}^{\pos}$.\footnote{
			Note that this is indeed a cone as the Cauchy-Schwarz inequality shows that if $A,B \in M_2(\CC_n)_{\herm}^{\pos}$ then $A+B \in M_2(\CC_n)_{\herm}^{\pos}$.}
		\item Two positive definite matrices $A$ and $B$ are \emph{homothetic} if and only if there exists some $k \in \RR_{>0}$ such that $A=kB$. 
	\end{enumerate} 
\end{definition}
Note that in the definition of positive definite it suffices to take either $a>0$ or $c>0$. 
Also note that the notation is such that if $A = \begin{psmallmatrix} a & b  \\ \overline{b} & c \end{psmallmatrix} \in M_2(\CC_n)_{\herm}^{\pos}$ then $b\in V_n$. 
In the future we will let $\overline{M_2(\CC_n)}_{\herm}^{\pos} = M_2(\CC_n)_{\herm}^{\pos}/\RR_{>0}^{\times}$ be the collection of positive definite Hermitian matrices up to homothety.

\begin{theorem}
	There is a well-defined left action of $\SL_2(\CC_n)$ on $M_2(\CC_n)_{\herm}^{\pos}$ given by $g\cdot A = gAg^{\dag}$. 
	This action respects homothety in the sense that $A$ and $B$ are homothetic if and only if $g\cdot A$ and $g\cdot B$ are for each $g \in \SL_2(\CC_n)$.
\end{theorem}
\begin{proof}
	This is \cite[Prop 3.3]{Elstrodt1988}.
\end{proof}

\subsection{Clifford-Hermitian Matrices and the Action of $\SL_2$}\label{sec:action-hermitian}
In section 3 of \cite{Elstrodt1988}, Elstrodt, Grunewald, and Mennicke give a $\SL_2(\CC_n)$-equivariant bijection between homothety classes of positive definite matrices and $\Hcal^{n+1}$. 
We will make use of this bijection and develop a further relationship which will help us give a formula for the boundary of the fundamental domain of $\Gamma$.
\begin{theorem}
	The map $\Phi:\overline{M_2(\CC_n)}_{\herm}^{\pos} \to \Hcal^{n+1}$ given by 
	$$\begin{pmatrix}
	a & b \\
	\overline{b} & c
	\end{pmatrix} \mapsto \dfrac{b}{c} + i_n \dfrac{\sqrt{ac-\vert b\vert^2}}{c}$$ is an $\SL_2(\CC_n)$-equivariant bijection.
	The inverse is given by  
	$$\Phi^{-1}(z+i_n \zeta) = \begin{pmatrix}
	c(\zeta^2 + \vert z \vert^2) & cz \\
	c \overline{z} & c
	\end{pmatrix}$$
	for some choice of $c>0$. 
	Here we are writing $x \in \Hcal^{n+1}$ as $x=z+i_n \zeta$ where $\zeta \in \RR_{>0}$ and $z \in V_n$. 
\end{theorem}
\begin{proof}
	This is \cite[Prop.~3.4]{Elstrodt1988} and \cite[Prop~3.5]{Elstrodt1988}. 
	They show that there exists an inverse map $\Psi: \Hcal^{n+1} \to M_2(\CC_n)_{\herm}^{\pos}$ which is the ``Hermitian form associated to $x$'' that we will see
        in Theorem~\ref{T:characterization-of-bubble-domain}, (\ref{I:proper-minimal-value}).
    The second part is implicit in the proof of Theorem~\ref{thm:hermitian-form-shape} which we will delay.
\end{proof}

We will need some enhancements to this map. 

\begin{definition}\label{def:matrix-to-form}
	Let $A \in M_2(\CC_n)_{\herm}^{\pos}$.  The \emph{Hermitian form} associated to $A$ is the map $q_A: \CC_n^2 \to \CC_n$ given by 
	\begin{equation}\label{E:quad-form}
	q_A(w) = w^{\dag} A w, \qquad w = \begin{pmatrix}u \\ v \end{pmatrix}.
	\end{equation}
	By Theorem~\ref{thm:hermitian-form-shape}, $q_A(w)\in\RR$ whenever $w$ satisfies \eqref{E:the-condition}.
\end{definition}

Under certain circumstances we can complete the square. 
This gives a description of the map $\Phi$ from \cite{Elstrodt1988} in terms of these parameters.
\begin{theorem}\label{thm:hermitian-form-shape}
	Let $A =\begin{psmallmatrix} a & b \\ \bar{b} & c \end{psmallmatrix}$ be a positive definite Hermitian matrix. 
	Let $w = \begin{pmatrix} u \\ v \end{pmatrix}$ and suppose that $u,v \in \CC_n^\mon \cup \lbrace 0 \rbrace$ are such that there exists a matrix 
	\begin{equation}\label{E:the-condition}
	\begin{pmatrix}
	u &  *\\
	v & *
	\end{pmatrix} \in \GL_2(\CC_n).
	\end{equation}
	\begin{enumerate}
		\item We may write the Hermitian form as
		\begin{equation}\label{E:expanded-form}
		q_A(w) = a \vert u \vert^2 + 2\Re(\overline{u} b v) + c \vert v \vert^2.
		\end{equation}
		\item With $w$ as above, if we let $z=b/c\in V_n$ ($c$ is real) and $\zeta=\sqrt{ac-\vert b\vert^2}/c$, then 
		$$ q_A(w) = c \left ( \vert \zeta u \vert^2 + \vert \bar{z}u+v \vert^2 \right), \qquad w = \begin{pmatrix} u \\ v \end{pmatrix}.$$ 
	\end{enumerate}
\end{theorem}
\begin{proof}
	\begin{enumerate}
		\item 
		Here we just have 
		\begin{align*}
		w^{\dag} A w &= \bar{u}(au+bv) + \bar{v}(\bar{b}u+cv)\\
		&= a \vert u \vert^2 + \bar{u}bv + \bar{v}\bar{b}u+c \vert v \vert^2 \\
		&=a \vert u \vert^2 + 2 \Re(\bar{u}bv) + c \vert v \vert^2.
		\end{align*}
		We used that $\vert \alpha \vert^2 = \alpha \overline{\alpha} = \overline{\alpha} \alpha$ for $\alpha \in \CC_n^\mon$, and that $\alpha+\overline{\alpha} = 2\Re(\alpha)$ for $\alpha=\bar{u}bv$. The last fact is clear if $u=0$. If $u\neq 0$, then $p=vu^{-1}\in V_n$ by \eqref{E:the-condition} (see the proof of (2)), so $\alpha+\overline{\alpha} = \bar{u}(bp+\bar{p}\bar{b})u = (bp+\bar{p}\bar{b})\vert u\vert^2$, because $bp+\bar{p}\bar{b}$ is real for $b,p\in V_n$. Hence $\alpha+\overline{\alpha}$ is real, so it equals $2\Re(\alpha)$.
\item We will show that $\bar{z}u+v \in \CC_n^\mon$, so that we may expand $\vert \bar{z}u+v\vert^2$ in terms of conjugates. 
If $u=0$, then $\bar{z}u+v=v\in\CC_n^\mon\cup\lbrace 0\rbrace$ and there is nothing to prove, so let $u\neq0$. Then $\bar{z}u+v=(\bar{z}+vu^{-1})u$. 
If $v=0$, then $vu^{-1}=0\in V_n$; otherwise the condition \eqref{E:the-condition} implies that $uv^{-1} \in V_n$ and hence that $vu^{-1}\in V_n$ since multiplicative inverses of Clifford vectors are Clifford vectors. This proves that $\bar{z}+vu^{-1}\in V_n$.
Since the product of a Clifford vector and a Clifford group element is a Clifford group element, we have shown that $\vert \bar{z}u+v\vert^2 = \overline{(\bar{z}u+v)}(\bar{z}u+v)=(\bar{z}u+v)\overline{(\bar{z}u+v)}$.

		We get 
		\begin{align*}
		\zeta^2\vert u\vert^2 + \vert \bar{z}u+v\vert^2 =& \zeta^2 \vert u\vert^2 + \vert z\vert^2 \vert u \vert^2 + \overline{\bar{z}u}v + \overline{v}\bar{z}u + \vert v \vert^2\\
		&=(\zeta^2 + \vert z \vert^2)\vert u \vert^2+2\Re(\bar{u}zv)+\vert v\vert^2.
		\end{align*}
		Now, multiplying the above expression by $c$ and matching terms with \eqref{E:expanded-form} we find $a/c=(\zeta^2+\vert z\vert^2)$ and $b/c=z$ give a solution.
		This gives $z=b/c$ and $\zeta^2=(ac-\vert b \vert^2)/c^2$ or $\zeta = \sqrt{ac-\vert b \vert^2}/c$.
	\end{enumerate}	
\end{proof}

From now on we will allow ourselves to conflate $\overline{M_2(\CC_n)}_{\herm}^{\pos}$ and $\Hcal^{n+1}$ using this $\PSL_2(\CC_n)$-equivariant bijection.

\subsection{Unimodular Pairs}
We remind the reader that we call $ab^{-1}$ left division and $b^{-1}a$ right division, following conventions in the literature (cf. Remark~\ref{rem:left-vs-right-division}).

We need an equivariance statement here that describes the relation between the
matrix of the same transformation in two different sets of coordinates. This is
analogous to the classical statement that an inner product space whose inner product
is given by the matrix $M$ in a fixed basis has the matrix $A^TMA$ when the basis
is changed by the action of $A$.  

We record this in the form of a lemma for future reference. 
\begin{lemma}\label{lem:qa-is-g-equivariant}
	For all $g\in \SL_2(\CC_n)$, all $A\in M_2(\CC_n)_{\herm}^{\pos}$, and all $w\in\CC_n^2$, we have  
	$q_A(g w)=q_{g^{\dag}Ag}(w)$. 
\end{lemma}
\begin{proof}
For $w^{\dag}$, a unimodular row vector, we have $q_A(w) = w^\dag A w$.
Then acting on the left of $w$ gives $q_A(g w) = (gw)^{\dag} A (gw) = w^{\dag} g^{\dag} A g w = q_{g^{\dag} A g}(w)$.
\end{proof}

\subsection{Cusps and Ideal Classes}\label{sec:cusps}
Consider $\Clf(q)$ for $q$ a positive definite quadratic form in $n-1$
variables. 
Let $K=\Clf(q)\otimes\QQ$. 
Let $\Ocal$ be an order in $K$. 
Let $\widehat{V} = \Vec(K) \cup \lbrace \infty \rbrace$. 
Let $\Gamma = \PSL_2(\Ocal)$.
The aim of this section is to prove that $\Gamma\backslash\widehat{V}$ is finite using ideal classes. 

The terminology of ideals is borrowed from \cite{Reiner1975}.
We build on some terminology from \S\ref{sec:maximal-orders}.
By a left \emph{$\Ocal$-lattice} $L$ we mean a finitely generated left $\Ocal$-module that is free as a $\ZZ$-module, as in \cite{Reiner1975}. 
When $L$ is a subset of a free left $K$-module $M$, we call it \emph{full} if $\mathbb{Q}L = M$. 
By a \emph{fractional ideal} in $K$ we mean a full left $\Ocal$-lattice in $M=K$. 
Since an isomorphism of fractional ideals induces an automorphism of $K$ as a left $K$-module, we have that $I\cong J$ if and only if there exists some $a \in U(K)$ such that $I=Ja$.

\begin{definition}
	Given two fractional ideals $I, J \subset K$, we declare them to be equivalent if and only if there exists some unit $a \in U(K)$ such that $I=Ja$. 
	We will write $I\sim J$ when this is the case. 
	The set of equivalence classes will be denoted by $\Cl(\Ocal)$ and called the \emph{left ideal class set}.
\end{definition}

We use the finiteness of the set of ideal classes. 
\begin{theorem}[Jordan-Zassenhaus {\cite[Theorem 26.4]{Reiner1975}}]
	If $\Ocal$ is an order
	in a semisimple $\mathbb{Q}$-algebra $K$, then there are only finitely many isomorphism classes of left $\Ocal$-lattices of any given $\mathbb{Z}$-rank. 
	As a consequence, $\Cl(\Ocal)$ is a finite set. 
\end{theorem}

This applies to our situation, as $q$ is nondegenerate so $\Clf(q)_{\QQ}=K$ is
semisimple. (If the arity $n-1$ of the quadratic form is even it is central simple.)

We record the following Lemma relating the notion of equivalence to generation of ideals.
\begin{lemma}\label{L:ideals}
Suppose that $\Ocal$ is closed under the Clifford conjugations and that $(c,d)\in \Ocal^2$ is unimodular. Then 
	\begin{enumerate}
		\item $\Ocal c^*+\Ocal d^*=\Ocal$ as left $\Ocal$-ideals;
		\item $c\Ocal^*+d\Ocal^*=\Ocal^*$ as right $\Ocal^*$-ideals;
		\item \label{I:vertical} $\Ocal\overline{c} + \Ocal \overline{d} = \Ocal$ as left $\Ocal$-ideals;
		\item $c\Ocal + d\Ocal = \Ocal$ as right $\Ocal$-ideals.
	\end{enumerate}
\end{lemma}
\begin{proof}
	Suppose that $\gamma=\begin{psmallmatrix} a& b\\c&d\end{psmallmatrix} \in \SL_2(\Ocal)$.
	Then $ad^*-bc^*=1$, which implies that $-b(c^*)+ad^* = 1$, which implies that $\Ocal = \Ocal c^*+\Ocal d^*$. 
	This can also seen by multiplying $\gamma^{-1}$ on the left by $\gamma$. 
	The second part follows from applying the $*$ operation.
	The third part follows from the existence of a left matrix inverse for $\begin{psmallmatrix} \bar{c} & * \\
	\bar{d} & * \end{psmallmatrix}$. 
	The fourth part follows from the existence of a right matrix inverse for a matrix of the form $\begin{psmallmatrix}
	c & d \\
	* & *
	\end{psmallmatrix}$.
\end{proof}

For $x\in\widehat{V}$, define $\chi_x$ and $L_x$ by
\begin{align*}
  \chi_x&=\begin{psmallmatrix}x\\1\end{psmallmatrix}, & L_x&=\Ocal x+\Ocal &&(x\in\Vec(K)),\\
  \chi_\infty&=\begin{psmallmatrix}1\\0\end{psmallmatrix}, & L_\infty&=\Ocal.
\end{align*}
Thus $L_x$ is the left ideal generated by the entries of $\chi_x$.
For a fractional ideal $L$, define its right dual by
\[
  L^{\vee}=\lbrace z\in K\colon Lz\subseteq\Ocal\rbrace.
\]
This is a right $\Ocal$-module, with $LL^{\vee}\subseteq\Ocal$. Moreover,
\begin{align*}
  (Lc)^{\vee}&=c^{-1}L^{\vee} &&(c\in U(K)),\\
  L_x^{\vee}&=\lbrace z\in\Ocal\colon xz\in\Ocal\rbrace &&(x\in\Vec(K)).
\end{align*}

\begin{definition}\label{def:cuspidal-class-set}
	Two fractional ideals $I, J\subset K$ are \emph{Clifford equivalent} if $I=Jc$ for some $c$ in the Clifford group $K^{\times}$ (Definition~\ref{defn:CliffordGroup}).
	The \emph{cuspidal class set} $\Cl^{\mathrm{cusp}}(\Ocal)$ is the set of Clifford equivalence classes of the ideals $L_x$ with $x\in\widehat{V}$.
\end{definition}
Figure~\ref{fig:cusps-ideal-classes} shows the map $x\mapsto[L_x]$ in the classical case $\Ocal=\ZZ[\sqrt{-5}]$.
\begin{figure}[htbp!]
\centering
\begin{tikzpicture}[scale=0.85,lab]
 \def\s{2.236}
 \foreach \a in {-2,...,2}{\foreach \b in {-1,0,1}{\fill (\a,\b*\s) circle (2.2pt);}}
 \foreach \a in {-2,...,1}{\foreach \b in {-1,0}{\fill[cuspC] (\a+0.5,\b*\s+\s/2) circle (2.4pt);}}
 \foreach \a in {-2,...,2}{\foreach \b in {-1,0}{\draw (\a,\b*\s+\s/2) circle (2pt);}}
 \draw[gray,->] (-2.7,0)--(2.7,0) node[right]{$\RR$}; \draw[gray,->] (0,-2.7)--(0,2.9) node[above]{$\sqrt{-5}\,\RR$};
 \node[cuspC,above right] at (0.5,\s/2){$\tfrac{1+\sqrt{-5}}2$};
 \node[right=8pt] at (3.6,0.9){$\infty$ and black, open points $\longmapsto [\Ocal]$};
 \node[right=8pt,cuspC] at (3.6,0){red points $\longmapsto [(2,1+\sqrt{-5})]$};
 \node[right=8pt] at (3.6,-0.9){$\Cl(\Ocal)=\{[\Ocal],[(2,1+\sqrt{-5})]\}$};
\end{tikzpicture}
\caption{Cusps and ideal classes (\S\ref{sec:cusps}) in the commutative case $\Ocal=\ZZ[\sqrt{-5}]$, $n=2$. A cusp $x$ goes to the class of $L_x=\Ocal x+\Ocal$. Black: points of $\Ocal$; open: points $\tfrac{\sqrt{-5}}2+\Ocal$, where $L_x$ is again principal; both lie in the orbit of $\infty$. Red: points $\tfrac{1+\sqrt{-5}}2+\Ocal$, where $L_x$ is the non-principal class. So $\PSL_2(\Ocal)$ has two cusps.}\label{fig:cusps-ideal-classes}
\end{figure}

Since $K^{\times}\subseteq U(K)$, there is a natural map $\Cl^{\mathrm{cusp}}(\Ocal)\to\Cl(\Ocal)$.
Equivalently, $\Cl^{\mathrm{cusp}}(\Ocal)$ is the set of Clifford equivalence classes of the left ideals $\Ocal\alpha+\Ocal\beta$, where $\begin{psmallmatrix}\alpha\\ \beta\end{psmallmatrix}$ runs over the first columns of elements of $\SL_2(K)$: such a column is $\chi_x u$ with $u\in K^{\times}$, where $x=\alpha\beta^{-1}$ and $u=\beta$ if $\beta\neq0$, and $x=\infty$ and $u=\alpha$ if $\beta=0$; conversely each $\chi_xu$ is such a column by Lemma~\ref{lem:completable} below.

\begin{definition}\label{def:completable-cusp}
	A cusp $x\in\Vec(K)$ is \emph{completable} (for $\Ocal$) if there is $t\in\Vec(\Ocal)$ such that $tx\in\Ocal$ and $x-xtx\in\Ocal$.
	The cusp $\infty$ is completable.
\end{definition}
For example, $x$ is completable if $x\in\Vec(\Ocal)$ (take $t=0$) or if $x^{-1}\in\Vec(\Ocal)$ (take $t=x^{-1}$).

\begin{lemma}\label{lem:completable}
	Let $\Ocal$ be closed under the Clifford conjugations.
	\begin{enumerate}
		\item \label{I:At} For $x,t\in\Vec(K)$, the matrix $A_t$ and its inverse are
		\[
		  A_t=\begin{pmatrix}x&xt-1\\1&t\end{pmatrix}\in\SL_2(K),
		  \qquad
		  A_t^{-1}=\begin{pmatrix}t&1-tx\\-1&x\end{pmatrix}.
		\]
		Every element of $\SL_2(K)$ with first column $\chi_x$ is of this form.
		\item \label{I:completing} A cusp $x\in\widehat{V}$ is completable if and only if some $A\in\SL_2(K)$ has first column $\chi_x$ and the first row of $A^{-1}$ has all its entries in $L_x^{\vee}$.
		\item \label{I:cocycle} For $\gamma\in\SL_2(\Ocal)$ and $x\in\widehat{V}$ we have $\gamma\chi_x=\chi_{\gamma x}u$ for some $u\in K^{\times}$, and then $L_x=L_{\gamma x}u$.
		\item \label{I:invariant} If $x$ is completable and $\gamma\in\SL_2(\Ocal)$, then $\gamma x$ is completable.
	\end{enumerate}
\end{lemma}
\begin{proof}
	\eqref{I:At} The entries of $A_t$ lie in the Clifford monoid: if $x\neq0$ then $xt-1=x(t-x^{-1})$ is a product of two vectors.
	We have $\Delta(A_t)=xt^*-(xt-1)=1$, $x(xt-1)^*=xtx-x\in\Vec(K)$ and $t1^*=t\in\Vec(K)$, so $A_t\in\SL_2(K)$ by Theorem~\ref{thm:gl2-short}, which also gives the inverse.
	If $\begin{psmallmatrix}x&b\\1&d\end{psmallmatrix}\in\SL_2(K)$, then $d=d1^*\in\Vec(K)$ and $1=xd^*-b$, so the matrix is $A_d$.

	\eqref{I:completing} For $x=\infty$ take $A=1$. For $x\in\Vec(K)$, by \eqref{I:At} we need $t\in\Vec(K)$ with $t,1-tx\in L_x^{\vee}$.
	Since $x=x^*$ and $t=t^*$, we have $xt\in\Ocal$ if and only if $tx=(xt)^*\in\Ocal$.
	So $t\in L_x^{\vee}$ if and only if $t\in\Vec(\Ocal)$ and $tx\in\Ocal$, and in that case $1-tx\in L_x^{\vee}$ if and only if $x-xtx\in\Ocal$.

	\eqref{I:cocycle} Let $\gamma=\begin{psmallmatrix}a&b\\c&d\end{psmallmatrix}$, so that $ad^*-bc^*=1$ and $dc^*\in\Vec(K)$.
	If $c\neq0$ then $c^{-1}d=c^{-1}(dc^*)(c^{-1})^*\in\Vec(K)$, so $(c^{-1}d)^*=c^{-1}d$, that is, $c^{-1}dc^*=d^*$.
	Let $x\in\Vec(K)$. If $cx+d\neq0$, then $\gamma\chi_x=\chi_{\gamma x}u$ with $u=cx+d$, which is $d$ if $c=0$ and $c(x+c^{-1}d)$ if $c\neq0$; in both cases $u$ is a nonzero element of the Clifford monoid, so $u\in K^{\times}$ by Proposition~\ref{prop:anisotropic}.
	If $cx+d=0$, then $c\neq0$, $\gamma x=\infty$ and $\gamma\chi_x=\chi_\infty u$ with $u=ax+b=b-ac^{-1}d=(bc^*-ad^*)(c^*)^{-1}=-(c^*)^{-1}$.
	Finally $\gamma\chi_\infty=\chi_{\gamma\infty}u$ with $u=c$ if $c\neq0$ and $u=a$ if $c=0$.
	Since $\gamma$ and $\gamma^{-1}=\begin{psmallmatrix}d^*&-b^*\\-c^*&a^*\end{psmallmatrix}$ have entries in $\Ocal$, the entries of $\gamma\chi_x$ generate the same left ideal $L_x$ as the entries of $\chi_x$; and the entries of $\chi_{\gamma x}u$ generate $L_{\gamma x}u$.

	\eqref{I:invariant} Choose $A\in\SL_2(K)$ with first column $\chi_x$ and first-row entries of $A^{-1}$ in $L_x^\vee$. Choose $u\in K^\times$ with $\gamma\chi_x=\chi_{\gamma x}u$, and put
	\[
	  A'=\gamma A\,\mathrm{diag}(u^{-1},u^*)\in\SL_2(K).
	\]
	Then $A'$ has first column $\gamma\chi_xu^{-1}=\chi_{\gamma x}$, and the first row of $A'^{-1}=\mathrm{diag}(u,(u^*)^{-1})A^{-1}\gamma^{-1}$ is $u$ times the first row of $A^{-1}\gamma^{-1}$.
	Its entries lie in $uL_x^{\vee}\Ocal\subseteq uL_x^{\vee}=(L_xu^{-1})^{\vee}=L_{\gamma x}^{\vee}$, so $\gamma x$ is completable by \eqref{I:completing}.
\end{proof}

\begin{theorem}\label{thm:finitely-many-cusp-classes}
  Let $K=\quat{-d_1,\ldots,-d_{n-1}}{\QQ}$. Suppose that $\Ocal$ is an order in $K$ closed under the Clifford conjugations.
  \begin{enumerate}
  \item \label{I:cusp-map} The rule $x\mapsto L_x$ induces a well-defined surjective map $\Gamma \backslash \widehat{V} \to \Cl^{\mathrm{cusp}}(\Ocal)$, and hence a well-defined map $\Gamma \backslash \widehat{V} \to \Cl(\Ocal)$.
  \item \label{I:cusp-inj} The completable cusps form a union of $\Gamma$-orbits that contains $\infty$, and two completable cusps $x$ and $y$ are $\Gamma$-equivalent if and only if $L_x$ and $L_y$ are Clifford equivalent.
  \item \label{I:cusp-bij} If every cusp is completable, then $x\mapsto L_x$ induces a bijection $\Gamma \backslash \widehat{V} \to \Cl^{\mathrm{cusp}}(\Ocal)$.
  This is the case if $n=2$ and $\Ocal$ is the ring of integers of the imaginary quadratic field $K$, and then $\Cl^{\mathrm{cusp}}(\Ocal)=\Cl(\Ocal)$ is the ideal class group of $K$.
  \end{enumerate}
\end{theorem}
\begin{proof}
  \eqref{I:cusp-map} By Lemma~\ref{lem:completable}\eqref{I:cocycle}, $L_{\gamma x}$ is Clifford equivalent to $L_x$ for all $\gamma\in\SL_2(\Ocal)$ and $x\in\widehat{V}$. The map is surjective by the definition of $\Cl^{\mathrm{cusp}}(\Ocal)$.

  \eqref{I:cusp-inj} The first claim is Lemma~\ref{lem:completable}\eqref{I:invariant}, and one direction of the second is \eqref{I:cusp-map}.
  Conversely, let $x$ and $y$ be completable with $L_x=L_yc$ for some $c\in K^{\times}$.
  By Lemma~\ref{lem:completable}\eqref{I:completing} there are $A,B\in\SL_2(K)$ with first columns $\chi_x$ and $\chi_y$ such that the entries of the first rows of $A^{-1}$ and $B^{-1}$ lie in $L_x^{\vee}$ and $L_y^{\vee}$, respectively.
  Put $B'=B\,\mathrm{diag}(c,(c^*)^{-1})$ and $\sigma=B'A^{-1}\in\SL_2(K)$.
  Then $\sigma\chi_x=B'\begin{psmallmatrix}1\\0\end{psmallmatrix}=\chi_yc$, so $\sigma(x)=y$.
  The entries of the first column $\chi_yc$ of $B'$ lie in $L_yc=L_x$.
  If $\begin{psmallmatrix}b'\\d'\end{psmallmatrix}$ is the second column of $B'$, then the first row of $B'^{-1}$ is $(d'^*,-b'^*)$ by Theorem~\ref{thm:gl2-short}; it is $c^{-1}$ times the first row of $B^{-1}$, so $b'^*,d'^*\in c^{-1}L_y^{\vee}=L_x^{\vee}$.
  Write $\chi_x=\begin{psmallmatrix}\alpha\\ \beta\end{psmallmatrix}$, so that $L_x=\Ocal\alpha+\Ocal\beta$ and the second row of $A^{-1}$ is $(-\beta^*,\alpha^*)$, and let $r$ be the first row of $A^{-1}$. Then
  $$\sigma=B'A^{-1}=(\chi_yc)\,r+\begin{psmallmatrix}b'\\d'\end{psmallmatrix}(-\beta^*,\ \alpha^*).$$
  The entries of $(\chi_yc)\,r$ lie in $L_xL_x^{\vee}\subseteq\Ocal$.
  The entries of the second matrix are $-(\beta e^*)^*$ and $(\alpha e^*)^*$ with $e\in\lbrace b',d'\rbrace$, and these lie in $\Ocal$ because $\alpha e^*,\beta e^*\in L_xL_x^{\vee}\subseteq\Ocal$ and $\Ocal^*=\Ocal$.
  Hence $\sigma\in\SL_2(K)\cap M_2(\Ocal)=\SL_2(\Ocal)$, and $x$ and $y$ are $\Gamma$-equivalent.

  \eqref{I:cusp-bij} The first claim follows from \eqref{I:cusp-map} and \eqref{I:cusp-inj}.
  Let $n=2$ and $\Ocal=\Ocal_K$. Then $\Vec(K)=K$, $\Vec(\Ocal)=\Ocal$ and $K^{\times}=K\setminus\lbrace0\rbrace$.
  By the proof of Lemma~\ref{lem:completable}\eqref{I:completing}, $x\in K$ is completable if and only if there is $t\in L_x^{\vee}$ with $1-tx\in L_x^{\vee}$, that is, if and only if $1\in L_x^{\vee}x+L_x^{\vee}=L_x^{\vee}L_x$.
  This holds because $L_x^{\vee}$ is the inverse of the invertible ideal $L_x$ of the Dedekind domain $\Ocal_K$.
  Finally, every fractional ideal of $\Ocal_K$ is $\Ocal\alpha+\Ocal\beta$ for some $\alpha$ and some $\beta\neq0$, and then it equals $L_{\alpha\beta^{-1}}\beta$; so $\Cl^{\mathrm{cusp}}(\Ocal_K)=\Cl(\Ocal_K)$.
\end{proof}

\begin{corollary}\label{cor:class-set-finite}
  For any $*$-stable order $\Ocal$, even if not closed under the other Clifford conjugations,
  the set $\Gamma \backslash \widehat{V}$ is finite.
\end{corollary}

\begin{proof}
  Suppose first that $\Ocal$ is closed under the Clifford conjugations.
  The group $\SL_2(K)$ acts transitively on $\widehat{V}$, and the stabilizer of $\infty$ is its subgroup $P$ of upper triangular matrices.
  So $\Gamma \backslash \widehat{V}$ is in bijection with $\SL_2(\Ocal)\backslash\SL_2(K)/P$, which is finite by \cite[Prop.~6.2]{Elstrodt1988}. In general,
  let $\Ocal_c = \Ocal \cap \Ocal^* \cap \bar \Ocal \cap \Ocal'$. Then
  $\Ocal_c$ is an order closed under the Clifford conjugations, and so
  there are only finitely many equivalence classes of cusps for $\Ocal_c$.
  Since $\Ocal_c$-equivalent cusps are $\Ocal$-equivalent, the more general
  statement follows.
\end{proof}

\begin{definition}\label{def:completable-local}
	Let $\Ocal$ be an order in $K$ and $p$ a prime. Write $\ZZ_{(p)}$ for the localization of $\ZZ$ at $p$ and set $\Ocal_{(p)}=\ZZ_{(p)}\Ocal$.
	A cusp $x\in\Vec(K)$ is \emph{completable at $p$} if there is $t\in\Vec(\Ocal_{(p)})$ such that
	\[
	  tx\in\Ocal_{(p)}\qquad\text{and}\qquad x-xtx\in\Ocal_{(p)}.
	\]
	The cusp $\infty$ is completable at $p$.
\end{definition}

\begin{lemma}\label{lem:completable-local}
Let $\Ocal$ be an order in $K$ closed under the Clifford conjugations. Let 
	$$\Vec(\Ocal)^{\#}=\lbrace u\in\Vec(K)\colon u\bar v+v\bar u\in\ZZ\text{ for all }v\in\Vec(\Ocal)\rbrace,$$ and let $p^{e_p}$ be the exponent of the $p$-primary part of the finite group $\Vec(\Ocal)^{\#}/\Vec(\Ocal)$.
	\begin{enumerate}
		\item \label{I:local} A cusp $x\in\Vec(K)$ is completable if and only if it is completable at every prime $p$. It is completable at every $p$ with $x\in\Vec(\Ocal_{(p)})$.
		\item \label{I:finite-check} Suppose that for every prime $p$ with $e_p>0$, every $x\in p^{-e_p}\Vec(\Ocal)$ is completable. Then every cusp is completable.
	\end{enumerate}
	Since translations by $\Vec(\Ocal)$ lie in $\Gamma$, the hypothesis of \eqref{I:finite-check} only has to be checked for the finitely many classes of $p^{-e_p}\Vec(\Ocal)$ modulo $\Vec(\Ocal)$.
\end{lemma}
\begin{proof}
	\eqref{I:local} Let $\Lambda=\lbrace t\in\Vec(\Ocal)\colon tx\in\Ocal\rbrace$. Then $x$ is completable if and only if $x$ lies in the lattice $M=x\Lambda x+\Ocal$, and it is completable at $p$ if and only if $x\in\ZZ_{(p)}M$, since localization commutes with finite sums and intersections of lattices and with intersection with $\Vec(K)$. An element of $K$ lies in $M$ if and only if it lies in $\ZZ_{(p)}M$ for every $p$. If $x\in\Vec(\Ocal_{(p)})$, then $t=0$ works.

	\eqref{I:finite-check} The proof of Lemma~\ref{lem:completable} applies to $\Ocal_{(p)}$ and $\SL_2(\Ocal_{(p)})=\SL_2(K)\cap M_2(\Ocal_{(p)})$, so completability at $p$ is invariant under $\SL_2(\Ocal_{(p)})$.
	Fix $x\in\Vec(K)$ and a prime $p$, put $e=e_p$, and let $k\geq0$ be minimal with $y=p^kx\in\Vec(\Ocal_{(p)})$. We claim that $x$ is $\SL_2(\Ocal_{(p)})$-equivalent to a cusp $x'$ with $p^ex'\in\Vec(\Ocal_{(p)})$. This is clear if $k\leq e$, so let $k>e$.
	For $u,v\in\Vec(\Ocal_{(p)})$ the elements $\nrd(u)=u\bar u$ and $\langle u,v\rangle=u\bar v+v\bar u$ lie in $\ZZ_{(p)}$, and $\nrd(u+v)=\nrd(u)+\langle u,v\rangle+\nrd(v)$.
	Since $p^{-1}y\notin\Vec(\Ocal_{(p)})$ and $p^e$ kills $\ZZ_{(p)}\Vec(\Ocal)^{\#}/\Vec(\Ocal_{(p)})$, we have $p^{-e-1}y\notin\ZZ_{(p)}\Vec(\Ocal)^{\#}$, so there is $w\in\Vec(\Ocal)$ with $\operatorname{ord}_p\langle y,w\rangle\leq e$.
	If $\operatorname{ord}_p\nrd(y)\leq k+e$, replace $w$ by $0$.
	In both cases $y'=y+p^kw$ satisfies $\operatorname{ord}_p\nrd(y')\leq k+e$, because $\operatorname{ord}_p(p^{2k}\nrd(w))\geq 2k>k+e$.
	Now $z=x+w=p^{-k}y'$ is not in $\Vec(\Ocal_{(p)})$, so $z\neq0$, and the element $S\begin{psmallmatrix}1&w\\0&1\end{psmallmatrix}$ of $\Gamma$, where $S=\begin{psmallmatrix}0&-1\\1&0\end{psmallmatrix}$, sends $x$ to $x'=-z^{-1}=-p^k\nrd(y')^{-1}\bar{y}'$.
	As $\bar y'\in\Vec(\Ocal_{(p)})$ and $\operatorname{ord}_p(p^k\nrd(y')^{-1})\geq-e$, we get $p^ex'\in\Vec(\Ocal_{(p)})$, proving the claim.
	Translating by an element of $\Vec(\Ocal_{(p)})$, we may assume $x'\in p^{-e}\Vec(\Ocal)$. If $e=0$ then $x'\in\Vec(\Ocal)$ is completable, and if $e>0$ it is completable by hypothesis. By \eqref{I:local}, $x'$ is completable at $p$, hence so is $x$. As $p$ was arbitrary, $x$ is completable by \eqref{I:local}.
\end{proof}

\begin{remark}\label{rem:cusp-map-not-injective}
	We explain the relation between cusps and ideal classes for the orders of this paper.
	
	\begin{enumerate}
	\item \emph{An invariant.} For $c\in K^{\times}$, right multiplication by $c$ on the $2^{n-1}$-dimensional $\QQ$-vector space $K$ is conjugate under $x\mapsto x^*$ to left multiplication by $c^*$, so by Proposition~\ref{prop:anisotropic} (with $W$ of rank $n-1$) its determinant has absolute value $\nrd(c^*)^{2^{n-2}}=(c\bar c)^{2^{n-2}}$.
	Hence the index $[L_x:\Ocal]$ modulo $2^{n-2}$-th powers of positive rational numbers depends only on the class of $L_x$ in $\Cl^{\mathrm{cusp}}(\Ocal)$, and so only on the $\Gamma$-orbit of $x$.
	In particular, if $x$ is equivalent to $\infty$ then $[L_x:\Ocal]$ is a $2^{n-2}$-th power.

	\item \emph{The map to $\Cl(\Ocal)$.} The map $\Gamma\backslash\widehat{V}\to\Cl(\Ocal)$ need not be injective, even for maximal orders.
	Let $\Ocal=\Ocal_{5,!}$ be the maximal order of \S\ref{sec:c5-oddball} and $s=(1+i_1+i_2+i_3)/2$.
	Then $[L_s:\Ocal]=2^4$ is not an $8$th power, so $s$ is not equivalent to $\infty$.
	But $L_s=\Ocal\lambda$ with $\lambda=\frac14(-3i_1-i_1i_2+i_3-i_2i_3+i_4+i_2i_4-i_1i_3i_4-i_1i_2i_3i_4)\in U(K)$, so $L_s$ and $L_\infty=\Ocal$ have the same class in $\Cl(\Ocal)$.
	The map need not be surjective either: by a \texttt{magma} computation the maximal order $B(-1,-1,-3)_0$ has two left ideal classes but only one cusp.
	This is why we work with $\Cl^{\mathrm{cusp}}(\Ocal)$ rather than $\Cl(\Ocal)$.

	\item \emph{Maximal orders.} For every maximal order $\Ocal$ in Table~\ref{tab:orders}, every cusp is completable, so Theorem~\ref{thm:finitely-many-cusp-classes} gives a bijection between $\Gamma\backslash\widehat{V}$ and $\Cl^{\mathrm{cusp}}(\Ocal)$.
	We checked this with Lemma~\ref{lem:completable-local}.
	For these orders $\Vec(\Ocal)^{\#}/\Vec(\Ocal)$ has exponent $2$ (the Hurwitz order, $\Ocal_4$ and the orders in $\quat{-1,-1,-1,-1}{\QQ}$) or $6$ (the orders in $\quat{-1,-3}{\QQ}$ and $\quat{-1,-1,-3}{\QQ}$), so $e_2,e_3\leq1$, and it suffices to check the cusps in $\frac12\Vec(\Ocal)$ and $\frac13\Vec(\Ocal)$ modulo $\Vec(\Ocal)$, at most $2^4+3^4$ of them for each order; a computer check shows that all of them are completable.
	For instance, $t=-1-i_4$ shows that $s=(1+i_1+i_2+i_3)/2$ is completable for $\Ocal_{5,!}$.
	For the orders in Table~\ref{tab:orders} with one cusp, it follows that every ideal $L_x$ is generated by an element of the Clifford group.
	For $\Ocal_{5,!}$ the classes of $L_\infty$ and $L_s$ in $\Cl^{\mathrm{cusp}}(\Ocal)$ are distinct, and so are those of $L_\infty$ and $L_{s'}$ for $A(-1,-1,-3)$ and $s'=(1+a_3)/2$, since $[L_{s'}:\Ocal]=4$ is not a $4$th power; so both orders have at least two cusps.

	\item \emph{Non-maximal orders.} For the order $\ZZ[i_1,i_2,i_3]$ of \S\ref{sec:1-1-1-order} the cusp $\zeta=(1+i_1+i_2+i_3)/2$ is not completable, and for $\ZZ[\sqrt{-3}]$ the cusp $(1+\sqrt{-3})/2$ of Example~\ref{ex:sl2-z3} is not completable.
	We do not know whether every cusp is completable for every maximal order closed under the Clifford conjugations.
	\end{enumerate}

	The finiteness of the set of cusps is a consequence of the reduction theory of arithmetic groups \cite[Proposition 15.6]{Borel2019}; see also \cite[Prop.~6.2]{Elstrodt1988}.
	
\end{remark}

\subsection{Characterization of Spheres}
See Figure~\ref{fig:characterization-of-spheres}.
\begin{figure}[htbp!]
\centering
\begin{tikzpicture}[scale=2.4,lab]
 \fill[bubbleC!8] (-0.2,0) rectangle (3.4,1.6);
 \fill[white] (2.4,0) arc[start angle=0,end angle=180,radius=1];
 \fill[cuspC!10] (2.4,0) arc[start angle=0,end angle=180,radius=1];
 \draw[bub] (2.4,0) arc[start angle=0,end angle=180,radius=1];
 \draw[bnd] (-0.2,0)--(3.4,0);
 \dotat[holeC]{1.4,0} \node[below] at (1.4,0){$-\bar u^{-1}\bar v$};
 \draw[covC,-{Stealth}] (1.4,0)--++(125:1) node[midway,right]{$1/|\bar u|$};
 \fill (0.3,1.2) circle (1.4pt) node[right]{$x$: $q_A(u,v)\ge 1$};
 \fill (1.75,0.3) circle (1.4pt) node[below=1pt,font=\scriptsize]{$x'$: $q_A(u,v)<1$};
 \draw[dotted] (0.3,1.2)--(0.3,0) node[below]{$z$}; \node[left] at (0.3,0.6){$\zeta$};
 \node at (3.1,1.4){$\Hcal^{n+1}$};
\end{tikzpicture}
\caption{Theorem~\ref{T:characterization-of-spheres}: for $x=z+i_n\zeta\in\Hcal^{n+1}$ with Hermitian matrix $A$, the condition $q_A(u,v)\ge q_A(0,1)=1$ holds exactly when $x$ lies outside or on the hemisphere (blue) of radius $1/|\bar u|$ centred at $-\bar u^{-1}\bar v\in V_n$. Schematic; vertical section through the centre.}\label{fig:characterization-of-spheres}
\end{figure}
\begin{theorem}\label{T:characterization-of-spheres}
	Let $u,v\in \CC_n^\mon$ with $u\neq 0$ be such that there exists a matrix
	$\begin{psmallmatrix}
	u &  *\\
	v & *
	\end{psmallmatrix} \in \GL_2(\CC_n)$.
	Let $x\in \Hcal^{n+1}$ and let $[A]\in \overline{M_2(\CC_n)}_{\herm}^{\pos}$ be its corresponding homothety class of positive definite matrices.
	
The condition $q_A(u,v)\geq q_A(0,1)$ on the matrix $A$ is equivalent to the condition that $x$ lies outside or on $B(-\bar{u}^{-1}\bar{v})$ where $-\bar{u}^{-1}\bar{v}\in \Vec(K)=\partial \Hcal^{n+1}$.
The radius of  $B(-\bar{u}^{-1}\bar{v})$ is  $1/\vert \bar{u} \vert$.
\end{theorem}
\begin{proof}
	$x = z+i_n\zeta\in \Hcal^{n+1}$ where $\zeta\in \RR_{>0}$ and $z\in V_n$.
	We write $q_A(u,v)$ as 
	$$q_A(u,v) = \vert \zeta u \vert^2 + \vert \bar{z}u+v \vert^2.$$
Here $A$ is normalized in its homothety class so that its lower right entry is $1$ (Theorem~\ref{thm:hermitian-form-shape}(2)); in particular, setting $(u,v)=(0,1)$ gives $q_A(0,1)=1$, and the condition on the matrix becomes $q_A(u,v)\geq 1$.
\begin{align*}
	\vert \zeta u \vert^2 + \vert \bar{z}u+v \vert^2\geq 1 &\iff \vert \zeta\vert^2 \vert u \vert^2 + \vert \bar{z} + vu^{-1} \vert^2 \vert u \vert^2 \geq 1\\
	&\iff \vert \bar{x} + vu^{-1} \vert  \geq 1/\vert u \vert \\
	&\iff \vert x + \bar{u}^{-1}\bar{v}\vert  \geq 1/\vert \bar{u} \vert 
	\end{align*}
	The second to last inequality is equivalent to $\vert x - (-\overline{vu^{-1}}) \vert \geq 1/\vert u \vert$. 
	Since $vu^{-1}$ is a Clifford group element we can use the formula $y^{-1} = \overline{y}/\vert y \vert^2$ or $\bar{y} = \vert y \vert^2 y^{-1}$. This implies $\vert y \vert^2 y^{-1} = \overline{y}$, and hence that $\overline{(vu^{-1})} = \vert vu^{-1} \vert^2 (vu^{-1})^{-1}=\vert v\vert^2 \vert u \vert^{-2} uv^{-1}=\vert v \vert^2 \bar{u}^{-1} \bar{v} \vert v \vert^{-2} = \bar{u}^{-1}\bar{v}$.
\end{proof}

\subsection{Bubbles $B_{(c,d)}$ and the Bubble Domain $B$}\label{sec:bubble-domain}
In virtue of the above definition we make the following definition. 

\begin{definition}
	Let $(\lambda,\mu)\in \Ocal^2$ with $\mu\neq 0$ be such that $(\mu,-\lambda)$ is right unimodular (Definition~\ref{def:right-unimodular}). 
	The \emph{bubble} $B_{(\lambda,\mu)}$ (or just $B(\mu^{-1}\lambda)$) in $\Hcal^{n+1}$ at $(\lambda,\mu)$ is a sphere of radius $1/\vert \mu \vert$ centered at $\mu^{-1}\lambda\in \Vec(K)$.
	Its equation is given by 
	\begin{equation}\label{E:above-sphere}
	\vert x-\mu^{-1}\lambda \vert = 1/\vert \mu \vert.
	\end{equation}
\end{definition}

The closed set $B\subset \Hcal^{n+1}$, consisting of elements no lower than any integral boundary sphere, plays an important role in defining our fundamental domain.

\begin{definition}\label{def:bubble-domain}
	The \emph{bubble domain} is the set $B\subset \Hcal^{n+1}$ defined by 
	\begin{equation}
	B=\lbrace x \in \Hcal^{n+1} \colon \forall (\lambda,\mu), \vert x -\mu^{-1}\lambda\vert \geq 1/\vert\mu\vert \rbrace.
	\end{equation}
	where $(\lambda,\mu)$ runs over the pairs of the previous definition.
\end{definition}
Figure~\ref{fig:bubble-domain} shows bubbles and the bubble domain for $\ZZ$ and in a slice for $\ZZ[\omega]$, $\omega=(1+\sqrt{-19})/2$.
\begin{figure}[htbp!]
\centering
\begin{tikzpicture}[scale=2.4,lab]
 \fill[bubbleC!8] (-1.25,0) rectangle (1.25,1.3);
 \begin{scope}\clip (-1.25,0) rectangle (1.25,1.3);
 \foreach \c/\r in {-1/1,0/1,1/1,-0.5/0.5,0.5/0.5,-0.333/0.333,0.333/0.333,-0.667/0.333,0.667/0.333,-0.25/0.25,0.25/0.25,-0.75/0.25,0.75/0.25}{
   \draw[bub] (\c+\r,0) arc[start angle=0,end angle=180,radius=\r];}
 \end{scope}
 \begin{scope}\clip (-1.25,0) rectangle (1.25,1.3);
 \foreach \c in {-1,0,1}{\fill[white,fill opacity=0.0] (\c,0) circle (1);}\end{scope}
 \draw[bnd] (-1.25,0)--(1.25,0);
 \foreach \c in {-1,0,1}{\dotat{\c,0}}
 \foreach \c in {-0.5,0.5,-0.333,0.333,-0.667,0.667,-0.25,0.25,-0.75,0.75}{\fill (\c,0) circle (0.9pt);}
 \node at (0,1.15) {$B$};
 \node[below] at (0,0){$0$}; \node[below] at (1,0){$1$}; \node[below] at (0.5,0){$\tfrac12$};
 \node at (0,-0.35){(a) $\Ocal=\ZZ$: bubbles with $|\mu|\le 4$};
\end{tikzpicture}\qquad
\begin{tikzpicture}[scale=1.7,lab]
\fill[bubbleC!8] (-0.3,0) rectangle (2.5,1.25);
 \draw[bub] (0.968,0) arc[start angle=0,end angle=180,radius=0.968];
 \draw[bub] (2.179+0.968,0) arc[start angle=0,end angle=180,radius=0.968];
 \draw[bub,line width=1.4pt] (1.0897+0.5,0) arc[start angle=0,end angle=180,radius=0.5];
 \draw[bnd] (-0.3,0)--(2.5,0);
 \node[bubbleC] at (0,1.05){$B(0)$}; \node[bubbleC] at (2.2,1.05){$B(\omega)$};
 \node[bubbleC,above] at (1.09,0.5){$B(\omega/2)$};
 \node[below] at (0,0){$0$}; \node[below] at (1.09,0){$\tfrac{\sqrt{19}}4$};\node[below] at (2.179,0){$\tfrac{\sqrt{19}}2$};
 \node at (1.1,-0.42){(b) $\ZZ[\omega]$, slice $x_0=\tfrac14$, axis $x_1$};
 \begin{scope}[shift={(0,0)}]\clip (-0.3,0) rectangle (2.5,1.3);\end{scope}
\end{tikzpicture}
\caption{Bubbles $B(\mu^{-1}\lambda)$ (blue) have radius $1/|\mu|$; the bubble domain $B$ (light blue) lies above all of them. (a) $\Ocal=\ZZ$, all bubbles with $|\mu|\le 4$: every bubble other than those at the integers lies under a unit bubble. (b) $\Ocal=\ZZ[\omega]$, $\omega=(1+\sqrt{-19})/2$, the vertical slice $x_0=1/4$ of $\Hcal^3$: the bubble $B(\omega/2)$ of radius $1/2$ reaches above $B(0)$ and $B(\omega)$, so it bounds $B$. Only the bubbles named in \S\ref{S:sqrt19} are drawn in (b).}\label{fig:bubble-domain}
\end{figure}
\begin{definition} \label{defn:fundamentalDomain} Let $E$ be a contractible space and let $\Gamma$ be a group
	acting on $E$ with discrete orbits. An {\em open fundamental domain}
	for the action of $\Gamma$ on $E$ is an open subset $D$ such that
        $E = \cup_{\gamma \in \Gamma} \gamma \bar D$
	and if $\gamma d \in D$ for $\gamma \in \Gamma$ and $d \in D$ then $\gamma = e_\Gamma$.
        We refer to the $\gamma D$ as {\em tiles}.
\end{definition}
Let $F$ be an open fundamental domain for the action of $\Gamma_{\infty}$ on $V_n=\partial \Hcal^{n+1}\setminus\lbrace\infty\rbrace$, chosen as in Theorem~\ref{T:fundamental-domain}.  
Then the fundamental domain for $\Gamma\subset \PSL_2(\CC_n)$ is the subset of the interior $B^{\circ}$ of $B$ which projects down to $F\subset V_n$ under the map $x_0+\cdots + x_{n-1}i_{n-1}+x_ni_n \mapsto x_0+\cdots + x_{n-1}i_{n-1}$. 
In a formula, the fundamental domain $D$ is given by
\begin{equation}\label{E:fundamental-domain}
D=\lbrace x_0+x_1i_1+\cdots + x_ni_n \in B^{\circ} \colon x_0 + x_1 i_1+ \cdots + x_{n-1}i_{n-1} \in F \rbrace.
\end{equation}

The idea, then, is to bring everything to the region $B$ by maximizing the $x_n$ value by inverting outside of these bubbles, and then moving these points back to the region that sits above the fundamental domain in the boundary lattice using translations by lattice elements.

We now give a characterization of being an element of $B$; this will give us a ``reduction theory'' to $B$.
To do this we introduce a definition.
\begin{definition}
A \emph{unimodular value} of a Clifford-Hermitian form $q_A$ is a value $q_A(\bar{\mu},\bar{\nu})$ with $(\mu,\nu)$ right unimodular (Definition~\ref{def:right-unimodular}).
\end{definition}
We now characterize the bubble domain $B$ as the collection of Hermitian forms which attain their minimum unimodular value at $(0,1)$. 
\begin{theorem}[Characterization of the Bubble Domain]\label{T:characterization-of-bubble-domain}
	The following are equivalent conditions on a point $x \in \Hcal^{n+1}$.
	\begin{enumerate}
		\item \label{I:in-B} $x\in B$.
		\item \label{I:proper-minimal-value} The Hermitian form associated to $x\in \Hcal^{n+1}$ has a minimal unimodular value at $(0,1)$. 
		\item \label{I:heights}  For all $\gamma \in \Gamma$ we have $x_n\geq (\gamma(x))_n.$
		\item \label{I:small-absolute-value} For all $\gamma\in \Gamma$ we have $\vert \gamma'(x)\vert \leq 1.$
		Here, $\gamma'(x)$ is the Jacobian derivative of the transformation $x
		\mapsto (ax+b)(cx+d)^{-1}.$
	\end{enumerate}
	Moreover, $x$ lies in the interior $B^{\circ}$ of $B$ if and only if $x_n> (\gamma(x))_n$ for all $\gamma\in\Gamma\setminus\Gamma_{\infty}$.
\end{theorem}
\begin{proof}
	\begin{itemize}
		\item We will first show that \eqref{I:heights} is equivalent to \eqref{I:in-B}.
		Let $g = \begin{psmallmatrix} a & b \\ c & d\end{psmallmatrix} \in \SL_2(\Ocal)$ be arbitrary.
		From the Magic Formula (Theorem~\ref{thm:magic}) we know that $(gx)_n = \vert cx+d\vert^{-2} x_n. $
		This holds unconditionally. 
		The condition about maximality of the $n$th component is equivalent to 
		$\vert cx+d \vert\geq 1$ for all $\begin{psmallmatrix} a & b \\ c & d\end{psmallmatrix} \in \SL_2(\Ocal)$.
If $c=0$, then $d$ is a unit, so $\vert d \vert = 1$ and the inequality holds with equality; here we use integrality: $\vert d\vert^2$ and $\vert d^{-1}\vert^2$ lie in $\Ocal\cap\RR_{>0}\subseteq\ZZ_{>0}$ and have product $1$. If $c\neq 0$, then it reads $1/\vert c\vert \leq \vert x+c^{-1}d\vert$. 
		By Lemma~\ref{L:proper-inputs} we are allowed to vary over matrices of a different shape, and these turn out to be the integral boundary spheres, hence equivalent to $x\in B$.
		
		\item To see that \eqref{I:proper-minimal-value} is equivalent to \eqref{I:in-B} we only need to use Theorem~\ref{T:characterization-of-spheres} and again Lemma~\ref{L:proper-inputs} to change the spheres coming from the minimal proper value condition being at $(0,1)$ to the general shape.
		
		\item We will show that \eqref{I:heights} and \eqref{I:small-absolute-value} are equivalent. 
		By the Magic Formula, again we have 
		$\gamma'(x): V_n \to V_n$ defined by $ \gamma'(x)u = \Delta(\gamma)(xc^*+d^*)^{-1}u(cx+d)^{-1}, $
		whenever $x$ and $\gamma(x)$ are not $\infty$.
		Since $\Delta(\gamma)=1$ we have $\vert \gamma'(x) \vert =\vert cx+d \vert^{-2}$. 
		The condition $\vert \gamma'(x)\vert \leq 1$ then again translates to $1 \leq \vert cx+d\vert^2$. 
\item For the last assertion, note first that only finitely many bubbles meet a compact set $C\subset\Hcal^{n+1}$. If $C\subset\lbrace x\colon x_n\geq\varepsilon,\ \vert x\vert\leq R\rbrace$ meets $B_{(\lambda,\mu)}$, then $1\leq\vert\mu\vert\leq1/\varepsilon$ (the lower bound because $\vert\mu\vert^2\in\Ocal\cap\RR_{>0}\subseteq\ZZ_{>0}$) and $\vert\lambda\vert=\vert\mu\vert\,\vert\mu^{-1}\lambda\vert\leq(R+1)/\varepsilon$, and $\Ocal$ is a lattice. If $x_n>(\gamma(x))_n$ for all $\gamma\notin\Gamma_\infty$, then by the first item $x$ lies strictly outside every bubble. Hence so does a neighbourhood of $x$, and $x\in B^{\circ}$.
		If $x_n=(\gamma(x))_n$ for some $\gamma\notin\Gamma_\infty$, then $x$ lies on a bubble, and points inside that bubble arbitrarily close to $x$ are not in $B$.
	\end{itemize}
\end{proof}

\subsection{Boundedness of Siegel Heights}

We now give a proof of boundedness of heights in orbits.
\begin{theorem}[Bounded Heights]\label{T:siegel-heights}
	For every $ x\in \Hcal^{n+1}$ we have $\sup\lbrace \gamma(x)_n \colon \gamma\in \Gamma \rbrace<\infty$, and the supremum is achieved.
\end{theorem}
\begin{proof}
	
	If $\gamma(x) = (ax+b)(cx+d)^{-1}$ then $\gamma(x)_n = x_n \vert cx+d\vert^{-2}$.
	This allows us to translate statements about $\sup \gamma(x)_n$ into
        statements about $\inf \vert cx+d \vert$.
        In particular, we need to show that there is no infinite decreasing sequence
        $\vert c_jx+d_j \vert$ with 
	$\begin{psmallmatrix}
	a_j & b_j \\ c_j & d_j
	\end{psmallmatrix} \in \SL_2(\Ocal).$
	This follows from the fact that $\Ocal \subset \CC_n$ is a lattice. Indeed, for fixed $x$ the map $(c,d)\mapsto cx+d$ from $\CC_n^2$ to $\CC_{n+1}$ is $\RR$-linear and injective, since $\CC_{n+1}=\CC_n\oplus\CC_n i_n$ and $x_n>0$. So $(c,d)\mapsto\vert cx+d\vert^2$ is a positive definite quadratic form on $\CC_n^2$, and it takes only finitely many values below any bound on the lattice $\Ocal^2$. 
	Another way to think about this is that there is no infinite decreasing sequence $\vert c_nx+d_n\vert$ given $(c_n,d_n) \in \Ocal^2$; hence, there is no infinite decreasing sequence with the additional condition that there exists some $a_n,b_n \in \Ocal$ such that  $\begin{psmallmatrix}
	a_n & b_n \\ c_n & d_n
	\end{psmallmatrix} \in \SL_2(\Ocal).$
\end{proof}

\subsection{The Fundamental Domain}\label{sec:fundamental-domain}

Figure~\ref{fig:fundamental-domain} shows $D$ for $\Ocal=\ZZ$ and for $\Ocal=\ZZ[i]$.
\begin{figure}[htbp!]
\centering
\begin{tikzpicture}[scale=1.5,lab]
 \begin{scope}\clip (-1.6,-0.05) rectangle (1.6,2.2);
 \fill[bubbleC!8] (-1.6,0) rectangle (1.6,2.2);
 \foreach \c in {-1,0,1}{\fill[white] (\c+1,0) arc[start angle=0,end angle=180,radius=1];}
 \foreach \c in {-2,-1,0,1,2}{\draw[bub] (\c+1,0) arc[start angle=0,end angle=180,radius=1];}
 \fill[domC!35] (-0.5,2.2)--(-0.5,0.866) arc[start angle=120,end angle=60,radius=1]--(0.5,2.2)--cycle;
 \draw[domC!80!black,line width=0.8pt] (-0.5,2.2)--(-0.5,0.866) (0.5,0.866)--(0.5,2.2);
 \end{scope}
 \draw[bnd] (-1.6,0)--(1.6,0);
 \draw[domC!80!black,line width=2.5pt] (-0.5,0)--(0.5,0) node[below,midway,black]{$F$};
 \node at (0,1.5) {$D$}; \node at (1.25,1.8) {$B$};
 \node[bubbleC] at (0,0.55) {$B(0)$};
 \foreach \c in {-1,0,1}{\dotat{\c,0}}
 \node at (0,-0.5){(a) $n=1$, $\Ocal=\ZZ$};
\end{tikzpicture}\qquad
\begin{tikzpicture}[scale=1.9,lab,x={(1cm,0cm)},y={(0.5cm,0.35cm)},z={(0cm,1cm)}]
 \draw[bnd] (-0.9,-0.8,0)--(0.9,-0.8,0)--(0.9,0.8,0)--(-0.9,0.8,0)--cycle;
 \draw[dom,fill=domC!45] (-0.5,-0.5,0)--(0.5,-0.5,0)--(0.5,0.5,0)--(-0.5,0.5,0)--cycle;
 \node at (0,0,0) {$F$};
\draw[bub,domain=-0.5:0.5,samples=30] plot (\x,-0.5,{sqrt(0.75-\x*\x)});
 \draw[bub,domain=-0.5:0.5,samples=30] plot (\x,0.5,{sqrt(0.75-\x*\x)});
 \draw[bub,domain=-0.5:0.5,samples=30] plot (0.5,\x,{sqrt(0.75-\x*\x)});
 \draw[bub,domain=-0.5:0.5,samples=30] plot (-0.5,\x,{sqrt(0.75-\x*\x)});
 \draw[bub,dotted,domain=-0.5:0.5,samples=30] plot (\x,0,{sqrt(1-\x*\x)});
 \draw[bub,dotted,domain=-0.5:0.5,samples=30] plot (0,\x,{sqrt(1-\x*\x)});
 \foreach \a/\b in {-0.5/-0.5,0.5/-0.5,0.5/0.5,-0.5/0.5}{
   \draw[domC!80!black,dashed] (\a,\b,0)--(\a,\b,0.707);
   \draw[domC!80!black,line width=0.8pt] (\a,\b,0.707)--(\a,\b,2.0);}
 \node at (0,0,1.55) {$D$};
 \node[bubbleC,right] at (0.5,0.5,0.75) {floor on $B(0)$};
 \node at (0,-0.9,-0.25) {(b) $n=2$, $\Ocal=\ZZ[i]$ (schematic)};
\end{tikzpicture}
\caption{The fundamental domain of Theorem~\ref{T:fundamental-domain}. (a) For $\Ocal=\ZZ$ acting on $\Hcal^2$: the bubble domain $B$ (light blue) is the region above all bubbles $B(\mu^{-1}\lambda)$; here the visible ones are the unit bubbles at the integers. The fundamental domain $D$ (orange) is the part of $B^\circ$ above $F=[-1/2,1/2]$. (b) For $\Ocal=\ZZ[i]$ acting on $\Hcal^3$: $D$ is the column over the square $F\subset V_2$ cut off below by $B(0)$; the walls continue to $\infty$. Schematic; (b) is drawn in oblique perspective.}\label{fig:fundamental-domain}
\end{figure}
\begin{theorem}\label{T:fundamental-domain}
	Let $F$ be an open fundamental domain for $\Gamma_{\infty}$ acting on $V_n$ which is the interior of a compact convex polytope (for example a Dirichlet domain for $\Gamma_\infty$ \cite[\S 6.6]{Ratcliffe2019}). 
	The open set $D \subset \Hcal^{n+1}$ given by 
	$$D = \lbrace x_0+x_1i_1+\cdots + x_ni_n \in B^{\circ} \colon x_0 + x_1 i_1+ \cdots + x_{n-1}i_{n-1} \in F \rbrace$$ 
	is an open fundamental domain for $\Gamma = \PSL_2(\Ocal)$ acting on $\Hcal^{n+1}$, and $\overline{D}=\lbrace x\in B\colon x_0+\cdots+x_{n-1}i_{n-1}\in\overline{F}\rbrace$.
	
\end{theorem}
\begin{proof}
	\begin{itemize}
\item Tiling Property: Write $\pi(x)=x_0+\cdots+x_{n-1}i_{n-1}$ and $h(x)=x_n$. By Lemma~\ref{lem:stabilizer-of-infinity}, $\Gamma_{\infty}$ acts by $x\mapsto x+s$ with $s\in\Vec(\Ocal)$ and by $x\mapsto txt^*$ with $t\in\Ocal^\times$. Here $\vert t\vert=1$: the diagonal entries $t$ and $t^{*-1}$ of $\sigma_t$ satisfy $\vert t\vert^2\,\vert t^{*-1}\vert^2=1$, and both factors lie in $\Ocal\cap\QQ_{>0}=\ZZ_{>0}$. These maps preserve $h$, are compatible with $\pi$, and preserve $B$, since $B$ is described by Theorem~\ref{T:characterization-of-bubble-domain}\eqref{I:heights} in terms of $\Gamma$-orbits.
		Let $x \in \Hcal^{n+1}$. 
		By Theorem~\ref{T:siegel-heights} there is $\gamma\in\Gamma$ such that $h(\gamma(x))$ is maximal on the orbit of $x$, so $\gamma(x)\in B$ by Theorem~\ref{T:characterization-of-bubble-domain}\eqref{I:heights}. 
		Since $\Gamma_\infty\overline{F}=V_n$, there is $t\in\Gamma_\infty$ with $t\gamma(x)\in B\cap\pi^{-1}(\overline{F})$.
		\item Closure and Connectedness: $B$ is an intersection of closed hyperbolic half-spaces (the exteriors of the bubbles) and has non-empty interior, so $B$ is convex and $B=\overline{B^{\circ}}$; also $\pi^{-1}(F)$ is open and convex. Hence $D$ is open and convex, in particular path connected. It is non-empty, since it contains the points of $\pi^{-1}(F)$ of height $>1$. So $\overline{D}=B\cap\pi^{-1}(\overline{F})$, and by the previous item $\Gamma\overline{D}=\Hcal^{n+1}$ and $\Gamma_\infty\overline{D}=B$.
		\item Finite Fixed Points in the Quotient: we will show that for all $x\in \Hcal^{n+1}$ there are only finitely many $\gamma$ such that $\overline{D} \owns \gamma(x)$.
		Having $\gamma_1(x)$ and $\gamma_2(x) $ in $\overline{D}$ is equivalent to having $x\in \overline{D}$ and some $\gamma(x)\in \overline{D}$.  
		By maximality of the representatives in the interior of $\overline{D}$ we must have that $x \in \partial D$.
		
		We claim that $\orb_{\Gamma}(x) \cap \overline{D}$ for any $x\in \Hcal^{n+1}$ has a height which is bounded above and below. 
		This will imply that 
		$$ \orb_{\Gamma}(x) \cap \overline{D} = \orb_{\Gamma}(x) \cap A$$
		for some compact set $A$ given by $A = \overline{D} \cap \lbrace x \in \Hcal^{n+1} \colon r\leq x_n \leq R \rbrace$ for $r,R\in \RR$. 
		By $\Gamma$ acting properly discontinuously (which only needs the discreteness of $\SL_2(\Ocal)$; see the Corollary to Theorem~\ref{thm:borel-harish-chandra}), we have that there are only finitely many $\gamma\in \Gamma$ such that $A \cap \gamma(A) \neq \emptyset$.
		Since every $\gamma(x) \in \overline{D}$ has the property that $\gamma(x)\in A$, we are done. 
		
		We now prove the bounds. 
		The upper bound on the height comes from boundedness of Siegel heights. 
		The lower bound uses the fact that $\overline{D} \subset \overline{B}$ and $x,\gamma(x) \in \overline{B}$ implies that $x_n=\gamma(x)_n$.
		This is Theorem~\ref{T:characterization-of-bubble-domain} item~\ref{I:heights}.
		
		\item Disjoint Translates of Interior: Suppose $x\in D$ and $\gamma(x)\in D$. Both points lie in $B$, so $h(x)\geq h(\gamma(x))\geq h(\gamma^{-1}\gamma(x))=h(x)$ by Theorem~\ref{T:characterization-of-bubble-domain}\eqref{I:heights}.
		Since $x\in B^{\circ}$, the last assertion of Theorem~\ref{T:characterization-of-bubble-domain} forces $\gamma\in\Gamma_{\infty}$, so that $\gamma(x) = x+v$ for $v\in \Vec(\Ocal)$ or $\gamma(x) = uxu^*$ for some $u \in \Ocal^{\times}$, or a composite of these. 
		Then $\pi(\gamma(x))$ is the image of $\pi(x)$ under $\gamma$, and both lie in $F$, a fundamental domain for $\Gamma_{\infty}$. So $\gamma$ acts trivially on $V_n$, and $\gamma=e_\Gamma$.
	\end{itemize}
\end{proof}

\begin{example}
	In the case of Bianchi groups, the group $\PSL_2(\ZZ[i])$ has extra stabilizers of $\infty$ of the form $\begin{psmallmatrix} i & 0 \\ 0 &-i \end{psmallmatrix}$ which will actually make the fundamental domain smaller. This gives extra automorphisms of the form $x \mapsto ixi$ where $x=x_0+x_1i+x_2j$. Similarly for $\PSL_2(\ZZ[\omega])$, where $\omega^2+\omega+1=0$.
\end{example}

\begin{theorem}\label{T:fundamental-domain-finiteness}
	Let $D=B^\circ\cap\pi^{-1}(F)$ be the open fundamental domain of Theorem~\ref{T:fundamental-domain}.
	Then the set $\Sigma=\lbrace \gamma\in\Gamma \colon \gamma\overline{D}\cap\overline{D}\neq\emptyset\rbrace$ is finite, and the tiles $\gamma\overline{D}$, $\gamma \in \Gamma$, form a locally finite family.
	In particular, the closed fundamental domain $\overline{D}$ intersects finitely many bubbles on its boundary, and $\overline{D}$ is a convex polyhedron with finitely many sides.
\end{theorem}
\begin{proof}
	Write $h(x)=x_n$ and $\pi(x)=x_0+\cdots+x_{n-1}i_{n-1}$.
	For $g=\begin{psmallmatrix}a&b\\c&d\end{psmallmatrix}\in\SL_2(\CC_n)$ we use the Magic Formula (Theorem~\ref{thm:magic}), $h(g(x))=h(x)/\vert cx+d\vert^{2}$.
	The two terms in $cx+d=(c\pi(x)+d)+c\,h(x)i_n$ are orthogonal, so
	\begin{equation}\label{E:cx+d}
	\vert cx+d\vert^2=\vert c\,\pi(x)+d\vert^2+\vert c\vert^2h(x)^2.
	\end{equation}
	If $c\neq0$ is an entry of an element of $\Gamma$, then $\vert c\vert^2=c\bar c\in\Ocal\cap\QQ_{>0}=\ZZ_{>0}$, so $\vert c\vert\geq1$.
	The group $\Gamma\subset\SL_2(\Ocal)$ is discrete and hence acts properly discontinuously on $\Hcal^{n+1}$ \cite[Theorem 5.3.5]{Ratcliffe2019}. \\

	\noindent \emph{Step 1 (points of bounded height).}
	Let $x\in\overline{D}$ and $\gamma\in\Gamma$ with $\gamma(x)\in\overline{D}$.
	Both points lie in $B$, so $h(\gamma(x))=h(x)$ by Theorem~\ref{T:characterization-of-bubble-domain}\eqref{I:heights}.
	If $\gamma\notin\Gamma_\infty$ this means $\vert cx+d\vert=1$, so $h(x)\leq 1/\vert c\vert\leq 1$ by \eqref{E:cx+d}.
For $\varepsilon>0$ the set $A_\varepsilon=\overline{D}\cap\lbrace \varepsilon\leq h\leq 1\rbrace$ is compact, so only finitely many $\gamma$ satisfy $\gamma(A_\varepsilon)\cap A_\varepsilon\neq\emptyset$; and only finitely many $\gamma\in\Gamma_\infty$ satisfy $\gamma(\overline{F})\cap\overline{F}\neq\emptyset$, because $\overline F$ is compact and $\Gamma_\infty$, a discrete group acting on $V_n$ by Euclidean isometries, acts properly discontinuously on $V_n$. 
Hence only finitely many $\gamma$ arise from points $x$ with $h(x)\geq\varepsilon$. \\

	\noindent \emph{Step 2 (reduction theory).}
	By \cite[Prop.~6.3]{Elstrodt1988} there are representatives $\sigma_1=1,\sigma_2,\ldots,\sigma_m\in\SL_2(K)$ of the double cosets $\Gamma\backslash\SL_2(K)/P$, where $P$ is the stabilizer of $\infty$; a compact set $\mathcal{K}\subset\Hcal^{n+1}$; a real number $r>0$; and compact sets $Q_j\subset V_n$, such that
	$\Hcal^{n+1}=\Gamma\bigl(\mathcal{K}\cup\bigcup_{j}\sigma_j\mathcal{S}_j\bigr)$ with $\mathcal{S}_j=\lbrace w\colon \pi(w)\in Q_j,\ h(w)\geq r\rbrace$.
	Here $Q_j$ is a fundamental parallelogram for the lattice $\Lambda_j$ of translations in $\Gamma_j'=\sigma_j^{-1}\Gamma\sigma_j$.
	Choose $\varepsilon_0>0$ smaller than $h$ on $\mathcal{K}\cup\mathcal{S}_1$.
	If $x\in\overline{D}$ and $h(x)<\varepsilon_0$, write $x=g^{-1}(y)$ with $g\in\Gamma$ and $y$ in this union. Since $x\in B$ we have $h(y)\leq h(x)<\varepsilon_0$, so $y\in\sigma_j\mathcal{S}_j$ for some $j\geq2$.
	That is, $x=\rho(w)$ with $\rho\in\Gamma\sigma_j$ and $h(w)\geq r$.

	Fix $j\geq 2$, and choose $N\geq1$ with $\Gamma\sigma_j\cup\Gamma'_j\subset M_2(N^{-1}\Ocal)$.
	Let $\begin{psmallmatrix}a&b\\0&a^{*-1}\end{psmallmatrix}$ lie in the stabilizer $\Gamma'_{j,\infty}$ of $\infty$ in $\Gamma'_j$. All its powers have nonzero diagonal entries in the lattice $N^{-1}\Ocal$, so $\vert a\vert=1$.
	Thus $\Gamma'_{j,\infty}$ preserves $h$ and acts on $V_n$ by Euclidean isometries. It contains the translations $\tau_\lambda$ by $\lambda\in\Lambda_j$. Let $R_j$ be the covering radius of $\Lambda_j$.
	For $\rho=\begin{psmallmatrix}A&*\\C&D\end{psmallmatrix}\in\Gamma\sigma_j$ we have $C\neq0$, since $\rho(\infty)$ is not $\Gamma$-equivalent to $\infty$.
	Put $H_\rho=\rho(\lbrace h\geq r\rbrace)$. This is a horoball at $\rho(\infty)$ that depends only on $\rho(\infty)$.
	For $\rho,\rho_0\in\Gamma\sigma_j$ and $\tau\in\Gamma'_j$ we have $\rho_0\tau\rho^{-1}\in\Gamma$. \\

	\noindent \emph{Step 3 (finitely many cusps).}
	Let $m_j$ be the minimum of $\vert C\vert$ over $\rho\in\Gamma\sigma_j$. It exists because $C\in N^{-1}\Ocal\setminus\lbrace0\rbrace$.
	Let $x=\rho(w)\in\overline{D}$ with $\rho\in\Gamma\sigma_j$ and $h(w)\geq r$.
	Choose $\rho_0\in\Gamma\sigma_j$ with lower row $(C_0,D_0)$ and $\vert C_0\vert=m_j$, and choose $\lambda\in\Lambda_j$ with $\vert\pi(w)+\lambda+C_0^{-1}D_0\vert\leq R_j$.
	Since $x$ has maximal height in its orbit and $\rho_0\tau_\lambda\rho^{-1}\in\Gamma$, \eqref{E:cx+d} gives
	$$\frac{1}{\vert C\vert^{2}h(w)}\;\geq\; h(x)\;\geq\; h(\rho_0\tau_\lambda(w))\;\geq\;\frac{h(w)}{m_j^2\,(R_j^2+h(w)^2)}.$$
	Hence $\vert C\vert\leq m_j(1+R_j^2/r^2)^{1/2}$.
	The horoball $H_\rho$ is the Euclidean ball tangent to $V_n$ at $\rho(\infty)=AC^{-1}$ of diameter $1/(r\vert C\vert^2)\leq 1/(rm_j^2)$, and $\pi(x)\in\overline{F}$, so $\vert AC^{-1}\vert$ is bounded.
	Since $A,C\in N^{-1}\Ocal$, with $\vert C\vert$ and $\vert A\vert=\vert AC^{-1}\vert\,\vert C\vert$ bounded, the cusp $\rho(\infty)$ lies in a finite set $S_j$.
	For each $c\in S_j$ fix $\rho_c\in\Gamma\sigma_j$ with $\rho_c(\infty)=c$ and lower row $(C_c,D_c)$, and put $p_c=-C_c^{-1}D_c\in V_n$. \\

	\noindent \emph{Step 4 (bounded links).}
	Let $c\in S_j$ and $x=\rho_c(w)\in\overline{D}$ with $h(w)\geq r$.
	For $\lambda\in\Lambda_j$, by \eqref{E:cx+d},
	$h(\rho_c\tau_\lambda(w))=h(w)/\bigl(\vert C_c\vert^2(\vert\pi(w)+\lambda-p_c\vert^2+h(w)^2)\bigr)$.
	Since $\rho_c\tau_\lambda\rho_c^{-1}\in\Gamma$ and $x$ has maximal height in its orbit, $\vert\pi(w)-p_c\vert\leq\vert\pi(w)+\lambda-p_c\vert$ for all $\lambda$, hence $\vert\pi(w)-p_c\vert\leq R_j$. \\

	\noindent \emph{Step 5 (conclusion).}
	Let $x\in\overline{D}$ and $\gamma\in\Gamma$ with $\gamma(x)\in\overline{D}$, and suppose $h(x)<\varepsilon_0$.
By Steps 2 and 3 we have $x=\rho(w)$ with $\rho\in\Gamma\sigma_j$, $c=\rho(\infty)\in S_j$ for some $j\geq2$ and $h(w)\geq r$. Then $\rho_c^{-1}\rho\in\Gamma'_{j,\infty}$, which preserves heights (Step 2), so replacing $w$ by $\rho_c^{-1}\rho(w)$ we may write $x=\rho_c(w)$ with $h(w)\geq r$. 
Then $\gamma(x)=(\gamma\rho_c)(w)\in\overline{D}$ with $\gamma\rho_c\in\Gamma\sigma_j$, so $\gamma(c)\in S_j$ by Step 3, and $\tau=\rho_{\gamma(c)}^{-1}\gamma\rho_c\in\Gamma'_{j,\infty}$.
	By Step 4, $\pi(w)$ and $\pi(\tau(w))$ lie within $R_j$ of $p_c$ and $p_{\gamma(c)}$ respectively.
	So the entries of $\tau$ are bounded, and they lie in $N^{-1}\Ocal$. Hence $\tau$, and with it $\gamma=\rho_{\gamma(c)}\tau\rho_c^{-1}$, lies in a finite set.
	Together with Step 1, this shows that $\Sigma$ is finite. \\

	\noindent \emph{Local finiteness.}
	Let $C\subset\lbrace \varepsilon\leq h\leq R\rbrace$ be compact, and let $y=\gamma(x)\in C$ with $x\in\overline{D}$.
	Then $h(x)\geq h(y)\geq\varepsilon$. If $\gamma\in\Gamma_\infty$ then $h(x)=h(y)\leq R$; otherwise $h(y)\leq 1/(\vert c\vert^2h(x))$ by \eqref{E:cx+d}, so $h(x)\leq1/\varepsilon$.
	So $x$ lies in the compact set $\overline{D}\cap\lbrace\varepsilon\leq h\leq\max(R,1/\varepsilon)\rbrace$, and only finitely many $\gamma$ occur. \\

	\noindent \emph{Bubbles and sides.}
	Suppose a bubble meets $\overline{D}$ at $x$. By the proof of Theorem~\ref{T:characterization-of-bubble-domain}, the bubble is $\lbrace y\colon\vert cy+d\vert=1\rbrace$ for some $\eta=\begin{psmallmatrix}a&b\\c&d\end{psmallmatrix}\in\Gamma\setminus\Gamma_\infty$. Then $h(\eta(x))=h(x)$, so $\eta(x)\in B$, and $t\eta(x)\in\overline{D}$ for some $t\in\Gamma_\infty$.
	Now $t\eta\in\Sigma$, and its lower row is a unit multiple of $(c,d)$, so the bubble is determined by $t\eta$. Hence only finitely many bubbles meet $\overline{D}$.
	Let $P$ be the intersection of $\pi^{-1}(\overline{F})$ with the closed exteriors of these bubbles.
If $p\in P\setminus\overline{D}$, the geodesic segment from a point of $D$ to $p$ leaves $\overline{D}$ through a bubble meeting $\overline{D}$ (only finitely many bubbles meet a compact neighbourhood of the exit point, by the last part of the proof of Theorem~\ref{T:characterization-of-bubble-domain}), and it stays inside that bubble afterwards, which contradicts $p\in P$. 
So $\overline{D}=P$ is a convex polyhedron with finitely many sides.
\end{proof}

\section{Generators and Relations for $\PSL_2(\Ocal)$}\label{sec:gens-and-rels}
There is a general philosophy that information about the finite presentation of a group that acts discretely on a topological space comes from the geometry of the fundamental domain of the action. 
These ideas go back to Poincar\'{e}, and many variations of it can be found in many textbooks. 
We couldn't find a reference that presented this material exactly as we wanted it, so in this section we cover it in a format that will work with our applications. 

A good general reference for the relationship between finite presentations and groups acting on topological spaces, for example,  is \cite{Bridson1999}.
The authors of \cite{Bridson1999} cover the ground, in some sense, three times: in chapters I.8 (starting on p.~131), II.12, and III.${\mathcal C}$, at an increasing degree of sophistication.  
	
\subsection{Setup}
Recall the definition of a fundamental domain (Defn.~\ref{defn:fundamentalDomain}).
In our application, $E$ will be $\Hcal^{n+1}$, on which $\Gamma$ acts properly discontinuously by isometries, and $\overline{D}$ will be a convex geodesic polyhedron \cite[Chapter~6]{Ratcliffe2019} whose tiles are locally finite: every compact subset of $\Hcal^{n+1}$ meets only finitely many tiles $\gamma\overline{D}$.
The {\em $k$-faces} are the faces of codimension $k$, and in particular the maximal proper faces are {\em facets} or \emph{sides}.

Form a graph $G=G(\Gamma,D)$ where the vertices are open fundamental domains (equivalently elements of $\Gamma$ or closed fundamental domains), and edges exist between $D$ and $D'$ if and only if the closed fundamental domains intersect in a facet.
\begin{definition}
	We will call the graph $G$ the \emph{tessellation graph}.
	By a \emph{basic step} we mean an edge in the graph.
	We will denote the neighbors of $G$ at $T$ by $N_G(T)$ and the basic steps from a $T$ by $E_G(T)$.
\end{definition}

Every basic step corresponds to an element of $\Gamma$ which takes the fundamental domain to its neighbor. The goal will be to convert walks in the graph $G$ to elements of the group $\Gamma$.
Every vertex of this graph has the same degree, namely the number of facets of $D$.

\begin{proposition}\label{prop:good-neighbours}
	Let $P$ be a point of a facet of a tile $T$ that is not in
	any face of codimension~$2$ of any tile. Then a sufficiently small neighborhood of $P$ contains
	points of exactly two tiles, $T$ and a tile $T'$, and $T\cap T'$ has dimension~$n$.
	The tile $T'$ may depend on $P$ (Remark~\ref{rem:subdivide-facets}).
\end{proposition}

\begin{proof} By local finiteness, a small ball $U$ about $P$ meets only
	tiles that contain $P$, and since $P$ lies in the relative interior of a facet of $T$,
	$U\cap T$ is a half-ball bounded by the hyperplane of that facet.
	Some tile $T'\neq T$ meets $U\setminus T$ in a nonempty open set.
	As the interiors of distinct tiles are disjoint, $P$ is not an interior point of $T'$,
	so, $P$ not being in a face of codimension~$2$ of $T'$, it lies in the relative interior of a facet of $T'$,
	and after shrinking $U$ the set $U\cap T'$ is also a half-ball.
	Disjointness of interiors forces the two half-balls to be the two sides of the same hyperplane.
	Hence $U\subset T\cup T'$, no third tile meets $U$, and $T\cap T'$ contains a
	neighborhood of $P$ in that hyperplane.
\end{proof}

\begin{remark}\label{rem:subdivide-facets}
  It is important to understand that a point may belong to a
  single facet of a fundamental domain but two or more facets of a translate
  of the fundamental domain.  For example, consider the plane $\RR^2$ with the
  standard action of $\ZZ^2$.  This does not happen if the fundamental domain
  is the obvious square $0 < x,y < 1$, but it does if the fundamental domain
  is taken to be a parallelogram with vertices
  $(0,0), (1,0), (-1/2,1), (1/2,1)$, for example.  In order to deal with this,
  we subdivide the facets along the codimension-$2$ intersections of facets
  of translates that meet on the boundary of the fundamental domain.  In this
  example we would break the facet joining $(0,0)$ to $(1,0)$ into
  $(0,0),(1/2,0)$ and $(1/2,0),(1,0)$, and similarly we would introduce
  $(0,1)$ between $(-1/2,1)$ and $(1/2,1)$.
  Precisely, the \emph{subdivided facets} of $\overline{D}$ are the sets
  $\overline{D}\cap\gamma\overline{D}$, $\gamma\in\Gamma$, of dimension~$n$.
  By Proposition~\ref{prop:good-neighbours}, every point of $\partial\overline{D}$
  outside the faces of codimension~$2$ of the tiles lies in exactly one of them,
  and $\overline{D}\cap\gamma\overline{D}$ borders only the tiles $\overline{D}$ and $\gamma\overline{D}$.
\end{remark}

\subsection{Finite Generation}
We now fix some notation.

\begin{notation}
\begin{itemize}
\item Let $F_1, \dots, F_m$ be the facets of $D$, subdivided in
  accordance with Remark~\ref{rem:subdivide-facets}.  
\item Let $N_G(D) = \lbrace N_1,\ldots, N_m \rbrace$ be the corresponding neighbors.  (The subdivision in the first step ensures that the number of neighbors
  is equal to the number of facets and also guarantees the uniqueness in the remaining steps.)
	\item Define $\gamma_i$ be the unique elements of $\Gamma$ such that $\gamma_i D = N_i.$
	\item Given a tile $T$, define $\gamma_T \in \Gamma$ to be the unique element such that $\gamma_T D = T.$
	\item Define $F_{i,T} = \gamma_T F_i.$
	\item Define the $N_{i,T}$ to be the second tiles associated to the $F_{i,T}$ in
	Proposition~\ref{prop:good-neighbours}.
\end{itemize}
\end{notation}

\begin{definition} A {\em walk} is a sequence of tiles $T_0, T_1, \dots, T_r$,
	such that $T_{i-1}$ and $T_i$ share a facet for all $1 \le i \le r$.
	The {\em route} of a walk is the sequence $a_1, \dots, a_r \in \lbrace 1, \ldots, m \rbrace$ such that
	$T_{i-1}$ and $T_i$ share the facet $F_{a_i,T_{i-1}}$ for all $i$.
\end{definition}

Consider the neighbors of $D$ in the graph $G$.
When $\overline{D}$ has finitely many subdivided facets we write $N_G(D) = \lbrace D_1,\ldots,D_r\rbrace$.
Each neighbor $D_i$ is a basic step away from $D$ and hence gives a well-defined element $\gamma_i$ and the basic steps of $D$ are  $E_G(D)=\lbrace \gamma_1,\gamma_2,\ldots,\gamma_r\rbrace.$
Note that $N_G(\gamma D) = \lbrace \gamma D_1,\ldots,\gamma D_r\rbrace$ and that the basic steps of $\gamma D$ are  $E_G(\gamma D)=\lbrace \gamma \gamma_1 \gamma^{-1}, \gamma \gamma_2 \gamma^{-1},\ldots, \gamma \gamma_r \gamma^{-1} \rbrace.$
The neighborhoods are pictured in Figure~\ref{F:neighbors}.
Using this notion of basic step we can now prove finite generation of $\Gamma$ and give the generators explicitly.
\begin{figure}[htbp!]
	\centering
	\begin{tikzpicture}[>=stealth]
  \begin{scope}
    \node[circle, draw, minimum size=0.9cm] (D) at (0,0) {$D$};
    \foreach \k in {1,...,8} {
      \pgfmathsetmacro{\a}{45*(\k-1)}
      \draw[->, thick] (\a:0.45) -- (\a:1.5);
      \node at (\a:1.9) {$\gamma_{\k}$};
    }
  \end{scope}
  \begin{scope}[shift={(6.2,0)}]
    \node[circle, draw, minimum size=0.9cm] (gD) at (0,0) {$\gamma D$};
    \foreach \k in {1,...,8} {
      \pgfmathsetmacro{\a}{45*(\k-1)}
      \draw[->, thick] (\a:0.45) -- (\a:1.5);
      \node at (\a:2.3) {$\gamma\gamma_{\k}\gamma^{-1}$};
    }
  \end{scope}
\end{tikzpicture}
	\caption{The basic steps of $D$ and of $\gamma D$.
		If the basic steps of $D$ are $\gamma_1,\dots,\gamma_r$, then those of $\gamma D$ are $\gamma\gamma_1\gamma^{-1},\dots,\gamma\gamma_r\gamma^{-1}$ (here $r=8$).
	}\label{F:neighbors}
\end{figure}
\begin{theorem}\label{thm:generators-polyhedron}
  Let $\Gamma$ act properly discontinuously by isometries on $\Hcal^{n+1}$, and let $D$ be an open fundamental domain whose closure is a convex geodesic polyhedron with locally finite translates.
  The group $\Gamma$ is generated by the collection of $\gamma\in \Gamma$ such that $\overline{D}\cap\gamma\overline{D}$ is a subdivided facet of $\overline{D}$ (Remark~\ref{rem:subdivide-facets}).
  In the notation above $\Gamma = \langle \gamma_1,\ldots, \gamma_r\rangle$.
\end{theorem}
\begin{proof}
	Let $\Gamma'$ be the subgroup generated by these elements, and put
	$U=\bigcup_{h\in\Gamma'}h\overline{D}$ and $V=\bigcup_{\gamma\in\Gamma\setminus\Gamma'}\gamma\overline{D}$.
	Both are closed because the tiles are locally finite, and $U\cup V=\Hcal^{n+1}$.
	Let $\Sigma$ be the union of the faces of codimension at least~$2$ of all tiles.
	If $h\in\Gamma'$, $\gamma\notin\Gamma'$ and $h\overline{D}\cap\gamma\overline{D}$ contains a point outside $\Sigma$,
	then by Proposition~\ref{prop:good-neighbours} the set $h\overline{D}\cap\gamma\overline{D}$, and hence
	$\overline{D}\cap h^{-1}\gamma\overline{D}$, has dimension~$n$, so $h^{-1}\gamma$ is one of the generators and $\gamma\in\Gamma'$, a contradiction.
	Thus $U\cap V\subset\Sigma$.
	Since $\Sigma$ is a locally finite union of closed convex sets of dimension at most $n-1$, it does not
	separate $\Hcal^{n+1}$, so $\Hcal^{n+1}\setminus\Sigma$ is connected. 
It is the disjoint union of the relatively closed
	sets $U\setminus\Sigma$, which contains $D$, and $V\setminus\Sigma$. Hence $V\setminus\Sigma=\emptyset$, and since every tile has
	nonempty interior, $V=\emptyset$, that is, $\Gamma'=\Gamma$.
	Concretely, a path from a point of $D$ to a point of $\gamma D$ avoiding $\Sigma$ crosses finitely many
	subdivided facets and so traces a walk in the tessellation graph; see Figure~\ref{F:walk}.
	Proposition~\ref{prop:walk} below turns its route into a word in $\gamma_1,\ldots,\gamma_r$.

\begin{figure}[htbp!]
	\centering
	\begin{tikzpicture}[>=stealth,
    tile/.style={circle, draw, thick, minimum size=0.75cm, inner sep=0pt, fill=white},
    nb/.style={circle, draw=black!35, minimum size=0.35cm, inner sep=0pt},
    wstep/.style={->, very thick},
    lab/.style={font=\small, fill=white, inner sep=1pt}]
  \begin{scope}
    \node[tile] (T0) at (0,0) {$D$};
    \node[tile] (T1) at (1.8,0.9) {$T_1$};
    \node[tile] (T2) at (3.6,0.4) {$T_2$};
    \node[tile] (T3) at (5.2,1.4) {$T_3$};
    \node[tile] (T4) at (6.4,2.9) {$T_4$};
    \foreach \t/\angs in {T0/{150,200,250,300}, T1/{60,110,160}, T2/{-40,-90,-140}, T3/{10,-35,-80,130}, T4/{-10,40,90,140,190}}
      \foreach \a in \angs { \draw[black!35] (\t.\a) -- ++(\a:0.55) node[nb, anchor=\a+180] {}; }
    \draw[wstep] (T0) -- (T1) node[midway, lab, above left] {$a_1$};
    \draw[wstep] (T1) -- (T2) node[midway, lab, above] {$a_2$};
    \draw[wstep] (T2) -- (T3) node[midway, lab, above left] {$a_3$};
    \draw[wstep] (T3) -- (T4) node[midway, lab, left] {$a_4$};
    \node[font=\small] at (3.2,-1.35) {$T_k=\gamma_{a_1}\cdots\gamma_{a_k}D$};
    \node at (3.2,-2.0) {(a)};
  \end{scope}
  \begin{scope}[shift={(9.6,0.8)}]
    \node[tile] (C0) at (0,-1) {$D$};
    \node[tile] (C1) at (1.5,0) {$T_1$};
    \node[tile] (C2) at (0.6,1.5) {$T_2$};
    \node[tile] (C3) at (-1.0,1.3) {$T_3$};
    \node[tile] (C4) at (-1.5,-0.2) {$T_4$};
    \draw[wstep] (C0) -- (C1) node[midway, lab, below right] {$a_1$};
    \draw[wstep] (C1) -- (C2) node[midway, lab, above right] {$a_2$};
    \draw[wstep] (C2) -- (C3) node[midway, lab, above] {$a_3$};
    \draw[wstep] (C3) -- (C4) node[midway, lab, left] {$a_4$};
    \draw[wstep] (C4) -- (C0) node[midway, lab, below left] {$a_5$};
    \node[font=\small] at (0,-1.9) {$\gamma_{a_1}\gamma_{a_2}\gamma_{a_3}\gamma_{a_4}\gamma_{a_5}=1$};
    \node at (0,-2.8) {(b)};
  \end{scope}
\end{tikzpicture}
	\caption{(a) A walk $D=T_0,T_1,\dots,T_4$ in the tessellation graph with route $a_1,\dots,a_4$; the faint nodes are other neighbours.
		By Proposition~\ref{prop:walk}, $T_k=\gamma_{a_1}\cdots\gamma_{a_k}D$.
		(b) A closed walk from $D$ back to $D$ gives a relation in $\Gamma$ (Corollary~\ref{cor:loop-relations}).
	}\label{F:walk}
\end{figure}
\end{proof}

\subsection{Finite Relations}\label{sec:finite-relations}

All of the relations in our group, aside from the simple ones that assert that
one generator is the inverse of another (or itself),
come from taking cycles in the tessellation graph.
This means we need an effective way of converting routes into group elements. 
The recipe is given by simply flipping the order of the route and multiplying the associated group elements. Here $\conj_g(x)=gxg^{-1}$. 
Before diving into a general proof, we write out two examples to clarify things. 
\begin{example}
	The route 12 first applies $\alpha_1=\gamma_1$ then applies $\alpha_2 = \conj_{\alpha_1}(\gamma_2)$ which gives $ \alpha_2\alpha_1 = \gamma_1\gamma_2\gamma_1^{-1}\cdot \gamma_1= \gamma_1 \gamma_2.$
\end{example}

\begin{example}
	Consider the route 123. 
	Then one has $\alpha_1=\gamma_1$ first, then $\alpha_2=\conj_{\alpha_1}(\gamma_2)$, then $\alpha_3 =\conj_{\alpha_2\alpha_1}(\gamma_3)$.
	This gives $\alpha_3\alpha_2\alpha_1$ being the element corresponding to the route $123$.
	Writing this out gives
	$ \alpha_3 \alpha_2 \alpha_1 = \conj_{\alpha_2\alpha_1}(\gamma_3) \conj_{\alpha_1}(\gamma_2) \alpha_1
= \gamma_1\gamma_2 \gamma_3 (\gamma_1\gamma_2)^{-1} \cdot \gamma_1 \gamma_2 \gamma_1^{-1} \cdot \gamma_1
	=\gamma_1\gamma_2 \gamma_3.$ \end{example}

We now give the general proof.
\begin{proposition}\label{prop:walk} Let $D = T_0, \dots, T_r$ be a walk with route $a_1, \dots, a_r$.
	Then $T_r = \gamma_{a_1} \dots \gamma_{a_r} D.$
\end{proposition}
\begin{proof} Induction on $r$.  For $r = 0$ there is nothing to do, and for
	$r = 1$ this is the definition.  
	
	So suppose it is true for $r = k$ and
	let us prove it for $r = k+1$.  
	Let $B = D \cup \gamma_{a_{k+1}} D$;
	let $w_k = \gamma_{a_1} \dots \gamma_{a_k}$. Since $w_k D = T_k$, it
	follows that $w_k B = T_k \cup N_{a_{k+1},T_k} = T_k \cup T_{k+1}$.
	By definition, $w_k(D) = T_k$, and $w_k(\gamma_{a_{k+1}} D) \cap w_k(D) = w_k(D \cap \gamma_{a_{k+1}} D) = w_k(F_{a_{k+1}}) = F_{a_{k+1},w_k}$ by how the neighbors were set up.  So $w_k(\gamma_{a_{k+1}} D)$ is the tile that shares $F_{a_{k+1},w_k}$ with $w_k(D)$. Again, by definition this is $T_{k+1}$.
	But since $w_k D = T_k$, it follows that $w_k \gamma_{a_{k+1}} D = T_{k+1}$. Since
	$w_k \gamma_{a_{k+1}} = w_{k+1}$, this is what we wanted.
	
Here is a second proof: let $\gamma$ be the element which takes $D$ to $T_{k+1}$. 
We have $\gamma = \alpha_2\alpha_1$ where by hypothesis we can write $\alpha_1= \gamma_{a_1}\cdots \gamma_{a_k}$ and $\alpha_2 = \conj_{\alpha_1}(\gamma_{a_{k+1}} )$.
We compute 
 $$ \gamma = \alpha_2\alpha_1 = \alpha_1\gamma_{a_{k+1}} \alpha_1^{-1} \alpha_1 = \alpha_1 \gamma_{a_{k+1}} = \gamma_{a_1}\gamma_{a_2} \cdots \gamma_{a_k} \gamma_{a_{k+1}},$$
which proves the result. 
\end{proof}

\begin{corollary}\label{cor:switch-relations} Suppose that the facet $F_k$ of
	$D$ is labeled $F_{k',N_k}$ as a facet of $N_k$.
	Then $\gamma_{k} \gamma_{k'} = 1$.
\end{corollary}
\begin{proof}
	The tile across $F_{k',N_k}=\gamma_kF_{k'}$ from $N_k=\gamma_kD$ is $N_{k',N_k}=\gamma_k\gamma_{k'}D$; it is also $D$.
	So $\gamma_k\gamma_{k'}D=D$, and $\gamma_k\gamma_{k'}=1$, since an element of $\Gamma$ mapping a point of $D$ into $D$
	is the identity (Definition~\ref{defn:fundamentalDomain}).
\end{proof}

\begin{corollary}\label{cor:loop-relations}
	Let $D = T_0, \dots, T_r = D$ be a walk with route $a_1, \dots, a_r$.
	Then $\gamma_{a_1} \dots \gamma_{a_r}$ is a relation in $\Gamma$.
\end{corollary}
\begin{proof}
	By Proposition~\ref{prop:walk}, $\gamma_{a_1}\cdots\gamma_{a_r}D=T_r=D$. An element of $\Gamma$ that maps $D$ to itself is the identity by Definition~\ref{defn:fundamentalDomain}, giving the relation.
\end{proof}

Now let $H$ be a $2$-face of $\overline{D}$, and let $P$ be a point of $H$ lying in no face of codimension~$3$ of any tile.
By local finiteness only finitely many tiles contain $P$.
By slicing near $P$ with a generic
$2$-plane we get an arrangement of cones in $\RR^2$, whose $1$-faces are
the restrictions of facets of tiles in the tiling. So we associate a walk
to $H$, starting at $D$, by crossing these facets in order once around $P$. Since this walk starts and ends
on the same tile, it gives a relation as in Corollary~\ref{cor:loop-relations}.
We call this the $H$-relation. (If the set of tiles containing $P$ varies along $H$, we subdivide $H$
as in Remark~\ref{rem:subdivide-facets} and take one relation for each piece.)
If $H'=\gamma H$ is another $2$-face of $\overline{D}$ with $\gamma\in\Gamma$, the $H'$-relation is a cyclic
permutation of the $H$-relation or of its inverse, so it suffices to take one $H$ from each $\Gamma$-orbit.
\begin{theorem}\label{thm:pres-gamma}
	$\Gamma$ is presented by the generators $\gamma_1,\gamma_2,\ldots$ of Theorem~\ref{thm:generators-polyhedron},
	the relations of Corollary~\ref{cor:switch-relations}, and the $H$-relations,
	one for each $2$-face $H$ of $\overline{D}$ up to the action of $\Gamma$.
	If only finitely many tiles meet $\overline{D}$, this presentation is finite.
\end{theorem}
\begin{proof}
	This is the presentation part of Poincar\'e's polyhedron theorem; compare
	\cite[Theorem~I.8.10]{Bridson1999}, \cite[Chapter~6]{Ratcliffe2019} and \cite{Swan1971}. We sketch the argument.
	The given relations hold by Corollaries~\ref{cor:switch-relations} and~\ref{cor:loop-relations}.
	Conversely, let $\gamma_{a_1}\cdots\gamma_{a_m}=1$. By Proposition~\ref{prop:walk} this is the route of a closed
	walk at $D$, which we realize by a loop $\ell$ based in $D$ that crosses subdivided facets transversally in the
	order of the route and avoids the union $\Sigma_2$ of the faces of codimension at least~$2$ of the tiles.
	Let $\Sigma_3\subset\Sigma_2$ be the union of the faces of codimension at least~$3$. It is a locally finite union of
	closed convex sets of codimension at least~$3$, so $\Hcal^{n+1}\setminus\Sigma_3$ is simply connected and, by general position,
	$\ell$ is null-homotopic through a disk in $\Hcal^{n+1}\setminus\Sigma_3$ meeting $\Sigma_2$ in finitely many points of $2$-faces.
	A homotopy in $\Hcal^{n+1}\setminus\Sigma_2$ changes the route only by inserting or deleting a crossing of a facet followed by the
	crossing back, which changes the word by a relation of Corollary~\ref{cor:switch-relations}.
	Passing the loop across a point of a $2$-face $\gamma_T H$ of a tile $T$ changes the word by inserting a conjugate of the
	$H$-relation or of its inverse. Hence $\gamma_{a_1}\cdots\gamma_{a_m}$ lies in the normal subgroup generated by the given relations.
	If only finitely many tiles meet $\overline{D}$, there are finitely many generators and finitely many $H$-relations.
\end{proof}

\section{Algorithms}\label{sec:algorithms}

In this section, we present what is needed for the computation of fundamental
domains for $\SL_2(\Ocal)$, where $\Ocal$ is a Clifford-stable order in a Clifford
algebra over $\QQ$. The results and algorithms of this section are not found in
\cite{Bygott1998,Lingham2005,Rahm2010} in exactly this form.
Nevertheless, the reader familiar
with these dissertations will immediately recognize our debt to their authors.

\subsection{Overview}
As stated in the introduction, we will explicitly compute the fundamental domains of $\PSL_2(\Ocal)$ for various orders $\Ocal$ in \texttt{magma}. 
The algorithm for producing the maximal orders is found in \S\ref{sec:maximal-orders}, essentially Algorithm~\ref{alg:p-maximal} after some discriminant considerations. 

As we have stated previously, many of the orders we found were Clifford-Euclidean.
The division algorithm (part of the definition of what it means for an order to be Clifford-Euclidean), Algorithm~\ref{alg:gcd}, gives the gcd $\gamma$ of two elements $\alpha,\beta \in \Ocal$, and Algorithm ~\ref{alg:gcd-coefs} tells us which $\lambda,\mu\in \Ocal$ given $\lambda \alpha + \mu \beta = \gamma$.
For example, this tells us that a one-sided ideal generated by two elements of $\Ocal^{\mon}$ to which the division algorithm applies is principal, and gives a generator. It does not make every one-sided ideal principal: $\Ocal_4$ is Clifford-Euclidean but not Clifford-principal (\S\ref{sec:maximal-c4}).

We have seen in \S\ref{sec:fundamental-domains} that the fundamental domain $D$ consists of the set of elements $x\in \Hcal^{n+1}$, which are above the boundary bubbles and project to the fundamental domain $F$ for the stabilizer of $\infty$, $\PSL_2(\Ocal)_{\infty}$. 
Since $D$ is convex (we take $F$ to be a convex polytope), determining the faces of $D$ is equivalent to determining the facets of $D$ and these fall into two different classes.

First, there are the sides coming from the fundamental domain  $F$ for $\PSL_2(\Ocal)_{\infty} \cong \Vec(\Ocal) \rtimes (\Ocal^{\times}/\lbrace \pm 1\rbrace)$.
This amounts to computing $\Ocal^{\times}$: Algorithm~\ref{alg:units} computes the group of Clifford units of an order. 
(We also used built-in algorithms in \texttt{magma} for dealing with lattices, which we do not discuss.)

Second, there are the sides coming from the bubbles that lie above $F$.
The simplest possible algorithm consists of producing many boundary bubbles, then only taking those which appear as a maximum above a point $y\in F$. 
We can exclude bubbles using Lemma~\ref{lem:singly-dominated}.
This is implementable, and we did do this at first, but it is not
conducive to a rigorous determination of the fundamental domain: we can
list all cusps in order, but we do not know a priori when to stop
(Remark~\ref{rem:stopping} gives an a posteriori criterion).
(Question~\ref{qstn:bubble-bound} in our open problems is exactly this
question.) Also, this produces an inefficient description of the fundamental
domain, since there can be bubbles that are dominated by a set of other
bubbles without being dominated by any of them individually.
It is more efficient to use linear programming.
This approach is given in Algorithm~\ref{alg:under-spheres}.
This algorithm also adds bubbles when the linear program is unbounded (Step~\ref{item:first-solve}).

Finally, this section describes an algorithm for computing generators of
$\SL_2(\Ocal)$ in certain cases.
First we do some abstract mathematics, then we concern ourselves with algorithms. 
Theorem~\ref{thm:generators-ocal} describes generators for $\PSL_2(\Ocal)$ and it suffices for each side $H$ of $\overline{D}$ to find some $g_H\in \Gamma$ such that $g(\overline{D}) \cap \overline{D} = H$.

Hence, to produce generators it suffices to find $g_H$ for each facet $H$ of $D$. 
Again we break this down into bubbles and non-bubble facets. 
This is easier for the sides corresponding to the sides of $F$.
One can translate by a generator of $\Vec(\Ocal)$ and then rotate by an element of $\Ocal^{\times}$.
Hence, it suffices to deal with the bubble case; for each bubble $H$ we need to find $g_H\in \Gamma$.
The general problem here seems tractable but would require investigation beyond the scope of this manuscript (see Question~\ref{qstn:crossing-bubbles}).
In the known cases where the domain is Clifford-Euclidean, we find that adding $x\mapsto -x^{-1}$ suffices. 
We suspect this holds for Clifford-Euclidean orders in general
(see Conjecture~\ref{conj:euclidean-bubbles}).
We are uncertain of how difficult this problem is, but it should be attempted in future investigations.
One approach, though perhaps not the most insightful,
would be to show there are finitely many Clifford-Euclidean orders and then compute the fundamental domains for all of those.

\subsection{Distance Lemma}\label{sec:distance-lemma}
  	The following concerns spheres in $V_n = \Vec(\CC_n)$ and the element 
  	 $$ h=\frac{1+i_1 + \cdots + i_{n-1}}{2}, $$
  	which is a deep hole for the standard cubic lattice $\ZZ^n = \Vec(\ZZ[i_1,i_2,\ldots,i_{n-1}])$.
Figure~\ref{fig:deep-hole} illustrates the lemma for $n=2$.
\begin{figure}[htbp!]
\centering
\begin{tikzpicture}[x=7cm,y=7cm,lab]
 \begin{scope}\clip (0,0) rectangle (0.5,0.5);
  \fill[covC!15] (0,0) circle (0.6); \fill[holeC!15] (0.5,0.5) circle (0.374);
  \draw[cov] (0,0) circle (0.6); \draw[holeC,line width=0.8pt,dashed] (0.5,0.5) circle (0.374);
 \end{scope}
 \draw[line width=0.8pt] (0,0) rectangle (0.5,0.5);
 \dotat{0,0} \node[below left] at (0,0){$0$};
 \dotat[holeC]{0.5,0.5} \node[right=3pt] at (0.5,0.5){$h=\tfrac{1+i_1}{2}$};
 \draw[covC,-{Stealth}] (0,0)--(0.6*0.8,0.6*0.6) node[pos=0.55,below right]{$r_1$};
 \draw[holeC,-{Stealth}] (0.5,0.5)--++(-0.374*0.5,-0.374*0.866) node[pos=0.6,right]{$r_2$};
 \node[below] at (0.25,0){$x_0$}; \node[left] at (0,0.25){$x_1$};
\end{tikzpicture}
\caption{Lemma~\ref{lem:deephole} drawn for $n=2$ with $r_1=0.6$, $r_2=\sqrt{0.14}$, so $r_1^2+r_2^2=n/4$: every point $a$ of the box $\{0\le a_j\le 1/2\}$ (black square) lies within $r_1$ of $0$ (green) or within $r_2$ of the deep hole $h$ (purple). The lemma is used for $n=4,5$; this picture is the reduced-dimension case.}\label{fig:deep-hole}
\end{figure}
  \begin{lemma} \label{lem:deephole}
  	Let $r_1$ and $r_2$ be non-negative real numbers such that
  	 $r_1^2+r_2^2=n/4$. 
   If $a =a_0+a_1i_1 + \cdots +a_{n-1}i_{n-1}$ with $0\leq a_j \leq 1/2$ for $j=0,1,2,\ldots,n-1$, then 
    $$ \vert a \vert \leq r_1 \quad \mbox{ or } \quad \vert a -h \vert \leq r_2.$$
  \end{lemma}
 \begin{proof} The minimum of $f(x) = x - 2x^2$ on $[0,1/2]$ is $0$.
    Thus we have
    $$r_1^2+r_2^2 =n/4 = \sum_{j=0}^{n-1} 1/4 \ge \sum_{i=0}^{n-1}( 1/4 - a_i + 2a_i^2 )= \sum_{i=0}^{n-1} \bigl((a_i-1/2)^2 + a_i^2\bigr) = \vert a-h\vert^2 + \vert a\vert^2.$$
    It follows that either $|a|^2 \le r_1^2$ or $|a-h|^2 \le r_2^2$.  Since
    $|a|, |a-h|, r_1, r_2 \ge 0$, the result follows.
 \end{proof}
 This is useful in determining distances to the point $h$ in various dimensions.
\begin{example}
	In what follows, $a \in V_n$ satisfies $0 \leq a_j \leq 1/2$ for all $j$.
	\begin{enumerate}
		\item If $n=4, r_1=1, r_2=0$, then either $a=h$ or $\vert a \vert \leq 1$.
		\item If $n=4$, we can also say that for every $a$ we have $\vert a \vert \leq1/\sqrt{2}$ or $\vert a-h \vert \leq 1/\sqrt{2}$. 
		\item If $n=5, r_1=1, r_2=1/2$, then $\vert a \vert \leq 1$ or $\vert a -h \vert \leq 1/2.$
		Another way of seeing the same numbers is that $\vert a \vert \leq 1/2$ or $\vert a -h \vert \leq 1.$
	\end{enumerate}
\end{example}

\subsection{Bubble Algorithm}\label{sec:bubbles}
We take this opportunity to state a simple lemma that allows us to exclude
many spheres from consideration.

Figure~\ref{fig:exclusion-lemma} shows the two cases of the lemma.
\begin{figure}[htbp!]
\centering
\begin{tikzpicture}[scale=1.9,lab]
 \draw[bub] (1.6,0) arc[start angle=0,end angle=180,radius=1.6];
 \draw[bub,fill=bubbleC!12] (1.3,0) arc[start angle=0,end angle=180,radius=0.6];
 \draw[bnd] (-1.8,0)--(1.8,0);
 \dotat{0,0} \node[below] at (0,0){$P$}; \dotat{0.7,0}\node[below] at (0.7,0){$Q$};
 \draw[covC,decorate,decoration={brace,amplitude=3pt,mirror}] (0,-0.25)--(0.7,-0.25) node[midway,below=2pt]{$d$};
 \node[bubbleC] at (-0.8,1.5){$H_P$}; \node[bubbleC] at (0.7,0.3){$H_Q$};
 \node at (0,-0.75){(a) $r\ge d(P,Q)+s$};
\end{tikzpicture}\qquad
\begin{tikzpicture}[scale=1.9,lab]
 \draw[bub] (1.6,0) arc[start angle=0,end angle=180,radius=1.6];
 \draw[bub,fill=cuspC!20] (2.2,0) arc[start angle=0,end angle=180,radius=0.75];
 \begin{scope}\clip (2.2,0) arc[start angle=0,end angle=180,radius=0.75]--cycle; \fill[white] (0,0) circle (1.6); \end{scope}
 \draw[bub] (2.2,0) arc[start angle=0,end angle=180,radius=0.75]; \draw[bub] (1.6,0) arc[start angle=0,end angle=90,radius=1.6];
 \draw[bnd] (-1.8,0)--(2.4,0);
 \dotat{0,0} \node[below] at (0,0){$P$}; \dotat{1.45,0}\node[below] at (1.45,0){$Q$};
 \node[bubbleC] at (-0.8,1.5){$H_P$}; \node[cuspC] at (2.15,0.85){exposed};
 \node at (0.1,-0.75){(b) $r< d(P,Q)+s$};
\end{tikzpicture}
\caption{Lemma~\ref{lem:singly-dominated}, vertical section through $P$ and $Q$. Hemispheres $H_P$, $H_Q$ (blue) of radii $r,s$. (a) If $r\ge d(P,Q)+s$, every point of $H_Q$ is on or beneath $H_P$. (b) Otherwise a part of $H_Q$ (red) lies above $H_P$.}\label{fig:exclusion-lemma}
\end{figure}
\begin{lemma}[Exclusion Lemma] \label{lem:singly-dominated}
  Let $P, Q$ be on the boundary of $\Hcal^{n}$, and consider hemispheres
  $H_P, H_Q$ of radius $r, s > 0$ and center $P, Q$ respectively.  
  Then every point of $H_Q$ is equal to or beneath a point of $H_P$
  if and only if $r \ge d(P,Q) + s$,
  where $d$ is the usual Euclidean distance.
\end{lemma}

\begin{proof} To prove ``if'', let $R$ be another point of the plane. We have
  $d(P,R) \le d(P,Q) + d(Q,R)$. The height of the point of $H_Q$ above $R$ is
  $\sqrt{s^2 - d(Q,R)^2}$, and the height of the point of $H_P$ above $R$ is
  \begin{align*}
    &\sqrt{r^2 - d(P,R)^2} \ge \sqrt{r^2 - (d(P,Q)+d(Q,R))^2} \ge
    \sqrt{(d(P,Q)+s)^2 - (d(P,Q)+d(Q,R))^2} \\
    &\quad = \sqrt{s^2 + 2sd(P,Q) - 2d(P,Q)d(Q,R) - d(Q,R)^2} \ge \sqrt{s^2 - d(Q,R)^2},
  \end{align*}
  assuming that $s \ge d(Q,R)$, which is a necessary condition for $H_Q$
  to contain a point lying above $R$.

  For ``only if'', take $R$ to be a point on the half-line from $P$ to $Q$ at a 
  distance slightly less than $s$ beyond $Q$ (if $P = Q$, any half-line from $P$ will do).
\end{proof}

More generally, we will need to consider the question of whether a hemisphere
$H_0$ lies under the union of hemispheres $H_1, \dots, H_n$. Inevitably, there
is no condition as simple as the one just given, but in practice the problem
can be solved rapidly by means of linear programming.

\begin{algorithm}\label{alg:under-spheres}
  For $i \in \{0,1,\dots,k\}$, let $H_i$ be a hemisphere with center
  $P_i \in \RR^n$ and radius $r_i$. We give a procedure that
  determines whether $H_0$ lies under $\cup_{i=1}^k H_i$.
  Let the coordinates on $\RR^n$ be $x_1, \dots, x_n$
  and let $P$ be the generic point $(x_1,\dots,x_n)$. We suppose that there
  are enough $P_i$ far away from $P_0$ in different directions that
  $r_i^2 - d(P_i,Q)^2 > r_0^2 - d(P_0,Q)^2$ for $d(P_0,Q)$ is large enough.

  \begin{enumerate}
  \item\label{item:find-ineqs}
    For each $i$, determine the linear inequality on the $x_i$ that is
    equivalent to $r_0^2-d(P_0,P)^2 \ge r_i^2-d(P_i,P)^2$. (This condition is
    linear because the $x_i^2$ terms are the same in $d(P_0,P)$ and $d(P_i,P)$.)
    Also, let $L$ be an empty list.
  \item\label{item:first-solve}
    Attempt to solve the linear program given by the inequalities of
    Step~\ref{item:find-ineqs}
    with any convenient objective function. If it is infeasible,
    then $H_0$ does lie under $\cup_{i=1}^k H_i$. If it is unbounded, then
    we need more $P_i$. (These are easily obtained in our context by
    translating by Clifford vectors in the order.) If it is solvable, add
    the solution to the list $L$.
  \item\label{item:find-0}
    Let $S$ be the set of constraints that hold with equality.
  \item\label{item:solve-new-obj}
    While $S$ is not empty, choose an element $S_j \in S$ and let it be the
    new objective function. Solve the linear programming problem and add
    the solution to the list $L$.
  \item\label{item:under}
    If the optimal value of the objective function is $0$, then
    $H_0$ does lie under $\cup_{i=1}^k H_i$.
  \item\label{item:not-under}
    If not, replace $S$ by the set of elements of $S$ that still hold
    with equality for the new optimal solution; this removes $S_j$ from $S$.
    If $S$ is not yet empty, return to Step~\ref{item:solve-new-obj}.
    If $S$ is empty, let $C_0 \subseteq \RR^n$ be the polyhedron defined by the
    inequalities of Step~\ref{item:find-ineqs}, and compute exactly
    $\delta = \min_{P \in C_0} d(P_0,P)^2$. (The minimum is attained at the
    orthogonal projection of $P_0$ onto the affine span of a face of $C_0$, so
    $\delta$ is the least value of $d(P_0,Q)^2$ over the finitely many such
    projections $Q$ that lie in $C_0$.)
    If $\delta < r_0^2$, then $H_0$ does not lie under $\cup_{i=1}^k H_i$;
    if $\delta \geq r_0^2$, then $H_0$ does lie under $\cup_{i=1}^k H_i$.
  \end{enumerate}

  The algorithm terminates: each pass through Step~\ref{item:not-under}
  removes $S_j$ from $S$, so at most $k+1$ linear programs are solved, and
  $\delta$ is found by a finite enumeration. If the closed discs in $\RR^n$
  bounded by $H_1, \dots, H_k$ cover $C_0$, then $\delta < r_0^2$ holds
  automatically once $S$ is empty. Indeed, if $P \in C_0$ lies in the disc
  of $H_j$, then $r_0^2-d(P_0,P)^2 \ge r_j^2-d(P_j,P)^2 \ge 0$, so $C_0$ lies
  in the closed disc of $H_0$, and the interior of $C_0$, which is nonempty,
  lies in the open disc. Our implementation (\texttt{can\_be\_pos.sage},
  called from \texttt{fundamental\_domain.mag}) performs only
  Steps~\ref{item:find-ineqs}--\ref{item:under}, declares $H_0$ needed as soon
  as $S$ is empty, and treats a failure of the first linear program,
  including unboundedness, as infeasibility. It is therefore correct when the
  discs cover $C_0$ and $C_0$ is bounded. When the discs do not cover $C_0$, it
  can keep a hemisphere that is not needed; this does not change the region
  above the kept hemispheres, since every hemisphere tested is a bubble.
  
  Let us show that if the algorithm
  terminates, then the answer is correct. If the algorithm terminates at Step
  \ref{item:first-solve}, it means that not all of the inequalities
  can hold, and so there is no point $P$ for which
  $r_0^2-d(P_0,P)^2 \ge r_i^2-d(P_i,P)^2$ for all $i > 0$, which means that
  there is no point above which hemisphere $H_0$ is higher than all of
  $H_1, \dots, H_k$. If the algorithm terminates at Step~\ref{item:under},
  then $r_0^2-d(P_0,P)^2 \ge r_i^2-d(P_i,P)^2$ for all $i$ implies that
  $r_0^2-d(P_0,P)^2 = r_j^2 -d(P_j,P)^2$, so again there is no point
  above which $H_0$ is strictly higher than all of $H_1, \dots, H_k$.

  If the algorithm terminates at Step~\ref{item:not-under}, then for
  every constraint we have an element of $L$ where that constraint is
  strictly satisfied, and all the rest are satisfied with possible
  equality. So the average $y$ of the elements of $L$ satisfies all
  inequalities strictly, and $C_0$ has nonempty interior. If $\delta < r_0^2$,
  choose $Q \in C_0$ with $d(P_0,Q) < r_0$. Every point of the segment from $Q$
  to $y$ other than $Q$ satisfies all inequalities strictly, and those close
  to $Q$ are at distance less than $r_0$ from $P_0$; above such a point $H_0$
  is strictly higher than all of $H_1, \dots, H_k$. If $\delta \geq r_0^2$, then
  every point $P$ with $d(P_0,P) < r_0$ lies outside $C_0$, so
  $r_i^2-d(P_i,P)^2 > r_0^2-d(P_0,P)^2$ for some $i$, and $H_0$ lies under
  $\cup_{i=1}^k H_i$.
  $\vartriangle$
\end{algorithm}

\begin{remark}\label{rem:stopping}
  Algorithm~\ref{alg:under-spheres} decides whether one hemisphere is needed
  given finitely many others; it does not say when all bubbles have been
  found. The following a posteriori criterion does. Let $\mathcal{B}$ be a
  finite set of bubbles, and let $h_{\mathcal{B}}$ be the minimum over $x \in \overline{F}$ of
  the largest height above $x$ of a $\Gamma_\infty$-translate of an element of
  $\mathcal{B}$. (Since the squared height of a hemisphere is a concave function on
  the region where that hemisphere is highest, $h_{\mathcal{B}}$ is attained at a vertex of
  this decomposition of $\overline{F}$ and can be computed exactly.) Suppose that
  $h_{\mathcal{B}}>0$ and that every bubble of radius greater than $h_{\mathcal{B}}$ with center in
  $\overline{F}$ is a $\Gamma_\infty$-translate of an element of $\mathcal{B}$ or lies
  under the union of these translates. Then the region above $\overline{F}$ and
  above all $\Gamma_\infty$-translates of elements of $\mathcal{B}$ is $\overline{D}$.
  Indeed, $\Gamma_\infty$ acts on $V_n$ by Euclidean isometries, so it preserves
  radii, and a bubble of radius at most $h_{\mathcal{B}}$ lies nowhere above this region.
  Since $B_g$ has radius $1/\vert c \vert$, only the finitely many bubbles with
  $\vert c \vert < 1/h_{\mathcal{B}}$ and center in $\overline{F}$ need to be tested. If $h_{\mathcal{B}} = 0$,
  for example when $\overline{F}$ contains singular cusps, the criterion does not
  apply, and a search shows only that no further bubble with $\vert c \vert$ up to
  the bound searched is needed.
\end{remark}

\subsection{Generators of $\SL_2(\Ocal)$ (theoretical)}

Let $g=\begin{pmatrix} a & b \\ c & d \end{pmatrix}$.
We will continue using the notation 
$$B_g = \lbrace x \in \Hcal^{n+1,\Sat} \colon \vert c x +d \vert = 1\rbrace = B(-c^{-1}d) = B_{1/\vert c \vert}(-c^{-1}d).$$ 
We will let $B_g^+$ denote the open interior of the sphere, not including the boundary. 
We will let $B_g^-$ denote the exterior. 
We will let $\overline{B}_g^+ = \lbrace x \in \Hcal^{n+1,\Sat} \colon \vert cx+d \vert \leq 1 \rbrace$ denote the interior with the closure and $\overline{B}_g^- = \lbrace x \in \Hcal^{n+1,\Sat} \colon \vert cx + d \vert \geq 1 \rbrace $. 
We will use similar notation for the spheres $B(x)$ where $x \in \Vec(K)$.

\begin{lemma}
Suppose $c \neq 0$. We have $g(B_g^-) = B_{g^{-1}}^+$. Moreover, $g^{-1}(\infty)$ is the center of $B_g$ and $g(\infty)$ is the center of $B_{g^{-1}}$.
\end{lemma}
\begin{proof}
	We will work with spheres of the form $B(\mu^{-1}\lambda)$ and take 
	$g = \begin{psmallmatrix}-\beta  & \alpha \\
	-\mu & \lambda
	\end{psmallmatrix} \in \SL_2(\Ocal)$
	so that $B(\mu^{-1}\lambda) = B_g$.
	A point $x$ lies in $g(\overline{B}_g^+)$ if and only if $y = g^{-1}x$ satisfies $\vert - \mu y + \lambda\vert \leq 1$.
	This implies that $g(\overline{B}_g^+)$ is defined by the inequality $\vert - \mu(a x + b)(c x+ d)^{-1} + \lambda\vert \leq 1,$ which simplifies to  $ \vert -\mu (a x + b) + \lambda(cx+d)\vert \leq \vert cx + d\vert $.
	Here we used
	$\begin{psmallmatrix}
	a & b \\
	c & d
	\end{psmallmatrix} = \begin{psmallmatrix}-\beta  & \alpha \\
	-\mu &  \lambda \end{psmallmatrix}^{-1} =
	\begin{psmallmatrix}
	\lambda^*&  -\alpha^*\\
	\mu^* & -\beta^*
	\end{psmallmatrix} $
	to avoid excessive superscripts.
	We now need to plug these in:
	$$\vert -\mu \lambda^* x + \mu \alpha^* + \lambda(\mu^* x-\beta^*)\vert \leq \vert \mu^* x - \beta^*\vert. $$
	The left hand side simplifies to  $ (\lambda \mu^* - \mu \lambda^*)x + \mu \alpha^* - \lambda \beta^* $.
	Since $cd^* = dc^*$ for $g$ (Proposition~\ref{prop:vn-equalities}), we have $\mu \lambda^* = \lambda\mu^*$, so the coefficient of $x$ is zero.
	Also, $\Delta(g) = -\beta \lambda^* + \alpha \mu^*$, hence the constant term is $\Delta(g)^*$, which is equal to $1$.
	This means the entire left hand side is just  $1$.
	The image of $\overline{B}_g^+$ under $g$ is then given by
	$$ 1 \leq \vert \mu^* x - \beta^*\vert, $$
	and $g^{-1}$ has bottom row $(\mu^*,-\beta^*)$. This implies that $g(\overline{B}_g^+)= \overline{B}_{g^{-1}}^-$; the same computation with equality gives $g(B_g) = B_{g^{-1}}$, and hence $g(B_g^-) = B_{g^{-1}}^+$.
	The claim is a direct computation using that $B_g$ is given by $\vert cx+d\vert =1$.
\end{proof}

In what follows, $B$ is the bubble domain of Definition~\ref{def:bubble-domain}. It is closed, so $\overline{B}=B$, and $B^{\circ}$ denotes its interior.\footnote{This notation follows Swan \cite{Swan1971}. }

\begin{lemma}\label{lem:side-pairing}
	If $B_g$ contains a side of $\overline{B}$, then $B_{g^{-1}}$ contains a side of $\overline{B}$.
\end{lemma}
\begin{proof}
By Theorem~\ref{T:characterization-of-bubble-domain},
	$B = \lbrace x \in \Hcal^{n+1} \colon \vert \gamma'(x) \vert \leq 1 \text{ for every } \gamma \in \Gamma \rbrace$,
	and $B^{\circ}$ is the set where the inequality is strict for every $\gamma\in\Gamma\setminus\Gamma_\infty$, that is, for every $\gamma$ with $c\neq 0$.
	Elements of $\Gamma_\infty$ have $\vert\gamma'(x)\vert=\vert d\vert^{-2}=1$ everywhere, so no strict inequality can be asked of them.
	For $\gamma$ with $c\neq0$, the sphere $B_\gamma$ is the locus $\vert\gamma'(x)\vert=\vert cx+d\vert^{-2}=1$.

	Take $x\in B\cap B_g$ and put $y=g(x)$.
	Since $\vert cx+d\vert=1$, the preceding lemma gives $y\in B_{g^{-1}}$, and the Magic Formula (Theorem~\ref{thm:magic}) gives $y_n=x_n$ and $\vert g'(x)\vert=1$.
	Let $\gamma\in\Gamma$. Since $\gamma g\in\Gamma$ and $x\in B$, we have $\vert(\gamma g)'(x)\vert\leq1$, and the chain rule gives $(\gamma g)'(x)=\gamma'(y)\, g'(x)$. Each of these derivatives is a scalar multiple of an orthogonal map (by the formula for $\gamma'(x)$ in the proof of Theorem~\ref{T:characterization-of-bubble-domain}), so their norms multiply and
	\[
	\vert\gamma'(y)\vert
	=\frac{\vert(\gamma g)'(x)\vert}{\vert g'(x)\vert}
	\leq1 .
	\]
	Equivalently, in terms of heights, $(\gamma(y))_n=((\gamma g)(x))_n\leq x_n=y_n$.
	Thus $y\in B$, and $g(B\cap B_g)\subseteq B\cap B_{g^{-1}}$.
	Applying the same argument to $g^{-1}$ gives the reverse inclusion, so $g(B\cap B_g)= B\cap B_{g^{-1}}$.
	Since $g$ is a homeomorphism of $\Hcal^{n+1}$ carrying $B_g$ onto $B_{g^{-1}}$, it carries a side of $B$ contained in $B_g$ onto a side of $B$ contained in $B_{g^{-1}}$.
\end{proof}

\begin{definition}
	Let $p\in V_n$, and fix $\overline{F}$, a compact closed fundamental domain for $\Gamma_{\infty}$.
	An \emph{additive reduction element} is an element $\sigma_p \in \Gamma_{\infty}$ such that $\sigma_p(p) \in \overline{F}$ (these are unique for points in the interior, and there are only finitely many on the boundary).
\end{definition}

The following Theorem tells us how the sides of the fundamental domain $D$ are related to generators of the group.
\begin{theorem}[Generators of $\PSL_2(\Ocal)$]\label{thm:generators-ocal}
	Suppose the sides of $D$ are given by the walls of $F$ together with bubbles $$B(p_1), \quad B(p_2), \quad  \ldots, \quad B(p_r).$$ 
	Let $g_j$ be such that $B_{g_j} = B(p_j)$ for $j=1,\ldots,r$ so that $p_j = g_j^{-1}(\infty)$.
	Let $q_j= g_j(\infty)$.
	The group $\PSL_2(\Ocal)$ is generated by the generators of $\PSL_2(\Ocal)_{\infty}$ together with 
	$$\sigma_{q_1}g_1, \quad \sigma_{q_2}g_2, \quad  \ldots \quad \sigma_{q_r}g_r. $$
\end{theorem}
\begin{proof}
	By Theorem~\ref{thm:generators-polyhedron}, $\Gamma$ is generated by the $\gamma$ for which $\overline{D}\cap\gamma\overline{D}$
	is a subdivided facet. Fix such a $\gamma$ and let $(c,d)$ be the bottom row of $\gamma^{-1}$.
	For $Q\in\overline{D}\cap\gamma\overline{D}$ both $Q$ and $\gamma^{-1}(Q)$ lie in $\overline{D}$, which is contained in the closed bubble domain,
	so they have the same height (Theorem~\ref{T:characterization-of-bubble-domain}\eqref{I:heights}); by the Magic Formula
	(Theorem~\ref{thm:magic}), $\vert cQ+d\vert=1$.
	If $c=0$, then $\gamma\in\Gamma_\infty$.
	If $c\neq 0$, then the $n$-dimensional set $\overline{D}\cap\gamma\overline{D}$ lies on the bubble $B_{\gamma^{-1}}$, which therefore carries
	a side of $D$, so $B_{\gamma^{-1}}=B(p_j)=B_{g_j}$ for some $j$. Comparing centres (the first lemma of this subsection),
	$\gamma(\infty)=p_j=g_j^{-1}(\infty)$, so $g_j\gamma\in\Gamma_\infty$ and $\gamma\in g_j^{-1}\Gamma_\infty$.
	Since $\sigma_{q_j}\in\Gamma_\infty$, every such $\gamma$ lies in the subgroup generated by $\Gamma_\infty$ and
	$\sigma_{q_1}g_1,\ldots,\sigma_{q_r}g_r$. In particular the choice of the $\sigma_{q_j}$ and of the $g_j$ does not matter.
\end{proof}

\subsection{Generators of $\SL_2(\Ocal)$ (algorithms)}\label{sec:tidy}
Throughout this subsection we fix
a Clifford algebra $K=\quat{-d_1,\dots,-d_n}{\QQ}$ and a Clifford-stable maximal order
$\Ocal$ in it. The essential issues are already visible in the case of
all $d_i$ equal to $1$, and the reader may prefer to simplify the notation
by restricting to that case.

Theorem~\ref{thm:generators-ocal} tells us that the generators of $\SL_2(\Ocal)$ come from the elements that allow us to cross the facets of our fundamental domains. 
Thus it remains to find the elements which allow us to cross the facets of our fundamental domain.

One element that is always present in all of our generating sets for
$\SL_2(\Ocal)$ is the matrix $\begin{psmallmatrix}0&1\\-1&0\end{psmallmatrix}$.
The utility of this one is obvious, even from an elementary point of view,
in view of the theory of continued fractions of rational numbers,  
but we also see it as giving the map that crosses the hemisphere $B(0)$,
which is always a facet of the boundary of the fundamental domain.

More generally, there are certain cusps $c \in \Vec(K) \cup \lbrace \infty \rbrace$ for which it is always possible to
write down a map that crosses the facet of the corresponding hemisphere.
In this section, we will describe the construction.  

\begin{definition}\label{def:tidy}
  Let $\lambda \in \Vec(\Ocal)$, $\mu \in \ZZ$, and suppose that
  $\bar\lambda \lambda \bmod \mu \in \{\pm 1\}$. Then the cusp
  $\lambda\mu^{-1}$ is called {\em tidy}. If $\lambda\mu^{-1}$ is
  tidy, then let $\varepsilon = 1$ if $\bar\lambda\lambda \equiv 1 \bmod \mu$
  and $\varepsilon = -1$ otherwise, and let $c$ be the integer such that
  $\bar\lambda\lambda - c\mu = \varepsilon$.
\end{definition}

\begin{remark} The fact that $\bar\lambda\lambda - c\mu = \varepsilon = \pm 1$
  shows that a tidy cusp is always unimodular.
\end{remark}

\begin{example} It is not difficult to give examples of tidy cusps. For
  example, in $\CC_4$ the cusp $(2+2i_1+2i_2+i_3)/6$ is tidy, because
  the norm of the numerator is $4+4+4+1=13\equiv 1 \bmod 6$. On the other
  hand, the cusp $(2+2i_1+2i_2+i_3)/5$ is not tidy.
\end{example}

Before describing how we can cross the hemisphere associated to a tidy cusp,
we introduce a bit of notation.

\begin{definition}\label{def:ms-tidy}
  If $s = \lambda \mu^{-1}$ is a tidy cusp, then let
  $s^\vee = -\varepsilon\bar\lambda \mu^{-1}$: it is also tidy. Further, let
  $H_s$ be the hemisphere with center $s$ and radius $1/\vert\mu\vert$, and
  similarly for $H_{s^\vee}$. Let $M_s$ be the matrix
  $\begin{psmallmatrix}-\bar\lambda&c \\ \mu&-\lambda\end{psmallmatrix}$ if
  $\varepsilon = 1$. If $\varepsilon = -1$,
  then change the sign of the first row and define $M_s$ to be
  $\begin{psmallmatrix}\bar\lambda&-c \\ \mu&-\lambda\end{psmallmatrix}$.
\end{definition}

The main point of this section is the following proposition and its
corollary:
\begin{proposition}\label{prop:cross-tidy-cusp}
  With notation as above, the matrix $M_s$ takes $H_s$ to $H_{s^\vee}$ and
  takes the interior of one hemisphere to the exterior of the other.
  More precisely, if $P$ is a point on $H_s$, let
  $P^\vee = s^\vee - \overline{(P - s)}$ (where we have extended $\bar{\cdot}$ to
  $\Hcal$ by letting it change the sign of the final coordinate, so that
  $P^\vee$ has the same height as $P$); then
  $M_s(P) = P^\vee$.
\end{proposition}

\begin{proof}
  Note first that this matrix satisfies the conditions of (\ref{E:conditions})
  and hence acts on the hyperbolic space.
  It is routine to verify the following equation:
  $$M_s = \begin{pmatrix}1&s^\vee\\0&1\end{pmatrix} \cdot
    \begin{pmatrix}0&-\mu^{-1}\\\mu&0\end{pmatrix} \cdot
      \begin{pmatrix}1&-\lambda \mu^{-1}\\0&1\end{pmatrix}.$$
  Considering the three factors in order, we have translation from
  $H_s$ to a hemisphere of radius $1/\vert\mu\vert$ centered at the origin,
  inversion in this hemisphere, and translation to $H_{s^\vee}$.
  The last statement is now clear, because inversion in a hemisphere
  centered at the origin takes a point $Q = y + t\,i_n$ to $-\bar y + t\,i_n$.
\end{proof}

\begin{remark}\label{rem:need-integers-for-sl2}
  This proposition is purely algebraic and is equally valid for
  $c, \mu \in \QQ$. However, if $c, \mu$ are not integers, we do not
  obtain an element of $\SL_2(\Ocal)$. We also note that
  $M_{s^\vee} = -M_s^{-1}$, so $M_{s^\vee}$ and $M_s^{-1}$ act in the same way.
\end{remark}

\begin{corollary}\label{cor:tidy-cusp}
  Suppose that $s$ is a tidy cusp such that part of $H_s$ is a facet $F$ of the
  boundary of the fundamental domain $D$.
  Then there is $T \in \Gamma_\infty$ such that the matrix $TM_s$ takes the
  copy $D_F$ of $D$ across $F$ to $D$.
\end{corollary}

\begin{proof} Let $P_s$ be a point in the relative interior of $F$. The
  hemisphere $H_s$ is $B_{M_s}$, and $M_s$ takes $H_s$ to $H_{s^\vee}$
  (Proposition~\ref{prop:cross-tidy-cusp}).
By the proof of Lemma~\ref{lem:side-pairing}, $M_s(P_s)$ lies in $H_{s^\vee} \cap \overline{B}$. Since
  $\overline{B}$ is invariant under $\Gamma_\infty$, there is $T \in \Gamma_\infty$
  such that $TM_s(P_s) \in \overline{D}$, and the corresponding part of
  $TH_{s^\vee}$ is a facet $F^\vee$ of the boundary of $D$.
  This done, if $P'_s$ is a point slightly beyond $P_s$ and in the
  copy of $D$ across $F$, then $TM_s(P'_s)$ is a point near $F^\vee$ but
  outside $TH_{s^\vee}$ and therefore in $D$. It follows that $TM_s$ takes
  one interior point of $D_F$ to a point of $D$; but this means it does
  the same for all such points.
\end{proof}

Here is a somewhat more general version of the construction of a reflection.
(We thank Daniel Martin for suggesting this idea to us.)
We note that the only point of $\Gamma\cdot\infty$ in $\overline{D}$ is $\infty$,
since every other point of $\Gamma\cdot\infty$ is the center of a bubble of positive radius.
\begin{proposition}\label{prop:copy-contains-s}
  Let $s$ be a unimodular cusp and let $D'$ be a translate of $D$ whose closure contains
  $s$. Let $M_s \in \SL_2(\Ocal)$ be a matrix taking $s$ to $\infty$.
  Fix $P \in D'$. Then there exists $T \in \Gamma_\infty$ such that
  $TM_s(P) \in D$.
\end{proposition}

\begin{proof}
  Since $M_s(s) = \infty$ and $s \in \overline{D'}$, it follows that
  $\infty \in M_s(\overline{D'})$. Let $T \in \SL_2(\Ocal)$ be such that
  $TM_s(D') = D$. Then $T(\infty)$ is a point of $\Gamma\cdot\infty$ in $\overline{D}$, so
  $T(\infty) = \infty$ and $T \in \Gamma_\infty$.
\end{proof}

In attempting to use Proposition~\ref{prop:copy-contains-s} to find the
generators of $\SL_2(\Ocal)$ from a description of the fundamental domain
as in Theorem~\ref{thm:generators-ocal}, there are
two difficulties. First, it may not be obvious how to write down
$M_s$ if $\Ocal$ does not admit a Euclidean algorithm and the cusp is not
given as $\lambda \mu^{-1}$ with $\mu \in \ZZ$ and
$\bar\lambda\lambda\equiv\pm1\bmod\mu$ (Definition~\ref{def:tidy}). Second, it is not necessarily clear how
to determine the cusp contained in a copy of $D$ adjacent to it.
It seems not to be true that if the facet of $D$ is a region $R$ of a
hemisphere with center $s$, then the cusp in the copy $D_R$ of $D$
lying across $R$ is necessarily $s$.

It appears that at least for some
maximal orders in $\CC_9$ the hemisphere above the untidy cusp,
$$\frac{3}{7}(1+i_1+i_2)+ \frac{2}{7}(i_3+i_4+i_5+i_6+i_7+i_8),$$
will contribute to the boundary of the fundamental
domain, and we do not know how to construct the corresponding reflection.

\subsection{Finite Index Subgroups and Passage to Suborders}\label{sec:finite-index-subgroups}

The purpose of this section is to explain how to pass to subgroups of finite index:
if $D$ is a fundamental domain for $\Gamma$ and $\Gamma' \subset \Gamma$ is of
finite index with representatives $g_1,\ldots,g_e \in \Gamma$ of the right cosets $\Gamma'\backslash\Gamma$, then 
$g_1 D \cup g_2 D \cup \cdots \cup g_e D$ is a fundamental domain for $\Gamma'$.  
In practice, this situation is actually quite complicated, and in \S\ref{sec:1-1-1} we pass to a finite-index suborder and pick up singular cusps.
Let $G$ act
transitively on $S$ and let $H \subseteq G$ be a subgroup of finite index.
(We formulate this abstractly, but in our applications $S=\Vec(K) \cup \lbrace \infty \rbrace$ and $G$ will be a group like $\PSL_2(\Ocal)$, and $H$ will be a finite index subgroup like $\PSL_2(\Ocal')$; then transitivity of $G$ on $S$ is the statement that $\Ocal$ is cuspidally principal (Definition~\ref{def:cpir}), which has to be checked.)

Let $x_1, \ldots, x_n$ be representatives for the orbits of $H$ on $S$.
Write $G_x$ (resp.~$H_x$) for the stabilizer of $x\in S$ in $G$ (resp.~$H$).
Let $c_{i,1},\ldots,c_{i,e_i}$ be coset representatives for $H_{x_i}\backslash G_{x_i}$.
Fix $x_0 \in S$ and let $g_1, \dots, g_n\in G$ be such that $g_i\cdot x_0 = x_i$ for $i=1,
\ldots n$. 

\begin{proposition}\label{prop:coset-reps-from-orbits} 
  The $c_{i,j}g_i$ form a set of coset representatives of $H\backslash G$:
  in other words, we have $G = \coprod_{i=1}^n \coprod_{j=1}^{e_i} H c_{i,j} g_i$.
  Furthermore, the index of $H$ in $G$ is $\sum_{i=1}^n [G_{x_i}:H_{x_i}]$.
\end{proposition}

\begin{proof} 
  First, we prove that distinct representatives give distinct cosets.
  Suppose that $Hc_{i,j}g_i = Hc_{i',j'}g_{i'}$. Then for some
  $h, h' \in H$ we have $hc_{i,j}g_i = h'c_{i',j'}g_{i'}$. Letting both sides
  act on $x_0$, we get that $hx_i = h' x_{i'}$ and so $i = i'$. But then
  $hc_{i,j} = h'c_{i,j'}$ and so $c_{i,j}c_{i,j'}^{-1} \in H \cap G_{x_i} = H_{x_i}$.
  Thus $j = j'$ by definition of the $c_{i}$.

  Now we show that every element of $G$ is in a coset of $H$
  represented by one of the $c_{i,j} g_i$.
  Fix $g \in G$; suppose that
  $gx_0$ is $H$-equivalent to $x_i$, and write $g = g'g_i$. Then
  we have $gx_0 = g'g_ix_0 = g'x_i$. Choose $h \in H$ with $hx_i = gx_0$:
  then $h(h^{-1}g')g_ix_0 = g'g_i x_0 = gx_0 = hx_i = hg_ix_0$. So
  $h$ takes $(h^{-1}g')g_ix_0$ and $g_ix_0$ to the same element of $S$,
  which means that $h^{-1}g' \in G_{x_i}$. Thus
  $h^{-1}g' = h'c_{i,j}$ for some $j$, and it follows that $g = hh'c_{i,j} g_i$,
  as desired.
  The last statement simply expresses the fact that the index of a subgroup is
  equal to the number of its cosets.
\end{proof}

Figure~\ref{fig:gamma0-2} shows the proposition for $P\Gamma_0(2)\subset\PSL_2(\ZZ)$.
\begin{figure}[htbp!]
\centering
\begin{tikzpicture}[scale=3.0,lab]
 \fill[domC!35] (-0.5,1.9)--(-0.5,0.866) arc[start angle=120,end angle=60,radius=1]--(0.5,1.9)--cycle;
 \fill[domC!20] (0.5,0.866) arc[start angle=60,end angle=120,radius=1] arc[start angle=60,end angle=0,radius=1] arc[start angle=180,end angle=120,radius=1];
 \fill[cuspC!15] (-0.5,0.2887)--(-0.5,0.866) arc[start angle=60,end angle=0,radius=1] arc[start angle=0,end angle=120,radius=0.3333];
 \draw[domC!80!black,line width=0.8pt] (-0.5,1.9)--(-0.5,0.866) (0.5,1.9)--(0.5,0.866);
 \draw[bub] (0.5,0.866) arc[start angle=60,end angle=120,radius=1];
 \draw[bub] (-0.5,0.866) arc[start angle=60,end angle=0,radius=1];
 \draw[bub] (0,0) arc[start angle=180,end angle=120,radius=1];
 \draw[bub] (0,0) arc[start angle=0,end angle=120,radius=0.3333];
 \draw[line width=0.8pt] (-0.5,0.2887)--(-0.5,0.866);
 \draw[bnd] (-1.0,0)--(1.0,0);
 \node at (0,1.4){$D$}; \node at (0,0.62){$SD$}; \node at (-0.36,0.42){$STD$};
 \fill[black] (0,0) circle (0.6pt) node[below]{$0$}; \node at (0,2.0){$\infty$};
 \node[below] at (-0.5,0){$-\tfrac12$}; \node[below] at (0.5,0){$\tfrac12$};
\end{tikzpicture}
\caption{Proposition~\ref{prop:coset-reps-from-orbits} for $H=P\Gamma_0(2)$ in $G=\PSL_2(\ZZ)$. The cusps $\infty$ and $0$ represent the $H$-orbits, with $[G_\infty:H_\infty]=1$ and $[G_0:H_0]=2$, so $[G:H]=3$ with coset representatives $1$, $S$, $ST$. The translates $D$, $SD$, $STD$ of the domain $D$ of $\PSL_2(\ZZ)$ together form a fundamental domain for $H$, with two tiles meeting at the cusp $0$.}\label{fig:gamma0-2}
\end{figure}

The following example shows that our proposition recovers the well-known example of the index of $\Gamma_0(p) \subset \SL_2(\ZZ)$.
\begin{example}\label{ex:gamma-0-p}
Let $p$ be prime, let $\Gamma_0(p) \subset \SL_2(\Z)$ be the usual
congruence subgroup consisting of matrices whose lower left entry is
divisible by $p$, and let $P\Gamma_0(p)\subset\PSL_2(\Z)$ be its image.
  The transitive action of $\PSL_2(\Z)$ on $\QQ \cup \lbrace \infty \rbrace$ is well-known; likewise,
the fact that $P\Gamma_0(p)$ acts on $\QQ \cup \lbrace \infty \rbrace$ with two orbits, represented by
  $\infty, 0$. The stabilizers of $\infty, 0$ in $\PSL_2(\Z)$ are generated by
  $\begin{psmallmatrix}1&1\\0&1\end{psmallmatrix}$ and $\begin{psmallmatrix}1&0\\1&1\end{psmallmatrix}$,
respectively. Thus the stabilizers in $P\Gamma_0(p)$ are of index $1$ and $p$.
  By Proposition~\ref{prop:coset-reps-from-orbits} we recover the well-known
fact that $[\PSL_2(\Z):P\Gamma_0(p)] = p+1$.
\end{example}

The following example was an important test case of Proposition~\ref{prop:coset-reps-from-orbits} for us.
\begin{example}\label{ex:sl2-z3}
  For an example more in keeping with the focus of this paper, let
  $\Ocal = \Z[\zeta_3]$, let $\Ocal' = \Z[\sqrt{-3}]$, and take
  $G = \PSL_2(\Ocal), H = \PSL_2(\Ocal')$. We will show that $[G:H] = 10$.
  Since $\Ocal$ is Euclidean, the group $G$ acts transitively on
  $\Q(\zeta_3) \cup \lbrace \infty \rbrace$.
  We claim that there are two orbits of cusps for $\Ocal'$,
  represented by $\infty$ and $\zeta_6 = \frac{1+\sqrt{-3}}{2}$.
  This follows from the
  Euclidean algorithm for $\Ocal$: every nonprincipal ideal of $\Ocal'$ is of the
  form $(2a,(1+\sqrt{-3})a)$ for some $a \in \Ocal$. The stabilizer of
  $\infty$ in $\SL_2(\Ocal)$ is generated by the three matrices
  $$S_1 = \begin{pmatrix}1&1\\0&1\end{pmatrix},\quad
    S_2 = \begin{pmatrix}1&\zeta_6\\0&1\end{pmatrix},\quad
    S_3 = \begin{pmatrix}\zeta_6&0\\0&\zeta_6^{-1}\end{pmatrix},$$
    with the relations $S_1S_2 = S_2S_1, S_3S_1S_3^{-1} = S_1^{-1} S_2,
    S_3S_2S_3^{-1} = S_1^{-1}, S_3^6 = 1$,
    while the stabilizer of $\zeta_6$ is conjugate to this by a matrix taking
    $\infty$ to $\zeta_6$, such as
    $\begin{psmallmatrix}\zeta_6&0\\1&\zeta_6^{-1}\end{psmallmatrix}$,
    and is therefore generated by
    $$T_1 = \begin{pmatrix}1-\zeta_6&\zeta_6-1\\-1&\zeta_6+1\end{pmatrix},\quad
    T_2 = \begin{pmatrix}-\zeta_6+2&-1\\-\zeta_6&\zeta_6\end{pmatrix},\quad
    T_3 = \begin{pmatrix}\zeta_6&0\\\zeta_6+1&-\zeta_6+1\end{pmatrix}$$
    with the same relations.
    We see that the subgroup of $\langle S_1, S_2, S_3 \rangle$ of matrices with
    entries in $\Ocal'$ is generated by $S_1, S_2^2, S_3^3$ and has index $6$,
    while the corresponding subgroup for the $T_i$ is generated by
    $T_1^2, T_2^2, T_3T_1$ and has index $4$.  (Both of these claims may
    be verified in \texttt{magma} by defining a finitely presented group
    using the relations above, checking that the given
    subgroups are in $\SL_2(\Ocal')$ and that they have the indices asserted,
    and verifying that the only coset representative in $\SL_2(\Ocal')$ is that
    of the identity.)
    Thus, according to Proposition~\ref{prop:coset-reps-from-orbits},
    the index is $10$.
\end{example}

\begin{remark} Although we do not need it for the example just given,
  we point out that presentations for $\SL_2(\ZZ[\zeta_3])$ and many other
  Bianchi groups may be found in \cite[Section 2.2]{FGT2010}.
\end{remark}
  
There are two sources of cusps for us. 
The first is arithmetic of ideal class sets, and the other is from singularities. 
These are discussed in the next two remarks (Remark~\ref{rem:cusps-come-from} and Remark~\ref{rem:cusps-from-singularities}).

Remark~\ref{rem:cusps-come-from} discusses cusps coming from ideal class sets and possible statements of Chebotarev density for orders of Clifford algebras.
\begin{remark}\label{rem:cusps-come-from}
  If $\Ocal$ is the maximal order in a number field $K$, then it is well-known
  \cite[Prop. 1.1]{vanderGeer1988}
  that the number of orbits for the action of $\PSL_2(\Ocal)$ on $K \cup \lbrace \infty \rbrace$ is
  $h_\Ocal$. (In this reference $K$ is assumed to be totally real, but that
  hypothesis is not used in the proof.) This is consistent with the fact
  that, by class field theory and the Chebotarev density theorem, the
  density of principal ideals among the prime ideals of $\Ocal$ is $1/h$.
  We do not know of a similar statement for maximal Clifford orders.

  \begin{samepage}
  To illustrate the situation, let us consider the
  maximal order $\Ocal$ in $\quat{-2,-3,-5}{\Q}$ generated by $a_1,a_2,a_3$
  and
  $$(1+a_2+a_3+a_2a_3)/2, \quad (a_1+a_1a_2a_3)/2, \quad (a_1a_2+a_1a_2a_3)/2,$$
  where $a_1^2=-2$, $a_2^2=-3$, and $a_3^2=-5$; its centre is $\ZZ[a_1a_2a_3]\cong\ZZ[\sqrt{30}]$.
  \end{samepage}
  Let $p>5$ be prime. If $\left(\frac{30}{p}\right) = 1$, then there are
  $2(p+1)$ maximal left ideals of $\Ocal$ of index $p^2$. These can never
  be generated by an element of the Clifford group, since the index of such
  a left ideal $\Ocal c$ is $(c\bar c)^4$, a fourth power (Remark~\ref{rem:cusp-map-not-injective}). On the other
  hand, there are $(p+1)^2+2(p^2+p+1)$ left ideals of index $p^4$, of which $(p+1)^2$
  arise from pairs of maximal left ideals in the quaternion algebras into
  which $\Ocal$ decomposes locally at $p$, and the other $2(p^2+p+1)$ from pairs
  consisting of a left ideal of index $p^4$ in one factor and the whole of the
  other factor; among these are $(p,1)$ and $(1,p)$.
  Exactly $(p+1)^2+2$ of these contain $p\Ocal$, namely the $(p+1)^2$ pairs and $(p,1)$, $(1,p)$.
  When $\left(\frac{30}{p}\right) = -1$, then there are no left
  ideals of index $p^2$ and $p^2+1$ of index $p^4$.
  
  To prove these statements, observe that if $\left(\frac{30}{p}\right) = 1$
  then $\Ocal$ decomposes over $\ZZ_p$ into a direct sum of two quaternion
  algebras (possibly split), while for $\left(\frac{30}{p}\right) = -1$
  it is a quaternion algebra over $\ZZ_p[\sqrt{30}]$, the unramified
  quadratic extension of $\ZZ_p$.
  
  Empirically, of the left ideals of index $p^4$ containing $p\Ocal$, the number
  of the form $\Ocal x$ with $x\in\Ocal$ is $26, 54, 74, 138, 154$ for
  $p = 7, 11, 13, 17, 19$, respectively (out of $66, 122, 198, 326, 402$),
  a proportion of roughly $2/5$ that slowly decreases for larger $p$
  (it is $322/902$ for $p=29$). It seems reasonable to guess that this is explained by a
  Chebotarev density theorem, perhaps by a positive answer to
  Question~\ref{conj:quaternion-chebotarev}.
\end{remark}

The next remark (Remark~\ref{rem:cusps-from-singularities}) talks about orders
``deep inside a maximal order''.
The intuition is that orders of large index inside a maximal order are singular and that these singularities produce singular points of the fundamental domain (note that there are two distinct meanings of the word ``singular'' in this sentence).
\begin{remark}\label{rem:cusps-from-singularities}
  Nonprincipal cusps also arise from the nonmaximality
  of orders. It is well-known, for example, that the ideal
  $(2,\sqrt{-3}+1) \subset \Z[\sqrt{-3}]$ is not principal but that the
  Picard group of $\Z[\sqrt{-3}]$ (the group of invertible fractional ideals)
  is trivial.  
  
  A similar phenomenon occurs in $\CC_4$. Let $\Ocal'$ be the
  nonmaximal order of $\CC_4$ generated by the Clifford units and let
  $\Ocal$ be the unique maximal order of $\CC_4$ containing $\Ocal'$
  (see \S\ref{sec:1-1-1}). Later we will show that $\Ocal$ is
  cuspidally principal in the sense of Definition~\ref{def:cpir}.
  The (left or right) ideal $(2,1+i_1+i_2+i_3)$ in $\Ocal'$ is not principal,
  but every ideal $I_p$ of index $p^4$ in $\Ocal'$, where $p$ is an odd prime,
  is principal. This can be proved by finding a generator $g_p$ of the extension
  of $I_p$ to $\Ocal$ and multiplying it by a unit of $\Ocal$ to obtain an
  element of $\Ocal'$, which will be a generator of $I_p$.
  More precisely, every left or right ideal of index $p^4$ becomes a free module
after tensoring with $\ZZ_p$ for any odd prime $p$, but $(2,1+i_1+i_2+i_3)$ remains
nonfree over $\ZZ_2$.

  Given a Clifford order $\Ocal$, by analogy with the classical situation we define
  the {\em Picard number} of $\Ocal$ to be the number of isomorphism classes of
  locally free $\Ocal$-modules of rank $1$. With this definition, we expect
  that the Picard number of $\Ocal'$ is still $1$.  
  On the other hand, we expect that the classical formula
  \cite[Prop. I.12.9]{Neukirch1999} expressing the class number
  of a quadratic order in terms of the class number of the maximal order, the unit
  index, and the decomposition of the level will have an analogue for Clifford orders,
  so that in all but a few exceptional cases (possibly including a small number of
  infinite families) the Picard number of a nonmaximal order
  will be larger than that of a maximal order containing it.
\end{remark}

\subsection{Changing the Order}
In this section we will indicate how to pass from one order to another
contained in it in practice.
In general, we feel that the most natural construction begins
from a maximal order, and we prefer to deduce a presentation for $\SL_2$
of a nonmaximal order from that of a maximal order rather than computing it
directly. Our approach is simply to find a lower bound for the index
by exhibiting distinct cosets and an upper bound by constructing matrices with
entries in the smaller order.

We begin with a basic fact of group theory. In the case where all the
subgroups are normal in $G$ this is part of the third isomorphism theorem,
but it is perhaps less familiar otherwise.

\begin{lemma}\label{lem:third} Let $G$ be a group with subgroups $H, K$
  such that $K \subseteq H$ and $K$ is normal in $G$. Then there is a
  canonical bijection $G/H \leftrightarrow (G/K)/(H/K)$.
\end{lemma}

\begin{proof} Simply take $aH$ to $a(H/K)$. Since $G/H$ and $(G/K)/(H/K)$
  have no natural group structure, there are no group operations to verify.
\end{proof}

We apply this in the following situation. Let $\Ocal$ be a $*$-stable
maximal order
in a Clifford algebra, let $R \subseteq \Ocal$ be another $*$-stable
order, and let
$\gamma_1, \dots, \gamma_n$ generate $\SL_2(\Ocal)$. Now let $G$ be the
free group on $n$ generators, let $K$ be the group of relations among the
$\gamma_i$ (that is, the kernel of the natural map $G \to \SL_2(\Ocal)$),
and let $H$ be the subgroup of $G$ consisting of words that map to matrices
in $\SL_2(R)$. Our goal is to find $[G/K:H/K]$, which by the lemma is equal
to $[G:H]$. We use the following algorithm:

\begin{algorithm}\label{alg:sl2-index}
  \begin{enumerate}
  \item Initialize two lists, one list $L_1$ of elements of $\SL_2(\Ocal)$ and
    one list $L_2$ of words in the generators of $\SL_2(\Ocal)$. Initially,
    both are empty.
  \item By day, search for elements $y \in \SL_2(\Ocal)$ such that
    $xy^{-1} \notin \SL_2(R)$ for all $x \in L_1$, and add them to $L_1$.
  \item By night, search for words in the generators whose product is in
    $\SL_2(R)$ and add them to $L_2$.
  \item Attempt to prove (using, say, $1$ second the first time and
    twice as much time each time as the previous)
    that the index of the group generated by $L_2$
    in $G$ is finite and of index equal to the length of $L_1$.
If this fails, return to the second step.
  \end{enumerate}
$\vartriangle$
\end{algorithm}

\begin{proposition}\label{prop:sl2-index-works}
  Algorithm~\ref{alg:sl2-index} is correct and always terminates.
\end{proposition}

\begin{proof} The length of $L_1$ is a lower bound for
  $[\SL_2(\Ocal):\SL_2(R)]$, while the index of the group generated by $L_2$
  is an upper bound, so if they are equal their common value must be equal
  to the index. Thus the algorithm cannot return an incorrect result.
  On the other hand, $[G:H] = [\SL_2(\Ocal):\SL_2(R)]$ is finite, since $\SL_2(R)$
  contains the principal congruence subgroup of level $m$ for any integer $m\geq 1$
  with $m\Ocal \subseteq R$, and then $H$ is finitely generated by Schreier's lemma.
  Provided the searches in the second and third steps eventually list every element
  of $\SL_2(\Ocal)$ and every word, $L_2$ eventually contains generators of $H$, and
  coset enumeration \cite{todd-coxeter} will eventually terminate and
  will give the index of $H$ in $G$.
\end{proof}

\begin{example}\label{ex:q-sqrt-minus3}
  Let $\Ocal = \ZZ[\frac{1+\sqrt{-3}}{2}]$ be the maximal order in
  $\QQ(\sqrt{-3})$ and let $R$ be the subring generated by $\sqrt{-3}$.
  For convenience, let $\zeta_6 = (1+\sqrt{-3})/2$. Having already
  proved that $\SL_2(R)$ has index~$10$ in $\SL_2(\Ocal)$, we now prove it
  again by means of Algorithm~\ref{alg:sl2-index}.

    In this case, the argument is essentially
    a slightly less motivated version of the first, but it
    is applicable in some situations where our control of the group theory
    is not as strong. First, we
consider $10$ matrices in distinct right cosets $\SL_2(R)M_i$ 
: $6$ upper triangular matrices,
  namely $M_1, M_2$ of the form $\begin{psmallmatrix}1&a\\0&1\end{psmallmatrix}$ for
  $a = 0, \zeta_6$ and $M_3, M_4,M_5,M_6$ given by
  $\begin{psmallmatrix}\zeta_6^i&j\\0&\zeta_6^{-i}\end{psmallmatrix}$ for $(i,j) \in \{1,2\} \times \{0,1\}$.
  In addition to these we have $M_7, M_8,M_9, M_{10}$ which are respectively
  $$\begin{pmatrix}-\zeta_6&-1 \\ 1&0 \end{pmatrix},\quad
  \begin{pmatrix}-1&\zeta_6\\ \zeta_6-1&0 \end{pmatrix},\quad
  \begin{pmatrix}-\zeta_6&1-\zeta_6 \\ -1&-1 \end{pmatrix},\quad
  \begin{pmatrix}1-\zeta_6&1-\zeta_6\\ -\zeta_6&0 \end{pmatrix}.$$

Since $M_i M_j^{-1} \notin \SL_2(R)$ for $i \ne j$,
  the index is at least $10$.

  On the other hand, since $\Ocal$ is Euclidean there is only one cusp, and
  hence $\SL_2(\Ocal)$ is generated by the stabilizer of infinity and
  $\begin{psmallmatrix}0&1\\-1&0\end{psmallmatrix}$. In turn, the stabilizer of infinity
  is generated by $\begin{psmallmatrix}1&i\\0&1\end{psmallmatrix}$ for an additive
  basis of $\Ocal$, say $1, \zeta_6$, and
$\begin{psmallmatrix}u&0\\0&u^{-1}\end{psmallmatrix}$ for generators $u$ of $\Ocal^\times$,
say $\zeta_6$.  The products of at most $7$ of these matrices generate
  a subgroup of index $10$ in the free group on $4$ generators.  Thus the
  index is $10$.
\end{example}
\section{The Cases $\QQ, \QQ(i)$ and $\QQ(\sqrt{-19})$}\label{sec:classical}

We briefly discuss the actions of $\PSL_2(\ZZ), \PSL_2(\ZZ[i])$ and $\PSL_2(\ZZ[\frac{1+\sqrt{-19}}{2}])$ since we have many new examples which are similar to these. 
There is nothing new in this section; it is included only for ease of reference and
so that the reader can compare the old cases to the new cases, as this is how we developed the theory.
\begin{itemize}
\item The example $\PSL_2(\ZZ)$ is the minimal example for our theory. 
\item The example $\Gamma=\PSL_2(\ZZ[i])$ is the minimal example with a nontrivial polyhedron determining the fundamental domain of $\Gamma_{\infty}$ on $V_n$.
The fundamental domain for $\Lambda=\Vec(\Ocal)$ is distinct from the fundamental domain for $\Gamma_{\infty}$.
\item The example $\PSL_2(\ZZ[\omega])$ for $\omega = (1+\sqrt{-19})/2$ is the minimal example, among fields of class number one, of a fundamental domain with a bubble $B(x)$, $x\in \Vec(K)$, of radius less than $1$. 
\end{itemize}

\begin{figure}[htbp!]
	\centering
	\begin{minipage}{0.47\textwidth}\centering\begin{tikzpicture}[scale=2.5, >=stealth]
  \pgfmathsetmacro{\h}{sqrt(19)/2}
  \pgfmathsetmacro{\hh}{sqrt(19)/4}
  \begin{scope}
    \clip (-0.5,0) rectangle (0.5,\h);
    \draw[thick, blue!70!black] (0,0) circle (1);
    \draw[thick, red!70!black] (0.5,\h) circle (1);
    \draw[thick, red!70!black, dashed] (-0.5,\h) circle (1);
    \draw[thick, green!45!black] (0.25,\hh) circle (0.5);
    \draw[thick, green!45!black, dashed] (-0.25,\hh) circle (0.5);
  \end{scope}
  \draw (-0.5,0) rectangle (0.5,\h);
  \draw[->] (-0.65,0) -- (0.7,0) node[right] {$x$};
  \draw[->] (0,-0.05) -- (0,\h+0.2) node[above] {$y$};
  \foreach \x in {-0.5,0.5} \draw (\x,0) -- (\x,-0.04) node[below, font=\scriptsize] {$\x$};
  \draw (-0.5,\h) -- (-0.54,\h) node[left, font=\scriptsize] {$\frac{\sqrt{19}}{2}$};
  \foreach \y in {1} \draw (-0.5,\y) -- (-0.54,\y) node[left, font=\scriptsize] {$\y$};
  \draw (-0.5,0.5) -- (-0.54,0.5) node[left, font=\scriptsize] {$0.5$};
\begin{scope}[font=\small]
    \draw[very thin, gray] (0.42,0.908) -- (0.62,0.80) node[right, black] {$B(0)$};
    \draw[very thin, gray] (0.45,1.548) -- (0.62,1.50) node[right, black] {$B(\omega/2)$};
    \draw[very thin, gray] (0.45,1.868) -- (0.62,1.95) node[right, black] {$B(\omega-1)$};
    \draw[very thin, gray] (-0.45,1.548) -- (-0.80,1.50) node[left, black] {$B(\tfrac{\omega-1}{2})$};
    \draw[very thin, gray] (-0.45,1.868) -- (-0.80,1.95) node[left, black] {$B(\omega)$};
  \end{scope}
\end{tikzpicture}\end{minipage}\hfill
	\begin{minipage}{0.5\textwidth}\centering\includegraphics[width=\linewidth]{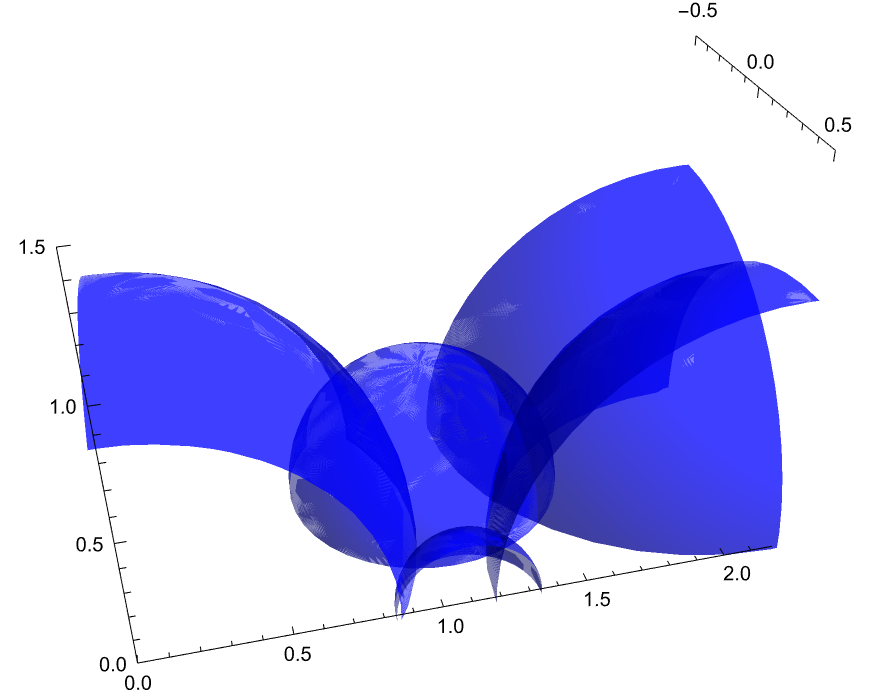}\end{minipage}
	\caption{The bubbles of $\QQ(\sqrt{-19})$ above the set of $x+iy$ with $-\frac12\le x\le\frac12$ and $0\le y\le\frac{\sqrt{19}}{2}$.
		Left: the circles where the five bubbles meet the boundary.
		Right: an oblique view of the bubbles in $\Hcal^3$.
	}\label{F:sqrt-19}
\end{figure}

Figure~\ref{F:sqrt-19} shows two images of the bubbles for $\PSL_2(\ZZ[\omega])$ where $\omega = (1+\sqrt{-19})/2$. 
The first image is where the spheres meet the boundary and is a collection of circles. 
The second image is a bird's-eye view of the bubbles. 
The first image is the type of image we will reproduce for $\Hcal^4$. 
Beyond $\Hcal^4$ it is only possible to do stereographic projections of certain polyhedra in $V_n$. 

\subsection{The Case of $\PSL_2(\ZZ)$}\label{sec:psl-2-z}
\label{ex:psl-2-z}
	For $\PSL_2(\ZZ)$ we let $E$ be the upper half-plane and
	$D$ the usual fundamental domain bounded by $x^2 + y^2 = 1$ and $x = \pm 1/2$.
	There is only one bubble $B(0)$, and $F = [-1/2,1/2]$ is the fundamental domain for $\Lambda=\ZZ$.
	
	For concreteness, let $F_1, F_2, F_3$ respectively be segments of the line
	$x = 1/2$, the circle $x^2 + y^2 = 1$, and the line $x = -1/2$. In terms of
	matrices we then have
	$$\tau_1 = \begin{pmatrix}1&1\cr0&1\cr\end{pmatrix},\quad
	S = \begin{pmatrix}0&1\cr-1&0\cr\end{pmatrix},\quad
	\tau_{-1} = \begin{pmatrix}1&-1\cr0&1\cr\end{pmatrix}.$$
The domain $D$, its sides and the side pairings are shown in Figure~\ref{fig:psl2z-domain}.
\begin{figure}[htbp!]
\centering
\begin{tikzpicture}[scale=2.6,lab]
 \fill[domC!35] (-0.5,2.1)--(-0.5,0.866) arc[start angle=120,end angle=60,radius=1]--(0.5,2.1)--cycle;
 \draw[gray!50] (1,0) arc[start angle=0,end angle=180,radius=1];
 \draw[gray!50] (-0.5,0.866)--(-0.5,0) (0.5,0.866)--(0.5,0);
 \draw[gray!50] (0,0) arc[start angle=0,end angle=180,radius=1] (2,0) arc[start angle=0,end angle=180,radius=1];
 \draw[domC!80!black,line width=1pt] (-0.5,2.1)--(-0.5,0.866) node[pos=0.4,left]{$F_3$};
 \draw[domC!80!black,line width=1pt] (0.5,2.1)--(0.5,0.866) node[pos=0.4,right]{$F_1$};
 \draw[bub,line width=1pt] (0.5,0.866) arc[start angle=60,end angle=120,radius=1];
 \node[bubbleC,below] at (0.2,0.98){$F_2$};
 \dotat[holeC]{0.5,0.866} \node[right] at (0.5,0.8){$\zeta_6$};
 \dotat[holeC]{-0.5,0.866} \node[left] at (-0.5,0.8){$\zeta_3$};
 \fill (0,1) circle (1.2pt) node[above]{$i$};
 \draw[-{Stealth},thick] (-0.5,1.75) to[bend left=35] node[above]{$\tau_1$} (0.5,1.75);
 \draw[-{Stealth},thick] (-0.25,1.06) to[bend left=50] node[above=-1pt]{\scriptsize $S$} (0.25,1.06);
 \draw[bnd] (-1.5,0)--(1.5,0);
 \draw[domC!80!black,line width=2.5pt] (-0.5,0)--(0.5,0) node[below,midway,black]{$F$};
 \foreach \c in {-1,0,1}{\dotat{\c,0}} \node[below] at (-1,0){$-1$};\node[below] at (1,0){$1$};
 \node at (0,1.42){$D$};
\end{tikzpicture}
\caption{The fundamental domain $D$ (orange) of $\PSL_2(\ZZ)$ from \S\ref{sec:psl-2-z}, bounded by $F_1=\{x=1/2\}$, $F_2=\{|x|=1\}$ (the bubble $B(0)$) and $F_3=\{x=-1/2\}$. The translation $\tau_1$ carries $F_3$ to $F_1$ and $S$ carries $F_2$ to itself, fixing $i$. The vertices $\zeta_3,\zeta_6$ are $\Gamma$-equivalent; the cycle at $\zeta_6$ gives the relation $(\tau_1S)^3=e$. Grey: neighbouring translates of $D$.}\label{fig:psl2z-domain}
\end{figure}
	
	The relations of Corollary~\ref{cor:switch-relations} are that
	$\tau_1 \tau_{-1} = \tau_{-1} \tau_1 = S^2 = e$. The fundamental
	domain itself has two vertices, namely $\zeta_6$ and $\zeta_3$, but these
	are $\Gamma$-equivalent so we need only consider one of them, say $\zeta_6$.
One easily sees from Figure~\ref{fig:psl2z-domain}
	that the associated word relation is $(\tau_1 S)^3$. We conclude
	that $\PSL_2(\ZZ) = \langle \tau_1, S \rangle/(S^2, (\tau_1 S)^3)$.
	Replacing $\tau_1$ by $\tau_1 S$ as a generator gives the usual
	presentation of $\PSL_2(\ZZ)$.

\subsection{The Case $\QQ(i)$}\label{ex:psl-2-z-i}
	The group $\Gamma=\PSL_2(\ZZ[i])$ acts on the hyperbolic $3$-space consisting of
	triples $(x,y,z) \in \RR^3$ with $z > 0$ in the usual way.
	
	\subsubsection{Fundamental Domain}
\begin{figure}[htbp!]
		\centering
		\begin{minipage}{0.45\textwidth}\centering\includegraphics[width=0.85\linewidth]{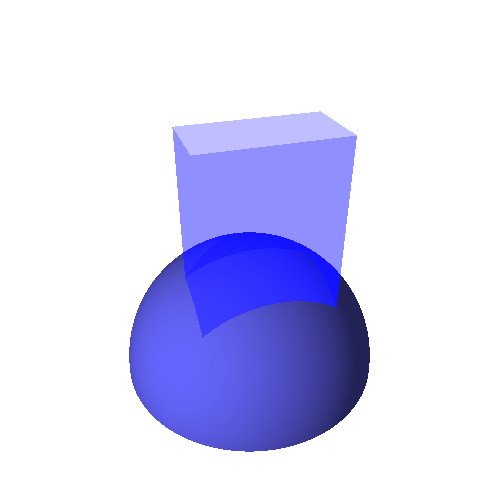}\end{minipage}\hfill
		\begin{minipage}{0.5\textwidth}\centering\begin{tikzpicture}[scale=2.0, >=stealth]
  \fill[black!12] (-0.5,0) rectangle (0.5,0.5);
  \draw[thick] (-0.5,0) rectangle (0.5,0.5);
  \draw[thin] (0,0) circle (1);
  \draw[->] (-1.25,0) -- (1.3,0) node[right] {$x_0$};
  \draw[->] (0,-1.2) -- (0,1.3) node[above] {$x_1$};
  \foreach \x in {-1,-0.5,0.5,1} \draw (\x,0.03) -- (\x,-0.03) node[below right=-1pt, font=\scriptsize, fill=white, inner sep=1pt] {$\x$};
  \foreach \y in {-1,0.5,1} \draw (0.03,\y) -- (-0.03,\y) node[above left=-1pt, font=\scriptsize, fill=white, inner sep=1pt] {$\y$};
  \node[font=\small] at (-45:1.28) {$|x|=1$};
\end{tikzpicture}\end{minipage}
		\caption{Left: the fundamental domain for $\Gamma=\PSL_2(\ZZ[i])$ acting on $\Hcal^3$: the region above the unit hemisphere over the rectangle $-\frac12\le x_0\le\frac12$, $0\le x_1\le\frac12$.
			Right: that rectangle, the fundamental domain for $\Gamma_{\infty}$ in $V_2$, together with the unit circle $|x|=1$ where the hemisphere meets $V_2$.
		}\label{F:h3-basic-fund}
	\end{figure}

	The fundamental domain for $\Ocal=\ZZ[i]$ acting on $\Hcal^3$ has a single sphere $B(0)$ of radius $1$ centered at zero above $F=\lbrace x_0 + i_1 x_1 \colon -1/2 < x_0 < 1/2, 0<x_1<1/2 \rbrace$.
	So there are five sides:
	$$ \lbrace x_0 = -1/2 \rbrace, \quad \lbrace x_0 = 1/2 \rbrace, \quad \lbrace x_1 = 0 \rbrace, \quad \lbrace x_1 =1/2 \rbrace, \quad \lbrace \vert x \vert =1 \rbrace.$$
	\cite[page 34]{Whitley1990}.
	An image of this fundamental domain is pictured in Figure~\ref{F:h3-basic-fund}.
	
\subsubsection{Generators}

The associated matrices are given by
$$\tau_1=\begin{pmatrix}1&1\\0&1\end{pmatrix},\quad 
\tau_{-1}=\begin{pmatrix}1&-1\\0&1\end{pmatrix}, \quad 
\gamma_3=\begin{pmatrix}i&1\\0&-i\end{pmatrix}, \quad
\pi_i=\begin{pmatrix}i&0\\0&-i\end{pmatrix},\quad
S=\begin{pmatrix}0&1\\-1&0\end{pmatrix}.$$
This may be checked by showing that they pair the respective facets
and take points inside the fundamental domain and close to them across
the facets. 

As before, Corollary~\ref{cor:switch-relations} confirms that
	$\tau_1 \tau_{-1} = \gamma_3^2 = \pi_i^2 = S^2 = e$.
	In addition, the $2$-faces (edges) induce the relations
	$$(\tau_1 \gamma_3)^2 = (\tau_1 \pi_i)^2 = (\tau_1 S)^3
= (\gamma_3 S)^3 = (\pi_i S)^2 = e.$$
In \cite[Sect. 2.3]{Sengun2012} we find the presentation
	$$\langle a,b,c,d|a^3,b^2,c^3,d^2,(ac)^2,(ad)^2,(bc)^2,(bd)^2 \rangle,$$
	which is transformed into ours by setting
	$a = \tau_1 S, b = S, c = S \gamma_3, d = S \pi_i$.

\subsection{The Case of $\QQ(\sqrt{-19})$} \label{S:sqrt19}
The case of the class number one field $K = \QQ(\sqrt{-19})$ with maximal order $\Ocal = \ZZ[\omega]$ where $\omega = \frac{1+\sqrt{-19}}{2}$ is a prototypical example of what we want to generalize. 
We let $\Gamma = \PSL_2(\Ocal)$.

Among fields of class number one, this is the first example in which there is a sphere whose radius is not
equal to $1$. The maximal orders in $\QQ(\sqrt{-m})$ for $m = 1, 2, 3, 7, 11$
are exactly the norm-Euclidean ones: every point of $\CC$ is at distance less than $1$ from $\Ocal$, and they give rise to simple fundamental domains where all
of the balls are centered at integers. The order $\Ocal$ of $\QQ(\sqrt{-19})$ has class number one but is not norm-Euclidean, and smaller spheres are needed.
The same happens for $m = 5, 6, 15$, for which the maximal order has class number $2$: the points $\sqrt{-5}/2$, $(1+\sqrt{-6})/2$ and $4\sqrt{-15}/15$ are at distance greater than $1$ from $\Ocal$, and they are the centers of bubbles of radius $1/2$, $1/2$ and $1/\sqrt{15}$, respectively.
All of these examples were investigated by Bianchi and Swan
\cite{Swan1971}. 

\subsubsection{Unit Group and Fundamental Domain for $\Gamma_{\infty}$}

The order $\Ocal$ has $\Ocal^{\times} = \lbrace \pm 1\rbrace$, so $\Gamma_{\infty} \cong \Vec(\Ocal) = \Ocal$ and the fundamental domain $F$ for the lattice $\Lambda = \Vec(\Ocal)$ is then $F=\lbrace x_0 + i x_1 \colon -1/2 \leq x_0 \leq 1/2, -\sqrt{19}/4 \leq x_1 \leq \sqrt{19}/4 \rbrace. $
The division by $4$ has to do with the congruence of 19 modulo $4$ and the shape of the basis of the $\Ocal$. 
By symmetry, it is enough to know the basis when $x_0\geq 0$ and $x_1\geq 0$, and we will often restrict to this case. 

\subsubsection{Fundamental Domain}
We show the spheres above the set
\[
  -\frac12\leq x_0\leq\frac12,\qquad
  0\leq x_1\leq\frac{\sqrt{19}}2.
\]
There are five bubbles involved:
 $$ B(0), \quad B(\omega), \quad B(\omega-1), \quad B(\omega/2), \quad B((\omega-1)/2).$$
The first three bubbles have radius $1$---the base points are elements of $\Ocal$. 
The remaining two bubbles have radius $1/2$. 
There are two images shown in Figure~\ref{F:sqrt-19}.
Figure~\ref{fig:sqrt19-top} shows the same bubbles from above.

\begin{figure}[htbp!]
\centering
\begin{tikzpicture}[scale=1.9,lab]
 \begin{scope}\clip (-0.5,0) rectangle (0.5,2.179);
  \fill[cuspC!30] (-0.5,0) rectangle (0.5,2.179);
  \foreach \p in {(0,0),(0.5,2.179),(-0.5,2.179),(1,0),(-1,0)}{\fill[bubbleC!12] \p circle (1); \draw[bub] \p circle (1);}
  \foreach \p in {(0.25,1.0897),(-0.25,1.0897)}{\draw[bub,dashed,line width=1pt] \p circle (0.5);}
 \end{scope}
 \draw[line width=0.8pt] (-0.5,0) rectangle (0.5,2.179);
 \fill (0,0) circle (1.3pt) node[below]{$0$};
 \fill (0.5,2.179) circle (1.3pt) node[above]{$\omega$}; \fill (-0.5,2.179) circle (1.3pt) node[above]{$\omega-1$};
 \fill[bubbleC] (0.25,1.0897) circle (1.2pt) node[above right=-1pt]{\scriptsize$\tfrac\omega2$};
 \fill[bubbleC] (-0.25,1.0897) circle (1.2pt) node[above left=-1pt]{\scriptsize$\tfrac{\omega-1}2$};
 \node[right] at (0.55,0.3){$x_0$ across, $x_1$ up};
\end{tikzpicture}
\caption{The case $\Ocal=\ZZ[\omega]$, $\omega=(1+\sqrt{-19})/2$ (\S\ref{S:sqrt19}), seen from above over $|x_0|\le\tfrac12$, $0\le x_1\le\sqrt{19}/2$. Blue disks: shadows of the unit bubbles $B(0)$, $B(\omega)$, $B(\omega-1)$ (and $B(\pm1)$). Red: points at distance at least $1$ from $\Ocal$; this is where $\Ocal$ fails to be norm-Euclidean. Dashed: the bubbles $B(\omega/2)$, $B((\omega-1)/2)$ of radius $\tfrac12$, which cover the red region.}\label{fig:sqrt19-top}
\end{figure}

The first shows the circles where the spheres meet the boundary $x_2=0$, and the other shows an oblique view of this selection of spheres. Note that when two circles intersect, their intersection points determine a line. 
Along these lines are exactly where the spheres have equal heights.
The segments where the two lines intersect are the projections of edges of the fundamental polyhedron---the edge is the edge of the two faces corresponding to the two bubbles. 

We remark that all of these cusps are tidy in the sense of \S\ref{sec:tidy}.

\section{The Case of $\quat{-1,-1}{\QQ}$} \label{S:hamilton}
For the quaternion algebra $\quat{-1,-1}{\QQ}$ we will use the usual quaternion notation with $i_1=i$ and $i_2=j$ with $ij=k$. 
We will consider two orders.
One is the order $\ZZ[i,j],$ which is also called the \emph{Lipschitz order} which has the $\ZZ[i,j]^{\times}= Q_8$, the quaternion group of order $8$. 
We also have the unique maximal order containing it called the \emph{Hurwitz Order}, which we denote by $\Ocal_3= \ZZ[i,j,\zeta],$ where 
$$\zeta =\frac{1+i+j+k}{2}.$$
Note that this makes $\Ocal_3 \subset \CC_3 \cong \RR^4$ a $D_4$ lattice while $\ZZ[i,j]$ is just a standard cubic lattice isomorphic to $\ZZ^4$.
On the other hand, in both of these cases we have $\Vec(O_3) = \Vec(\ZZ[i,j])\cong \ZZ^3$, which is just the standard cubic lattice.

\begin{proposition}\label{prop:lipschitz-euclidean} The Lipschitz order $\ZZ[i,j]$ and Hurwitz order $\Ocal_3$ are Clifford-Euclidean for the norm.
\end{proposition}
\begin{proof} 
	Both orders are $*$-stable, and $\Vec(\ZZ[i,j]) = \Vec(\Ocal_3) \cong \ZZ^3$ has covering radius $\sqrt{3}/2 < 1$, so this follows from Theorem~\ref{thm:small-covering-radius-euclidean}.
\end{proof}

\begin{remark}
It is important to note that $\PSL_2(\ZZ[i,j])$ can have class number 1 without $\ZZ[i,j]$ being a left principal ideal domain.  
The point is that, although $ \ZZ[i,j]\cdot (1 + i + j + k)+\ZZ[i,j]\cdot 2$ is not left principal, every cusp of the form $b^{-1}a$ where $a, b$ are Clifford vectors, i.e., of the form $r + s i+ t j$, is principal, that is, in the orbit of $\infty$.
\end{remark}

\begin{lemma}
	Let $\Ocal = \ZZ[i_1,\ldots,i_{n-1}]$.
	Let $\Gamma = \SL_2(\Ocal)$ act on $\Hcal^{n+1}$.
	The set $F\subset V_n$, defined as the set of $x=x_0 + x_1 i_1 + \cdots + x_{n-1} i_{n-1}$ such that  
	$$ -\frac{1}{2} < x_0 < \frac{1}{2}, \qquad 0 < x_j < \frac{1}{2}, \quad 1\leq j \leq n-1, $$
	is an open fundamental domain for the action of $\Gamma_{\infty}$ on $V_n$. 
\end{lemma}

\begin{proof}
	Clearly, $F_0 = \lbrace x \in V_n \colon -1/2 < x_j < 1/2\rbrace$ is an
	open fundamental domain for $\Vec(\Ocal)$. Every element of $\Ocal^{\times}$ has integer coordinates and norm $1$, so it is of the form $\pm i_{j_1}\cdots i_{j_r}$; hence $\Gamma_{\infty}$ acts on $V_n$ by translations by $\Vec(\Ocal)$ composed with products of the maps $T_j$ below. 
	Let $T_j(x) = i_j x i_j$. One has $T_j(i_j) = -i_j$, $T_j(i_k)=i_k$ for $0<k\neq j$, $T_j(1) = -1$.
	So $T_j(x)$ flips the sign of the $x_0$ and $x_j$ component.
	The map $T_j$ is a $180^{\circ}$ rotation in the $x_0x_j$-plane.
	The map $T_j$ has any number of fundamental domains in the $x_0x_j$-plane, but we can choose to always take $x_0$ to be long and $x_j$ to be short so that $-1/2<x_0<1/2$ and $0<x_j<1/2$. 
	We then intersect all of the fundamental domains for $T_j$ to get the result. 
	A picture of the fundamental domain is given in Figure~\ref{F:h3-stab-fund-dom}.
\begin{figure}[htbp!]
		\centering
		\begin{tikzpicture}[x={(1.5cm,0cm)}, y={(0.75cm,0.45cm)}, z={(0cm,1.5cm)},
    line join=round, >=stealth, lab/.style={font=\small}]
\newcommand{\hidden}{\draw[dashed, black!60] (-0.5,0.5,-0.5) -- (0.5,0.5,-0.5);
    \draw[dashed, black!60] (-0.5,0.5,-0.5) -- (-0.5,-0.5,-0.5);
    \draw[dashed, black!60] (-0.5,0.5,-0.5) -- (-0.5,0.5,0.5);}
  \newcommand{\visible}{\draw (-0.5,-0.5,-0.5) -- (0.5,-0.5,-0.5) -- (0.5,-0.5,0.5) -- (-0.5,-0.5,0.5) -- cycle;
    \draw (0.5,-0.5,-0.5) -- (0.5,0.5,-0.5) -- (0.5,0.5,0.5) -- (0.5,-0.5,0.5);
    \draw (-0.5,-0.5,0.5) -- (-0.5,0.5,0.5) -- (0.5,0.5,0.5);}
\newcommand{\subbox}[6]{\fill[black!22] (#1,#3,#5) -- (#2,#3,#5) -- (#2,#3,#6) -- (#1,#3,#6) -- cycle;
    \fill[black!12] (#1,#3,#6) -- (#2,#3,#6) -- (#2,#4,#6) -- (#1,#4,#6) -- cycle;
    \fill[black!32] (#2,#3,#5) -- (#2,#4,#5) -- (#2,#4,#6) -- (#2,#3,#6) -- cycle;
    \draw[thick] (#1,#3,#5) -- (#2,#3,#5) -- (#2,#3,#6) -- (#1,#3,#6) -- cycle;
    \draw[thick] (#1,#3,#6) -- (#1,#4,#6) -- (#2,#4,#6) -- (#2,#3,#6);
    \draw[thick] (#2,#3,#5) -- (#2,#4,#5) -- (#2,#4,#6);}
\begin{scope}[shift={(0.75,0,1.9)}]
    \hidden \visible
    \draw[->, very thin] (-0.85,-0.5,0.85) node[left, lab] {$x_2=\tfrac12$} -- (-0.1,0,0.5);
    \draw[->, very thin, dashed] (-0.85,-0.5,-0.95) node[left, lab] {$x_2=-\tfrac12$} -- (-0.1,0,-0.5);
    \draw[->, very thin, dashed] (-1.15,-0.5,0.1) node[left, lab] {$x_0=-\tfrac12$} -- (-0.5,-0.1,0.1);
    \draw[->, very thin] (1.15,0,0.1) node[right, lab] {$x_0=\tfrac12$} -- (0.5,0,0.1);
    \draw[->, very thin] (0.35,-0.5,-0.95) node[below right, lab] {$x_1=\tfrac12$ (front)} -- (0.15,-0.5,-0.2);
    \draw[->, very thin, dashed] (0.85,0.5,0.95) node[right, lab] {$x_1=-\tfrac12$ (back)} -- (0.1,0.5,0.3);
  \end{scope}
\begin{scope}[shift={(-0.9,0,-0.6)}]
    \hidden \subbox{-0.5}{0.5}{-0.5}{0.5}{0}{0.5} \visible
    \node[lab] at (0,-0.5,-0.85) {$T_2$: $0<x_2<\tfrac12$};
  \end{scope}
  \node at (0.35,0,-0.6) {$\cap$};
  \begin{scope}[shift={(1.6,0,-0.6)}]
    \hidden \subbox{-0.5}{0.5}{-0.5}{0}{-0.5}{0.5} \visible
    \node[lab] at (0,-0.5,-0.85) {$T_1$: $0<x_1<\tfrac12$};
  \end{scope}
  \node at (2.85,0,-0.6) {$=$};
  \begin{scope}[shift={(4.1,0,-0.6)}]
    \hidden \subbox{-0.5}{0.5}{-0.5}{0}{0}{0.5} \visible
    \node[lab] at (0,-0.5,-0.85) {both};
  \end{scope}
\end{tikzpicture}
		\caption{The construction of the fundamental domain for $\PSL_2(\ZZ[i,j])_{\infty}$ in $V_3$.
			Top: the box $[-\frac12,\frac12]^3$ with its faces labelled.
			Bottom: the fundamental domain for $T_2$ ($0<x_2<\frac12$) intersected with that for $T_1$ ($0<x_1<\frac12$) gives the fundamental domain $-\frac12<x_0<\frac12$, $0<x_1,x_2<\frac12$.
		}\label{F:h3-stab-fund-dom}
	\end{figure}
\end{proof}

\subsection{The Case of the Hurwitz Order in $\quat{-1,-1}{\QQ}$} \label{S:hurwitz}

There are some standard actions of the Hurwitz order $\Ocal_3^{\times}$ that are well-known. 
For example, this group acts on the tetrahedron embedded in the $i,j,k$ space. 
We will do something different, acting on the space $V_3=\RR+i\RR+j\RR \subset \HH$ rather than the usual subspace $\RR i + \RR j + \RR k \subset \HH$.

\subsubsection{Clifford Units}
	The Clifford group of units $\Ocal_3^{\times}$ has order $24$ and is generated by $i,j$ and $\zeta$, a 6th root of unity:
	 $$ \Ocal_3^{\times} = \langle i,j, \zeta\rangle, \quad \zeta = \frac{1+i+j+ij}{2}.$$
	 Note that while $\zeta$ is a Clifford group element, it can only be written as a product of Clifford vectors from $K=(-1,-1/\QQ)$ and not $\Vec(\Ocal_3)$.
	 We have $\zeta = (1+i)(1+j)/2,$ but we can't write $\zeta$ as a product of elements from $\Vec(\Ocal_3) \cap \Ocal_3^{\times}$.
	 We denote the subgroup generated by $\Vec(\Ocal_3) \cap \Ocal_3^{\times}$ by $\Vec(\Ocal_3)^{\times}$. 
	 This group is just the quaternion group
	  $$\Vec(\Ocal_3)^{\times} = Q_8 = \langle i,j\rangle,$$
	  and its image in $\Ocal_3^{\times}/\lbrace \pm 1 \rbrace \subset \SO(V_3)$ is the group $C_2^2$ generated by $\pi_i$ and $\pi_j$ which act as 
	 $$\pi_{i}(x+yi+zj) = (-x-yi+zj), \quad \pi_{j}(x+yi+zj) = -x+yi-zj.$$

	  Note that $\Vec(\Ocal_3) \cong \ZZ^3$, and that this group is actually the subgroup $\Weyl(\ZZ^3)^+$
          of the Weyl group $\Weyl(\ZZ^3)$ of the root system $\Phi= \lbrace \pm 1,\pm i,\pm j\rbrace$ in the Euclidean space $V_3\cong \RR^3$ with inner product $\langle x, y\rangle = 2x\cdot y$.

	 Here $C_2$ denotes a cyclic group of order 2.
	 
See Figure~\ref{fig:cubic3} for the Dynkin diagram of the simply laced root system of the root lattice $\ZZ^3$.
	 The full Weyl group would have all of the reflections $r_1,r_i,r_j$ where $r_{i_a}$ reflects the basis element $i_a$ and fixes $i_b$ for $a\neq b$.
	
	Under the isomorphism $C_2^3 \cong \FF_2^3$ we can view this subgroup of the Weyl group as a code which is sometimes useful (see e.g. \cite[p.~1]{Dupuy2023}).
	 Note that this group-code is not to be confused with our doubly even lattice-codes.
	 The group generated by $\pi_i$ and $\pi_j$ corresponds to the group generated by $101$ and $110$ in $\FF_2^3$. 
	 Here $101$ corresponds to $\pi_j$, and $110$ corresponds to $\pi_i$.
\begin{figure}[htbp!]
	 	\centering
	 	\begin{tikzpicture}[lab]
 \node[dnode,label={[dlab]below:$1$}] at (0,0) {};
 \node[dnode,label={[dlab]below:$i$}] at (\dstep,0) {};
 \node[dnode,label={[dlab]below:$j$}] at (2*\dstep,0) {};
\end{tikzpicture}
	 	\caption{The Dynkin diagram of the root system $\{\pm1,\pm i,\pm j\}$ of $\ZZ^3$, of type $A_1\times A_1\times A_1$.}\label{fig:cubic3}
	 \end{figure}
	
	The full quotient $\Ocal_3^{\times}/\lbrace \pm 1 \rbrace$ is isomorphic to $A_4$.
	We have $A_4 \cong C_2^2 \rtimes C_3$, and indeed the element $\zeta \in \Ocal^{\times}_3 \setminus \Vec(\Ocal_3)^{\times}$ comes from an automorphism of the root system $\Phi$ outside of $\Weyl(\Phi)$. 
	The special non-Clifford vector $\zeta = (1+i+j+k)/2$ acts by a signed cyclic permutation of the coordinates
	$$ \pi_{\zeta}(x+yi+zj) = -y-zi+xj.$$
	One sees that $\Aut(\Dynkin(\ZZ^3))=S_3$ and $C_3 = \langle \pi_{\zeta}\rangle \cong A_3$, which maps onto the subgroup of even permutations in $\Aut(\Dynkin(\ZZ^3))$.
	
\subsubsection{Fundamental Domain for $\Gamma_{\infty}$}
A fundamental domain for $\Gamma_{\infty}$ is pictured in Figure~\ref{fig:trace2-1-1}.
\begin{theorem}
	The fundamental domain for $\Gamma_{\infty}$ in $V_3$ can be described as the collection $F$ of $x+yi+z j$ where $z \in [0,1/2]$ and $0 \leq x \leq z$ and $0 \leq \vert y \vert \leq z$. 
\end{theorem}
\begin{proof}
The space of Clifford vectors is 3-dimensional, and we can translate to make all 3 coordinates have an absolute value less than $1/2$ (uniquely except on the boundary). 
We break up the positive octant $0 \leq x, y, z \leq 1/2$ into 6 pieces depending on the order of $x, y,z$ so one of the pieces is $0 \leq y \leq x \leq z \leq 1/2$, etc. We break up the others by taking the image of this triangulation by the Klein four-group maps generated by $\pi_{i_1}$ and $\pi_{i_2}$. 
So we have a set of 48 tetrahedra that is preserved by the action of $A_4$. The fundamental domain is a set of orbit representatives, ideally one such that the closures meet nicely in facets. One choice is the four tetrahedra:
$$\{0 \leq y \leq x \leq z \leq 1/2\}, \quad \{0 \leq x \leq y \leq z \leq 1/2\}$$
$$ \{0 \leq -y \leq x \leq z \leq 1/2\} \quad  \{0 \leq x \leq -y \leq z \leq 1/2\}.$$
One then checks directly that the union of these is precisely the region we have described.
\end{proof}

\begin{figure}[htbp!]
	\centering
	\begin{tikzpicture}[x={(-0.5299cm,-0.3177cm)}, y={(0.8480cm,-0.1985cm)}, z={(0cm,0.9272cm)}, scale=2.6, line join=round, >=stealth,
    c/.style={circle, fill, inner sep=1.2pt}]
  \coordinate (O) at (0,0,0);
  \coordinate (A) at (0,-0.5,0.5);
  \coordinate (B) at (0.5,-0.5,0.5);
  \coordinate (C) at (0.5,0.5,0.5);
  \coordinate (E) at (0,0.5,0.5);
\draw[dashed, black!50] (0,0,0) -- (-1.1,0,0);
  \draw[dashed, black!50] (0,0,0) -- (0,-1.15,0);
  \draw[dashed, black!50] (0,0,0) -- (0,0,-1.15);
\fill[orange!60, opacity=0.55] (O) -- (A) -- (E) -- cycle;
  \fill[orange!70, opacity=0.55] (O) -- (A) -- (B) -- cycle;
  \fill[orange!80!black, opacity=0.45] (O) -- (B) -- (C) -- cycle;
  \fill[orange!50, opacity=0.55] (O) -- (C) -- (E) -- cycle;
  \fill[orange!35, opacity=0.65] (A) -- (B) -- (C) -- (E) -- cycle;
  \draw[thick] (A) -- (B) -- (C) -- (E) -- cycle;
  \draw[thick] (O) -- (A) (O) -- (B) (O) -- (C) (O) -- (E);
\draw[->] (0,0,0) -- (1.35,0,0) node[anchor=north east] {$x$};
  \draw[->] (0,0,0) -- (0,1.35,0) node[anchor=north west] {$y$};
  \draw[->] (0,0,0) -- (0,0,1.35) node[anchor=south] {$z$};
\node[c, label={[font=\small]below:$0$}] at (0,0,0) {};
  \node[c, label={[font=\small]below:$1$}] at (1,0,0) {};
  \node[c, label={[font=\small]below:$i$}] at (0,1,0) {};
  \node[c, label={[font=\small]below:$-i$}] at (0,-1,0) {};
  \node[c, label={[font=\small]right:$j$}] at (0,0,1) {};
  \node[c, label={[font=\small]right:$-j$}] at (0,0,-1) {};
\end{tikzpicture}
	\caption{The fundamental domain $F=\{x+yi+zj : 0\le z\le\frac12,\ 0\le x\le z,\ |y|\le z\}$ of $\Gamma_{\infty}$ for $\Gamma=\PSL_2(\Ocal_3)$, a square pyramid with apex $0$.
		Also pictured are the centres $0,1,\pm i,\pm j$ of the bubbles $B(0),B(1),B(\pm i),B(\pm j)$.
	}\label{fig:trace2-1-1}
\end{figure}

\subsubsection{Fundamental Domain for $\Gamma$}

\begin{theorem}\label{thm:one-bubble-o3}
The fundamental domain for $\PSL_2(\Ocal_3)$ acting on $\Hcal^4$ is the region above $F$ with the single bubble $B(0)$. 
\end{theorem}
\begin{proof}
	The region $F$ has a convenient description as $\vert y \vert \leq z$ and $0\leq x\leq z$ for $z\in [0,1/2]$. This is pictured above with centers for the bubbles $B(0),B(1),B(i),B(-i),B(j),B(-j)$. 
	These are the only bubbles that could matter. 
	
	The 4-cell bordering $B(0)$ and $B(1)$ sits above a 3-cell in $\lbrace x=1/2\rbrace$ which is on the boundary of $F$. This is the same for the other ones. They project to $\lbrace x=\pm 1/2 \rbrace$, $\lbrace y = \pm 1/2\rbrace$ and $\lbrace z = \pm 1/2\rbrace$. 
	None of these are on the interior of $F$ and hence we can omit them.

We will now check the spheres of other radii. 
First, we consider $B((1+i)/2)$. Since $(1+i)/2 = 1/(1-i)$ in lowest terms,
the radius is $1/\sqrt{2}$.
Again, consider a point $(1/2+a,1/2+b,c)$, where $a, b < 0$.  
The point above this in $B(0)$ has fourth coordinate squaring to $1/2 - (a+b) - (a^2+b^2+c^2)$, and for $B((1+i)/2)$ it is $1/2 - (a^2+b^2+c^2)$: also smaller.

After that we need to look at $(1+i+j)/3 = 1/(1-i-j)$. The radius squared is $1/3$. We have to be more careful, because there is less symmetry.  
Nevertheless, we claim that a sphere is not needed here.
Consider a point lying under this hemisphere whose coordinates are $$(a+1/3,b+1/3,c+1/3,\sqrt{1/3-(a^2+b^2+c^2)}).$$  
We first compare to the hemisphere centered at the origin, which includes the point $$(a+1/3,b+1/3,c+1/3,\sqrt{1-((a+1/3)^2+(b+1/3)^2+(c+1/3)^2)}).$$  
Squaring the last coordinates and expanding, we have to compare $1/3 - (a^2+b^2+c^2)$ to $2/3 - 2/3(a+b+c) - (a^2+b^2+c^2)$.  
After cancelling the square terms, we find that the point on the basic hemisphere centered at $0$ is at least as high if and only if $1/2 \ge a+b+c$.  
This does not follow from the assumption that $a^2+b^2+c^2 \le 1/3$, so
we consider other hemispheres as well, following Algorithm~\ref{alg:under-spheres}.
In particular, with the hemisphere with center $(1,1,1)$
we compare $1/3 - (a^2+b^2+c^2)$ to 
$$1 - (2/3-a)^2 - (2/3-b)^2 - (2/3-c)^2 = 1 - 4/3 + 4/3(a+b+c) - (a^2+b^2+c^2).$$  
The right-hand side is greater or equal if and only if $a+b+c \ge 1/2$.
Of course, either this condition or the previous condition $1/2 \ge a+b+c$
must hold. (In terms of our linear programming, although
$1/2 \ge a+b+c, a+b+c \ge 1/2$ is a feasible set of inequalities, it is
not possible for both of them to hold strictly.)
It follows that the hemisphere based at $(1+i+j)/3$ is dominated by the pair of hemispheres with centers $0$ and $1+i+j$.  
By symmetry it is unnecessary to consider $(\pm 1 \pm i \pm j)/3$.

If the denominator has norm 4, the radius is $1/2$. Either the center is
one of $1/2, i/2, j/2$
and Lemma~\ref{lem:singly-dominated} applies with $P = 0$
or it is $(\pm 1 \pm i \pm j)/2$. It suffices to treat the positive signs;
since $(1+i+j)/2$ is in lowest terms the radius is $1/2$.
Consider a point $(1/2+a,1/2+b,1/2+c) \in {\mathbb R}^3$;
first, suppose that $a, b, c < 0$.  
Then the fourth coordinate of the point of $B(0)$ above it squares to
$1/4 - (a+b+c) - (a^2+b^2+c^2)$, and of $B((1+i+j)/2)$ to
$1/4 - (a^2+b^2+c^2)$, which is smaller. Similarly, with a different set
of signs we would choose $x, y, z \in \{0,1\}$ such that $x-1/2, y-1/2, z-1/2$
have the same signs as $a,b,c$ and use $B(x+yi+zj)$ in place of $B(0)$.

If the denominator has norm 5, we can take the center to be
$(x+yi+zj)/5$ where $5|(x^2+y^2+z^2)$ and $x, y, z$ are $0, 1, 2$. Up
to symmetry there is only one possibility, and the distance from
$(2+i)/5$ to the origin plus the radius $1/\sqrt 5$ is less than $1$,
so Lemma~\ref{lem:singly-dominated} applies.

If the denominator has norm 6, the center is $(x+yi+zj)/6$ where
$6|(x^2+y^2+z^2)$ and $x, y, z \in \{0, 1, 2, 3\}$. We can assume that
$(x,y,z) = 1$, so essentially the only possibility is $(2+i+j)/6$.
Again, Lemma~\ref{lem:singly-dominated} shows that this sphere is not needed.

If the denominator has norm 7, the center is $(x+yi+zj)/7$ where
$7|(x^2+y^2+z^2)$ and $x,y,z \in \{0,1,2,3\}$.
The only choice is $(3+2i+j)/7$, and again this is taken care of by
Lemma~\ref{lem:singly-dominated}.

Finally, we claim that no cusp whose denominator has norm $\ge 8$ can ever
be required. Indeed, every point of $F$ satisfies $x^2+y^2+z^2\leq 3/4$, so the point of $B(0)$ above it has height at least $\sqrt{1-3/4}=1/2$.
A hemisphere of radius at most $1/\sqrt{8}<1/2$ has height less than $1/2$ everywhere, so above $F$ it lies below $B(0)$
and contributes no side of the fundamental domain.
\end{proof}

\subsubsection{Generators}
		The generators are found using Theorem~\ref{thm:generators-ocal}, where there is a generator for each side of the fundamental domain $D$.
		They come from the bubbles and the generators of $\PSL_2(\Ocal_3)_{\infty}$.
		The group $\PSL_2(\Ocal_3)$ is generated by 
		$$S, \qquad \tau_1,\quad \tau_{i}, \quad \tau_j, \qquad \pi_{\zeta}, \quad \pi_{i}, \quad \pi_{j}.$$
	
\subsection{The Case of $\quat{-1,-1}{\ZZ}$} \label{sec:lipschitz}
This was first carried out in \cite[Theorem 8]{Maclachlan1989}.
A calculation similar to the calculation of the fundamental domain of $\ZZ[i]$ acting on $\Hcal^3$ can be performed for $\ZZ[i,j]$ acting on $\Hcal^4$.
The group of units here is the well-known quaternion group 
$$ \ZZ[i,j]^{\times}  =Q_8 = \langle \pm1 ,\pm i, \pm j, \pm k \rangle.$$

\subsubsection{Fundamental Domain}
The fundamental domain for $\Gamma' = \PSL_2(\ZZ[i,j])$ acting on $\mathcal{H}^4$ has as boundaries 
	$$\lbrace x = -1/2 \rbrace, \quad \lbrace x = 1/2 \rbrace, \quad \lbrace y = 0 \rbrace, \quad \lbrace y=1/2 \rbrace$$
	$$  \lbrace z=0 \rbrace, \quad \lbrace z =1/2 \rbrace, \quad B(0)=\lbrace \vert x \vert =1 \rbrace.$$
The box $F$ and the heights of $B(0)$ over its corners are shown in Figure~\ref{fig:lipschitz-domain}.
\begin{figure}[htbp!]
\centering
\begin{tikzpicture}[scale=3.2,lab,x={(1cm,0cm)},y={(0.55cm,0.4cm)},z={(0cm,1cm)}]
 \foreach \a in {-0.5,0.5}{\foreach \b in {0,0.5}{\foreach \c in {0,0.5}{\coordinate (p\a\b\c) at (\a,\b,\c);}}}
 \fill[domC!25] (-0.5,0,0)--(0.5,0,0)--(0.5,0,0.5)--(-0.5,0,0.5)--cycle;
 \fill[domC!35] (0.5,0,0)--(0.5,0.5,0)--(0.5,0.5,0.5)--(0.5,0,0.5)--cycle;
 \fill[domC!15] (-0.5,0,0.5)--(0.5,0,0.5)--(0.5,0.5,0.5)--(-0.5,0.5,0.5)--cycle;
 \draw[domC!80!black,line width=0.8pt] (-0.5,0,0)--(0.5,0,0)--(0.5,0.5,0)--(0.5,0.5,0.5)--(-0.5,0.5,0.5)--(-0.5,0,0.5)--cycle
   (-0.5,0,0.5)--(0.5,0,0.5)--(0.5,0.5,0.5) (0.5,0,0)--(0.5,0,0.5);
 \draw[domC!80!black,dashed] (-0.5,0,0)--(-0.5,0.5,0)--(0.5,0.5,0) (-0.5,0.5,0)--(-0.5,0.5,0.5);
 \fill (0,0,0) circle (0.5pt) node[below]{$0$: height $1$};
 \node[below left] at (-0.5,0,0){$\tfrac{\sqrt3}{2}$}; \node[below right] at (0.5,0,0){$\tfrac{\sqrt3}{2}$};
 \node[above right] at (0.5,0.5,0.5){$\tfrac12$}; \node[above left] at (-0.5,0.5,0.5){$\tfrac12$};
 \node[left] at (-0.5,0,0.5){$\tfrac1{\sqrt2}$};
 \draw[->,gray] (0.7,0,0)--(0.95,0,0) node[right]{$x_0$};
 \draw[->,gray] (0.7,0,0)--(0.7,0.35,0) node[right]{$x_1$};
 \draw[->,gray] (0.7,0,0)--(0.7,0,0.3) node[above]{$x_2$};
\end{tikzpicture}
\caption{The fundamental domain for $\PSL_2(\ZZ[i,j])$ acting on $\Hcal^4$ (\S\ref{sec:lipschitz}). Orange: the box $F=\{|x_0|\le\tfrac12,\ 0\le x_1,x_2\le\tfrac12\}\subset V_3$. The domain $D$ is the region above $F$ and outside $B(0)$; numbers are the height $\sqrt{1-|x|^2}$ of $B(0)$ at the corners of $F$. The neighbouring bubbles $B(\pm1)$ meet $B(0)$ over the faces $x_0=\pm\tfrac12$.}\label{fig:lipschitz-domain}
\end{figure}
Note that we don't include $B(1)$ and $B(-1)$ (or other translates by the cubic lattice) because they meet $B(0)$ at $x=1/2$ and $x=-1/2$, respectively.
The image of the fundamental domain for $\Gamma_{\infty}$ acting on $V_3$ was already pictured in Figure~\ref{F:h3-stab-fund-dom}.

\begin{remark}
	This example is deceptively simple. 
	The same naive example in higher dimensions becomes wildly complicated.
	We cannot expect the domain bounded by
	$|x_i| \le 1/2$ and $\sum_{i=1}^n x_i^2 + y^2 = 1$ to be fundamental for
	$n \geq 4$, because $\sum_{i=1}^n x_i^2$ can be at least $1$
	(for $n=4$ the corner $(1+i_1+i_2+i_3)/2$ lies on $\vert x \vert = 1$). Thus we would not expect a set of generators analogous to those
	above to be sufficient to generate $\PSL_2(\ZZ[i_1,i_2,i_3])$.
\end{remark}
	
\subsubsection{Generators}\label{sec:lipschitz-generators}
The generators for $\PSL_2(\ZZ[i,j])$ are the matrices 
	$$
	\tau_1 = \begin{pmatrix}
	1 & 1 \\
	0 & 1 \\
	\end{pmatrix},
	\tau_1^{-1} =\begin{pmatrix}
	1 & -1 \\
	0 & 1 \\
	\end{pmatrix},
	\gamma_3=\begin{pmatrix}
	i & 1 \\
	0 & -i \\
	\end{pmatrix},
	\pi_i=\begin{pmatrix}
	i & 0 \\
	0 & -i \\
	\end{pmatrix},$$
	$$
	\gamma_5=\begin{pmatrix}
	j & 1 \\
	0 & -j \\
	\end{pmatrix},
	\pi_j=\begin{pmatrix}
	j & 0 \\
	0 & -j \\
	\end{pmatrix},
	S=\begin{pmatrix}
	0 & 1 \\
	-1 & 0 \\
	\end{pmatrix}.
	$$
Besides the switch relations $S^2=\gamma_3^2=\pi_i^2=\gamma_5^2=\pi_j^2=1$ of Corollary~\ref{cor:switch-relations}, the relations, which we computed using our algorithm, are rather complicated but given by 
	\begin{equation*}
\begin{aligned}
& \tau_1\gamma_3\tau_1\gamma_3 = 1 & \quad & \tau_1\pi_i\tau_1\pi_i = 1 & \quad & \tau_1\gamma_5\tau_1\gamma_5 = 1 \\
& \tau_1\pi_j\tau_1\pi_j = 1 & \quad & \tau_1S\tau_1S\tau_1S = 1 & \quad & \tau_1^{-1}\gamma_3\tau_1^{-1}\gamma_3 = 1 \\
& \tau_1^{-1}\pi_i\tau_1^{-1}\pi_i = 1 & \quad & \tau_1^{-1}\gamma_5\tau_1^{-1}\gamma_5 = 1 & \quad & \tau_1^{-1}\pi_j\tau_1^{-1}\pi_j = 1 \\
& \tau_1^{-1}S\tau_1^{-1}S\tau_1^{-1}S = 1 & \quad & \gamma_3\gamma_5\gamma_3\gamma_5 = 1 & \quad & \gamma_3\pi_j\gamma_3\pi_j = 1 \\
& \gamma_3S\gamma_3S\gamma_3S = 1 & \quad & \pi_i\gamma_5\pi_i\gamma_5 = 1 & \quad & \pi_i\pi_j\pi_i\pi_j = 1 \\
& \pi_iS\pi_iS = 1 & \quad & \gamma_5S\gamma_5S\gamma_5S = 1 & \quad & \pi_jS\pi_jS = 1. \\
\end{aligned}
	\end{equation*}
Writing $\tau_v=\begin{psmallmatrix}1&v\\0&1\end{psmallmatrix}$, one checks that $\pi_i=S^{-1}\tau_iS^{-1}\tau_i^{-1}S^{-1}\tau_i$, $\pi_j=S^{-1}\tau_jS^{-1}\tau_j^{-1}S^{-1}\tau_j$, $\gamma_3=\pi_i\tau_i^{-1}$ and $\gamma_5=\pi_j\tau_j^{-1}$. Hence the group generated by the translations $\tau_1,\tau_i,\tau_j$ and $S$, which is of the kind considered in \cite{Krausshar2004}, is all of $\PSL_2(\ZZ[i,j])$.
\section{The Case of  $\quat{-1,-3}{\QQ}$} \label{S:quat13}
We will view $\quat{-1,-3}{\QQ}$ as sitting inside Hamilton's quaternions, and so $\quat{-1,-3}{\QQ} = \QQ[i,\sqrt{3}j]$.
Inside $\quat{-1,-3}{\QQ}$ there are two maximal orders $\Ocal(-1,-3)_1$ and $\Ocal(-1,-3)_2$ which contain the elements $(1+\sqrt{3}j)/2$ and $(1+\sqrt{3}k)/2$, respectively. 
Note that one has a sixth root of unity which is a Clifford vector while the other does not. These orders are conjugate to each other by $1+i$. 
Also observe that we cannot generate a single order with both roots of unity in it, because $(1+\sqrt{3}j)/2 - (1+\sqrt{3}k)/2$ is not integral; it
has minimal polynomial $x^2 + 3/2$.

\subsection{The Case of  $\Ocal(-1,-3)_2$}
\begin{figure}[htbp!]
		\begin{center}
			\includegraphics[scale=0.33]{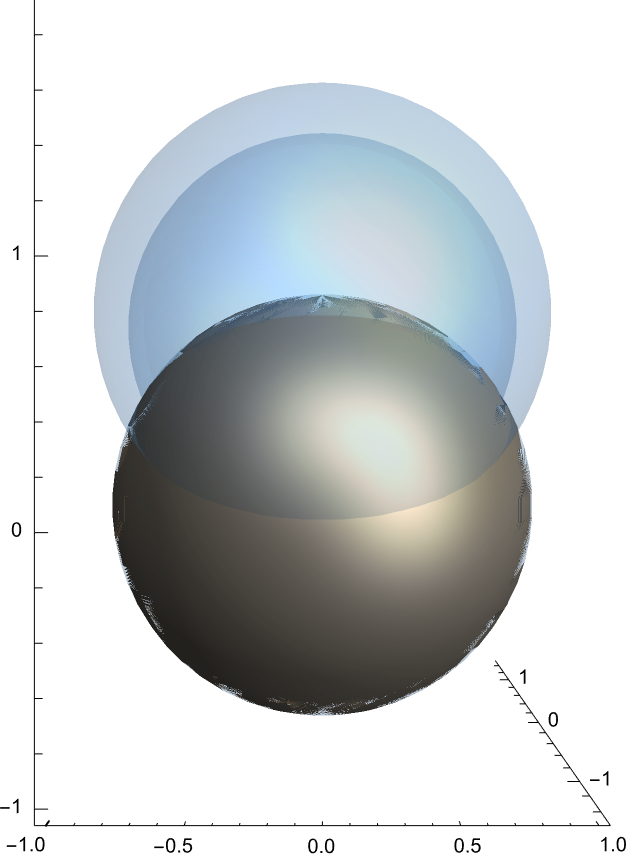} \qquad
			\includegraphics[scale=0.33]{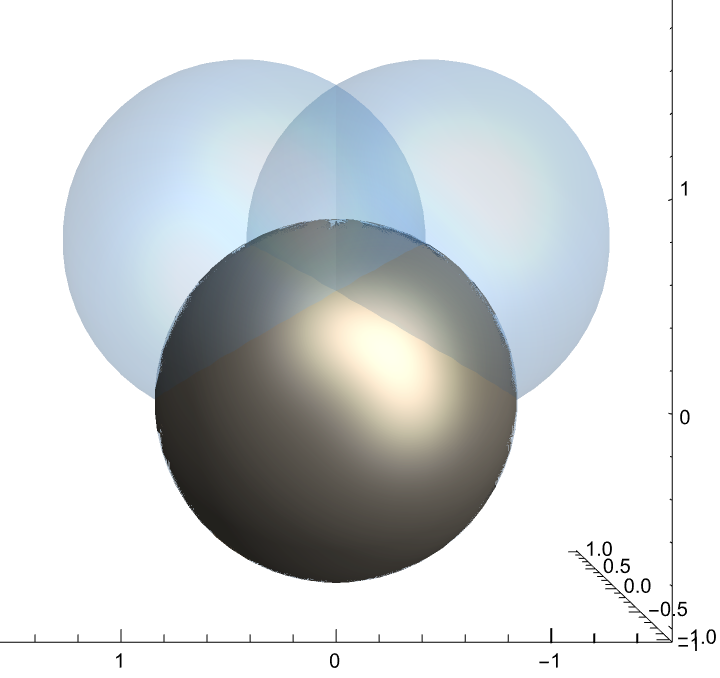} \qquad
			\includegraphics[scale=0.33]{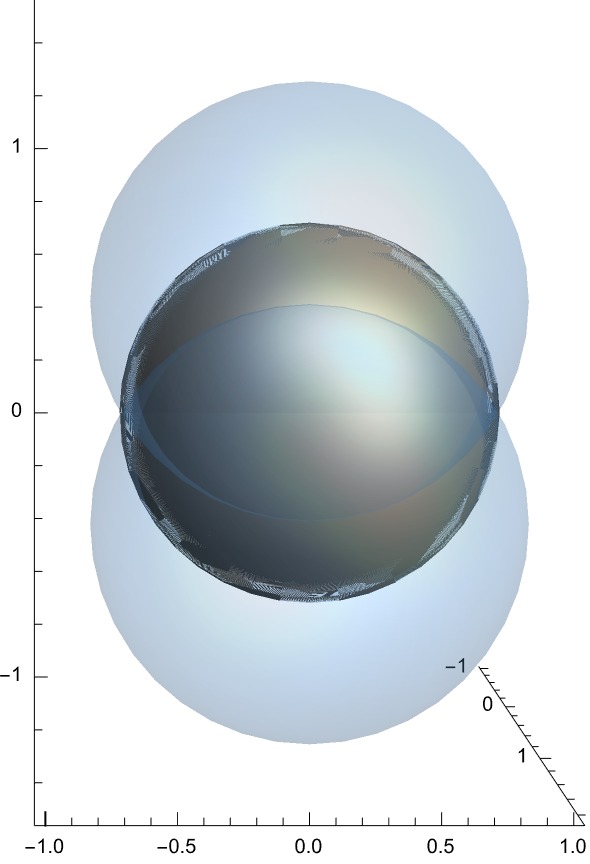} 
		\end{center}
		\begin{center}
			\includegraphics[scale=0.33]{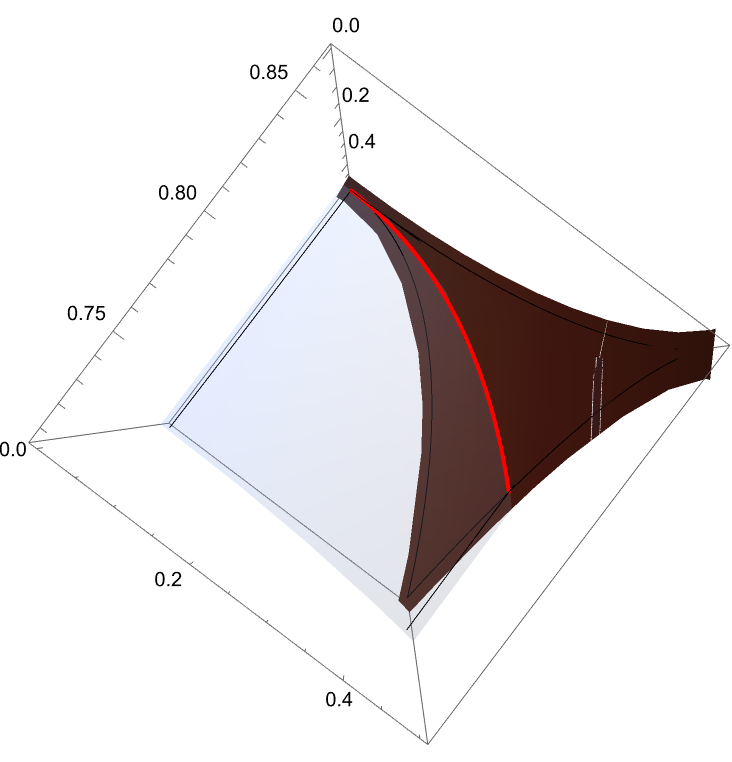} \qquad
			\includegraphics[scale=0.33]{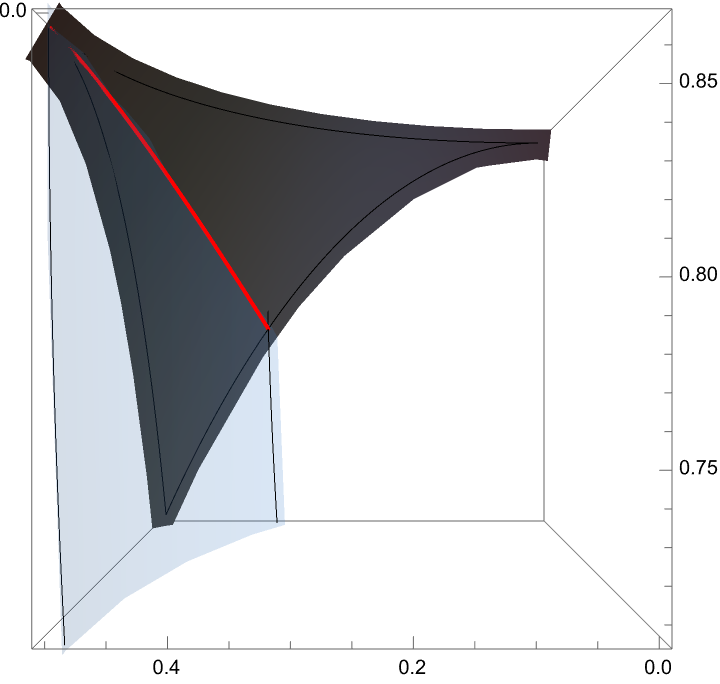} \qquad
			\includegraphics[scale=0.33]{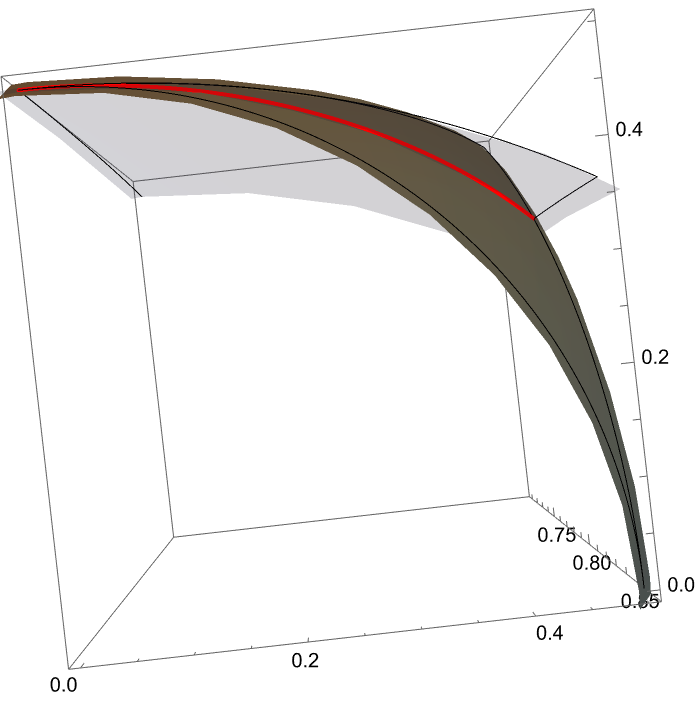} 
		\end{center}
		\begin{center}
			\includegraphics[scale=0.33]{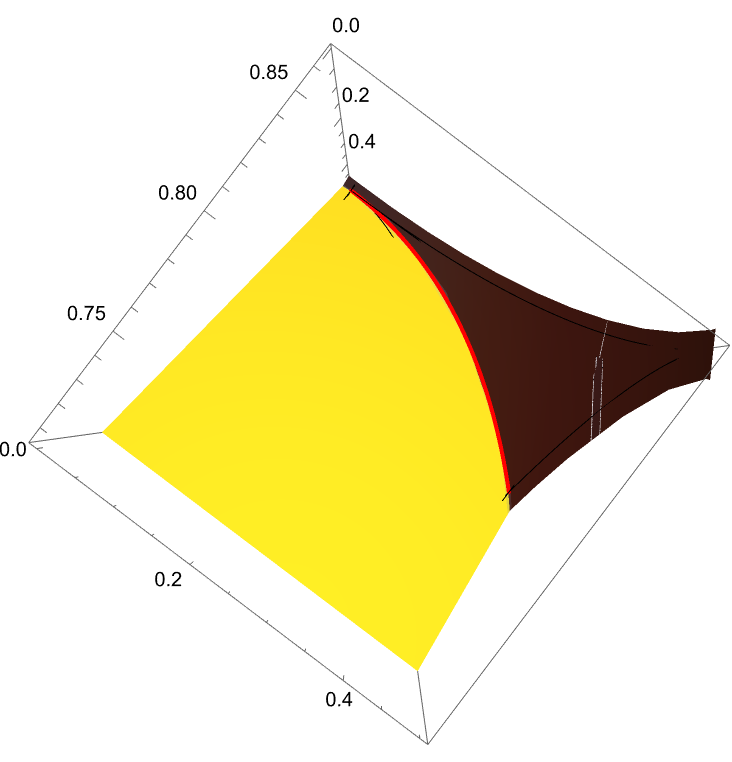} \qquad
			\includegraphics[scale=0.33]{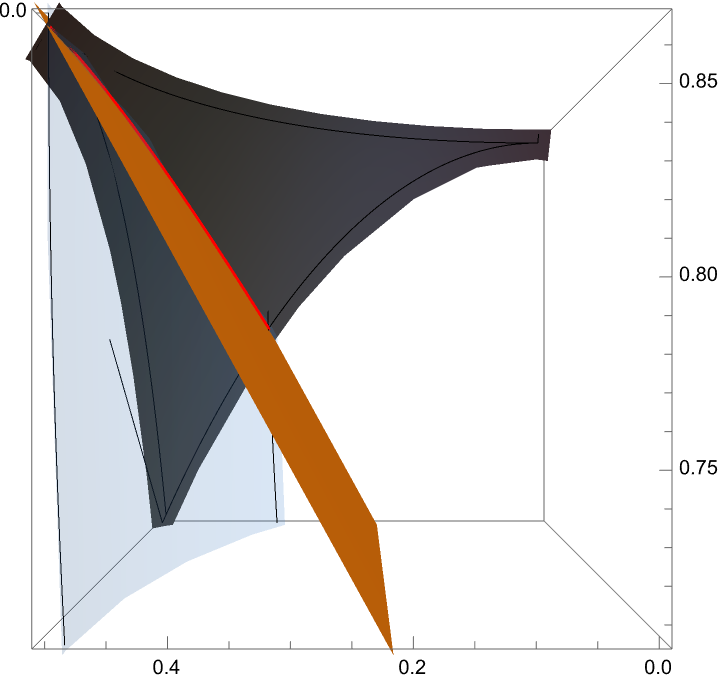} \qquad
			\includegraphics[scale=0.33]{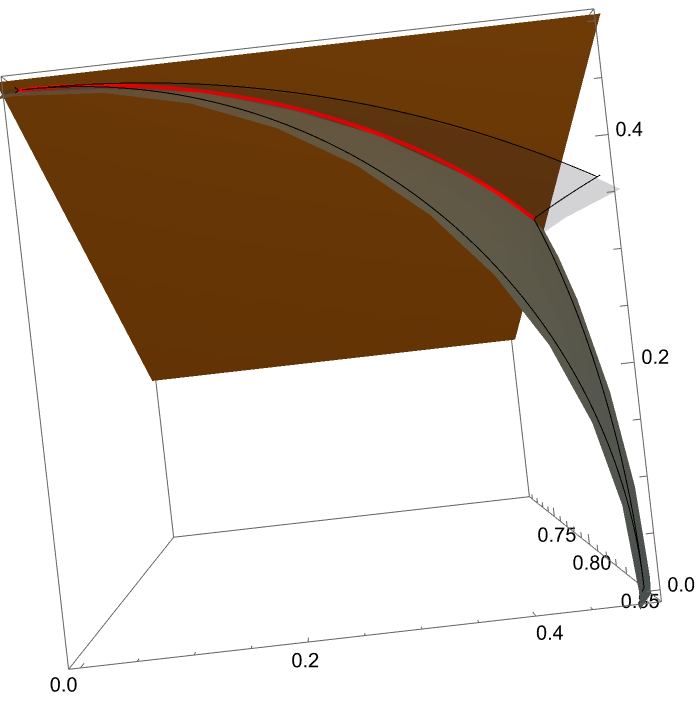} 
		\end{center}
		\caption{
			The first row shows the collection of spheres that are above to the region with $-1/2\leq x_j \leq 1/2$.
			The three potential spheres $B(0), B(\frac{i_1+\sqrt{3}i_2}{2})$, and $B(\frac{-i_1+\sqrt{3}i_2}{2})$ are to be used as bubbles for the maximal order of $\quat{-1,-3}{\QQ}$ where $\partial \Hcal^4 = V_3$. 
			These spheres are pictured from the front, left, and top, respectively, in each column.
			On the second and third rows we restrict to the region of $x = x_0 + i_1 x_1 +\sqrt{3}i_2 x_2$ with $ 0\leq x_j \leq 1/2$ (so the fundamental domain of $\Gamma_{\infty}$ is not included in this picture).
			The middle row shows the where the sphere $B(\frac{i_1+\sqrt{3}i_2}{2})$ intersects $B(1)$ in a red arc.
			The bottom row displays the plane containing the red line. 
			There is a region in $V_3$ in this plane containing the red line on the interior of both spheres where the height of $B(0)$ and $B(\frac{i_1+\sqrt{3}i_2}{2})$ in $\Hcal^4$ are equal.
			This is analogous to the chord determined by two circles---along this finite chord the two spheres sitting above the circles in the plane will have the same height. 
			Note the radial symmetry in the picture from the left (middle column) where the $x_1x_2$-axis is easier to see.
		}
		\label{fig:plot-1-3}
	\end{figure}

We will work with $\Ocal=\Ocal(-1,-3)_2.$
Everything is the same for $\Ocal(-1,-3)_1$: with $h=(1+i)/\sqrt{2}$, conjugation by $\diag(h,h^{-1})$, which acts on $\Hcal^4$ by the rotation $\pi_{1+i}$ ($1\mapsto i$, $i\mapsto -1$, $j\mapsto j$),
carries $\SL_2(\Ocal(-1,-3)_1)$ onto $\SL_2(\Ocal(-1,-3)_2)$.
This order contains the elements 
$$ \zeta = \frac{1+\sqrt{3}k}{2}, \quad \zeta^6=1, \qquad J = \frac{i+\sqrt{3}j}{2}, \quad J^2=-1.$$
The element $J$ may appear strange, since it acts like $\sqrt{-1}$ but looks like a primitive 6th root of unity in the $yz$-plane.
After running our algorithm for finding bubbles above the region with 
	$$ \Omega =\lbrace x+iy+jz \colon -1/2\leq x, y \leq 1/2, 0\leq z \leq \sqrt{3}/2 \rbrace$$
we found three potential spheres,
	$$B(0), \quad B(J),\quad B(J-i).$$
By symmetry we plotted these in Figure~\ref{fig:plot-1-3} (read the caption).
We will reduce the size of the region $\Omega$ and show that only $B(0)$ is actually needed.

\subsubsection{Unit Group}
The Clifford unit group of $\Ocal^{\times}$ is given by 
$$\Ocal^{\times} = \lbrace \pm 1 , \frac{\pm 1 \pm \sqrt{3}k}{2}, \frac{\pm i \pm \sqrt{3}j}{2}, \pm i \rbrace = \langle \zeta, J \rangle.$$ 
This group is generated by $\zeta$ and $J$, which have order $6$ and $4$, respectively. 
The elements $J$, $J\zeta$, and $\zeta J$ are all square roots of $-1$ and $J\zeta=i$. There are $1,1,2,6,2$ units of order $1,2,3,4,6$, respectively. 

The group modulo $\lbrace \pm 1 \rbrace$, which acts by conjugation, also has a nice description:
$$\Ocal^{\times}/\lbrace \pm 1\rbrace  = \langle \overline{\zeta}, \overline{J} \rangle \cong D_3.$$
Note that $\zeta J \zeta J = -1$ implies that $\Ocal^{\times}/\lbrace \pm 1 \rbrace \cong D_3$ with generators $\overline{\zeta}$ and $\overline{J}$.

\subsubsection{Fundamental Domain for $\Gamma_{\infty}$}
	For $a\in V_3$ the map $\pi_{\zeta}(a) = \zeta a \zeta^*$
	acts like a $2\pi/3$-rotation in the $yz$-plane and acts trivially in the $x$-plane.
	The map $\pi_{J}$ is given by 
	$$\pi_J(1) = -1, \quad \pi_J(i) = i/2-\sqrt{3}j/2,\quad \pi_J(j) = -\sqrt{3}i/2-j/2.$$ 
	This is the reflection in the $yz$-plane across the line
        $y + \sqrt{3} z = 0$ with a change in sign in the $x$-coordinate. 
	Note that this sign change in the $x$-axis converts the reflection
        into an orientation-preserving map, which can be thought of as a
        rotation (as usual, the compositum of two reflections is a rotation).

	One can see that the representation $\pi: D_3=\Ocal^{\times}/\lbrace \pm 1 \rbrace$ is reducible with $\RR \subset V_3$ being just the sign and $D_3$ acting on $\RR i + \RR j \subset V_3$ via the usual dihedral action.
	
	The fundamental domain in the $yi+zj$-plane is pictured in Figure~\ref{fig:stained-glass}. 
	This was found with a little trial and error and came from modifying a Dirichlet domain for $\Im(\Vec(\Ocal))$ to be invariant under the $D_3$ action.
	If we let $H_{12}$ be the inner honeycomb in Figure~\ref{fig:stained-glass}, which is not the fundamental domain, we can construct the ``Allen wrench'' 
	   $$A=[-1/2,1/2]\times H_{12}, $$
	 and observe that the Allen wrench is a fundamental domain for $\Vec(\Ocal)$. 
	 It is also invariant under the dihedral action with reflections turning around the Allen wrench. 
	 Then the wedge of the honeycomb (or ``orange slice'')  $W$, pictured with the fat lines in the stained glass picture (Figure~\ref{fig:stained-glass}) together with the extrusion $]-1/2,1/2[$, gives the fundamental domain for $\Gamma_{\infty}$,
         $$F =]-1/2,1/2[\times W.$$

\subsubsection{Fundamental Domain for $\Gamma$}
	The open fundamental domain for $\PSL_2(\Ocal(-1,-3)_2)$ acting on $\Hcal^4$ is
	$$D = \lbrace x+i_1y+i_2z+i_3w \colon  \vert x+i_1y+i_2z+i_3w\vert > 1 ,x+i_1y+i_2z \in F\rbrace.$$ 
	We are claiming that the only bubble that is actually used is $B(0)$.
	The bubbles $B(J)$ and $B(J-i)$ are not new: since $J, J-i\in\Vec(\Ocal)$, they are translates of $B(0)$ by elements of $\Gamma_{\infty}$.
	The region $F$ lies in the Allen wrench $A$, the Voronoi cell of $\Vec(\Ocal)$, on which $x^2+y^2+z^2\leq 1/4+1/3=7/12$. Hence above $F$ the points of $B(0)$ have height at least $\sqrt{5/12}$, and $B(0)$ is the highest of the unit hemispheres centered at points of $\Vec(\Ocal)$.
	A bubble $B(\mu^{-1}\lambda)$ has radius $1/\vert\mu\vert$, so only bubbles with $\vert\mu\vert^2\leq 2$ can reach height $\sqrt{5/12}$, and those with $\vert\mu\vert^2=1$ are translates of $B(0)$.
	
	For each of the finitely many bubbles with $\vert\mu\vert^2=2$ that meet the region above $A$, rising strictly above $B(0)$ at a point $v$ is a strict linear inequality in $v$; checking it at the vertices of $A$ shows that none of them rises above $B(0)$. The closest ones, for example those centered at $(\pm1\pm i)/2$, are tangent to $B(0)$ from below. 
	So we do not need any bubble other than $B(0)$, but only by the skin of our teeth.
\subsubsection{Generators}
	The generators are found using Theorem~\ref{thm:generators-ocal}. 
	The group $\SL_2(\Ocal)$ is generated by 
	$$ S=\begin{pmatrix}
	  0 & -1 \\
	  1 & 0
	\end{pmatrix}, \  
	\tau_J=\begin{pmatrix}
	1 & J \\
	0 & 1
	\end{pmatrix}, 
	\  
	\tau_i=\begin{pmatrix}
	1 & i \\
	0  & 1
	\end{pmatrix},\ 
	\tau_1=\begin{pmatrix}
	1 & 1\\
	0 & 1
	\end{pmatrix}, \ \pi_\zeta=\begin{pmatrix}
		\zeta & 0 \\
		0 & \zeta
	\end{pmatrix}, 
	\ \pi_J=\begin{pmatrix}
J & 0\\
	0 & -J
	\end{pmatrix}.
$$
The generator $S$ is for the bubble $B(0)$, the generators
$\tau_J,\tau_i,\tau_1$ are for the walls corresponding to
$\Vec(\Ocal)$, and $\pi_{\zeta}$ and $\pi_J$ are the extra symmetries
for the other walls corresponding to vertices of $F$.  Concretely, let
the fundamental domain depicted in Figure~\ref{fig:stained-glass} be $F_0$,
and for $i \in \ZZ/6\ZZ$ let $F_i$ be its counterclockwise rotation by
an angle of $2\pi i/6$.  Then $\pi_{\zeta}(F_i) = F_{i+2}$, while
$\pi_J(F_i) = F_{3-i}$.  In particular, the matrices $\pi_\zeta \pi_J$ and
$\pi_J \pi_\zeta$ give reflections across two of the thick black boundary
walls of $F_0$; for the third we may use the composition of $\pi_J$ with
translation by $J$.  One readily checks that
$\langle \pi_\zeta, \pi_J \rangle = \langle \pi_\zeta \pi_J, \pi_J \pi_\zeta \rangle$,
establishing the correctness of the given set of generators.  In fact,
we do not need $\tau_i$, since it does not correspond to any wall of the
fundamental domain.  A more basic set of generators would be
$\tau_1, \pi_J \pi_\zeta, \pi_\zeta \pi_J, \tau_J \pi_J, S$.  (The fundamental
domain has six walls, but $\tau_1^{-1}$, which crosses the wall at $x = -1/2$,
may be omitted, since $\tau_1$ is already present.  As expected, we can
express $\tau_i$ in terms of these generators: it is
$(\tau_J \pi_J) (\pi_\zeta \pi_J) (\tau_J \pi_J) (\pi_J \pi_\zeta)$.)

Relations arise from pairs of intersecting walls of the fundamental domain
as in Section~\ref{sec:finite-relations}.
Thus the product of $\tau_1$ with any of
$\pi_J \pi_\zeta, \pi_\zeta \pi_J, \tau_J \pi_J$, or of $\pi_J \pi_\zeta$ or
$\pi_\zeta \pi_J$ with $S$, has order $2$; the product of
any two of $\pi_J \pi_\zeta, \pi_\zeta \pi_J, \tau_J \pi_J$ has order $3$;
and likewise $(\tau_1 S)^3 = (\tau_J \pi_J S)^3 = 1$; together with the switch relations
$S^2 = (\pi_J \pi_\zeta)^2 = (\pi_\zeta \pi_J)^2 = (\tau_J \pi_J)^2 = 1$, these are all of
the relations.

Finally, we remark that for $\Ocal = \quat{-1,-3}{\ZZ}$, the analogue of the
Lipschitz order in this context, the group
$\SL_2(\Ocal)$ has index $15$ in $\SL_2(\Ocal(-1,-3)_2)$.  This can be proved by
Proposition~\ref{prop:coset-reps-from-orbits} or directly from the
presentation.

\section{The Case of  $\quat{-1,-1,-1}{\QQ}$} \label{sec:1-1-1}
By computation,
the integral Clifford algebra $\ZZ[i_1,i_2, i_3]$ is contained in a unique maximal order $\Ocal_4$.
The order $\Ocal_4$ is generated over $\ZZ[i_1,i_2,i_3]$ as an associative algebra by the elements $(1+i_{123})/2$ and $\zeta=(1+i_1+i_2+i_3)/2$.
Thus the code associated to this order is spanned by $1111$.
However, the order $\Ocal_4$ is not generated by any single element over
$\ZZ[i_1,i_2,i_3]$: to prove this, it suffices to consider one representative
for each coset of the additive group $\Ocal_4/\ZZ[i_1,i_2,i_3]$, whose order
is $64$. It is generated by $(1+i_{123})/2$ and $(1+i_1+i_2+i_{12})/2$.

The lattice $\Lambda = \Vec(\Ocal_4)$ is $\frac{1}{2}\Lambda_\code$,
which is the dual $D_4^{*}$ of the checkerboard lattice $D_4$; it is isometric to $\frac{1}{\sqrt{2}}D_4$.
This is the analog of the Hurwitz quaternions.

We remind the reader that the $D_4$ lattice in our presentation also has an $F_4$ root system; root lattices with short and long vectors can have two distinct root systems, and furthermore such lattices can be root lattices for both root systems. 
The ADE lattices are classified by their \emph{simply laced} Dynkin diagrams (no double arrows), which in this case is $D_4$.

\subsection{The case of the unique maximal order $\Ocal_4$ in $\quat{-1,-1,-1}{\QQ}$} \label{sec:maximal-c4}
This is an example where $\PSL_2(\Ocal_4)_{\infty}$ is very interesting but the bubbles are not so interesting. 
The order $\Ocal_4$ is our example of a Clifford-Euclidean order, hence cuspidally principal (Corollary~\ref{cor:Euclidean-implies-unimodular}), which is not Clifford-principal: not every right ideal generated by elements of $\Ocal_4^\mon$ is generated by a single element of $\Ocal_4^\mon$.

\begin{theorem}
	The order $\Ocal_4$ is Clifford-Euclidean but not 
	Clifford-principal.
\end{theorem} 
\begin{proof}
	It follows from Theorem~\ref{thm:clifford-euclidean-orders}
	that $\Ocal_4$ is Clifford-Euclidean.
	On the other hand, we claim that $\Ocal_4$
	is not Clifford-principal. Indeed, let
	$a = (-3 + i_1 + i_2 - i_3)(i_1), b = (-2-i_1-i_3)(1+i_1)$: both belong to
	$\Ocal_4^\mon$.
	The index in $\Ocal_4$
	of the right ideal $aR + bR$ is
	$9 \cdot 4^4$, which is not a fourth power. Therefore, this ideal is not
	generated by an element of $\Ocal_4^\mon$.
\end{proof}

Let $\Gamma = \PSL_2(\Ocal_4)$. 
The fundamental domain is clearly contained in the fundamental domain $F'$  for $\Gamma'_{\infty}$ where $\Gamma' = \PSL_2(\ZZ[i_1,i_2,i_3])$, which has three 180-degree rotations in the $x_0x_1$, $x_0x_2$ and $x_0x_3$ planes giving 
 $$ -\frac{1}{2}\leq x_0 \leq \frac{1}{2}, \quad 0 \leq x_1 \leq \frac{1}{2}, \quad 0 \leq x_2 \leq \frac{1}{2}, \quad 0 \leq x_3 \leq \frac{1}{2}.$$
We can figure out which spheres lie over $F'$.

Figure~\ref{fig:o4-spheres-slice} illustrates the lemma in the plane through $0$, $1$ and $\zeta$.
\begin{figure}[htbp!]
\centering
\begin{tikzpicture}[scale=1.6,lab]
 \begin{scope}\clip (-1.3,-0.6) rectangle (2.3,1.6);
 \foreach \a in {-2,...,3}{\foreach \b in {-1,0,1,2}{
   \draw[bub,fill=bubbleC!8,fill opacity=0.6] ({\a+0.5*\b},{0.866*\b}) circle (1);}}
 \foreach \a in {-2,...,3}{\foreach \b in {-1,0,1,2}{\fill ({\a+0.5*\b},{0.866*\b}) circle (1.5pt);}}
 \end{scope}
 \draw[holeC,line width=1pt,dashed] (0.5,0.2887) circle (0.7071);
 \fill[holeC] (0.5,0.2887) circle (1.5pt);
 \node[below left] at (0,0){$0$}; \node[below right] at (1,0){$1$}; \node[above] at (0.5,0.95){$\zeta$};
 \draw[gray] (0,0)--(1.5,0.866);
 \node at (0.5,-0.85){(a) top view of the slice};
\end{tikzpicture}\qquad
\begin{tikzpicture}[xscale=2.6,yscale=2.6,lab]
 \draw[bub,domain=0:1,samples=60] plot (\x,{sqrt(1-\x*\x)});
 \draw[bub,domain=0:1.05,samples=80] plot (\x,{sqrt(max(0,1-(\x*\x-1.732*\x+1)))});
 \draw[holeC,dashed,line width=1pt,domain=-0.13:1.28,samples=80] plot (\x,{sqrt(max(0,0.5-(\x-0.5774)^2))});
 \draw[bnd] (-0.1,0)--(1.1,0);
 \fill (0,0) circle (1.2pt) node[below]{$0$}; \fill[holeC] (0.5774,0) circle (1.2pt) node[below]{deep point};
 \node[below] at (0.866,-0.12){\scriptsize edge midpoint};
 \node at (0.55,-0.42){(b) profile along the grey line};
\end{tikzpicture}
\caption{Lemma~\ref{lem:O4-spheres} in the plane through $0$, $1$, $\zeta=(1+i_1+i_2+i_3)/2$, where $\Vec(\Ocal_4)$ is a triangular lattice. Blue: unit hemispheres $B(\nu)$, $\nu\in\Vec(\Ocal_4)$. Purple: a hemisphere of radius $1/\sqrt2$, the largest radius allowed by the lemma, centred at the point of the slice farthest from the lattice. (b) Heights along the grey line from $0$: the purple curve stays below the blue. Two-dimensional slice of $V_4$; the purple hemisphere is a test case, not a specific bubble.}\label{fig:o4-spheres-slice}
\end{figure}
\begin{lemma}\label{lem:O4-spheres}
	Suppose that $a = \mu^{-1}\lambda\in \Vec(K)$ with $(\mu,\lambda)$ unimodular and $a\notin\Vec(\Ocal_4)$.
	Then $B(a)$ is contained in the union $\bigcup_{\nu\in\Vec(\Ocal_4)}B(\nu)$ of the unit hemispheres centered at the points of $\Vec(\Ocal_4)$, which are the $\Gamma_\infty$-translates of $B(0)$.
\end{lemma}
\begin{proof}
	The radius of $B(a)$ is $r=1/\vert\mu\vert$, and $\vert\mu\vert^2=\mu\bar\mu$ is a positive integer.
	If $\vert\mu\vert^2=1$, then $\mu$ is a Clifford unit of $\Ocal_4$ and $a\in\Ocal_4\cap\Vec(K)=\Vec(\Ocal_4)$, which we have excluded; so $r^2\le 1/2$.
	The squared covering radius of $\Vec(\Ocal_4)=\ZZ^4+\ZZ\zeta$ is $1/2$: the Voronoi cell of $0$ is the 24-cell $\{x:\vert x_j\vert\le 1/2,\ \sum_j\pm x_j\le 1\}$, whose $24$ vertices, the permutations of $(\pm\frac12,\pm\frac12,0,0)$, all have squared norm $1/2$.
	Let $(x,t)$ be a point under $B(a)$, so that $\vert x-a\vert^2+t^2<r^2$, and let $\nu\in\Vec(\Ocal_4)$ be a point nearest to $x$.
	Then $\vert x-\nu\vert^2+t^2<\frac12+r^2\le 1$, so $(x,t)$ lies under $B(\nu)$.
\end{proof}

\begin{remark}
In the case of $\PSL_2(\ZZ[i_1,i_2,i_3])$ the element $\zeta$ will now be a deep hole for the lattice and is a \emph{singular cusp}. We cannot place a sphere there. 
We deal with this in \S\ref{sec:1-1-1-order}.
\end{remark}

\subsubsection{Unit Group}

The group $\Ocal_4^\times$ of Clifford units has order $576$ and is given by 
$$\Ocal_4^{\times} = \langle i_1,i_2,i_3,\zeta, \alpha \rangle, \qquad \zeta = \frac{1+i_1+i_2+i_3}{2}, \quad \alpha=\frac{1 + i_{12}-i_{13}+i_{23}}{2}.$$
The subgroup $\langle i_1,i_2,i_3\rangle=\langle i_1,i_2\rangle\times\langle i_{123}\rangle\cong Q_8\times C_2$ has order 16, and the maps $\pi_{i_a}$ act by $\pi_{i_a}(1)=-1$, $\pi_{i_a}(i_a)=-i_a$, $\pi_{i_a}(i_b)=i_b$, and should be viewed as a rotation of angle $\pi$ in the $x_0x_a$-plane.
The subgroup $\Vec(\Ocal_4)^{\times}$, defined to be the group generated by the Clifford units of $\Ocal_4$ that lie in $\Vec(\Ocal_4)$, is $\langle i_1, i_2,i_3, \zeta\rangle$, which has order $192$. 
The group $\lbrace \pi_v \colon v\in \Vec(\Ocal_4)^{\times}\rbrace \cong \Vec(\Ocal_4)^{\times}/\lbrace \pm 1 \rbrace$ that describes the action of the integral
vectors that are Clifford units on the space of Clifford vectors has order $96$ and is a subgroup of the positive orientation Weyl group $\Weyl(D_4)^+$, which has index two in $\Weyl(D_4)$ (these are the transformations that have determinant~$1$).
The Weyl group of $D_4$ is isomorphic to $C_2^3\rtimes S_4$ (in general $\Weyl(D_n)$ is $C_2^{n-1}\rtimes S_n$ where $S_n$ permutes the factors and $C_2^{n-1}$ is viewed as an $n-1$-dimensional $\FF_2$-vector subspace of $\FF_2^n$ ) and $\Weyl(D_4)^+ \cong C_2^3\rtimes A_4$. 
Since $96 = 2^3\cdot (4!/2)$, it follows that
\begin{equation}\label{eqn:weyl-d4plus}
  \langle \pi_{i_1}, \pi_{i_2}, \pi_{i_3}, \pi_{\zeta} \rangle = \Weyl(D_4)^+,
\end{equation}
with the elements $\pi_{i_1}, \pi_{i_2}, \pi_{i_3}$ being the three generators of $C_2^3$ and $\pi_{\zeta}$ acting as a rotation by an angle of $2\pi/3$. 

The quotient group $\Ocal_4^{\times}/\lbrace \pm 1 \rbrace$ has order
$576/2 = 288 = 3\cdot 96$. 
The element $\pi_{\alpha}$ acts as a cyclic permutation of $i_1,i_2,i_3$, and this is actually an automorphism of the Dynkin diagram (so an automorphism of the root system which isn't in the group generated by reflections in the roots).
One has $\Aut(\Dynkin(D_4))=S_3$ (see Figure~\ref{fig:d4}); this has to do
with the famous {\em triality automorphism} of $\SO_8$.
\begin{figure}[htbp!]
	\centering
	\begin{tikzpicture}[lab]
 \node[dnode,label={[dlab]right:$-\tfrac12(1+i_1+i_2+i_3)$}] (c) at (0,0) {};
 \node[dnode,label={[dlab]below:$i_1$}] (a) at ($(c)+(180:\dstep)$) {};
 \node[dnode,label={[dlab]right:$i_2$}] (b) at ($(c)+(60:\dstep)$) {};
 \node[dnode,label={[dlab]right:$i_3$}] (d) at ($(c)+(-60:\dstep)$) {};
 \draw[dedge] (c)--(a) (c)--(b) (c)--(d);
 \node at (0,-1.75){(a) $D_4$};
\end{tikzpicture}\qquad
\begin{tikzpicture}[lab]
 \node[dnode,label={[dlab]below:$1-i_1$}] (n0) at (0,0) {};
 \node[dnode,label={[dlab]below:$i_1-i_2$}] (n1) at (\dstep,0) {};
 \node[dnode,label={[dlab]below:$i_2$}] (n2) at (2*\dstep,0) {};
 \node[dnode,label={[dlab]right:$-\tfrac12(1+i_1+i_2+i_3)$}] (n3) at (3*\dstep,0) {};
 \draw[dedge] (n0)--(n1) (n2)--(n3); \dyndouble{n1}{n2}
 \node at (2*\dstep,-1.0){(b) $F_4$};
\end{tikzpicture}
	\caption{The lattice $\Vec(\Ocal_4)$ is a root lattice for both the $D_4$ and $F_4$ root systems.
		(a) The Dynkin diagram of $D_4$ with simple roots in $\Vec(\Ocal_4)$; all of them are short.
		(b) The Dynkin diagram of $F_4$; the double bond points from the long roots to the short roots.
	}\label{fig:d4}
\end{figure}
This element of order $3$ describes $\Ocal_4^{\times}/\lbrace \pm 1 \rbrace$ as an extension of $\Weyl(D_4)^+$ by $C_3$, so rather than an action on just $D_4$, it is perhaps best to think about this as a representation to $\Aut(\Phi)$ where $\Phi$ is the root system.
The group $\Weyl(F_4)$ has order $1152$, which is $4 \cdot 288$.  We now
explain how $288$ arises as the order of $\Ocal^{\times}_4/\lbrace \pm 1\rbrace$.

The $D_4$-lattice $\Vec(\Ocal_4)$ also has an $F_4$-root system. 
Now, when we speak of $D_4$ we will speak of the root system $\Phi_0 \subset \Vec(\Ocal_4)$ of this lattice and not the lattice itself. 
The $D_4$-root system is given by $\pm$ coordinate vectors and all vectors with entries $\pm 1/2$,
so $24$ in all.
The $F_4$-root system is the $D_4$-root system together with all vectors with two
entries $0$ and two $\pm 1$.  

We have $\Aut(D_4) = \Aut(F_4) = \Weyl(F_4)$, which we will just call $G$. 
This is a Schl\"afli  group and hence is the symmetry group of a polytope in $C_{24}$ (see \S\ref{sec:polyhedra}). 
We have $\Weyl(D_4) \subset \Aut(D_4)$ being a normal subgroup generated by reflections,
and $\Aut(\Phi)/\Weyl(\Phi)$ is the set of
graph automorphisms of the Dynkin diagram for the $D_4$-root system (Figure~\ref{fig:d4}).
This is the group $S_3$.
This means the triality map $\pi_{\alpha}$ is actually in $\Weyl(F_4)^+$, so that the image of $\Ocal_4^\times$ has index $2$ in $\Weyl(F_4)^+$ and index $4$ in $\Weyl(F_4)$. Since $\Weyl(F_4)$ is the symmetry group of the 24-cell $C_{24}$, which is the Voronoi cell of $0$ in $\Vec(\Ocal_4)$, a fundamental domain for $\Gamma_{\infty}$ in $V_4$ is made of four of the $1152$ congruent simplices into which the mirrors of $\Weyl(F_4)$ cut $C_{24}$ (see the theorem below).

\subsection{Coxeter Groups and Classification of Polyhedra}\label{sec:polyhedra}
In this section, we give a brief introduction to Coxeter groups and their
relation to the classification of polyhedra.
The reader unfamiliar with Coxeter groups may wish to consult a standard
reference such as \cite{BB2005}. Here we only state a few basic definitions
and results.

A \emph{Coxeter group} \cite[1.1]{BB2005}
is a pair consisting of a group $G$ and a set of
involutions $S = \{s_1, \dots, s_n\}$ that generate $G$, such that the
relations are all of the form $(s_is_j)^{m_{ij}} = 1$. Given a Coxeter group,
the associated \emph{Coxeter diagram} or \emph{Coxeter graph}
is the graph whose vertex set is in bijection with $S$ and that has
an edge labeled $m_{ij}$ between $v_i$ and $v_j$ if $m_{ij} > 2$.

In particular, two vertices are adjacent if and only if the corresponding
generators do not commute.
Every group of isometries of $\RR^n$ or $\Hcal^n$ generated by
reflections is a Coxeter group,
so the distinguished generators are often called {\em reflections}.

A \emph{Schl\"afli symbol} is a sequence $\lbrace m_{0,1}, m_{1,2},
m_{2,3},\ldots, m_{n-1,n}\rbrace$ of integers $\ge 3$ that encodes a
Coxeter 
group with generators $s_0,\ldots,s_n$
and relations $(s_is_{i+1})^{m_{i,i+1}}=1$ for $i=0,1,\ldots,n-1$ and
$(s_is_j)^2 = 1$ for $|i-j| > 1$. For example, the symbol $\lbrace
3,4,3\rbrace$ is the group generated by $s_0,s_1,s_2,s_3$ where
$$(s_0s_1)^{3}=(s_1s_2)^4=(s_2s_3)^3=(s_0s_2)^2 = (s_0s_3)^2 = (s_1s_3)^2 = 1.$$
Schl\"{a}fli symbols are in bijection with Coxeter diagrams whose underlying graph is a path.
See Figure~\ref{fig:coxeter-diagrams}.
\begin{figure}[htbp!]
\centering
\begin{tikzpicture}[lab]
 \node[dnode,label={[dlab]below:$s_0$}] (a) at (0,0) {};
 \node[dnode,label={[dlab]below:$s_1$}] (b) at (\dstep,0) {};
 \node[dnode,label={[dlab]below:$s_2$}] (c) at (2*\dstep,0) {};
 \node (d) at (3*\dstep,0) {$\cdots$};
 \node[dnode,label={[dlab]below:$s_n$}] (e) at (4.6*\dstep,0) {};
 \draw[dedge] (a)--node[dm]{$m_{0,1}$}(b)--node[dm]{$m_{1,2}$}(c)--(d)--node[dm]{$m_{n-1,n}$}(e);
 \node at (2*\dstep,-0.9){(a) $\{m_{0,1},m_{1,2},\dots,m_{n-1,n}\}$};
\end{tikzpicture}\qquad
\begin{tikzpicture}[lab]
 \foreach \k in {0,...,3}{\node[dnode,label={[dlab]below:$s_\k$}] (n\k) at (\k*\dstep,0) {};}
 \draw[dedge] (n0)--(n1) (n2)--(n3); \draw[dedge] (n1)--node[dm]{$4$}(n2);
 \node at (1.5*\dstep,-0.9){(b) $\{3,4,3\}$, $\Weyl(F_4)$};
\end{tikzpicture}
\caption{(a) A Coxeter diagram whose graph is a path, labelled by the Schl\"afli symbol $\{m_{0,1},\dots,m_{n-1,n}\}$: generators $s_i$, with $(s_is_{i+1})^{m_{i,i+1}}=1$; unjoined vertices commute. Edges with $m=3$ are left unlabelled. (b) The symbol $\{3,4,3\}$ of the 24-cell, whose isometry group is $\Weyl(F_4)$.}\label{fig:coxeter-diagrams}
\end{figure}
	
\begin{theorem}

  The isometry group of a regular polyhedron $P$ is a Coxeter group whose
  underlying graph is a path (and can hence be described by a Schl\"afli
  symbol) \cite{Coxeter1973}. Regular polyhedra up to similarity and duality
  are classified by their associated isometry group. 
\end{theorem}

\subsubsection{Fundamental Domain for $\Gamma_{\infty}$}
In what follows we make use of the bijection between polyhedra and their Coxeter group. 
See \S\ref{sec:polyhedra}.
See Figure~\ref{fig:24-cell}(a).
\begin{figure}[htbp!]
\centering
\begin{tikzpicture}[scale=1.0,lab]
 \draw[cellE] (1.954,1.954)--(2.669,0.715);
\draw[cellE] (1.954,1.954)--(2.669,-0.715);
\draw[cellE] (1.954,1.954)--(1.239,0.715);
\draw[cellE] (1.954,1.954)--(1.239,-0.715);
\draw[cellE] (1.954,1.954)--(0.715,2.669);
\draw[cellE] (1.954,1.954)--(0.715,1.239);
\draw[cellE] (1.954,1.954)--(-0.715,2.669);
\draw[cellE] (1.954,1.954)--(-0.715,1.239);
\draw[cellE] (2.669,0.715)--(2.669,-0.715);
\draw[cellE] (2.669,0.715)--(1.239,0.715);
\draw[cellE] (2.669,0.715)--(1.954,-1.954);
\draw[cellE] (2.669,0.715)--(0.715,2.669);
\draw[cellE] (2.669,0.715)--(1.430,-0.000);
\draw[cellE] (2.669,0.715)--(0.000,1.430);
\draw[cellE] (2.669,0.715)--(0.715,-1.239);
\draw[cellE] (2.669,-0.715)--(1.239,-0.715);
\draw[cellE] (2.669,-0.715)--(1.954,-1.954);
\draw[cellE] (2.669,-0.715)--(0.715,1.239);
\draw[cellE] (2.669,-0.715)--(1.430,-0.000);
\draw[cellE] (2.669,-0.715)--(-0.000,-1.430);
\draw[cellE] (2.669,-0.715)--(0.715,-2.669);
\draw[cellE] (1.239,0.715)--(1.239,-0.715);
\draw[cellE] (1.239,0.715)--(1.954,-1.954);
\draw[cellE] (1.239,0.715)--(-0.715,2.669);
\draw[cellE] (1.239,0.715)--(0.000,1.430);
\draw[cellE] (1.239,0.715)--(-1.430,0.000);
\draw[cellE] (1.239,0.715)--(-0.715,-1.239);
\draw[cellE] (1.239,-0.715)--(1.954,-1.954);
\draw[cellE] (1.239,-0.715)--(-0.715,1.239);
\draw[cellE] (1.239,-0.715)--(-0.000,-1.430);
\draw[cellE] (1.239,-0.715)--(-1.430,0.000);
\draw[cellE] (1.239,-0.715)--(-0.715,-2.669);
\draw[cellE] (1.954,-1.954)--(0.715,-1.239);
\draw[cellE] (1.954,-1.954)--(0.715,-2.669);
\draw[cellE] (1.954,-1.954)--(-0.715,-1.239);
\draw[cellE] (1.954,-1.954)--(-0.715,-2.669);
\draw[cellE] (0.715,2.669)--(0.715,1.239);
\draw[cellE] (0.715,2.669)--(-0.715,2.669);
\draw[cellE] (0.715,2.669)--(1.430,-0.000);
\draw[cellE] (0.715,2.669)--(0.000,1.430);
\draw[cellE] (0.715,2.669)--(-1.954,1.954);
\draw[cellE] (0.715,2.669)--(-1.239,0.715);
\draw[cellE] (0.715,1.239)--(-0.715,1.239);
\draw[cellE] (0.715,1.239)--(1.430,-0.000);
\draw[cellE] (0.715,1.239)--(-0.000,-1.430);
\draw[cellE] (0.715,1.239)--(-1.954,1.954);
\draw[cellE] (0.715,1.239)--(-1.239,-0.715);
\draw[cellE] (-0.715,2.669)--(-0.715,1.239);
\draw[cellE] (-0.715,2.669)--(0.000,1.430);
\draw[cellE] (-0.715,2.669)--(-1.430,0.000);
\draw[cellE] (-0.715,2.669)--(-1.954,1.954);
\draw[cellE] (-0.715,2.669)--(-2.669,0.715);
\draw[cellE] (-0.715,1.239)--(-0.000,-1.430);
\draw[cellE] (-0.715,1.239)--(-1.430,0.000);
\draw[cellE] (-0.715,1.239)--(-1.954,1.954);
\draw[cellE] (-0.715,1.239)--(-2.669,-0.715);
\draw[cellE] (1.430,-0.000)--(0.715,-1.239);
\draw[cellE] (1.430,-0.000)--(0.715,-2.669);
\draw[cellE] (1.430,-0.000)--(-1.239,0.715);
\draw[cellE] (1.430,-0.000)--(-1.239,-0.715);
\draw[cellE] (0.000,1.430)--(0.715,-1.239);
\draw[cellE] (0.000,1.430)--(-0.715,-1.239);
\draw[cellE] (0.000,1.430)--(-1.239,0.715);
\draw[cellE] (0.000,1.430)--(-2.669,0.715);
\draw[cellE] (-0.000,-1.430)--(0.715,-2.669);
\draw[cellE] (-0.000,-1.430)--(-0.715,-2.669);
\draw[cellE] (-0.000,-1.430)--(-1.239,-0.715);
\draw[cellE] (-0.000,-1.430)--(-2.669,-0.715);
\draw[cellE] (-1.430,0.000)--(-0.715,-1.239);
\draw[cellE] (-1.430,0.000)--(-0.715,-2.669);
\draw[cellE] (-1.430,0.000)--(-2.669,0.715);
\draw[cellE] (-1.430,0.000)--(-2.669,-0.715);
\draw[cellE] (0.715,-1.239)--(0.715,-2.669);
\draw[cellE] (0.715,-1.239)--(-0.715,-1.239);
\draw[cellE] (0.715,-1.239)--(-1.239,0.715);
\draw[cellE] (0.715,-1.239)--(-1.954,-1.954);
\draw[cellE] (0.715,-2.669)--(-0.715,-2.669);
\draw[cellE] (0.715,-2.669)--(-1.239,-0.715);
\draw[cellE] (0.715,-2.669)--(-1.954,-1.954);
\draw[cellE] (-0.715,-1.239)--(-0.715,-2.669);
\draw[cellE] (-0.715,-1.239)--(-2.669,0.715);
\draw[cellE] (-0.715,-1.239)--(-1.954,-1.954);
\draw[cellE] (-0.715,-2.669)--(-2.669,-0.715);
\draw[cellE] (-0.715,-2.669)--(-1.954,-1.954);
\draw[cellE] (-1.954,1.954)--(-1.239,0.715);
\draw[cellE] (-1.954,1.954)--(-1.239,-0.715);
\draw[cellE] (-1.954,1.954)--(-2.669,0.715);
\draw[cellE] (-1.954,1.954)--(-2.669,-0.715);
\draw[cellE] (-1.239,0.715)--(-1.239,-0.715);
\draw[cellE] (-1.239,0.715)--(-2.669,0.715);
\draw[cellE] (-1.239,0.715)--(-1.954,-1.954);
\draw[cellE] (-1.239,-0.715)--(-2.669,-0.715);
\draw[cellE] (-1.239,-0.715)--(-1.954,-1.954);
\draw[cellE] (-2.669,0.715)--(-2.669,-0.715);
\draw[cellE] (-2.669,0.715)--(-1.954,-1.954);
\draw[cellE] (-2.669,-0.715)--(-1.954,-1.954);
\fill (1.954,1.954) circle (1.1pt);
\fill (2.669,0.715) circle (1.1pt);
\fill (2.669,-0.715) circle (1.1pt);
\fill (1.239,0.715) circle (1.1pt);
\fill (1.239,-0.715) circle (1.1pt);
\fill (1.954,-1.954) circle (1.1pt);
\fill (0.715,2.669) circle (1.1pt);
\fill (0.715,1.239) circle (1.1pt);
\fill (-0.715,2.669) circle (1.1pt);
\fill (-0.715,1.239) circle (1.1pt);
\fill (1.430,-0.000) circle (1.1pt);
\fill (0.000,1.430) circle (1.1pt);
\fill (-0.000,-1.430) circle (1.1pt);
\fill (-1.430,0.000) circle (1.1pt);
\fill (0.715,-1.239) circle (1.1pt);
\fill (0.715,-2.669) circle (1.1pt);
\fill (-0.715,-1.239) circle (1.1pt);
\fill (-0.715,-2.669) circle (1.1pt);
\fill (-1.954,1.954) circle (1.1pt);
\fill (-1.239,0.715) circle (1.1pt);
\fill (-1.239,-0.715) circle (1.1pt);
\fill (-2.669,0.715) circle (1.1pt);
\fill (-2.669,-0.715) circle (1.1pt);
\fill (-1.954,-1.954) circle (1.1pt);
 \node at (0,-3.0){(a) 24-cell, $F_4$ Coxeter-plane projection};
\end{tikzpicture}\qquad
\begin{tikzpicture}[lab]
 \node[dnode] (1) at (0,0){}; \node[dnode] (2) at (\dstep,0){}; \node[dnode] (3) at (2*\dstep,0){};
 \node[dnode] (4) at ($(3)+(60:\dstep)$){}; \node[dnode] (5) at ($(3)+(-60:\dstep)$){};
 \draw[dedge] (1)--(2)--(3)--(4) (3)--(5);
 \node at (1.2*\dstep,-1.9){(b) $D_5$};
\end{tikzpicture}
\caption{(a) The 24-cell $C_{24}=\{x\in V_4: |x_j|\le\tfrac12,\ \sum_j\pm x_j\le1\}$, the Voronoi cell of $0$ in $\Vec(\Ocal_4)$ (\S\ref{sec:maximal-c4}), projected onto the Coxeter plane of $F_4$ (drawn from the $24$ points $\pm e_a\pm e_b$, a rotated and rescaled copy of $C_{24}$; edges join vertices at distance $\sqrt2$). Since $\Vec(\Ocal_{5,0})=\ZZ\oplus\Lambda'$ with $\Lambda'$ a copy of $\Vec(\Ocal_4)$, the Voronoi cell for $\Ocal_{5,0}$ is $[-\tfrac12,\tfrac12]\times C_{24}$. (b) The Dynkin diagram of $D_5$; the Clifford units of $\Ocal_{5,!}$ act on $V_5$ through $\Weyl(D_5)^+$. Projection computed numerically.}\label{fig:24-cell}
\end{figure}
The \emph{24-cell} is the unique polyhedron with Schl\"afli symbol $\lbrace 3,4,3\rbrace.$
Its isometry group is $\Weyl(F_4)$.
The 24-cell also goes by the names  \emph{icositetrachoron}, \emph{octaplex}, \emph{icosatetrahedroid}, \emph{octacube}, \emph{hyperdiamond} or \emph{polyoctahedron}.

\begin{theorem}
	 The fundamental domain for $\PSL_2(\Ocal_4)_{\infty}$ acting on $V_4$ is a union of four chambers of the affine reflection group $\Vec(\Ocal_4)\rtimes\Weyl(F_4)$, that is, of four of the $1152$ congruent simplices into which the mirrors of $\Weyl(F_4)$ cut the 24-cell $C_{24}$. Its volume is $1/576$, which is $1/288$ of the volume of $C_{24}$.
\end{theorem}
\begin{proof}
The 24-cell $C_{24}=\{x\in V_4:\vert x_j\vert\le 1/2,\ \sum_j\pm x_j\le 1\}$ is the Voronoi cell of $0$ in $\Vec(\Ocal_4)$, so it is a fundamental domain for the translations, and its symmetry group $\Weyl(F_4)$ preserves $\Vec(\Ocal_4)$.
Let $U$ be the image of $\pi:\Ocal_4^{\times}/\lbrace \pm 1\rbrace \to \GL(V_4)$ given by $\pi_u(x) = uxu^*$, so that $\PSL_2(\Ocal_4)_\infty$ acts on $V_4$ as $\Vec(\Ocal_4)\rtimes U$.
We showed via a \texttt{magma} computation that $U$ has index $4$ in $\Weyl(F_4)$.

The group $U$ intersects the group of signed permutation matrices in a subgroup of order $96$ (we can apply an even permutation to the coordinates and independently change the sign of an even number of them).
\begin{samepage}
For example, the unit
\[
  \zeta=\frac{1+i_1+i_2+i_3}{2}
\]
gives a matrix that acts in a more complicated way.
\end{samepage}
One can check that this group of order $288$ is a subgroup of the Coxeter group $\Weyl(F_4)$.
We can extend this group by the diagonal matrix $\diag(1,1,1,-1)$ and the permutation matrix $(3,4)$.
These together with $U$ generate a group isomorphic to $\Weyl(F_4)$.
So if $L$ is a fundamental domain for $\Vec(\Ocal_4)\rtimes\Weyl(F_4)$ acting on $V_4$, then the union of the images of $L$ under representatives $g$ of the cosets $Ug$ of $U$ in $\Weyl(F_4)$ is a fundamental domain for $\Vec(\Ocal_4)\rtimes U$.
In particular we may take the identity, the two elements $\diag(1,1,1,-1)$ and
$(3,4)$ mentioned above, and their product, and the domain $L$ of the lemma below; this gives Corollary~\ref{cor:O4-F}. Since $L$ has volume $1/2304$, the union has volume $4/2304=1/576$, while $C_{24}$ has volume $\operatorname{covol}\Vec(\Ocal_4)=1/2$.
\end{proof}

A fundamental domain for the group of signed permutations of order 384,
which is the Weyl group $\Weyl(B_4)$ acting on the unit cube centered at
$0$, is the subset
$$\lbrace (a,b,c,d) \colon 0 \leq d \leq c\leq b \leq a \leq 1/2
\rbrace.$$
Our group is an extension of this, and all the rows of
elements of our group either have one $\pm 1$ and three $0$ entries or
all four entries $\pm 1/2$. Knowing the inequalities just given, we
can determine all inequalities formed by dot products with vectors of
$\pm 1$ except for two. For example, we know that $a + c \geq b + d$;
however, the inequalities $a+d \geq b+c$ and  $a \geq b+c+d$ are
undetermined.
	\begin{lemma}
		The locus $L$ defined by
		$$L=\lbrace (a,b,c,d) \colon 0 \leq d \leq c \leq b \leq  a \leq 1/2, a  \geq  b + c  + d \rbrace $$
		is a fundamental domain for the affine group $\Vec(\Ocal_4)\rtimes\Weyl(F_4)$ acting on $V_4$; equivalently, it is a fundamental domain for $\Weyl(F_4)$ acting on the 24-cell $C_{24}$. It is the $4$-simplex with vertices $0$, $(\frac12,0,0,0)$, $(\frac12,\frac12,0,0)$, $(\frac12,\frac14,\frac14,0)$, $(\frac12,\frac16,\frac16,\frac16)$, and its volume is $1/2304$.
	\end{lemma}
	\begin{proof}
	The vectors $e_2-e_3$, $e_3-e_4$, $e_4$, $\frac12(e_1-e_2-e_3-e_4)$ form a base of the $F_4$ root system described above \cite[\S12.1]{Humphreys1978}, so the cone $\{0\le d\le c\le b,\ b+c+d\le a\}$ is the closure of a Weyl chamber, which is a fundamental domain for $\Weyl(F_4)$ acting on $V_4$ \cite[\S10.3]{Humphreys1978}.
	On this cone $\vert x_j\vert\le a$ for all $j$ and $\sum_j\vert x_j\vert\le 2a$, so its intersection with $C_{24}$ is its intersection with $\{a\le 1/2\}$, which is $L$.
	Since $\Weyl(F_4)$ preserves $C_{24}$, and $C_{24}$ is a fundamental domain for the translations by $\Vec(\Ocal_4)$, the first two claims follow. The vertices and the volume are a direct computation.
	\end{proof}

\begin{corollary}\label{cor:O4-F}
	Let $\sigma$ be the permutation matrix given by $\sigma(a,b,c,d) = (a,b,d,c)$. 
	Let $\tau$ be the transformation $\tau(a,b,c,d) = (a,b,c,-d)$.
	The fundamental domain for $\PSL_2(\Ocal_4)_{\infty}$ is then 
	$$ F = L \cup \sigma(L) \cup \tau(L) \cup \sigma\tau(L).$$
\end{corollary}

\begin{remark} We can describe $F$ by the inequalities
  $1/2 \ge a \ge b \ge \pm c; b \ge \pm d; a \ge b + c \pm d; a \ge b - c + d; c + d \ge 0$.
  In fact $a \ge b$ is unnecessary, since it is implied by
  $a \ge b+c+d$ and $c+d \ge 0$.  Nor do we need $b \ge -c$ or $b \ge -d$,
  which follow from $b \ge d$ (respectively $b \ge c$) and $c+d \ge 0$.
  So in fact the only inequalities needed are
  $$1/2 \ge a, \quad b \ge c, \quad b \ge d, \quad c+d \ge 0, \quad a \ge b+c\pm d, \quad a \ge b - c + d.$$
\end{remark}

\subsubsection{Fundamental Domain for $\Gamma$}
By Lemma~\ref{lem:O4-spheres}, every hemisphere lies in the union of the unit hemispheres $B(\nu)$, $\nu\in\Vec(\Ocal_4)$. Since $F$ lies in the Voronoi cell $C_{24}$ of $0$, a point above $F$ that lies outside $B(0)$ lies outside every $B(\nu)$. So the fundamental domain is the region above $F$ and outside $B(0)$; it is bounded by the walls above the walls of $F$ together with $B(0)$.

\subsubsection{Generators for $\Gamma$} \label{sec:gen-gamma-o4}
The group $\Gamma= \SL_2(\Ocal_4)$ is generated by the generators of the unit group, the inversion in the single sphere, and the generators of the lattice $\Vec(\Ocal_4)$.
This implies that $\Gamma$ is generated by
\begin{equation}\label{eqn:o4-generators}
 S, \quad \pi_{i_1}, \pi_{i_2}, \pi_{i_3}, \pi_{\zeta},\pi_{\alpha} \quad \tau_1, \tau_{i_1}, \tau_{i_2}, \tau_{\zeta}.
\end{equation}

Alternatively, we may describe a presentation of the group based on the
walls of the fundamental domain of $\Gamma_\infty$ as above.  With one
exception, if any single
inequality from those defining $F$ is violated, we still remain inside
the fundamental domain of $\Vec \Ocal_4$, so the element of $\PSL_2$
that maps the fundamental domain to its translate is always conjugation by
a unit.  The exception is $1/2 \ge a$, for which the appropriate transformation
first negates $a$ and $d$ (if $c>d$) or $a$ and $c$ (if $c<d$) and then translates by $1$.  We thus give a table of
the affine transformations and corresponding elements of $\PSL_2$.  These
elements together with $\begin{psmallmatrix}0&1\\-1&0\end{psmallmatrix}$
generate the group.  The ``Map'' column indicates the 
affine transformation that maps the adjacent translate of $F$ back to $F$, evaluated at a general point
$(a,b,c,d)$.  If $u$ is a Clifford unit we write $\pi_u$ for the matrix
$\begin{psmallmatrix}u&0\\0&(u^*)^{-1}\end{psmallmatrix}$.  The group is then
generated by $S$ and the entries in the ``Element'' column.

\begin{table}[h]
  \begin{center}
    \small
    \renewcommand{\arraystretch}{1.2}\begin{tabular}{|c|c|c|c|}\hline
      Inequality & Map & Element & $z_k$ \\ \hline \hline
      $1/2 \ge a$, $c > d$ & $(1-a,b,c,-d)$ & $T_0 = \tau_1 \pi_{z_0}$ & $z_0 = i_3$ \\ \hline
      $1/2 \ge a$, $c < d$ & $(1-a,b,-c,d)$ & $T_1 = \tau_1 \pi_{z_1}$ & $z_1 = i_2$ \\ \hline
$b \ge c$ & $(a,c,d,b)$ & $\pi_{z_2}$ & $z_2 = (1-i_{12}+i_{13}-i_{23})/2$ \\ \hline
$b \ge d$ & $(a,d,b,c)$ & $\pi_{z_3}$ & $z_3 = (1+i_{12}-i_{13}+i_{23})/2$ \\ \hline
$c + d \ge 0$ & $(a,b,-c,-d)$ & $\pi_{z_4}$ & $z_4 = i_{23}$ \\ \hline
      $a \ge b+c+d$, $c>d$ & \begin{tabular}{@{}l@{}}$\tfrac12(a+b+c+d,\,a+b-c-d,$\\ $\phantom{\tfrac12(}a-b+c-d,\,-a+b+c-d)$\end{tabular} & $\pi_{z_5}$ & $z_5 = (1-i_3+i_{13}+i_{23})/2$ \\ \hline
      $a \ge b+c+d$, $c<d$ & \begin{tabular}{@{}l@{}}$\tfrac12(a+b+c+d,\,a+b-c-d,$\\ $\phantom{\tfrac12(}{-a}+b-c+d,\,a-b-c+d)$\end{tabular} & $\pi_{z_6}$ & $z_6 = (1-i_2+i_{12}-i_{23})/2$ \\ \hline
      $a \ge b+c-d$ & \begin{tabular}{@{}l@{}}$\tfrac12(a+b+c-d,\,a+b-c+d,$\\ $\phantom{\tfrac12(}a-b+c+d,\,a-b-c-d)$\end{tabular} & $\pi_{z_7}$ & $z_7 = (1+i_3-i_{13}-i_{23})/2$ \\ \hline
      $a \ge b-c+d$ & \begin{tabular}{@{}l@{}}$\tfrac12(a+b-c+d,\,a+b+c-d,$\\ $\phantom{\tfrac12(}a-b-c-d,\,a-b+c+d)$\end{tabular} & $\pi_{z_8}$ & $z_8 = (1+i_2-i_{12}+i_{23})/2$ \\ \hline
    \end{tabular}
  \end{center}
  \caption{Facets of the fundamental domain $F$ of $\Gamma_\infty$ (Corollary~\ref{cor:O4-F}), where
    $\Gamma = \PSL_2(\Ocal_4)$, and the elements of $\Gamma$ that map the
    adjacent images of $F$ to $F$.  Note that there are two distinct
    translates of the fundamental domain adjacent across the hyperplanes
  $a = b+c+d$ and $a = 1/2$.}
  \label{tab:gens-psl2-o4}
\end{table}

\subsubsection{Relations for $\Gamma$}\label{sec:rel-gamma-o4}
Using the generators given above, we may give a complete presentation for
$\PSL_2(\Ocal_4)$.  To our knowledge, this is the first time that such a
thing has been done for an order in a Clifford algebra with more than $2$
imaginary units.  The group generated by $\pi_{z_i}$ for $2 \le i \le 8$
is finite (being the image of the group $\Ocal_4^\times$ of Clifford units) and so its
relations may easily be described.  For the remaining relations, we note
that the unit hemisphere, like the hyperplane $a=b+c+d$, separates the
fundamental domain from two of its translates, which are again separated by
the hyperplane $c = d$.  The corresponding generators are given by
$S_0=S\pi_{z_0}$ and $S_1=S\pi_{z_1}$, where $S=\begin{psmallmatrix}0&1\\-1&0\end{psmallmatrix}$.  For purposes of generating the group this is
of no importance, since the $\pi_{z_i}$ from Table~\ref{tab:gens-psl2-o4}
already generate the image of $\Ocal_4^\times$, but these
additional generators are needed to determine the relations in a sensible way,
as described in Section~\ref{sec:gens-and-rels}.
There are $20$ relations that involve one or more of the
$S_i, T_i$, where $T_0,T_1$ are as in Table~\ref{tab:gens-psl2-o4} and $\pi_k$ denotes $\pi_{z_k}$.  These are as follows:

\begin{center}
  \begin{tabular}{lllll}
    $S_1 \pi_4 S_0$ & $(T_0 S_0)^3$ & $ (\pi_2 S_1 \pi_3 S_0)^2$ & $\pi_4 S_1 S_0$ & $(\pi_5 S_0)^2$ \\
    $ (\pi_7 S_0)^2$ & $(T_1 S_1)^3$ & $(\pi_3 S_0 \pi_2 S_1)^2$ & $\pi_4 S_0 S_1$ & $(\pi_6 S_1)^2$ \\
    $ (\pi_8 S_1)^2$ & $ T_1 \pi_4 T_0$ & $(\pi_2 T_1 \pi_3 T_0)^2$ & $ \pi_4 T_1 T_0$ & $ (\pi_5 T_0)^3$ \\
    $(\pi_7 T_0)^3$ & $(\pi_3 T_0 \pi_2 T_1)^2$ & $\pi_4 T_0 T_1$ & $(\pi_6 T_1)^3$ & $(\pi_8 T_1)^3$ \\
  \end{tabular}
\end{center}

Magma rapidly reduces the presentation to one on the four generators
$S_0, T_0, \pi_2, \pi_5$ of order $2, 2, 3, 3$ respectively and satisfying
the additional relations 

\begin{equation*}
  \begin{aligned}
    &    (S_0\pi_5^{-1})^2,\quad (T_0\pi_5^{-1})^3,\quad (T_0S_0)^3,\quad (\pi_2\pi_5^{-1}\pi_2^{-1}\pi_5)^2,\quad S_0 \pi_2 \pi_5 \pi_2^{-1} \pi_5S_0 \pi_2\pi_5^{-1} \pi_2^{-1}\\
    & \quad (T_0 \pi_2 \pi_5 \pi_2^{-1}\pi_5^{-1})^2,\quad (\pi_2 \pi_5 \pi_2^{-1} \pi_5)^3,\quad \pi_2\pi_5^{-1}S_0\pi_2\pi_5^{-1}\pi_2S_0 \pi_2^{-1} \pi_5 \pi_2^{-1}S_0 \pi_5 \pi_2^{-1}S_0 \\
    & \qquad \pi_2 T_0 \pi_5\pi_2\pi_5^{-1} \pi_2 T_0 \pi_2 T_0 \pi_2 \pi_5 \pi_2^{-1}\pi_5^{-1} \pi_2^{-1} T_0. \\
  \end{aligned}
\end{equation*}

\subsection{The Case of $\quat{-1,-1,-1}{\ZZ}$} \label{sec:1-1-1-order}
The example of $\PSL_2(\ZZ[i_1,i_2,i_3])$ acting on $\Hcal^5$ is very different from the behavior of $\PSL_2(\ZZ[i_1,i_2])$ acting on $\Hcal^4$ (which was dealt with in \cite{Maclachlan1989}) and $\PSL_2(\ZZ[i_1])$ acting on $\Hcal^3$ (which is classical). 
What is interesting about $\ZZ[i_1,i_2,i_3] \subset \Ocal_4$ is that passing to the subgroup $\PSL_2( \ZZ[i_1,i_2,i_3]) \subset \PSL_2(\Ocal_4)$ causes
the cusp at $\zeta = (1+i_1+i_2+i_3)/2$ to become inequivalent to $\infty$
since this element is no longer in our order.
This issue is dealt with abstractly in \S\ref{sec:finite-index-subgroups} where we dealt with finite index subgroups. 

For the rest of this section we use the notation
 $$ \Ocal = \ZZ[i_1,i_2,i_3].$$

The following was performed in \texttt{magma}. 
Following Algorithm~\ref{alg:sl2-index}, to find generators for $\SL_2(\ZZ[i_1,i_2,i_3])$ we took words of length up to $6$ in the $10$ generators of $\SL_2(\Ocal_4)$ (equation \eqref{eqn:o4-generators}) with a sampling of other words. 
Call this group $\Gamma$. 

In order to give a more enlightening proof, we use the method of
Proposition~\ref{prop:coset-reps-from-orbits}.
We retain our notation and let $\Gamma_x=\Stab_{\Gamma}(x)$ for a group $\Gamma$ acting on a space containing an element $x$. 
Figure~\ref{fig:z-i1i2i3-cusps} shows the two cusp classes along the line $\RR\zeta$.
\begin{figure}[htbp!]
\centering
\begin{tikzpicture}[scale=1.5,lab]
 \foreach \k in {-2,...,3}{\draw[bub] (\k+1,0) arc[start angle=0,end angle=180,radius=1];}
 \draw[bnd] (-2.4,0)--(3.4,0) node[right]{$\RR\zeta$};
 \foreach \k in {-2,0,2}{\fill (\k,0) circle (2.2pt);}
 \foreach \k in {-1,1,3}{\fill[cuspC] (\k,0) circle (2.4pt);}
 \node[below] at (0,0){$0$}; \node[below,cuspC] at (1,0){$\zeta$}; \node[below] at (2,0){$2\zeta$};
 \node[below] at (-1,0){$-\zeta$};\node[below] at (3,0){$3\zeta$};
 \draw[->] (0.3,1.5)--(0.3,1.9) node[above]{$\infty$};
 \node[right] at (4.1,0.9){black: $H\infty$};
 \node[right,cuspC] at (4.1,0.6){red: $H\zeta$};
\end{tikzpicture}
\caption{Cusps on the line $\RR\zeta\subset V_4$, $\zeta=(1+i_1+i_2+i_3)/2$, with the unit bubbles centred there (blue). For $\PSL_2(\Ocal_4)$ every point $k\zeta$ is a cusp equivalent to $\infty$. For the subgroup $H=\PSL_2(\ZZ[i_1,i_2,i_3])$, the points $2k\zeta\in\ZZ^4$ remain in the orbit of $\infty$ (black) and the points $(2k+1)\zeta\in\zeta+\ZZ^4$ lie in the orbit of $\zeta$ (red); see Proposition~\ref{prop:index-120}. One-dimensional slice.}\label{fig:z-i1i2i3-cusps}
\end{figure}
\begin{proposition}\label{prop:index-120}
\begin{enumerate}
	\item $\SL_2(\ZZ[i_1,i_2,i_3])_0\backslash\SL_2(\Ocal_4)_0$ has representatives $$\begin{pmatrix}v&0\\sv&{v^{-1}}^*\end{pmatrix}$$ 
	where $v, s$ run over coset representatives for $\Ocal_4^\times/\Ocal^\times$ (note that $\Ocal^{\times} = \langle i_1,i_2,i_3\rangle$ and $\Ocal_4^{\times} = \langle i_1,i_2,i_3,\zeta,\alpha\rangle$)
	and $\Vec(\Ocal_4)/\Vec(\Ocal)=\lbrace [0],[\zeta]\rbrace$, respectively.
	\item $\SL_2(\ZZ[i_1,i_2,i_3])_{\zeta}\backslash\SL_2(\Ocal_4)_{\zeta}$ has representatives which are the $\tau_{\zeta}$ conjugates of $$\begin{pmatrix}w&0\\tw&{w^{-1}}^*\end{pmatrix}$$ 

	Here $w$ and $t$ run over coset representatives for, respectively,
	\begin{align*}
	  \Ocal_4^\times/\pi^{-1}(\Weyl(D_4)^+)&=\langle[\alpha]\rangle,\\
	  \Vec(\Ocal_4)/2\Vec(\Ocal_4)&\cong\FF_2^4.
	\end{align*}
	\item $\SL_2(\ZZ[i_1,i_2,i_3])$ has index 120 in $\SL_2(\Ocal_4)$, and we can give explicit coset representatives.

\end{enumerate}
\end{proposition}
\begin{proof}
  We retain our notation $\Ocal=\ZZ[i_1,i_2,i_3]$, and write $G=\SL_2(\Ocal_4)$ and $H=\SL_2(\Ocal)$.
  We first show that $H$ has exactly two orbits on $\Vec(K)\cup\{\infty\}$, represented by $0$ and $\zeta$.
  Since $\Ocal_4$ is Clifford-Euclidean, $G$ acts transitively on $\Vec(K)\cup\{\infty\}$ (Corollary~\ref{cor:Euclidean-implies-unimodular}).
  For a cusp $v\neq\infty$ choose $g=\begin{psmallmatrix}a&b\\c&d\end{psmallmatrix}\in G$ with $g\infty=v=ac^{-1}$, and put $N(v)=c\bar c$.
  This is a positive integer that does not depend on $g$, because any other choice has first column $(au,cu)$ with $u\in\Ocal_4^\times$.
  If $q\in\Vec(\Ocal)=\ZZ^4$, then $S\tau_{-q}g$ has first column $(c,-(v-q)c)$, so $N(S\tau_{-q}v)=\vert v-q\vert^2N(v)$.
  The covering radius of $\ZZ^4$ is $1$, and a nearest $q$ has $\vert v-q\vert=1$ only if $v-q\in\zeta+\ZZ^4\subset\Vec(\Ocal_4)$.
  So if $v\notin\Vec(\Ocal_4)$ we can decrease $N(v)$ by an element of $H$; hence every cusp is $H$-equivalent to $\infty=S0$ or to a cusp with $N(v)=1$, that is, to a point of $\Vec(\Ocal_4)=\ZZ^4\cup(\zeta+\ZZ^4)$, and so to $0$ or to $\zeta$.
  Conversely, $\zeta$ is not $H$-equivalent to $\infty$: the elements of $G$ taking $\infty$ to $\zeta$ are those with first column $(\zeta u,u)$, $u\in\Ocal_4^\times$, and if $u$ and $\zeta u$ were both in $\Ocal$ then so would be $\zeta=(\zeta u)\bar u$.
  So the cusps $0$ and $\zeta$ represent the two $H$-orbits.
\begin{enumerate}
\item 
The stabilizer of $0$ in $\SL_2(\Ocal_4)$ is given by the matrices of the
form $\begin{psmallmatrix}u&0\\ru&(u^{-1})^*\end{psmallmatrix}$, where
$u$ is a Clifford unit and $r$ is a Clifford vector in $\Ocal_4$.
To calculate the index of the subgroup of $\SL_2(\Ocal_4)_0$ consisting
of matrices with entries in $\SL_2(\Ocal)$, let
$\Gamma' \subset \SL_2(\Ocal_4)_0$ be the subgroup of matrices with
top left entry in $\Ocal$. By the same stabilizer calculation used in Lemma~\ref{lem:third}, we have
$[\SL_2(\Ocal_4)_0:\Gamma'] = [\Ocal_4^\times:\Ocal^\times] = 576/16 = 36$
and $[\Gamma':\SL_2(\Ocal)_0] = [\Vec(\Ocal_4):\Vec(\Ocal)] = 2$.
Thus $[\SL_2(\Ocal_4)_0:\SL_2(\Ocal)_0] = 72$.  In fact this argument
shows that the set of matrices of the form
$$\begin{pmatrix}v&0\\sv&{v^{-1}}^*\end{pmatrix}$$ 
  is a set of coset representatives for $\SL_2(\Ocal)_0$ in $\SL_2(\Ocal_4)_0$,
where $v, s$ run over coset representatives for $\Ocal_4^\times/\Ocal^\times$
and $\Vec(\Ocal_4)/\Vec(\Ocal)$ respectively.
\item 
  Similarly, the stabilizer of $\zeta$ in $\Ocal_4$ is the conjugate of
the stabilizer subgroup $\SL_2(\Ocal_4)_0$ by a matrix, such as
$\begin{psmallmatrix}1&-\zeta\\0&1\end{psmallmatrix}$, that takes
$\zeta$ to $0$. This conjugate is
\begin{equation}\label{eqn:conjugate}
\SL_2(\Ocal_4)_{\zeta} = \left\{\begin{pmatrix}u+\zeta ru& -u\zeta - \zeta ru \zeta + \zeta {u^{-1}}^* \\ ru & -ru\zeta + {u^{-1}}^*\end{pmatrix}: u \in \Ocal_4^\times, r \in \Vec(\Ocal_4)\right\},
\end{equation}
and we need to determine the index of the subgroup consisting of elements all
of whose entries are in $\Ocal$. Let $M_{u,r}$ be the matrix above. We
first observe that, if $r'-r$ is such that $\zeta^i (r'-r)u \zeta^j$ has
integral coefficients
for all
$u \in \Ocal_4^\times$, $i,j \in \{0,1\}$,
then $M_{u,r} \in \SL_2(\Ocal)_{\zeta}$ if and only if
$M_{u,r'} \in \SL_2(\Ocal)_{\zeta}$. Indeed, all entries of $M_{u,r'}-M_{u,r}$
are of the form $\zeta^i (r'-r)u \zeta^j$ as indicated. One calculates that
this holds if and only if $r \in 2 \Vec(\Ocal_4)$.

Let $\Gamma_{1,\zeta}$ be the subgroup of $\SL_2(\Ocal)_{\zeta}$ defined by 
$$\Gamma_{1,\zeta}=\lbrace \begin{pmatrix}1+\zeta r& -\zeta - \zeta r \zeta + \zeta  \\ r & -r\zeta + 1\end{pmatrix}  \colon r \in 2 \Vec(\Ocal_4) \rbrace.$$
This is the subgroup where we have set $u=1$ in equation \eqref{eqn:conjugate}.
Since $\Gamma_{1,\zeta}\subseteq\SL_2(\Ocal)_{\zeta}\subseteq\SL_2(\Ocal_4)_{\zeta}$, multiplicativity of the index gives
\[
  [\SL_2(\Ocal_4)_{\zeta}:\SL_2(\Ocal)_{\zeta}]
  =\frac{[\SL_2(\Ocal_4)_{\zeta}:\Gamma_{1,\zeta}]}
         {[\SL_2(\Ocal)_{\zeta}:\Gamma_{1,\zeta}]}.
\]
We will compute the two indices on the right.

A set of coset representatives for
$\Gamma_{1,\zeta}$ in $\SL_2(\Ocal_4)_{\zeta}$ is obtained by letting $u$ range
over $\Ocal_4^\times$ and $r$ over $\Vec(\Ocal_4)/2\Vec(\Ocal_4)$ and hence
has order $576 \cdot 16$.

  This is not so many that one could not list them all on a computer
and determine directly
which ones lie in $\SL_2(\Ocal)_{\zeta}$, but it is easier to argue as follows.
Fix $u$ and consider the different $[r] \in \Vec(\Ocal_4)/2\Vec(\Ocal_4)$.  
If $r, r'$ both give matrices in $\SL_2(\Ocal)_\zeta$, then
$(r'-r)u, (r'-r)u\zeta$ both belong to $\Ocal$; 
it follows that
$(r'-r)u \in 2 \Ocal_4$, so $r, r'$ represent the same coset. 
In other words,
for each $u$ there is at most one possible choice of
$[r] \in \Vec(\Ocal_4)/2 \Vec(\Ocal_4)$.

Take a set of Clifford units generating the subgroup $\pi^{-1}(\Weyl(D_4)^+) = \langle i_1,i_2,i_3,\zeta\rangle$ of index $3$ in $\Ocal_4^\times$ from display \eqref{eqn:weyl-d4plus}, where $\pi:\Ocal_4^{\times} \to \Aut(V_4)$ is the usual representation.
We find that for every element of this subgroup there is in fact one
possible choice of $[r]  \in \Vec(\Ocal_4)/2 \Vec(\Ocal_4)$ giving an
element of $\SL_2(\Ocal)$.
Since there are Clifford units
for which there is no such choice (any unit not in $\frac{1}{2} \Ocal$, for
example), the index is exactly $3$.
Thus $192$ of the $576 \cdot 16$
choices are valid, giving an index of $48$.  
As before, we may take
the conjugates by $\tau_{\zeta}=\begin{pmatrix}1&\zeta\\0&1\end{pmatrix}$ of the
$$\begin{pmatrix}w&0\\tw&{w^{-1}}^*\end{pmatrix}$$ 
as coset
representatives, where $w, t$ run over coset representatives for the
index-$3$ subgroup of $\Ocal_4^\times$ and $\Vec(\Ocal_4)/2\Vec(\Ocal_4)$,
respectively.

\item 
We conclude by applying Proposition~\ref{prop:coset-reps-from-orbits},
finding that $[\SL_2(\Ocal_4):\SL_2(\Ocal)] = 72 + 48 = 120.$
Proposition~\ref{prop:coset-reps-from-orbits} also gives a recipe for finding the coset representatives from the cosets of the stabilizer subgroups.
\end{enumerate}
\end{proof}

\begin{remark} The index $[\SL_2(\Ocal_4):\SL_2(\Ocal)]$ can in principle
  be determined from the presentation given in Section~\ref{sec:rel-gamma-o4},
  but this is a somewhat painful computation.  The group of integral matrices
  was not generated by the words of length at most $12$ in the $4$ generators that we examined.
\end{remark}

\begin{remark}
	The lattice $\Lambda_{\zeta} = 2\Vec(\Ocal_4)$ defines the torus $V_4/\Lambda_{\zeta}$, and $f^{-1}([\zeta])$ is the quotient of this torus by a finite group, where 
	$$f: \PSL_2(\ZZ[i_1,i_2,i_3])\backslash \Hcal^{5,\BS} \to \PSL_2(\ZZ[i_1,i_2,i_3]) \backslash \Hcal^{5,\Sat}$$ 
	is the map from the Borel--Serre compactification to the Satake compactification.
\end{remark}

\subsubsection{Relation to Known Reflection Groups}\label{sec:relation-to-known}
Abstractly, we can make some contact with the non-algebraic arithmetic group $\Gamma_5$ defined in \cite[\S 7, pg 298]{Ratcliffe2019}; more generally, this discussion holds for certain $\Gamma_{n+1}$ and $\PSL_2(\ZZ[i_1,\ldots,i_{n-1}])$ for all $n\geq 2$.
Recall that the nonalgebraic group $\PO(1,n+1) \cong \Isom(\Hcal^{n+1})$ preserves the $y>0$ sheet of the real quadric $-y^2+x_0^2+\cdots +x_n^2 =-1$. It contains the transformation $(1,-1_{n+1})$, of determinant $(-1)^{n+1}$, given by $(y,x)\mapsto (y,-x)$.
The groups $\Gamma_{n+1}$ are defined by
$$ \Gamma_{n+1} = \PO(1,n+1) \cap \O_{1,n+1}(\ZZ). $$
The group $\Gamma_{5} \subset \Isom(\Hcal^5)$ is a non-compact 5-simplex reflection group and hence has 6 generators.
A similar result holds for $\Gamma_i$ with $2 \le i \le 9$ but is false for
larger $i$.
Generalized polytopes are covered in \cite[p. 266]{Ratcliffe2019}.

\begin{question}\label{prop:known-groups}
  Let $G_1, G_2$ be the images of $\PSL_2(\ZZ[i_1,\ldots,i_{n-1}])$
  and $\Gamma_{n+1}$ in $\Isom(\Hcal^{n+1})$. For which $n\geq 2$ is $G_1$ conjugate to a subgroup of index $2$ in $G_2$?
\end{question}

\begin{remark}
  Both groups are arithmetic, and $G_1$ is commensurable with a conjugate of $G_2$: by the proof of Theorem~\ref{thm:arithmeticity}, the image of $G_1$ is commensurable with $\SO_Q(\ZZ)\cap\O_Q(\RR)^{\circ}$ for $Q=x^2-yz+x_1^2+\cdots+x_{n-1}^2$, and the substitution $y=a+b$, $z=a-b$ identifies the index-$2$ sublattice $\{y\equiv z \bmod 2\}$ of $(\ZZ^{n+2},Q)$ with the standard Lorentzian lattice defining $\Gamma_{n+1}$.
  The spin sequence of Theorem~\ref{thm:short-exact-sequence} shows that the image of $G_1$ has index $2$ in the image of $\PSpin_Q(\ZZ)$, but it does not compute the index in $G_2$: for $n\geq 2$ the map $\PSpin_Q(\ZZ)\to\SO_Q(\ZZ)$ is not surjective (Remark~\ref{rem:spin-sequence-at-2}), and $\Gamma_{n+1}$ is defined by a different lattice.
  For $n=1$ the index is not $2$: $\Gamma_2$ contains the reflection group of the triangle with angles $\pi/2,\pi/4,0$, of covolume $\pi/4$, while $\PSL_2(\ZZ)$ has covolume $\pi/3$.
\end{remark}

In order to understand the relation between $\SO$ and $\PO$, we recall some
facts about the outer automorphisms of $\O_{1,n+1}(\RR)$. 
Much of what is written here is known to experts but is not written down formally anywhere according to
\cite[\href{https://mathoverflow.net/questions/235758/automorphism-group-of-real-orthogonal-lie-groups}{MO235758}]{mathoverflow-autgp}, which we follow.

Let $\O(p,q)$ denote the real Lie group 
$$\O(p,q) = \O_{p,q}(\RR)$$ 
where $\O_{p,q}$ is the group scheme over $\ZZ$ associated to
the standard indefinite quadratic form of signature $(p,q)$.
We first define the characters $\det_p, \det_q: \O(p,q) \to \{\pm 1\}$.
The subgroup $\O(p) \times O(q)$ is maximal compact in $\O(p,q)$. Define
$\det_p, \det_q$ to be the maps taking $(M_p,M_q) \in \O(p) \times O(q)$ to
$\det M_p, \det M_q$, respectively. Then $(\det_p,\det_q)$ is a surjective
homomorphism with connected kernel to a discrete group, so it is the
map to the component group.
It is known \cite[Chapter XV, Theorem 3.1]{Hochschild1965}
that a real Lie group with finite component group is smoothly homeomorphic
to the product of a maximal compact subgroup by a Euclidean space.
Hence the component groups are equal and the
homomorphism $(\det_p,\det_q): \O(p) \times O(q) \to \{\pm 1\}^2$ extends to
$\O(p,q)$. We again refer to the components as $\det_p, \det_q$; the other
two elements of the group are reasonably called $1, \det$.

\begin{proposition}\label{prop:easy-auts} Let $\chi: \O(p,q) \to \{\pm 1\}$ be a character and define
  $\mu_{\chi}: \O(p,q) \to \O(p,q)$ by $\mu_\chi(g) = \chi(g) g$. If
  $\chi(-I_{p+q}) = 1$, then $\mu_{\chi}$ is an automorphism of $\O(p,q)$.
\end{proposition}

\begin{proof} Since $\chi(g)\in\{\pm1\}$ is a central scalar, $\mu_\chi$ is a homomorphism: $\mu_\chi(gh)=\chi(g)\chi(h)gh=\mu_\chi(g)\mu_\chi(h)$.
  Moreover $\mu_\chi(\mu_\chi(g))=\chi(\chi(g)I_{p+q})\,\chi(g)g=g$, because $\chi(\pm I_{p+q})=1$ and $\chi(g)^2=1$.
  So $\mu_\chi$ is an involutive automorphism.
\end{proof}

\begin{corollary}\label{cor:outer-auts}
Let $p,q\ge 1$. We have $$\operatorname{Out}( \O(p,q)) \supseteq \begin{cases}
\lbrace 1, [\mu_{\det}]\rbrace & \mbox{ $p,q$ odd } \\
\lbrace 1, [\mu_{\det_p}] \rbrace, & \mbox{ $p$ even, $q$ odd } \\
\lbrace 1, [\mu_{\det_q}] \rbrace, & \mbox{ $p$ odd, $q$ even } \\
\lbrace 1, [\mu_{\det_p}], [\mu_{\det_q}], [\mu_{\det}] \rbrace, & \mbox{$p,q$ even}
\end{cases}$$
where the angle brackets denote equivalence classes of outer automorphisms.
\end{corollary}

\begin{proof} By Proposition~\ref{prop:easy-auts} these are automorphisms. Let $\chi\ne 1$ be one of the listed characters and let $r$ be a reflection, in a vector $v$ of suitable sign, with $\chi(r)=-1$. Then $\mu_\chi(r)=-r$. The $(-1)$-eigenspace of $r$ is the line $\RR v$ and that of $-r$ is $v^\perp$, so an isometry conjugating $r$ to $-r$ would map $\RR v$ isometrically onto $v^\perp$; this is impossible, by dimension if $p+q\ge 3$ and because $v^\perp$ and $\RR v$ have opposite signs if $p=q=1$. So $\mu_\chi$ is not inner, and since $\mu_\chi\mu_{\chi'}=\mu_{\chi\chi'}$ the listed classes are distinct.
\end{proof}

\begin{remark} It is proved in the MathOverflow thread cited above that these
  containments are equalities. However, we do not need this.
\end{remark}

\begin{proposition}\label{prop:so-po-same-n-odd}
  For $n$ odd, the groups $\SO(1,n+1)$ and $\PO(1,n+1)$ are isomorphic but not conjugate
  as subgroups of $\O(1,n+1)$.
\end{proposition}

\begin{proof}
  They are not conjugate, because they are distinct subgroups of index $2$ in $O(1,n+1)$,
  and subgroups of index $2$ are normal.
  They are isomorphic, because they are exchanged by $\mu_{\det_{n+1}}$.
\end{proof}

\begin{proposition}\label{prop:so-po-different-n-even}
  For $n$ even and positive, the groups $\SO(1,n+1)$ and $\PO(1,n+1)$ are not isomorphic.
\end{proposition}

\begin{remark}
  Before proving this, we remark that the outer automorphism $\mu_{\det{}}$
  in this case preserves $\SO(1,n+1)$, and this is already enough to prove that
  they are not isomorphic as subgroups of $\O(1,n+1)$.
\end{remark}

\begin{proof} 
  Note first that $\SO(1,n+1)$ has a nontrivial center generated by $-1$.
  On the other hand, we show that the center of $\PO(1,n+1)$ is trivial.
  Let $D$ be the diagonal matrix with entries $1, -1, -1, \dots, -1$. This
  belongs to $\PO(1,n+1)$, so
  $Z_{\PO(1,n+1)} \subseteq Z_{\PO(1,n+1)}(D) = \{1\} \times \O(n+1)$. On the other hand,
  the center of $\O(n+1)$ is $\pm 1$ and so
  $Z_{\PO(1,n+1)} \subseteq \{1,D\}$. It is easy to see that $D \notin Z_{\PO(1,n+1)}$
  (again, for $n > 0$): for $n = 2$ we have the matrix
  $$\begin{pmatrix}2&1&1&1\\1&0&1&1\\1&1&0&1\\1&1&1&0\end{pmatrix} \in \PO(1,n+1)$$
  that does not commute with $D$, and for larger $n$ we can use a diagonal block matrix
  whose blocks are this and $I_{n-2}$.
\end{proof}

\begin{remark}
  When $n+1$ is odd, $\SO(1,n+1)$ is preserved by every $\mu_\chi$.
  Indeed, if $g \in \SO(1,n+1)$, then $\det(\mu_{\chi}(g)) = \det( \chi(g)g) = \chi(g)^{n+2}\det(g)=\det(g)$, so $\mu_\chi(g)\in\SO(1,n+1)$. 
\end{remark}

\section{The Case of $\quat{-1,-1,-3}{\QQ}$} \label{sec:1-1-3}
In this section we consider the Clifford algebra $\quat{-1,-1,-3}{\QQ}$.
Throughout we will write $a_3$ for $\sqrt 3 i_3$.  
We write the standard basis of $\quat{-1,-1,-3}{\QQ}$ as
\[
\begin{aligned}
  v_1&=1, & v_2&=i_1, & v_3&=i_2, & v_4&=i_1i_2,\\
  v_5&=a_3, & v_6&=i_1a_3, & v_7&=i_2a_3, & v_8&=i_1i_2a_3.
\end{aligned}
\]

This is the order of generators chosen by \texttt{magma}; it is lexicographic
on the reversed bit strings, but is not ordered by the degree of elements
in the tensor algebra.

The Clifford order $\ZZ[i_1,i_2,a_3]= \quat{-1,-1,-3}{\ZZ}$ generated by the $i_1,i_2,a_3 = a_3$ and is contained in $4$ maximal orders,
with index $16$ in each case: 
$$ A(-1,-1,-3), \qquad B(-1,-1,-3)_0, \quad B(-1,-1,-3)_1, \quad B(-1,-1,-3)_2.$$
It will be convenient to introduce some special elements
$$\quad J_1=\frac{i_1 + a_3}{2}, \quad J_2=\frac{i_2 + a_3}{2}$$
$$\zeta_0=\frac{1 + a_3}{2}, \quad \zeta_1=\frac{1 + \sqrt{3}i_{13}}{2}, \quad \zeta_2=\frac{1 + \sqrt{3}i_{23}}{2}.$$ 
The elements satisfy $J_1^2=J_2^2=-1$ and $\zeta_0^6=\zeta_1^6=\zeta_2^6=1$.
Each of the orders $B(-1,-1,-3)_j$ is generated by its root of unity
$\zeta_j$ for $0 \le j \le 2$. In formulas, we have
$$B(-1,-1,-3)_j = \ZZ[i_1,i_2,a_3][\zeta_j].$$
Note that the element $\zeta_0$ is a Clifford vector, unlike $\zeta_1$ and $\zeta_2$.
More precisely,
\begin{align*}
  J_1&=a_3+i_1\zeta_1\in B(-1,-1,-3)_1,\\
  J_2&=a_3+i_2\zeta_2\in B(-1,-1,-3)_2.
\end{align*}
Thus these orders also contain Clifford vectors with nonintegral coordinates.

The orders $B(-1,-1,-3)_j$ are Clifford-conjugate.
Put $i_0=1$. Conjugation by $i_j+i_k$ sends
$B(-1,-1,-3)_j$ to $B(-1,-1,-3)_k$.

We can give a description for $A(-1,-1,-3)$ in a similar style with generators in terms of special named elements
$$A(-1,-1,-3)=\ZZ[i_1,i_2,a_3][\alpha,\beta,\gamma], \quad \alpha=\frac{1+i_{12} + \sqrt{3}(i_{13} + i_{23})}{2}, $$
$$\beta=\frac{i_1+i_{12}+\sqrt{3}(i_{13}+i_{123})}{2}, \quad \gamma = \frac{1+i_{1}+i_{2}+i_{12}}{2}.$$ 
The order $A(-1,-1,-3)$ does not contain any Clifford vectors with
nonintegral components and (as we will soon see) is not conjugate to the
others, so $\Vec(A(-1,-1,-3)) = \Vec(\ZZ[i_1,i_2,a_3])$.\footnote{
This is consistent with Remark~\ref{rem:need-n-gt-3}, since $\quat{-1,-1,-3}{\QQ}$ is a
totally definite quaternion algebra over its centre $\Q(\sqrt 3)$.  A similar argument would lead to the same conclusion if
$3$ were replaced by any positive integer congruent to $3 \bmod 4$.}
We prepare to prove this by determining the condition on an element
of $\quat{-1,-1,-3}{\QQ}$ to have rational square.
\begin{lemma}\label{lem:rat-square}
  Let $x = \sum_{j=1}^8 c_j v_j$. Then $x^2 \in \Q$ if and only if one of the
  following holds:
  \begin{enumerate}
  \item All $c_j$ except $c_1$ are $0$.
  \item All $c_j$ except $c_8$ are $0$.
  \item We have $c_1 = c_8 = 0$ and $c_2 c_7 - c_3 c_6 + c_4 c_5 = 0$.
  \end{enumerate}
  In the last case we have $x^2 = -(\sum_{j=2}^4 c_j^2 + 3 \sum_{k=5}^7 c_k^2)$.
\end{lemma}

\begin{proof} In $\quat{-1,-1,-3}{\QQ} \otimes_\QQ \QQ(c_1, \dots, c_8)$ we have
  $\sum_{j=1}^8 c_j v_j$ as a generic element whose
  square has coefficients in $\QQ(c_1, \dots, c_8)$. If $x \in \quat{-1,-1,-3}{\QQ}$,
then $x^2 \in \Q$ if and only if all coefficients of the square other than
  the coefficient of $1$ vanish on the coefficients of $x$.
  This defines a subscheme of $\PP^7(\Q)$.
  Computation in \texttt{magma} reveals that it has four irreducible components,
  three corresponding to the cases above and one defined by equations
  including $c_j^2 - 3c_{9-j}^2 = 0$ for $1 \le j \le 4$ and therefore devoid
  of rational points. The last assertion is an easy calculation.  
\end{proof}

For $i \in \{i_1,i_2,i_1i_2\}$, one checks that
$\frac{i + i_1 i_2 a_3}{2}$ has minimal polynomial $x^4 - x^2 + 1$ and lies in
$B(-1,-1,-3)_2$, $B(-1,-1,-3)_1$, $B(-1,-1,-3)_0$, respectively, so that
$\Z[\zeta_{12}]$ embeds into each $B(-1,-1,-3)_j$.
We will prove that $\ZZ[\zeta_{12}]$ does not embed into $A(-1,-1,-3)$.
Note that both types of maximal order have Clifford unit group of order
$24$, though with only one of the two groups having an element of
order $12$ they cannot be isomorphic.

\begin{lemma}\label{lem:roots-of-unity}
	Let $W_1, I_1$ be the subsets of $A(-1,-1,-3)$ of elements satisfying
	the equations $x^2+x+1 = 0, x^2+1 = 0$, respectively. Then $\#W_1 = 8$ and
	$\#I_1 = 6$, and no element of $W_1$ commutes with an element of $I_1$.
\end{lemma}

\begin{proof} If $w \in W_1$, then $2w+1 \in A(-1,-1,-3)$ and it
	satisfies the equation $x^2 + 3 = 0$. Thus we must be in the last case of
Lemma~\ref{lem:rat-square}, so that $\sum_{j=2}^4 c_j^2+3\sum_{k=5}^7 c_k^2=3$.
	The ring $A(-1,-1,-3)$ is a lattice containing $\ZZ[i_1,i_2,a_3]$ with some finite index $m$, so the coordinates $c_j$ of its elements lie in $\frac{1}{m}\ZZ$, and together with this equation only finitely many $(c_j)$ are possible.
	Among the integral ones, the possibilities are
	$c_2, c_3, c_4 \in \{\pm 1\}, c_5 = c_6 = c_7 = 0$ and for one of
	$c_5, c_6, c_7$ to be $\pm 1$ and all the other $c_i$ to be $0$. In the
	first case we obtain $8$ elements of $W_1$ from the choices of sign; the
	second does not give elements of $A$.
	There are also nonintegral solutions of the equation and of the condition of Lemma~\ref{lem:rat-square}, for example $c_2=\frac{3}{2}$, $c_5=\frac{1}{2}$ and all other $c_j=0$; a direct check against a $\ZZ$-basis of $A(-1,-1,-3)$ shows that none of them is of the form $2w+1$ with $w\in A(-1,-1,-3)$.
	The argument for $I_1$ is
	similar but simpler, the elements of $I_1$ being
	$\pm i_1, \pm i_2, \pm i_1 i_2$; again the nonintegral solutions, such as $c_2=c_5=\frac{1}{2}$, are excluded by the same check. It is now routine to verify the last
	claim.
\end{proof}

\begin{proposition}\label{prop:no-zeta-12}
	The maximal orders $B(-1,-1,-3)_1$ and $A(-1,-1,-3)$ are not isomorphic, and
	therefore are not Clifford conjugate.
\end{proposition}
\begin{proof} As already stated,
	the element $z = (i_2 + i_1 i_2 a_3)/2$ of $B(-1,-1,-3)_1$
	has minimal polynomial $x^4 - x^2 + 1$. As a result, the elements
	$z^4$ and $z^3$ have the minimal polynomials $x^2+x+1$ and $x^2+1$,
	respectively, and commute with each other. It follows that $A(-1,-1,-3)$ is
	not isomorphic to $B(-1,-1,-3)_1$, by Lemma~\ref{lem:roots-of-unity}.
\end{proof}

We are now ready to consider the fundamental domains of the respective orders.

\subsection{The Case of $B(-1,-1,-3)$}\label{sec:113-first}
  We begin with $B(-1,-1,-3)_0$, which is simpler because it is Clifford-Euclidean.
  This follows from Theorem~\ref{thm:small-covering-radius-euclidean}, the
  squared covering radius of $\Vec(B(-1,-1,-3)_0)$ being $5/6<1$. One can compute it from the
  fact that the lattice is the orthogonal direct sum of the hexagonal lattice
  with $\ZZ^2$; the squared covering radius of an orthogonal direct sum of lattices is
  the sum of the squared covering radii of the factors, here $1/3+1/2$.
  Alternatively, it is an immediate computation in \texttt{magma}.
  
\subsubsection{Clifford Unit Group}
The group of orthogonal transformations induced by the Clifford units has order $12$ and is isomorphic to the dihedral group of order $12$.
In terms of explicit generators we have
$$B(-1,-1,-3)_0^{\times} = \langle \alpha, i_1 \rangle, \quad \alpha = \frac{-i_{12} + i_{12}a_3}{2},$$
a group of order $24$.
Note that $i_2$ is contained in this set as $\alpha^3 = -i_{12}$, so $\alpha^3i_1 = -i_2$.

We have $r=\pi_{\alpha}$ of order $6$, acting by rotation by $2\pi/3$ in the
plane spanned by $1$ and $a_3$ and by $-1$ on the span of $i_1,i_2$.
The transformation $s = \pi_{i_1}$ is order $2$ and acts as $\pi_{i_1}(x_0+x_1i_1+x_2i_2+x_3a_3) = -x_0-x_1i_1+x_2i_2+x_3a_3$.
One checks that
$i_1\alpha i_1 = -\alpha^{-1}$ in $B(-1,-1,-3)_0$, so the actions on the Clifford vectors satisfy
$srs = r^{-1}$.  

\subsubsection{Fundamental Domain of $\Gamma_{\infty}$}
See Figure~\ref{fig:113-domains}(a),(b).
\begin{figure}[htbp!]
\centering
\begin{tikzpicture}[scale=2.4,lab]
 \draw[line width=0.8pt] (0,0.577)--(0.5,0.289)--(0.5,-0.289)--(0,-0.577)--(-0.5,-0.289)--(-0.5,0.289)--cycle;
 \fill[domC!35,draw=domC!80!black,line width=0.8pt] (0,0)--(0.5,0.289)--(0,0.577)--cycle;
 \foreach \p in {(0,0),(1,0),(-1,0),(0.5,0.866),(-0.5,0.866),(0.5,-0.866),(-0.5,-0.866)}{\fill \p circle (1.3pt);}
 \node[below] at (0,0){$0$}; \node[right] at (1,0){$1$}; \node[above] at (0.5,0.866){$\zeta_0$};
 \node[right] at (0.5,0.289){$\tfrac{3+a_3}6$}; \node[above left] at (0,0.577){$\tfrac{a_3}3$};
 \node at (0,-1.15){(a) $B_0$: $x_0a_3$-plane};
\end{tikzpicture}\quad
\begin{tikzpicture}[scale=2.4,lab]
 \draw[gray] (-0.5,-0.5) rectangle (0.5,0.5);
 \fill[domC!35,draw=domC!80!black,line width=0.8pt] (0,-0.5) rectangle (0.5,0.5);
 \fill (0,0) circle (1.2pt);
 \node[below] at (0.25,-0.5){$0\le x_1\le\tfrac12$}; \node[rotate=90,above] at (0.5,0){$|x_2|\le\tfrac12$};
 \node at (0,-1.0){(b) $B_0$: $x_1x_2$-plane};
\end{tikzpicture}\quad
\begin{tikzpicture}[scale=2.6,lab,x={(1cm,0cm)},y={(0.5cm,0.38cm)},z={(0cm,1cm)}]
\coordinate (O) at (0,0,0);
 \coordinate (a) at (0.5,0,-0.5); \coordinate (b) at (0.5,0.5,-0.5); \coordinate (c) at (0.5,0.5,0.5); \coordinate (d) at (0.5,0,0.5);
 \fill[domC!30] (O)--(a)--(d)--cycle; \fill[domC!45] (a)--(b)--(c)--(d)--cycle;
 \draw[domC!80!black,line width=0.8pt] (O)--(a)--(b)--(c)--(d)--cycle (O)--(c) (a)--(d);
 \draw[domC!80!black,dashed] (O)--(b);
 \fill (O) circle (1.2pt) node[left]{$0$};
 \draw[->,gray] (0.7,0,-0.5)--(0.95,0,-0.5) node[right]{$x_0$};
 \draw[->,gray] (0.7,0,-0.5)--(0.7,0.3,-0.5) node[right]{$x_1$};
 \draw[->,gray] (0.7,0,-0.5)--(0.7,0,-0.25) node[above]{$x_2$};
 \node at (0.3,0,-0.95){(c) $A$: $x_0x_1x_2$-space};
\end{tikzpicture}
\caption{Fundamental domains for $\Gamma_\infty$ in \S\ref{sec:1-1-3}, with $a_3=\sqrt3 i_3$. (a), (b) For $B(-1,-1,-3)_0$: the product of the triangle with vertices $0$, $(3+a_3)/6$, $a_3/3$ (orange) inside the hexagonal Voronoi cell of $\ZZ+\ZZ\zeta_0$ in the $x_0a_3$-plane, and the rectangle $0\le x_1\le\tfrac12$, $|x_2|\le\tfrac12$. (c) For $A(-1,-1,-3)$: the pyramid $x_0\ge|x_1|$, $x_0\ge|x_2|$, $x_1\ge0$, $x_0\le\tfrac12$, times the interval $|x_3|\le\tfrac12$ for the coefficient of $a_3$.}\label{fig:113-domains}
\end{figure}
The stabilizer of $\infty$ on the Clifford vectors contains the translations
by $1, i_1, i_2, \zeta_0$.  

A fundamental domain for these can be described by
the inequalities $|x_1|, |x_2| \le 1/2$ together with those that put
$x_0+x_3a_3$ in the regular hexagon with center $0$, inradius $1/2$ and vertices $\pm a_3/3$, $\pm(3\pm a_3)/6$, which is the Voronoi cell of $0$ in $\ZZ+\ZZ\zeta_0$.

The rotations $r^2,r^4$ (rotations by $\pm 2\pi/3$ of the plane of $1,a_3$, trivial on $x_1,x_2$) and $s$ (which acts on that plane by $x_0\mapsto -x_0$) generate a dihedral group of order $6$ acting on the hexagon, and they can be used to put $x_0+x_3a_3$ into
one of the six triangles making up the hexagon, namely the one with vertices
$0, (3 + a_3)/6, a_3/3$.

The element $r^3=\pi_{i_{12}}$ fixes this plane and negates $x_1$ and $x_2$, so we can then require $x_1 \ge 0$. Thus the fundamental domain is described
by the three inequalities giving the triangle with vertices $0, (3+a_3)/6, a_3/3$ in the $x_0a_3$-plane and by the inequalities $0 \le x_1 \le 1/2$,
 $-1/2 \le x_2 \le 1/2$ in the $x_1x_2$-plane.

\subsubsection{Fundamental Domain and Generators for $\Gamma$}
The only necessary hemisphere is the unit
hemisphere with center at the origin $B(0)$.
Indeed, let $B(a)$ be a hemisphere with $a\notin\Vec(B(-1,-1,-3)_0)$; its radius $r$ satisfies $r^{-2}=N$ for an integer $N\ge 2$.
If $N\ge 6$, then $r^2\le 1/6$, and the argument of Lemma~\ref{lem:O4-spheres}, with squared covering radius $5/6$, shows that $B(a)$ lies in the union of the unit hemispheres $B(\nu)$, $\nu\in\Vec(B(-1,-1,-3)_0)$.
For $2\le N\le 5$ the center $a$ lies in $\frac1N\Vec(B(-1,-1,-3)_0)$, and for each of the finitely many classes of such $a$ modulo $\Vec(B(-1,-1,-3)_0)$ an exact computation gives points $\nu_k\in\Vec(B(-1,-1,-3)_0)$ and weights $\lambda_k\ge 0$ with $\sum_k\lambda_k=1$, $\sum_k\lambda_k\nu_k=a$ and $\sum_k\lambda_k\vert\nu_k-a\vert^2\le 1-r^2$.
Since $\sum_k\lambda_k(1-\vert x-\nu_k\vert^2)-(r^2-\vert x-a\vert^2)=1-r^2-\sum_k\lambda_k\vert\nu_k-a\vert^2\ge 0$ for every $x$, again $B(a)\subseteq\bigcup_k B(\nu_k)$.
As the fundamental domain of $\Gamma_\infty$ lies in the Voronoi cell of $0$, only $B(0)$ meets the region above it.
As usual, the reflection across this hemisphere is given by the matrix
$S = \binom{0  \ -1}{ 1 \ \ 0 }$, so this together with the rotations $\pi_{i_1}$, $\pi_{\alpha}$ and translations 
generate $\Gamma = \PSL_2(B(-1,-1,-3)_0)$. Here is a full list:
$$ \tau_{1}, \tau_{i_1},\tau_{i_2}, \tau_{a_3}, \quad S, \quad \pi_{i_1}, \quad \pi_{\alpha}. $$
Our algorithm tells us that $S\tau_1$ has order $3$, while $S\tau_vS\tau_{-v}$ has order $3$ for
$v \in \{i_1,i_2,a_3\}$, and $S\tau_{\zeta_0}S\tau_{1-\zeta_0}$ also has order $3$.
In addition, we find that $S\pi_{\alpha}$ and $S\pi_{i_1}$ both have order $2$.

\subsection{The Case of $A(-1,-1,-3)$}\label{sec:113-second}
The order $A(-1,-1,-3)$ is not Clifford-Euclidean for the norm $N(x)=x\bar x$: for example, there is no
way to divide $1+i_1+a_3$ by $2$ with a smaller remainder.
    
    \subsubsection{Clifford Unit Group, Fundamental Domain of $\Gamma_{\infty}$, and Fundamental Domain of $\Gamma$ }
    The group $A(-1,-1,-3)^\times$ of Clifford units has order $24$, and its image under $\pi$, of order $12$, is isomorphic to $A_4$; it acts on the first
    three
    coordinates by cyclic permutations with an even number of negative signs
    (and trivially on the fourth coordinate, the coefficient of $a_3$).

    Every Clifford unit has coefficient $0$ for all generators of the
    Clifford algebra involving $a_3 = \sqrt{3}i_3$; the inclusion of the Hurwitz order
    (see \S\ref{S:hurwitz}) in
    the Hamilton quaternions into $A(-1,-1,-3)$ induces a bijection on
    groups of Clifford units.

    Therefore, a fundamental domain for $\Gamma_\infty$
is defined by $x_0 \ge |x_1|$, $x_0\ge |x_2|$, $x_1\ge 0$ and, as always, $|x_0|,|x_1|,|x_2| \le 1/2$, with the coefficient of $a_3$ between $-1/2$ and $1/2$. See Figure~\ref{fig:113-domains}(c).
    Since the pair $(\mu,\lambda)=(2,1+i_1+\sqrt{3}i_3)$ is unimodular\footnote{It is both left and right unimodular because it is invariant under $*$ and the order is invariant under $*$}, it is necessary to introduce a hemisphere $\mu^{-1}\lambda = (1+i_1+\sqrt{3}i_3)/2$ centered there in addition to the one based at the origin.
    
    Up to elements of the stabilizer of $\infty$, it appears that these
    are the only two bubbles that are needed; we have not proved this. The cusp $s=(1+i_1+a_3)/2$ is
    tidy (here $\bar\lambda\lambda=5\equiv 1 \bmod 2$), so we may use the matrix $M_s$ of Definition~\ref{def:ms-tidy} for it.
    Assuming that no other bubbles are needed, the generators are those of the stabilizer of infinity together
    with
    $$S=M_0 = \begin{pmatrix}0&1\\-1&0\end{pmatrix},\quad
      M_s = \begin{pmatrix}(-1+i_1+a_3)&2 \\ 2 & (-1-i_1-a_3)\end{pmatrix}.$$
    Let the stabilizer of infinity be generated by $\tau_1, \tau_{i_1}, \tau_{i_2}, \tau_{a_3}$,
    where $\tau_a$ is the matrix representing translation by $a$, and the
    additional elements $\pi_u = \begin{psmallmatrix}u&0\\0&(u^{-1})^*\end{psmallmatrix}$,
    where $u \in A(-1,-1,-3)^\times$.

    We do not determine the relations of $\SL_2(A(-1,-1,-3))$ here. The groups
    $\SL_2(A(-1,-1,-3))$ and $\SL_2(B(-1,-1,-3)_0)$ both contain $\SL_2(\ZZ[i_1,i_2,a_3])$, but neither
    contains the other: for the Hurwitz unit $\gamma=(1+i_1+i_2+i_{12})/2$ we have
    $\pi_\gamma\in\SL_2(A(-1,-1,-3))\setminus\SL_2(B(-1,-1,-3)_0)$, while
    $\pi_\alpha\in\SL_2(B(-1,-1,-3)_0)\setminus\SL_2(A(-1,-1,-3))$ for the unit $\alpha$ of \S\ref{sec:113-first}.
    So a presentation of one of these groups does not directly give one of the other.
    The matrix $M_s$, however, has entries in $\ZZ[i_1,i_2,a_3]$, so it also lies in $\SL_2(B(-1,-1,-3)_0)$,
    and we can express it in terms of the generators of \S\ref{sec:113-first}. It
    takes $\infty$ to $(-1+i_1+a_3)/2$. We
    precompose with translation by $(1-a_3)/2$ to reach $i_1/2$, then by
    $S$ to get to $2i_1$, then translation by $-2i_1$ and $S$ again to
    return to $\infty$. The relation
    $$M_s = -\tau_{(-1+a_3)/2}S\tau_{2i_1}M_0\pi_{i_1}\tau_{(-1-a_3)/2}$$ in $\SL_2(B(-1,-1,-3)_0)$ is
    then apparent from inspection of the matrix
    $S\tau_{-2i_1}S\tau_{(1-a_3)/2}M_s$.

\section{The Case of $\quat{-1,-1,-1,-1}{\QQ}$} \label{sec:1-1-1-1} 
The integral Clifford algebra  $\ZZ[i_1,i_2,i_3,i_4]$ is contained in $6$
maximal orders.  

Five of these have code of dimension $1$,
spanned by one of the five binary vectors of length $5$ and Hamming weight
$4$. These will be denoted by $\Ocal_{5,i}$ according to the position of the
zero: for example, $\Ocal_{5,2}$ is the order whose code is spanned by
$11011$, or equivalently that contains $(i_0+i_1+i_3+i_4)/2$.
Although the sum of four $i_1,i_2,i_3,i_4$ divided by $2$ generates an order
over the integral Clifford algebra that is not maximal, the
maximal order containing it is unique.

The code associated to the remaining maximal order is trivial; since this
order is special in various ways, we will denote it by $\Ocal_{5,!}$.

All six of the maximal orders are conjugate and hence isomorphic as abstract associative algebras: $(1+i_j)^{-1}\Ocal_{5,0}(1+i_j)=\Ocal_{5,j}$ for $1\le j\le 4$, and $\Ocal_{5,0}$ is conjugate to $\Ocal_{5,!}$ (Proposition~\ref{rem:not-invariant}).
Still, among these, all have distinct Clifford vectors, and there are
two distinct classes up to conjugacy by elements of the Clifford monoid.
One of the conjugacy classes is $\{\Ocal_{5,!}\}$, while
the other is $\{\Ocal_{5,i}: 0 \le i \le 4\}$.
To refer to the larger class
in the abstract, when it is unnecessary to distinguish
among the conjugates, we will call them simply $\Ocal_{[5,1,4]}$.

\begin{definition}\label{def:clifford-conjugacy} Let $\Ocal, \Ocal'$ be
  orders in a Clifford algebra $\Clf$. If $x^{-1} \Ocal x = \Ocal'$ for
  some $x \in \Clf^\mon$, then $\Ocal$ and $\Ocal'$ are
  {\em Clifford conjugate}.
\end{definition}

\begin{remark} In contrast to the situation for Clifford algebras with
  $3$ imaginary units, the maximal order of $\CC_5$ is unique up to
  conjugacy, as we proved in Proposition~\ref{prop:maximal-orders-unique}.
  However, the Clifford monoid is no longer a conjugation invariant and the
  order can interact with the Clifford vectors in many different ways, so
  we expect that there are infinitely many maximal orders up to Clifford conjugacy.
  Under the assumption that $x \in \Clf^\mon$, an element $u \in \Clf$ is
  a Clifford unit if and only if $x^{-1} u x$ is, so
  two Clifford conjugate orders have isomorphic groups of Clifford units.
  This is not true for orders that are conjugate but not Clifford conjugate
  (Proposition~\ref{rem:not-invariant}).
  We will not consider such examples further in this paper.
\end{remark}

\begin{proposition}\label{rem:not-invariant}
  There exist maximal orders $\Ocal, \Ocal' \subset \CC_5$ that are isomorphic,
  but not by any
  automorphism of the Clifford algebra that preserves the set of
  Clifford vectors, and for which 
  $\Ocal^\times$ and $\Ocal'^\times$ are not isomorphic.
\end{proposition}
\begin{proof}  
  This holds for $\Ocal= \Ocal_{5,0}$ and $\Ocal' = \Ocal_{5,!}$, for example.
  We compute that the groups of Clifford units have orders $\vert \Ocal_{[5,1,4]}^{\times}\vert = 1152$ and
  $\vert\Ocal_{5,!}^{\times}\vert=1920$.
Conjugation by an element of the Clifford monoid $\CC_5^\mon$ would preserve the structure of
the group of Clifford units, so this shows that $\Ocal_{[5,1,4]}$ and
$\Ocal_{5,!}$ are not conjugate by such an element. On the other hand,
the first part of Proposition~\ref{prop:maximal-orders-unique} shows that
they are conjugate. Indeed, a short computer search finds that
$v\Ocal_{5,0}v^{-1} = \Ocal_{5,!}$ for $v = 1 + i_2 + i_1i_2 + i_2i_3i_4$.
\end{proof}

\begin{remark} Note also that $\Ocal_{5,0}$ is Clifford-Euclidean by
  Theorem~\ref{thm:clifford-euclidean-orders}, while $\Ocal_{5,!}$ is not Clifford-Euclidean for the norm $x\bar x$,
  since it is impossible to divide $\sum_{j=0}^3 i_j$ by $2$ in $\Ocal_{5,!}$
  and obtain a smaller remainder. (In $\Ocal_{5,0}$ we would have
  $\sum_{j=0}^3 i_j = 2(\sum_{j=1}^4 i_j/2) + (i_0 - i_4)$.)
\end{remark}

\subsection{The Oddball Maximal Order in $\quat{-1,-1,-1,-1}{\ZZ}$} \label{sec:c5-oddball}
\subsubsection{Clifford Unit Group}

The group $\Ocal_{5,!}^\times$ of Clifford units is of order $1920$, and its
center is $\{-1,1\}$. This group acts with kernel $\pm 1$
on the boundary $\RR^5$ of hyperbolic space by linear
transformations; indeed, the action is by the group of signed permutation
matrices of determinant $1$ with even underlying permutation. Thus
$\Ocal_{5,!}^\times/\pm 1 \cong \Weyl(D_5)^+$, where $\Weyl(D_5)^+$ (Figure~\ref{fig:24-cell}(b)) refers to the 
subgroup of the Weyl group consisting of products of an even number of
reflections.

\subsubsection{Fundamental Domain of $\Gamma_{\infty}$}
This allows us to describe 
the fundamental domain of the stabilizer of $\infty$.

\begin{proposition}\label{prop:fd-inf-oddball}
  The fundamental domain of $\Gamma_\infty$ is bounded by the following
  hyperplanes:
  \begin{enumerate}
  \item  $x_i < 1/2$ for $0 \le i \le 4$;
  \item  $x_i > 0$ for $0 \le i \le 3$;
  \item $x_4 > -1/2$;
  \item $x_i > x_{i+1}$ for $0 \le i \le 2$;
  \item $x_2 > x_4, x_2 > -x_4$.
  \end{enumerate}
\end{proposition}

\begin{proof}
  Given a point $(x_0,\dots,x_4)$ in $\RR^5$, as usual we can
  translate so that $-1/2 \le x_i \le 1/2$ for $0 \le i \le 4$, and then
  we can apply an even permutation such that
  $|x_0| \ge |x_1| \ge |x_2| \ge |x_3|, |x_4|$ and a sign change after which
  $x_0, x_1, x_2, x_3 \ge 0$. Generically these are unique.
\end{proof}

\subsubsection{Fundamental Domain of $\Gamma$}
In order to determine the fundamental domain of $\Gamma$, we must study the
cusps. First, we note that the right ideal
$(1+i_1+i_2+i_3)\Ocal_{5,!} + 2\Ocal_{5,!}$ is of index $2^{12}$, which is not an
$8$th power, so this ideal cannot be generated by any single element of the
Clifford monoid $\Ocal_{5,!}^{\mon}$. (It is, however, generated by an element of
$\Ocal_{5,!}$ outside the Clifford monoid.)

We conjecture that this is essentially the only failure of $\Ocal_{5,!}$ to
be cuspidally principal. To be exact:
\begin{conjecture}\label{conj:o5-not-principal}
  Let $x, y \in \Vec(\Ocal_{5,!})$ such that $x\Ocal_{5,!} + y\Ocal_{5,!}$ is
  not generated by an element of $\Ocal_{5,!}^\mon$. Then for some $r \in \Ocal_{5,!}^\mon$, some $w\in\Vec(\Ocal_{5,!})$ and some
  $S \subset \{0,1,2,3,4\}$ with $|S| = 4$ we have
  $\{x,y\} = \{2r,r(\sum_{j \in S} i_j+2w)\}$.
\end{conjecture}

When the right ideal $I = xR + yR$ is generated by a single element of
$\Ocal_{5,!}^\mon$, we can find a generator by
finding the shortest vector in the lattice $I$ relative to the Clifford norm $x \bar x$,
which in this case coincides with the Euclidean norm. We have verified that
up to denominator $10$ every hemisphere is either a translate of that centered
at $0$ or $\sigma_\pm/2$, where $\sigma_\pm=1+i_1+i_2+i_3\pm i_4$, or is dominated by the union of such hemispheres.

\begin{proposition}\label{prop:zero-or-halves} Every point in the closure of the
  fundamental domain lies strictly under either
  the unit hemisphere or $B((1+i_1+i_2+i_3\pm i_4)/2,1/2)$, except for the
  equivalent points
  $(1/2,1/2,1/2,0,\pm 1/2)$ and $(1/2,1/2,1/2,1/2,0)$ which are on two or
  three of these.
\end{proposition}

\begin{proof} This is much like Lemma~\ref{lem:deephole}.
  Let the point be $(x_0,\dots,x_4)$. By symmetry we take $x_4 \ge 0$.
  Suppose that $\sum_{i=0}^4 x_i^2 \ge 1$ and $\sum_{i=0}^4 (1/2-x_i)^2 \ge 1/4$.
  Adding these two inequalities, we find $\sum_{i=0}^4 2x_i^2 - x_i \ge 0$.
  However, the maximum of $2x^2-x$ on $[0,1/2]$ is $0$, achieved at both
  endpoints. Thus all $x_i$ are $0$ or $1/2$. A simple check shows that
  only $(1/2,1/2,1/2,0,1/2)$ and $(1/2,1/2,1/2,1/2,0)$ are in the closure
  of the fundamental domain and satisfy both inequalities.
\end{proof}
  
\subsubsection{Generators}

Let $\sigma=\sigma_+$. The matrix $\binom{3 \ \ \bar \sigma}{ \sigma \ \ 2}$ belongs to
$\SL_2(\Ocal_{5,!})$, so $(\sigma,2)$ is unimodular.
If, as our computations up to denominator $10$ indicate, the hemispheres $B(0)$ and $B(\sigma_\pm/2)$ and their $\Gamma_\infty$-translates are the only ones needed, then the generators of $\SL_2(\Ocal_{5,!})$
consist of the generators of $\Gamma_\infty$ and two additional generators,
one for each of the spheres $B(0)$ and $B(\sigma/2)$. As usual,
the first of these is given by the matrix
$S=\binom{0 \ \ 1}{-1 \ 0}$. The second will be based on
the matrix $M_\sigma = \binom{3 \ \ \bar{\sigma}}{ \sigma \ \ 2}$
that we used above to show that $\sigma/2$ is a unimodular cusp.

\begin{proposition}\label{invert-cusp-o5}
  Let $\sigma=\sigma_+=1+i_1+i_2+i_3+i_4$, $\sigma_-=\sigma-2i_4$ and $s=\sigma/2$. The cusp $s$ is tidy, and the matrix
  $M_s=\begin{psmallmatrix}-\bar\sigma&2\\2&-\sigma\end{psmallmatrix}$ of Definition~\ref{def:ms-tidy}, the matrix
  $\pi_{i_4}=\begin{psmallmatrix}i_4&0\\0&-i_4\end{psmallmatrix}$ and $\tau_1$ lie in $\SL_2(\Ocal_{5,!})$.
  The element $\gamma_1=\tau_1\pi_{i_4}\tau_1M_s$ maps the hemisphere $B(\sigma_+/2)$ onto $B(\sigma_-/2)$, acting on it by
  $(x_0,\dots,x_4,t)\mapsto(x_0,x_1,x_2,x_3,-x_4,t)$, and it maps the inside of $B(\sigma_+/2)$ to the outside of $B(\sigma_-/2)$.
\end{proposition}

\begin{proof} Here $\bar\sigma\sigma-2\cdot 2=1$, so $s$ is tidy, and all entries lie in $\ZZ[i_1,\ldots,i_4]$.
  By Proposition~\ref{prop:cross-tidy-cusp}, $M_s$ maps $B(s)$ onto $B(s^\vee)$, where $s^\vee=-\bar\sigma/2=s-1$, and exchanges inside and outside;
  $\tau_1$ maps $B(s-1)$ back to $B(s)$; and $\tau_1\pi_{i_4}$, which acts by $(x_0,\dots,x_4,t)\mapsto(1-x_0,x_1,x_2,x_3,-x_4,t)$, maps $B(\sigma_+/2)$ onto $B(\sigma_-/2)$.
  A direct computation gives the stated formula for $\gamma_1$ on $B(\sigma_+/2)$.
\end{proof}

\begin{remark} Let $\bar D$ be the standard fundamental domain. By
  Proposition~\ref{prop:zero-or-halves},
  it has boundary components that are part of the hemispheres
  $B(\sigma/2)$ and $B(\sigma/2-i_4)$. Let $\bar D_1, \bar D_2$ be the
  translates of $\bar D$ adjacent to it across these components. We have
  just described the element $\gamma_1$ of $\SL_2(\Ocal_{5,!})$ taking $\bar D$
  to $\bar D_1$; the reader might expect an additional generator $\gamma_2$
  to be necessary to take $\bar D$ to $\bar D_2$. However, it can be checked
  that $\gamma_2 = \gamma_1^{-1}$, or equivalently that
  $\gamma_1(\bar D_2) = \bar D$. One way to verify this is to determine
  $\gamma_2$ directly, either from the formula of
  Proposition~\ref{prop:cross-tidy-cusp} or analogously to our description
  of $\gamma_1$.
\end{remark}

\subsection{The Case of $[5,1,4]$ Orders in $\quat{-1,-1,-1,-1}{\ZZ}$} \label{sec:c5-reasonable}

The five $[5,1,4]$ orders in $\quat{-1,-1,-1,-1}{\ZZ}$ are Clifford-conjugate,
so it suffices to consider one of them. For concreteness, we consider
$\Ocal_{5,0}$, the maximal order containing $\sum_{j=1}^4 i_j/2$.  

\subsubsection{Clifford Unit Group}

As in \S\ref{sec:1-1-1}, the group $\Ocal_{5,0}^\times$ of Clifford units is of order $1152$ and the
group of matrices giving its action on the space of Clifford
vectors is of order $576$. Its action on Euclidean space is again connected
with the $24$-cell (Figure~\ref{fig:24-cell}(a)).

\subsubsection{Fundamental Domain of $\Gamma_{\infty}$}

Similarly to what we saw in \S\ref{sec:1-1-1}, we have:

\begin{proposition} A fundamental domain for the action of $\Gamma_\infty$
  on $\Vec(\CC_5)$ is $\{(x_0,\dots,x_4): 0\le x_0\le 1/2,\ (x_1,x_2,x_3,x_4)\in F\}$, where $F$ is the domain of Corollary~\ref{cor:O4-F} in the coordinates $(a,b,c,d)=(x_1,x_2,x_3,x_4)$.
\end{proposition}

\begin{proof} We have $\Vec(\Ocal_{5,0})=\ZZ\oplus\Lambda'$, where $\Lambda'=\ZZ^4+\ZZ\frac{i_1+i_2+i_3+i_4}{2}$ corresponds to $\Vec(\Ocal_4)$. See Figure~\ref{fig:a1-d4}. A computation shows that the Clifford units act by $(x_0,y)\mapsto(\det(B)x_0,By)$, where $B$ runs over a group of order $576$ whose elements with $\det B=1$ form the image of $\Ocal_4^\times$ from \S\ref{sec:1-1-1}. Using translation by $1$ and an element with $\det B=-1$ we may assume $0\le x_0\le 1/2$, and then Corollary~\ref{cor:O4-F} applies to $y$.
\end{proof}
\begin{figure}[htbp!]
\centering
\begin{tikzpicture}[lab]
 \node[dnode,label={[dlab]below:$1$}] (e) at (-2.2*\dstep,0) {};
 \node[dnode,label={[dlab]right:$-\tfrac12(i_1+i_2+i_3+i_4)$}] (c) at (0,0) {};
 \node[dnode,label={[dlab]below:$i_2$}] (a) at ($(c)+(180:\dstep)$) {};
 \node[dnode,label={[dlab]right:$i_3$}] (b) at ($(c)+(60:\dstep)$) {};
 \node[dnode,label={[dlab]right:$i_4$}] (d) at ($(c)+(-60:\dstep)$) {};
 \draw[dedge] (c)--(a) (c)--(b) (c)--(d);
\end{tikzpicture}
\caption{The Dynkin diagram of $A_1\oplus D_4$ with simple roots in $\Vec(\Ocal_{5,0})=\ZZ\oplus\Lambda'$, where $\Lambda'$ is a copy of $\Vec(\Ocal_4)$.}\label{fig:a1-d4}
\end{figure}

\subsubsection{Fundamental Domain of $\Gamma$}

\begin{proposition} The only hemisphere needed in the construction of the
  fundamental domain is $B(0)$.
\end{proposition}

\begin{proof} The squared covering radius of $\Vec(\Ocal_{5,0})=\ZZ\oplus\Lambda'$ is $\frac14+\frac12=\frac34$, and the fundamental domain of $\Gamma_\infty$ lies in the Voronoi cell of $0$.
  Let $B(a)$ be a hemisphere with $a\notin\Vec(\Ocal_{5,0})$; its radius $r$ satisfies $r^{-2}=N$ for an integer $N\ge 2$.
  If $N\ge 4$, then $r^2\le\frac14$ and the argument of Lemma~\ref{lem:O4-spheres} shows that $B(a)$ lies in the union of the unit hemispheres centered at the points of $\Vec(\Ocal_{5,0})$.
  For $N\in\{2,3\}$, the center $a$ lies in $\frac1N\Vec(\Ocal_{5,0})$, and the exact computation described in \S\ref{sec:113-first} ($15$ classes of $a$ for $N=2$ and $80$ for $N=3$) shows the same.
  Consequently, only $B(0)$ meets the region above the fundamental domain of $\Gamma_\infty$.
\end{proof}

Since the squared covering radius of $\Vec(\Ocal_{5,0})$ is $3/4$,
the order is Clifford-Euclidean by Theorem~\ref{thm:clifford-euclidean-orders},
and the Euclidean division can be performed by solving a closest-vector
problem in this lattice, which is easy.

Thus the fundamental domain of $\Gamma$ is the region above the fundamental domain of $\Gamma_\infty$ and outside $B(0)$.

\subsubsection{Generators}
According to the description just given, the group is generated by
translations, the Clifford units, and the generator associated to $B(0)$. As usual, the generator for $B(0)$ is $z \to -1/z$.

\section{Additional Questions}\label{sec:open-problems}

We have collected a list of questions. 
\subsection{Class Number Problems}

In the case of imaginary quadratic fields, Gauss gave a series of conjectures concerning the behavior of class numbers for imaginary quadratic fields.
All of these questions have analogs for our orders. 
We start with some basic ones.

\begin{question}
  Are there finitely many Clifford-Euclidean orders? Is there an effective
  algorithm to determine them all, and can it be made practical?
  Same questions for cuspidally principal orders.
\end{question}

\begin{question}
  How does the class number of an order grow as a function of the
  dimension and discriminant? In particular, make the results of
  \cite[Chapters 26--29]{Voight2020} concrete in this situation,
  especially Theorem 26.2.3 and Main Theorem 29.10.1. Are there special
  families of Clifford algebras (for example, in which all the
  generators $i_j$ square to $-1$ or $-3$) that behave differently from the general case?
\end{question}

\begin{question}
Is there an order in the rational Clifford algebra of a positive definite quadratic form that is weakly Clifford-Euclidean but not Clifford-Euclidean? 
\end{question}

\subsection{Bubble Algorithmic Issues}
We have observed experimentally that all of the Clifford-Euclidean orders seem to have fundamental domains which require only a single bubble $B(0)$, up to the action of $\Gamma_{\infty}$.
We don't know of any \emph{Euclidean} examples with a bubble of radius not equal to one. 
Note that our non-Euclidean example $\Ocal_{5,!}$ (\S\ref{sec:c5-oddball}) needs two spheres, $B(0)$ and $B(\sigma/2)$.
\begin{conjecture}\label{conj:euclidean-bubbles}
	If $\Ocal$ is a Clifford-Euclidean order, then all boundary bubbles have radius 1 and there exists a fundamental domain whose bubbles are all $\Gamma_{\infty}$-translates of $B(0)$.
\end{conjecture}

In the course of our computations of the boundary spheres, we needed to know when to stop testing bubbles. 
In the Bianchi case there are bounds on the size of the bubbles in terms of the discriminant of the order \cite{Swan1971,Rahm2010}.
Does something like this work in the positive definite Clifford setting?
We do not know the answer even for cuspidally principal orders, where an argument along the lines of \cite{Swan1971} may apply; orders with singular cusps seem harder.
\begin{question}\label{qstn:bubble-bound}
Can one obtain an a priori bound on the curvature of the bubbles of $\PSL_2(\Ocal)$ that depends on the discriminant of $\Ocal$?
\end{question}

A set of generators for $\PSL_2(\Ocal)$ can be obtained as the set of elements that send the fundamental domain to each of its neighbors.
Determining how to cross the bubbles of the fundamental domain in general seems quite tricky.
This has been determined in the ``tidy'' case: see \S\ref{sec:tidy}.
How to do this in general is unclear.
If $B=B_g$ for a known $g\in\SL_2(\Ocal)$, then $g$ carries $B_g$ to $B_{g^{-1}}$; the difficulty is to find such a $g$, that is, to complete the pair $(c,d)$ defining $B$ to a matrix satisfying the conditions defining $\SL_2(\Ocal)$.
\begin{question}\label{qstn:crossing-bubbles}
Given a bubble $B$ in the fundamental domain $D$, what is the general procedure for finding a transformation $\gamma$ such that $\gamma(D)\cap D$ meets in the side defined by $B$? 
\end{question}

\subsection{Issues with Orders Closed Under Involutions}

We cannot enumerate the maximal orders of $\CC_n$ containing the
Clifford order for $n > 5$, but we can determine a single maximal
order containing a given order.  
Calculations for $\ZZ[i_1,\ldots,i_{n-1}]$ suggest the following conjecture:
\begin{conjecture}\label{conj:fixed-order-exists}
There always exists a maximal order containing $\ZZ[i_1,\ldots,i_{n-1}]$ which is closed under the involutions.
\end{conjecture}
This question makes sense for general quadratic forms but we have no data on these questions.

\begin{question}\label{q:most-orders-not-fixed}
Choose a maximal order $\Ocal$ of $\CC_n$ containing the Clifford order
uniformly at random. Does the probability of $\Ocal$ being fixed
by any of the involutions tend to $0$ as $n \to \infty$?
\end{question}

It might be more natural to choose $\Ocal$ with probability proportional
to $1/\#\Aut(\Ocal)$, and it is unclear what difference this will make for questions of this type.

\subsection{Algebraic Number Theory}
The following question has a positive answer if the Clifford algebra is an
	imaginary quadratic field, where it is the classical statement that the
	ideals of norm at most $X$ are equidistributed among the ideal classes
	as $X\to\infty$ (for prime ideals this is a case of the Chebotarev
	density theorem). It also has a positive answer if the Clifford algebra
	is a quaternion algebra, by counting lattice points in each class.
	It has some empirical support more generally.
	This should be compared to Remark~\ref{rem:cusps-come-from} where we discuss where cusps come from.
\begin{question}\label{conj:quaternion-chebotarev}
	Let $K$ be a rational Clifford algebra of a positive definite quadratic form and let
	$\Ocal_1, \dots, \Ocal_m$ be representatives for the set of isomorphism
	classes of maximal orders of $K$. Let ${\mathcal L}_1, \dots, {\mathcal L}_k$
	be the elements of the left ideal class set $\Cl(\Ocal_1)$, and define
	$r: \{1,\dots,k\} \to \{1,\dots,m\}$ so that the right
	order of elements preserving an ideal $L_i \in {\mathcal L}_i$ by right
	multiplication is isomorphic to $\Ocal_{r(i)}$. When the left ideals of $\Ocal_1$
	contained in $\Ocal_1$ are ordered by their index in $\Ocal_1$, are they
	distributed among the ${\mathcal L}_i$ in proportion to
	$1/\#\Ocal_{r(i)}^\times$?
\end{question}

\subsection{Lattices and Doubly Even Codes}
There are a number of outstanding problems related to the correspondence between lattices and codes which remains open.
Given a code $C$, form $\Ocal:=\ZZ[\frac{I \cdot c}{2}: c\in C]$. 
\begin{question}
When is $\Vec(\Ocal) = \frac{1}{2}\Lambda_C$?
\end{question}
We propose the following statement.
\begin{conjecture}\label{conj:generate-order}
	Let $\code$ be a doubly even code of length $n$. 
	Let $\Lambda_C$ be the inverse image of $C$ under the natural map $\ZZ + \ZZ i_1 + \cdots + \ZZ i_{n-1} \to \FF_2 + \FF_2 i_1 + \cdots + \FF_2 i_{n-1}$.
	Then $\ZZ[\frac{1}{2}\Lambda_C ]$ is an order in $\quat{(-1)^{n-1}}{\QQ}$.
\end{conjecture}

Certainly the $\sum_j c_j i_j/2$ are integral; the problem is to prove that
their products remain so. In the commutative setting this would be a
standard fact, but it is not clear here. Another way to say this
is that, if Conjecture~\ref{conj:generate-order} holds and $\Vec(\ZZ[\frac{1}{2}\Lambda_C])=\frac{1}{2}\Lambda_C$, then this construction shows that the image of the map of
Lemma~\ref{lem:order-to-code} includes all maximal codes.

\subsection{Arithmetic Hyperbolic Torsion}
An arithmetic hyperbolic reflection group is an arithmetic group $\Gamma \subset \PO_{1,n+1}(\RR)$ that is generated by reflections. 
It is a theorem of Vinberg \cite{Vinberg1984} that these can only exist in $\Hcal^{m}$ for $m<30$; for non-cocompact groups Esselmann \cite{Esselmann1996} showed that $m\leq 21$, and $m=21$ is attained by an example of Borcherds (see \cite{Belolipetsky2016} for a survey). 
An interesting question is to relax the condition on reflections and seek $\Gamma$ generated by torsion elements. 
Given that $\PSL_2(\ZZ)$ is generated by an element of order 2 and an element of order 3, one might ask if orders exist with this property. 
A systematic investigation of the groups $\SL_2(\Ocal)$ and their torsion or homological torsion is certainly worthy of future investigation. See for example \cite{Rahm2013} for the Bianchi case.

\subsection{Bott Periodicity and Even Unimodular Lattices}

We have observed that there is a special Clifford-Euclidean order $\Ocal_{H(8,4)}$ with $\Vec(\Ocal_{H(8,4)})=\frac{1}{2}\Lambda_{H(8,4)}$ (Theorem~\ref{thm:clifford-euclidean-orders}), which contains $\ZZ[\frac{c\cdot I}{2}\colon c\in H(8,4)]$ corresponding to a doubly even code, whose lattice $\Vec(\Ocal_{H(8,4)})$ is similar to $E_8$.
We suspect that there exist maximal orders corresponding to the Niemeier lattices that come from doubly even self-dual codes of length $24$, for instance to the Niemeier lattice with root system $24A_1$ (corresponding to the Golay code $G_{24}$); these are undoubtedly interesting objects. The Leech lattice does not occur in this way: if $\Ocal\supseteq\ZZ[i_1,\ldots,i_{23}]$, then by Lemma~\ref{lem:order-to-code} the lattice $\Vec(\Ocal)=\frac{1}{2}\Lambda_C$ has minimal squared length $1$ and covolume $2^{-\dim C}\geq 2^{-12}$, whereas the rescaling of the Leech lattice with minimal squared length $1$ has covolume $2^{-24}$. 
\begin{question}
	Do there exist maximal orders $\Ocal$ with $\Vec(\Ocal)$ being one of these special lattices in dimension 24?
\end{question}
The answer to the above question may simply be ``no''. 
Matrix algebras over Euclidean rings are Euclidean.
Does something similar hold for the condition of Clifford-Euclideanity in light of the Bott periodicity isomorphisms?
Can we say anything about the ``tameness'' of Niemeier orders, supposing they exist? 
\begin{question}
	Supposing that $\Ocal_{G_{24}}$ exists, what tameness properties does it have? 
\end{question}

Positive definite even unimodular lattices can only exist in dimension $n$ when $n \equiv 0 \mod 8$. 
It is an interesting question to ask if there is a ``geometric reason'' for this. 
See for example \cite[\href{https://math.stackexchange.com/questions/2058249/geometric-reason-why-even-unimodular-positive-definite-lattices-exist-only-in-di}{MSE 2058249}]{mathse-posdef}.
\begin{question}
	Can Bott periodicity and the geometry of Clifford orders be used to prove that positive definite even unimodular lattices can only exist in dimension $n$ when $n\equiv 0 \mod 8$? More generally, what is the connection between Bott periodicity and the dimensions of even unimodular lattices?
\end{question}

\subsection{Modular Symbols}
The first two authors are developing a theory of modular symbols and will report on this in a separate paper.

What is unclear to the authors is to what extent the Satake construction given in \S\ref{sec:satake-abelian} depends on the parameters chosen. 
In what sense are these locally symmetric spaces ``moduli of Satake abelian varieties'' or ``moduli of Hodge structures'', if at all? 
Is there a connection between the Satake construction and the Bianchi-Humbert interpretation of $\Hcal^{n+1}$ as positive definite Clifford-Hermitian forms up to scalars?

\bibliographystyle{alpha}
\bibliography{cliff-revised}

\begin{bibdiv}
\begin{biblist}

\bib{AtiyahBottShapiro1964}{article}{
      author={Atiyah, M.~F.},
      author={Bott, R.},
      author={Shapiro, A.},
       title={Clifford modules},
        date={1964},
     journal={Topology},
      volume={3},
      number={suppl. 1},
       pages={3\ndash 38},
}

\bib{Ahlfors1984}{article}{
      author={Ahlfors, Lars~V.},
       title={Old and new in {M}\"{o}bius groups},
        date={1984},
        ISSN={0066-1953},
     journal={Ann. Acad. Sci. Fenn. Ser. A I Math.},
      volume={9},
       pages={93\ndash 105},
         url={https://doi-org.ezproxy.uvm.edu/10.5186/aasfm.1984.0901},
      review={\MR{752394}},
}

\bib{Ahlfors1985}{incollection}{
      author={Ahlfors, Lars~V.},
       title={M\"{o}bius transformations and {C}lifford numbers},
        date={1985},
   booktitle={Differential geometry and complex analysis},
   publisher={Springer, Berlin},
       pages={65\ndash 73},
      review={\MR{780036}},
}

\bib{Alahmadi2014}{article}{
      author={Alahmadi, Adel},
      author={Jain, S.~K.},
      author={Lam, T.~Y.},
      author={Leroy, A.},
       title={Euclidean pairs and quasi-{E}uclidean rings},
        date={2014},
        ISSN={0021-8693,1090-266X},
     journal={J. Algebra},
      volume={406},
       pages={154\ndash 170},
         url={https://doi.org/10.1016/j.jalgebra.2014.02.009},
      review={\MR{3188333}},
}

\bib{Adem2007}{book}{
      author={Adem, Alejandro},
      author={Leida, Johann},
      author={Ruan, Yongbin},
       title={Orbifolds and stringy topology},
      series={Cambridge Tracts in Mathematics},
   publisher={Cambridge University Press, Cambridge},
        date={2007},
      volume={171},
        ISBN={978-0-521-87004-7; 0-521-87004-6},
         url={https://doi.org/10.1017/CBO9780511543081},
      review={\MR{2359514}},
}

\bib{AnthropicClaudeCode2026}{misc}{
      author={{Anthropic}},
       title={Claude code},
         how={\url{https://claude.com/claude-code}},
        note={Version 2.1.319, desktop app; model \texttt{claude-opus-5-5}},
}

\bib{AssmusPless1983}{article}{
      author={Assmus, Edward F.~{J}r.},
      author={Pless, Vera},
       title={On the covering radius of extremal self-dual codes},
        date={1983},
     journal={IEEE Trans. Inform. Theory},
      volume={29},
       pages={359\ndash 363},
}

\bib{Atiyah1967}{book}{
      author={Atiyah, M.~F.},
       title={K-theory},
   publisher={W. A. Benjamin, New York--Amsterdam},
        date={1967},
}

\bib{Auel2009}{book}{
      author={Auel, Asher},
       title={Cohomological invariants of line bundle-valued symmetric bilinear
  forms},
   publisher={ProQuest LLC, Ann Arbor, MI},
        date={2009},
        ISBN={978-1109-22433-7},
  url={http://gateway.proquest.com.ezproxy.uvm.edu/openurl?url_ver=Z39.88-2004&rft_val_fmt=info:ofi/fmt:kev:mtx:dissertation&res_dat=xri:pqdiss&rft_dat=xri:pqdiss:3363244},
        note={Thesis (Ph.D.)--University of Pennsylvania},
      review={\MR{2717684}},
}

\bib{Bass1974}{article}{
      author={Bass, Hyman},
       title={Clifford algebras and spinor norms over a commutative ring},
        date={1974},
        ISSN={0002-9327},
     journal={Amer. J. Math.},
      volume={96},
       pages={156\ndash 206},
         url={https://doi.org/10.2307/2373586},
      review={\MR{360645}},
}

\bib{BB2005}{book}{
      author={Bj\"orner, Anders},
      author={Brenti, Francesco},
       title={Combinatorics of {C}oxeter groups},
      series={Graduate Texts in Mathematics},
   publisher={Springer},
        date={2005},
      volume={231},
        ISBN={978-3540-442387},
}

\bib{Bulla2010}{article}{
      author={Bulla, E.},
      author={Constales, D.},
      author={Krau\ss{}har, R.~S.},
      author={Ryan, John},
       title={Dirac type operators for {S}pin manifolds associated to
  congruence subgroups of generalized modular groups},
        date={2010},
        ISSN={0075-4102,1435-5345},
     journal={J. Reine Angew. Math.},
      volume={643},
       pages={1\ndash 19},
         url={https://doi.org/10.1515/CRELLE.2010.042},
      review={\MR{2658187}},
}

\bib{Magma1997}{article}{
      author={Bosma, Wieb},
      author={Cannon, John},
      author={Playoust, Catherine},
       title={The {M}agma algebra system. {I}. {T}he user language},
        date={1997},
     journal={J. Symbolic Comput.},
      volume={24},
      number={3-4},
       pages={235\ndash 265},
        note={Version V2.29-2 used},
}

\bib{Behrens2002}{article}{
      author={Behrens, Mark~J.},
       title={A new proof of the {B}ott periodicity theorem},
        date={2002},
        ISSN={0166-8641,1879-3207},
     journal={Topology Appl.},
      volume={119},
      number={2},
       pages={167\ndash 183},
         url={https://doi.org/10.1016/S0166-8641(01)00060-8},
      review={\MR{1886093}},
}

\bib{Belolipetsky2016}{article}{
      author={Belolipetsky, Mikhail},
       title={Arithmetic hyperbolic reflection groups},
        date={2016},
        ISSN={0273-0979},
     journal={Bull. Amer. Math. Soc. (N.S.)},
      volume={53},
      number={3},
       pages={437\ndash 475},
         url={https://doi-org.ezproxy.uvm.edu/10.1090/bull/1530},
      review={\MR{3501796}},
}

\bib{Bridson1999}{book}{
      author={Bridson, Martin~R.},
      author={Haefliger, Andr{\'e}},
       title={Metric spaces of non-positive curvature},
    language={English},
      series={Grundlehren Math. Wiss.},
   publisher={Berlin: Springer},
        date={1999},
      volume={319},
        ISBN={3-540-64324-9},
}

\bib{Borel2006}{book}{
      author={Borel, Armand},
      author={Ji, Lizhen},
       title={Compactifications of symmetric and locally symmetric spaces},
      series={Mathematics: Theory \& Applications},
   publisher={Birkh\"{a}user Boston, Inc., Boston, MA},
        date={2006},
        ISBN={978-0-8176-3247-2; 0-8176-3247-6},
      review={\MR{2189882}},
}

\bib{Bosch1990}{book}{
      author={Bosch, Siegfried},
      author={L\"{u}tkebohmert, Werner},
      author={Raynaud, Michel},
       title={N\'{e}ron models},
      series={Ergebnisse der Mathematik und ihrer Grenzgebiete (3) [Results in
  Mathematics and Related Areas (3)]},
   publisher={Springer-Verlag, Berlin},
        date={1990},
      volume={21},
        ISBN={3-540-50587-3},
         url={https://doi.org/10.1007/978-3-642-51438-8},
      review={\MR{1045822}},
}

\bib{Boothby1986}{book}{
      author={Boothby, William~M.},
       title={An introduction to differentiable manifolds and {R}iemannian
  geometry},
     edition={Second},
      series={Pure and Applied Mathematics},
   publisher={Academic Press, Inc., Orlando, FL},
        date={1986},
      volume={120},
        ISBN={0-12-116052-1; 0-12-116053-X},
      review={\MR{861409}},
}

\bib{Borel2019}{book}{
      author={Borel, Armand},
       title={Introduction to arithmetic groups},
      series={University Lecture Series},
   publisher={American Mathematical Society, Providence, RI},
        date={2019},
      volume={73},
        ISBN={978-1-4704-5231-5},
         url={https://doi.org/10.1090/ulect/073},
        note={Translated from the 1969 French original [MR0244260] by Lam
  Laurent Pham, Edited and with a preface by Dave Witte Morris},
      review={\MR{3970984}},
}

\bib{Brungs1973}{article}{
      author={Brungs, H.~H.},
       title={Left {E}uclidean rings},
        date={1973},
        ISSN={0030-8730,1945-5844},
     journal={Pacific J. Math.},
      volume={45},
       pages={27\ndash 33},
         url={http://projecteuclid.org/euclid.pjm/1102947704},
      review={\MR{316490}},
}

\bib{Borel1973}{article}{
      author={Borel, A.},
      author={Serre, J.-P.},
       title={Corners and arithmetic groups},
        date={1973},
        ISSN={0010-2571,1420-8946},
     journal={Comment. Math. Helv.},
      volume={48},
       pages={436\ndash 491},
         url={https://doi.org/10.1007/BF02566134},
      review={\MR{387495}},
}

\bib{Bygott1998}{thesis}{
      author={Bygott, Jeremy},
       title={Modular forms and modular symbols over imaginary quadratic
  fields},
        type={Ph.D. Thesis},
        date={1998},
}

\bib{Constales2013}{article}{
      author={Constales, D.},
      author={Grob, D.},
      author={Krausshar, R.~S.},
       title={A new class of hypercomplex analytic cusp forms},
        date={2013},
        ISSN={0002-9947,1088-6850},
     journal={Trans. Amer. Math. Soc.},
      volume={365},
      number={2},
       pages={811\ndash 835},
         url={https://doi.org/10.1090/S0002-9947-2012-05613-3},
      review={\MR{2995374}},
}

\bib{CliffLean2026}{misc}{
      author={{Clifford-Bianchi Lean project}},
       title={Lean companion to {The Basic Theory of Clifford-Bianchi Groups
  for Hyperbolic n-Space}},
         how={\url{https://github.com/tdupu/cliff-lean}},
        date={2026},
        note={Public repository, revision
  \texttt{19c129c40d331d53695b2d12bb0aa6806a0a9b44}},
}

\bib{Conrad2012}{misc}{
      author={Conrad, Brian},
       title={Clifford algebras},
         how={Course notes for Math 210C, available at
  \url{https://math.stanford.edu/~conrad/210CPage/handouts/clifford.pdf}},
        date={2012},
        note={Accessed: 2025-01-07},
}

\bib{Conrad2014}{incollection}{
      author={Conrad, Brian},
       title={Reductive group schemes},
        date={2014},
   booktitle={Autour des sch\'{e}mas en groupes. {V}ol. {I}},
      series={Panor. Synth\`eses},
      volume={42/43},
   publisher={Soc. Math. France, Paris},
       pages={93\ndash 444},
      review={\MR{3362641}},
}

\bib{Coxeter1973}{book}{
      author={Coxeter, H. S.~M.},
       title={Regular polytopes},
     edition={Third},
   publisher={Dover Publications, New York},
        date={1973},
}

\bib{Cremona1984}{article}{
      author={Cremona, J.~E.},
       title={Hyperbolic tessellations, modular symbols, and elliptic curves
  over complex quadratic fields},
        date={1984},
        ISSN={0010-437X,1570-5846},
     journal={Compositio Math.},
      volume={51},
      number={3},
       pages={275\ndash 324},
         url={http://www.numdam.org/item?id=CM_1984__51_3_275_0},
      review={\MR{743014}},
}

\bib{SPLAG}{book}{
      author={Conway, J.~H.},
      author={Sloane, N. J.~A.},
       title={Sphere packings, lattices and groups},
    language={English},
     edition={3rd ed.},
      series={Grundlehren Math. Wiss.},
   publisher={New York, NY: Springer},
        date={1999},
      volume={290},
        ISBN={0-387-98585-9},
        note={{With} additional contributions by {E}. {Bannai}, {R}. {E}.
  {Borcherds}, {J}. {Leech}, {S}. {P}. {Norton}, {A}. {M}. {Odlyzko}, {R}. {A}.
  {Parker}, {L}. {Queen} and {B}. {B}. {Venkov}.},
}

\bib{Cremona1994}{article}{
      author={Cremona, J.~E.},
      author={Whitley, E.},
       title={Periods of cusp forms and elliptic curves over imaginary
  quadratic fields},
        date={1994},
        ISSN={0025-5718,1088-6842},
     journal={Math. Comp.},
      volume={62},
      number={205},
       pages={407\ndash 429},
         url={https://doi.org/10.2307/2153419},
      review={\MR{1185241}},
}

\bib{Dupuy2023}{article}{
      author={Dupuy, Taylor},
      author={Kedlaya, Kiran~S.},
      author={Zureick-Brown, David},
       title={Angle ranks of abelian varieties},
        date={2023},
     journal={Mathematische Annalen},
       pages={1\ndash 17},
}

\bib{MagmaCliffordAlgebras2026}{misc}{
      author={Dupuy, Taylor},
      author={Logan, Adam},
       title={magma-clifford-algebras: a {M}agma package for {C}lifford
  algebras and {C}lifford--{B}ianchi groups},
         how={\url{https://github.com/tdupu/magma-clifford-algebras}},
        date={2026},
        note={Public repository, revision
  \texttt{51cfdeb0135c19d2f660ca985a8865bf91d3f65d}},
}

\bib{DHIL2}{unpublished}{
      author={Dupuy, Taylor},
      author={Logan, Adam},
       title={Results on {{${\mathrm{SL}}_2$}} of orders in quaternion
  algebras},
        note={In preparation},
}

\bib{Lean4}{inproceedings}{
      author={de~Moura, Leonardo},
      author={Ullrich, Sebastian},
       title={The {Lean} 4 theorem prover and programming language},
        date={2021},
   booktitle={Automated deduction -- cade 28},
      series={Lecture Notes in Computer Science},
      volume={12699},
   publisher={Springer},
       pages={625\ndash 635},
        note={Version used: v4.35.0},
}

\bib{AgentSkillsPublic}{misc}{
      author={Dupuy, Taylor},
       title={agent-skills-public: reusable workflows for mathematical writing
  and verification},
         how={\url{https://github.com/tdupu/agent-skills-public}},
        note={Public repository, revision
  \texttt{8d13efb5e66306aeae6e3b948bc463e7e2db1b65}},
}

\bib{Elstrodt1987}{article}{
      author={Elstrodt, J.},
      author={Grunewald, F.},
      author={Mennicke, J.},
       title={Vahlen's group of {C}lifford matrices and spin-groups},
        date={1987},
        ISSN={0025-5874},
     journal={Math. Z.},
      volume={196},
      number={3},
       pages={369\ndash 390},
         url={https://doi-org.ezproxy.uvm.edu/10.1007/BF01200359},
      review={\MR{913663}},
}

\bib{Elstrodt1988}{article}{
      author={Elstrodt, J.},
      author={Grunewald, F.},
      author={Mennicke, J.},
       title={Arithmetic applications of the hyperbolic lattice point theorem},
        date={1988},
        ISSN={0024-6115},
     journal={Proc. London Math. Soc. (3)},
      volume={57},
      number={2},
       pages={239\ndash 283},
         url={https://doi-org.ezproxy.uvm.edu/10.1112/plms/s3-57.2.239},
      review={\MR{950591}},
}

\bib{Elstrodt1990}{article}{
      author={Elstrodt, J.},
      author={Grunewald, F.},
      author={Mennicke, J.},
       title={Kloosterman sums for {C}lifford algebras and a lower bound for
  the positive eigenvalues of the {L}aplacian for congruence subgroups acting
  on hyperbolic spaces},
        date={1990},
        ISSN={0020-9910},
     journal={Invent. Math.},
      volume={101},
      number={3},
       pages={641\ndash 685},
         url={https://doi-org.ezproxy.uvm.edu/10.1007/BF01231519},
      review={\MR{1062799}},
}

\bib{EKM2008}{book}{
      author={Elman, Richard},
      author={Karpenko, Nikita},
      author={Merkurjev, Alexander},
       title={The algebraic and geometric theory of quadratic forms},
      series={Colloquium Publications},
   publisher={American Mathematical Society, Providence, RI},
        date={2008},
      volume={56},
}

\bib{Esselmann1996}{article}{
      author={Esselmann, Frank},
       title={{\"U}ber die maximale {D}imension von {L}orentz-{G}ittern mit
  coendlicher {S}piegelungsgruppe},
        date={1996},
     journal={J. Number Theory},
      volume={61},
      number={1},
       pages={103\ndash 144},
}

\bib{FGT2010}{article}{
      author={Finis, Tobias},
      author={Grunewald, Fritz},
      author={Tirao, Paulo},
       title={{The Cohomology of Lattices in ${\SL}(2,\mathbb{C})$}},
        date={2010},
     journal={Experimental Mathematics},
      volume={19},
      number={1},
       pages={29 \ndash  63},
}

\bib{Fpylll}{misc}{
      author={{FPLLL Development Team}},
       title={{fpylll}: a python interface for the {fplll} lattice-reduction
  library},
         how={\url{https://github.com/fplll/fpylll}},
        note={Version used not recorded},
}

\bib{GAP4}{misc}{
      author={{GAP Group}},
       title={{GAP} -- groups, algorithms, and programming},
         how={\url{https://www.gap-system.org}},
        note={Version 4.14.0},
}

\bib{Grob2015}{article}{
      author={Grob, D.},
      author={Krau\ss{}har, R.~S.},
       title={A {S}elberg trace formula for hypercomplex analytic cusp forms},
        date={2015},
        ISSN={0022-314X,1096-1658},
     journal={J. Number Theory},
      volume={148},
       pages={398\ndash 428},
         url={https://doi.org/10.1016/j.jnt.2014.09.002},
      review={\MR{3283187}},
}

\bib{Hahn2004}{incollection}{
      author={Hahn, Alexander},
       title={The {C}lifford algebra in the theory of algebras, quadratic
  forms, and classical groups},
        date={2004},
   booktitle={Clifford algebras ({C}ookeville, {TN}, 2002)},
      series={Prog. Math. Phys.},
      volume={34},
   publisher={Birkh\"{a}user Boston, Boston, MA},
       pages={305\ndash 322},
      review={\MR{2025987}},
}

\bib{Hatcher2002}{book}{
      author={Hatcher, Allen},
       title={Algebraic topology},
   publisher={Cambridge University Press, Cambridge},
        date={2002},
        ISBN={0-521-79160-X; 0-521-79540-0},
      review={\MR{1867354}},
}

\bib{Hochschild1965}{book}{
      author={Hochschild, G.},
       title={The structure of {L}ie groups},
   publisher={Holden-Day, San Francisco},
        date={1965},
}

\bib{Humphreys1978}{book}{
      author={Humphreys, James~E.},
       title={Introduction to {L}ie algebras and representation theory},
      series={Graduate Texts in Mathematics},
   publisher={Springer-Verlag, New York-Berlin},
        date={1978},
      volume={9},
        ISBN={0-387-90053-5},
        note={Second printing, revised},
      review={\MR{499562}},
}

\bib{Iga2021}{article}{
      author={Iga, Kevin},
       title={Adinkras: graphs of clifford algebra representations,
  supersymmetry, and codes},
        date={2021},
     journal={arXiv preprint arXiv:2110.01665},
}

\bib{IvanyosRonyai}{article}{
      author={Ivanyos, G{\'a}bor},
      author={R{\'o}nyai, Lajos},
       title={Finding maximal orders in semisimple algebras over
  {{\(\mathbb{Q}\)}}},
    language={English},
        date={1993},
        ISSN={1016-3328},
     journal={Comput. Complexity},
      volume={3},
      number={3},
       pages={245\ndash 261},
}

\bib{Ji2002}{article}{
      author={Ji, L.},
      author={MacPherson, R.},
       title={Geometry of compactifications of locally symmetric spaces},
        date={2002},
        ISSN={0373-0956,1777-5310},
     journal={Ann. Inst. Fourier (Grenoble)},
      volume={52},
      number={2},
       pages={457\ndash 559},
         url={http://aif.cedram.org/item?id=AIF_2002__52_2_457_0},
      review={\MR{1906482}},
}

\bib{Joyce2012}{incollection}{
      author={Joyce, Dominic},
       title={On manifolds with corners},
        date={2012},
   booktitle={Advances in geometric analysis},
      series={Adv. Lect. Math. (ALM)},
      volume={21},
   publisher={Int. Press, Somerville, MA},
       pages={225\ndash 258},
      review={\MR{3077259}},
}

\bib{Joyce2014}{incollection}{
      author={Joyce, Dominic},
       title={An introduction to d-manifolds and derived differential
  geometry},
        date={2014},
   booktitle={Moduli spaces},
      series={London Math. Soc. Lecture Note Ser.},
      volume={411},
   publisher={Cambridge Univ. Press, Cambridge},
       pages={230\ndash 281},
      review={\MR{3221297}},
}

\bib{Knus1991}{book}{
      author={Knus, Max-Albert},
       title={Quadratic and {H}ermitian forms over rings},
      series={Grundlehren der mathematischen Wissenschaften [Fundamental
  Principles of Mathematical Sciences]},
   publisher={Springer-Verlag, Berlin},
        date={1991},
      volume={294},
        ISBN={3-540-52117-8},
         url={https://doi-org.ezproxy.uvm.edu/10.1007/978-3-642-75401-2},
        note={With a foreword by I. Bertuccioni},
      review={\MR{1096299}},
}

\bib{Krausshar2004}{book}{
      author={Krausshar, Rolf~S\"{o}ren},
       title={Generalized analytic automorphic forms in hypercomplex spaces},
      series={Frontiers in Mathematics},
   publisher={Birkh\"{a}user Verlag, Basel},
        date={2004},
        ISBN={3-7643-7059-9},
         url={https://doi-org.ezproxy.uvm.edu/10.1007/b95203},
      review={\MR{2037622}},
}

\bib{Lee2013SmoothManifolds}{book}{
      author={Lee, John~M.},
       title={Introduction to smooth manifolds},
     edition={2},
      series={Graduate Texts in Mathematics},
   publisher={Springer},
        date={2013},
      volume={218},
}

\bib{Lingham2005}{thesis}{
      author={Lingham, Mark},
       title={Modular forms and elliptic curves over imaginary quadratic
  fields},
        type={Ph.D. Thesis},
        date={2005},
}

\bib{Lounesto2001}{book}{
      author={Lounesto, Pertti},
       title={Clifford algebras and spinors},
     edition={Second},
      series={London Mathematical Society Lecture Note Series},
   publisher={Cambridge University Press, Cambridge},
        date={2001},
      volume={286},
        ISBN={0-521-00551-5},
         url={https://doi-org.ezproxy.uvm.edu/10.1017/CBO9780511526022},
      review={\MR{1834977}},
}

\bib{Li1987}{article}{
      author={Li, Jian-Shu},
      author={Piatetski-Shapiro, I.},
      author={Sarnak, P.},
       title={Poincar\'{e} series for {${\rm SO}(n,1)$}},
        date={1987},
        ISSN={0253-4142,0973-7685},
     journal={Proc. Indian Acad. Sci. Math. Sci.},
      volume={97},
      number={1-3},
       pages={231\ndash 237},
         url={https://doi.org/10.1007/BF02837825},
      review={\MR{983616}},
}

\bib{Maass1949}{inproceedings}{
      author={Maass, Hans},
       title={Automorphe {F}unktionen von mehreren {V}er{\"a}nderlichen und
  {D}irichletsche {R}eihen},
organization={Springer},
        date={1949},
   booktitle={Abhandlungen aus dem {M}athematischen {S}eminar der
  {U}niversit\"at {H}amburg},
      volume={16},
       pages={72\ndash 100},
}

\bib{mathoverflow-grothendieck-topology}{misc}{
      author={{Mathoverflow authors}, The},
       title={Mathoverflow},
         how={\url{https://mathoverflow.net}},
        date={2010--2019},
}

\bib{mathse-posdef}{misc}{
      author={{Mathematics Stack Exchange authors}, The},
       title={Mathematics {S}tack {E}xchange},
         how={\url{https://math.stackexchange.com}},
        date={2016},
}

\bib{mathoverflow-autgp}{misc}{
      author={{Mathoverflow authors}, The},
       title={Mathoverflow},
         how={\url{https://mathoverflow.net}},
        date={2016},
}

\bib{math-stack-exchange-left-translation}{misc}{
      author={{Mathematics Stack Exchange authors}, The},
       title={Mathematics {S}tack {E}xchange},
         how={\url{https://math.stackexchange.com}},
        date={2019},
}

\bib{Mathlib2020}{inproceedings}{
      author={{Mathlib Community}},
       title={The {Lean} mathematical library},
        date={2020},
   booktitle={Proceedings of the 9th acm sigplan international conference on
  certified programs and proofs (cpp 2020)},
   publisher={ACM},
       pages={367\ndash 381},
        note={Revision used:
  \texttt{074063781836922cb02701cd12943c9f51731afa}},
}

\bib{McInroy2016}{article}{
      author={McInroy, Justin},
       title={Vahlen groups defined over commutative rings},
        date={2016},
        ISSN={0025-5874},
     journal={Math. Z.},
      volume={284},
      number={3-4},
       pages={901\ndash 917},
         url={https://doi-org.ezproxy.uvm.edu/10.1007/s00209-016-1678-x},
      review={\MR{3563259}},
}

\bib{Milnor1963}{book}{
      author={Milnor, J.},
       title={Morse theory},
      series={Annals of Mathematics Studies},
   publisher={Princeton University Press, Princeton, NJ},
        date={1963},
      volume={51},
}

\bib{Miller2024}{misc}{
      author={Miller, Robert~L.},
       title={{D}oubly-{E}ven {C}odes},
         how={\url{https://rlmill.github.io/de_codes/}},
        note={Accessed April 24, 2024.},
}

\bib{Morris2015}{book}{
      author={Morris, Dave~Witte},
       title={Introduction to arithmetic groups},
   publisher={Deductive Press, [place of publication not identified]},
        date={2015},
        ISBN={978-0-9865716-0-2; 978-0-9865716-1-9},
      review={\MR{3307755}},
}

\bib{Maclachlan2003}{book}{
      author={Maclachlan, Colin},
      author={Reid, Alan~W.},
       title={The arithmetic of hyperbolic 3-manifolds},
      series={Graduate Texts in Mathematics},
   publisher={Springer-Verlag, New York},
        date={2003},
      volume={219},
        ISBN={0-387-98386-4},
         url={https://doi-org.ezproxy.uvm.edu/10.1007/978-1-4757-6720-9},
      review={\MR{1937957}},
}

\bib{Moosa2010}{article}{
      author={Moosa, Rahim},
      author={Scanlon, Thomas},
       title={Jet and prolongation spaces},
        date={2010},
        ISSN={1474-7480,1475-3030},
     journal={J. Inst. Math. Jussieu},
      volume={9},
      number={2},
       pages={391\ndash 430},
         url={https://doi.org/10.1017/S1474748010000010},
      review={\MR{2602031}},
}

\bib{Maclachlan1989}{article}{
      author={Maclachlan, C.},
      author={Waterman, P.~L.},
      author={Wielenberg, N.~J.},
       title={Higher-dimensional analogues of the modular and {P}icard groups},
        date={1989},
        ISSN={0002-9947},
     journal={Trans. Amer. Math. Soc.},
      volume={312},
      number={2},
       pages={739\ndash 753},
         url={https://doi-org.ezproxy.uvm.edu/10.2307/2001009},
      review={\MR{965301}},
}

\bib{Neukirch1999}{book}{
      author={Neukirch, J\"urgen},
       title={Algebraic number theory},
      series={Grundlehren Math. Wiss.},
   publisher={Springer},
        date={1999},
      volume={322},
        ISBN={978-3-642-08473-7},
         url={https://doi.org/10.1007/978-3-662-0398},
        note={Translated from the 1992 German original by Norbert Schappacher},
}

\bib{NLabQF}{misc}{
      author={{nLab authors}},
       title={Quadratic form},
         how={\url{https://ncatlab.org/nlab/show/quadratic+form}},
        date={2024},
        note={\url{https://ncatlab.org/nlab/show/quadratic+form}, Accessed July
  18, 2024.},
}

\bib{OpenAICodex2026}{misc}{
      author={{OpenAI}},
       title={Codex},
         how={\url{https://openai.com/codex/}},
        note={Model identifier not recorded},
}

\bib{Porteous1995}{book}{
      author={Porteous, Ian~R.},
       title={Clifford algebras and the classical groups},
      series={Cambridge Studies in Advanced Mathematics},
   publisher={Cambridge University Press, Cambridge},
        date={1995},
      volume={50},
        ISBN={0-521-55177-3},
         url={https://doi-org.ezproxy.uvm.edu/10.1017/CBO9780511470912},
      review={\MR{1369094}},
}

\bib{Platonov1994}{book}{
      author={Platonov, Vladimir},
      author={Rapinchuk, Andrei},
       title={Algebraic groups and number theory},
      series={Pure and Applied Mathematics},
   publisher={Academic Press, Inc., Boston, MA},
        date={1994},
      volume={139},
        ISBN={0-12-558180-7},
        note={Translated from the 1991 Russian original by Rachel Rowen},
      review={\MR{1278263}},
}

\bib{Python2026}{misc}{
      author={{Python Software Foundation}},
       title={Python},
         how={\url{https://www.python.org/}},
        note={Version 3.13.7},
}

\bib{Rahm2010}{thesis}{
      author={Rahm, Alexander},
       title={({C}o)homologies and {K}-theory of {B}ianchi groups using
  computational geometric models},
        type={Ph.D. Thesis},
        date={2010},
}

\bib{Rahm2013}{article}{
      author={Rahm, Alexander~D.},
       title={The homological torsion of {$\rm{PSL}_2$} of the imaginary
  quadratic integers},
        date={2013},
        ISSN={0002-9947,1088-6850},
     journal={Trans. Amer. Math. Soc.},
      volume={365},
      number={3},
       pages={1603\ndash 1635},
         url={https://doi.org/10.1090/S0002-9947-2012-05690-X},
      review={\MR{3003276}},
}

\bib{Ratcliffe2019}{book}{
      author={Ratcliffe, John~G.},
       title={Foundations of hyperbolic manifolds},
      series={Graduate Texts in Mathematics},
   publisher={Springer, Cham},
        date={2019},
      volume={149},
        ISBN={978-3-030-31597-9; 978-3-030-31596-2},
         url={https://doi-org.ezproxy.uvm.edu/10.1007/978-3-030-31597-9},
        note={Third edition [of 1299730]},
      review={\MR{4221225}},
}

\bib{Reiner1975}{book}{
      author={Reiner, I.},
       title={Maximal orders},
   publisher={Academic Press [A subsidiary of Harcourt Brace Jovanovich,
  Publishers], London-New York},
        date={1975},
        note={London Mathematical Society Monographs, No. 5},
      review={\MR{0393100}},
}

\bib{SageMath}{misc}{
      author={{Sage Developers}},
       title={{SageMath}, the sage mathematics software system},
         how={\url{https://www.sagemath.org}},
        note={Version 10.8},
}

\bib{Satake1966}{article}{
      author={Satake, I.},
       title={Clifford algebras and families of abelian varieties},
        date={1966},
        ISSN={0027-7630},
     journal={Nagoya Math. J.},
      volume={27},
       pages={435\ndash 446},
         url={http://projecteuclid.org.ezproxy.uvm.edu/euclid.nmj/1118801764},
      review={\MR{210716}},
}

\bib{Segal2003}{book}{
      author={Segal, Sanford~L.},
       title={Mathematicians under the {N}azis},
   publisher={Princeton University Press},
        date={2003},
}

\bib{Selberg1965}{incollection}{
      author={Selberg, Atle},
       title={On the estimation of {F}ourier coefficients of modular forms},
        date={1965},
   booktitle={Theory of numbers},
      editor={Whiteman, Albert~Leon},
      series={Proc. Sympos. Pure Math.},
      volume={VIII},
   publisher={Amer. Math. Soc., Providence, RI},
       pages={1\ndash 15},
      review={\MR{182610}},
}

\bib{Sengun2012}{incollection}{
      author={{\c{S}}eng{\"u}n, Mehmet~Haluk},
       title={Arithmetic aspects of {B}ianchi groups},
    language={English},
        date={2014},
   booktitle={Computations with modular forms. proceedings of a summer school
  and conference, {H}eidelberg, {G}ermany, {A}ugust--{S}eptember 2011},
   publisher={Cham: Springer},
       pages={279\ndash 315},
}

\bib{Serre1979LocalFields}{book}{
      author={Serre, Jean-Pierre},
       title={Local fields},
      series={Graduate Texts in Mathematics},
   publisher={Springer},
        date={1979},
      volume={67},
}

\bib{Sheydvasser2019}{article}{
      author={Sheydvasser, Arseniy},
       title={Quaternion orders and sphere packings},
        date={2019},
     journal={Journal of Number Theory},
      volume={204},
       pages={41\ndash 98},
}

\bib{Shimura2004}{book}{
      author={Shimura, Goro},
       title={Arithmetic and analytic theories of quadratic forms and
  {C}lifford groups},
      series={Mathematical Surveys and Monographs},
   publisher={American Mathematical Society, Providence, RI},
        date={2004},
      volume={109},
        ISBN={0-8218-3573-4},
         url={https://doi-org.ezproxy.uvm.edu/10.1090/surv/109},
      review={\MR{2027702}},
}

\bib{Solomon2020}{article}{
      author={Solomon, Jake~P.},
      author={Tukachinsky, Sara~B.},
       title={Differential forms on orbifolds with corners},
        date={2020},
     journal={arXiv preprint arXiv:2011.10030},
}

\bib{stacks-project}{misc}{
      author={{Stacks project authors}, The},
       title={The stacks project},
         how={\url{https://stacks.math.columbia.edu}},
        date={2022},
}

\bib{Swan1971}{article}{
      author={Swan, Richard~G.},
       title={Generators and relations for certain special linear groups},
        date={1971},
        ISSN={0001-8708},
     journal={Advances in Math.},
      volume={6},
       pages={1\ndash 77 (1971)},
         url={https://doi-org.ezproxy.uvm.edu/10.1016/0001-8708(71)90027-2},
      review={\MR{284516}},
}

\bib{TauCeti2026}{misc}{
      author={{TauCeti}},
       title={Tau {Ceti}},
         how={\url{https://github.com/tdupu/TauCeti}},
        date={2026},
        note={Revision \texttt{85ed098}},
}

\bib{todd-coxeter}{article}{
      author={Todd, J.~A.},
      author={Coxeter, H. S.~M.},
       title={A practical method for enumerating cosets of a finite abstract
  group},
        date={1936},
     journal={Proc. Edinb. Math. Soc.},
      volume={5},
       pages={26\ndash 34},
}

\bib{Vahlen1902}{article}{
      author={Vahlen, K.~Th.},
       title={Ueber {B}ewegungen und complexe {Z}ahlen},
        date={1902},
     journal={Mathematische Annalen},
      volume={55},
       pages={585\ndash 593},
         url={http://eudml.org/doc/158052},
}

\bib{Vakil2023}{misc}{
      author={Vakil, Ravi},
       title={The rising sea: foundations of algebraic geometry},
  how={\url{http://math.stanford.edu/~vakil/216blog/FOAGjul3123public.pdf}},
        date={2023},
        note={Version of July 31, 2023, accessed December 18, 2023},
}

\bib{vanderGeer1988}{book}{
      author={van~der Geer, Gerard},
       title={Hilbert modular surfaces},
      series={Ergebnisse der Mathematik und ihrer Grenzgebiete},
   publisher={Springer},
        date={1988},
      volume={3.16},
}

\bib{Vinberg1984}{article}{
      author={Vinberg, {\`E}.~B.},
       title={Absence of crystallographic groups of reflections in
  {L}obachevski{\u\i} spaces of large dimension},
        date={1984},
     journal={Trudy Moskov. Mat. Obshch.},
      volume={47},
       pages={68\ndash 102},
        note={In Russian},
      review={\MR{774946}},
}

\bib{Voight2020}{misc}{
      author={Voight, John},
       title={Quaternion algebras},
        date={2020},
        note={Version from May 2020},
}

\bib{Vulakh1993}{article}{
      author={Vulakh, L.~Ya.},
       title={Higher-dimensional analogues of {F}uchsian subgroups of {${\rm
  PSL}(2,\germ o)$}},
        date={1993},
        ISSN={0002-9947},
     journal={Trans. Amer. Math. Soc.},
      volume={337},
      number={2},
       pages={947\ndash 963},
         url={https://doi-org.ezproxy.uvm.edu/10.2307/2154251},
      review={\MR{1145965}},
}

\bib{Vulakh1995}{article}{
      author={Vulakh, L.~Ya.},
       title={Diophantine approximation in {${\bf R}^n$}},
        date={1995},
        ISSN={0002-9947},
     journal={Trans. Amer. Math. Soc.},
      volume={347},
      number={2},
       pages={573\ndash 585},
         url={https://doi-org.ezproxy.uvm.edu/10.2307/2154902},
      review={\MR{1276937}},
}

\bib{Vulakh1999}{article}{
      author={Vulakh, L.~Ya.},
       title={Farey polytopes and continued fractions associated with discrete
  hyperbolic groups},
        date={1999},
        ISSN={0002-9947},
     journal={Trans. Amer. Math. Soc.},
      volume={351},
      number={6},
       pages={2295\ndash 2323},
         url={https://doi-org.ezproxy.uvm.edu/10.1090/S0002-9947-99-02151-0},
      review={\MR{1467477}},
}

\bib{Waterman1993}{article}{
      author={Waterman, P.~L.},
       title={M\"{o}bius transformations in several dimensions},
        date={1993},
        ISSN={0001-8708},
     journal={Adv. Math.},
      volume={101},
      number={1},
       pages={87\ndash 113},
         url={https://doi-org.ezproxy.uvm.edu/10.1006/aima.1993.1043},
      review={\MR{1239454}},
}

\bib{Whitley1990}{thesis}{
      author={Whitley, Elise},
       title={Modular forms and elliptic curves over imaginary quadratic
  fields},
        type={Ph.D. Thesis},
        date={1990},
}

\end{biblist}
\end{bibdiv}

\appendix

\section{Topological Bott Periodicity}\label{A:classical-bott}

There are two types of periodicity which have the name ``Bott Periodicity'': one about periodicity of $K$-theory, which we call Topological Bott Periodicity, and one about periodicity of Clifford algebras. 
The aim of this section is to explain the relationship between these two phenomena for the uninitiated.

To help some make some statements about Clifford algebras more simple, we will use the notation $R(n)=M_n(R),$ for an associative algebra $R$.
Two notions of Bott periodicity are displayed in Table~\ref{tab:bott-periodicity}.
Figure~\ref{fig:bott-clock} arranges the same data as a clock.

\begin{table}[h]
	\begin{center}
		\begin{tabular}{|c|cccccccccc|}
			\hline	$n$ & $1$ & $2$ & $3$ & $4$ & $5$ & $6$ & $7$ & $8$ & $9$\\  \hline
			$\CC_n$ & $\RR$ & $\CC$ & $\HH$ &$\HH\oplus \HH$ &$ \HH(2)$ &$\CC(4)$ &$\RR(8)$ &$\RR(8)\oplus \RR(8)$ &$\RR(16)$\\ \hline
			$X_n^{\infty}$ & - & $\O/\U$ & $\U/\Sp$ &$\BSp\times \ZZ$& $\Sp$ &$\Sp/\U$ &$\U/\O$ &$\BO\times \ZZ$ &$ \O$\\
			\hline
		\end{tabular}
		
	\end{center}
	\caption{A table of Bott periodicity.}\label{tab:bott-periodicity}
\end{table}

A prototypical example of Topological Bott Periodicity is the statement $\KU^n(X) \cong \KU^{n+2}(X),$
where $X$ is a topological space and $\KU$ is the complex $K$-theory pertaining to complex vector bundles. 
A prototype for periodicity of Clifford algebras is the statement 
$ C_{p+1,q+1} =C_{p,q}(2).$
The two statements are connected via the construction of spectra of the $K$-Theory, which are certain topological spaces which represent $K$-Theory. 
These spaces are constructed via consideration of vector spaces with Clifford multiplication and how to prolong them to larger and larger Clifford algebras together with Morita equivalence.

We will focus on real $K$-theory $\KO^n(X)$ for topological spaces $X$ and periodicity of $\CC_n$ in this survey section. 
This section is only intended to give context for non-experts and will not be used elsewhere.
None of its statements is new, and we do not prove them; see \cite{Atiyah1967} for $\KO$-theory and its representability, \cite[Part~IV]{Milnor1963} for Bott periodicity, and \cite{AtiyahBottShapiro1964} and the addendum to \cite{Behrens2002} for the relation with Clifford modules.

Recall that $\KO^0(X)$ is the Grothendieck group of real vector bundles on $X$.
The higher topological $K$-groups are defined by suspension.\footnote{
	The (reduced) suspension of a pointed topological space $X$ is $\Sigma X = S^1 \wedge X$.
	Here $\wedge$ denotes the smash product.
	We should also recall that $\Sigma$ has a right adjoint $\Omega$ which is the based loop space $\operatorname{Hom}_*(S^1,X)$ of base-point-preserving maps, which has the structure of a topological space using the compact open topology (generally smash and hom are adjoint).
	This is the adjunction $\Sigma \dashv \Omega$.}
\begin{definition}\label{T:real-K-theory}
	$KO^n(X) = KO(\Sigma^nX)$
\end{definition}
\begin{theorem}
	The $KO^0$ is represented by $\BO\times \ZZ$ where $\BO = \varinjlim_n \operatorname{Gr}(n,\infty)$ is a limit of real Grassmannians of $n$-dimensional real subspaces inside $\RR^{\infty}$, and $\ZZ$ is given the discrete topology.
	More precisely,
	\begin{equation}
	\KO^0(X) = [X,BO\times \ZZ] \quad (X \text{ compact Hausdorff}),
	\end{equation} 
	where square brackets denotes homotopy classes of morphisms of topological spaces.
\end{theorem}
The space $BO$ is the classifying space for the infinite orthogonal group $O =O_{\infty}(\RR)=\varinjlim \O_n(\RR)$.
Maps to it give real vector bundles.
We note that some authors use the notation $KO_0=BO\times \ZZ$, which is consistent with what follows.

\begin{definition}
	The \emph{real $K$-theory spectrum} is the collection of spaces 
	$(KO_n)_{n\geq 0}=(\Omega^n(BO\times \ZZ))_{n\geq 0}$. 
\end{definition}

\begin{corollary}
	$\KO^n(X) = [X, KO_n].$
\end{corollary}
\begin{proof}
	This follows from the adjunction $\Sigma \dashv \Omega$ and the definition of real $K$-theory $\KO^n(X) = [\Sigma^nX, KO_0] = [X,\Omega^nKO_0] = [X,KO_n]$.
\end{proof}
The classical Bott periodicity in this context is given below.
\begin{theorem}\label{E:real-bott}
	The real Bott periodicity theorem states that $\KO^n(X) = \KO^{n+8}(X)$.
	Equivalently, there is a weak equivalence $\Omega^8KO_n \approx KO_n$.
\end{theorem}

The remainder of the section is devoted to explaining how these statements are related to the periodicity of Clifford algebras. We follow \S2 of the addendum to \cite{Behrens2002}.

The connection between the periodicities of real $\KO^n(X)$ and periodicity of Clifford algebras $\CC_n$ comes from an explicit construction of the spectrum for real $K$-theory by considerations of $\CC_n$-multiplication structures on real inner product spaces. 
The spectrum will be naturally shifted by $8$. 
There are others that fall out of this like $\KSp^{n+4}(X) = \KO^n(X)$. 
Here $\KSp$ is symplectic $K$-theory and $\KSp^0$ is represented by $\BSp\times \ZZ$ in a similar fashion, so the symplectic $K$-theory spectrum and the real $K$-theory spectrum will essentially be the same thing. 

When we prove the generalization of the statement $\CC_{n+8} = M_{16}(\CC_n)$, what is relevant is Morita equivalence.
\begin{definition}
	We recall that two associative algebras $A$ and $B$ are \emph{Morita equivalent} if and only if their categories of left modules $\Mod_A$ and $\Mod_B$ are equivalent. 
\end{definition}
In the application of periodicity we use the fundamental example.
\begin{example}
	Matrices over a ring $R$ are always Morita equivalent to $R$, so the spaces $\CC_n$ and $\CC_{n+8}$ are Morita equivalent.
	The map in this case is given below 
	\begin{equation}\label{E:morita-equivalence}
	\Mod_{\CC_n} \xrightarrow{\sim} \Mod_{M_{16}(\CC_n)} \cong \Mod_{\CC_{n+8}}, \quad W \mapsto W^{\oplus 16}.
	\end{equation}
\end{example}

We will now be dealing with spaces with $\CC_n$-multiplication.
\begin{definition}
	Let $W$ be a real inner product space of finite dimension.
	A \emph{$\CC_n$-multiplication} is a collection $J_1,J_2,\ldots,J_{n-1}$ of isometries of $W$ such that the rule 
	$$ i_a \cdot w = J_a(w), \quad 0<a<n, $$
	gives $W$ the structure of a $\CC_n$-module. 
	We will denote such objects by $(W, (J_1,J_2,\ldots,J_{n-1}))$, and when the multiplications are clear simply by $W$.
	We denote the category of such objects as $\Hilb^{\CC_n}_{\RR}$.
\end{definition}

Note that given a inner product space $W$, giving it a $\CC_1$-multiplication does nothing, giving it a $\CC_2$-structure makes it a complex vector space (where $J_1$ is an almost complex structure) with inner product, giving it a $\CC_3$-structure makes it a quaternionic vector space (a left $\HH$-module) with compatible inner product, and so on. The act of extending these structures means choosing $J_1,J_2,J_3,\ldots$ in a sequence so that $J_n$ will be compatible with the previous $J_1,\ldots, J_{n-1}$ in a way that it extends the $\CC_n$-module structure to a $\CC_{n+1}$-module structure.

Let $W\in \Hilb_{\RR}^{\CC_n}$ with $J_1,\ldots,J_{n-1}$ giving the $\CC_n$-structure.
We can define now a new set of isometries in the orthogonal group of the inner product space $W$ to be the collection of isometries $f:W\to W$, which are equivariant with respect to $J_1,\ldots,J_{n-1}$. 
We could call this group $\O( (W,(J_1,\ldots,J_{n-1})))$. 
To keep things short we will keep $W$ fixed for now and let 
$$ G_n =  \O( (W,(J_1,\ldots,J_{n-1})))$$
be the group of isometries preserving the $\CC_n$-multiplication.

As we increase $n$ the conditions on these groups become more restrictive, so $G_n \supset G_{n+1}$. 

Fix $(W,(J_1,\ldots,J_{n-1}))\in\Hilb_{\RR}^{\CC_n}$.
The choices of $J_n$ that extend the $\CC_n$-multiplication to a
$\CC_{n+1}$-multiplication form a space $X_{n+1}(W)$, abbreviated $X_{n+1}$.

In notation
$$ X_{n+1} = \lbrace J_n \colon \mbox{$J_n$  prolongs  $J_1,\ldots,J_{n-1}$ }\rbrace.$$
Thus $X_{n+1}$ is the space of $\CC_{n+1}$-multiplications extending a fixed $\CC_n$-multiplication; for example $X_2$ is the space of orthogonal complex structures on $W$.
We care about a limiting version of these spaces, which we will call
$X_n^{\infty}$. 
These are obtained by taking a colimit over the spaces $W$ where the colimit  $W_{\infty}$ has the property that every irreducible object of $\Hilb^{\CC_n}_{\RR}$ appears as a direct summand of $W_{\infty}$ infinitely many times.

\begin{theorem}\label{T:spectrum}
	The collection of spaces $(X_n^{\infty})_{n\geq 2}$ form a spectrum:
	$\Omega X_n^{\infty} \approx X_{n+1}^{\infty}$.
Moreover, we have $X_4^{\infty} = \BSp\times \ZZ$ and $X_8^{\infty} = \BO\times \ZZ$, which means for every compact Hausdorff space $Y$ we have 
	$$\KO^n(Y) = [Y,X_{n+8}^{\infty}], \quad \KSp^n(Y) =[Y,X_{n+4}^{\infty}].$$
\end{theorem}

We now want to explain the above representations of real and symplectic $K$-theory and why they give rise to Bott periodicity (this, as perhaps expected, has to do with Morita equivalence). 
\begin{lemma}
	Suppose $W$ admits at least one $J_n$ prolonging $J_1,\ldots, J_{n-1}$ (and hence defining $G_{n+1}$).
	Then the map
	\begin{equation}
	G_n/G_{n+1} \to X_{n+1}, \quad g \mapsto g J_n g^{-1},
	\end{equation}
	is a homeomorphism onto the $G_n$-orbit of $J_n$. This orbit is all of $X_{n+1}$ when $\CC_{n+1}$ is simple; when $\CC_{n+1}\cong A\oplus A$ (that is, $n+1\equiv 0 \bmod 4$), $X_{n+1}$ is a disjoint union of such orbits.
\end{lemma}

We now sketch the proof of Theorem~\ref{T:spectrum} and how Bott periodicity follows. 
If we return to the finite dimension spaces $X_n=X_n(W)$, we see that they fit into a fibration
$$ X_{n+1} \to E_n \to X_n,$$
where $E_n$ is a total space.
We have very good control of $G_n/G_{n+1}$ (we know exactly what these are) and all of this fits together well in a limit giving space $X^{\infty}_n$, which we understand, and a fibration
\begin{equation}\label{E:fibration}
X_{n+1}^{\infty} \to E_n^{\infty} \to X_n^{\infty}.
\end{equation}
In this, $E_n^{\infty}$ is contractible (see \cite[addendum Theorem 3.1]{Behrens2002}). 
Rotating this sequence gives 
$$ \Omega E_n^{\infty} \to \Omega X_n^{\infty} \to X_{n+1}^{\infty} \to E_n^{\infty} $$
and the two end terms are contractible giving a weak equivalence 
$\Omega X_n^{\infty} \approx X_{n+1}^{\infty}$.

At the finite level the Morita equivalence \eqref{E:morita-equivalence} will induce things like 
$$ X_n(W) \cong X_{n+8}(W^{\oplus 16}),$$
which comes from Clifford periodicity.
The increase in dimension washes away in the limit to give $X_n^{\infty} \approx X_{n+8}^{\infty}$, which is what proves Theorem~\ref{E:real-bott}.
It is also an explicit description of the groups and quotients $G_n(W)/G_{n+1}(W)$ that allows us to prove $X_4^{\infty} = \BSp\times \ZZ$ and $X_8^{\infty} = \BO\times \ZZ$ and relate the periodicity induced by Morita equivalence to $K$-theory.

\section{Group Schemes}
\label{sec:group-schemes}
Let $G$ be a group scheme over $\QQ$.
The goal of this appendix is to explain why the group $G(\QQ)$ and the group scheme $G$ are not the same thing and why $G(\ZZ)$ is not well-defined. 

	We remind the reader that group schemes over a commutative ring $R$ are group objects in the category of schemes over $R$, which will we denote by $\GrpSch_R$ \cite[Section 7.6.4]{Vakil2023}.
	Being a group object means that morphisms to this object in the category naturally (functorially) have the structure of a group---the hom functor with the object plugged into the right entry is a functor to not just sets, but groups.
	In our case this means these are schemes equipped with morphisms $\mu:G\times G\to G$, $S:G\to G$, $e: \Spec(R) \to G$ such that for every $R$-algebra $R'$ its points $G(R')$ have the structure of a group in a functorial way, i.e., $G(R')$ is an object in the category of groups. If $R'\to R''$ is a morphism of $R$-algebras then we have a group homomorphism $G(R') \to G(R'')$. 
	Note that $G(R)$ is shorthand for the set of morphisms $\Spec(R) \to G$.
	
	With the exception of abelian varieties, all group schemes in the paper are affine, which means they are of the form $G= \Spec(A)$ for some $R$-algebra $A$, which has the structure of a Hopf algebra (the Hopf algebra structure is the encoding of the group operations in the language of rings so all the axioms are ``opposite'' of ring axioms---so, for example, a multiplication map $G\times G\to G$ turns into a comultiplication map $A \to A\otimes A$). 
	The set $G(R')$ is in bijection with the set of $R$-algebra homomorphisms $A\to R'$; we usually think of $G(R')$ as entries in some matrix. 
	When $G=\Spec(A) \in \GrpSch_R$ and $R'$ is an $R$-algebra, we let $G_{R'} = \Spec(A \otimes_R R') \in \GrpSch_{R'}$. 
	So for example, $\GG_m = \GG_{m,\ZZ} = \Spec( \ZZ[x,x^{-1}])$ and $\GG_{m,\QQ} = \Spec( \QQ[x,x^{-1}])$ and $\GL_2 = \Spec \ZZ[x_{11},x_{12},x_{21},x_{22},y]/\langle (x_{11}x_{22}-x_{12}x_{21})y-1\rangle$.

        There are many examples to show that $G(\QQ)$ does not determine $G$.  For example, suppose that $E, E'$ are two elliptic curves over $\QQ$ whose Mordell-Weil groups are isomorphic but that are not isogenous, and that $n$ is an integer such that the quadratic twists of $E, E'$ by $n$ have different ranks.  Then $E, E'$ define group schemes over $\QQ$ with $E(\QQ) \cong E'(\QQ)$, but $E(\QQ(\sqrt{n}))$ and $E'(\QQ(\sqrt{n}))$ have different ranks as abelian groups, so they are certainly not isomorphic.
        
	The following provides an example where the set of $\ZZ$-points of a $\QQ$-group scheme depends on the choice of $\ZZ$-group scheme and therefore fails to be well-defined.
	
\begin{example} Fix $n \in \ZZ$.  The functor on $\ZZ$-algebras taking
		$R$ to the multiplicative group of elements of norm $1$ in $R[t]/(t^2-n)$
		is a group scheme $S_n$ represented by $\Spec(\ZZ[x,y]/(y^2 - nx^2 - 1))$.
		Now let $a > 1$ and consider $S_n, S_{na^2}$. There is a map
		$S_{na^2}(\ZZ) \to S_n(\ZZ)$ induced by the algebra map
		taking $x, y$ to $ax, y$; this becomes an isomorphism over $\QQ$, but it
		is not over $\ZZ$ in general.
		If $n = 2, a = 3$, for example, then $(2,3) \in S_2(\ZZ)$ is
		not the image of any point of $S_{18}(\ZZ)$. Even if the $\ZZ$-points happen
		to be the same, the functors will still be different. For example,
		with $n = 2, a = 2$ we have $S_2(\ZZ) = S_8(\ZZ)$ as sets, because
every unit of $\ZZ[\sqrt 2]$ of norm $1$ is of the form $\pm(3 + 2 \sqrt 2)^k$ and hence belongs to $\ZZ[\sqrt 8]$. This is not at all true for general
		rings, however; in the ring $\ZZ[\sqrt 3,\sqrt 2]$ the element
		$\sqrt 3 - \sqrt 2$ is a unit, and it does not belong to
		$\ZZ[\sqrt 3,\sqrt 8]$, so $S_2(\ZZ[\sqrt 3]) \ne S_8(\ZZ[\sqrt 3])$.
\end{example}

\section{AI Methods}\label{app:ai-methods}

\subsection{Typesetting}
Agents using Claude Code with \texttt{claude-opus-5-5} edited the \LaTeX{}
source and prepared the figures in TikZ. Of the $34$ agent-drawn figures in
this version, $13$ redraw earlier figures and $21$ are new. They illustrate
the inclusions among orders, Bott periodicity, codes and lattices, M\"obius
maps, horoballs, cusp classes, and fundamental-domain constructions. The
redrawings include domain boundaries, horocycles, and Coxeter diagrams.
Drawing them exposed errors in a root sign, a face label, an angle range,
and the description of geodesics meeting horocycles.

\subsection{Agent management}
The workflow instructions have a versioned public source in the
\texttt{agent-skills-public} repository~\cite[\href{https://github.com/tdupu/agent-skills-public/blob/8d13efb5e66306aeae6e3b948bc463e7e2db1b65/README.md}{README}]{AgentSkillsPublic}.
The repository AI policy was established with
\agentskill{new-repo-ai-policy}. Task provenance was recorded using
\agentskill{update-ai-usage}. Token counts and gaps in those counts were
recorded using \agentskill{update-tokens}.

\subsection{Natural language proofs}
The \agentskill{frontier-dump} workflow produced repair plans for findings
that required new arguments. \agentskill{triage-referee-report} linked
findings to proposed source changes, and \agentskill{revise} was used to
apply and check revisions. Agents in Claude Code with
\texttt{claude-opus-5-5} and in OpenAI Codex (model not recorded) also
drafted proof text while porting arguments to Lean. The authors reviewed
the text and revised the arguments included here.

\subsection{Adversarial reviews}
Two whole-manuscript reviews used \agentskill{referee-report} in Claude Code
with \texttt{claude-opus-5-5}. Each assigned nine agents to ranges of the paper.
They checked proofs against their statements, tested small cases for
counterexamples, compared the introduction with the body, and checked
citations. In the first review, another agent attempted to refute reported
errors. Supporting \agentskill{grill-with-docs}, \agentskill{grilling},
and \agentskill{domain-modeling} workflows helped examine disputed steps.
A later proofreading pass examined notation, missing cases, and agreement
between statements and proofs. The authors reviewed the resulting changes;
the agent reviews themselves are not human verification.

\subsection{AI assisted computations}
Review agents wrote exact-rational scripts in
Python~\cite[3.13.7]{Python2026}. They used GAP~\cite[4.14.0]{GAP4}
for finite group calculations. SageMath~\cite[10.8]{SageMath} supplied linear programming and
polyhedral calculations. The \texttt{fpylll} library~\cite{Fpylll}
(version not recorded) was used to enumerate short lattice vectors. Agents
also used Magma's \texttt{MaximalOrder} command~\cite{Magma1997} to test
selected orders. Monte Carlo orbit counts in the review records are
approximate rather than exact.

An agent reran an exact certificate in Magma~\cite[V2.29-2]{Magma1997},
using a development version of the authors' Magma package. A public version
is available in \texttt{magma-clifford-algebras}~\cite[\href{https://github.com/tdupu/magma-clifford-algebras/tree/51cfdeb0135c19d2f660ca985a8865bf91d3f65d}{repository, commit 51cfdeb}]{MagmaCliffordAlgebras2026};
that public revision does not yet contain the certificate code. The certificate passed for a presentation
of $\SL_2(\Ocal_4)$ with $10$ generators and $75$ relators, derived from a
fundamental domain different from the one used in this manuscript. Modulo
$-1$, that presentation and the one in this manuscript have the same
abelianization and the same numbers of subgroups at each index up to $10$.
These comparisons do not establish that the presented groups are
isomorphic. The Magma computations in the mathematical body were carried
out by the human authors.

\subsection{AI statement generation}
This discussion, the abstract notice, and the appendix were drafted and
revised in Claude Code with \texttt{claude-opus-5-5}, using the
\agentskill{write-materials-and-methods} workflow. A separate agent
checked the disclosure with \agentskill{latex-ai-statement}; its factual
findings were used in revisions. The authors reviewed this account and
take responsibility for its final wording.

\subsection{Proofreading and adversarial review}
The \agentskill{referee-report} findings and the Lean port prompted
mathematical corrections to the Clifford $\SL_2$ and Spin statements,
the classification of maximal orders, arithmeticity, fundamental domains,
cusp classes, and group presentations. Smaller corrections addressed signs
and labels in examples and figures, notation, citations, cross-references,
and explanatory wording. The authors reviewed and revised the changes
included in this version. The review reports and test scripts are retained
in the repository.

\subsection{Literature review}
Agents suggested some bibliography entries while checking claims and
references. The provenance records distinguish these suggestions from
references selected by the authors. The authors reviewed the bibliography
for this version.

\subsection{Formalization}
Agents used Lean~4~\cite{Lean4}, Mathlib~\cite{Mathlib2020}, and
Tau~Ceti~\cite{TauCeti2026}. The formalization source is in the public
\texttt{cliff-lean} repository~\cite[\href{https://github.com/tdupu/cliff-lean/tree/19c129c40d331d53695b2d12bb0aa6806a0a9b44}{repository, commit 19c129c}]{CliffLean2026}.
The \agentskill{using-leanpowers} route used
the \agentskill{lean-workflow} pipeline: \agentskill{lean-align-statement},
\agentskill{lean-prove}, and \agentskill{lean-verify}. Statement alignment
compared each Lean target with the corresponding manuscript claim.

The inventory for the automated revision contains $314$ theorem-type
environments, including remarks. All $209$ items selected for the Lean port
have checked conditional proofs; $105$ items, mostly expository or
computational, were outside its scope. The proofs use only Lean's standard
axioms and explicit hypotheses for background results that we expect to
formalize in future work. The manuscript-facing dependencies are the
commensurability of images of arithmetic groups under rational
morphisms~\cite[Prop.~9.4.5]{Morris2015}
(cf.~\cite[Thm.~4.1]{Platonov1994}), the
Borel--Harish-Chandra finite-covolume theorem~\cite[Thm.~4.14]{Platonov1994},
results on reduced norms, normal ideals, Eichler conjugacy, and local
indices~\cite[\S\S~9, 22--24, 31, 34--35]{Reiner1975}, the
local-to-global construction of maximal orders~\cite[Thm.~5.3 and
Cor.~11.6]{Reiner1975}, local-field results on unramified quadratic
extensions~\cite[Thm.~12.8]{Reiner1975} and~\cite{Serre1979LocalFields},
and the
identification of fppf $\mathbb{G}_m$-torsors with the Picard
group~\href{https://stacks.math.columbia.edu/tag/03P8}{(Stacks Project,
Tag~03P8)}. Auxiliary Lean modules also assume the homogeneous-space
manifold and orbit theorems for Lie group actions~\cite[Ch.~21]{Lee2013SmoothManifolds},
but the manuscript-facing statements do not depend on them.
The written proofs in this paper do not depend on the Lean port.

\subsection{Record keeping and token accounting}\label{app:tokens}
The repository retains provenance and token records for AI-assisted
tasks. Token counts were not recorded for every task, and no costs were
calculated. The available counts are:

\begin{center}
\small\setlength{\tabcolsep}{4pt}\begin{tabular}{llll}
\hline
Task & Model / harness & Tokens & Cost (USD)\\
\hline
Satake work & not recorded / Codex & not recorded & unpriced\\
Initial Lean pass & not recorded & not recorded & unpriced\\
Source editing & Claude Opus 5.5 / not recorded & not recorded & unpriced\\
Adversarial review & \texttt{claude-opus-5-5} / Claude Code & 2{,}287{,}174$^{*}$ & unpriced\\
Review triage & \texttt{claude-opus-5-5} / Claude Code & not recorded & unpriced\\
Statement drafting & \texttt{claude-opus-5-5} / Claude Code & not recorded & unpriced\\
Lean formalization & \texttt{claude-opus-5-5} / Claude Code & 1{,}871{,}637$^{*}$ & unpriced\\
\hline
\end{tabular}
\end{center}
\noindent$^{*}$Each starred count covers the recorded subagents for that
task, but excludes the coordinating session and web searches. It has no
breakdown by token class. The cost is unpriced because no rate was applied.

\end{document}